\documentclass{amsart}
\usepackage{amssymb,amsmath}
\usepackage{graphicx}
\usepackage{enumitem,array}
\usepackage{cite}
\usepackage{xcolor}
\usepackage{caption}

\usepackage{float}
\usepackage{tikz}
\usetikzlibrary{arrows.meta}
\usetikzlibrary{patterns}
\usepackage{geometry}
\allowdisplaybreaks
\usepackage[expansion=false]{microtype}
\makeatletter
\def\thm@space@setup{\thm@preskip=3pt \thm@postskip=3pt }
\def\subsection{\@startsection{subsection}{2}%
\z@{.25\linespacing\@plus.35\linespacing}{-.5em}%
{\normalfont\bfseries}}
\def\subsubsection{\@startsection{subsubsection}{3}%
\z@{.25\linespacing\@plus.35\linespacing}{-.5em}%
{\normalfont\itshape}}
\makeatother
\setlist{topsep=2pt,itemsep=2pt,parsep=0pt,partopsep=0pt}
\setlist[enumerate,1]{label={\rm(\roman*)},ref={\roman*}}
\setlist[enumerate,2]{label={\rm(\alph*)},ref={\alph*}}
\setlist[enumerate,3]{label={\rm(\arabic*)},ref={\arabic*}}
\selectfont
\usepackage[title,titletoc]{appendix}

\usepackage[colorlinks=true,plainpages=true,pdfpagelabels,hypertexnames=true,pdfstartview={FitV null},linkcolor=magenta,citecolor=red,urlcolor=black]{hyperref}

\newcommand{\boundarydataarrows}{%
\foreach \X/\Y in {
-0.38/0.45,
0.30/0.45,
-0.65/0.00,
0.05/0.00,
0.65/0.00,
-0.38/-0.45,
0.30/-0.45
}{
\draw[
->,
line width=1.00pt,
blue!45!black
]
({\X-0.23},\Y)--++(0.46,0);
}%
}

\newtheorem{defn}{Definition}[section]
\newtheorem{definition}[defn]{Definition}
\newtheorem{lemma}[defn]{Lemma}
\newtheorem{proposition}[defn]{Proposition}
\newtheorem{theorem}[defn]{Theorem}
\newtheorem{remark}[defn]{Remark}
\numberwithin{equation}{section}

\newcommand{\ep}{\varepsilon}
\newcommand{\R}{\mathbb{R}}
\newcommand{\norm}[1]{\left\Vert#1\right\Vert}
\DeclareMathOperator{\dist}{dist}
\DeclareMathOperator{\supp}{supp}
\DeclareMathOperator{\sgn}{sgn}

\begin{document}
\raggedbottom

\setlength{\abovedisplayskip}{2pt plus 1pt minus 1pt}
\setlength{\belowdisplayskip}{2pt plus 1pt minus 1pt}
\setlength{\abovedisplayshortskip}{0pt plus 1pt}
\setlength{\belowdisplayshortskip}{1pt plus 1pt minus 1pt}
\title[Separation points and pressure gradients of the Stokes system]{Separation points and pressure gradients for the nonstationary Stokes system with localized boundary data in the half-space}

\author{Kyungkeun Kang}
\author{Chanhong Min}
\thanks{Affiliation of both authors: Department of Mathematics, Yonsei University, 50 Yonsei-ro, Seodaemun-gu, Seoul 03722, Republic of Korea. Email: kkang@yonsei.ac.kr (K. Kang), chopinmin1@yonsei.ac.kr (C. Min).}
\thanks{Corresponding author: Chanhong Min (chopinmin1@yonsei.ac.kr).}
\date{}

\makeatletter
\@namedef{subjclassname@2020}{%
	\textup{2020} Mathematics Subject Classification}
\makeatother

\begin{abstract}
We study separation points of the nonstationary Stokes system in a half-space with zero initial data and nonnegative, space-time separable boundary data supported in the unit ball and in $[0,1]$, respectively. We call $(x',t)\in\mathbb{R}^{n-1}\times(0,\infty)$, with $x'$ outside the spatial support, a separation point if one velocity component takes both signs arbitrarily near $(x',0,t)$. The main ingredient of our analysis is an explicit representation of velocity and pressure involving confluent hypergeometric functions.

For normal boundary data that are unimodal and log-concave in time, the $i$th tangential velocity component has a unique separation time where $|x_i|>1$. For sufficiently large $|x'|$, the normal velocity component has two separation times, one before and one after $t=1$. Along the $x_i$-axis, the separation time for the tangential velocity component eventually increases toward the unique zero of the temporal data's Caputo half-derivative. Both separation times for the normal velocity component eventually increase radially, the earlier approaching the zero of the Caputo half-derivative and the later growing quadratically in $|x'|$.

For tangential boundary data that are unimodal in time, with $n\geq3$, a separation point of the parallel velocity component exists at $|x'|>1$ precisely when the corresponding Riesz transform's derivative in the data direction is positive, with uniqueness under log-concavity. For the remaining velocity components, we establish unique separation points on axis-parallel lines far from the origin and derive their asymptotics.

For radial spatial data, we determine the signs and zero sets of tangential pressure derivatives at the boundary, thereby locating regions of favorable and adverse pressure gradients.
\end{abstract}
\maketitle

\section{Introduction}

In this paper, we study the nonstationary Stokes system in the half-space with zero initial data and nonzero boundary data.
\begin{align}\label{0308-01}
	\begin{cases}
		\partial_{t}w-\Delta w+\nabla p=0,\,\,\,\nabla\cdot w=0, & \quad \text{in } \R_+^n \times (0,\infty), \\
		w(x',0,t)=g(x',t), & \quad \text{in } \Sigma\times (0,\infty),\\
		w(x,0)=0, & \quad \text{in } \R_+^n,
	\end{cases}
\end{align}
where we assume \(n\geq2\) and write \(\R^n_+=\{x=(x',x_n)\in\R^n:x_n>0\}\) and \(\Sigma=\R^{n-1}\).
Here \(w: \R_+^n \times (0,\infty)\rightarrow \R^n\) denotes the velocity field, \(p:\R_+^n\times(0,\infty)\rightarrow \R\) denotes the corresponding pressure, and \(g:\Sigma\times(0,\infty)\rightarrow \R^n\) denotes the boundary data prescribed on the plane \(\{x_n=0\}\).
The solution formula for \eqref{0308-01} with \(g=(g_1,\ldots,g_n)\) is given by \cite{Sol231,KangLaiLaiTsai2023}
\begin{align}\label{0330-01}
	w_i(x,t)=\sum_{j=1}^n\displaystyle\int_{0}^t\int_{\Sigma}K_{ij}(x-y',t-s)g_j(y',s)\,dy'ds,\quad 1\leq i\leq n,
\end{align}
where \(x-y':=(x'-y', x_n)\) and \(K_{ij}\) is the Golovkin tensor given by
\begin{align}\label{0308-02}
	K_{ij}(x,t)= -2\delta_{ij}\partial_{x_n}\Gamma(x,t) + 2\delta_{jn}\delta(t) \partial_{x_i}N(x) + 4L_{ij}(x,t),\quad 1\leq i,j\leq n.
\end{align}
Here \(\delta_{ij}\) is the Kronecker delta, and \(\delta\) is the Dirac mass at time zero.
The heat kernel \(\Gamma\) and the fundamental solution \(N\) of the Laplace equation are given by
\begin{align*}
	\Gamma(x,t) & =
	\begin{cases}
		\displaystyle \frac{1}{(4\pi t)^{n/2}} e^{-\frac{|x|^2}{4t}}, & \quad t>0,\\[6pt]
		0,& \quad t\leq 0,
	\end{cases}
	\qquad N(x) =
	\begin{cases}
		\displaystyle -\frac{1}{(n-2)\omega_{n-1}|x|^{n-2}}, & \quad n\geq 3, \\[6pt]
		\displaystyle \frac{1}{2\pi}\ln |x|, & \quad n=2,
	\end{cases}
\end{align*}
where \(\omega_{n-1} = \frac{2\pi^{\frac{n}{2}}}{\Gamma\left(\frac{n}{2}\right)}\) is the surface area of the unit sphere \(\mathbb S^{n-1}\subset\R^n\), and \(\Gamma\) denotes the classical Gamma function.
With this sign convention, \(\Delta N=\delta_0\) in the sense of distributions.
The kernel \(L_{ij}\) is given by
\begin{align}\label{0419-02}
	\displaystyle L_{ij}(x,t) :=\partial_{x_j}\int_{\Sigma\times(0,x_n)}\partial_{z_n}\Gamma(z,t)\partial_{x_i}N(x-z)\,dz.
\end{align}

Throughout this paper we factor the heat kernel as \(\Gamma(x,t)=\Gamma'(x',t)\Gamma_1(x_n,t)\), where
\begin{align*}
	\Gamma'(x',t):=(4\pi t)^{-\frac{n-1}{2}}e^{-\frac{|x'|^2}{4t}}, \qquad \Gamma_1(x_n,t):=(4\pi t)^{-\frac12}e^{-\frac{x_n^2}{4t}},\quad t>0,
\end{align*}
so that \(\Gamma'\) is the \((n-1)\)-dimensional heat kernel.

For \(1\leq i\leq n-1\), the \(i\)th \((n-1)\)-dimensional Riesz transform \(R_i'\) is defined by
\begin{align}\label{0419-021}
	R_i'f(x')=2\,\mathrm{p.v.} \int_{\Sigma}\partial_iN(x'-y',0)f(y')\,dy' = \frac{\Gamma\left( \frac{n}{2}\right)}{\pi^{n/2}}\,\mathrm{p.v.}\int_{\Sigma}\frac{x_i-y_i}{|x'-y'|^n}f(y')\,dy',
\end{align}
where \(\mathrm{p.v.}\)
denotes the Cauchy principal value.

Throughout this paper we consider boundary data that are localized in space and time and point in a single fixed coordinate direction.
That is, we fix \(1\leq j\leq n\) and take
\begin{align}\label{0324-13}
	g(x',t)=g_j(x',t)\mathbf e_j,\qquad g_j(x',t)=\psi(x')\phi(t),
\end{align}
where \(\psi:\Sigma\to\R\) is nonnegative and supported in \(\overline{B_1'}=\{x'\in \Sigma : |x'|\leq 1\}\), and \(\phi\) is nonnegative and supported in \([0,1]\).
The precise assumptions on \(\psi\) and \(\phi\) are {\rm\textbf{(S)}} and {\rm\textbf{(T1)}}--{\rm\textbf{(T5)}} in Subsection \ref{260002} below.

\subsection{Known results and motivation}

Boundary-layer separation is an important way in which viscosity changes a flow near a solid boundary.
Prandtl's boundary layer theory \cite{Prandtl1904} describes the thin viscous region next to the boundary.
In two dimensions, with \(x\) the coordinate along the wall and \(y>0\) the distance from it, the velocity \((u,v)\) in this region satisfies, after the standard boundary-layer rescaling, the Prandtl equations
\begin{equation*}
	\begin{cases}
		\partial_tu+u\partial_xu+v\partial_yu-\partial_y^2u=-\partial_xp^E, & y>0,\\
		\partial_xu+\partial_yv=0, & y>0,\\
		u=v=0, & y=0,\\
		u\to u^E, & y\to\infty,
	\end{cases}
\end{equation*}
where \(u^E(t,x)\) and \(p^E(t,x)\) are the tangential velocity and the pressure of the outer Euler flow at the wall, so that \(\partial_tu^E+u^E\partial_xu^E=-\partial_xp^E\).
The steady equations are obtained by dropping \(\partial_tu\).

In the classical steady two-dimensional setting, separation is associated with the vanishing of the wall shear stress and the onset of reversed flow near the wall.
Goldstein \cite{Goldstein1948} analyzed the singular behavior of steady boundary-layer solutions near separation, while Lighthill \cite{Lighthill1963} emphasized the role of vorticity production at the wall.
Dalibard and Masmoudi \cite{DalibardMasmoudi2019} proved separation for the stationary Prandtl equation under an adverse pressure gradient and justified the Goldstein singularity for a suitable class of boundary data.
Further results on separation and the local behavior of steady Prandtl solutions were obtained by Shen, Wang, and Zhang \cite{ShenWangZhang2021}.
Iyer and Masmoudi \cite{IyerMasmoudiReversal} constructed stationary Prandtl flows exhibiting reversal and recirculation by treating the change of sign of the tangential velocity through a mixed-type free-boundary formulation.

For unsteady or higher-dimensional flows, vanishing wall shear alone does not characterize separation, and additional geometric and kinematic information is needed \cite{PerryFairlie1975,SearsTelionis1975,TobakPeake1982}.
The singularity observed numerically by van Dommelen and Shen \cite{VanDommelenShen1980} was studied through a Lagrangian description by van Dommelen and Cowley \cite{VanDommelenCowley1990}.
Kukavica, Vicol, and Wang \cite{KukavicaVicolWang2017} rigorously established finite-time blowup of the boundary-layer thickness for the unsteady Prandtl equations.
For the inviscid two-dimensional Prandtl system, Collot, Ghoul, and Masmoudi \cite{CollotGhoulMasmoudi2021} analyzed singularity formation and unsteady separation, including self-similar solutions and the generic appearance of a particular separation profile.
For the viscous unsteady Prandtl system, Collot, Ghoul, Ibrahim, and Masmoudi \cite{CollotGhoulIbrahimMasmoudi2022} identified a stable blowup pattern for a reduced one-dimensional equation governing a tangential derivative of the velocity along a transversal axis.
They also studied the persistence of analyticity near the axis when this derivative, restricted to the axis, follows the same blowup pattern.

The well-posedness of the Prandtl equations has been established in several distinct settings.
The classical monotonicity theory is developed in \cite{OleinikSamokhin1999}.
Under Oleinik's monotonicity assumption, Masmoudi and Wong \cite{MasmoudiWong2015} proved local existence and uniqueness in weighted Sobolev spaces by energy estimates exploiting a nonlinear cancellation.
Sammartino and Caflisch \cite{SammartinoCaflisch1998a,SammartinoCaflisch1998b} established analytic solutions and justified the associated inviscid limit.
Beyond analyticity, G\'erard-Varet and Masmoudi \cite{GerardVaretMasmoudi2015} proved local well-posedness in tangential Gevrey class \(7/4\) for a class of nonmonotone initial data.
G\'erard-Varet, Maekawa, and Masmoudi \cite{GerardVaretMaekawaMasmoudi2018} established Gevrey stability of Prandtl expansions for two-dimensional Navier--Stokes flows under monotonicity and concavity assumptions on the boundary-layer profile.
Related criteria and results for the inviscid limit were obtained in \cite{Kato1984,Kelliher2007,Kelliher2008,Maekawa2014}.

Steady Prandtl expansions have also been justified for various outer flows and boundaries, with corresponding regularity and stability results \cite{GuoIyer2021,GuoIyer2023,GuoNguyen2017,Iyer2017RotatingDisk,Iyer2019MovingI,Iyer2019MovingII,Iyer2019NonShear,Iyer2020Blasius,Iyer2020MovingIII,IyerMasmoudi2021,IyerMasmoudi2026}.
Complementary results include finite-time blowup \cite{EEngquist1997}, ill-posedness \cite{GerardVaretDormy2010,GuoNguyen2011}, and instability of Navier--Stokes boundary layers and their sublayers \cite{BianGrenier2024,Grenier2000,GrenierGuoNguyen2016a,GrenierGuoNguyen2016b,GrenierNguyen2018,GrenierNguyen2019,GrenierNguyen2024}.
These developments motivate a quantitative study of how boundary data, viscous diffusion, and nonlocal pressure determine separation and reversal.
The present paper investigates these mechanisms within the nonstationary Stokes system by determining separation points, their temporal and spatial behavior, and the signs of the associated pressure gradients.

A closely related line of research concerns boundary-induced
singularities arising from localized boundary data for the
Stokes and Navier--Stokes systems in the half-space.
Kang \cite{Kang2005} showed that compactly supported,
H\"older-continuous boundary flux can produce unbounded
normal derivatives for the Stokes system, even when the
solution remains bounded and has finite energy.
More precisely, he constructed a weak solution $u$ satisfying
\[
\|\partial_{x_n} u\|_{L^\infty(Q_{1/2}^{+})} = \infty,
\qquad
\|u\|_{L^\infty(Q_1^{+})}
+
\|\nabla u\|_{L^2(Q_1^{+})} < \infty,
\]
where \(B_r^+:=\{x\in\R^n_+:|x|<r\}\) and \(Q_r^+:=B_r^+\times(1-r^2,1)\) for \(0<r\leq1\).
Subsequently, Kang, Lai, Lai, and Tsai
\cite{KangLaiLaiTsai2022} demonstrated that analogous
boundary-gradient singularities can also occur for
finite-energy solutions of the Navier--Stokes equations.
Their analysis uses sharp pointwise bounds for the Green tensor of the nonstationary Stokes system \cite{KangLaiLaiTsai2023},
and such pointwise estimates were later used to obtain \(L^1\)--\(L^q\) estimates and mild solutions with pointwise decay \cite{KangLaiLaiTsai2025}.
Chang and Kang studied local regularity near the boundary.
In \cite{ChangKang2023}, they constructed bounded solutions whose normal derivatives are singular in every local \(L^p\), \(1<p<\infty\), and established criteria for boundedness and H\"older regularity of gradients.
They studied further singular behavior near the boundary in \cite{ChangKang2025JMFM,ChangKang2026NA}.
Kang and Min \cite{KM2024} showed that, near the boundary and away from the support of tangential boundary data, the regularity of the solution depends on the temporal regularity of the boundary data.
They also showed that normal and tangential boundary data produce boundary singularities of different orders, namely unbounded first and second normal derivatives, respectively.

The works \cite{CKM2026,CK2026} are the most directly related to the present paper.
Chang, Kang, and Min \cite{CKM2026} constructed space-time localized boundary data that produce flow reversal in the Stokes system \eqref{0308-01}.
They showed that both the tangential velocity components and the normal component change sign at least once.
For normal boundary data, Chang and Kang \cite{CK2026} related separation of a tangential velocity component to a temporal singular integral \(M\) proportional to the Caputo half-derivative of the temporal boundary data.
In their setting for the temporal boundary data, if \(M<0\) on an interval \((a,b)\subset(1,2)\), separation occurs at a time in \((1,a)\) at each sufficiently distant point of a suitable cone \cite[Assumption~1.4 and Theorem~1.7(i)]{CK2026}.
For \(\phi(1)>0\), they also gave a sufficient condition for \(M(t)>2^{-1/2}\phi(1)\) throughout \(1<t<2\) \cite[Lemma~6.1(iii)]{CK2026}.
They analyzed the motion of the separation point when it exists and showed that the pressure gradient is adverse.
They also constructed Navier--Stokes solutions with the same qualitative behavior by a perturbation argument.

Motivated by \cite{CK2026}, our first objective is to provide a more precise and quantitative description of separation points and their dynamics. 
For the tangential velocity component, Chang and Kang proved uniqueness of the separation time in \((1,2)\) under the regularity and initial compatibility hypotheses of \cite[Theorem~1.9]{CK2026}, together with the additional conditions \(M(2)<0\) and \(M'(t)<-e_0<0\) on \((1,2)\).
In the present paper, log-concavity of the temporal boundary data alone gives this uniqueness, for every H\"older exponent \(0<\beta<1\) in {\rm\textbf{(T1)}} (Theorem~\ref{thm-0324-01}).
We further investigate the temporal structure of separation points for both tangential and normal velocity components and derive far-field asymptotic expansions that provide explicit convergence and divergence rates in terms of the boundary data.
Our second objective is to study the normal velocity component and to extend the analysis to all velocity components for tangential boundary data.
For normal boundary data, Chang and Kang also proved that the normal velocity component has a unique separation time in \((1,2)\) at sufficiently distant boundary points \cite[Theorem~1.9(ii)]{CK2026}.
Under our definition and assumptions, we obtain two separation times at each sufficiently distant boundary point.
The first occurs before \(t=1\), while the boundary data are still nonzero, as in \cite[Theorem~1.9(ii)]{CK2026}.
The second occurs after \(t=1\), when the prescribed boundary data have vanished, and so lies outside the interval considered there.
We derive the asymptotic behavior of both times.
For both normal and tangential boundary data, we determine the signs of the tangential pressure derivatives and identify regions of favorable and adverse pressure gradients relative to the separation points.

The quantitative analysis in this paper rests on the fact that, via the solution formula \eqref{0330-01}, the boundary derivatives of the velocity and the tangential pressure derivatives reduce to explicit space-time integrals with kernels involving confluent hypergeometric functions.
By Proposition~\ref{prop-260046}{\rm(i)}, \(\partial_{x_n}w_i\) vanishes at separation points of \(w_i\) for \(i<n\), and \(\partial_{x_n}^2w_n\) vanishes at separation points of \(w_n\).
Conversely, a point \((x',t)\) with \(x'\notin\supp\psi\) is a separation point if this derivative, as a function on \(\Sigma\times(0,\infty)\), takes both signs in every neighborhood of \((x',t)\).
A detailed analysis of the confluent hypergeometric functions and their zeros allows us to establish the existence, uniqueness, and asymptotic behavior of the separation points, as well as the signs of the pressure derivatives.
Related special-function representations occur in the half-plane vortex analysis of Dalibard and Gallay \cite{DalibardGallay2026} and in the viscous vortex-ring analysis of Fukumoto and Kaplanski \cite{FukumotoKaplanski2008,Fukumoto2010}.

\vspace{\baselineskip}

\textbf{Overview.}\quad
We now summarize the main results of this paper, using the notation of Subsection~\ref{subsect-notation}.
We derive formulas for the first normal derivatives of the tangential velocity components and the second normal derivative of the normal component at the boundary, and use their signs to study separation points outside the support of the spatial boundary data.
We also derive formulas for the tangential pressure derivatives and determine their signs for both normal and tangential boundary data. These derivatives are initially represented as space-time convolutions whose kernels are themselves defined by spatial convolutions over \(\Sigma\).
The key step is to evaluate the convolutions defining these kernels explicitly, which gives exact formulas in terms of confluent hypergeometric functions.

\vspace{\baselineskip}
\textbf{Normal boundary data.}
\quad
For localized normal boundary data, we study separation points of the tangential and normal velocity components, together with the signs of the associated tangential pressure derivatives.

\textbf{Normal boundary data: tangential velocity components.}
\quad
We first consider the tangential velocity components in Theorem~\ref{thm-0324-01}.
At every boundary point \(x'\) with \(|x_i|>1\), we establish the existence of a separation point \((x',t_i^*(x'))\) (see Definition~\ref{260000}) with \(t_i^*(x')\) in the time interval where the temporal boundary data decrease, and we exclude separation while they initially increase and after the boundary data vanish.
We also estimate how fast \(\partial_{x_n}w_i(x',0,\cdot)\) vanishes near \(t_i^*(x')\), with rates depending on the H\"older exponent \(\beta\) of \(\phi'\).
For every \(0<\beta\leq\frac12\), Proposition~\ref{prop-0512-01} constructs temporal boundary data for which the upper bound is attained on both sides of a separation time.
The rate is \(h^{\beta+1/2}\) for \(\beta<1/2\) and \(h\ln(1/h)\) for \(\beta=1/2\), where \(h\) is the distance to that time.
If \(\phi\) is log-concave on \((0,1)\), the separation point is unique for every \(0<\beta<1\).
Under {\rm\textbf{(T1)}}--{\rm\textbf{(T4)}}, \(t_i^*(r\mathbf e_i')\) and \(t_i^*(-r\mathbf e_i')\) are strictly increasing in \(r\) for sufficiently large \(r\), and
\begin{align*}
	t_i^*(x_i\mathbf e_i')=t_\infty+\frac{nT_1(t_\infty)}{T_0'(t_\infty)}x_i^{-2}+O(|x_i|^{-1-2\beta})\qquad\text{as }|x_i|\to\infty,
\end{align*}
where \(T_k\) is defined in \eqref{0323-02} and \(t_\infty\), defined in \eqref{260004}, is the unique zero in \((0,1)\) of the Caputo half-derivative of \(\phi\).

\textbf{Normal boundary data: normal velocity component.}
\quad
In contrast, the normal velocity component has two distinct separation points \((x',t_n^*(x'))\) and \((x',\widetilde{t_n^*}(x'))\) at every sufficiently distant boundary point \(x'\), one before and one after the boundary data vanish (see Theorem~\ref{thm-0507-01}).
For all nonzero temporal boundary data satisfying {\rm\textbf{(T1)}}, the later separation time \(\widetilde{t_n^*}(x')\) is a simple zero of \(\partial_{x_n}^2w_n(x',0,\cdot)\), is unique at sufficiently distant points, increases strictly along each ray there, and satisfies
\begin{align*}
	\widetilde{t_n^*}(x')=\frac{|x'|^2}{4z_1}+O(|x'|)\qquad\text{as }|x'|\to\infty,
\end{align*}
where \(z_1\) is the positive zero of a confluent hypergeometric function.
We also estimate how fast \(\partial_{x_n}^2w_n(x',0,\cdot)\) vanishes near the earlier separation time \(t_n^*(x')\).
Log-concavity gives uniqueness of \(t_n^*(x')\) at sufficiently distant points for every \(0<\beta<1\).
Under {\rm\textbf{(T1)}}--{\rm\textbf{(T4)}}, it increases strictly along each ray at large distances and satisfies
\begin{align*}
	t_\infty-t_n^*(x')=-\frac{3nT_1(t_\infty)}{T_0'(t_\infty)}|x'|^{-2}+O(|x'|^{-1-2\beta}).
\end{align*}

\textbf{Normal boundary data: tangential pressure derivatives.}
\quad
We further investigate the tangential pressure derivatives in Theorem~\ref{thm-260019}.
For radial spatial boundary data satisfying {\rm\textbf{(S)}} and temporal boundary data satisfying {\rm\textbf{(T1)}} and {\rm\textbf{(T2)}}, we obtain a continuous radius \(r_{\rm nor}(t)\) for \(0<t<t_1\).
The sign of \(\partial_{x_i}p(x',0,t)\) agrees with that of \(x_i\) for \(1<|x'|\leq r_{\rm nor}(t)-1\) and is opposite for \(|x'|\geq r_{\rm nor}(t)+1\).
This radius tends to zero as \(t\to0^+\) precisely when \(\phi(t)/\phi'(t)\to0\), a condition implied by log-concavity of \(\phi\) on \((0,1)\).
At the other endpoint,
\begin{align*}
	r_{\rm nor}(t)\sim\left(\frac{n\phi(t)}{\phi'(t)}\right)^{\frac12}\qquad\text{as }t\to t_1^-.
\end{align*}
For \(t<t_1\) sufficiently close to \(t_1\), the derivative changes sign across the sphere of radius \(\rho_{\rm nor}(t)\) on every ray with \(x_i\ne0\).
On such a ray, \(|x'|=\rho_{\rm nor}(t)\) is the only zero with \(|x'|>1\).
Moreover, \(\rho_{\rm nor}(t)^2-r_{\rm nor}(t)^2\) tends to \(\frac{n+3}{n-1}\frac{M_2}{M_0}>0\) as \(t\to t_1^-\), where \(M_0,M_2\) are defined in \eqref{260003}.
For \(t\geq t_1\), \(\partial_{x_i}p\) has the sign of \(x_i\) throughout \(\{|x'|>1\}\).

\vspace{\baselineskip}

\textbf{Tangential boundary data.}
\quad
For localized tangential boundary data, we study separation points of the tangential and normal velocity components, together with the signs of the associated tangential pressure derivatives.
The behavior of the separation points depends on the direction of the velocity component relative to the boundary data.
Without loss of generality, we take the boundary data in the \(x_1\)-direction.

\textbf{Tangential boundary data: velocity component in the direction of the boundary data.}
\quad
Suppose that \(\phi\) satisfies {\rm\textbf{(T1)}} and {\rm\textbf{(T2)}}.
For the velocity component \(w_1\) in the direction of the boundary data, Theorem~\ref{thm-260010} shows that, when \(n\geq3\), separation occurs at some time \(0<t<1\) at a given boundary point \(x'\) with \(|x'|>1\) if and only if
\begin{align*}
	\mathcal S(x'):=\partial_{x_1}R_1'\psi(x')>0.
\end{align*}
If \(\phi\) is also log-concave on \((0,1)\), the separation time is unique at every \(x'\) with \(|x'|>1\) and \(\mathcal S(x')>0\).
Separation is excluded throughout \(\{|x'|>1\}\) in two dimensions and, in every dimension, for \(t\geq1\).
For radial \(\psi\), it is also excluded on the \(x_1\)-axis at \(x'=x_1\mathbf e_1'\) with \(|x_1|>1\).

\textbf{Tangential boundary data: velocity component in the direction of the boundary data for radial \(\psi\).}
\quad
Now suppose that \(n\geq3\) and that the spatial boundary data are radial.
For every \(t_0\in(0,1)\), there is an \(R_{t_0}>1\) such that the radius \(\rho_*(x_1,t)\) of the separation points \((x_1,\rho_*(x_1,t)\omega'')\), \(\omega''\in\mathbb S^{n-3}\), is unique whenever \(|x_1|>R_{t_0}\) and \(0<t\leq t_0\).
This radius is strictly increasing in \(|x_1|\) and satisfies
\begin{align*}
	\rho_*(x_1,t)=\sqrt{n-1}\,|x_1|+\frac{M_2}{2(n-1)^{3/2}M_0|x_1|}+O(|x_1|^{-3})\qquad\text{as }|x_1|\to\infty,
\end{align*}
uniformly for \(0<t\leq t_0\), where \(M_k\) is defined in \eqref{260003}.
We also derive the asymptotics of the separation times at points of the hyperplane \(x_1=0\) as \(|x''|\to\infty\), under additional smoothness and nondegeneracy assumptions on the temporal boundary data.

\textbf{Tangential boundary data: tangential velocity component in a different direction.}
\quad
For a tangential velocity component \(w_2\) in a direction different from that of the boundary data, Theorem~\ref{thm-260014} gives, in dimensions \(n\geq3\), separation points whose spatial coordinates are independent of time throughout \(0<t\leq1\).
At each \(0<t<1\), this separation point is unique on each line parallel to the \(x_1\)-axis with \(|x_2|>1\), and on each line parallel to the \(x_2\)-axis with \(|x_1|>1\), once the line is sufficiently far from the corresponding axis.
As \(|x'|\to\infty\) along rays, these coordinates approach \(m_1/M_0\) and \(m_2/M_0\), respectively, with an error of order \(|x'|^{-1}\).
If \(\int_\Sigma(y_1-m_1/M_0)y_2\psi(y')\,dy'\neq0\), we also establish their strict monotonicity for sufficiently large \(|x'|\).
For radial spatial boundary data, a point \((x',t)\) with either \(0<t<1\) and \(|x'|>1\), or \(t\geq1\) and \(|x'|\geq\sqrt2\), is a separation point if and only if \(x_1x_2=0\).

\textbf{Tangential boundary data: normal velocity component.}
\quad
Theorem~\ref{thm-260017} describes the separation points of the normal velocity component \(w_n\).
For \(n\geq3\), at every positive time, we establish the existence of a separation point on each line parallel to the \(x_1\)-axis whose distance from the origin exceeds one.
For sufficiently large distances, this point is unique, and, by Proposition~\ref{prop-260127}, its \(x_1\)-coordinate \(x_1^*(x'',t)\) is continuously differentiable in \(t\).
The same proposition shows that, for \(t\geq3/2\), the second normal derivative has exactly one simple zero on every line parallel to the \(x_1\)-axis, and that this zero gives a separation point whenever it lies outside \(\supp\psi\).
For each fixed \(0<t<1\),
\begin{align*}
	x_1^*(x'',t)=\frac{m_1}{M_0}+O(|x''|^{-1})\qquad\text{as }|x''|\to\infty.
\end{align*}
After the boundary data vanish, the limit as \(|x''|\to\infty\) along a fixed ray may instead be determined by the intersection of \(\supp\psi\) with a supporting hyperplane orthogonal to the ray, under an additional geometric condition.
On the other hand, for each fixed \(x''\) with \(|x''|\) sufficiently large,
\begin{align*}
	\lim_{t\to\infty}x_1^*(x'',t)=\frac{m_1}{M_0}.
\end{align*}
Thus the limits of \(x_1^*(x'',t)\) as \(|x''|\to\infty\) and as \(t\to\infty\) need not agree.
The example following Theorem~\ref{thm-260017} realizes the two values \(0\) and \(1/4\).

\textbf{Tangential boundary data: pressure derivative in the direction of the boundary data.}
\quad
We also analyze the pressure derivative \(\partial_{x_1}p\) in the direction of the boundary data in Theorem~\ref{thm-260021}.
For radial spatial boundary data and temporal boundary data that initially increase as in {\rm\textbf{(T2)}}, this derivative is negative on the \(x_1\)-axis for \(|x_1|\geq1\) at all sufficiently small positive times.
In dimensions \(n\geq3\), for large \(|x'|\), the sign changes across the cones
\begin{align*}
	|x''|=\sqrt{n-1}\,|x_1|\quad(t\neq t_\infty),\qquad |x''|=\sqrt{n+1}\,|x_1|\quad(t=t_\infty),
\end{align*}
where the case \(t=t_\infty\) requires the spatial boundary data to be radially nonincreasing.
For radially nonincreasing spatial boundary data and temporal boundary data that are concave on \((t_1,t_\infty)\), we prove that, shortly before \(t_\infty\), the zero \(\rho(t)>1\) of \(\partial_{x_1}p\) on the positive \(x_1\)-axis is unique and satisfies
\begin{align*}
	\lim_{t\to t_\infty^-}(t_\infty-t)\rho(t)^2=-\frac{n(n+1)T_1(t_\infty)}{(n-1)T_0'(t_\infty)}.
\end{align*}
For sufficiently large times, it suffices that \(\psi\) satisfies {\rm\textbf{(S)}} and is radial, and that \(\phi\) satisfies {\rm\textbf{(T1)}} and \(\phi\not\equiv0\).
On the positive \(x_1\)-axis with \(x_1\geq1\), the unique zero is simple and satisfies
\begin{align*}
	\rho(t)=2\sqrt{z_7t}+O_n(t^{-1/2})\qquad\text{as }t\to\infty,
\end{align*}
where \(z_7\) is the positive zero of a linear combination of confluent hypergeometric functions.

\textbf{Tangential boundary data: zero curves of the pressure derivative.}
\quad
For \(n\geq3\) and radial \(\psi\), fix \(t>0\) with \(\norm{\phi}_{L^1(0,t)}>0\).
Under {\rm\textbf{(T1)}}, the pressure derivative changes sign at a positive \(x_1\) on every line parallel to the \(x_1\)-axis with \(|x''|\geq1\).
For large \(|x''|\), this zero is simple and unique.
As \(|x''|\to\infty\), this zero is asymptotic to \(|x''|/\sqrt{n-1}\) if \(T_0(t)\ne0\) and to \(|x''|/\sqrt{n+1}\) if \(T_0(t)=0\).
Under the additional hypotheses of Theorem~\ref{thm-260021}{\rm(vi)(b)}, shortly before \(t_\infty\), the largest positive zero \(\rho(|x''|,t)\) exists on every line parallel to the \(x_1\)-axis, including the axis itself, decreases in \(|x''|\) near the axis, and increases for large \(|x''|\).

\textbf{Tangential boundary data: pressure derivative in a different direction.}
\quad
Finally, Theorem~\ref{thm-260027} addresses the pressure derivative \(\partial_{x_2}p\) in a tangential direction different from that of the boundary data.
For \(n\geq3\) and radial spatial boundary data, which again need not be radially nonincreasing, we show that, throughout the region \(|x'|\geq\sqrt2\), its sign is opposite to that of \(x_1x_2\) for \(0<t\leq t_1\) and agrees with it for \(t\geq t_\infty\).
For \(t_1<t<t_\infty\), there is a continuous radius \(r_{\rm tan}(t)\).
In the region \(|x'|\geq\sqrt2\), the sign of \(\partial_{x_2}p(x',0,t)\) agrees with that of \(x_1x_2\) for \(|x'|\leq r_{\rm tan}(t)-1\) and is opposite for \(|x'|\geq r_{\rm tan}(t)+1\).
For \(t<t_\infty\) sufficiently close to \(t_\infty\), each ray with \(x_1x_2\ne0\) has a unique sign-changing zero in \(|x'|\geq\sqrt2\), at radius \(\rho_{\rm tan}(t)\), where
\begin{align*}
	\rho_{\rm tan}(t)\sim\left(-\frac{(n+2)T_1(t_\infty)}{T_0'(t_\infty)}\right)^{\frac12}(t_\infty-t)^{-\frac12}\qquad\text{as }t\to t_\infty^-.
\end{align*}
Together with the signs of the corresponding velocity derivatives near the boundary, these pressure results identify regions of favorable and adverse pressure gradients and describe their evolution in time.

\subsection{Definition of the separation point}

\begin{definition}\label{260000}
	Let \(1\leq j\leq n\), let \(w\) be the solution of the Stokes system \eqref{0308-01} given by \eqref{0330-01}, with the boundary data \(g(x',t)=g_j(x',t)\mathbf e_j\) as in \eqref{0324-13}, and let \(1\leq i\leq n\).
	A point \((x',t_*)\in (\Sigma\setminus \supp\psi)\times (0,\infty)\) is a separation point of \(w_i\) if there exist \(\varepsilon_*, M_*>0\) such that, for any \(0<\varepsilon<\varepsilon_*\) and \(0<M<M_*\), the sets \(S_{i,\varepsilon, M}^+, S_{i,\varepsilon, M}^-\subset \R_+^n\times (0,\infty)\) given by
\begin{align*}
		S_{i,\varepsilon, M}^+:= Q_{\varepsilon, M}(x', t_*) \cap \{(z,s): w_i(z,s)>0\}, \\
		S_{i,\varepsilon, M}^-:= Q_{\varepsilon, M}(x', t_*) \cap \{(z,s): w_i(z,s)<0\}
	\end{align*}
	both have positive Lebesgue measure, where \(Q_{\varepsilon, M}(x',t):= \{(z,s)\in\R_+^n \times (0,\infty): |z'-x'|<M, \, 0<z_n<\varepsilon,\, |s-t|<M^2\}\).
	In this case, \(t_*\) is a separation time of \(w_i\) at \(x'\).

	Moreover, a point \((x', t_*)\in (\Sigma\setminus \supp\psi)\times(0,\infty)\) is a separation point of \(w\) if there exists \(1\leq i\leq n\) such that \((x', t_*)\) is a separation point of \(w_i\).
\end{definition}

\begin{figure}[H]
	\centering
	\begin{tikzpicture}[x=1.12cm,y=0.70cm,font=\small,>=stealth]

		\fill[gray!7]
		(-4.40,0)--(4.75,0)--(7.445,-3.85)--(-1.705,-3.85)--cycle;
		\draw[gray!55,line width=.6pt]
		(-4.40,0)--(4.75,0)--(7.445,-3.85)--(-1.705,-3.85)--cycle;

		\node[gray!60!black] at (-0.85,-3.52)
		{$\Sigma\times(0,\infty)$};

		\fill[
			pattern=north east lines,
			pattern color=gray!60
		]
		(-4.40,0)
		--(4.75,0)
		--(5.184,-.62)
		--(-3.966,-.62)
		--cycle;

		\draw[dashed,gray!65,line width=.55pt]
		(-4.40,0)
		--(4.75,0)
		--(5.184,-.62)
		--(-3.966,-.62)
		--cycle;

		\node[font=\small,inner sep=1pt]
		at (-2.10,-.34)
		{$\supp\psi\times(0,\infty)$};

		\draw[->,line width=.85pt]
		(-4.40,0)--(5.15,0) node[right] {$t$};

		\draw[->,line width=.85pt]
		(-4.40,0)--(-4.40,2.05) node[above] {$x_n$};

		\draw[->,line width=.85pt]
		(-4.40,0)--(-1.495,-4.15) node[below right] {$x'$};

		\begin{scope}[shift={(0,-1.30)}]
			\coordinate (A)at (-0.35,-0.15);
			\coordinate (B)at ( 3.05,-0.15);

			\coordinate (At) at (-0.35, 2.52);
			\coordinate (Bt) at ( 3.05, 2.52);

			\coordinate (D)at ( 1.45,-2.03);
			\coordinate (C)at ( 4.85,-2.03);

			\coordinate (Dt) at ( 1.45, 0.64);
			\coordinate (Ct) at ( 4.85, 0.64);

			\coordinate (E)at (0.30,-0.15);
			\coordinate (Et) at (0.75,2.52);

			\coordinate (F)at ( 3.85,-2.03);
			\coordinate (Ft) at ( 4.35, 0.64);
			\coordinate (P)at ( 2.18,-1.13);

			\fill[blue!18,opacity=.72]
			(A)--(E)
			.. controls (.50,.62) and (.60,1.72) .. (Et)
			--(At)--cycle;

			\fill[red!18,opacity=.72]
			(E)--(B)--(Bt)--(Et)
			.. controls (.60,1.72) and (.50,.62) .. cycle;

			\fill[blue!22,opacity=.72]
			(D)--(F)
			.. controls (3.63,-.80) and (3.82,.15) .. (Ft)
			--(Dt)--cycle;

			\fill[red!22,opacity=.72]
			(F)--(C)--(Ct)--(Ft)
			.. controls (3.82,.15) and (3.63,-.80) .. cycle;

			\fill[blue!15,opacity=.72]
			(At)--(Et)
			.. controls (1.55,2.38) and (3.20,1.35) .. (Ft)
			--(Dt)--cycle;

			\fill[red!15,opacity=.72]
			(Et)--(Bt)--(Ct)--(Ft)
			.. controls (3.20,1.35) and (1.55,2.38) .. cycle;

			\fill[blue!14,opacity=.68]
			(A)--(D)--(Dt)--(At)--cycle;

			\fill[red!14,opacity=.68]
			(B)--(C)--(Ct)--(Bt)--cycle;

			\shade[
				left color=gray!5,
				right color=gray!38,
				shading angle=18,
				opacity=.90
			]
			(E)
			.. controls (.55,-.18) and (1.35,-.82) .. (P)
			.. controls (2.70,-1.48) and (3.38,-1.93) .. (F)
			.. controls (3.98,-1.10) and (4.22,-.15) .. (Ft)
			.. controls (3.73,1.43) and (1.78,2.37) .. (Et)
			.. controls (.70,1.66) and (.50,.55) .. cycle;

			\draw[gray!70!black,line width=.8pt]
			(E)
			.. controls (.55,-.18) and (1.35,-.82) .. (P)
			.. controls (2.70,-1.48) and (3.38,-1.93) .. (F)
			.. controls (3.98,-1.10) and (4.22,-.15) .. (Ft)
			.. controls (3.73,1.43) and (1.78,2.37) .. (Et)
			.. controls (.70,1.66) and (.50,.55) .. cycle;

			\draw[gray!62,line width=.4pt]
			(.57,.36)
			.. controls (1.40,.08) and (2.63,-.70) ..
			(3.58,-.95);

			\draw[gray!62,line width=.4pt]
			(.63,1.06)
			.. controls (1.52,.75) and (2.80,-.04) ..
			(3.72,-.27);

			\draw[gray!62,line width=.4pt]
			(.69,1.77)
			.. controls (1.62,1.43) and (3.00,.61) ..
			(3.86,.39);

			\draw[dashed,line width=.7pt] (A)--(B) (B)--(Bt) (B)--(C);
			\draw[line width=.85pt] (A)--(At) (A)--(D) (C)--(D)  (C)--(Ct);
			\draw[line width=.85pt] (At)--(Bt)--(Ct)--(Dt)--cycle;
			\draw[line width=.85pt] (Dt)--(D);

			\node[blue!65!black,align=center,font=\large,anchor=east] (pluslab)
			at (-1.55,2.10)
			{$S_{i,\varepsilon,M}^{+}$\\[3pt]$w_i>0$};
			\draw[->,blue!65!black,line width=.8pt]
			(pluslab.east) .. controls (-0.95,2.05) and (-0.45,1.75) .. (.15,1.60);

			\node[red!65!black,align=center,font=\large,anchor=west] (minuslab)
			at (5.35,2.10)
			{$S_{i,\varepsilon,M}^{-}$\\[3pt]$w_i<0$};
			\draw[->,red!65!black,line width=.8pt]
			(minuslab.west) .. controls (4.80,2.05) and (4.10,1.75) .. (3.55,1.55);

			\node[inner sep=2pt] at (2.02,3.20)
			{$Q_{\varepsilon,M}(x',t_*)$};

			\fill (P) circle (1.8pt);

			\node[anchor=north,inner sep=1.5pt] (sep) at (2.40,-1.20)
			{$\bigl(x',t_*\bigr)$};

			\draw[<->] (1.45,-2.27)--(4.85,-2.27)
			node[pos=.52,below,yshift=-6pt,inner sep=1pt] {$|s-t_*|<M^2$};

			\draw[<->] (-0.72,-0.15)--(1.08,-2.03)
			node[midway,below,rotate=-33,inner sep=2pt]
			{$|z'-x'|<M$};

			\draw[<->] (5.12,-2.03)--(5.12,.64)
			node[midway,right,inner sep=1pt] {$0<z_n<\varepsilon$};

		\end{scope}

	\end{tikzpicture}
	\caption{A schematic illustration of a separation point given in Definition \ref{260000}}
	\label{260001}
\end{figure}
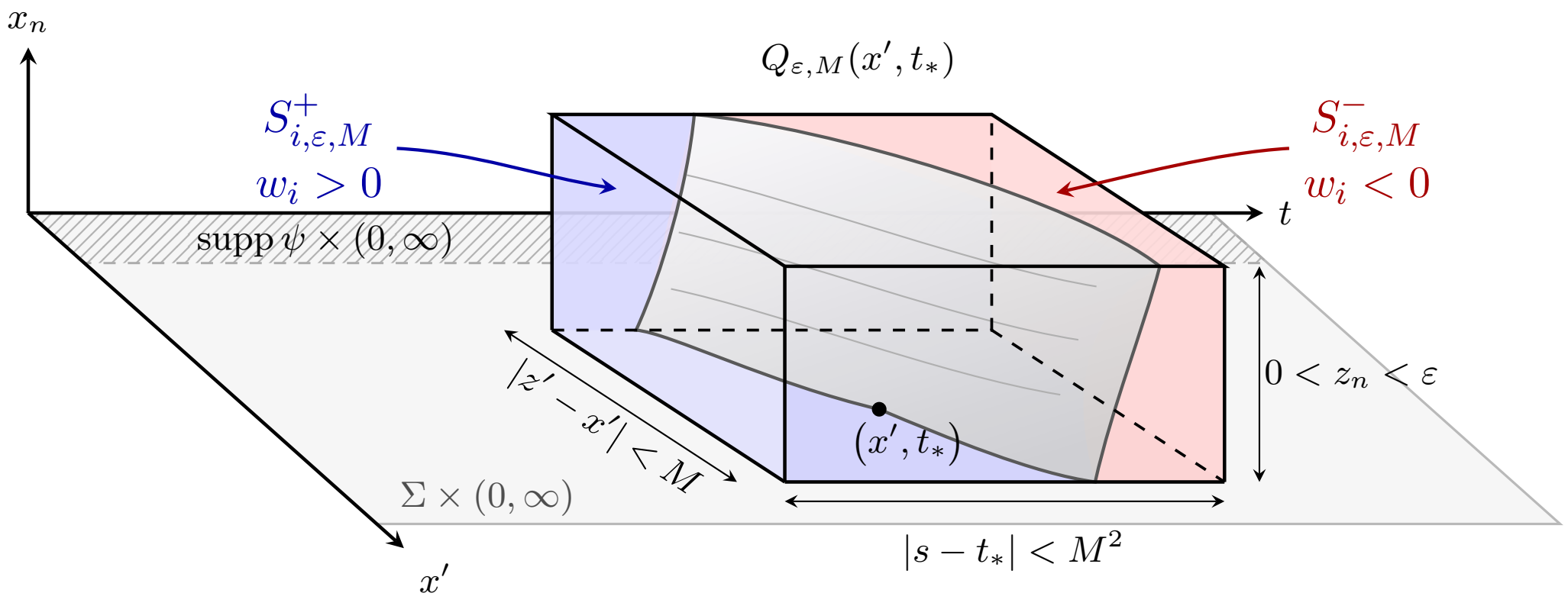

This definition involves only the values of the single component \(w_i\) near \((x',0,t_*)\).
Since \(Q_{\varepsilon,M}\) increases with \(\varepsilon\) and \(M\), the condition in Definition \ref{260000} holds for some \(\varepsilon_*,M_*>0\) if and only if it holds for all \(\varepsilon,M>0\).
Figure~\ref{260001} illustrates the definition.
The two signs may be obtained by varying the tangential position, time, or both.
In particular, separation points may persist over an interval of time due to spatial symmetry, without any temporal sign change at a fixed spatial location.
The individual theorems state when the relevant normal derivative at the boundary has a zero in time at a fixed \(x'\) and when this zero is unique.
Proposition \ref{prop-260046} gives the precise criterion in terms of normal derivatives at the
boundary.

\begin{remark}
If \(x'\notin\supp\psi\), a separation point in the sense of \cite[Definition~1.1]{CK2026} is one in the sense of Definition~\ref{260000}, since \(w_i\) is continuous in the interior.
The converse fails.
For boundary data in the \(\mathbf e_1\)-direction with radial \(\psi\), \(\phi\) satisfying {\rm\textbf{(T1)}} and {\rm\textbf{(T2)}}, and \(n\geq3\), reflection across \(x_2=0\) makes \(w_2\) odd in \(x_2\), so \(w_2\) vanishes on \(\{x_2=0\}\). Nevertheless, by Theorem~\ref{thm-260014}{\rm(ii)}, \((r\mathbf e_1',t)\) with \(r>1\) is a separation point of \(w_2\) for every \(0<t<1\).
\end{remark}
\begin{remark}[Physical meaning]\label{rmk-physical-meaning}
By Proposition~\ref{prop-260046}, separation points in the sense of Definition~\ref{260000} can be studied through the sign of the leading normal derivative of \(w_i\) at the boundary.
For \(i<n\), this is the wall shear stress \(\partial_{x_n}w_i(x',0,t)\), and a separation point is a point where it changes sign, which is the classical criterion of vanishing wall shear \cite{Prandtl1904,Goldstein1948}.
For \(w_n\), \(\partial_{x_n}w_n=0\) on \(\Sigma\setminus\supp\psi\), and a separation point is a point where the near-wall normal velocity changes sign, as on attachment and detachment lines \cite{PerryFairlie1975,TobakPeake1982}.
For \(w_n\), and for \(w_i\) with \(i\ne j\), the term separation point is used in this generalized sense.
\end{remark}

\subsection{Notation}\label{subsect-notation}
\begin{itemize}
	\item \(\{\mathbf e_1',\dots,\mathbf e_{n-1}'\}\) is the standard basis of \(\Sigma\).
	\item For \(1\leq\ell\leq n-1\), let \((\mathbf e_\ell')^\perp=\{z'\in\Sigma:z'_\ell=0\}\) and write \(x^{(\ell)}:=x'-x_\ell\mathbf e_\ell'\in(\mathbf e_\ell')^\perp\). In particular, \(x''=x^{(1)}=(0,x_2,\ldots,x_{n-1})\).
	\item We write \((s,x''):=s\mathbf e_1'+x''\) for \(s\in\R\), so that \(x'=(x_1,x'')\).
	\item We write \(\mathbb S^{n-2}\) for the unit sphere in \(\Sigma\) and, for \(n\geq3\), identify \(\mathbb S^{n-3}\) with the unit sphere in \((\mathbf e_1')^\perp\).
	Points of \(\mathbb S^{n-2}\) and \(\mathbb S^{n-3}\) are denoted by \(\omega'\) and \(\omega''\), respectively.
	\item For \(r>0\), \(B_r'\) denotes the open ball in \(\Sigma\) of radius \(r\) centered at the origin.
\item \(\nabla'\) and \(\Delta'\) denote the gradient and the Laplacian in the tangential variable \(x'\), respectively.
\item For \(k\geq0\) and \(1\leq\ell\leq n-1\),
\begin{align}\label{260003}
M_k:=\norm{|y'|^k\psi}_{L^1(\Sigma)}, \qquad m_\ell:=\int_\Sigma y_\ell\psi(y')\,dy'.
\end{align}
We also write \(\mathbf m:=(m_1,\dots,m_{n-1})\).
\item For integers \(k \geq 0\),
\begin{align}\label{0323-02}
	T_k(t) := \int_0^t (t-s)^{\frac{2k-1}{2}}\phi'(s)\,ds.
\end{align}
Up to normalization, \(T_0\) is the Caputo half-derivative of \(\phi\), and \(T_1\) is the Riemann--Liouville half-integral of \(\phi\).
\end{itemize}

\subsection{Main results}\label{260002}

We first list the assumptions on the boundary data \(\psi\) and \(\phi\).
\begin{enumerate}
	\item[\textbf{(S)}] \(\psi\in C_c^{2,\alpha}(\Sigma)\) for some \(\alpha\in(0,1)\), \(\supp\psi\subset\overline{B_1'}\), and \(\psi\) is nonnegative and not identically zero.

	\item[\textbf{(T1)}] \(\phi\in C^{1,\beta}([0,\infty))\) for some \(\beta\in(0,1)\), and \(\phi\) is nonnegative.
	Moreover, \(\phi(0)=0\) and \(\phi(t)=0\) for \(t\geq1\).

	\item[\textbf{(T2)}] There exists \(t_1\in(0,1)\) such that \(\phi'(s)>0\) for \(0<s<t_1\) and \(\phi'(s)<0\) for \(t_1<s<1\).

	\item[\textbf{(T3)}] The H\"older exponent \(\beta\) in {\rm\textbf{(T1)}} satisfies \(\beta\in\left(\frac12,1\right)\).

	\item[\textbf{(T4)}] \(\phi\) is strictly log-concave on \((0,1)\).
\end{enumerate}
Whenever {\rm\textbf{(T2)}} holds and \(T_0\) has a zero in \((t_1,1)\), we denote its first zero by
\begin{align}\label{260004}
	t_\infty:=\inf\{t\in(t_1,1):T_0(t)=0\}.
\end{align}
We will also need the following assumption on \(\phi\).
\begin{enumerate}
	\item[\textbf{(T5)}] \(\phi\) is concave on \((t_1,t_\infty)\).
\end{enumerate}

Assumptions {\rm\textbf{(T1)}} and {\rm\textbf{(T2)}} imply \(\phi(t)>0\) for every \(0<t<1\).
The boundedness of \(\phi'\) also gives
\begin{align}\label{0511-02}
	\int_0^t \frac{|\phi'(s)|}{\sqrt{t-s}}\,ds < \infty
\end{align}
for all \(t>0\).

\begin{remark}\label{rmk-0730-03}
	Under {\rm\textbf{(T1)}} and {\rm\textbf{(T2)}}, \(T_0\) is continuous, \(T_0>0\) on \((0,t_1]\), and \(T_0(1)< (1-t_1)^{-1/2}\int_0^1\phi'(s)\,ds=0\), so \(t_\infty\in(t_1,1)\) and \(T_0(t_\infty)=0\).
	If \(\phi\) is log-concave on \((0,1)\), then \(t_\infty\) is the only positive zero of \(T_0\).
	Indeed, extending \(\phi\) by zero for negative arguments, integration by parts in \eqref{0323-02} gives
	\begin{align*}
		\frac{T_0(t)}{\phi(t)}
		=\frac12\int_0^\infty
		\left(1-\frac{\phi(t-s)}{\phi(t)}\right)\frac{ds}{s^{3/2}},\qquad 0<t<1,
	\end{align*}
	and, for \(0<t<u<1\), log-concavity gives \(\phi(u-s)/\phi(u)\geq\phi(t-s)/\phi(t)\) for \(0<s<t\), strictly for \(t<s<u\).
	Hence \(T_0/\phi\) is strictly decreasing on \((0,1)\), while \(T_0(t)=-\frac12\int_0^1\phi(s)(t-s)^{-3/2}\,ds<0\) for \(t>1\).
	Under {\rm\textbf{(T1)}}--{\rm\textbf{(T4)}}, Lemma~\ref{lem-0325-01}{\rm(ii)} with \(K(u)=u^{-1/2}\) and \(\sigma=\frac12\) gives \(T_0'(t_\infty)<0\).
\end{remark}

\begin{remark}
The assumptions {\rm\textbf{(T1)}}--{\rm\textbf{(T5)}} are compatible.
For example, \(\phi(t)=t^2(1-t)^2\) on \([0,1]\), extended by zero, satisfies {\rm\textbf{(T1)}}--{\rm\textbf{(T4)}} with \(t_1=1/2\), and
\begin{align}\label{example:abel-polynomial}
T_0(t)=\frac8{105}t^{3/2}(35-84t+48t^2),\qquad t_\infty=\frac{21-\sqrt{21}}{24},\quad 0<t<1.
\end{align}
Since \(\phi''<0\) on \((1/2,t_\infty)\), it also satisfies {\rm\textbf{(T5)}}. See Remark~\ref{rmk-T1-relaxed} for a comment on the regularity required in {\rm\textbf{(T1)}} near \(t=1\).
\end{remark}

We state the first main theorem concerning the separation points of the tangential velocity components \(w_i\) with normal boundary data.
\begin{theorem}\label{thm-0324-01}\textbf{(Tangential velocity component with normal boundary data)}
	Let \(n\geq2\), \(1\leq i\leq n-1\), and let \(w\) solve the Stokes system \eqref{0308-01} with boundary data \eqref{0324-13}, where \(j=n\) and \(\psi\) satisfies {\rm\textbf{(S)}}.
	Let \(x'\in\Sigma\) satisfy \(|x_i|>1\).
	Then the following hold.
	\begin{enumerate}
		\item\label{item-0429-11} \textbf{(Nonexistence and existence)} Suppose that \(\phi\) satisfies {\rm\textbf{(T1)}}.
If \(0<t<1\), \(\phi\) is nondecreasing on \([0,t]\), and \(\phi(t)>0\), then \((x',t)\) is not a separation point of \(w_i\).
		Moreover, if \(t\geq1\), then \((x',t)\) is not a separation point of \(w_i\).

		On the other hand, if \(\phi\) also satisfies {\rm\textbf{(T2)}}, then there exists \(t_i^*(x')\in(t_1,1)\) such that \((x',t_i^*(x'))\) is a separation point of \(w_i\).
			\item\label{item-0429-12} \textbf{(Temporal vanishing rates)} Suppose that \(\phi\) satisfies {\rm\textbf{(T1)}} and {\rm\textbf{(T2)}}.
			Let \(\beta\in(0,1)\) be the H\"older exponent of \(\phi'\), and let \(t_i^*=t_i^*(x')\) be given in {\rm(i)}.
		Then, for \(h:=|t-t_i^*|>0\) sufficiently small,
		\begin{align*}
			|\partial_{x_n}w_i(x',0,t)| \leq \frac{C|x'|^2}{(|x'|-1)^{n+1}}\begin{cases}h^{\beta+\frac12}, & \beta<\frac12, \\
h|\ln h|,& \beta=\frac12, \\
h^\sigma,&\beta>\frac12,
			\end{cases}
		\end{align*}
		where \(C>0\) is independent of \(x'\) and, for \(\beta>1/2\), \(\sigma=1\) if \(\partial_t\partial_{x_n}w_i(x',0,t_i^*)\neq0\) and \(\sigma=\beta+\frac12\) otherwise.
		Moreover, if \(0<\beta<\frac12\) and
		\begin{align}\label{260005}
			\phi'(s)-\phi'(t_i^*)=-c_0(s-t_i^*)|s-t_i^*|^{\beta-1} +o(|s-t_i^*|^\beta)\qquad\text{as }s\to t_i^*,
		\end{align}
		for some \(c_0>0\), then
		\begin{align*}
			\lim_{\varepsilon\to0^+}\varepsilon^{-\beta-\frac12} \partial_{x_n}w_i(x',0,t_i^*+\varepsilon) =-C_\beta c_0\int_\Sigma \frac{x_i-y_i}{|x'-y'|^n}\psi(y')\,dy',
		\end{align*}
		and
		\begin{align*}
			\lim_{\varepsilon\to0^-}|\varepsilon|^{-\beta-\frac12} \partial_{x_n}w_i(x',0,t_i^*+\varepsilon) =\widetilde C_\beta c_0\int_\Sigma \frac{x_i-y_i}{|x'-y'|^n}\psi(y')\,dy',
		\end{align*}
		where \(C_\beta,\widetilde C_\beta>0\) depend only on \(n\) and \(\beta\). 
		\item\label{item-0511-01}\textbf{(Uniqueness)} Suppose that \(\phi\) satisfies {\rm\textbf{(T1)}} and {\rm\textbf{(T2)}} and is log-concave on \((0,1)\).
		Then the separation time \(t_i^*(x')\) is unique.
		Moreover, \(t_1<t_i^*(x')<t_\infty\) and
		\begin{align*}
			\sgn\bigl(\partial_{x_n}w_i(x',0,t)\bigr) =\begin{cases} \sgn(x_i),& 0<t<t_i^*(x'), \\
				-\sgn(x_i), & t>t_i^*(x').
			\end{cases}
		\end{align*}

		\item\label{item-0429-13}\textbf{(Monotonicity and asymptotics)} Suppose that \(\phi\) satisfies {\rm\textbf{(T1)}}--{\rm\textbf{(T4)}}.
		Then \(t_i^*(r\mathbf e_i')\) and \(t_i^*(-r\mathbf e_i')\) are strictly increasing in \(r\) for all sufficiently large \(r\).
		Moreover,
		\begin{align}\label{260006}
			t_i^*(x_i\mathbf e_i')=t_\infty +\frac{nT_1(t_\infty)}{T_0'(t_\infty)}x_i^{-2}+O(|x_i|^{-1-2\beta})
		\end{align}
		as \(|x_i|\to\infty\).
	\end{enumerate}
\end{theorem}
\begin{figure}[H]
	\centering
	\begin{tikzpicture}[x=1.1cm,y=2.15cm,>=latex,
			every node/.style={font=\footnotesize},
			declare function={sep(\x)=0.72-0.52/(\x*\x);}]

		\draw[->] (0,0) -- (6.15,0) node[below right] {\(x_1\)};
		\draw[->] (0,0) -- (0,1.08) node[above right,inner sep=1pt] {\(t\)};

		\draw[densely dotted,black!55] (0,1) -- (5.95,1);
		\node[left] at (0,1) {\(1\)};
		\draw[dashed,black!65] (0,0.72) -- (5.95,0.72);
		\node[left] at (0,0.72) {\(t_\infty\)};

		\draw[very thick]
		plot[domain=1.2:5.8,samples=80,smooth] (\x,{sep(\x)});
		\node[below right] at (2.55,{sep(2.55)})
		{\(t=t_1^*(x_1\mathbf e_1')\)};
		\node[font=\Large] at (1.65,0.25) {\(+\)};
		\node[font=\Large] at (1.65,0.86) {\(-\)};

		\draw (1,-0.018) -- (1,0.018);
		\node[below] at (1,0) {\(1\)};

		\begin{scope}[shift={(6.45,0.92)}]
			\draw[very thick] (0,0) -- (0.65,0);
			\node[right] at (0.78,0) {separation points};
			\draw[dashed,black!65] (0,-0.20) -- (0.65,-0.20);
			\node[right] at (0.78,-0.20) {asymptote};
		\end{scope}

	\end{tikzpicture}
	\caption{A schematic illustration of separation points and signs of \(\partial_{x_n}w_1\) on the positive \(x_1\)-axis when \(n=3\)}
	\label{260007}
\end{figure}
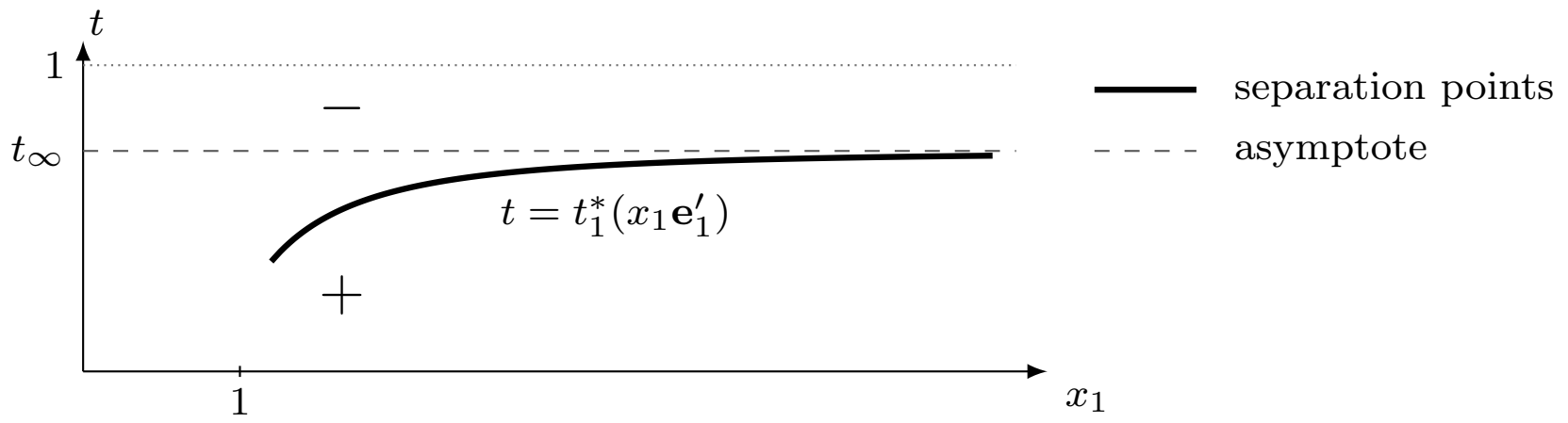

Theorem \ref{thm-0324-01}{\rm(i)} follows from Propositions \ref{prop-0730-01} and \ref{prop-0504-03}.
Theorem \ref{thm-0324-01}{\rm(ii)} is proved in Proposition \ref{prop-0512-02}, and Theorem \ref{thm-0324-01}{\rm(iii)} follows from Proposition \ref{prop-0324-02} and Lemma \ref{lem-0828-01}.
Theorem \ref{thm-0324-01}{\rm(iv)} follows from Propositions \ref{prop-0421-02} and \ref{prop-0731-11}.
Figure~\ref{260007} illustrates Theorem \ref{thm-0324-01} for \(n=3\).

We now state the second main theorem concerning the separation points of the normal velocity component \(w_n\) with normal boundary data.

\begin{theorem}\label{thm-0507-01}\textbf{(Normal velocity component with normal boundary data)}
	Let \(n\geq2\) and let \(w\) solve the Stokes system \eqref{0308-01} with boundary data \eqref{0324-13}, where \(j=n\) and \(\psi\) satisfies {\rm\textbf{(S)}}.
	Then the following hold.
	\begin{enumerate}
		\item \textbf{(Nonexistence and existence)} Suppose that \(\phi\) satisfies {\rm\textbf{(T1)}}.
		If \(0<t<1\), \(\phi\) is nondecreasing on \([0,t]\), and \(\phi(t)>0\), then \((x',t)\) is not a separation point of \(w_n\) whenever \(|x'|\geq1+C\sqrt{t}\), for some \(C=C(n)>0\).
		
		On the other hand, if \(\phi\not\equiv0\), then there is a separation point \((x',\widetilde{t_n^*}(x'))\) of \(w_n\) with \(\widetilde{t_n^*}(x')>1\) for sufficiently large \(|x'|\).
		If \(\phi\) also satisfies {\rm\textbf{(T2)}}, then there is a separation point \((x',t_n^*(x'))\) of \(w_n\) with \(0<t_n^*(x')<1\) for sufficiently large \(|x'|\).
		\item\label{item-0507-12} \textbf{(Temporal vanishing rates)} Suppose that \(\phi\) satisfies {\rm\textbf{(T1)}} and \(\phi\not\equiv0\).
		Then, for \(|x'|\) sufficiently large and \(\widetilde h:=|t-\widetilde{t_n^*}(x')|>0\) sufficiently small,
		\begin{align*}
		    c_n M_0\|\phi\|_{L^1(0,1)}|x'|^{-n-5}\widetilde h
		    \leq |\partial_{x_n}^2w_n(x',0,t)|
		    \leq C_n M_0\|\phi\|_{L^1(0,1)}|x'|^{-n-5}\widetilde h,
		\end{align*}
		where \(0<c_n\leq C_n\) depend only on \(n\).

		Suppose, in addition, that \(\phi\) satisfies {\rm\textbf{(T2)}}. Let \(\beta\in(0,1)\) be the H\"older exponent of \(\phi'\), and let \(t_n^*=t_n^*(x')\) be given in {\rm(i)}.
		Then, for \(h:=|t-t_n^*|>0\) sufficiently small,
		\begin{align*}
		    |\partial_{x_n}^2w_n(x',0,t)|\leq\frac{C|x'|^2}{(|x'|-1)^{n+2}}
		    \begin{cases}
		        h^{\beta+1/2}, & \beta<1/2,\\
		        h|\ln h|, & \beta=1/2,\\
		        h^\sigma, & \beta>1/2,
		    \end{cases}
		\end{align*}
		where \(C>0\) is independent of \(x'\) and, for \(\beta>1/2\), \(\sigma=1\) if \(\partial_t\partial_{x_n}^2w_n(x',0,t_n^*)\ne0\) and \(\sigma=\beta+1/2\) otherwise.
		Moreover, if \(0<\beta<\frac12\) and \eqref{260005} holds at \(t_n^*\), then
		\begin{align*}
			\lim_{\varepsilon\to0^+}\varepsilon^{-\beta-\frac12} \partial_{x_n}^2w_n(x',0,t_n^*+\varepsilon) =-C_\beta c_0\int_\Sigma\frac{\psi(y')}{|x'-y'|^n}\,dy',
		\end{align*}
		and
		\begin{align*}
			\lim_{\varepsilon\to0^-}|\varepsilon|^{-\beta-\frac12} \partial_{x_n}^2w_n(x',0,t_n^*+\varepsilon) =\widetilde C_\beta c_0\int_\Sigma\frac{\psi(y')}{|x'-y'|^n}\,dy',
		\end{align*}
		where \(C_\beta\) and \(\widetilde C_\beta\) are the constants in Theorem~\ref{thm-0324-01}{\rm(ii)}.
		\item\label{item-0507-14}\textbf{(Uniqueness)} Suppose that \(\phi\) satisfies {\rm\textbf{(T1)}} and \(\phi\not\equiv0\). Then \(\widetilde{t_n^*}(x')\) given in {\rm(i)} is the unique zero of \(\partial_{x_n}^2w_n(x',0,\cdot)\) in \((1,\infty)\) for sufficiently large \(|x'|\).
		Suppose, in addition, that \(\phi\) satisfies {\rm\textbf{(T2)}} and is log-concave on \((0,1)\). Then \(t_n^*(x')\) is the unique zero of \(\partial_{x_n}^2w_n(x',0,\cdot)\) in \((0,1)\) for sufficiently large \(|x'|\). Moreover,
		\begin{align*}
			\sgn\bigl(\partial_{x_n}^2w_n(x',0,t)\bigr) =\begin{cases} 1,& 0<t<t_n^*(x'),\\
				-1, & t_n^*(x')<t<\widetilde{t_n^*}(x'), \\
				1,& t>\widetilde{t_n^*}(x').
			\end{cases}
		\end{align*}
		\item\label{item-0507-13}\textbf{(Monotonicity and asymptotics)} 		
		Suppose that \(\phi\) satisfies {\rm\textbf{(T1)}} and \(\phi\not\equiv0\). Then \(r\mapsto\widetilde{t_n^*}(r\omega')\) is strictly increasing for sufficiently large \(r\), uniformly in \(\omega'\in\mathbb{S}^{n-2}\), and
		\begin{align}
			\widetilde{t_n^*}(x')&=\tau_0|x'|^2+O(|x'|),\label{0512-20}
		\end{align}
		for some \(\tau_0\) depending only on \(n\).

		Suppose, in addition, that \(\phi\) satisfies {\rm\textbf{(T2)}}--{\rm\textbf{(T4)}}. Then \(r\mapsto t_n^*(r\omega')\) is strictly increasing for sufficiently large \(r\), uniformly in \(\omega'\in \mathbb S^{n-2}\), and
		\begin{align}
			t_\infty-t_n^*(x')&=-\frac{3nT_1(t_\infty)}{T_0'(t_\infty)}|x'|^{-2}+O(|x'|^{-1-2\beta}).\label{0512-19}
		\end{align}
	\end{enumerate}
\end{theorem}
\begin{figure}[H]
	\centering
	\begin{tikzpicture}[x=0.95cm,y=0.52cm,>=latex,
			every node/.style={font=\footnotesize},
			declare function={small(\x)=0.68-0.42/(\x*\x);
					large(\x)=0.18*\x*\x+0.35;
					parabola(\x)=0.18*\x*\x;}]

		\draw[->] (0,0) -- (5.25,0) node[below right] {\(\lvert x'\rvert\)};
		\draw[->] (0,0) -- (0,4.75) node[left] {\(t\)};

		\draw[densely dotted,black!55] (0,1) -- (5.05,1);
		\node[above left,inner ysep=0pt] at (0,1) {\(1\)};
		\draw[dashed,black!65] (0,0.68) -- (5.05,0.68);
		\node[below left,inner ysep=0pt] at (0,0.68) {\(t_\infty\)};

		\draw[very thick]
		plot[domain=1.25:4.95,samples=80,smooth] (\x,{small(\x)});
		\node[right] at (5.55,0.50)
		{\(t=t_n^*(x')\)};

		\draw[densely dotted,black!55]
		plot[domain=2.2:4.9,samples=80,smooth] (\x,{parabola(\x)});
		\draw[very thick]
		plot[domain=2.0:4.85,samples=80,smooth] (\x,{large(\x)});
		\node[below right,inner sep=1pt] at (3.2,{large(3.2)})
		{\(t=\widetilde{t_n^*}(x')\)};
		\node[font=\normalsize] at (1.65,0.30) {\(+\)};
		\node[font=\Large] at (3.65,1.35) {\(-\)};
		\node[font=\Large] at (3.65,3.45) {\(+\)};

		\draw (1,-0.08) -- (1,0.08);
		\node[below] at (1,0) {\(1\)};

		\begin{scope}[shift={(5.55,4.25)},yscale=2]
			\draw[very thick] (0,0) -- (0.65,0);
			\node[right] at (0.78,0) {separation points};
			\draw[dashed,black!65] (0,-0.3) -- (0.65,-0.3);
			\node[right] at (0.78,-0.3) {\(t=t_\infty\)};
			\draw[densely dotted,black!55] (0,-0.6) -- (0.65,-0.6);
			\node[right] at (0.78,-0.6) {\(t=\tau_0|x'|^2\)};
		\end{scope}

	\end{tikzpicture}
	\caption{A schematic illustration of separation points of \(w_n\) and signs of \(\partial_{x_n}^2w_n\) along a fixed ray when \(n=3\)}
	\label{260008}
\end{figure}
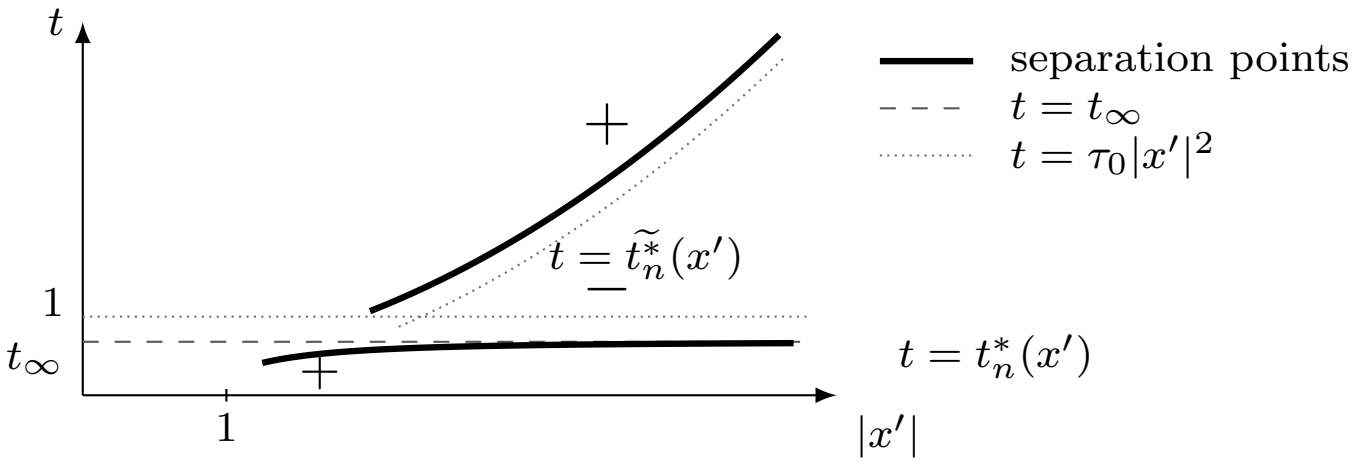

The nonexistence and existence statements in Theorem \ref{thm-0507-01}{\rm(i)} are proved in Propositions \ref{prop-0730-20} and \ref{prop-0504-02}, respectively.
Theorem \ref{thm-0507-01}{\rm(ii)} is proved in Proposition \ref{prop-0512-01}, and Theorem \ref{thm-0507-01}{\rm(iii)} in Propositions \ref{prop-0504-02} and \ref{prop-0505-01}, while Theorem \ref{thm-0507-01}{\rm(iv)} follows from Propositions \ref{prop-0514-01} and \ref{prop-0421-4}.
Figure~\ref{260008} illustrates Theorem \ref{thm-0507-01} for \(n=3\).

We now state the third main theorem concerning the separation points of the tangential velocity component in the direction of the tangential boundary data.
If \(n\geq3\), let \(\omega'\in\mathbb S^{n-2}\) satisfy \(\omega'\cdot\mathbf e_i'=0\).
Set
\begin{align}\label{260009}
	\mathcal S(x'):=\partial_{x_i}R_j'\psi(x').
\end{align}

\begin{theorem}
	\label{thm-260010}\textbf{(Tangential velocity component in the direction of the tangential boundary data)}
	Let \(n\geq2\), \(1\leq i=j\leq n-1\), and let \(w\) solve the Stokes system \eqref{0308-01} with boundary data \eqref{0324-13}, where \(\psi\) satisfies {\rm\textbf{(S)}}.
	Then, with \(\mathcal S\) given in \eqref{260009}, the following hold for \(|x'|>1\).
	\begin{enumerate}
		\item \textbf{(Nonexistence)} Suppose that \(\phi\) satisfies {\rm\textbf{(T1)}} and {\rm\textbf{(T2)}}.
		Then \((x',t)\) is not a separation point of \(w_i\) if one of the following holds.
		\begin{enumerate}
			\item\(\mathcal S(x')\leq0\) and \(t>0\).
			\item \(t\geq1\).
			\item \(n=2\) and \(t>0\).
			\item \(n\geq 3\), \(x'=x_i\mathbf{e}_i'\), and \(\psi\) is radial.
		\end{enumerate}

		\item \textbf{(Existence)} Let \(n\geq 3\) and suppose that \(\phi\) satisfies {\rm\textbf{(T1)}} and {\rm\textbf{(T2)}}.
		Then the following hold.
		\begin{enumerate}
			\item If \(x'\) satisfies \(\mathcal S(x')>0\), then there exists \(t_i^*(x')\in(0,1)\) such that \((x',t_i^*(x'))\) is a separation point of \(w_i\).
			If, in addition, \(\psi\) is radial and \(|x^{(i)}|>1\), then there exists \(t_i^*(x^{(i)})\in(0,1)\), depending only on \(|x^{(i)}|\), such that \((x^{(i)},t_i^*(x^{(i)}))\) is a separation point of \(w_i\).
			\item If \(\psi\) is radial, then, for every \(t_0\in(0,1)\), there exists \(R_{t_0}=R_{t_0}(n,\psi,\phi)>1\) such that, for every \(|x_i|>R_{t_0}\) and \(0<t\leq t_0\), there exists \(\rho_*(x_i,t)>0\) such that, for every \(\omega'\in\mathbb S^{n-2}\) with \(\omega'\cdot\mathbf e_i'=0\), \((x_i\mathbf e_i'+\rho_*(x_i,t)\omega',t)\) is a separation point of \(w_i\).
		\end{enumerate}

		\item \textbf{(Uniqueness)} Let \(n\geq 3\) and suppose that \(\phi\) satisfies {\rm\textbf{(T1)}} and {\rm\textbf{(T2)}}.
		If \(\psi\) is radial, then, for every \(t_0\in(0,1)\), there exists \(\overline{R}_{t_0}= \overline{R}_{t_0}(n,\psi,\phi)\geq R_{t_0}\) such that, if \(|x_i|>\overline{R}_{t_0}\) and \(0<t\leq t_0\), then \(\rho_*(x_i,t)\) given in {\rm(ii)(b)} is unique.
		If \(\phi\) is log-concave on \((0,1)\), then \(t_i^*(x')\) and \(t_i^*(x^{(i)})\) given in {\rm(ii)(a)} are unique for every \(|x'|>1\) and \(|x^{(i)}|>1\), respectively.

		\item \textbf{(Monotonicity and asymptotics)} Let \(n\geq 3\) and suppose that \(\phi\) satisfies {\rm\textbf{(T1)}} and {\rm\textbf{(T2)}}.
		If \(\psi\) is radial, then, for every \(t_0\in(0,1)\), there exists \(\widetilde{R}_{t_0}=\widetilde{R}_{t_0}(n,\psi,\phi)\geq \overline{R}_{t_0}\) such that, for every \(0<t\leq t_0\), \(\rho_*(x_i,t)\) given in {\rm(ii)(b)} is strictly increasing in \(|x_i|\) for \(|x_i|>\widetilde{R}_{t_0}\), and the following expansion holds uniformly for \(0<t\leq t_0\) as \(|x_i|\to\infty\):
		\begin{align*}
			\rho_*(x_i,t)=\sqrt{n-1}\,|x_i|+\frac{M_2}{2(n-1)^{3/2}M_0|x_i|}+O\bigl(|x_i|^{-3}\bigr).
		\end{align*}
		If, in addition, \(\phi\in C^N([0,1])\) for some \(N\geq1\), some derivative of \(\phi\) of order at most \(N\) is nonzero at each of \(0\) and \(1\), and \(k_1\) is the least positive integer such that \(\phi^{(k_1)}(1)\neq0\), then every \(t_i^*(x^{(i)})\) given in {\rm(ii)(a)} satisfies
		\begin{align*}
			\ln(1-t_i^*(x^{(i)}))=-\frac{|x^{(i)}|^2}{4k_1}+O(|x^{(i)}|).
		\end{align*}

		\item \textbf{(Limits of positivity regions)} Suppose that \(\phi\) satisfies {\rm\textbf{(T1)}} and {\rm\textbf{(T2)}}.
		For \(t>0\), the indicator functions of
		\begin{align*}
			\Omega_t:=\bigl\{x'\in\Sigma:\ |x'|>1,\ \partial_{x_n}w_i(x',0,t)>0\bigr\}, \qquad \Omega_0:=\bigl\{x'\in\Sigma:\ |x'|>1,\ \mathcal S(x')\leq0\bigr\}
		\end{align*}
		satisfy
		\begin{align*}
			\mathbf 1_{\Omega_t}\to\mathbf 1_{\Omega_0}\quad\text{pointwise on } \{|x'|>1\}\text{ as }t\to0^+, \qquad \mathbf 1_{\Omega_t}\to1\quad\text{pointwise on } \{|x'|>1\}\text{ as }t\to1^-.
		\end{align*}
	\end{enumerate}
\end{theorem}
\begin{figure}[H]
	\centering
	\begin{tikzpicture}[scale=0.72,>=latex,every node/.style={font=\footnotesize}]
		\begin{scope}[shift={(-4.2,0)}]
			\fill[black!12] (0,0)--(2.8,2.35)--(2.8,-2.35)--cycle;
			\fill[black!12] (0,0)--(-2.8,2.35)--(-2.8,-2.35)--cycle;
			\fill[white] (0,0) circle (1.35);
			\draw[densely dotted,black!65] (0,0) circle (1.35);
			\draw[dashed,black!60] (-2.8,2.35)--(2.8,-2.35);
			\draw[dashed,black!60] (-2.8,-2.35)--(2.8,2.35);
			\draw[->] (-3,0)--(3.05,0) node[below] {\(x_1\)};
			\draw[->] (0,-2.75)--(0,2.85) node[left] {\(x_2\)};
			\fill[black!18] (0,0) circle (1);
			\draw (0,0) circle (1);

			\boundarydataarrows

			\node[font=\scriptsize] at (0,0.3) {$\overline{B_1'}$};
			\node[below,font=\scriptsize] at (1.35,-0.05) {\(r_\varepsilon\)};
			\foreach \sx in {-1,1}{\foreach \sy in {-1,1}{
							\node[font=\Large] at ({\sx*2.15},{\sy*0.58}) {\(+\)};}}
			\node[font=\scriptsize] at (0,-3.15) {\(|x_2|\leq(\sqrt{n-1}-\varepsilon)|x_1|\)};
		\end{scope}
	\end{tikzpicture}
	\caption{A schematic illustration, for \(n=3\), of the region \(\{|x'|\geq r_\varepsilon,\ |x^{(i)}|\leq(\sqrt{n-1}-\varepsilon)|x_i|\}\) in which \(w_i\) has no separation point, with \(r_\varepsilon=n/\varepsilon\)}
	\label{260011}
\end{figure}
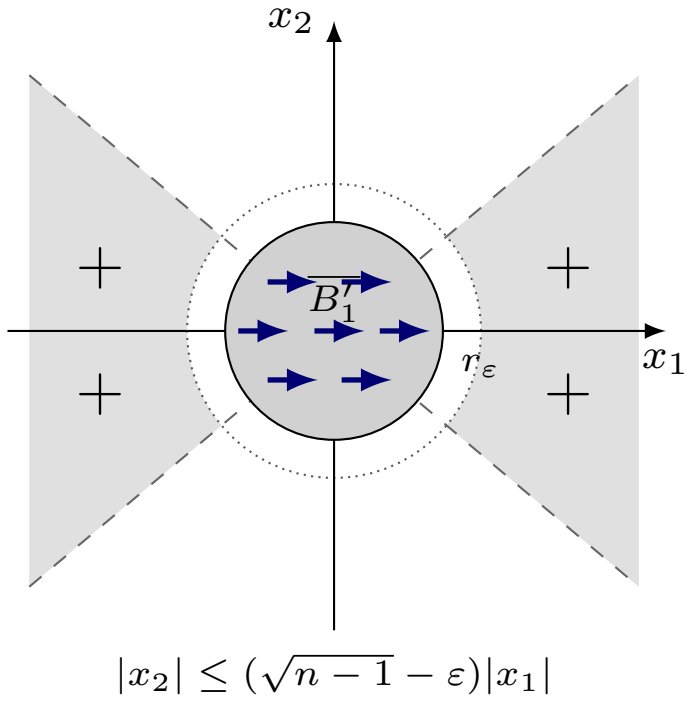

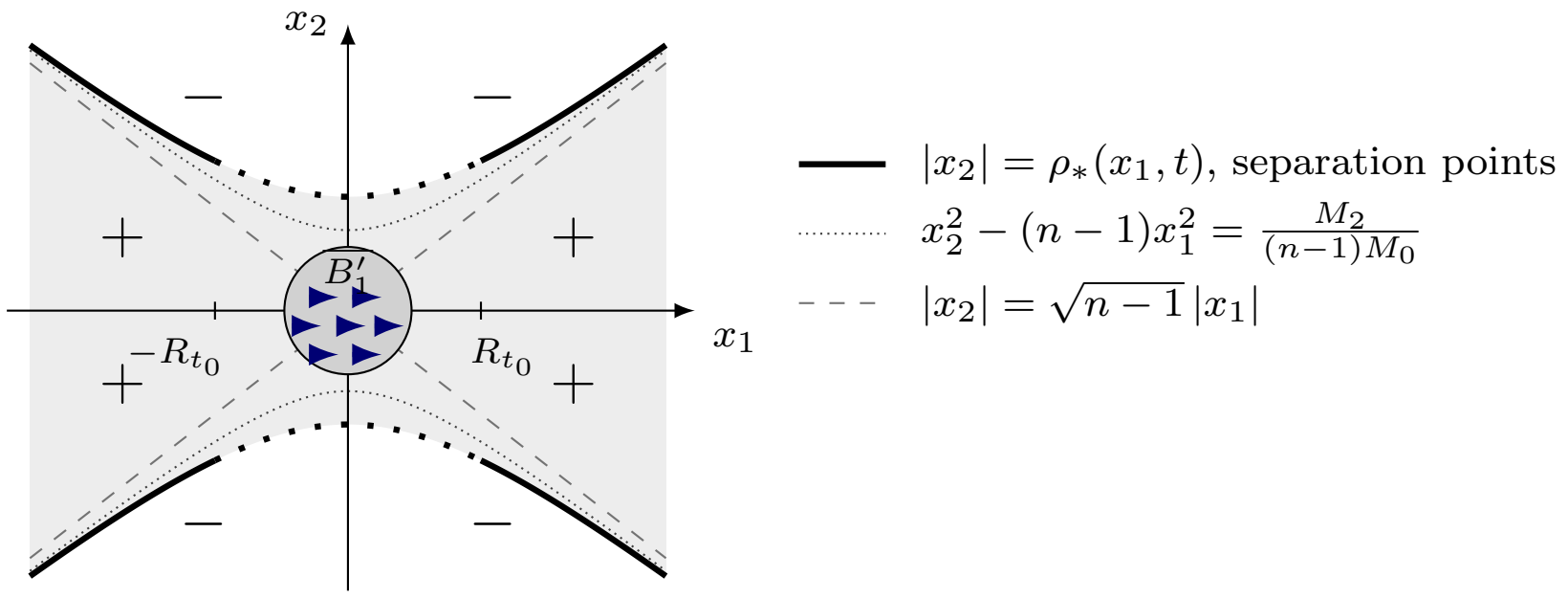
\begin{figure}[H]
	\centering
	\begin{tikzpicture}[scale=0.8,yscale=0.55,>=latex,
			every node/.style={font=\footnotesize},
			declare function={rho(\x)=sqrt(2*\x*\x+1.6+1.6/(1+0.18*\x*\x));
					hyp(\x)=sqrt(2*\x*\x+1.6);}]

		\begin{scope}
			\clip (-2.75,-4.25) rectangle (2.75,4.25);
			\fill[black!7]
			plot[domain=-2.85:2.85,samples=90,smooth] (\x,{-rho(\x)})
			-- plot[domain=2.85:-2.85,samples=90,smooth] (\x,{rho(\x)})
			-- cycle;
		\end{scope}

		\draw[dashed,black!55] (-2.75,3.889) -- (2.75,-3.889);
		\draw[dashed,black!55] (-2.75,-3.889) -- (2.75,3.889);

		\draw[densely dotted,black!75]
		plot[domain=-2.75:2.75,samples=70,smooth] (\x,{hyp(\x)});
		\draw[densely dotted,black!75]
		plot[domain=-2.75:2.75,samples=70,smooth] (\x,{-hyp(\x)});

		\foreach \s in {1,-1}{
				\draw[very thick] plot[domain=1.15:2.75,samples=55,smooth] (\x,{\s*rho(\x)});
				\draw[very thick] plot[domain=-2.75:-1.15,samples=55,smooth] (\x,{\s*rho(\x)});
				\draw[very thick,loosely dotted]
				plot[domain=-1.15:1.15,samples=40,smooth] (\x,{\s*rho(\x)});
			}

		\draw[->] (-2.95,0) -- (3.0,0) node[below right] {\(x_1\)};
		\draw[->] (0,-4.4) -- (0,4.5) node[left] {\(x_2\)};

		\fill[black!18] (0,0) ellipse [x radius=.55,y radius=1];
		\draw (0,0) ellipse [x radius=.55,y radius=1];

		\begin{scope}[shift={(0,-0.24)},xscale=.55]
			\boundarydataarrows
		\end{scope}

		\node[font=\scriptsize] at (0,0.62) {$\overline{B_1'}$};

		\foreach \s in {1,-1}{\draw (\s*1.15,-0.14) -- (\s*1.15,0.14);}
		\node[below right,font=\scriptsize] at (0.9,-0.15) {\(R_{t_0}\)};
		\node[below left,font=\scriptsize] at (-0.9,-0.15) {\(-R_{t_0}\)};

		\foreach \sx in {-1,1}{\foreach \sy in {-1,1}{
						\node[font=\Large] at ({\sx*1.25},{\sy*3.35}) {\(-\)};
						\node[font=\Large] at ({\sx*1.95},{\sy*1.15}) {\(+\)};}}

		\begin{scope}[shift={(3.9,2.3)},yscale=1.45455]
			\draw[very thick] (0,0) -- (0.75,0);
			\node[right] at (0.9,0) {\(|x_2|=\rho_*(x_1,t)\), separation points};
			\draw[densely dotted,black!75] (0,-0.75) -- (0.75,-0.75);
			\node[right] at (0.9,-0.75) {\(x_2^2-(n-1)x_1^2=\frac{M_2}{(n-1)M_0}\)};
			\draw[dashed,black!55] (0,-1.5) -- (0.75,-1.5);
			\node[right] at (0.9,-1.5) {\(|x_2|=\sqrt{n-1}\,|x_1|\)};
		\end{scope}

	\end{tikzpicture}
	\caption{A schematic illustration of separation points of \(w_1\) for radial \(\psi\), \(n=3\), and fixed \(t\in(0,1)\)}
	\label{260012}
\end{figure}

Theorem \ref{thm-260010}{\rm(i)} is proved in Proposition \ref{prop-260061}.
Theorem \ref{thm-260010}{\rm(ii)(a)} is proved in Proposition \ref{prop-260072}{\rm(i)} and {\rm(iii)}, and Theorem \ref{thm-260010}{\rm(ii)(b)} in Proposition \ref{prop-260072}{\rm(ii)}.
Theorem \ref{thm-260010}{\rm(iii)} is proved in Proposition \ref{prop-260081}, {\rm(iv)} in Propositions \ref{prop-260088} and \ref{prop-260086}, and {\rm(v)} in Proposition \ref{prop-0703-08}.
Figures~\ref{260011} and~\ref{260012} illustrate Theorem \ref{thm-260010} for \(n=3\).

\begin{remark}
The conditions \(\mathcal S(x')\leq0\) and \(\mathcal S(x')>0\) in Theorem \ref{thm-260010}{\rm(i)(a)} and {\rm(ii)(a)} hold on regions independent of \(\psi\): since \(\max_{0\leq\alpha\leq1}(\sqrt{1-\alpha^2}+\sqrt{n-1}\,\alpha)=\sqrt n\),
	\begin{align}\label{260070}
	\mathcal S(x')\leq0 \quad\text{if}\quad \sqrt{n-1}\,|x_i|-|x^{(i)}|\geq\sqrt n, \qquad \mathcal S(x')>0 \quad\text{if}\quad |x^{(i)}|-\sqrt{n-1}\,|x_i|\geq\sqrt n.
	\end{align}
In particular, if \(0<\varepsilon<\sqrt{n-1}\) and \(|x^{(i)}|\leq(\sqrt{n-1}-\varepsilon)|x_i|\), then \(|x'|\leq\sqrt n\,|x_i|\) and \(\sqrt{n-1}\,|x_i|-|x^{(i)}|\geq\varepsilon|x_i|\), so \(\mathcal S(x')\leq0\) on the region in Figure~\ref{260011} with \(r_\varepsilon=n/\varepsilon\).
\end{remark}

We now state the fourth main theorem concerning the separation points of a tangential velocity component in a direction different from that of the tangential boundary data.
We take \(i=2\) and \(j=1\).
We write \(x^\flat=x^{(2)}=x'-x_2\mathbf e_2'\in(\mathbf e_2')^\perp\), with the same convention for \(y^\flat\). The map \(\sigma\) exchanges the first two coordinates.
Then
\begin{align}\label{260013}
	x'=(x_1,x'')=x_2\mathbf e_2'+x^\flat.
\end{align}

\begin{theorem}
	\label{thm-260014}\textbf{(Tangential velocity component in a direction different from that of the tangential boundary data)}
	Let \(n\geq3\) and let \(w\) solve the Stokes system \eqref{0308-01} with boundary data \eqref{0324-13}, where \(j=1\), \(\psi\) satisfies {\rm\textbf{(S)}}, and \(\phi\) satisfies {\rm\textbf{(T1)}} and {\rm\textbf{(T2)}}.
	Then the following hold for \(w_2\).

	\begin{enumerate}
		\item \textbf{(Nonexistence)} If \(|x_1|\geq 1\), \(|x_2|\geq 1\), and \(t>0\), then \((x',t)\) is not a separation point of \(w_2\).

		\item \textbf{(Existence)} For each fixed \(x''\in(\mathbf e_1')^\perp\) with \(|x_2|>1\) and every \(0<t\leq1\), there exists \(x_1^*(x'')\in(-1,1)\), independent of \(t\), such that \((x_1^*(x'')\mathbf{e}_1'+x'',t)\) is a separation point of \(w_2\).
		
		Similarly, for each fixed \(x^\flat\in(\mathbf e_2')^\perp\) with \(|x_1|>1\) and every \(0<t\leq1\), there exists \(x_2^*(x^\flat)\in(-1,1)\), independent of  \(t\), such that \((x_2^*(x^\flat)\mathbf e_2'+x^\flat,t)\) is a separation point of \(w_2\).

		If, in addition, \(\psi\) is radial and if either \(0<t<1\) and \(|x'|>1\), or \(t\geq1\) and \(|x'|\geq\sqrt2\), then \((x',t)\) is a separation point of \(w_2\) if and only if \(x_1x_2=0\).
		Moreover, if \(|x'|>1\) and \(x_1x_2=0\), then \((x',1)\) is a separation point of \(w_2\).

		\item \textbf{(Uniqueness)} If \(|x_2|>1\) and \(|x''|>1+2\sqrt{n+1}\), then, for every \(0<t<1\), \((x_1^*(x'')\mathbf{e}_1'+x'',t)\) is the unique separation point of \(w_2\) on the line \(\{x_1\mathbf e_1'+x'':x_1\in\R\}\).
		Moreover,
		\begin{align*}
			\sgn\bigl(\partial_{x_n}w_2(x',0,t)\bigr) =\sgn(x_2)\,\sgn\bigl(x_1-x_1^*(x'')\bigr).
		\end{align*}
		Similarly, if \(|x_1|>1\) and \(|x^\flat|>1+2\sqrt{n+1}\), then, for every \(0<t<1\), \((x_2^*(x^\flat)\mathbf e_2'+x^\flat,t)\) is the unique separation point of \(w_2\) on the line \(\{x_2\mathbf e_2'+x^\flat:x_2\in\R\}\).
		Moreover,
		\begin{align*}
		\sgn\bigl(\partial_{x_n}w_2(x',0,t)\bigr) =\sgn(x_1)\sgn\bigl(x_2-x_2^*(x^\flat)\bigr).
		\end{align*}
		\item \textbf{(Monotonicity and asymptotics)} If
		\begin{align*}
			\displaystyle \int_\Sigma\left(y_1-\frac{m_1}{M_0}\right)y_2\psi(y')\,dy'\neq0,
		\end{align*}
		then \(x_1^*(x_2\mathbf e_2')\) and \(x_2^*(x_1\mathbf e_1')\) are either both strictly increasing or both strictly decreasing on each of \((R_1,\infty)\) and \((-\infty,-R_1)\), for some sufficiently large \(R_1\).
		For each fixed \(\omega''\in\mathbb S^{n-3}\) with \(\omega_2''\neq0\), as \(r\to\infty\),
		\begin{align*}
			x_1^*(r\omega'')=\frac{m_1}{M_0} +O(r^{-1}), \qquad x_2^*(r\sigma\omega'')=\frac{m_2}{M_0} +O(r^{-1}).
		\end{align*}
	\end{enumerate}
\end{theorem}

\begin{figure}[H]
	\centering
	\begin{tikzpicture}[scale=1.0,yscale=0.62,>=latex,
			every node/.style={font=\footnotesize}]

		\fill[black!4] (-1,-3.2) rectangle (1,3.2);
		\fill[black!4] (-3.2,-1) rectangle (3.2,1);

		\fill[pattern=north east lines,pattern color=black!45,even odd rule]
		(-1,-1) rectangle (1,1)(0,0) circle (1);

		\draw[densely dotted,black!60] (-1,-3.2) -- (-1,3.2);
		\draw[densely dotted,black!60] ( 1,-3.2) -- ( 1,3.2);
		\draw[densely dotted,black!60] (-3.2,-1) -- (3.2,-1);
		\draw[densely dotted,black!60] (-3.2, 1) -- (3.2, 1);

		\draw[dashed,black!60] (0.35,-3.2) -- (0.35,3.2);
		\draw[dashed,black!60] (-3.2,-0.2) -- (3.2,-0.2);

		\draw[->] (-3.45,0) -- (3.55,0) node[below right] {\(x_1\)};
		\draw[->] (0,-3.45) -- (0,3.6) node[left] {\(x_2\)};

		\draw[very thick]
		plot[domain=1.03:3.2,samples=60,smooth] ({0.35+0.48/\x},\x);
		\draw[very thick]
		plot[domain=-3.2:-1.03,samples=60,smooth] ({0.35+0.48/\x},\x);

		\draw[very thick,black!55]
		plot[domain=1.03:3.2,samples=60,smooth] (\x,{-0.2+0.48/\x});
		\draw[very thick,black!55]
		plot[domain=-3.2:-1.03,samples=60,smooth] (\x,{-0.2+0.48/\x});

		\fill[black!18] (0,0) ellipse [x radius=.62,y radius=1];
		\draw (0,0) ellipse [x radius=.62,y radius=1];

		\begin{scope}[shift={(0,-0.17)},xscale=.62]
			\boundarydataarrows
		\end{scope}

		\node[font=\scriptsize] at (0,0.65) {$\overline{B_1'}$};

		\node[font=\Large] at (2.35,2.35) {\(+\)};
		\node[font=\Large] at (2.35,0.55) {\(+\)};
		\node[font=\Large] at (0.80,2.35) {\(+\)};
		\node[font=\Large] at (-2.35,2.35) {\(-\)};
		\node[font=\Large] at (-0.55,2.35) {\(-\)};
		\node[font=\Large] at (-2.35,0.55) {\(-\)};
		\node[font=\Large] at (-2.35,-2.35) {\(+\)};

		\node[font=\Large] at (-2.35,-0.75) {\(+\)};
		\node[font=\Large] at (2.35,-0.55) {\(-\)};
		\node[font=\Large] at (2.35,-2.35) {\(-\)};
		\node[font=\Large] at (0.80,-2.35) {\(-\)};
		\node[font=\Large] at (-0.55,-2.35) {\(+\)};

		\node[above right,font=\scriptsize,inner sep=1pt] at (1.08,0.35){\(1\)};
		\node[above left,font=\scriptsize,inner sep=1pt] at (-1.08,0.35) {\(-1\)};
		\node[above left,font=\scriptsize]at (0,1){\(1\)};
		\node[below left,font=\scriptsize]at (0,-1) {\(-1\)};

		\begin{scope}[shift={(4.35,2.2)}]
			\draw[very thick] (0,0) -- (0.75,0);
			\node[right] at (0.9,0) {\(x_1=x_1^*(x'')\) on \(\{|x_2|>1\}\)};
			\draw[very thick,black!55] (0,-0.7) -- (0.75,-0.7);
			\node[right] at (0.9,-0.7) {\(x_2=x_2^*(x^\flat)\) on \(\{|x_1|>1\}\)};
			\draw[dashed,black!60] (0,-1.4) -- (0.75,-1.4);
			\node[right] at (0.9,-1.4) {\(x_1=m_1/M_0\), \(x_2=m_2/M_0\)};
			\draw[densely dotted,black!60] (0,-2.1) -- (0.75,-2.1);
			\node[right] at (0.9,-2.1) {\(|x_1|=1\), \ \(|x_2|=1\)};
		\end{scope}

	\end{tikzpicture}
	\caption{A schematic illustration of separation points of \(w_2\) for \(n=3\) and fixed \(t\in(0,1)\)}
	\label{260015}
\end{figure}
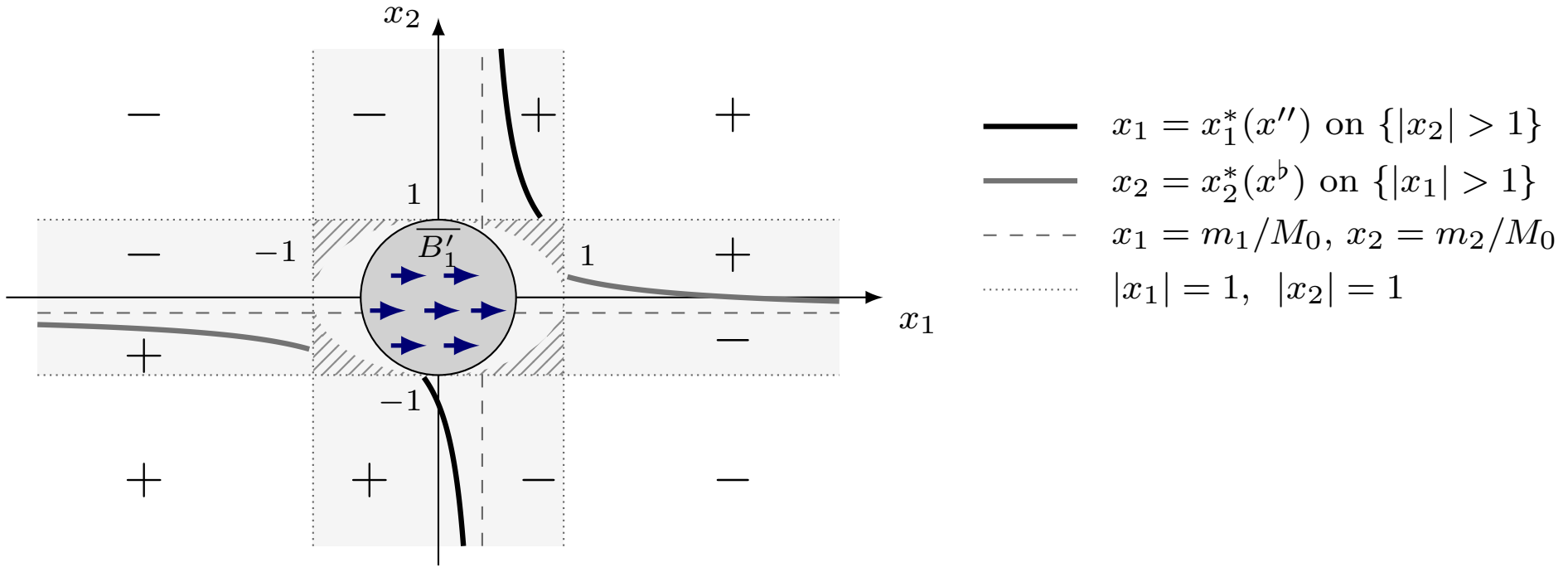

\begin{figure}[H]
	\centering
	\begin{tikzpicture}[
			scale=0.92,
			xscale=0.75,
			yscale=0.50,
			>=stealth,
			font=\small,
			declare function={
					rho(\x)=sqrt(2*\x*\x+1.6+1.6/(1+0.18*\x*\x));
					xone(\y)=0.35+0.48/\y;
					xtwo(\x)=-0.2+0.48/\x;
				}
		]

		\def\Rt{1.15}

		\begin{scope}
			\clip (-3.15,-4.35) rectangle (3.15,4.35);

			\fill[black!6]
			plot[domain=-2.95:2.95,samples=120,smooth]
			(\x,{-rho(\x)})
			--
			plot[domain=2.95:-2.95,samples=120,smooth]
			(\x,{rho(\x)})
			-- cycle;
		\end{scope}

		\begin{scope}[yscale=1.5]
		\fill[
			pattern=north east lines,
			pattern color=black!45,
			even odd rule
		]
		(-1,-1) rectangle (1,1)
		(0,0) circle (1);

		\draw[densely dotted,black!60,line width=.65pt]
		(-1,-1) rectangle (1,1);

		\fill[black!16] (0,0) circle (1);
		\draw[line width=.75pt] (0,0) circle (1);
		\end{scope}

		\draw[->,line width=.8pt]
		(-3.25,0)--(3.32,0)
		node[below right] {$x_1$};

		\draw[->,line width=.8pt]
		(0,-4.48)--(0,4.55)
		node[left] {$x_2$};

		\draw[line width=1.6pt]
		plot[domain=\Rt:2.95,samples=80,smooth]
		(\x,{rho(\x)});

		\draw[line width=1.6pt]
		plot[domain=-2.95:-\Rt,samples=80,smooth]
		(\x,{rho(\x)});

		\draw[line width=1.6pt,loosely dotted,black!65]
		plot[domain=-\Rt:\Rt,samples=50,smooth]
		(\x,{rho(\x)});

		\draw[line width=1.6pt]
		plot[domain=\Rt:2.95,samples=80,smooth]
		(\x,{-rho(\x)});

		\draw[line width=1.6pt]
		plot[domain=-2.95:-\Rt,samples=80,smooth]
		(\x,{-rho(\x)});

		\draw[line width=1.6pt,loosely dotted,black!65]
		plot[domain=-\Rt:\Rt,samples=50,smooth]
		(\x,{-rho(\x)});

		\draw[line width=1.6pt,black!55]
		plot[domain=1.03:4.25,samples=120,smooth]
		({xone(\x)},\x);

		\draw[line width=1.6pt,black!55]
		plot[domain=-4.25:-1.03,samples=120,smooth]
		({xone(\x)},\x);

		\draw[line width=1.6pt,black!55]
		plot[domain=1.03:3.05,samples=90,smooth]
		(\x,{xtwo(\x)});

		\draw[line width=1.6pt,black!55]
		plot[domain=-3.05:-1.03,samples=90,smooth]
		(\x,{xtwo(\x)});

		\def\arrowlength{0.23}

		\foreach \i in {-4.5,-3.5,-2.5,-1.5,-0.5,0.5,1.5,2.5,3.5,4.5}{
				\foreach \j in {0,...,16}{

						\pgfmathsetmacro{\xx}{0.56*\i}
						\pgfmathsetmacro{\yy}{-4.08+0.51*\j}
						\pgfmathtruncatemacro{\insupp}{
							ifthenelse(sqrt(\xx*\xx+(\yy/1.5)*(\yy/1.5))<=1,1,0)
						}

						\pgfmathtruncatemacro{\inC}{
							ifthenelse(
							abs(\xx)<=1,
							ifthenelse(
							abs(\yy)<=1.5,
							ifthenelse(sqrt(\xx*\xx+(\yy/1.5)*(\yy/1.5))>1,1,0),
							0
							),
							0
							)
						}

						\ifnum\insupp=0
							\ifnum\inC=0

								\pgfmathsetmacro{\rhoSq}{
									2*\xx*\xx
									+1.6
									+1.6/(1+0.18*\xx*\xx)
								}

								\pgfmathsetmacro{\uone}{
									(\rhoSq-\yy*\yy)/4.0
								}

								\pgfmathsetmacro{\Avec}{
									\xx*\yy - 0.35*\yy - 0.48
								}

								\pgfmathsetmacro{\Bvec}{
									\xx*\yy + 0.20*\xx - 0.48
								}

								\pgfmathsetmacro{\blendraw}{
									max(0,min(1,(abs(\yy)-0.96)/0.06))
								}

								\pgfmathsetmacro{\blend}{
									\blendraw*\blendraw*(3-2*\blendraw)
								}

								\pgfmathsetmacro{\utwo}{
									(
									\blend*\Avec
									+(1-\blend)*\Bvec
									)/4.0
								}

								\pgfmathsetmacro{\vmag}{
									sqrt(\uone*\uone+\utwo*\utwo)+0.000001
								}

								\pgfmathsetmacro{\dx}{
									\arrowlength*\uone/\vmag
								}

								\pgfmathsetmacro{\dy}{
									\arrowlength*\utwo/\vmag
								}

								\pgfmathsetmacro{\xa}{\xx-\dx/2}
								\pgfmathsetmacro{\ya}{\yy-\dy/2}
								\pgfmathsetmacro{\xb}{\xx+\dx/2}
								\pgfmathsetmacro{\yb}{\yy+\dy/2}

								\draw[
									->,
									line width=.55pt
								]
								(\xa,\ya)--(\xb,\yb);

							\fi
						\fi
					}
			}

		\begin{scope}[yscale=1.5]
		\foreach \X/\Y in {
				-0.38/0.45,
				0.30/0.45,
				-0.65/0.00,
				0.05/0.00,
				0.65/0.00,
				-0.38/-0.60,
				0.30/-0.60
			}{
				\draw[
					->,
					line width=1.00pt,
					blue!45!black
				]
				({\X-0.23},\Y)--++(0.46,0);
			}

		\node[below right,font=\scriptsize,fill=white,inner sep=.5pt]
		at (1,0) {$1$};

		\node[below left,font=\scriptsize,fill=white,inner sep=.5pt]
		at (-1,0) {$-1$};

		\node[above,font=\scriptsize,fill=white,inner sep=.5pt]
		at (0,1) {$1$};

		\node[below,font=\scriptsize,fill=white,inner sep=.5pt]
		at (0,-1) {$-1$};

		\node[font=\scriptsize,fill=black!16,inner sep=0pt] at (0,0.73) {$\overline{B_1'}$};
		\node[font=\scriptsize,fill=black!16,inner sep=0pt] at (0,-0.30) {$O$};
		\end{scope}

	\end{tikzpicture}

	\caption{A schematic illustration of the velocity field near the boundary for tangential boundary data}
	\label{fig:tangential_boundary_combined}
\end{figure}
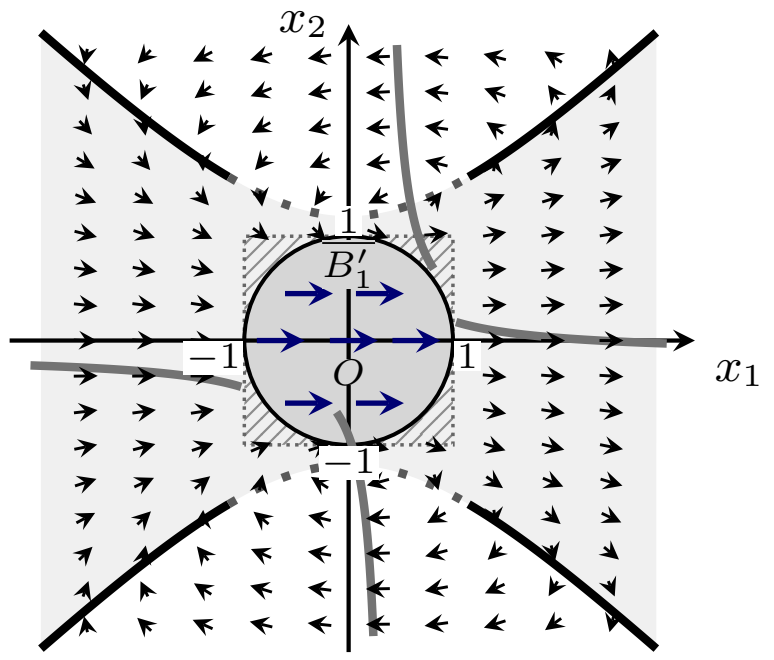

Theorem \ref{thm-260014}{\rm(i)} is proved in Proposition \ref{prop-260097}.
The statements in Theorem \ref{thm-260014}{\rm(ii)} for general and for radial \(\psi\) are proved in Propositions \ref{prop-260101} and \ref{prop-0708-10}, respectively.
Theorem \ref{thm-260014}{\rm(iii)} is proved in Proposition \ref{prop-260103}, while Theorem \ref{thm-260014}{\rm(iv)} is proved in Proposition \ref{prop-260107}.
Figure~\ref{260015} illustrates Theorem \ref{thm-260014} for \(n=3\), and Figure~\ref{fig:tangential_boundary_combined} sketches the velocity field near the boundary for tangential boundary data.

We now state the fifth main theorem concerning the separation points of the normal velocity component with tangential boundary data.
Since the divergence-free condition gives \(\partial_{x_n}w_n(x',0,t)=-\phi(t)\partial_{x_j}\psi(x')=0\) for \(|x'|>1\), we study \(\partial_{x_n}^2w_n(x',0,t)\) instead.
\begin{theorem}
	\label{thm-260017}\textbf{(Normal velocity component with tangential boundary data)}
		Let \(n\geq2\) and let \(w\) solve the Stokes system \eqref{0308-01} with boundary data \eqref{0324-13}, where \(j=1\), \(\psi\) satisfies {\rm\textbf{(S)}}, and \(\phi\) satisfies {\rm\textbf{(T1)}} and {\rm\textbf{(T2)}}.
	Then the following hold for \(|x'|>1\).
	\begin{enumerate}
		\item \textbf{(Nonexistence)} If \(|x_1|\geq 1\) and \(t>0\), then \((x',t)\) is not a separation point of \(w_n\).
		\item \textbf{(Existence)} Let \(n\geq 3\).
		For each \(t>0\) and each fixed \(x''\in(\mathbf e_1')^\perp\) with \(|x''|>1\), there exists  \(x_1^*(x'',t)\in(-1,1)\) such that \((x_1^*(x'',t)\mathbf e_1'+x'',t)\) is a separation point of \(w_n\).
		If, in addition, \(\psi\) is radial, then \((x',t)\) is a separation point of \(w_n\) if and only if \(x_1=0\).

		\item \textbf{(Uniqueness)} Let \(n\geq 3\).
		If \(|x''|\) is sufficiently large, then, for every \(t>0\), \((x_1^*(x'',t)\mathbf e_1'+x'',t)\) is the unique separation point of \(w_n\) on the line \(\{x_1\mathbf e_1'+x'':x_1\in\R\}\).
		Moreover,
		\begin{align*}
			\sgn(\partial_{x_n}^2w_n(x',0,t)) =\sgn(x_1-x_1^*(x'',t)).
		\end{align*}
		\item \textbf{(Asymptotics)} Let \(n\geq 3\).
		\begin{enumerate}
			\item There exists \(R\geq2\), depending only on \(n,\psi,\phi\), such that, for \(|x''|>R\) and every \(0<t<1\),
			\begin{align*}
				\left|x_1^*(x'',t)-\frac{m_1}{M_0}\right| \leq\frac{n+2}{|x''|} +C_n\frac{\phi(t_1)}{\phi(\max\{t,t_1\})}e^{-c|x''|^2},
			\end{align*}
			and
			\begin{align*}
				\left|x_1^*(x'',t)-\frac{m_1}{M_0} -\frac{n+2}{M_0|x''|^2}x''\cdot\int_\Sigma \left(y_1-\frac{m_1}{M_0}\right)y''\psi(y')\,dy'\right| \leq\frac{C_{n,\psi}}{|x''|^2} +C_n\frac{\phi(t_1)}{\phi(\max\{t,t_1\})}e^{-c|x''|^2},
			\end{align*}
			where \(c>0\) is an absolute constant.
			\item Let \(t\geq1\) and \(\omega''\in\mathbb S^{n-3}\), and set \(F_{\omega''}:=\left\{y'\in\supp\psi:\omega''\cdot y'' =\max_{z'\in\supp\psi}\omega''\cdot z''\right\}\).
			If \(\{y_1:y'\in F_{\omega''}\}=\{p_1\}\) for some \(p_1\), then
			\begin{align*}
				\lim_{r\to\infty}x_1^*(r\omega'',t)=p_1.
			\end{align*}
			\item For \(|x''|>R\) and every \(t\geq1\),
			\begin{align*}
				\left|x_1^*(x'',t)-\frac{m_1}{M_0}\right| \leq2\min\left\{1,\frac{1+|x''|}{t}\right\}.
			\end{align*}
		\end{enumerate}
	\end{enumerate}
\end{theorem}

For an example with unequal iterated limits, take \(n\geq3\), \(\phi\) from \eqref{example:abel-polynomial}, and a nonnegative smooth radial function \(b\) with support \(\overline{B'_{1/8}}\) and positive integral.
Set
\begin{align*}
	\psi(y')=b(y'-\tfrac12\mathbf e_2')+b(y'-\tfrac12\mathbf e_1').
\end{align*}
Both supports lie in \(B'_1\).
The two translates have equal mass, so \(m_1/M_0=1/4\).
In the direction \(\omega''=\mathbf e_2'\), the maximum of \(y_2\) on the support is attained only at \(\tfrac58\mathbf e_2'\), whose first coordinate is zero.
Thus Theorem \ref{thm-260017}{\rm(iv)(b)} and {\rm(iv)(c)} give, respectively,
\begin{align*}
	\lim_{t\to\infty}\lim_{r\to\infty}x_1^*(r\mathbf e_2',t)=0,
	\qquad
	\lim_{r\to\infty}\lim_{t\to\infty}x_1^*(r\mathbf e_2',t)=\tfrac14.
\end{align*}

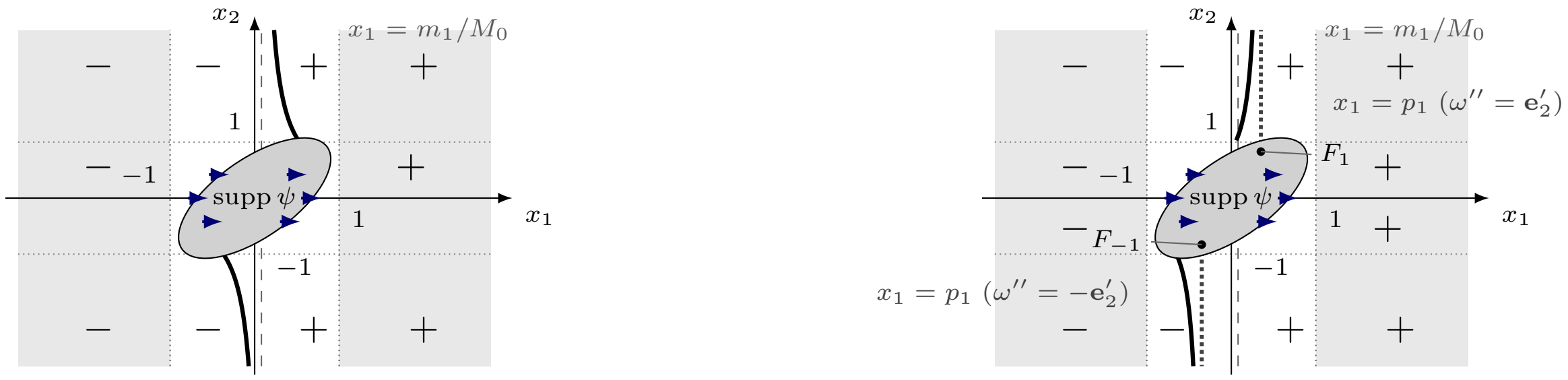
\begin{figure}[H]
	\centering
	\begin{minipage}[t]{0.49\textwidth}
		\centering
		\begin{tikzpicture}[x=.83cm,y=.55cm,>=latex,
				every node/.style={font=\scriptsize}]
			\fill[black!9] (1,-3) rectangle (2.8,3);
			\fill[black!9] (-2.8,-3) rectangle (-1,3);
			\draw[densely dotted,black!60] (-1,-3) -- (-1,3);
			\draw[densely dotted,black!60] (1,-3) -- (1,3);
			\draw[densely dotted,black!45] (-2.8,1) -- (2.8,1);
			\draw[densely dotted,black!45] (-2.8,-1) -- (2.8,-1);
			\draw[dashed,black!60] (0.08,-3) -- (0.08,3);
			\draw[->] (-2.95,0) -- (3.05,0) node[below right] {\(x_1\)};
			\draw[->] (0,-3.15) -- (0,3.25) node[left] {\(x_2\)};
			\draw[very thick]
			plot[domain=1.03:3,samples=60,smooth] ({0.08+0.45/\x},\x);
			\draw[very thick]
			plot[domain=-3:-1.03,samples=60,smooth] ({0.08+0.45/\x},\x);
			\begin{scope}[rotate=35]
				\fill[black!18] (0,0) ellipse (1.05 and 0.70);
				\draw (0,0) ellipse (1.05 and 0.70);
			\end{scope}
			\foreach \X/\Y in {-0.42/0.42,0.50/0.42,-0.67/0,0.67/0,
					-0.50/-0.42,0.42/-0.42}{
					\draw[->,line width=.8pt,blue!45!black]
					({\X-0.12},\Y)--++(0.24,0);
				}
			\node at (0,0) {$\supp\psi$};
			\node[font=\large] at (2,2.35) {\(+\)};
			\node[font=\large] at (-1.85,2.35) {\(-\)};
			\node[font=\large] at (2,-2.35) {\(+\)};
			\node[font=\large] at (-1.85,-2.35) {\(-\)};
			\node[font=\large] at (-0.55,2.35) {\(-\)};
			\node[font=\large] at (0.70,2.35) {\(+\)};
			\node[font=\large] at (-0.55,-2.35) {\(-\)};
			\node[font=\large] at (0.70,-2.35) {\(+\)};
			\node[font=\large] at (-1.85,0.55) {\(-\)};
			\node[font=\large] at (1.85,0.55) {\(+\)};
			\node[below right] at (1,0) {\(1\)};
			\node[above left] at (-1,0) {\(-1\)};
			\node[above left] at (0,1) {\(1\)};
			\node[below right,inner sep=1pt] at (0.20,-1) {\(-1\)};
			\node[above right,black!65] at (0.95,2.55) {\(x_1=m_1/M_0\)};
		\end{tikzpicture}

		\textnormal{\rm(a) \(0<t<1\)}
	\end{minipage}\hfill
	\begin{minipage}[t]{0.49\textwidth}
		\centering
		\begin{tikzpicture}[x=.83cm,y=.55cm,>=latex,
				every node/.style={font=\scriptsize}]
			\fill[black!9] (1,-3) rectangle (2.8,3);
			\fill[black!9] (-2.8,-3) rectangle (-1,3);
			\draw[densely dotted,black!60] (-1,-3) -- (-1,3);
			\draw[densely dotted,black!60] (1,-3) -- (1,3);
			\draw[densely dotted,black!45] (-2.8,1) -- (2.8,1);
			\draw[densely dotted,black!45] (-2.8,-1) -- (2.8,-1);
			\draw[dashed,black!60] (0.08,-3) -- (0.08,3);
			\draw[->] (-2.95,0) -- (3.05,0) node[below right] {\(x_1\)};
			\draw[->] (0,-3.15) -- (0,3.25) node[left] {\(x_2\)};
			\draw[very thick]
			plot[domain=1.03:3,samples=60,smooth] ({0.35-0.30/\x},\x);
			\draw[very thick]
			plot[domain=-3:-1.03,samples=60,smooth] ({-0.35+0.30/\x},\x);
			\draw[densely dotted,black!75,line width=1.1pt] (0.35,0.83) -- (0.35,3);
			\draw[densely dotted,black!75,line width=1.1pt] (-0.35,-0.83) -- (-0.35,-3);
			\begin{scope}[rotate=35]
				\fill[black!18] (0,0) ellipse (1.05 and 0.70);
				\draw (0,0) ellipse (1.05 and 0.70);
			\end{scope}
			\foreach \X/\Y in {-0.42/0.42,0.50/0.42,-0.67/0,0.67/0,
					-0.50/-0.42,0.42/-0.42}{
					\draw[->,line width=.8pt,blue!45!black]
					({\X-0.12},\Y)--++(0.24,0);
				}
			\node at (0,0) {$\supp\psi$};
			\node[font=\large] at (2,2.35) {\(+\)};
			\node[font=\large] at (-1.85,2.35) {\(-\)};
			\node[font=\large] at (2,-2.35) {\(+\)};
			\node[font=\large] at (-1.85,-2.35) {\(-\)};
			\node[font=\large] at (-0.7,2.35) {\(-\)};
			\node[font=\large] at (0.7,2.35) {\(+\)};
			\node[font=\large] at (-0.7,-2.35) {\(-\)};
			\node[font=\large] at (0.7,-2.35) {\(+\)};
			\node[font=\large] at (-1.85,0.55) {\(-\)};
			\node[font=\large] at (1.85,-0.55) {\(+\)};
			\node[font=\large] at (-1.85,-0.55) {\(-\)};
			\node[font=\large] at (1.85,0.55) {\(+\)};
			\node[below right] at (1,0) {\(1\)};
			\node[above left] at (-1,0) {\(-1\)};
			\node[above left] at (0,1) {\(1\)};
			\node[below right,inner sep=1pt] at (0.20,-1) {\(-1\)};
			\fill (0.35,0.83) circle (1.2pt);
			\fill (-0.35,-0.83) circle (1.2pt);
			\draw[thin,black!60] (0.35,0.83)--(0.95,0.78);
			\node[right,inner sep=1pt] at (1,0.78) {\(F_1\)};
			\draw[thin,black!60] (-0.35,-0.83)--(-0.95,-0.78);
			\node[left,inner sep=1pt] at (-1,-0.78) {\(F_{-1}\)};
			\node[right,black!75] at (1.05,1.7)
			{\(x_1=p_1\; (\omega''=\mathbf e_2')\)};
			\node[left,black!75] at (-1.05,-1.7)
			{\(x_1=p_1\; (\omega''=-\mathbf e_2')\)};
			\node[above right,black!65,align=left] at (0.95,2.55)
			{\(x_1=m_1/M_0\)};
		\end{tikzpicture}

		\textnormal{\rm(b) \(t>1\)}
	\end{minipage}
	\caption{A schematic illustration of zeros of \(\partial_{x_n}^2w_n(\cdot,0,t)\) and their far-field asymptotic behavior for \(n=3\)}
	\label{260018}
\end{figure}

Theorem \ref{thm-260017}{\rm(i)} is proved in Proposition \ref{prop-260121}.
The statements in Theorem \ref{thm-260017}{\rm(ii)} for general and for radial \(\psi\) are proved in Propositions \ref{prop-260125} and \ref{prop-260123}, respectively.
Theorem \ref{thm-260017}{\rm(iii)} is proved in Proposition \ref{prop-260127}.
Theorem \ref{thm-260017}{\rm(iv)(a)}, {\rm(iv)(b)}, and {\rm(iv)(c)} follow from Proposition \ref{prop-260131}{\rm(i)}, {\rm(ii)}, and {\rm(iii)}, respectively.
Figure~\ref{260018} illustrates Theorem \ref{thm-260017} for \(n=3\).

The remaining main theorems concern the signs of the tangential pressure derivatives for normal and tangential boundary data.
Favorable and adverse pressure gradients play an important role in boundary-layer separation.
For a chosen tangential direction, the pressure gradient is called favorable when it points opposite to the near-wall tangential velocity and adverse when it points in the same direction.
Outside \(\supp\psi\), the boundary velocity vanishes.
Hence, where \(\partial_{x_n}w_i(x',0,t)\neq0\), the gradient is favorable when \(\partial_{x_i}p(x',0,t)\partial_{x_n}w_i(x',0,t)<0\).
It is adverse when this product is positive.

We now state the sixth main theorem concerning the tangential pressure derivatives for normal boundary data.

\begin{theorem}\label{thm-260019}
	\textbf{(Pressure derivative with normal boundary data)} Let \(n\geq2\), \(1\leq i\leq n-1\), and let \(p\) be the pressure of the Stokes system \eqref{0308-01} with boundary data \eqref{0324-13}, where \(j=n\).
	Suppose that \(\psi\) satisfies {\rm\textbf{(S)}} and is radial, and that \(\phi\) satisfies {\rm\textbf{(T1)}} and {\rm\textbf{(T2)}}.
	Then the following hold.
	\begin{enumerate}
		\item For \(0<t<t_1\), there exists \(r_{\rm nor}(t)>0\) such that
		\begin{align*}
			\sgn(\partial_{x_i}p(x',0,t)) =
			\begin{cases}
				\sgn(x_i),& 1<|x'|\leq r_{\rm nor}(t)-1, \\
				-\sgn(x_i), & |x'|\geq r_{\rm nor}(t)+1,
			\end{cases}
		\end{align*}
		where the first region is empty if \(r_{\rm nor}(t)\leq 2\).
		Moreover, \(r_{\rm nor}\) is continuous on \((0,t_1)\), and
		\begin{align*}
			n<\frac{\phi'(t)r_{\rm nor}(t)^2}{\phi(t)}<2n \qquad(0<t<t_1),
		\end{align*}
		and
		\begin{align}\label{261009}
			r_{\rm nor}(t) =\left(\frac{n\phi(t)}{\phi'(t)}\right)^{\frac12}(1+o(1)) \qquad\text{as }t\to t_1^-.
		\end{align}
		If \(\phi\) is log-concave on \((0,1)\), then \(r_{\rm nor}(t)\to0\) as \(t\to0^+\).
		\item For \(t<t_1\) sufficiently close to \(t_1\), there exists a unique \(\rho_{\rm nor}(t)>1\) such that
		\begin{align*}
			\sgn(\partial_{x_i}p(x',0,t))=-\sgn(x_i)\sgn(|x'|-\rho_{\rm nor}(t)),\qquad |x'|>1.
		\end{align*}
		Moreover, \(|\rho_{\rm nor}(t)-r_{\rm nor}(t)|<1\), and as \(t\to t_1^-\),
		\begin{align*}
			\rho_{\rm nor}(t)^2-r_{\rm nor}(t)^2
			=\frac{n+3}{n-1}\frac{M_2}{M_0}+O\bigl(r_{\rm nor}(t)^{-2}\bigr).
		\end{align*}
		\item If \(|x'|>1\) and \(t\geq t_1\), then \(\sgn(\partial_{x_i}p(x',0,t))=\sgn(x_i)\).
	\end{enumerate}
\end{theorem}
\begin{figure}[H]
	\centering
	\begin{tikzpicture}[scale=0.82,>=latex,
			every node/.style={font=\footnotesize}]
		\begin{scope}[shift={(3.6,0)}]
			\draw[->] (-2.45,0)--(2.55,0) node[below right] {\(x_1\)};
			\draw[->] (0,-2.45)--(0,2.55) node[left] {\(x_2\)};
			\fill[black!18] (0,0) circle (1);
			\draw (0,0) circle (1);
			\node[font=\scriptsize] at (0,0) {$\overline{B_1'}$};
			\draw[very thick] (0,-2.3)--(0,-1.02);
			\draw[very thick] (0,1.02)--(0,2.3);
			\node[font=\Large] at (1.65,1.45) {\(+\)};
			\node[font=\Large] at (-1.65,1.45) {\(-\)};
			\node[font=\Large] at (1.65,-1.45) {\(+\)};
			\node[font=\Large] at (-1.65,-1.45) {\(-\)};
			\node at (0,-2.85) {\(t\geq t_1\)};
		\end{scope}
		\begin{scope}[shift={(-3.6,0)}]
			\draw[->] (-2.45,0)--(2.55,0) node[below right] {\(x_1\)};
			\draw[->] (0,-2.45)--(0,2.55) node[left] {\(x_2\)};
			\fill[black!18] (0,0) circle (1);
			\draw (0,0) circle (1);
			\node[font=\scriptsize] at (0,0) {$\overline{B_1'}$};
			\draw[very thick] (0,-2.3)--(0,-1.02);
			\draw[very thick] (0,1.02)--(0,2.3);
			\fill[black!7,even odd rule] (0,0) circle (2.25) (0,0) circle (1.95);
			\draw[very thick,densely dashed] (0,0) circle (1.95);
			\draw[very thick,densely dashed] (0,0) circle (2.25);
			\node[font=\scriptsize,inner sep=1pt,anchor=west] at (2.45,0.75)
			{\(r_{\rm nor}(t)-1\)};
			\node[font=\scriptsize,inner sep=1pt,anchor=west] at (2.30,-0.75)
			{\(r_{\rm nor}(t)+1\)};
			\node[font=\Large] at (1.20,0.45) {\(+\)};
			\node[font=\Large] at (-1.20,0.45) {\(-\)};
			\node[font=\Large] at (2.15,1.55) {\(-\)};
			\node[font=\Large] at (-2.15,1.55) {\(+\)};
			\node[font=\Large] at (1.20,-0.45) {\(+\)};
			\node[font=\Large] at (-1.20,-0.45) {\(-\)};
			\node[font=\Large] at (2.15,-1.55) {\(-\)};
			\node[font=\Large] at (-2.15,-1.55) {\(+\)};
			\node at (0,-2.85) {\(\substack{t<t_1\text{ sufficiently close}\\\text{to }t_1}\)};
		\end{scope}
	\end{tikzpicture}
	\caption{A schematic illustration of signs of \(\partial_{x_1}p(\cdot,0,t)\) with normal boundary data for \(n=3\) and fixed \(t\). The central disk is the closed unit ball containing \(\supp\psi\).}
	\label{260020}
\end{figure}
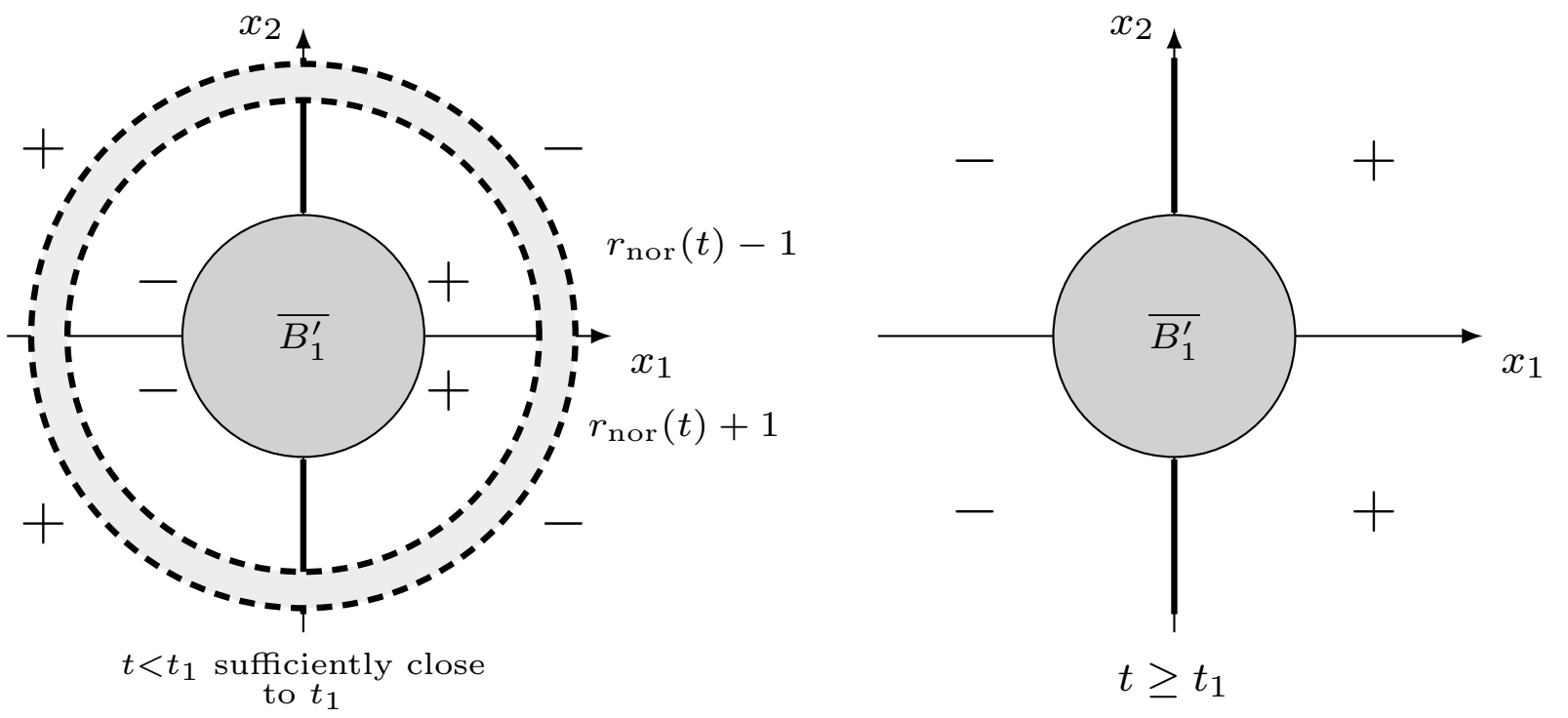
Theorem \ref{thm-260019}{\rm(i)}--{\rm(iii)} are proved in Proposition \ref{prop-260147}{\rm(i)}--{\rm(iii)}, respectively.
Figure~\ref{260020} illustrates Theorem \ref{thm-260019} for \(n=3\).

We now state the seventh main theorem concerning the pressure derivative in the direction of the tangential boundary data.
\begin{theorem}\label{thm-260021}
	\textbf{(Pressure derivative in the direction of the tangential boundary data)} Let \(n\geq2\), \(1\leq i=j\leq n-1\), and let \(p\) be the pressure of the Stokes system \eqref{0308-01} with boundary data \eqref{0324-13}.
	Suppose that \(\psi\) satisfies {\rm\textbf{(S)}} and is radial, and that \(\phi\) satisfies {\rm\textbf{(T1)}} and \(\phi\not\equiv0\).
	Then the following hold.
	\begin{enumerate}
		\item\label{1013-110}\textbf{(Signs for sufficiently large \(|x'|\))} Suppose that \(\phi\) satisfies {\rm\textbf{(T2)}}--{\rm\textbf{(T4)}}. Then the following hold.
		\begin{enumerate}
			\item\textbf{(On the \(x_i\)-axis)} For every fixed \(t>0\) with \(t\ne t_\infty\), there exists \(R_{\rm ax}(t)>1\) such that, for \(|r|>R_{\rm ax}(t)\), 
			\begin{align*}
				\sgn\bigl(\partial_{x_i}p(r\mathbf e_i',0,t)\bigr) =-\sgn(T_0(t)) =
				\begin{cases}
					-1, & 0<t<t_\infty, \\
					1, & t>t_\infty.
				\end{cases}
			\end{align*}
			\item\textbf{(For general \(x'\))} If \(n\geq3\), \(t\in(0,\infty)\setminus\{t_\infty\}\), and \(\varepsilon\in(0,\sqrt{n-1})\), then there exists \(R_\varepsilon(t)>1\) such that, for \(|x'|>R_\varepsilon(t)\),
			\begin{align*}
				\sgn(T_0(t))\partial_{x_i}p(x',0,t)
				\begin{cases}
					>0, & |x^{(i)}|\geq(\sqrt{n-1}+\varepsilon)|x_i|, \\
					<0, & |x^{(i)}|\leq(\sqrt{n-1}-\varepsilon)|x_i|.
				\end{cases}
			\end{align*}
		\end{enumerate}

		\item\label{260311}\textbf{(For sufficiently small \(t\))} 
		Suppose that \(\phi\) satisfies {\rm\textbf{(T2)}}. Then there exists \(t_0=t_0(n,\psi,\phi)\in(0,t_1)\) such that, for \(0<t<t_0\),
		\begin{align*}
			\partial_{x_i}p(r\mathbf e_i',0,t)<0,\qquad |r|\geq1.
		\end{align*}

		\item\label{260310}\textbf{(The time \(t_\infty\))} Suppose that \(\psi\) is radially nonincreasing and that \(\phi\) satisfies {\rm\textbf{(T2)}}--{\rm\textbf{(T4)}}. Then the following hold.
		\begin{enumerate}
			\item\textbf{(On the \(x_i\)-axis)} \(\partial_{x_i}p(r\mathbf e_i',0,t_\infty)>0\) for all sufficiently large \(|r|\).
			If \(\phi\) also satisfies {\rm\textbf{(T5)}}, then
			\begin{align*}
				\partial_{x_i}p(r\mathbf e_i',0,t_\infty)>0, \qquad |r|\geq1.
			\end{align*}
			\item\textbf{(For general \(x'\))} If \(n\geq3\) and \(\varepsilon\in(0,\sqrt{n+1})\), then there exists \(R_\varepsilon(t_\infty)>1\) such that, for \(|x'|>R_\varepsilon(t_\infty)\),
			\begin{align*}
				\partial_{x_i}p(x',0,t_\infty)
				\begin{cases}
					<0, & |x^{(i)}|\geq(\sqrt{n+1}+\varepsilon)|x_i|, \\
					>0, & |x^{(i)}|\leq(\sqrt{n+1}-\varepsilon)|x_i|.
				\end{cases}
			\end{align*}
		\end{enumerate}

		\item\label{261017}\textbf{(Before \(t_\infty\))} Suppose that \(\psi\) is radially nonincreasing and that \(\phi\) satisfies {\rm\textbf{(T2)}}--{\rm\textbf{(T5)}}.
		Then there exists \(\delta_\infty=\delta_\infty(n,\psi,\phi)\in(0,t_\infty-t_1)\) such that, for \(t_\infty-\delta_\infty<t<t_\infty\), there is a unique \(\rho(t)>1\) such that \(\partial_{x_i}p(\pm\rho(t)\mathbf e_i',0,t)=0\).
		Moreover,
		\begin{align}\label{261040}
			\sgn\bigl(\partial_{x_i}p(r\mathbf e_i',0,t)\bigr) =\sgn(\rho(t)-|r|),\qquad |r|\geq1,
		\end{align}
		and
		\begin{align*}
			\lim_{t\to t_\infty^-}(t_\infty-t)\rho(t)^2 =-\frac{n(n+1)T_1(t_\infty)} {(n-1)T_0'(t_\infty)}.
		\end{align*}

		\item\label{260312}\textbf{(Zeros on the \(x_i\)-axis at later times)}
		\begin{enumerate}
			\item\textbf{(The last zero in time at \(r=1\))} Suppose that \(\psi\) is radially nonincreasing and that \(\phi\) satisfies {\rm\textbf{(T2)}}--{\rm\textbf{(T5)}}, and define
		\begin{align}\label{0911-23}
			\overline t:=\sup\left\{t>0: \partial_{x_i}p(\mathbf e_i',0,t)=0\right\}.
		\end{align}
		Then \(\overline t\) is well-defined, \(t_\infty<\overline t<\infty\), and
		\begin{align*}
			\partial_{x_i}p(\pm\mathbf e_i',0,t)<0 \quad\text{for }t>\overline t, \qquad \partial_{x_i}p(\pm\mathbf e_i',0,\overline t)=0.
		\end{align*}
			For every \(t>\overline t\), the zeros of \(r\mapsto\partial_{x_i}p(r\mathbf e_i',0,t)\) on \((1,\infty)\) have a minimum \(\rho_-(t)\) and a maximum \(\rho_+(t)\), which satisfy
			\begin{align*}
				\partial_{x_i}p(r\mathbf e_i',0,t)<0 \quad\text{for }1\leq|r|<\rho_-(t),\qquad \partial_{x_i}p(r\mathbf e_i',0,t)>0 \quad\text{for }|r|>\rho_+(t),
			\end{align*}
			and every zero of \(r\mapsto\partial_{x_i}p(r\mathbf e_i',0,t)\) with \(|r|\geq1\) belongs to \([-\rho_+(t),-\rho_-(t)]\cup[\rho_-(t),\rho_+(t)]\).
			\item\textbf{(For sufficiently large \(t\))} There exists \(\widehat t=\widehat t(n)>1\) such that, for every \(t\geq\widehat t\), \(r\mapsto\partial_{x_i}p(r\mathbf e_i',0,t)\) has exactly one zero \(r=\rho(t)\) on \([1,\infty)\).
			It satisfies \(\rho(t)>1\) and \(\partial_r\partial_{x_i}p(r\mathbf e_i',0,t)>0\) at \(r=\rho(t)\).
			Moreover,
			\begin{align*}
				\sgn\bigl(\partial_{x_i}p(r\mathbf e_i',0,t)\bigr) =\sgn\bigl(|r|-\rho(t)\bigr),\qquad |r|\geq1,
			\end{align*}
			and
			\begin{align*}
				\bigl|\rho(t)-c_n\sqrt{t}\bigr|\leq C_nt^{-1/2},\qquad t\geq\widehat t,
			\end{align*}
			where \(c_n,C_n>0\) depend only on \(n\).
			Finally, under the hypotheses of {\rm(v)(a)}, \(\overline t<\widehat t\) and \(\rho_-(t)=\rho_+(t)=\rho(t)\) for \(t\geq\widehat t\).
		\end{enumerate}

		\item\label{261018}\textbf{(Zero curve for \(n\geq3\))} Let \(n\geq3\).
		\begin{enumerate}
			\item\textbf{(For fixed \(t>0\))} Fix \(x^{(i)}\in(\mathbf e_i')^\perp\) and suppose that \(\norm{\phi}_{L^1(0,t)}>0\).
			If \(|x^{(i)}|\geq1\), then \(x_i\mapsto\partial_{x_i}p(x_i\mathbf e_i'+x^{(i)},0,t)\) changes sign at some \(x_i>0\).
			Moreover, there exists \(R_{\rm off}(t)>1\) such that, if \(|x^{(i)}|>R_{\rm off}(t)\), then this positive zero is unique and simple, and \(x_i=\rho(|x^{(i)}|,t)\) is a strictly increasing \(C^1\) function of \(|x^{(i)}|\).
			As \(|x^{(i)}|\to\infty\),
			\begin{align*}
				\rho(|x^{(i)}|,t)=
				\begin{cases}
					|x^{(i)}|/\sqrt{n-1}+O_t(|x^{(i)}|^{-1}),&T_0(t)\ne0,\\
					|x^{(i)}|/\sqrt{n+1}+O_t(|x^{(i)}|^{-1}),&T_0(t)=0.
				\end{cases}
			\end{align*}
			\item\textbf{(For \(t<t_\infty\) sufficiently close to \(t_\infty\))} Suppose that \(\psi\) is radially nonincreasing and that \(\phi\) satisfies {\rm\textbf{(T2)}}--{\rm\textbf{(T5)}}.
			Then, for \(t<t_\infty\) sufficiently close to \(t_\infty\) and every \(x^{(i)}\in(\mathbf e_i')^\perp\), \(x_i\mapsto\partial_{x_i}p(x_i\mathbf e_i'+x^{(i)},0,t)\) has exactly one zero \(x_i>0\) with \(|x_i\mathbf e_i'+x^{(i)}|\geq\rho(t)\), where \(\rho(t)\) is given in {\rm(iv)}.
			This zero is the largest positive zero, and it is the zero \(\rho(|x^{(i)}|,t)\) given in {\rm(vi)(a)} when \(|x^{(i)}|>R_{\rm off}(t)\).
			This zero is simple for every \(x^{(i)}\).
			These zeros leave every bounded subset of \(\Sigma\) as \(t\to t_\infty^-\).
		\end{enumerate}
	\end{enumerate}
\end{theorem}

\begin{figure}[H]
	\centering
	\begin{tikzpicture}[x=1.08cm,y=.60cm,>=latex,
			every node/.style={font=\footnotesize}]
		\fill[black!5] (-1,0) rectangle (1,6.25);
		\draw[densely dashed,black!45] (-1,0)--(-1,6.25);
		\draw[densely dashed,black!45] (1,0)--(1,6.25);
		\draw[densely dashed,black!55] (-4.15,0)--(-4.15,6.25) (4.15,0)--(4.15,6.25);
		\foreach \y in {0.62,2.35,3.20,4.18,5.20}
			\draw[dashed,black!60] (-5.60,\y)--(5.60,\y);

		\draw[->] (-5.75,0)--(5.95,0) node[below right] {\(x_1\)};
		\draw[->] (0,0)--(0,6.55) node[left] {\(t\)};

		\foreach \x/\lab in {-1/{-1},1/{1}}{
		\draw (\x,-0.06)--(\x,0.06);
		\node[below] at (\x,0) {\(\lab\)};
		}
		\foreach \y/\lab in {0.62/{t_0},2.35/{t_\infty-\delta_\infty},3.20/{t_\infty},4.18/{\overline t},5.20/{\widehat t}}{
		\draw (-0.06,\y)--(0.06,\y);
		\node[left,font=\scriptsize,fill=black!5,inner sep=1pt] at (0,\y) {\(\lab\)};
		}
		\node[font=\large] at (-2.05,0.30) {\(-\)};
		\node[font=\large] at (2.05,0.30) {\(-\)};
		\node[font=\large] at (-4.75,1.55) {\(-\)};
		\node[font=\large] at (4.75,1.55) {\(-\)};

		\draw[very thick]
		plot[domain=1.60:5.55,samples=80,smooth]
		(\x,{3.20-2.176/(\x*\x)});
		\draw[very thick]
		plot[domain=-5.55:-1.60,samples=80,smooth]
		(\x,{3.20-2.176/(\x*\x)});
		\draw[very thick,dashed]
		(1.18,1.485)..controls(1.20,1.66)and(1.28,2.01)..(1.60,2.35);
		\draw[very thick,dashed]
		(-1.18,1.485)..controls(-1.20,1.66)and(-1.28,2.01)..(-1.60,2.35);
		\draw[fill=white] (1.60,2.35) circle (1.7pt);
		\draw[fill=white] (-1.60,2.35) circle (1.7pt);
		\node at (3.40,2.62) {\(x_1=\rho(t)\)};
		\node at (-3.40,2.62) {\(x_1=-\rho(t)\)};
		\node[font=\large] at (-4.75,2.58) {\(-\)};
		\node[font=\large] at (-1.30,2.70) {\(+\)};
		\node[font=\large] at (1.30,2.70) {\(+\)};
		\node[font=\large] at (4.75,2.58) {\(-\)};
		\node[font=\large] at (-1.30,3.60) {\(+\)};
		\node[font=\large] at (1.30,3.60) {\(+\)};
		\node[font=\large] at (-4.75,3.72) {\(+\)};
		\node[font=\large] at (4.75,3.72) {\(+\)};
		\node[fill=white,inner sep=1.5pt] at (3.00,3.60)
		{no zeros when \(t=t_\infty\)};

		\draw[very thick,densely dotted]
		(1.82,4.18)..controls(2.10,4.66)and(2.22,4.90)..(2.42,5.20);
		\draw[very thick,densely dotted]
		(-1.82,4.18)..controls(-2.10,4.66)and(-2.22,4.90)..(-2.42,5.20);
		\draw[very thick]
		(2.42,5.20)..controls(2.62,5.50)and(2.90,5.86)..(3.18,6.20);
		\draw[very thick]
		(-2.42,5.20)..controls(-2.62,5.50)and(-2.90,5.86)..(-3.18,6.20);
		\node[above] at (3.18,6.20) {\(x_1=\rho(t)\)};
		\node[above] at (-3.18,6.20) {\(x_1=-\rho(t)\)};
		\node[font=\large] at (-4.75,5.68) {\(+\)};
		\node[font=\large] at (-1.25,5.68) {\(-\)};
		\node[font=\large] at (1.25,5.68) {\(-\)};
		\node[font=\large] at (4.75,5.68) {\(+\)};
		\node[font=\scriptsize,align=center,fill=white,inner sep=1pt] at (-5.00,6.18)
		{\(|x_1|\) sufficiently\\[-2pt]large};
		\node[font=\scriptsize,align=center,fill=white,inner sep=1pt] at (5.00,6.18)
		{\(|x_1|\) sufficiently\\[-2pt]large};

		\node[font=\scriptsize,fill=white,inner sep=1pt] at (3.55,0.31)
		{no zeros for \(0<t<t_0\)};
	\end{tikzpicture}
	\caption{A schematic illustration of the zeros and signs of \(x_1\mapsto\partial_{x_1}p(x_1\mathbf e_1',0,t)\) in the \((x_1,t)\)-plane}
	\label{261300}
\end{figure}
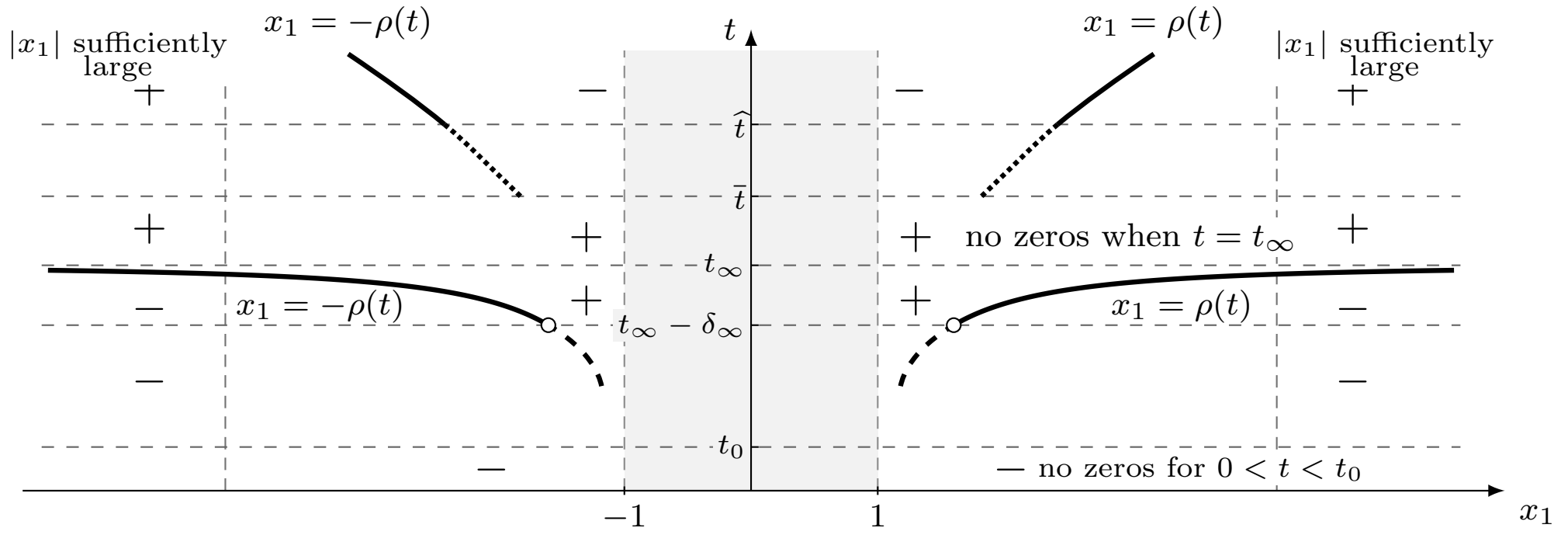

\begin{figure}[!htbp]
	\centering
	\begin{tikzpicture}[x=.80cm,y=.50cm,>=latex,
			every node/.style={font=\scriptsize}]
		\begin{scope}[shift={(-2.95,3.15)}]
			\draw[->] (-2.55,0)--(2.65,0) node[below right] {\(x_1\)};
			\draw[->] (0,-2.30)--(0,2.40) node[left] {\(x_2\)};
			\fill[black!18] (0,0) ellipse [x radius=.5625,y radius=.90];
			\draw (0,0) ellipse [x radius=.5625,y radius=.90];
			\begin{scope}[scale=.84,xscale=.625,yscale=1.35]
				\boundarydataarrows
			\end{scope}
			\node[font=\tiny] at (0,.23) {\(\overline{B_1'}\)};
			\draw[densely dashed,black!55]
			(0.78,1.10)--(1.56,2.20) (-0.78,1.10)--(-1.56,2.20)
			(0.78,-1.10)--(1.56,-2.20) (-0.78,-1.10)--(-1.56,-2.20);
			\node[font=\tiny,anchor=south east]
			at (2.55,2.58) {\(|x_2|=\sqrt{n-1}|x_1|\)};
			\draw[very thick] plot[domain=1.10:2.20,variable=\s]
			({sqrt(.45+\s*\s/2)},\s);
			\draw[very thick] plot[domain=-2.20:-1.10,variable=\s]
			({sqrt(.45+\s*\s/2)},\s);
			\draw[very thick] plot[domain=1.10:2.20,variable=\s]
			({-sqrt(.45+\s*\s/2)},\s);
			\draw[very thick] plot[domain=-2.20:-1.10,variable=\s]
			({-sqrt(.45+\s*\s/2)},\s);
			\foreach \sx/\sy in {1/1,1/-1,-1/1,-1/-1}{
					\draw[very thick,dashed]
					({\sx*.875},{\sy*.210})
					..controls({\sx*.79},{\sy*.42})and({\sx*.9093},{\sy*.88})
					..({\sx*1.0271},{\sy*1.10});
				}
			\foreach \sx/\sy in {1.0271/1.10,1.0271/-1.10,-1.0271/1.10,-1.0271/-1.10}
			\draw[fill=white] (\sx,\sy) circle (1.5pt);
			\foreach \sx in {-1,1}{
					\node[font=\large] at ({\sx*.36},1.52) {\(+\)};
					\node[font=\large] at ({\sx*.36},-1.52) {\(+\)};
					\node[font=\large] at ({\sx*2.14},.55) {\(-\)};
					\node[font=\large] at ({\sx*2.14},-.55) {\(-\)};
				}
			\node at (0,-2.66) {\({\rm(a)}\ 0<t<t_0\)};
		\end{scope}

		\begin{scope}[shift={(2.95,3.15)}]
			\draw[->] (-2.55,0)--(2.65,0) node[below right] {\(x_1\)};
			\draw[->] (0,-2.30)--(0,2.40) node[left] {\(x_2\)};
			\fill[black!18] (0,0) ellipse [x radius=.5625,y radius=.90];
			\draw (0,0) ellipse [x radius=.5625,y radius=.90];
			\begin{scope}[scale=.84,xscale=.625,yscale=1.35]
				\boundarydataarrows
			\end{scope}
			\node[font=\tiny] at (0,.23) {\(\overline{B_1'}\)};
			\draw[densely dashed,black!55]
			(0.78,1.10)--(1.56,2.20) (-0.78,1.10)--(-1.56,2.20)
			(0.78,-1.10)--(1.56,-2.20) (-0.78,-1.10)--(-1.56,-2.20);
			\node[font=\tiny,anchor=south east]
			at (2.55,2.58) {\(|x_2|=\sqrt{n-1}|x_1|\)};
			\foreach \sx/\sy in {1/1,1/-1,-1/1,-1/-1}{
					\draw[very thick]
					({\sx*1.204},0)
					..controls({\sx*1.204},{\sy*.40})and({\sx*1.05},{\sy*.62})..({\sx*1.04},{\sy*.95})
					..controls({\sx*1.03},{\sy*1.32})and({\sx*1.38},{\sy*1.90})..({\sx*1.58},{\sy*2.20});
				}
			\fill (1.204,0) circle (1.6pt);
			\fill (-1.204,0) circle (1.6pt);
			\node[above right] at (1.204,0) {\(\rho(t)\)};
			\node[above left] at (-1.204,0) {\(-\rho(t)\)};
			\foreach \sx in {-1,1}{
					\node[font=\large] at ({\sx*2.15},.90) {\(-\)};
					\node[font=\large] at ({\sx*2.15},-.52) {\(-\)};
				}
			\node at (0,-2.66)
			{\({\rm(b)}\ t_\infty-\delta_\infty<t<t_\infty\)};
		\end{scope}

		\begin{scope}[shift={(-2.95,-3.15)}]
			\draw[->] (-2.55,0)--(2.65,0) node[below right] {\(x_1\)};
			\draw[->] (0,-2.30)--(0,2.40) node[left] {\(x_2\)};
			\fill[black!18] (0,0) ellipse [x radius=.5625,y radius=.90];
			\draw (0,0) ellipse [x radius=.5625,y radius=.90];
			\begin{scope}[scale=.84,xscale=.625,yscale=1.35]
				\boundarydataarrows
			\end{scope}
			\node[font=\tiny] at (0,.23) {\(\overline{B_1'}\)};
			\draw[densely dashed,black!55]
			(0.55,1.10)--(1.10,2.20) (-0.55,1.10)--(-1.10,2.20)
			(0.55,-1.10)--(1.10,-2.20) (-0.55,-1.10)--(-1.10,-2.20);
			\node[font=\tiny,anchor=south east]
			at (2.55,2.58) {\(|x_2|=\sqrt{n+1}|x_1|\)};
			\foreach \sx/\sy in {1/1,1/-1,-1/1,-1/-1}{
					\draw[very thick]
					({\sx*.76},{\sy*1.28})
					..controls({\sx*.84},{\sy*1.55})and({\sx*1.04},{\sy*2.03})
					..({\sx*1.12},{\sy*2.20});
				}
			\foreach \sx/\sy in {1/1,1/-1,-1/1,-1/-1}{
					\draw[very thick,dashed]
					({\sx*.65},{\sy*.62})
					..controls({\sx*.59},{\sy*.78})and({\sx*.704},{\sy*1.091})
					..({\sx*.76},{\sy*1.28});
				}
			\foreach \sx/\sy in {.76/1.28,.76/-1.28,-.76/1.28,-.76/-1.28}
			\draw[fill=white] (\sx,\sy) circle (1.5pt);
			\foreach \sx in {-1,1}{
					\node[font=\large] at ({\sx*.34},1.54) {\(-\)};
					\node[font=\large] at ({\sx*.34},-1.54) {\(-\)};
					\node[font=\large] at ({\sx*2.12},.54) {\(+\)};
					\node[font=\large] at ({\sx*2.12},-.54) {\(+\)};
				}
			\node at (0,-2.66) {\({\rm(c)}\ t=t_\infty\)};
		\end{scope}

		\begin{scope}[shift={(2.95,-3.15)}]
			\draw[->] (-2.55,0)--(2.65,0) node[below right] {\(x_1\)};
			\draw[->] (0,-2.30)--(0,2.40) node[left] {\(x_2\)};
			\fill[black!18] (0,0) ellipse [x radius=.5625,y radius=.90];
			\draw (0,0) ellipse [x radius=.5625,y radius=.90];
			\begin{scope}[scale=.84,xscale=.625,yscale=1.35]
				\boundarydataarrows
			\end{scope}
			\node[font=\tiny] at (0,.23) {\(\overline{B_1'}\)};
			\draw[densely dashed,black!55]
			(0.78,1.10)--(1.56,2.20) (-0.78,1.10)--(-1.56,2.20)
			(0.78,-1.10)--(1.56,-2.20) (-0.78,-1.10)--(-1.56,-2.20);
			\node[font=\tiny,anchor=south east]
			at (2.55,2.58) {\(|x_2|=\sqrt{n-1}|x_1|\)};
			\foreach \sx/\sy in {1/1,1/-1,-1/1,-1/-1}{
					\draw[very thick]
					({\sx*1.48},{\sy*1.10})
					..controls({\sx*1.48},{\sy*1.45})and({\sx*1.78},{\sy*2.00})
					..({\sx*1.92},{\sy*2.20});
				}
			\foreach \sx in {-1,1}{
					\draw[very thick,densely dashed]
					({\sx*1.34},0)
					..controls({\sx*1.34},.55)and({\sx*1.48},.72)
					..({\sx*1.48},1.10);
					\draw[very thick,densely dashed]
					({\sx*1.34},0)
					..controls({\sx*1.34},-.55)and({\sx*1.48},-.72)
					..({\sx*1.48},-1.10);
				}
			\foreach \sx/\sy in {1.48/1.10,1.48/-1.10,-1.48/1.10,-1.48/-1.10}
			\draw[fill=white] (\sx,\sy) circle (1.5pt);
			\foreach \sx in {-1.34,-1.08,1.08,1.34} \fill (\sx,0) circle (1.6pt);
			\draw[thin,black!60] (1.08,-.06)--(.65,-1.10);
			\draw[thin,black!60] (1.34,-.06)--(2.02,-.96);
			\node[font=\scriptsize] at (.50,-1.40) {\(\rho_-(t)\)};
			\node[font=\scriptsize] at (2.15,-1.15) {\(\rho_+(t)\)};
			\foreach \sx in {-1,1}{
					\node[font=\large] at ({\sx*.34},1.70) {\(-\)};
					\node[font=\large] at ({\sx*.34},-1.95) {\(-\)};
					\node[font=\large] at ({\sx*1.12},.78) {\(-\)};
					\node[font=\large] at ({\sx*1.12},-.50) {\(-\)};
					\node[font=\large] at ({\sx*2.18},.82) {\(+\)};
					\node[font=\large] at ({\sx*2.18},-.58) {\(+\)};
				}
			\node at (0,-2.66) {\({\rm(d)}\ t>\overline t\)};
		\end{scope}
	\end{tikzpicture}
	\caption{A schematic illustration of the zeros of \(\partial_{x_1}p(\cdot,0,t)\) for \(j=1\) and \(n=3\)}
	\label{1013-111}
\end{figure}
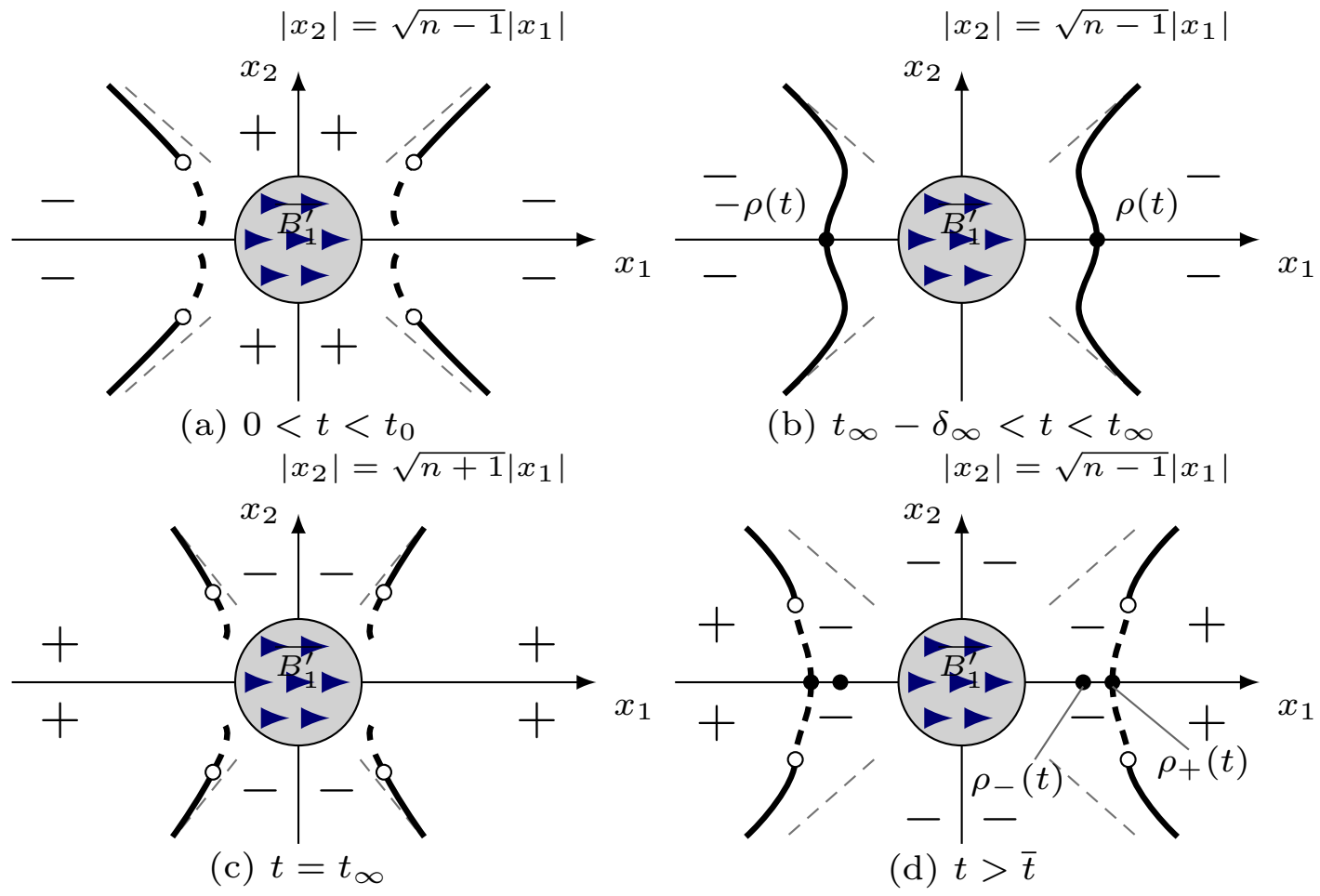

For \(n\geq3\), Theorem \ref{thm-260021}{\rm(i)}--{\rm(iv)} are proved in Proposition \ref{prop-261020}{\rm(i)}--{\rm(iv)}, respectively, and Theorem \ref{thm-260021}{\rm(v)(a)} follows from Proposition \ref{prop-261020}{\rm(iii)} and {\rm(v)}.
For \(n=2\), Theorem \ref{thm-260021}{\rm(i)} and {\rm(ii)} follow from Proposition \ref{prop-260154}{\rm(i)} and {\rm(iii)}, respectively, {\rm(iii)(a)} from Proposition \ref{prop-260154}{\rm(i)} and {\rm(ii)}, {\rm(iv)} from Proposition \ref{prop-260154}{\rm(iv)}, and Theorem \ref{thm-260021}{\rm(v)(a)} from Proposition \ref{prop-260154}{\rm(iii)} and {\rm(v)}.
For every \(n\geq2\), Theorem \ref{thm-260021}{\rm(v)(b)} follows from Lemma \ref{lem-261080}{\rm(iii)}--{\rm(iv)} and Propositions \ref{prop-260154}{\rm(v)} and \ref{prop-261020}{\rm(v)}.
Theorem \ref{thm-260021}{\rm(vi)} is proved in Proposition \ref{prop-261022}, and \eqref{261036} gives the coefficient of \(|x^{(i)}|^{-1}\) in {\rm(vi)(a)} when \(T_0(t)\ne0\).
Figures~\ref{261300} and~\ref{1013-111} illustrate Theorem \ref{thm-260021}.

We now state the eighth main theorem concerning the pressure derivative in a direction different from that of the tangential boundary data.
\begin{theorem}\label{thm-260027}
	\textbf{(Pressure derivative in a direction different from that of the tangential boundary data)} Let \(n\geq3\), \(1\leq i\neq j\leq n-1\), and let \(p\) be the pressure of the Stokes system \eqref{0308-01} with boundary data \eqref{0324-13}.
	Suppose that \(\psi\) satisfies {\rm\textbf{(S)}} and is radial, and that \(\phi\) satisfies {\rm\textbf{(T1)}} and \(\phi\not\equiv0\).
	Then the following hold.
	\begin{enumerate}
		\item Suppose that \(\phi\) satisfies {\rm\textbf{(T2)}}. If \(|x'|\geq\sqrt2\) and \(0<t\leq t_1\), then \(\sgn(\partial_{x_i}p(x',0,t))=-\sgn(x_ix_j)\).
		\item Suppose that \(\phi\) satisfies {\rm\textbf{(T2)}}--{\rm\textbf{(T4)}}. If \(t_1<t<t_\infty\), then there exists \(r_{\rm tan}(t)>0\) such that
		\begin{align*}
			\sgn(\partial_{x_i}p(x',0,t)) =
			\begin{cases}
				\sgn(x_ix_j),& \sqrt2\leq |x'|\leq r_{\rm tan}(t)-1, \\
				-\sgn(x_ix_j), & |x'|\geq\max\{\sqrt2,r_{\rm tan}(t)+1\},
			\end{cases}
		\end{align*}
		where the first region is empty if \(r_{\rm tan}(t)<1+\sqrt2\).
		Moreover, \(r_{\rm tan}(t)\) is continuous and satisfies
		\begin{align}\label{261010}
			r_{\rm tan}(t) =\left(-\frac{(n+2)T_1(t_\infty)}{T_0'(t_\infty)}\right)^{\frac12} (t_\infty-t)^{-\frac12}(1+o(1)) \qquad\text{as }t\to t_\infty^-.
		\end{align}
		\item Suppose that \(\phi\) satisfies {\rm\textbf{(T2)}}--{\rm\textbf{(T4)}}. Then, for \(t<t_\infty\) sufficiently close to \(t_\infty\), there exists a unique \(\rho_{\rm tan}(t)>\sqrt2\) such that
		\begin{align*}
			\sgn(\partial_{x_i}p(x',0,t))=-\sgn(x_ix_j)\sgn(|x'|-\rho_{\rm tan}(t)),\qquad |x'|\geq\sqrt2.
		\end{align*}
		Along each ray with \(x_ix_j\ne0\), this zero is simple.
		As \(t\to t_\infty^-\), \(\rho_{\rm tan}(t)-r_{\rm tan}(t)=O(r_{\rm tan}(t)^{-1})\), and
		\begin{align*}
			\rho_{\rm tan}(t)=\left(-\frac{(n+2)T_1(t_\infty)}{T_0'(t_\infty)}\right)^{\frac12}(t_\infty-t)^{-\frac12}(1+o(1)) \qquad\text{as }t\to t_\infty^-.
		\end{align*}
		\item If \(|x'|\geq\sqrt2\) and \(t\geq1\), then \(\sgn(\partial_{x_i}p(x',0,t))=\sgn(x_ix_j)\).
		Moreover, if \(\phi\) satisfies {\rm\textbf{(T2)}}--{\rm\textbf{(T4)}}, then the same sign holds for \(t_\infty\leq t<1\).
	\end{enumerate}
\end{theorem}
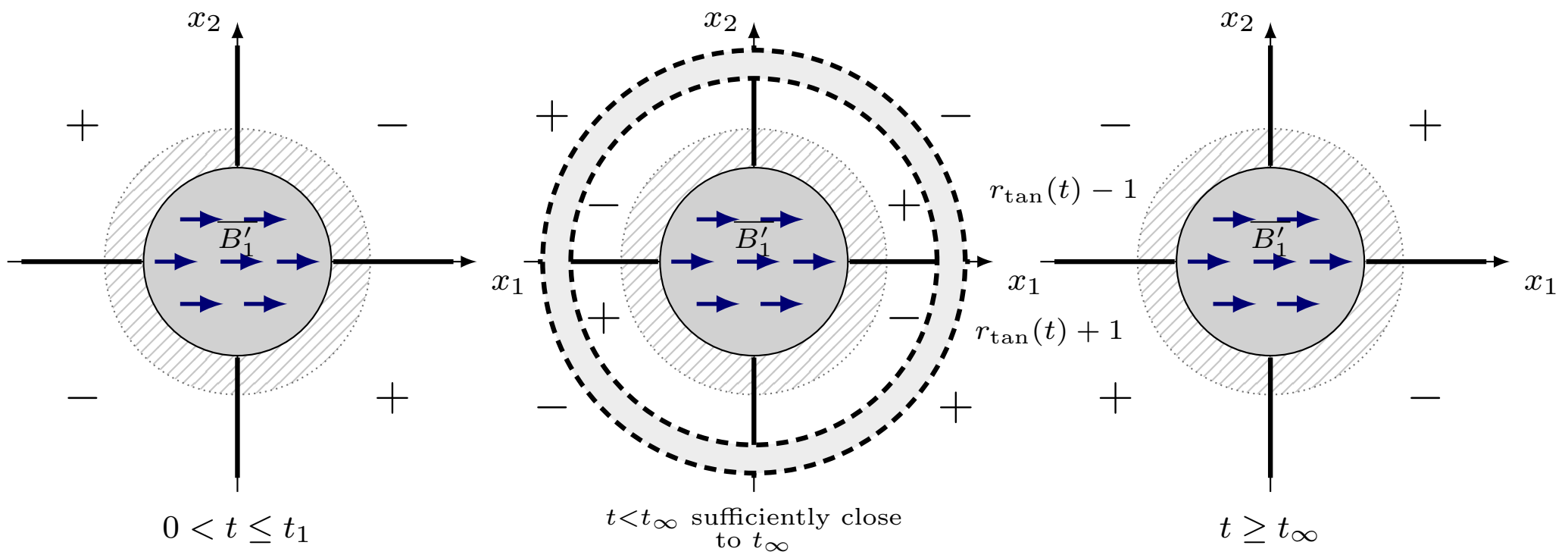
\begin{figure}[ht]
	\centering
	\begin{tikzpicture}[scale=0.82,>=latex,
			every node/.style={font=\footnotesize}]
		\foreach \x/\lab in {-5.5/{0<t\leq t_1},0/{\substack{t<t_\infty\text{ sufficiently close}\\\text{to }t_\infty}},5.5/{t\geq t_\infty}}{
		\begin{scope}[shift={(\x,0)}]
			\draw[->] (-2.45,0)--(2.55,0) node[below right] {\(x_1\)};
			\draw[->] (0,-2.45)--(0,2.55) node[left] {\(x_2\)};
			\fill[pattern=north east lines,pattern color=black!22,even odd rule]
(0,0) circle ({sqrt(2)}) (0,0) circle (1);
\draw[densely dotted,black!55] (0,0) circle ({sqrt(2)});
\fill[black!18] (0,0) circle (1);
			\draw (0,0) circle (1);

			\boundarydataarrows

			\node[font=\scriptsize] at (0,0.25) {\(\overline{B_1'}\)};
			\draw[very thick] (-2.3,0)--(-1.02,0);
			\draw[very thick] (1.02,0)--(2.3,0);
			\draw[very thick] (0,-2.3)--(0,-1.02);
			\draw[very thick] (0,1.02)--(0,2.3);
			\node at (0,-2.85) {\(\lab\)};
		\end{scope}
		}
		\begin{scope}[shift={(-5.5,0)}]
			\node[font=\Large] at (1.65,1.45) {\(-\)};
			\node[font=\Large] at (-1.65,1.45) {\(+\)};
			\node[font=\Large] at (-1.65,-1.45) {\(-\)};
			\node[font=\Large] at (1.65,-1.45) {\(+\)};
		\end{scope}
		\begin{scope}[shift={(0,0)}]
			\fill[black!7,even odd rule] (0,0) circle (2.25) (0,0) circle (1.95);
			\draw[very thick,densely dashed] (0,0) circle (1.95);
			\draw[very thick,densely dashed] (0,0) circle (2.25);
			\node[font=\scriptsize,inner sep=1pt,anchor=west] at (2.45,0.75)
			{\(r_{\rm tan}(t)-1\)};
			\node[font=\scriptsize,inner sep=1pt,anchor=west] at (2.30,-0.75)
			{\(r_{\rm tan}(t)+1\)};
			\node[font=\Large] at (1.60,0.60) {\(+\)};
			\node[font=\Large] at (-1.60,0.60) {\(-\)};
			\node[font=\Large] at (-1.60,-0.60) {\(+\)};
			\node[font=\Large] at (1.60,-0.60) {\(-\)};
			\node[font=\Large] at (2.15,1.55) {\(-\)};
			\node[font=\Large] at (-2.15,1.55) {\(+\)};
			\node[font=\Large] at (-2.15,-1.55) {\(-\)};
			\node[font=\Large] at (2.15,-1.55) {\(+\)};
		\end{scope}
		\begin{scope}[shift={(5.5,0)}]
			\node[font=\Large] at (1.65,1.45) {\(+\)};
			\node[font=\Large] at (-1.65,1.45) {\(-\)};
			\node[font=\Large] at (-1.65,-1.45) {\(+\)};
			\node[font=\Large] at (1.65,-1.45) {\(-\)};
		\end{scope}
	\end{tikzpicture}
		\caption{A schematic illustration of signs of \(\partial_{x_2}p(\cdot,0,t)\) for \(i=2\), \(j=1\), and \(n=3\), with \(t\) fixed. The central disk contains \(\supp\psi\), and the dotted circle has radius \(\sqrt{2}\).}
	\label{260028}
\end{figure}
Theorem \ref{thm-260027}{\rm(i)}--{\rm(iv)} are proved in Proposition \ref{prop-0831-01}{\rm(i)}--{\rm(iv)}, respectively.
Figure~\ref{260028} illustrates Theorem \ref{thm-260027} for \(n=3\).
\begin{remark}
	Take \(i=2\), \(j=1\), and radial \(\psi\).
	By Proposition \ref{prop-0708-10}, \(\sgn(\partial_{x_n}w_2(x',0,t))=\sgn(x_1x_2)\) for \(0<t<1\) and \(|x'|>1\).
	Hence, for \(x_1x_2\neq0\) and \(|x'|\geq\sqrt2\), Theorem \ref{thm-260027} gives a favorable pressure gradient for \(0<t\leq t_1\) and an adverse one for \(t_\infty\leq t<1\), and, for \(t_1<t<t_\infty\), an adverse one for \(|x'|\leq r_{\rm tan}(t)-1\) and a favorable one for \(|x'|\geq r_{\rm tan}(t)+1\).
	For \(t\geq1\), \eqref{1007-08} gives \(\partial_{x_n}w_2=0\) and \(\partial_{x_n}^2w_2=\partial_{x_2}p\) at the boundary, so the pressure gradient is adverse wherever \(\partial_{x_2}p(x',0,t)\ne0\).
\end{remark}
\begin{remark}\label{rmk-T1-relaxed}
The regularity of \(\phi\) at \(t=1\) required in {\rm\textbf{(T1)}} is assumed for simplicity.
In {\rm\textbf{(T1)}}, one could instead require only that \(\phi\) be continuous on \([0,\infty)\) and belong to \(C^{1,\beta}([0,T])\) for every \(T\in(0,1)\), which allows \(\phi'\) to be unbounded near \(t=1\).
For example, for \(0<a<1\), the function \(\phi(t)=t^2(1-t)^a\) on \([0,1]\), extended by zero, satisfies this weaker condition and {\rm\textbf{(T2)}}--{\rm\textbf{(T4)}}, but not {\rm\textbf{(T1)}}.
We expect that many of our results extend to such temporal boundary data, possibly under additional conditions for the assertions at \(t=1\) itself, but we do not pursue this here.
\end{remark}
\subsection{Relation to flow reversal}
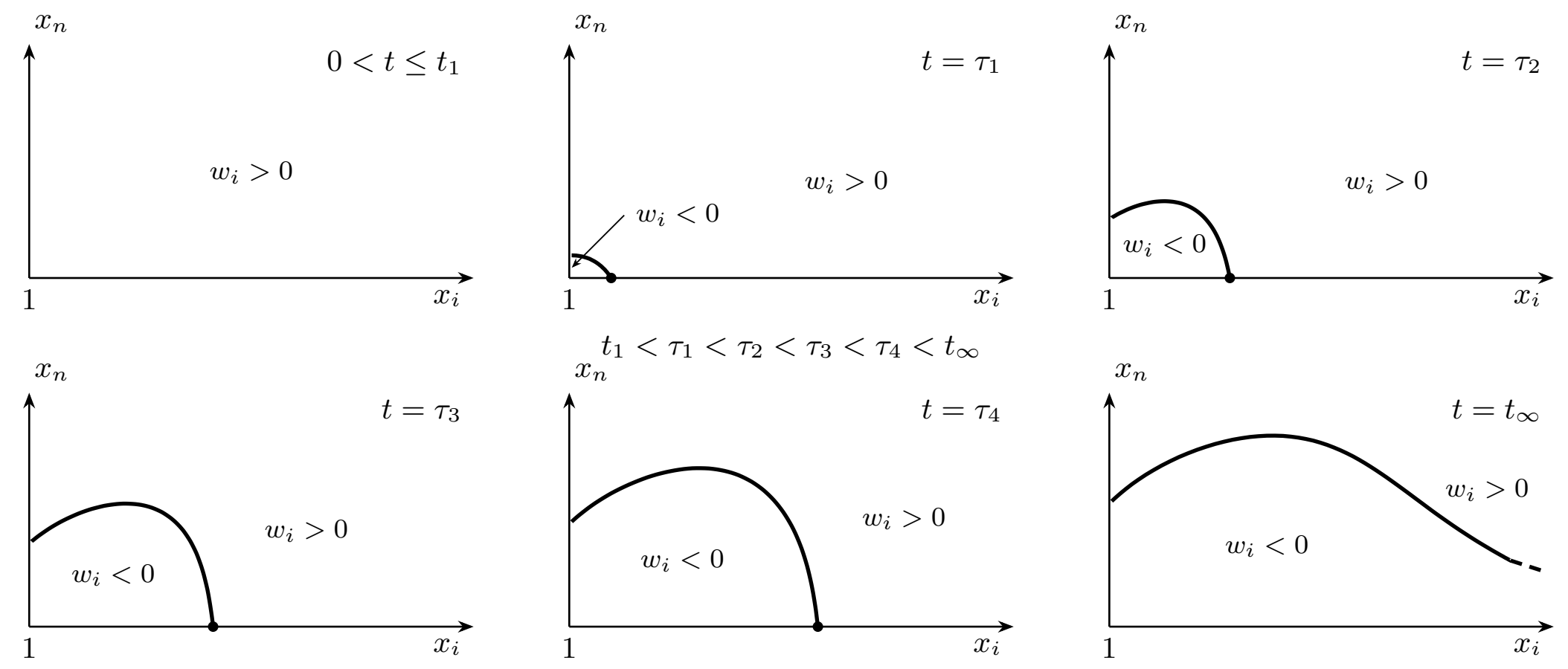
\begin{figure}[H]
\centering
	\begin{tikzpicture}[
		scale=0.85,
		transform shape,
		axis/.style={-{Stealth[length=4pt]}, line width=0.5pt},
		zero/.style={line width=1pt},
		td/.style={circle, fill, inner sep=1.1pt},
		lab/.style={font=\small},
		sgn/.style={font=\footnotesize},
		lead/.style={line width=0.35pt}]
        \node[lab] at (7.975,-0.72) {$t_1<\tau_1<\tau_2<\tau_3<\tau_4<t_\infty$};
		\begin{scope}[shift={(0.000,0.000)}]
			\draw[axis] (0,0) -- (4.650,0) node[below left, lab] {$x_i$};
			\draw[axis] (0,0) -- (0,2.450) node[above right, lab, xshift=-2pt] {$x_n$};
			\node[below, lab] at (0,0) {$1$};
			\node[lab, anchor=north east] at (4.650,2.500) {$0<t\leq t_1$};
			\node[sgn] at (2.325,1.103) {$w_i>0$};
		\end{scope}
		\begin{scope}[shift={(5.650,0.000)}]
			\draw[axis] (0,0) -- (4.650,0) node[below left, lab] {$x_i$};
			\draw[axis] (0,0) -- (0,2.450) node[above right, lab, xshift=-2pt] {$x_n$};
			\node[below, lab] at (0,0) {$1$};
			\node[lab, anchor=north east] at (4.650,2.500) {$t=\tau_1$};
			\draw[zero] plot[smooth, tension=0.35] coordinates {(0.025,0.237) (0.025,0.237) (0.027,0.237) (0.030,0.237) (0.033,0.237) (0.038,0.237) (0.044,0.236) (0.050,0.236) (0.058,0.236) (0.066,0.235) (0.075,0.235) (0.085,0.233) (0.096,0.232) (0.107,0.230) (0.119,0.228) (0.131,0.226) (0.144,0.223) (0.157,0.219) (0.170,0.215) (0.184,0.211) (0.198,0.206) (0.212,0.201) (0.226,0.194) (0.239,0.187) (0.253,0.180) (0.267,0.172) (0.280,0.164) (0.293,0.155) (0.306,0.146) (0.318,0.136) (0.331,0.126) (0.342,0.116) (0.353,0.105) (0.364,0.095) (0.374,0.085) (0.383,0.074) (0.392,0.064) (0.400,0.054) (0.407,0.045) (0.414,0.037) (0.420,0.028) (0.425,0.021) (0.429,0.015) (0.433,0.010) (0.436,0.005) (0.438,0.002) (0.439,0.001) (0.439,0.000)};
			\node[td] at (0.439,0) {};
			\node[sgn] at (2.900,1.000) {$w_i>0$};
			\node[sgn, anchor=west] (n1) at (0.575,0.657) {$w_i<0$};
			\draw[lead, -{Stealth[length=2.5pt]}] (n1.west) -- (0.025,0.107);
		\end{scope}
		\begin{scope}[shift={(11.300,0.000)}]
			\draw[axis] (0,0) -- (4.650,0) node[below left, lab] {$x_i$};
			\draw[axis] (0,0) -- (0,2.450) node[above right, lab, xshift=-2pt] {$x_n$};
			\node[below, lab] at (0,0) {$1$};
			\node[lab, anchor=north east] at (4.650,2.500) {$t=\tau_2$};
			\draw[zero] plot[smooth, tension=0.35] coordinates {(0.025,0.631) (0.027,0.633) (0.033,0.636) (0.044,0.643) (0.059,0.651) (0.077,0.662) (0.099,0.674) (0.125,0.687) (0.153,0.700) (0.184,0.715) (0.218,0.729) (0.254,0.743) (0.292,0.756) (0.331,0.768) (0.372,0.778) (0.413,0.787) (0.455,0.795) (0.497,0.800) (0.540,0.803) (0.582,0.804) (0.624,0.802) (0.666,0.798) (0.707,0.791) (0.747,0.782) (0.786,0.769) (0.824,0.754) (0.861,0.735) (0.896,0.714) (0.931,0.689) (0.963,0.662) (0.995,0.631) (1.024,0.598) (1.052,0.562) (1.079,0.523) (1.104,0.481) (1.127,0.437) (1.148,0.391) (1.167,0.344) (1.185,0.295) (1.201,0.247) (1.215,0.200) (1.227,0.154) (1.237,0.112) (1.246,0.074) (1.252,0.043) (1.257,0.019) (1.260,0.005) (1.261,0.000)};
			\node[td] at (1.261,0) {};
			\node[sgn] at (2.900,1.000) {$w_i>0$};
			\node[sgn] at (0.581,0.338) {$w_i<0$};
		\end{scope}
		\begin{scope}[shift={(0.000,-3.650)}]
			\draw[axis] (0,0) -- (4.650,0) node[below left, lab] {$x_i$};
			\draw[axis] (0,0) -- (0,2.450) node[above right, lab, xshift=-2pt] {$x_n$};
			\node[below, lab] at (0,0) {$1$};
			\node[lab, anchor=north east] at (4.650,2.500) {$t=\tau_3$};
			\draw[zero] plot[smooth, tension=0.35] coordinates {(0.025,0.892) (0.029,0.896) (0.042,0.906) (0.063,0.923) (0.091,0.945) (0.127,0.971) (0.169,1.000) (0.217,1.032) (0.270,1.064) (0.327,1.096) (0.388,1.127) (0.451,1.156) (0.516,1.184) (0.583,1.209) (0.650,1.231) (0.718,1.249) (0.785,1.265) (0.853,1.276) (0.919,1.284) (0.984,1.288) (1.047,1.287) (1.109,1.283) (1.169,1.275) (1.228,1.262) (1.284,1.245) (1.338,1.223) (1.390,1.198) (1.439,1.167) (1.486,1.133) (1.531,1.094) (1.574,1.050) (1.614,1.003) (1.651,0.950) (1.687,0.893) (1.719,0.832) (1.750,0.767) (1.777,0.697) (1.803,0.623) (1.826,0.546) (1.846,0.467) (1.865,0.386) (1.880,0.305) (1.893,0.227) (1.904,0.155) (1.913,0.092) (1.919,0.042) (1.922,0.011) (1.924,0.000)};
			\node[td] at (1.924,0) {};
			\node[sgn] at (2.900,1.000) {$w_i>0$};
			\node[sgn] at (0.879,0.537) {$w_i<0$};
		\end{scope}
		\begin{scope}[shift={(5.650,-3.650)}]
			\draw[axis] (0,0) -- (4.650,0) node[below left, lab] {$x_i$};
			\draw[axis] (0,0) -- (0,2.450) node[above right, lab, xshift=-2pt] {$x_n$};
			\node[below, lab] at (0,0) {$1$};
			\node[lab, anchor=north east] at (4.650,2.500) {$t=\tau_4$};
			\draw[zero] plot[smooth, tension=0.35] coordinates {(0.025,1.098) (0.032,1.105) (0.055,1.126) (0.092,1.158) (0.142,1.198) (0.203,1.245) (0.274,1.295) (0.354,1.346) (0.439,1.396) (0.530,1.443) (0.624,1.487) (0.720,1.527) (0.816,1.563) (0.913,1.593) (1.009,1.617) (1.104,1.636) (1.197,1.650) (1.288,1.658) (1.376,1.659) (1.461,1.656) (1.544,1.646) (1.623,1.631) (1.699,1.611) (1.773,1.585) (1.842,1.553) (1.909,1.517) (1.972,1.475) (2.032,1.429) (2.089,1.377) (2.143,1.321) (2.193,1.261) (2.241,1.196) (2.285,1.127) (2.327,1.054) (2.365,0.977) (2.400,0.896) (2.433,0.811) (2.462,0.723) (2.489,0.632) (2.512,0.540) (2.533,0.446) (2.551,0.353) (2.566,0.263) (2.579,0.179) (2.588,0.107) (2.595,0.049) (2.599,0.013) (2.601,0.000)};
			\node[td] at (2.601,0) {};
			\node[sgn] at (3.501,1.130) {$w_i>0$};
			\node[sgn] at (1.184,0.692) {$w_i<0$};
		\end{scope}
		\begin{scope}[shift={(11.300,-3.650)}]
			\draw[axis] (0,0) -- (4.650,0) node[below left, lab] {$x_i$};
			\draw[axis] (0,0) -- (0,2.450) node[above right, lab, xshift=-2pt] {$x_n$};
			\node[below, lab] at (0,0) {$1$};
			\node[lab, anchor=north east] at (4.650,2.500) {$t=t_\infty$};
			\draw[zero] plot[smooth, tension=0.35] coordinates {(0.025,1.313) (0.050,1.337) (0.123,1.401) (0.236,1.491) (0.381,1.590) (0.548,1.687) (0.729,1.776) (0.916,1.851) (1.105,1.912) (1.291,1.957) (1.473,1.986) (1.648,1.999) (1.817,1.997) (1.978,1.979) (2.132,1.948) (2.278,1.903) (2.417,1.847) (2.549,1.783) (2.673,1.712) (2.791,1.637) (2.903,1.561) (3.008,1.486) (3.108,1.413) (3.202,1.343) (3.291,1.277) (3.374,1.215) (3.453,1.158) (3.527,1.106) (3.596,1.058) (3.661,1.014) (3.722,0.974) (3.779,0.937) (3.831,0.904) (3.880,0.874) (3.925,0.847) (3.967,0.822) (4.004,0.800) (4.039,0.780) (4.069,0.763) (4.097,0.747) (4.121,0.734) (4.141,0.722) (4.159,0.713) (4.173,0.705) (4.184,0.699) (4.192,0.695) (4.197,0.692) (4.198,0.692)};
			\draw[zero, densely dashed] (4.198,0.692) -- (4.570,0.572);
			\node[sgn] at (3.953,1.430) {$w_i>0$};
			\node[sgn] at (1.648,0.840) {$w_i<0$};
		\end{scope}
	\end{tikzpicture}
	\caption{A schematic illustration of flow reversal near the boundary}\label{fig-flow-reversal}
\end{figure}

Chang, Kang, and Min \cite{CKM2026} established interior flow reversal for normal boundary data with \(\phi(s)=(1-s)^a\) on \((\frac12,1)\), where \(-1<a\leq1\).
For the boundary data considered here, Proposition \ref{prop-260046} gives the expansions in the normal coordinate
\begin{align}\label{intro:boundary-taylor}
w_i(x,t)=\partial_{x_n}w_i(x',0,t)x_n+o(x_n)\quad(i<n),\qquad w_n(x,t)=\frac{1}{2}\partial_{x_n}^2w_n(x',0,t)x_n^2+o(x_n^2),
\end{align}
uniformly on compact subsets of \((\Sigma\setminus\supp\psi)\times(0,\infty)\) as \(x_n\to0^+\).
Thus a nonzero normal derivative in \eqref{intro:boundary-taylor} determines the sign of the corresponding velocity component sufficiently close to the boundary.
The sign changes proved in this paper therefore describe separation in the sense of Definition \ref{260000}, and Figure~\ref{fig-flow-reversal} sketches the resulting flow reversal near the boundary.
We leave the study of interior flow reversal for the boundary data considered here to future work.

\subsection{Organization of the paper}
This paper is organized as follows.
Section \ref{260029} contains preliminary results used throughout the paper.
Section \ref{260049} treats normal boundary data, with Subsections \ref{260050} and \ref{260056} proving Theorems \ref{thm-0324-01} and \ref{thm-0507-01}, respectively.
Section \ref{secttan} treats tangential boundary data.
Subsections \ref{260058}, \ref{260092}, and \ref{260114} consider a tangential velocity component in the direction of the boundary data, a tangential velocity component in a different direction, and the normal velocity component, respectively.
They prove Theorems \ref{thm-260010}, \ref{thm-260014}, and \ref{thm-260017}.
Section \ref{260135} derives the boundary pressure formulas and proves Theorems \ref{thm-260019}, \ref{thm-260021}, and \ref{thm-260027}.
The appendices collect auxiliary results on Bessel functions, parabolic H\"{o}lder estimates, and boundary limits for Newtonian and heat potentials, together with the proofs of the preliminary results stated in Subsections \ref{260031}--\ref{subsect-0827-01}.

\section{Preliminaries}\label{260029}

\subsection{Confluent hypergeometric function \texorpdfstring{\({}_1F_1\)}{1F1}}

For parameters \(a,b\in\mathbb{C}\) with \(b\notin\{0,-1,-2,\dots\}\), the confluent hypergeometric function of the first kind is denoted by \({}_1F_1(a,b;z)\).
It is the unique solution of Kummer's equation
\begin{equation}\label{0327-03}
	zF''(z)+(b-z)F'(z)-aF(z)=0
\end{equation}
that is analytic in a neighborhood of \(z=0\) and satisfies \(F(0)=1\).
This solution is given by
\begin{equation}\label{260030}
	{}_1F_1(a,b;z) = \sum_{k=0}^\infty \frac{(a)_k}{(b)_k}\frac{z^k}{k!},
\end{equation}
where
\begin{equation*}
	(a)_k =
	\begin{cases}
		1,& k=0,\\
		a(a+1)\cdots(a+k-1), & k\geq 1.
	\end{cases}
\end{equation*}
The series converges for every \(z\in\mathbb C\).
The following properties are from \cite[Chapter 13]{Abram} and \cite[\S\S13.2--13.4, 13.7]{DLMF}.
\begin{proposition}\label{prop-0324-10}
	The following hold.
	\begin{enumerate}
		\item For \(b>a>0\),
		\begin{align}\label{0324-14}
			{}_1F_1(a,b;z) = \frac{\Gamma(b)}{\Gamma(a)\Gamma(b-a)}\int_0^1 e^{zt}t^{a-1}(1-t)^{b-a-1}\,dt.
		\end{align}
		\item Kummer's transformation is
		\begin{align}\label{0327-04}
			{}_1F_1(a,b;z)=e^z{}_1F_1(b-a,b;-z).
		\end{align}
		\item\label{item-0324-04} For \(z\in\mathbb C\) and parameters for which both sides are defined, the following hold, where \(a>0\), \(b>1\), \(c\geq0\), and \(s>0\) in \eqref{0911-21}.
		\begin{align}
			& \frac{d}{dz}{}_1F_1(a,b;z) = \frac{a}{b}{}_1F_1(a+1, b+1;z),\label{0327-07}\\
			& z\frac{d}{dz}{}_1F_1(a,b;z) + (b-1){}_1F_1(a,b;z) = (b-1){}_1F_1(a,b-1;z),\label{0423-10}\\
			& z\frac{d}{dz} {}_1F_1(a,b;z) + a{}_1F_1(a,b;z) = a {}_1F_1(a+1,b;z),\label{0327-01}\\
			& b{}_1F_1(a,b;z)=(b-a){}_1F_1(a,b+1;z) + a {}_1F_1 (a+1, b+1; z),\label{0327-02}\\
			& {}_1F_1(a,b;z) = {}_1F_1(a, b+1;z) + \frac{az}{b(b+1)}{}_1F_1(a+1, b+2;z),\label{0506-04}\\
			& \frac{d}{ds}\left[s^{-(b-1)}{}_1F_1\left(a,b;-\frac{c}{s}\right)\right] =-(b-1)s^{-b}{}_1F_1\left(a,b-1;-\frac{c}{s}\right).\label{0911-21}
		\end{align}
		\item\label{item-0504-01} Let \(a,b\in\R\) satisfy \(a,b,b-a\notin\{0,-1,-2,\ldots\}\).
		As \(z\to+\infty\),
		\begin{equation}\label{0327-10}
			\begin{aligned}
				{}_1F_1(a,b;z)& =\frac{\Gamma(b)}{\Gamma(a)}e^zz^{a-b} \left(1+O(z^{-1})\right),\\
				{}_1F_1(a,b;-z) & =\frac{\Gamma(b)}{\Gamma(b-a)}z^{-a} \left(1+O(z^{-1})\right).
			\end{aligned}
		\end{equation}
		More precisely, there exists \(z_0(a,b)>0\) such that, for \(z\geq z_0(a,b)\),
		\begin{align}\label{0504-06}
			{}_1F_1(a,b;-z) = \frac{\Gamma(b)}{\Gamma(b-a)}z^{-a}\left(1 + \frac{a(a-b+1)}{z} + \frac{a(a+1)(a-b+1)(a-b+2)}{2z^2}\right) + R_{a,b}(z),
		\end{align}
		where \(|R_{a,b}(z)|\leq C_{a,b}z^{-a-3}\) for some \(C_{a,b}>0\).
		\item\label{item-0514-01} For \(a,b>0\) and \(z>0\),
		\begin{align}\label{0514-02}
			|{}_1F_1(a,b;-z)|\leq \frac{C}{(1+z)^a},
		\end{align}
		where \(C>0\) depends only on \(a\) and \(b\).
	\end{enumerate}
\end{proposition}
\subsection{Some properties of the kernels involving the confluent hypergeometric functions}\label{260031}

For \(u\geq0\), define
\begin{align}
	 & Q(u):= \frac{\sqrt{\pi}}{2\Gamma\left(\frac{n+3}{2}\right)}u^{\frac{n}{2}}{}_1F_1\left(\frac{n}{2}, \frac{n+3}{2};-u\right),\label{0827-10}\\
	 & P(u):= {}_1F_1\left(\frac{n}{2}, \frac{n+3}{2};-u\right),\quad L(u):= \ln P(u),\label{0827-11}\\
	 & \mathcal{K}(u):= {}_1F_1\left(\frac{n}{2},\frac{n+1}{2}; - u\right) - \frac{2}{n+1}{}_1F_1\left(\frac{n}{2}, \frac{n+3}{2}; - u \right),\label{0827-12} \\
	 & \mathcal{H}(u):= \frac{n}{2}\mathcal{K}(u) + u \mathcal{K}'(u).\label{0827-13}
\end{align}
The integral representation \eqref{0324-14} shows that \(P(u)>0\) for \(u\geq0\), so \(L\) is well-defined.
Whenever \(\mathcal{K}(u)<0\), set
\begin{align}\label{0827-14}
	\mathcal{L}(u):=\ln (-\mathcal{K}(u)),\qquad \mathcal{M}(u):=-(u\mathcal{L}'(u))'.
\end{align}
Lemmas \ref{lem-0424-01}--\ref{lem-260032} are proved in Appendix \ref{260211}.
\begin{lemma}\label{lem-0424-01}
	Let \(n\geq 2\) and let \(Q\) be given in \eqref{0827-10}.
	Then the following hold.
	\begin{enumerate}
		\item\label{item-0428-01} \(Q(0)=0\), \(Q\) is strictly increasing on \((0,\infty)\), and \(Q(u)\to1\) as \(u\to\infty\).
    \item For \(u>0\),
    \begin{align}\label{0322-112}
        Q(u)< \min \left\{1, \frac{\sqrt{\pi}}{2\Gamma\left(\frac{n+3}{2}\right)}u^{\frac{n}{2}}\right\}.
    \end{align}
    \item For \(u>0\),
	\begin{align}\label{0322-12}
		\left|Q(u)-1+ \frac{n}{4u}\right|\leq \frac{n(n+2)}{8u^2}.
	\end{align}
	\end{enumerate}
\end{lemma}

\begin{lemma}\label{lem-0324-01}
	Let \(n\geq2\) and let \(P\) and \(L\) be given in \eqref{0827-11}. Then the following hold.
	\begin{enumerate}
		\item\label{item-0326-02} \(P(0)=1\), \(P\) is strictly decreasing on \([0,\infty)\), and
		\begin{align*}
			\frac{n+1}{2}P(u)+uP'(u)>0,\qquad u>0.
		\end{align*}
		\item\label{item-0326-03} For \(u\geq0\),
		\begin{align}\label{0529-03}
			1\leq (1 + 2u)^{\frac{n}{2}}P(u)\leq \frac{2^{\frac{n+2}{2}}}{\sqrt{\pi}}\Gamma\left(\frac{n+3}{2}\right).
		\end{align}
		\item\label{item-0327-01} For \(u>0\),
		\begin{align}\label{0529-02}
			-\frac{n}{1+2u}<L'(u)< -\frac{n}{(n+3)(1+2u)}.
		\end{align}
		\item\label{item-0327-02} For \(u>0\),
		\begin{align}\label{0602-01}
			0<L''(u)<\frac{2n}{(1+2u)^2}.
		\end{align}
		\item\label{item-0324-02} \(uL'(u)\) is strictly decreasing on \((0,\infty)\).
		For \(a_n:=\frac{(n+3)(n+5)}{6}\) and \(u>0\),
		\begin{align}\label{0529-04}
			- \frac{n(n+3)(n+5)^2}{36\left(u+a_n\right)^2}< (uL'(u))'< - \frac{n}{4\left(u+a_n\right)^2}.
		\end{align}
		In particular, for all \(u,v>0\),
		\begin{align}\label{0731-32}
			\frac{(uL'(u))'}{(vL'(v))'}\leq \frac{l_n}{9n}\,\frac{\left(v+a_n\right)^2}{\left(u+a_n\right)^2}, \qquad l_n:=n(n+3)(n+5)^2.
		\end{align}

		\item\label{item-0326-01}\(u^2L''(u)\) is strictly increasing on \((0,\infty)\).
	\end{enumerate}
\end{lemma}

\begin{lemma}\label{lem-0618-01}
		Let \(n\geq2\) and let \(\mathcal{K}\) and \(\mathcal{H}\) be given in \eqref{0827-12} and \eqref{0827-13}, respectively.
	Then there exists \(C>0\), depending only on \(n\), such that
	\begin{align*}
		|\mathcal{K}(u)|\leq \frac{C}{(1+u)^{\frac{n}{2}}},\qquad |\mathcal{H}(u)|\leq \frac{C}{(1+u)^{\frac{n+2}{2}}},\qquad u>0.
	\end{align*}
\end{lemma}

We use the following positive zeros, with \(\mathcal K\) given in \eqref{0827-12}.
\begin{enumerate}
	\item \(z_1\) is the unique positive zero of \({}_1F_1(\frac n2,\frac{n-1}{2};-z)=e^{-z}{}_1F_1(-\frac12,\frac{n-1}{2};z)=\frac{2}{n-1}\left(\frac{n+1}{2}\mathcal K(z)+z\mathcal K'(z)\right)\).
  \item \(z_2\) is the unique positive zero of \({}_1F_1(\frac{n+2}{2},\frac{n+1}{2};-z) =e^{-z}{}_1F_1(-\frac12,\frac{n+1}{2};z)\).
  \item \(z_3\) is the unique positive zero of \(\mathcal K(z)\).
  \item \(z_4\) is the unique positive zero of \(\frac n2\mathcal K(z)+z\mathcal K'(z) =\frac n2{}_1F_1(\frac{n+2}{2},\frac{n+1}{2};-z) -\frac n{n+1}{}_1F_1(\frac{n+2}{2},\frac{n+3}{2};-z)\).
  \item \(z_5\) is the larger positive zero of \({}_1F_1(\frac{n+2}{2},\frac{n-1}{2};-z) =e^{-z}{}_1F_1(-\frac32,\frac{n-1}{2};z)\).
  \item \(z_6\) is the largest positive zero of \( \frac{d}{dz}\left(z^{\frac{n+1}{2}} {}_1F_1\left(\frac n2,\frac{n-1}{2};-z\right)\right)\).
	\item \(z_7\) is the unique positive zero of \(n\,{}_1F_1(\frac{n+2}{2},\frac{n+1}{2};-z)-(n-1)\,{}_1F_1(\frac n2,\frac{n+1}{2};-z) =e^{-z}\bigl[n\,{}_1F_1(-\frac12,\frac{n+1}{2};z)-(n-1)\,{}_1F_1(\frac12,\frac{n+1}{2};z)\bigr]\).
\end{enumerate}
The first equality in the definition of \(z_1\) and the equalities in the definitions of \(z_2\), \(z_5\), and \(z_7\) follow from \eqref{0327-04}.
The last equality in the definition of \(z_1\) follows from \eqref{0515-01}.
Theorem \ref{thm-0430-01} gives the uniqueness of \(z_1,z_2\) and shows that \({}_1F_1(\frac{n+2}{2},\frac{n-1}{2};-z)\), which defines \(z_5\), has exactly two positive zeros.
The uniqueness of \(z_7\) is proved in Lemma \ref{lem-261080}{\rm(i)}.
The following lemmas give the remaining properties of these zeros.

\begin{lemma}\label{lem-0730-50}
	Let \(n\geq2\).
	Then \(z_1\) and \(z_2\) satisfy
	\begin{align}\label{0701-05}
		0<\sqrt{z_2}-\sqrt{z_1}<1.
	\end{align}
\end{lemma}

The uniqueness of the zeros \(z_3,z_4\) follows from the next lemma.
It is proved in Appendix \ref{260216}.
\begin{lemma}\label{lem-0514-01}
	Let \(n\geq2\) and let \(\mathcal{K}\) be given in \eqref{0827-12}.
	Then there exist unique \(z_3,z_4>0\), depending only on \(n\), such that the following hold.
	\begin{enumerate}
		\item \(\mathcal{K}(z)>0\) for \(0<z<z_3\), and \(\mathcal{K}(z)<0\) for \(z>z_3\).
\item \(\frac{n+1}{2}\mathcal{K}(z)+z\mathcal{K}'(z)>0\) for \(0<z<z_1\), and \(\frac{n+1}{2}\mathcal{K}(z)+z\mathcal{K}'(z)<0\) for \(z>z_1\).
\item \(\frac n2\mathcal{K}(z)+z\mathcal{K}'(z)>0\) for \(0<z<z_4\), and \(\frac n2\mathcal{K}(z)+z\mathcal{K}'(z)<0\) for \(z>z_4\).
	\end{enumerate}
	Moreover,
	\begin{align}\label{0623-02}
		z_4<z_1<z_3.
	\end{align}
\end{lemma}

\begin{lemma}\label{lem-0520-01}
	Let \(n\geq2\) and let \(\mathcal{K}\) be given in \eqref{0827-12}.
	For \(z>z_3\), let \(\mathcal{L}(z)\) and \(\mathcal{M}(z)\) be given in \eqref{0827-14}.
	Then there exists \(Z_0>z_3\), depending only on \(n\), such that the following hold.
	\begin{enumerate}
		\item For \(z\geq Z_0\), \(\mathcal{L}'(z)<0\), and \(\mathcal{M}(z)\) is positive and strictly decreasing.
		Moreover, 
		\begin{align}\label{0520-31}
			\frac{3n}{4z^2}\leq \mathcal{M}(z)\leq \frac{3n}{2z^2}.
		\end{align}
		\item For \(\zeta_2>\zeta_1>Z_0\),
		\begin{align}\label{0520-30}
			\frac{\mathcal{K}(\zeta_1)}{\mathcal{K}(\zeta_2)}<\left(\frac{\zeta_2}{\zeta_1}\right)^{\frac{n}{2}}.
		\end{align}
		\item There exists \(C_{\mathcal{L}}>0\), depending only on \(n\), such that, for every \(z>Z_0\), 
		\begin{align}\label{0520-40}
			|\mathcal{L}'(z)|+2 z|\mathcal{L}''(z)|\leq \frac{C_{\mathcal{L}}}{z}.
		\end{align}
	\end{enumerate}
\end{lemma}

\begin{lemma}\label{lem-260032}
	Let \(n\geq2\) and
	\begin{align*}
		F(z):=z^{\frac{n+1}{2}}{}_1F_1\left(\frac n2,\frac{n-1}{2};-z\right), \qquad z>0.
	\end{align*}
	Then \(F'\) has a largest positive zero \(z_6\), and \(F'(z)<0\) for \(z>z_6\).
	Moreover, there exists \(C=C(n)>0\) such that
	\begin{align*}
	    |F'(z)|\leq Cz^{-1/2},\qquad z>0.
	\end{align*}
\end{lemma}

\subsection{Uniqueness of zeros for Volterra-type operators}\label{260039}
For functions on positive times, we denote their time convolution by
\begin{align*}
	(K*f)(t):=\int_0^t K(t-s)f(s)\,ds
\end{align*}
whenever the integral converges.
When a lemma assumes bounds on \(K\) and \(K'\), we write these as
\begin{align}\label{0829-01}
	|K(t)| \leq \frac{C}{t^{\sigma}}, \quad |K'(t)|\leq \frac{C}{t^{1+\sigma}},
\end{align}
on the time interval specified in each statement. Lemmas \ref{lem-0424-20}--\ref{lem-0504-01} are proved in Appendix \ref{app-volterra}.

\begin{lemma}\label{lem-0424-20}
	Let \(K\in C^1(0,1)\) satisfy \eqref{0829-01} on \((0,1)\) for some \(\sigma\in(0,1)\).
	Let \(f\in C^{\alpha}([0,1])\) with \(\alpha\in(\sigma,1)\).
	Then \(K*f\in C^{1,\alpha-\sigma}_{\mathrm{loc}}(0,1)\), and
	\begin{align}\label{0325-02}
		\frac{d}{dt}(K*f)(t)=K(t)f(t) + \int_0^t K'(t-s)(f(s)-f(t))\,ds.
	\end{align}
	Suppose that \(K=K(r,u)\in C^1(J\times(0,1))\) for an open interval \(J\) and that, for each compact \(I\subset J\),
	\begin{align*}
		|K(r,u)|+|\partial_rK(r,u)|\leq C_Iu^{-\sigma},\qquad
		|\partial_uK(r,u)|\leq C_Iu^{-1-\sigma}
		\qquad(r\in I,\ 0<u<1).
	\end{align*}
	Then \((r,t)\mapsto(K(r,\cdot)*f)(t)\) is jointly \(C^1\).
\end{lemma}

\begin{lemma}\label{lem-0325-01}
	Let \(K\in C^1(0,\infty)\) be positive and strictly decreasing, and satisfy \eqref{0829-01} on \((0,1)\) for some \(\sigma\in(0,1)\).
	Let \(\phi\in C^1([0,1])\) be log-concave on \((0,1)\) with \(\phi(0)=0\) and \(\phi>0\) on \((0,1)\). Then the following hold.
	\begin{enumerate}
		\item \((K*\phi')/\phi\) is strictly decreasing on \((0,1)\). In particular, \(K*\phi'\) has at most one zero.
		\item If \(\phi'\in C^\alpha([0,1])\) for some \(\sigma<\alpha<1\), then \((K*\phi')'(t_*)<0\) at each zero \(t_*\in(0,1)\).
	\end{enumerate}
\end{lemma}

\begin{lemma}\label{lem-0506-01}
Let \(K_1\in C^1(0,1)\cap L^1(0,1)\) and \(K_2\in C(0,1)\cap L^1(0,1)\) with \(K_1'<0\) and \(K_2>0\) on \((0,1)\), and suppose that \(K_1\) satisfies \eqref{0829-01} on \((0,1)\) for some \(\sigma\in(0,1)\).
Let \(\phi\in C^1([0,1])\) be log-concave on \((0,1)\) with \(\phi(0)=0\) and \(\phi>0\) on \((0,1)\), and define, for \(t\in(0,1)\),
	\begin{equation}
		(T\phi)(t):=(K_1*\phi')(t)-(K_2*\phi)(t).
	\end{equation}
	Then the following hold.
	\begin{enumerate}
		\item \((T\phi)/\phi\) is strictly decreasing on \((0,1)\). In particular, \(T\phi\) has at most one zero in \((0,1)\), with positive values before any zero and negative values after it.
		\item If \(\phi'\in C^\alpha([0,1])\) for some \(\sigma<\alpha<1\), then \((T\phi)'(t_*)<0\) at each zero \(t_*\in(0,1)\).
	\end{enumerate}
\end{lemma}

\begin{lemma}\label{lem-0504-01}
	Let \(K\in C^1(0,\infty)\) have exactly one zero \(\tau_0>0\), with \(K(\tau)<0\) for \(0<\tau<\tau_0\) and \(K(\tau)>0\) for \(\tau>\tau_0\).
	Let \(\phi\in C^1([0,1])\) be log-concave on \((0,1)\) with \(\phi(0)=\phi(1)=0\) and \(\phi>0\) on \((0,1)\), extended by zero for \(s>1\).
	Then, for \(t>1\),
	\begin{equation}
		(K*\phi)(t)=\int_0^1 K(t-s)\phi(s)\,ds,
	\end{equation}
	and \(K*\phi\) has at most one zero in \((1,\infty)\).
	Moreover, if \((K*\phi)(t_*)=0\) for some \(t_*>1\), then \((K*\phi)'(t_*)>0\).
\end{lemma}

\subsection{Some integral identities}\label{subsect-0827-01}
Proposition \ref{prop-0322-001} and Lemma \ref{lem-0422-03} are proved in Appendix \ref{260240}.
\begin{proposition}\label{prop-0322-001}
	Let \(n\geq2\) and \(1\leq i\leq n-1\).
	For every \(x'\in\Sigma\) and \(t>0\),
	\begin{align}\label{0320-01}
		\mathrm{p.v.}\int_{\Sigma}\Gamma'(x'-z',t)\partial_{z_i}N(z',0)\,dz' = - \frac{\Gamma\left( \frac{n}{2}\right) }{2\Gamma\left( \frac{n+1}{2}\right)} \sqrt{t}\partial_{x_i} \Gamma'(x',t){}_1F_1\left( \frac{1}{2}, \frac{n+1}{2}; \frac{|x'|^2}{4t}\right).
	\end{align}
\end{proposition}

We also use the following formula in spherical coordinates. For \(n\geq3\) and every \(G\in L^1((-1,1); (1-\xi^2)^{\frac{n-4}{2}}d\xi)\),
\begin{align}\label{0907-01}
	\int_{\mathbb S^{n-2}}G(\omega_1)\,d\omega =|\mathbb S^{n-3}|\int_{-1}^1 G(\xi)(1-\xi^2)^{\frac{n-4}{2}}\,d\xi.
\end{align}

\begin{lemma}\label{lem-0422-03}
	Let \(n\geq2\) and \(0\leq s<r\).
	Then
	\begin{align}\label{0422-07}
		\int_{\mathbb S^{n-2}} \frac{d\sigma(\omega)} {(r^2+s^2-2rs\omega_1)^{\frac{n-2}{2}}} = |\mathbb S^{n-2}|\,r^{-(n-2)} \sum_{k=0}^\infty \frac{ \left(\frac{n-2}{2}\right)_k \left(\frac12\right)_k }{ k!\left(\frac{n-1}{2}\right)_k }\left(\frac{s}{r}\right)^{2k}.
	\end{align}
\end{lemma}

\subsection{Properties of the kernel \texorpdfstring{\(L_{ni}\)}{Lni}}
The next lemma gives properties of the kernel \(L_{ni}\).
\begin{lemma}\label{lem-260040}
		Let \(n\geq2\), \(1\leq i\leq n-1\), and let \(L_{ni}\) be given in \eqref{0419-02}.
	Then, for every \(x_n>0\) and \(t>0\), the following hold.
	\begin{enumerate}
		\item\label{260041}
		\begin{align}\label{0708-22}
			\int_{\Sigma}L_{ni}(x,t)\,dx'=0.
		\end{align}
		\item\label{260042}
		\begin{align}\label{0708-17}
			\int_{\Sigma}|x'|\,|L_{ni}(x,t)|\,dx' \leq \frac{C}{\sqrt{t}},
		\end{align}
		where \(C>0\) depends only on \(n\).
		\item\label{260043} For each fixed \(t>0\),
		\begin{align}\label{0708-18}
			\lim_{x_n\to0^+} \int_{\Sigma}|x'|\,|L_{ni}(x,t)|\,dx' = 0.
		\end{align}
	\end{enumerate}
\end{lemma}
\begin{proof}
First consider the convolution
\begin{align*}
    A_N(x,t)=\int_\Sigma\Gamma'(y',t)\partial_{x_n}N(x'-y',x_n)\,dy',
    \qquad x\in\R_+^n,\quad t>0.
\end{align*}
We claim that its tangential derivative satisfies
\begin{align}\label{0708-35}
	|\partial_{x_i}A_N(x,t)| \leq C \frac{x_n+\sqrt{t}}{(|x|^2+t)^{\frac{n+1}{2}}},
\end{align}
where \(C>0\) depends only on \(n\).
Using the representation \(\partial_{x_n}N(x)=-\int_0^\infty\partial_{x_n}\Gamma(x,s)\,ds\), we obtain
\begin{align}\label{0827-02}
	\partial_{x_i}A_N(x,t) & = -\int_0^\infty \partial_{x_n}\Gamma_1(x_n,s)\partial_{x_i}\Gamma'(x',t+s)\,ds.
\end{align}

Since \(|\partial_{x_i}\Gamma'(x',t+s)| \leq C|x'|(t+s)^{-\frac{n+1}{2}} e^{-\frac{|x'|^2}{4(t+s)}}\), \eqref{0827-02} gives
\begin{align}\label{0827-03}
	|\partial_{x_i}A_N(x,t)| & \leq Cx_n|x'|\int_0^\infty \frac{e^{-\frac{x_n^2}{4s}-\frac{|x'|^2}{4(t+s)}}} {s^{\frac32}(t+s)^{\frac{n+1}{2}}}\,ds =Cx_n|x'|(I_1+I_2),
\end{align}
where \(I_1\) and \(I_2\) denote the integrals over \((0,t)\) and \((t,\infty)\), respectively.

On \(0<s<t\), we have \(t\leq t+s\leq2t\), and the factor involving \(x_n\) satisfies \(e^{-x_n^2/(4s)}\leq e^{-x_n^2/(8t)}e^{-x_n^2/(8s)}\).
Hence
\begin{align}\label{0827-04}
	I_1 & \leq e^{-\frac{|x|^2}{8t}}\int_0^t \frac{e^{-\frac{x_n^2}{8s}}}{s^{\frac32}(t+s)^{\frac{n+1}{2}}}\,ds\leq \frac{e^{-\frac{|x|^2}{8t}}}{t^{\frac{n+1}{2}}}\int_0^t \frac{e^{-\frac{x_n^2}{8s}}}{s^{\frac32}}\,ds\leq \frac{Ce^{-\frac{|x|^2}{8t}}}{x_nt^{\frac{n+1}{2}}}\leq \frac{C\sqrt{t}}{x_n|x|(|x|^2+t)^{\frac{n+1}{2}}},
\end{align}
where we have used \(e^{-r^2/8}\leq\frac{C}{r(1+r^2)^{(n+1)/2}}\) with \(r=|x|/\sqrt t\).

On \(s\geq t\), we have \(s\leq t+s\leq2s\), so the product of the two exponential factors in \eqref{0827-03} is bounded by \(e^{-|x|^2/(8s)}\).
Setting \(u=1/s\), we obtain
\begin{align}\label{0827-05}
	I_2 & \leq \int_t^\infty \frac{e^{-\frac{|x|^2}{8s}}}{s^{\frac32}(s+t)^{\frac{n+1}{2}}}\,ds = \int_0^{\frac{1}{t}}\frac{u^{\frac{n}{2}}e^{-\frac{|x|^2u}{8}}}{(1+tu)^{\frac{n+1}{2}}}\,du \leq \int_0^{\frac{1}{t}}u^{\frac{n}{2}}e^{-\frac{|x|^2u}{8}}\,du\leq \frac{C}{(|x|^2+t)^{\frac{n+2}{2}}}.
\end{align}
Using \(|x'|\leq|x|\leq(|x|^2+t)^{1/2}\) then gives \eqref{0708-35}.

\underline{Proof of (i).}\quad
By \eqref{0419-02},
\begin{align}\label{0708-23}
	L_{ni}(x,t) & = \partial_{x_i} \int_0^{x_n} \partial_{y_n}\Gamma_1(y_n,t) \int_{\Sigma} \Gamma'(y',t)\partial_{x_n}N(x-y)\,dy'dy_n\nonumber \\
	 & =\int_0^{x_n} \partial_{y_n}\Gamma_1(y_n,t) \partial_{x_i}A_N(x',x_n-y_n,t) \,dy_n.
\end{align}
Integrating \eqref{0708-35} in polar coordinates in \(x'\) gives \(\norm{\partial_{x_i}A_N(\cdot,x_n,t)}_{L^1(\Sigma)} \leq C/\sqrt t\), so Fubini's theorem permits interchanging the \(x'\)- and \(y_n\)-integrals in \eqref{0708-23}.
For the inner spatial integral, move the derivative onto \(\Gamma'\) in \(A_N\) to obtain
\begin{align*}
		\int_{\Sigma} \partial_{x_i}A_N(x,t)\,dx' = \left( \int_{\Sigma} \partial_{x_i}\Gamma'(x',t)\,dx' \right) \left( \int_{\Sigma}\partial_{x_n}N(y',x_n)\,dy' \right) = 0.
\end{align*}
Here, \(\partial_{x_i}\Gamma'(\cdot,t)\) has integral zero by oddness.
Substituting this cancellation into \eqref{0708-23} gives \eqref{0708-22}.

\underline{Proof of (ii).}\quad
For \(0<y_n<x_n\), \eqref{0708-35} and polar coordinates give
\begin{align}\label{0827-06}
	\int_{\Sigma}|x'||\partial_{x_i}A_N(x',x_n-y_n,t)|\,dx' \leq C\frac{x_n-y_n+\sqrt t}{\sqrt{(x_n-y_n)^2+t}}\leq C.
\end{align}
Inserting this uniform first-moment bound into \eqref{0708-23} gives
\begin{align}\label{0708-19}
	\int_{\Sigma}|x'|\,|L_{ni}(x,t)|\,dx' & \leq \int_0^{x_n} |\partial_{y_n}\Gamma_1(y_n,t)| \int_{\Sigma}|x'|\, |\partial_{x_i}A_N(x',x_n-y_n,t)| \,dx'dy_n \leq C \int_0^{x_n} |\partial_{y_n}\Gamma_1(y_n,t)|\,dy_n\leq \frac{C}{\sqrt{t}},
\end{align}
which proves \eqref{260042}.

\underline{Proof of (iii).}\quad
For fixed \(t>0\), \(y_n\mapsto|\partial_{y_n}\Gamma_1(y_n,t)|\) is continuous at \(y_n=0\).
Hence the last integral in \eqref{0708-19}, taken over \((0,x_n)\), tends to zero as \(x_n\to0^+\).
This proves \eqref{260043}.
\end{proof}
\subsection{Real analyticity}
\begin{lemma}\label{lem-260044}
	Let \(n\geq2\), let \(E\subset\Sigma\) be bounded and measurable, let \(g\in L^1(E)\), and let \(x''\in(\mathbf e_1')^\perp\).
	Let \(I\subset\R\) be an open interval such that \(x'\notin\overline E\) whenever \(x_1\in I\).
	Then the following hold.
	\begin{enumerate}
		\item\label{1006-08} For every \(a>0\) and \(m\in\{0,1\}\),
		\begin{align*}
			\int_E\frac{(x_1-y_1)^m}{|x'-y'|^{2a}}g(y')\,dy'
		\end{align*}
		is real-analytic in \(x_1\in I\).
		\item\label{1006-09} For every \(T>0\), \(b\in\R\), \(m\in\{0,1\}\), and \(q\in L^\infty(0,T)\),
	\begin{align*}
			\int_0^T q(s)s^{-b}\int_E(x_1-y_1)^m e^{-\frac{|x'-y'|^2}{4s}}g(y')\,dy'ds
	\end{align*}
		is real-analytic in \(x_1\in I\).
		\item\label{1006-10} A real-analytic function on an open interval which is not identically zero has only finitely many zeros in each compact subinterval.
		\item\label{1006-11} Let \(T>0\), \(0<a<b\), \(c<a+1\), \(m\in\{0,1\}\), \(1\leq i\leq n-1\), and \(q\in C([0,T])\).
		Then
		\begin{align*}
			\int_E(x_i-y_i)^m g(y')\int_0^t\frac{q(t-s)}{s^c} {}_1F_1\left(a,b;-\frac{|x'-y'|^2}{4s}\right)\,dsdy'
		\end{align*}
		is continuous on \((\Sigma\setminus\overline E)\times(0,T]\) and real-analytic in \(x_1\in I\) for each fixed \(t\in(0,T]\).
		It is also real-analytic in \(r\) along every line \(x'=z'+r\omega'\), with \(z'\in\Sigma\) and \(\omega'\in\mathbb S^{n-2}\), on each open interval where \(z'+r\omega'\notin\overline E\).
		\item\label{1006-12} Let \(T>0\), \(a,b,c>0\), and \(q\in L^1(0,T)\).
		For each fixed \(x'\in\Sigma\),
		\begin{align*}
			\int_E g(y')\int_0^T\frac{q(s)}{(t-s)^c} {}_1F_1\left(a,b;-\frac{|x'-y'|^2}{4(t-s)}\right)\,dsdy'
		\end{align*}
		is real-analytic in \(t>T\).
		\item\label{261032} Let \(0<a<b\), \(1<c<a+2\), let \(P\) be a polynomial on \(\Sigma\), and suppose that \(q\in L^\infty(0,\infty)\) satisfies \(|q(s)|\leq C\min\{s,1\}\).
		Then
		\begin{align*}
			\int_E P(x'-y')g(y')\int_0^\infty\frac{q(s)}{s^c} {}_1F_1\left(a,b;-\frac{|x'-y'|^2}{4s}\right)\,dsdy'
		\end{align*}
		is real-analytic in \(r\) along every line \(x'=z'+r\omega'\) on each open interval where \(z'+r\omega'\notin\overline E\).
	\end{enumerate}
\end{lemma}
\begin{proof}
Fix a point \(x_0'\) of the given line with \(d=\dist(x_0',E)>0\), and write \(x_0'+h\omega'\) for nearby points, with \(\omega'=\mathbf e_1'\) in {\rm(i)} and {\rm(ii)}.
For \(|h|\leq d/8\),
\begin{align}\label{1006-14}
2|h||\omega'\cdot(x_0'-y')|+h^2 & \leq\tfrac12|x_0'-y'|^2,\qquad y'\in E.
\end{align}
In {\rm(i)}, {\rm(ii)}, {\rm(iv)}, and {\rm(vi)}, expand the kernel at \(x_0'+h\omega'\) in powers of \(h\) by the binomial and exponential series, after the integral representation \eqref{0324-14} in {\rm(iv)} and {\rm(vi)}.
By \eqref{1006-14}, the sums of the absolute values of the terms are bounded by \(C|x_0'-y'|^{-2a}\), \(Ce^{-d^2/(8s)}\), and \(Ce^{-d^2v/(8s)}\), respectively, since \(E\) is bounded.
The integral representation \eqref{0324-14} and \eqref{0514-02} give
\begin{align}\label{1006-13}
\frac1{s^c}\int_0^1e^{-\frac{d^2v}{8s}} v^{a-1}(1-v)^{b-a-1}\,dv \leq C(a,b)s^{-c}\left(1+\frac{d^2}{s}\right)^{-a} \leq C(a,b)d^{-2a}s^{a-c},
\end{align}
which is integrable near \(s=0\) under the stated conditions on \(c\).
In {\rm(vi)}, \(|q(s)|\leq C\min\{s,1\}\) and \(c>1\) give integrability at infinity.
Since \(g\in L^1(E)\), integrating term by term gives convergent power series in \(h\), and the same majorants and dominated convergence give the joint continuity in {\rm(iv)}.
For {\rm(v)}, fix \(t_0>T\) and \(|h|\leq(t_0-T)/2\).
The series \eqref{260030} and the binomial series expand the kernel at \(t_0+h\), with the sum of the absolute values of its terms bounded by
\begin{align*}
\left(\frac2{t_0-T}\right)^c {}_1F_1\left(a,b; \frac{(|x'|+\sup_{y'\in E}|y'|)^2}{2(t_0-T)}\right),
\end{align*}
uniformly for \(0\leq s\leq T\) and \(y'\in E\).
Finally, if zeros of a real-analytic function accumulate at a point of \(I\), all its derivatives vanish there, which gives {\rm(iii)}.
\end{proof}

\subsection{Boundary derivatives and separation points}
For the tangential velocity component in the direction of tangential boundary data, that is, for \(i=j<n\), and for \(x'\notin\supp\psi\), define
\begin{align}\label{260045}
T_{\rm van}(x'):=\bigl\{t>0:\partial_{x_n}w_i(x',0,t)=0\bigr\}.
\end{align}

The next proposition relates boundary derivatives of \(w_i\) to Definition \ref{260000}.
\begin{proposition}\label{prop-260046}
	Let \(1\leq i,j\leq n\), and let \(w\) solve the Stokes system \eqref{0308-01} with boundary data \eqref{0324-13}, where \(\psi\) satisfies {\rm\textbf{(S)}} and \(\phi\) satisfies {\rm\textbf{(T1)}}.
	Then the following hold.
	\begin{enumerate}
		\item\label{260047} If \((x_*',t_*)\) is a separation point of \(w_i\), then
		\begin{align}\label{1021-01}
			\partial_{x_n}^{\ell}w_i(x_*',0,t_*)=0,
		\end{align}
		where \(\ell=1\) if \(i<n\) and \(\ell=2\) if \(i=n\).
		Conversely, \((x_*',t_*)\in(\Sigma\setminus\supp\psi)\times(0,\infty)\) is a separation point of \(w_i\) if \(\partial_{x_n}^{\ell}w_i(\cdot,0,\cdot)\) takes both signs in every neighborhood of \((x_*',t_*)\).
		\item\label{260048} Let \(\ell\) be as in {\rm(i)}, let \(x'\notin\supp\psi\) and \(t_i^*(x')>0\), and suppose that
		\begin{align*}
			\partial_{x_n}^{\ell}w_i(x',0,t_i^*(x'))=0,\quad 
	\partial_{x_n}^{\ell}w_i(x',0,s)\,\partial_{x_n}^{\ell}w_i(x',0,t)<0 \quad (0<s<t_i^*(x')<t).
		\end{align*}
		Then \((x',t)\) is a separation point of \(w_i\) if and only if \(t=t_i^*(x')\).
		\item\label{1007-06} Fix \(t>0\) and a positive integer \(\ell\), let \(J\subset(0,\infty)\) be an open interval containing \(t\), and let \(U\subset\Sigma\setminus\supp\psi\) be open.
		Suppose that, for some \(h_0>0\), each \(\partial_{x_n}^m w_i\), \(1\leq m\leq\ell\), extends continuously to \(U\times[0,h_0)\times J\), and that
		\begin{align}\label{1021-04}
			\partial_{x_n}^m w_i(x',0,s)=0 \qquad(x'\in U,\ s\in J,\ 1\leq m<\ell).
		\end{align}
		Suppose also that \(f\in C(U)\) takes both signs in every neighborhood in \(U\) of each of its zeros and satisfies
		\begin{align}\label{1021-05}
			\sgn\bigl(\partial_{x_n}^{\ell}w_i(x',0,t)\bigr)=\sgn f(x')\qquad(x'\in U).
		\end{align}
		Then, for \(x'\in U\), \((x',t)\) is a separation point of \(w_i\) if and only if \(f(x')=0\).
	\end{enumerate}
\end{proposition}
\begin{proof}
Fix compact sets \(U\Subset\Sigma\setminus\supp\psi\) and \(J\Subset(0,\infty)\), choose \(T>\sup J\) and \(0<\gamma<\alpha\), and set \(d=\dist(U,\supp\psi)>0\).
We first show that \(\partial_{x_n}^{\ell}w_i\) extends continuously to \(U\times\{0\}\times J\) and that, for \(0<x_n\leq1\),
\begin{align}\label{260263}
	\sup_{(z',s)\in U\times J} \left| \frac{\ell!\,w_i(z',x_n,s)}{x_n^{\ell}} -\partial_{x_n}^{\ell}w_i(z',0,s) \right|\leq Cx_n^\gamma.
\end{align}

The H\"older bound for the Riesz transform used in Lemma \ref{lem-260265} gives
\begin{align}\label{260250}
	 & \sum_{|\nu|\leq1}\left( \|\partial_{x'}^\nu\psi\,\phi\|_{C^{1+\gamma,\frac{1+\gamma}{2}}(\Sigma\times[0,T])} +\sum_{k=1}^{n-1} \|R_k'(\partial_{x'}^\nu\psi)\,\phi\|_{C^{1+\gamma,\frac{1+\gamma}{2}}(\Sigma\times[0,T])} \right)\leq C_T\|\psi\|_{C^{2,\alpha}(\Sigma)} \|\phi\|_{C^1([0,T])}.
\end{align}
These functions vanish at \(s=0\), since \(\phi(0)=0\).
Applying \cite[Lemma~5.1]{CK2026} with \(k=1\) to \(w\) and to each \(\partial_{x_k}w\) gives
\begin{align}\label{260251}
	 & \|w\|_{C^{1+\gamma,\frac{1+\gamma}{2}}(\overline{\R_+^n}\times[0,T])} +\sum_{k=1}^{n-1} \|\partial_{x_k}w\|_{C^{1+\gamma,\frac{1+\gamma}{2}}(\overline{\R_+^n}\times[0,T])}\leq C_T\|\psi\|_{C^{2,\alpha}(\Sigma)} \|\phi\|_{C^1([0,T])}.
\end{align}
Thus \(\partial_{x_n}w_i\) for \(i<n\) and \(\partial_{x_n}\partial_{x_k}w_k\) for \(k<n\) extend continuously to the boundary.
By incompressibility,
\begin{align}\label{260252}
	\partial_{x_n}w_n=-\sum_{k=1}^{n-1}\partial_{x_k}w_k, \qquad \partial_{x_n}^2w_n=-\sum_{k=1}^{n-1} \partial_{x_n}\partial_{x_k}w_k,
\end{align}
so \(\partial_{x_n}^2w_n\) also extends continuously to the boundary, and \(\partial_{x_n}w_n(z',0,s)=0\) for \(z'\in U\), since \(U\) is disjoint from \(\supp\psi\).
The H\"older bound in \eqref{260251} gives
\begin{align}\label{260262}
	\sup_{(z',s)\in U\times J} \left|\partial_{x_n}^{\ell}w_i(z',x_n,s) -\partial_{x_n}^{\ell}w_i(z',0,s)\right| \leq Cx_n^\gamma.
\end{align}
Since \(w_i(z',0,s)=0\) on \(U\times J\),
\begin{align}\label{1007-19}
	w_i(z',x_n,s)=\int_0^{x_n}\partial_{x_n}w_i(z',r,s)\,dr\quad(i<n), \qquad w_n(z',x_n,s)=\int_0^{x_n}(x_n-r)\partial_{x_n}^2w_n(z',r,s)\,dr,
\end{align}
and \eqref{260262} gives \eqref{260263}.

We next identify these boundary values for tangential boundary data, that is, for \(j<n\).
Let \(v\) solve \eqref{0330-01} with boundary data \(f\phi\,\mathbf e_j\), where \(f=\psi\) or \(f=\partial_{x_k}\psi\) with \(k<n\).
Then \(v=w\) and \(v=\partial_{x_k}w\), respectively.

Extend \(\phi\) by zero to negative times and, for \(f\in C_c(\Sigma)\), let
\begin{align}\label{0703-05}
	k_f(y,t) := \int_0^t\int_\Sigma \partial_{y_n}\Gamma(y-z',t-s) f(z')\phi(s)\,dz'ds,\qquad k:=k_\psi.
\end{align}
For \(i<n\), Lemma~\ref{lem-260265} with exponent \(\gamma\) and \(\supp\partial_{x_j}f\subset\overline{B'_1}\) gives
\begin{align}\label{260253}
	\|R_i'(\partial_{x_j}f)\|_{C^\gamma(\Sigma)} \leq C\|f\|_{C^{1,\gamma}(\Sigma)} \leq C\|\psi\|_{C^{2,\alpha}(\Sigma)}.
\end{align}
Since the extension of \(\phi\) is Lipschitz and \(\int_\Sigma |y'|^\gamma\Gamma'(y',u)\,dy'=C_{n,\gamma}u^{\gamma/2}\), \eqref{260253} gives
\begin{align}\label{260254}
	\sup_{\substack{x'\in\Sigma \\
	0\leq s\leq T}} \left| \phi(s-u)\bigl(\Gamma'(\cdot,u)*_\Sigma R_i'(\partial_{x_j}f)\bigr)(x') -\phi(s)R_i'(\partial_{x_j}f)(x') \right| \leq C\min\{u^{\gamma/2},1\}.
\end{align}
The Fourier multiplier identity \(R_i'\partial_{x_j}(\Gamma'(\cdot,u)*_\Sigma f)=\Gamma'(\cdot,u)*_\Sigma R_i'(\partial_{x_j}f)\) and \(\int_0^\infty\partial_{x_n}\Gamma_1(x_n,u)\,du=-\frac12\) give
\begin{align*}
	R_i'\partial_{x_j}k_f(x,s)&+\frac12\phi(s)R_i'(\partial_{x_j}f)(x')\\
	 & =\int_0^\infty\partial_{x_n}\Gamma_1(x_n,u) \left[ \phi(s-u)\bigl(\Gamma'(\cdot,u)*_\Sigma R_i'(\partial_{x_j}f)\bigr)(x')-\phi(s)R_i'(\partial_{x_j}f)(x') \right]\,du.
\end{align*}
By \eqref{260254} and the scaling \(u=x_n^2v\),
\begin{align}
	\sup_{0\leq s\leq T} \left\|R_i'\partial_{x_j}k_f(\cdot,x_n,s) +\frac12\phi(s)R_i'(\partial_{x_j}f)\right\|_{L^\infty} & \leq Cx_n\int_0^\infty u^{-3/2}e^{-x_n^2/(4u)} \min\{u^{\gamma/2},1\}\,du \leq Cx_n^\gamma,\label{260255}\\
	\sup_{0\leq s\leq T} \left\|R_i'\partial_{x_j}\partial_{x_n}k_f(\cdot,x_n,s)\right\|_{L^\infty} & \leq C\int_0^\infty u^{-3/2}(1+x_n^2/u)e^{-x_n^2/(4u)} \min\{u^{\gamma/2},1\}\,du \leq Cx_n^{\gamma-1}.\label{260256}
\end{align}
The same estimate as in \eqref{260255}, with \(f\) in place of \(R_i'(\partial_{x_j}f)\), gives \(\sup_{0\leq s\leq T}\|k_f(\cdot,x_n,s)+f\phi(s)/2\|_{L^\infty}\leq Cx_n^\gamma\).

Moving the tangential derivative in \eqref{0419-02} onto \(f\) in \eqref{0330-01}, and using \(\partial_{x_j}k_f(\cdot,x_n-r,s)\in L^1\cap L^\infty\subset L^2(\Sigma)\), we obtain from \eqref{261031}
\begin{align}\label{260257}
	v_i(x,s) & =-2\delta_{ij}k_f(x,s)+4\int_0^{x_n}\int_\Sigma\partial_{y_n}N(y',r) R_i'\partial_{x_j}k_f(x'-y',x_n-r,s)\,dy'dr.
\end{align}
By Fubini's theorem and the fundamental theorem of calculus,
\begin{align}\label{260258}
	\partial_{x_n}v_i(x,s) & =-2\delta_{ij}\partial_{x_n}k_f(x,s) -2\phi(s)\int_\Sigma\partial_{y_n}N(y',x_n) R_i'(\partial_{x_j}f)(x'-y')\,dy'\nonumber\\*
	 & \quad+4\int_0^{x_n}\int_\Sigma\partial_{y_n}N(y',r) R_i'\partial_{x_j}\partial_{x_n} k_f(x'-y',x_n-r,s)\,dy'dr,
\end{align}
and scaling the Poisson kernel gives
\begin{align}\label{260259}
	\left\|2\partial_{y_n}N(\cdot,x_n)*_\Sigma R_i'(\partial_{x_j}f) -R_i'(\partial_{x_j}f)\right\|_{L^\infty} & \leq [R_i'(\partial_{x_j}f)]_{C^\gamma} \int_\Sigma2\partial_{y_n}N(y',x_n)|y'|^\gamma\,dy'\leq Cx_n^\gamma.
\end{align}

Since \(|z'-y'|\geq d\) for \(z'\in U\) and \(y'\in\supp f\), differentiating \eqref{0703-05} gives the continuous boundary value
\begin{align*}
	\partial_{x_n}k_f(z',0,s) =\int_0^s\int_\Sigma \partial_{x_n}^2\Gamma(z'-y',0,s-u)f(y')\phi(u)\,dy'du.
\end{align*}
Since \(\partial_{x_n}^3\Gamma(z'-y',0,u)=0\), Taylor's theorem and the heat kernel bound give
\begin{align}\label{260260}
	\sup_{(z',s)\in U\times J} \left|\partial_{x_n}k_f(z',x_n,s)-\partial_{x_n}k_f(z',0,s)\right| & \leq Cx_n^2\|f\|_{L^1}\|\phi\|_{L^\infty(0,1)} \int_0^T u^{-(n+4)/2}e^{-d^2/(Cu)}\,du \leq Cx_n^2.
\end{align}
The three terms in \eqref{260258} now give
\begin{align}\label{260261}
	\sup_{(z',s)\in U\times J} \left|\partial_{x_n}v_i(z',x_n,s) +2\delta_{ij}\partial_{x_n}k_f(z',0,s) +\phi(s)R_i'(\partial_{x_j}f)(z')\right| \leq Cx_n^\gamma.
\end{align}
The boundary value in \eqref{260261} is continuous in \((z',s)\).
With \(f=\psi\) and \(f=\partial_{x_k}\psi\), it gives \(\partial_{x_n}w_i\) and \(\partial_{x_n}\partial_{x_k}w_i\) on \(U\times\{0\}\times J\) for \(i,k<n\), and \eqref{260252} then gives \(\partial_{x_n}^2w_n(z',0,s)\).

\underline{Proof of (i).}\quad
Suppose that \(\partial_{x_n}^{\ell}w_i(x_*',0,t_*)\neq0\).
By continuity, there are \(c,r_0>0\) such that \(\partial_{x_n}^{\ell}w_i(z',0,s)\) has a fixed sign and absolute value at least \(c\) on \(\{|z'-x_*'|\leq r_0,\ |s-t_*|\leq r_0^2\}\subset(\Sigma\setminus\supp\psi)\times(0,\infty)\).
By \eqref{260263}, there is \(h_0>0\) such that \(w_i(z',x_n,s)\) has the same sign for \(|z'-x_*'|<r_0\), \(|s-t_*|<r_0^2\), and \(0<x_n<h_0\).
Hence one of \(S_{i,\varepsilon,M}^+\) and \(S_{i,\varepsilon,M}^-\) is empty for \(0<\varepsilon<\min\{\varepsilon_*,h_0\}\) and \(0<M<\min\{M_*,r_0\}\), where \(\varepsilon_*\) and \(M_*\) are the constants in Definition \ref{260000}.
This contradiction proves \eqref{1021-01}.

For the converse, choose \(M_*>0\) such that \(\{|z'-x_*'|\leq M_*,\ |s-t_*|\leq M_*^2\}\subset(\Sigma\setminus\supp\psi)\times(0,\infty)\), and fix \(0<M<M_*\).
By hypothesis and continuity, there are open sets \(V_\pm\) with closures in \(\{|z'-x_*'|<M,\ |s-t_*|<M^2\}\) and constants \(c_\pm>0\) such that \(\pm\partial_{x_n}^{\ell}w_i(z',0,s)\geq c_\pm\) on \(V_\pm\).
By \eqref{260263} with \(U=\overline{B_{M_*}'(x_*')}\) and \(J=[t_*-M_*^2,t_*+M_*^2]\), there is \(h_0\in(0,1]\) such that \(\pm w_i(z',x_n,s)\geq c_\pm x_n^{\ell}/(2\ell!)\) on \(V_\pm\) for \(0<x_n<h_0\).
Hence, for \(0<\varepsilon<1\), \(S_{i,\varepsilon,M}^\pm\) contains \(\{(z',x_n,s):(z',s)\in V_\pm,\ 0<x_n<\min\{\varepsilon,h_0\}\}\), whose measure is \(|V_\pm|\min\{\varepsilon,h_0\}>0\).
Thus \((x_*',t_*)\) is a separation point of \(w_i\).

\underline{Proof of (ii).}\quad
By hypothesis, \(\partial_{x_n}^{\ell}w_i(x',0,\cdot)\) takes both signs in every neighborhood of \(t_i^*(x')\), and \(t_i^*(x')\) is its only zero.
Hence {\rm(ii)} follows from {\rm(i)}.

\underline{Proof of (iii).}\quad
Fix \(x_*'\in U\), and choose an open ball \(B'\Subset U\) containing \(x_*'\), a compact interval \(J'\Subset J\) containing \(t\) in its interior, and \(0<h_1<h_0\).
By the boundary condition and \eqref{1021-04}, \(\partial_{x_n}^m w_i(z',0,s)=0\) for \(z'\in B'\), \(s\in J'\), and \(0\leq m<\ell\).
Taylor's theorem in \(x_n\) and the uniform continuity of \(\partial_{x_n}^{\ell}w_i\) on \(\overline{B'}\times[0,h_1]\times J'\) give
\begin{align}\label{1021-02}
	\sup_{(z',s)\in B'\times J'} \left|\frac{\ell!w_i(z',x_n,s)}{x_n^\ell} -\partial_{x_n}^{\ell}w_i(z',0,s)\right|\longrightarrow0 \qquad\text{as }x_n\to0^+.
\end{align}
The proof of {\rm(i)}, with \eqref{1021-02} in place of \eqref{260263}, shows that \((x_*',t)\) is a separation point if \(\partial_{x_n}^\ell w_i(\cdot,0,\cdot)\) takes both signs in every neighborhood of \((x_*',t)\), and is not one if \(\partial_{x_n}^\ell w_i(x_*',0,t)\neq0\).
By \eqref{1021-05} and the hypothesis on \(f\), the first case occurs if \(f(x_*')=0\), and the second if \(f(x_*')\neq0\).
\end{proof}

	\section{Separation points for normal boundary data}\label{260049}
In this section, we prove Theorems \ref{thm-0324-01} and \ref{thm-0507-01}.
Throughout this section the boundary data are \emph{normal}, that is, we take \(j=n\) in \eqref{0324-13}, so that \( g(x',t)=\psi(x')\phi(t)\mathbf e_n. \)

\subsection{Tangential velocity component}\label{260050}
\subsubsection{Formula for the normal derivatives of \(w_i\)}
\begin{proposition}\label{prop-0421-01}
		Let \(1\leq i\leq n-1\), and let \(w\) solve the Stokes system \eqref{0308-01} with boundary data \eqref{0324-13}, where \(j=n\).
Suppose that \(\psi\) satisfies {\rm\textbf{(S)}}, and that \(\phi\) satisfies {\rm\textbf{(T1)}}.
	Then, for every \(x'\in \Sigma\) with \(|x'|>1\) and \(t\geq1\),
	\begin{align}\label{0319-10}
		\begin{split}
			\partial_{x_n}w_i(x',0,t) & = -\frac{\Gamma\left( \frac{n}{2}\right)}{2^{n+1}\pi^{\frac{n}{2}}\Gamma\left(\frac{n+1}{2}\right)}\int_{\Sigma}(x_i-y_i)\psi(y')\\
			 & \quad \times\int_0^1\frac{\phi(s)}{(t-s)^{\frac{n+3}{2}}} {}_1F_1\left( \frac{n}{2}, \frac{n+1}{2}; -\frac{|x'-y'|^2}{4(t-s)}\right)\,dsdy'.
		\end{split}
	\end{align}
	For every \(x'\in\Sigma\) with \(|x'|>1\) and \(t>0\),
	\begin{align}\label{0322-04}
		\begin{split}
			\partial_{x_n}w_i(x',0,t) & =\frac{\Gamma\left(\frac{n}{2}\right)}{2^{n+1}\pi^{\frac{n}{2}}\Gamma\left(\frac{n+3}{2}\right)}\int_{\Sigma}(x_i-y_i)\psi(y')\\
			 & \quad \times\int_0^t \frac{\phi'(s)}{(t-s)^{\frac{n+1}{2}}}{}_1F_1\left( \frac{n}{2}, \frac{n+3}{2};- \frac{|x'-y'|^2}{4(t-s)}\right)\,dsdy'.
		\end{split}
	\end{align}
\end{proposition}
\begin{proof}
	We recall the following formula from \cite[Lemma 3.1]{CK2026}.
	\begin{align}\label{0324-10}
		\begin{split}
			\partial_{x_n}w_i(x',0,t) & =-\frac{1}{\sqrt{4\pi}}\int_{-\infty}^t \frac{\phi(s)}{(t-s)^{\frac32}}\int_{\Sigma}\Gamma'(x'-y',t-s)(R_i'\psi(y')-R_i'\psi(x'))\,dy'ds \\
			& \quad - \frac{1}{\sqrt{4\pi}} R_i'\psi(x') \int_{-\infty}^t \frac{\phi(s)-\phi(t)}{(t-s)^{\frac32}}\,ds,
		\end{split}
	\end{align}
	where we extend \(\phi\) by zero to \((-\infty,0)\).
	For \(s<t\), the equality
	\begin{align*}
		\Gamma'(\cdot,t-s)*_\Sigma R_i'\psi =R_i'\bigl(\Gamma'(\cdot,t-s)*_\Sigma\psi\bigr)
	\end{align*}
	holds.
	Together with \(\int_\Sigma\Gamma'(x'-y',t-s)\,dy'=1\), this gives
	\begin{align*}
		\int_{\Sigma}\Gamma'(x'-y',t-s) & (R_i'\psi(y')-R_i'\psi(x'))\,dy' = -R_i'\psi(x')+ 2\int_{\Sigma}\psi(y')\mathrm{p.v.}\int_{\Sigma}\partial_{z_i}N(z',0)\Gamma'(x'-y'-z',t-s)\,dz'dy' \\
		& = -R_i'\psi(x') - \frac{\Gamma\left(\frac{n}{2}\right)}{\Gamma\left( \frac{n+1}{2}\right)}\sqrt{t-s} \int_{\Sigma}\partial_{x_i}\Gamma'(x'-y',t-s){}_1F_1\left( \frac12, \frac{n+1}{2};\frac{|x'-y'|^2}{4(t-s)}\right)\psi(y')\,dy'\\
		 & =-\frac{\Gamma\left(\frac{n}{2}\right)}{\pi^{\frac{n}{2}}}\int_{\Sigma}\frac{x_i-y_i}{|x'-y'|^n}F\left(\frac{|x'-y'|^2}{4(t-s)}\right)\psi(y')\,dy'.
	\end{align*}
The inner integral in \(z'\) is understood in the principal-value sense of Proposition \ref{prop-0322-001}, which gives the second equality.
	The third follows from \eqref{0419-021} and Kummer's transformation \eqref{0327-04}, with
	\begin{align}\label{0322-03}
		F(u) & :=1-\frac{\sqrt{\pi}}{\Gamma\left(\frac{n+1}{2}\right)}u^{\frac{n}{2}}e^{-u}{}_1F_1\left(\frac{1}{2},\frac{n+1}{2}; u\right)\nonumber \\
		& =1-\frac{\sqrt{\pi}}{\Gamma\left(\frac{n+1}{2}\right)}u^{\frac{n}{2}}{}_1F_1\left(\frac{n}{2},\frac{n+1}{2};-u\right),\qquad u>0.
	\end{align}

	Substituting the last identity and \eqref{0419-021} into \eqref{0324-10} gives
	\begin{align}\label{0322-02}
		\partial_{x_n}w_i(x',0,t) & = \frac{\Gamma\left(\frac{n}{2}\right)}{2\pi^{\frac{n+1}{2}}}\int_{-\infty}^t \frac{\phi(s)}{(t-s)^{\frac32}}\int_{\Sigma}\frac{x_i-y_i}{|x'-y'|^n}F\left(\frac{|x'-y'|^2}{4(t-s)}\right)\psi(y')\,dy'ds\nonumber\\
		& \quad + \frac{\Gamma\left(\frac{n}{2}\right)}{2\pi^{\frac{n+1}{2}}}\int_{\Sigma}\frac{x_i-y_i}{|x'-y'|^n}\psi(y')\,dy'\int_{-\infty}^t \frac{\phi(t)-\phi(s)}{(t-s)^{\frac32}}\,ds\nonumber\\
		& = \frac{\Gamma\left(\frac{n}{2}\right)}{2\pi^{\frac{n+1}{2}}}\int_{\Sigma}\frac{x_i-y_i}{|x'-y'|^n}\psi(y')\int_{-\infty}^t \frac{1}{(t-s)^{\frac32}}\left( \phi(t)-\phi(s) + F\left(\frac{|x'-y'|^2}{4(t-s)}\right) \phi(s)\right)\,dsdy'.
	\end{align}

	First suppose that \(t\geq1\).
	Using \(\phi(t)=0\), and \(\phi(s)=0\) for \(s>1\), we find from \eqref{0322-02} that
	\begin{align}\label{0322-10}
		\partial_{x_n}w_i(x',0,t) & = - \frac{\Gamma\left(\frac{n}{2}\right)}{2\pi^{\frac{n+1}{2}}}\int_{\Sigma}\frac{x_i-y_i}{|x'-y'|^n}\psi(y')\int_{0}^1 \frac{1}{(t-s)^{\frac32}}\left( 1- F\left(\frac{|x'-y'|^2 }{4(t-s)}\right)\right)\phi(s)\,dsdy'.
	\end{align}
	By \eqref{0322-03},
	\begin{align*}
		\frac{x_i-y_i}{|x'-y'|^n}\left(1-F\left(\frac{|x'-y'|^2}{4(t-s)}\right)\right) =\frac{\sqrt\pi}{2^n\Gamma\left(\frac{n+1}{2}\right)}\frac{x_i-y_i}{(t-s)^{\frac n2}} {}_1F_1\left(\frac n2,\frac{n+1}{2};-\frac{|x'-y'|^2}{4(t-s)}\right).
	\end{align*}
Substituting this into \eqref{0322-10} gives \eqref{0319-10} for every \(t\geq1\); at \(t=1\) the integral converges absolutely by {\rm\textbf{(T1)}}.

	We now derive \eqref{0322-04} for arbitrary \(t>0\).
	We separate the time integral on the right-hand side of \eqref{0322-02} as follows.
	\begin{align*}
		\int_{-\infty}^t \frac{1}{(t-s)^{\frac32}} & \left( \phi(t)-\phi(s) + F\left(\frac{|x'-y'|^2}{4(t-s)}\right) \phi(s)\right)\,ds\\
		& = \int_{-\infty}^t \frac{\phi(t)-\phi(s)}{(t-s)^{\frac32}}\,ds + \int_{-\infty}^t\frac{1}{(t-s)^{\frac32}}F\left( \frac{|x'-y'|^2}{4(t-s)}\right) \phi(s)\,ds=: I_1(t) + I_2(x'-y',t).
	\end{align*}
	By \eqref{0511-02}, integration by parts gives
	\begin{align}\label{0322-07}
		I_1(t) & =\int_{0}^\infty\frac{\phi(t)-\phi(t-s)}{s^{\frac32}}\,ds= \lim_{s\rightarrow 0^+}\frac{2(\phi(t) - \phi(t-s))}{s^{\frac12}} + 2\int_0^\infty \frac{\phi'(t-s)}{s^{\frac12}}\,ds = 2\int_0^t \frac{\phi'(s)}{\sqrt{t-s}}\,ds.
	\end{align}

	For \(I_2\), we first define
	\begin{align}\label{0322-05}
		\begin{split}
			V(x'-y',t) & :=\int_0^t \frac{1}{s^{\frac32}}F\left(\frac{|x'-y'|^2}{4s}\right)\,ds\\
		& = \frac{2}{\sqrt{t}}\left[\frac{\sqrt{\pi}}{2\Gamma\left(\frac{n+3}{2}\right)}\left(\frac{|x'-y'|^2}{4t}\right)^{\frac{n}{2}}{}_1F_1\left(\frac{n}{2}, \frac{n+3}{2}; -\frac{|x'-y'|^2}{4t}\right)-1\right],
		\end{split}
	\end{align}
	where \eqref{0423-10} was used to compute the integral.
	Since \(\phi(s)=0\) for \(s\leq0\), integration by parts gives
	\begin{align}\label{0322-06}
		\begin{split}
			I_2(x'-y',t) & =\int_{-\infty}^t \partial_t V(x'-y',t-s)\phi(s)\,ds
			=\int_0^t \partial_s V(x'-y',s)\phi(t-s)\,ds\\
			& = \int_0^t V(x'-y',s)\phi'(t-s)\,ds.
		\end{split}
	\end{align}

	Adding \eqref{0322-07} and \eqref{0322-06}, and then using \eqref{0322-05}, gives
	\begin{align}\label{0322-08}
		\begin{split}
			\int_{-\infty}^t \frac{1}{(t-s)^{\frac32}} & \left( \phi(t)-\phi(s) + F\left(\frac{|x'-y'|^2}{4(t-s)}\right) \phi(s)\right)\,ds\\
			& = \frac{\sqrt{\pi}}{\Gamma\left(\frac{n+3}{2}\right)}\int_0^t \frac{\phi'(s)}{\sqrt{t-s}}\left( \frac{|x'-y'|^2}{4(t-s)}\right)^{\frac{n}{2}}{}_1F_1\left( \frac{n}{2}, \frac{n+3}{2};- \frac{|x'-y'|^2}{4(t-s)}\right)\,ds.
		\end{split}
	\end{align}
	Then \eqref{0322-04} follows from \eqref{0322-02} and \eqref{0322-08}.
\end{proof}

\subsubsection{Nonexistence of separation points}\label{subsect-0619-01}
The next proposition excludes separation points outside \(\{|x_i|\leq1\}\) while \(\phi\) is increasing and after the boundary data vanish.
This proves the nonexistence statement in Theorem \ref{thm-0324-01}{\rm(i)}.
\begin{proposition}\label{prop-0730-01}
	Let \(1\leq i \leq n-1\), and let \(w\) solve the Stokes system \eqref{0308-01} with boundary data \eqref{0324-13}, where \(j=n\).
	Suppose that \(\psi\) satisfies {\rm\textbf{(S)}}, and that \(\phi\) satisfies {\rm\textbf{(T1)}}.
	Then, for every \(x'\in \Sigma\) with \(|x_i|>1\), the following hold.
	\begin{enumerate}
		\item If \(0<t<1\) and \(\phi\) is nondecreasing on \([0,t]\), then \(\sgn(x_i)\partial_{x_n}w_i(x',0,t)\geq0\).
		Moreover, if \(\phi(t)>0\), then \(\sgn(x_i)\partial_{x_n}w_i(x',0,t)>0\).
    \item If \(t\geq1\), then \(\sgn(x_i)\partial_{x_n}w_i(x',0,t)\leq0\).
    Moreover, if \(\phi\not\equiv0\), then \(\sgn(x_i)\partial_{x_n}w_i(x',0,t)<0\).
	\end{enumerate}
\end{proposition}
\begin{proof}
	The integral representation \eqref{0324-14} shows that both confluent hypergeometric functions in \eqref{0319-10} and \eqref{0322-04} are positive, and \(\sgn(x_i)(x_i-y_i)>0\) on \(\supp\psi\).
	
	\underline{Proof of (i).}\quad
	Since \(\phi'(s)\geq0\) for \(0<s<t\), \eqref{0322-04} gives \(\sgn(x_i)\partial_{x_n}w_i(x',0,t)\geq0\), and the integral is nonzero when \(\phi(t)=\int_0^t\phi'(s)\,ds>0\).

	\underline{Proof of (ii).}\quad
	Since \(\phi\geq0\), \eqref{0319-10} gives \(\sgn(x_i)\partial_{x_n}w_i(x',0,t)\leq0\).
	The integral is nonzero when \(\phi\not\equiv0\).
	Proposition \ref{prop-260046}{\rm(i)} then excludes a separation point, including at \(t=1\).
\end{proof}

\subsubsection{Existence of separation points}\label{subsect-0619-001}
To study separation points of \(w_i\), we introduce \(G(x',t)\) for \(x'\notin\supp\psi\) and \(t>0\) as follows.
\begin{align}\label{0320-27}
	G(x',t):= \int_0^t \frac{\phi'(s)}{\sqrt{t-s}}\int_{\Sigma} \frac{x_i-y_i}{|x'-y'|^n}\psi(y') Q\left(\frac{|x'-y'|^2}{4(t-s)} \right)\,dy'ds,
\end{align}
where \(Q\) is given in \eqref{0827-10}.
Then \eqref{0322-04} gives, for every \(|x'|>1\) and \(t>0\),
\begin{align}\label{0630-07}
	\partial_{x_n}w_i(x',0,t) = \frac{\Gamma\left(\frac{n}{2}\right)}{\pi^{\frac{n+1}{2}}}G(x',t).
\end{align}

We also introduce the kernel
\begin{align}\label{0325-06}
	K(x',t) := \frac{1}{\sqrt{t}}\int_\Sigma \frac{x_i-y_i}{|x'-y'|^n}\psi(y')Q\left( \frac{|x'-y'|^2}{4t}\right)\,dy',
\end{align}
so that by \eqref{0320-27},
\begin{align}\label{0324-15}
	G(x',t) = \int_0^t \phi'(s)K(x',t-s)\,ds.
\end{align}

The next proposition gives a separation point of \(w_i\) at some time in \((t_1,1)\) for each \(x'\) with \(|x_i|>1\).
This proves the existence statement in Theorem \ref{thm-0324-01}{\rm(i)}.
\begin{proposition}\label{prop-0504-03}
	Let \(1\leq i\leq n-1\), and let \(w\) solve the Stokes system \eqref{0308-01} with boundary data \eqref{0324-13}, where \(j=n\).
		Suppose that \(\psi\) satisfies {\rm\textbf{(S)}}, and that \(\phi\) satisfies {\rm\textbf{(T1)}} and {\rm\textbf{(T2)}}.
	Then, for each \(x'\in\Sigma\) with \(|x_i|>1\), there exists \(t_i^*(x')\in(t_1,1)\) such that \((x',t_i^*(x'))\) is a separation point of \(w_i\).
\end{proposition}
\begin{lemma}\label{lem-0424-02}
	Let \(1\leq i\leq n-1\), suppose that \(\psi\) satisfies {\rm\textbf{(S)}}, and let \(K\) be given in \eqref{0325-06}.
	Then, for \(x'\in \Sigma\) with \(|x_i|>1\), the following hold.
	\begin{enumerate}
		\item\label{item-0828-01} \(\sgn(x_i)K(x',t)\) is positive and strictly decreasing in \(t\in(0,\infty)\).
		\item\label{item-0828-02} \(K\) satisfies
		\begin{align}\label{0424-10}
			|K(x',t)|\leq C\min\left\{\frac{|x_i|}{\sqrt{t}(|x'|-1)^n}, \frac{|x_i|}{t^{\frac{n+1}{2}}}\right\} M_0.
		\end{align}
		\item \(\partial_t K\) satisfies
		\begin{align}\label{0424-12}
			|\partial_tK(x',t)|\leq C\min\left\{\frac{|x_i|}{t^{\frac32}(|x'|-1)^n}, \frac{|x_i|}{t^{\frac{n+3}{2}}}\right\}M_0.
		\end{align}
	\end{enumerate}
	Here \(C>0\) depends only on \(n\).
\end{lemma}
\begin{proof}
	\underline{Proof of (i).}\quad
	By Lemma \ref{lem-0424-01}{\rm(i)}, \(Q(u)>0\) and \(Q'(u)\geq 0\) for \(u\geq 0\).

	Differentiating \eqref{0325-06} in \(t\) gives
	\begin{align}\label{0914-01}
		\partial_tK(x',t)=-\frac1{t^{3/2}}\int_\Sigma \frac{x_i-y_i}{|x'-y'|^n}\psi(y')\left(\frac12 Q\left( \frac{|x'-y'|^2}{4t}\right)+ \frac{|x'-y'|^2}{4t}Q'\left( \frac{|x'-y'|^2}{4t}\right)\right)\,dy'.
	\end{align}
	Since \(\sgn(x_i)(x_i-y_i)>0\) on \(\supp\psi\), \eqref{0325-06} and \eqref{0914-01} show that \(\sgn(x_i)K(x',t)\) is positive and strictly decreasing in \(t\).

	\underline{Proof of (ii).}\quad
	For \eqref{0424-10}, \eqref{0322-112} of Lemma \ref{lem-0424-01} gives
	\begin{align*}
		|K(x',t)| \leq\min\left\{ \frac{1}{\sqrt{t}}\int_{\Sigma}\frac{|x_i-y_i|}{|x'-y'|^n}\psi(y')\,dy', \frac{C}{t^{\frac{n+1}{2}}}\int_{\Sigma}|x_i-y_i|\psi(y')\,dy'\right\},
	\end{align*}
	and \(|x_i-y_i|\leq2|x_i|\) and \(|x'-y'|\geq|x'|-1\) prove \eqref{0424-10}.

	\underline{Proof of (iii).}\quad
For \eqref{0424-12}, \eqref{0827-10}--\eqref{0827-11} give
	\begin{align*}
		\frac12Q(u)+uQ'(u) =\frac{\sqrt\pi}{2\Gamma\left(\frac{n+3}{2}\right)}u^{\frac n2} \left(\frac{n+1}{2}P(u)+uP'(u)\right),\quad u>0.
	\end{align*}
	Then \eqref{0327-07} and \eqref{0514-02} give
	\begin{align*}
		\left|\frac12Q(u)+uQ'(u)\right| \leq Cu^{\frac n2}(1+u)^{-\frac n2} \leq C\min\{u^{\frac n2},1\}.
	\end{align*}
		Substituting this estimate into \eqref{0914-01} and using the estimates from \eqref{0424-10} gives \eqref{0424-12}.
\end{proof}

\begin{proof}[Proof of Proposition \ref{prop-0504-03}]
Relabeling or reflecting the tangential coordinates transforms \(\psi\) into a function that still satisfies {\rm\textbf{(S)}}.
		It therefore suffices to consider \(i=1\) and \(x_1>1\).

	Since \(\phi'>0\) on \((0,t_1)\) by {\rm\textbf{(T2)}} and \(K(x',\cdot)>0\) by Lemma \ref{lem-0424-02}{\rm(i)}, \eqref{0324-15} gives \(G(x',t)>0\) for every \(0<t\leq t_1\).
At the other endpoint, Proposition \ref{prop-0730-01}{\rm(ii)} and \eqref{0630-07} give \(G(x',1)<0\).

	By \eqref{0322-04}, \eqref{0630-07}, and Lemma \ref{lem-260044}{\rm(iv)}, applied with \(g=\psi\), \(q=\phi'\), \(m=1\), \(a=\frac n2\), \(b=\frac{n+3}{2}\), and \(c=\frac{n+1}{2}\), \(G\) is jointly continuous for \(x_1>1\), \(t>0\), and real-analytic in \(x_1>1\) for fixed \(x''\) and \(t\).

For fixed \(x''\) and \(t>0\), Lemma \ref{lem-0424-01}{\rm(iii)} and \eqref{0323-02} give, as \(r\to\infty\),
	\begin{align}\label{1006-03}
		G(r,x'',t) & =T_0(t)\int_\Sigma\frac{r-y_1}{|(r,x'')-y'|^n}\psi(y')\,dy' -nT_1(t)\int_\Sigma\frac{r-y_1}{|(r,x'')-y'|^{n+2}}\psi(y')\,dy' +O(r^{-n-3}).
	\end{align}
	Integration by parts gives \(T_1(t)=\frac12\int_0^t\phi(s)(t-s)^{-\frac12}\,ds>0\).
	Consequently, \eqref{1006-03} yields
	\begin{align*}
		\lim_{r\to\infty}r^{n-1}G(r,x'',t)=M_0T_0(t)\neq0 \quad\text{if }T_0(t)\neq0, \qquad \lim_{r\to\infty}r^{n+1}G(r,x'',t)=-nM_0T_1(t)<0 \quad\text{if }T_0(t)=0.
	\end{align*}
	Hence, by Lemma \ref{lem-260044}{\rm(iii)}, \(G(r,x'',t)\) cannot vanish on any open interval of \(r>1\).

	We now prove that there exists \(t_*\in(t_1,1)\) such that every neighborhood of \((x_1,t_*)\) in the \((r,s)\)-plane contains points with \(G(r,x'',s)>0\) and points with \(G(r,x'',s)<0\).
	Otherwise, the relatively open sets
	\begin{align*}
		A_\pm:=\{t\in[t_1,1]:\ \pm G(r,x'',s)\geq0 \text{ in some neighborhood of }(x_1,t)\}
	\end{align*}
	cover \([t_1,1]\).
	They are disjoint, since \(t\in A_{+}\cap A_{-}\) would imply \(G(r,x'',t)=0\) on an open interval containing \(x_1\).
	The joint continuity of \(G\) and \(G(x',t_1)>0>G(x',1)\) then give \(t_1\in A_+\) and \(1\in A_-\), contradicting the connectedness of \([t_1,1]\).
	The continuity of \(G\) gives \(G(x',t_*)=0\), and \eqref{0630-07} yields \(\partial_{x_n}w_1(x',0,t_*)=0\).
	Proposition \ref{prop-260046}{\rm(i)}, with \eqref{0630-07} and the two signs of \(G\) just proved, shows that \((x',t_*)\) is a separation point.
	Take \(t_1^*(x')=t_*\).
\end{proof}

\subsubsection{Temporal vanishing rates of \(\partial_{x_n}w_i\) near a zero}
The next proposition gives the temporal vanishing rates of \(\partial_{x_n}w_i\) near any such zero and the leading coefficients when \(0<\beta<\frac12\).
This proves Theorem \ref{thm-0324-01}{\rm(ii)}.

\begin{proposition}\label{prop-0512-02}
	Let \(1\leq i\leq n-1\), and let \(w\) solve the Stokes system \eqref{0308-01} with boundary data \eqref{0324-13}, where \(j=n\).
		Suppose that \(\psi\) satisfies {\rm\textbf{(S)}}, and that \(\phi\) satisfies {\rm\textbf{(T1)}} and {\rm\textbf{(T2)}}.
		Let \(\beta\in(0,1)\) be the H\"older exponent of \(\phi'\) in {\rm\textbf{(T1)}}.
		Let \(x'\in\Sigma\) satisfy \(|x_i|>1\), and let \(t_i^*=t_i^*(x')\in(t_1,1)\) be a zero of \(\partial_{x_n}w_i(x',0,\cdot)\).
	Then the following hold.
	\begin{enumerate}
		\item\label{item-0801-01} If \(h=|t-t_i^*|\) satisfies \(0<h<\min\{t_i^*,1-t_i^*\}\), then
		\begin{align}\label{0919-01}
			|\partial_{x_n}w_i(x',0,t)| \leq C\frac{|x'|^2}{(|x'|-1)^{n+1}}
			\begin{cases}
				h^{\beta+\frac12}, & 0<\beta<\frac12, \\
				h|\ln h|,& \beta=\frac12,\\
				h^\sigma,& \frac12<\beta<1
			\end{cases},
		\end{align}
		where, for \(\beta>1/2\), \(\sigma=1\) if \(\partial_t\partial_{x_n}w_i(x',0,t_i^*)\neq0\) and \(\sigma=\beta+\frac12\) otherwise, and \(C>0\) depends only on \(n\), \(\beta\), \(M_0\), \(t_1\), and \(\norm{\phi'}_{C^\beta}\).

		\item Suppose that \(0<\beta<\frac12\) and that, for some \(c_0>0\),
		\begin{align}\label{0730-70}
			\phi'(s)-\phi'(t_i^*)=-c_0(s-t_i^*)|s-t_i^*|^{\beta-1}+o\left(|s-t_i^*|^{\beta}\right), \qquad s\rightarrow t_i^*.
		\end{align}
		Then
		\begin{align}\label{0730-71}
			\lim_{\ep\rightarrow0^+}\ep^{-\beta-\frac12}\,\partial_{x_n}w_i(x',0,t_i^*+\ep) =-\frac{\Gamma\left(\frac n2\right)}{\pi^{\frac{n+1}{2}}}\,c_0 \frac{\Gamma(\beta+1)\Gamma(\frac12)}{\Gamma(\beta+\frac32)} \bigl(1+\sec(\pi\beta)\bigr) \int_{\Sigma}\frac{x_i-y_i}{|x'-y'|^{n}}\psi(y')\,dy',
		\end{align}
		and
		\begin{align}\label{0730-71a}
			\lim_{\ep\rightarrow0^-}|\ep|^{-\beta-\frac12}\,\partial_{x_n}w_i(x',0,t_i^*+\ep) =\frac{\Gamma\left(\frac n2\right)}{\pi^{\frac{n+1}{2}}}\,c_0 \frac{\Gamma(\beta+1)\Gamma(\frac12)}{\Gamma(\beta+\frac32)}\tan(\pi\beta) \int_{\Sigma}\frac{x_i-y_i}{|x'-y'|^{n}}\psi(y')\,dy'.
		\end{align}
		Thus the exponent \(\beta+\frac12\) in {\rm(i)} is sharp, and \(\partial_{x_n}w_i(x',0,t)\) is not differentiable at \(t_i^*\).
	\end{enumerate}
\end{proposition}

\begin{proof}
	\underline{Proof of (i).}\quad
	By \eqref{0320-27},
	\begin{align}\label{0629-01}
		G(x',t) = \int_0^t \frac{\phi'(s)}{\sqrt{t-s}}W(x',t-s)\,ds,\qquad W(x',t) := \int_{\Sigma}\frac{x_i-y_i}{|x'-y'|^n}\psi(y') Q\left( \frac{|x'-y'|^2}{4t}\right)\,dy'.
	\end{align}
	We may assume \(x_i>1\); for \(x_i<-1\), replace \(G\) and \(W\) by \(-G\) and \(-W\).
	By \eqref{0827-10}, \eqref{0327-01}, and \eqref{0514-02}, with \(u=|x'-y'|^2/(4t)\),
	\begin{align*}
		\partial_t^kW(x',t)=(-1)^k\frac{\sqrt\pi\left(\frac n2\right)_k} {2^{n+1}\Gamma\left(\frac{n+3}{2}\right)t^{\frac n2+k}} \int_\Sigma(x_i-y_i)\psi(y') {}_1F_1\left(\frac n2+k,\frac{n+3}{2};-u\right)\,dy',\qquad 0\leq k\leq3.
	\end{align*}
	Since \(|x_i-y_i|\leq|x'-y'|\) and \(|x'-y'|\geq|x'|-1\), \eqref{0514-02}, Lemma \ref{lem-0424-01}{\rm(i)}, and dominated convergence give \(W\in C^2([0,1])\) and
	\begin{align}
		 & W(x',0)=\lim_{t\rightarrow 0}W(x',t) = \int_{\Sigma}\frac{x_i-y_i}{|x'-y'|^n}\psi(y')\,dy',\label{0428-11}\\
		 & \norm{\partial_t^k W(x',\cdot)}_{L_t^\infty((0,1))}\leq \frac{CM_0}{(|x'|-1)^{n-1+2k}}, \quad 0\leq k\leq 3.\label{0428-12}
	\end{align}

	Let \(\varepsilon=t-t_i^*\).
	By \eqref{0629-01} and \eqref{0428-12}, \(|G(x',t)|\leq CM_0\norm{\phi'}_{L^\infty(0,1)}(|x'|-1)^{1-n}\) for \(0<t<1\), which gives \eqref{0919-01} for \(|\varepsilon|\geq t_1/2\), since \(t_i^*>t_1\), \(|\varepsilon|<1/2\), and \((|x'|-1)^{1-n}\leq |x'|^2(|x'|-1)^{-n-1}\).
	Let \(0<\varepsilon<t_1/2\); negative \(\varepsilon\) is analogous.
	Then
	\begin{align}\label{0920-10}
		G(x',t_i^*+\varepsilon)-G(x', t_i^*) & = \int_0^{t_i^*}\phi'(s)\frac{W(x',t_i^*+\varepsilon-s)-W(x',t_i^*-s)}{\sqrt{t_i^*+\varepsilon-s}}\,ds\nonumber\\
		 & \quad + \int_0^{t_i^*}\phi'(s)W(x',t_i^*-s)\left( \frac{1}{\sqrt{t_i^*+\varepsilon-s}}-\frac{1}{\sqrt{t_i^*-s}}\right)\,ds\nonumber\\
		 & \quad + \int_{t_i^*}^{t_i^*+\varepsilon} \frac{\phi'(s)W(x',t_i^*+\varepsilon-s)}{\sqrt{t_i^*+\varepsilon-s}}\,ds=:A+B+C,
	\end{align}
	and \eqref{0428-12} gives \(|A|\leq C\varepsilon M_0\norm{\phi'}_{L^\infty(0,1)}(|x'|-1)^{-n-1}\).
	Subtracting and adding \(\phi'(t_i^*)W(x',0)\) in \(B\) and \(C\) gives
	\begin{align}\label{0920-20}
		B+C & =\int_0^{t_i^*}(\phi'(s)W(x',t_i^*-s)-\phi'(t_i^*)W(x',0)) \left(\frac{1}{\sqrt{t_i^*+\varepsilon-s}}-\frac{1}{\sqrt{t_i^*-s}}\right)\,ds \nonumber \\
		 & \quad+\int_{t_i^*}^{t_i^*+\varepsilon}\frac{\phi'(s)W(x',t_i^*+\varepsilon-s) - \phi'(t_i^*)W(x',0) }{\sqrt{t_i^*+\varepsilon-s}}\,ds \nonumber\\
		 & \quad + 2\phi'(t_i^*)W(x',0)(\sqrt{t_i^*+\varepsilon}-\sqrt{t_i^*})=:I_1+I_2+I_3.
	\end{align}
	By \eqref{0428-12}, \(|\phi'(s)W(x',\tau)-\phi'(t_i^*)W(x',0)|\leq CM_0\norm{\phi'}_{C^\beta}|x'|^2(|x'|-1)^{-n-1}(|s-t_i^*|^\beta+\tau)\) for \(0<\tau<1\).
	In \(I_1\), \(\tau=t_i^*-s=:r\), and \(r^{-1/2}-(r+\varepsilon)^{-1/2}\leq\min\{r^{-1/2},\varepsilon r^{-3/2}/2\}\) gives
	\begin{align*}
		|I_1|\leq CM_0\norm{\phi'}_{C^\beta }\frac{|x'|^2}{(|x'|-1)^{n+1}}\left(\varepsilon^{\beta+\frac12}\mathbf{1}_{\beta< \frac12} + \ep |\ln \ep|\mathbf{1}_{\beta=\frac12} + \ep \mathbf{1}_{\beta>\frac12}\right).
	\end{align*}
	In \(I_2\), \(|s-t_i^*|\) and \(\tau\) are at most \(\varepsilon\), so, with \eqref{0428-12},
	\begin{align*}
		|I_2|\leq CM_0\norm{\phi'}_{C^\beta}\frac{|x'|^2\varepsilon^{\beta+\frac12}}{(|x'|-1)^{n+1}},\qquad |I_3|\leq CM_0\norm{\phi'}_{L^\infty(0,1)}\frac{\varepsilon}{(|x'|-1)^{n-1}}.
	\end{align*}
	Collecting these bounds, since \(G(x',t_i^*)=0\), \eqref{0630-07} gives \eqref{0919-01} for \(\beta\leq1/2\), and for \(\beta>1/2\) with \(\sigma=1\).

	Let \(\beta>1/2\) and \(\partial_tG(x',t_i^*)=0\).
	Lemma \ref{lem-0424-20}, applied with \(K(t)=W(x',t)/\sqrt{t}\) and \(f=\phi'\), gives \(G(x',\cdot)\in C^{1,\beta-\frac12}([t_1/2,1))\), since \eqref{0428-12} implies \eqref{0829-01} with \(\sigma=\frac12\), and its proof bounds the H\"older seminorm of \(\partial_tG(x',\cdot)\) there by \(CM_0\norm{\phi'}_{C^\beta}|x'|^2(|x'|-1)^{-n-1}\).
	Hence
	\begin{align*}
		|G(x',t)|=|G(x',t)-G(x',t_i^*)-\partial_tG(x',t_i^*)(t-t_i^*)| \leq CM_0\norm{\phi'}_{C^\beta}\frac{|x'|^2}{(|x'|-1)^{n+1}}h^{\beta+\frac12},
	\end{align*}
	and \eqref{0630-07} proves \eqref{0919-01}.

	\underline{Proof of (ii).}\quad
	By \eqref{0428-12}, \(u\mapsto(W(x',u)-W(x',0))/\sqrt u\) is bounded by \(C\sqrt u\) and has derivative bounded by \(C/\sqrt u\) on \((0,1)\).
	Hence replacing \(W(x',t-s)\) by \(W(x',0)\) in \eqref{0629-01} changes \(G(x',t+h)-G(x',t)\) by \(O(h)\) for \(0<t<t+h<1\).
	Since \(G(x',t_i^*)=0\) and \(2\phi'(t_i^*)W(x',0)(\sqrt{t_i^*+\ep}-\sqrt{t_i^*})=O(|\ep|)\), this gives
	\begin{align}\label{0808-22}
		G(x',t_i^*+\ep) =W(x',0)\left[ \int_0^{t_i^*+\ep}\frac{\phi'(s)-\phi'(t_i^*)}{\sqrt{t_i^*+\ep-s}}\,ds -\int_0^{t_i^*}\frac{\phi'(s)-\phi'(t_i^*)}{\sqrt{t_i^*-s}}\,ds \right]+O(|\ep|).
	\end{align}
	Given \(\eta>0\), \eqref{0730-70} provides \(0<d<\min\{t_i^*,1-t_i^*\}\) such that
	\begin{align}\label{260051}
		\left|\phi'(t_i^*+u)-\phi'(t_i^*)+c_0u|u|^{\beta-1}\right| \leq\eta|u|^\beta,\qquad 0<|u|<d.
	\end{align}

	Let \(0<\ep<d\).
	Splitting the first integral in \eqref{0808-22} at \(t_i^*\), substituting \(r=t_i^*-s\) and \(u=s-t_i^*\), and applying \eqref{260051} on \((0,d)\) give
	\begin{align*}
		 & \left|G(x',t_i^*+\ep)+c_0W(x',0) \left[\int_0^{t_i^*}r^\beta \left(\frac1{\sqrt r}-\frac1{\sqrt{r+\ep}}\right)\,dr +\int_0^\ep\frac{u^\beta}{\sqrt{\ep-u}}\,du\right]\right| \leq C(\beta)|W(x',0)|\eta\ep^{\beta+\frac12}+C_d\ep.
	\end{align*}
	Extending the first integral to \((0,\infty)\) changes it by \(O(\ep)\), and substituting \(r=\ep v\) and using the identity \(v^{-1/2}-(1+v)^{-1/2}=\frac12\int_0^1(v+s)^{-3/2}\,ds\), beta integrals \eqref{0324-14} at \(z=0\), and the recurrence and reflection formulas \cite[(5.5.1), (5.5.3)]{DLMF} gives
	\begin{align*}
		\int_0^{t_i^*}r^\beta \left(\frac1{\sqrt r}-\frac1{\sqrt{r+\ep}}\right)\,dr =\frac{\Gamma(\beta+1)\Gamma(\frac12)} {\Gamma(\beta+\frac32)}\sec(\pi\beta)\,\ep^{\beta+\frac12}+O(\ep),\qquad
		\int_0^\ep\frac{u^\beta}{\sqrt{\ep-u}}\,du =\frac{\Gamma(\beta+1)\Gamma(\frac12)} {\Gamma(\beta+\frac32)}\ep^{\beta+\frac12}.
	\end{align*}
	Dividing by \(\ep^{\beta+\frac12}\) and letting \(\ep\to0^+\) and then \(\eta\to0^+\), we obtain
	\begin{align*}
		\lim_{\ep\to0^+}\ep^{-\beta-\frac12}G(x',t_i^*+\ep) =-c_0W(x',0)\frac{\Gamma(\beta+1)\Gamma(\frac12)} {\Gamma(\beta+\frac32)}\bigl(1+\sec(\pi\beta)\bigr),
	\end{align*}
	and \eqref{0630-07} and \eqref{0428-11} prove \eqref{0730-71}.

	Let \(0<h<d/2\).
	Substituting \(r=t_i^*-s\) in \eqref{0808-22} and applying \eqref{260051} give
	\begin{align*}
		\left|G(x',t_i^*-h)-c_0W(x',0) \left[\int_h^{t_i^*}\frac{r^\beta}{\sqrt{r-h}}\,dr -\int_0^{t_i^*}r^{\beta-\frac12}\,dr\right]\right| \leq C(\beta)|W(x',0)|\eta h^{\beta+\frac12}+C_dh.
	\end{align*}
	Substituting \(r=u+h\) and \(u=hv\) and using the identity \((1+v)^\beta-v^\beta=\beta\int_0^1(v+s)^{\beta-1}\,ds\) and beta integrals gives
	\begin{align*}
		\int_h^{t_i^*}\frac{r^\beta}{\sqrt{r-h}}\,dr -\int_0^{t_i^*}r^{\beta-\frac12}\,dr =\frac{\Gamma(\beta+1)\Gamma(\frac12)} {\Gamma(\beta+\frac32)}\tan(\pi\beta)\,h^{\beta+\frac12}+O(h),
	\end{align*}
	and, as above, \eqref{0630-07} and \eqref{0428-11} give \eqref{0730-71a}.

	Finally, \(x_i-y_i\) has the sign of \(x_i\) on \(\supp\psi\), so \(W(x',0)\neq0\) by \eqref{0428-11}.
	This proves the sharpness of the exponent, and, since \(\beta<\frac12\), \(\ep^{-1}G(x',t_i^*+\ep)\) is unbounded as \(\ep\to0^+\).
\end{proof}
\begin{remark}\label{rmk-critical-abel}
Suppose that \eqref{0730-70} holds with \(\beta=1/2\), and write \(t_*=t_i^*(x')\).
For \(h>0\), replace \(\phi'(s)-\phi'(t_*)\) in \eqref{0808-22} by its leading term in \eqref{0730-70} and rescale \(s-t_*\) by \(h\). For \(\ep=h\) and \(\ep=-h\), this gives the integrands \(1-\sqrt{v/(v+1)}\) and \(\sqrt{(v+1)/v}-1\), up to constant factors. They are integrable near zero and equal \(1/(2v)+O(v^{-2})\) as \(v\to\infty\), so their integrals up to \(t_*/h\) equal \(\tfrac12\ln(1/h)+O(1)\).
The parts over \((t_*,t_*+h)\) and \((t_*-h,t_*)\) are \(O(h)\), and, splitting the error in \eqref{260051} at \(d\) as in the proof, its contribution is at most \(C|W(x',0)|\eta h(1+\ln(d/h))+O_d(h)\) for \(0<h<d/2\).
Dividing by \(h\ln(1/h)\) and letting \(h\downarrow0\) and then \(\eta\downarrow0\) give, as \(\ep\to0\) from either side,
\begin{align*}
	G(x',t_*+\ep)=-\tfrac12 c_0 W(x',0)\,\ep\ln(1/|\ep|)
	+o\!\left(|\ep|\ln(1/|\ep|)\right).
\end{align*}
By \eqref{0630-07} and \eqref{0428-11}, the logarithmic rate of \(\partial_{x_n}w_i\) is attained on both sides, and \(\partial_t\partial_{x_n}w_i(x',0,t_*)\) does not exist.
\end{remark}

\subsubsection{Uniqueness of the zero}\label{subsect-0829-01}
The next proposition establishes uniqueness of the zero of \(\partial_{x_n}w_i(x',0,\cdot)\), shows that the corresponding point is a separation point, and determines the sign before and after the zero.
\begin{proposition}\label{prop-0324-02}
	Let \(1\leq i\leq n-1\), and let \(w\) solve the Stokes system \eqref{0308-01} with boundary data \eqref{0324-13}, where \(j=n\).
	Suppose that \(\psi\) satisfies {\rm\textbf{(S)}}, and that \(\phi\) satisfies {\rm\textbf{(T1)}} and {\rm\textbf{(T2)}} and is log-concave on \((0,1)\).
	Then, for each \(x'\in \Sigma\) with \(|x_i|>1\), \(\partial_{x_n}w_i(x',0,\cdot)\) has exactly one zero \(t_i^*(x')\), and \((x',t_i^*(x'))\) is a separation point of \(w_i\).
	Moreover,
	\begin{align}\label{0929-30}
		\sgn\bigl(\partial_{x_n}w_i(x',0,t)\bigr) =-\sgn(x_i)\sgn(t-t_i^*(x')).
	\end{align}
\end{proposition}
\begin{proof}
	Fix \(x'\in \Sigma\) with \(|x_i|>1\).
	We apply Lemma \ref{lem-0325-01} to \(\phi\), with the kernel \(\sgn(x_i)K(x',\cdot)\) from \eqref{0325-06}.
	By \eqref{0914-01}, \(\sgn(x_i)\partial_tK(x',t)<0\) for every \(t>0\), and \eqref{0424-10} and \eqref{0424-12} give \eqref{0829-01} with \(\sigma=1/2\).

		Since \(\phi\) is positive and log-concave on \((0,1)\), Lemma~\ref{lem-0325-01} makes \(\sgn(x_i)G(x',t)/\phi(t)\) strictly decreasing there. Hence \(G(x',\cdot)\) has at most one zero in \((0,1)\).
	Proposition \ref{prop-0504-03} gives existence in \((0,1)\).

	For every \(t\geq1\), Proposition \ref{prop-0730-01}{\rm(ii)} and \eqref{0630-07} give \(\sgn(x_i)G(x',t)<0\).
	Under the additional assumption {\rm\textbf{(T3)}}, Lemma~\ref{lem-0325-01}{\rm(ii)} gives \(\sgn(x_i)\partial_tG(x',t_i^*(x'))<0\).
		Thus \eqref{0929-30} follows from \eqref{0630-07}, \eqref{0324-15}, the fact that \(\sgn(x_i)G(x',t)>0\) for \(0<t\leq t_1\), and the uniqueness of the zero.
\end{proof}
The next lemma gives the upper bound \(t_i^*(x')<t_\infty\). Together with Proposition~\ref{prop-0324-02}, it proves Theorem~\ref{thm-0324-01}{\rm(iii)}.
\begin{lemma}\label{lem-0828-01}
Suppose that \(\psi\) satisfies {\rm\textbf{(S)}}, and that \(\phi\) satisfies {\rm\textbf{(T1)}} and {\rm\textbf{(T2)}} and is log-concave on \((0,1)\).
	Let \(1\leq i\leq n-1\), and let \(t_\infty\) and \(t_i^*(x')\) be the zeros of \(T_0(t)\) and \(G(x',t)\), respectively.
	Then \(t_i^*(x')<t_\infty\) for every \(x'\in \Sigma\) with \(|x_i|>1\).
\end{lemma}
\begin{proof}
	Set \(\widetilde G(x',t)=\sgn(x_i)G(x',t)\) and \(\widetilde W(x',t)=\sgn(x_i)W(x',t)\).
	It suffices to show that \(\widetilde G(x',t_\infty)<0\).
	Indeed, since \(\widetilde G(x',t)>0\) for \(0<t<t_1\) by \eqref{0320-27}, the intermediate value theorem then gives a zero in \((t_1,t_\infty)\), and Proposition \ref{prop-0324-02} identifies it with \(t_i^*(x')\).

	By \eqref{260004},
	\begin{align}\label{0324-17}
		\int_0^{t_\infty}\frac{\phi'(s)}{\sqrt{t_\infty-s}}\,ds=0.
	\end{align}
	Since \(Q\) is strictly increasing, \(\widetilde W(x',\cdot)\) is positive and strictly decreasing.
	Split the integral defining \(\widetilde G\) at \(t_1\).
	For \(0<s<t_1\), \(\phi'(s)>0\) and \(\widetilde W(x',t_\infty-s)<\widetilde W(x',t_\infty-t_1)\), so
	\begin{align*}
		\int_0^{t_1}\frac{\phi'(s)}{\sqrt{t_\infty-s}} \widetilde W(x',t_\infty-s)\,ds < & \ \widetilde W(x',t_\infty-t_1) \int_0^{t_1}\frac{\phi'(s)}{\sqrt{t_\infty-s}}\,ds.
	\end{align*}
	For \(t_1<s<t_\infty\), \(\phi'(s)<0\) and \(\widetilde W(x',t_\infty-s)>\widetilde W(x',t_\infty-t_1)\), so
	\begin{align*}
		\int_{t_1}^{t_\infty}\frac{\phi'(s)}{\sqrt{t_\infty-s}} \widetilde W(x',t_\infty-s)\,ds < & \ \widetilde W(x',t_\infty-t_1) \int_{t_1}^{t_\infty}\frac{\phi'(s)}{\sqrt{t_\infty-s}}\,ds.
	\end{align*}
	Adding these inequalities and using \eqref{0324-17} gives \(\widetilde G(x',t_\infty)<0\).
\end{proof}

	\subsubsection{Monotonicity of the separation point}\label{subsect-0829-02}
The next proposition shows that \(t_i^*(r\mathbf e_i')\) and \(t_i^*(-r\mathbf e_i')\) are strictly increasing in \(r\) for sufficiently large \(r\).
This proves the monotonicity statement in Theorem \ref{thm-0324-01}{\rm(iv)}.

\begin{proposition}\label{prop-0421-02}
	Let \(1\leq i\leq n-1\), and let \(w\) solve the Stokes system \eqref{0308-01} with boundary data \eqref{0324-13}, where \(j=n\).
	Suppose that \(\psi\) satisfies {\rm\textbf{(S)}}, and that \(\phi\) satisfies {\rm\textbf{(T1)}}--{\rm\textbf{(T4)}}.
	Then there exists \(C=C(n)>0\) such that \(t_i^*(r\mathbf e_i')\) and \( t_i^*(-r\mathbf e_i')\) are strictly increasing in \(r\) for
	\begin{align}\label{0511-01}
		r > 1+\max\left\{ \left(C\frac{M_2}{M_0}\right)^{\frac{1}{n+8}},\left(C\frac{M_2}{M_0}\right)^{\frac{1}{2}}\right\}.
	\end{align}
\end{proposition}
To prove Proposition \ref{prop-0421-02}, we first prove the following lemma.
\begin{lemma}\label{lem-0326-01}
	Let \(1\leq i\leq n-1\) and \(M>0\), and let \(K\) be given in \eqref{0325-06}.
	Suppose that \(\psi\) satisfies {\rm\textbf{(S)}} and \(\frac{M_2}{M_0}\leq M\).
	Then there exists \(C=C(n)>0\) such that, for \(x_i>1+\max\{(CM)^{\frac{1}{n+8}},(CM)^{\frac12}\}\) and \(t>0\), \(K(x_i\mathbf e_i',t)\) satisfies
	\begin{align}\label{0326-11}
		\partial_{t}\partial_{x_i}\ln K(x_i\mathbf e_i',t)>0.
	\end{align}
\end{lemma}
\begin{proof}
	In this proof, we write \(x'=x_i\mathbf e_i'\).
	We define
	\begin{align*}
		H(z',t):= \frac{\gamma_n z_i}{t^{\frac{n+1}{2}}}P\left( \frac{|z'|^2}{4t}\right), \qquad h(x',y',t) := H(x'-y',t)\psi(y'),
	\end{align*}
	where \(P\) is given in \eqref{0827-11} and \(\gamma_n=\frac{\sqrt\pi}{2^{n+1}\Gamma(\frac{n+3}{2})}\).
	Then
	\begin{align*}
		K(x',t) = \int_{\Sigma}h(x',y',t)\,dy' = \int_{\Sigma}H(x'-y',t)\psi(y')\,dy'.
	\end{align*}
	For \(z_i>0\) define the logarithmic derivatives
	\begin{align*}
		f(z',t)=\partial_{z_i}\ln H(z',t),\qquad g(z',t)=\partial_t\ln H(z',t).
	\end{align*}
	Since \(\psi\) vanishes outside \(\overline{B_1'}\), we may integrate over \(B_1'\).
	There \(x_i-y_i\geq x_i-1>0\), so \(H(x'-y',t)>0\). 
	We first establish the identity
	\begin{align}\label{0326-01}
		\partial_t\partial_{x_i}\ln K(x',t) & = \frac{\displaystyle \int_{B_1'}\partial_t \partial_{x_i}\ln H(x'-y',t) h(x',y',t)\,dy'}{K(x',t)}\nonumber\\
		 & +\frac{\displaystyle \iint_{(B_1')^2}(f(x'-y',t)-f(x'-z',t))(g(x'-y',t)-g(x'-z',t))h(x',y',t)h(x',z',t)\,dy'dz'}{2K(x',t)^2}.
	\end{align}
	To verify \eqref{0326-01}, note that
	\begin{align*}
		\partial_{x_i} \ln K(x',t)
		&= \frac{\displaystyle\int_{B_1'}f(x'-y',t) h(x',y',t)\,dy'}{K(x',t)},
	\end{align*}
	and
	\begin{align*}
		\partial_th(x',y',t) = \partial_t H(x'-y',t)\psi(y') = g(x'-y',t)h(x',y',t).
	\end{align*}

	With \(u=\frac{|x'-y'|^2}{4t}\) and \(L=\ln P\),
	\begin{align}\label{0326-02}
		\partial_t \partial_{x_i} \ln H(x'-y',t) = \partial_t \partial_{x_i}L(u)= \partial_t \left( \frac{x_i-y_i}{2t}L'(u)\right) = -\frac{x_i-y_i}{2t^2}(L'(u) + u L''(u)).
	\end{align}
	Lemma \ref{lem-0324-01}{\rm(v)} gives \(\partial_t\partial_{x_i}\ln H(x'-y',t)>0\) for \(x_i>1\) and \(|y'|\leq1\).

	We next bound the first term on the right-hand side of \eqref{0326-01}.
	Set \(\Phi(u)=-(uL'(u))'\).
	Then Lemma \ref{lem-0324-01}{\rm(v)} gives \(\Phi(u)>0\), and
	\begin{align*}
		\Phi'(u) = -(uL'(u))'' = -uL'''(u) -2L''(u) = -\frac{1}{u}(u^2L''(u))',
	\end{align*}
	so Lemma \ref{lem-0324-01}{\rm(vi)} shows that \(\Phi\) is strictly decreasing.
	For \(|y'|\leq1\), \(x'=x_i\mathbf e_i'\), and \(x_i>1\), \eqref{0326-02} gives
	\begin{align*}
		\partial_t \partial_{x_i}\ln H(x'-y',t) = \frac{x_i-y_i}{2t^2}\Phi(u) \geq \frac{x_i-1}{2t^2} \Phi\left( \frac{(x_i+1)^2}{4t}\right).
	\end{align*}
	Thus
	\begin{align}\label{0326-03}
		\begin{split}
			\int_{B_1'}\partial_t \partial_{x_i} \ln H(x'-y',t) h(x',y',t)\,dy' & \geq \frac{x_i-1}{2t^2}\Phi\left( \frac{(x_i+1)^2}{4t}\right) \int_{B_1'}h(x',y',t)\,dy'.
		\end{split}
	\end{align}
	We next bound the second term on the right-hand side of \eqref{0326-01}.
	By the mean value theorem,
	\begin{align*}
		|(f(x'-y',t)-f(x'-z',t)) (g(x'-y',t)-g(x'-z',t))| 
		\leq |y'-z'|^2 M_f(x',t) M_g(x',t),
	\end{align*}
	where \(\nabla'\) denotes the gradient with respect to the spatial argument, \(M_f(x',t):=\sup_{|y'|\leq1}|\nabla'f(x'-y',t)|\), and \(M_g(x',t):=\sup_{|y'|\leq1}|\nabla'g(x'-y',t)|\).
	Hence
	\begin{align}\label{0326-04}
		\begin{split}
			 & \left|\iint_{(B_1')^2}(f(x'-y',t)-f(x'-z',t))(g(x'-y',t)-g(x'-z',t))h(x',y',t)h(x',z',t)\,dy'dz' \right| \\
			 & \quad \leq M_f(x',t) M_g(x',t) \iint_{(B_1')^2}|y'-z'|^2h(x',y',t)h(x',z',t)\,dy'dz'\\
			 & \quad \leq 2M_f(x',t) M_g(x',t)K(x',t) \int_{B_1'}|y'|^2H(x'-y',t) \psi(y')\,dy'.
		\end{split}
	\end{align}
	To estimate the remaining factor in \eqref{0326-04} relative to \(K(x',t)\), Lemma \ref{lem-0324-01}{\rm(i)} gives
	\begin{align}\label{0326-05}
		\frac{\gamma_n(x_i-1)}{t^{\frac{n+1}{2}}} P\left( \frac{(x_i+1)^2}{4t}\right) \leq H(x'-y',t) \leq \frac{\gamma_n(x_i+1)}{t^{\frac{n+1}{2}}}P\left( \frac{(x_i-1)^2}{4t}\right).
	\end{align}

	We next estimate \(M_f\) and \(M_g\).
For \(1\leq j\leq n-1\), differentiation gives
	\begin{align*}
		\partial_{y_j}f(x'-y',t) = \frac{\delta_{ij}}{(x_i-y_i)^2} -\frac{(x_i-y_i)(x_j-y_j)}{4t^2} L''(u)-\frac{\delta_{ij}}{2t}L'(u).
	\end{align*}
	Thus, by \eqref{0602-01} and \eqref{0529-02},
	\begin{align*}
		|\nabla' f(x'-y',t) |
		 & \leq \frac{1}{(x_i-y_i)^2} + \frac{3n}{2t+|x'-y'|^2}\leq \frac{3n+1}{(x_i-y_i)^2}\leq \frac{3n+1}{(x_i-1)^2}.
	\end{align*}

	Similarly, for \(M_g\),
	\begin{align*}
		\partial_{y_j}g(x'-y',t) = \partial_{y_j}\left( -\frac{n+1}{2t}- \frac{1}{t}uL'(u)\right) = -\frac{x_j-y_j}{2t^2}\Phi(u).
	\end{align*}
	Thus \eqref{0731-32} gives
	\begin{align*}
						|\nabla' g(x'-y',t)| & = \frac{|x'-y'|}{2t^2}\Phi(u)\leq \frac{l_n}{18n}\frac{(x_i+1)^5}{(x_i-1)^4t^2}\Phi\left( \frac{(x_i+1)^2}{4t}\right).
	\end{align*}

	Thus \eqref{0326-04}, \eqref{0326-05}, and \eqref{0529-03} give
	\begin{align*}
		 & \frac{\displaystyle \left|\iint_{(B_1')^2}(f(x'-y',t)-f(x'-z',t))(g(x'-y',t)-g(x'-z',t))h(x',y',t)h(x',z',t)\,dy'dz'\right|}{K(x',t)^2}\\
		 & \leq \frac{2^{\frac{n+4}{2}}}{\sqrt{\pi}}\Gamma\left(\frac{n+3}{2}\right) M_f(x',t) M_g(x',t) \frac{M_2}{M_0} \frac{(x_i+1)^{n+1}}{(x_i-1)^{n+1}}.
	\end{align*}
	Thus
	\begin{align}\label{0326-07}
		\begin{split}
			 & \frac{\displaystyle \left| \iint_{(B_1')^2}(f(x'-y',t)-f(x'-z',t))(g(x'-y',t)-g(x'-z',t))h(x',y',t)h(x',z',t)\,dy'dz' \right|}{K(x',t)^2} \\
			 & \leq C_{\mathrm{cov}}(n)\frac{(x_i+1)^{n+6}}{(x_i-1)^{n+7}t^2}\Phi\left( \frac{(x_i + 1)^2}{4t}\right) \frac{M_2}{M_0},
		\end{split}
	\end{align}
	where \(C_{\mathrm{cov}}(n) = \frac{2^{\frac{n+2}{2}}(3n+1)l_n}{9n\sqrt{\pi}}\Gamma\left(\frac{n+3}{2}\right)\).
	Then \eqref{0326-01}, \eqref{0326-03}, and \eqref{0326-07} give
	\begin{align*}
		\partial_t\partial_{x_i}\ln K(x',t) 
		 & \geq \frac{x_i-1}{2t^2}\Phi\left(\frac{(x_i+1)^2}{4t}\right) - \frac{C_{\mathrm{cov}}(n)}{2} \frac{(x_i+1)^{n+6}}{(x_i-1)^{n+7}t^{2}}\Phi\left(\frac{(x_i+1)^2}{4t}\right) \frac{M_2}{M_0}.
	\end{align*}
	Thus \(\partial_t \partial_{x_i}\ln K(x',t)>0\) if
	\begin{align*}
		\frac{M_2}{M_0}<\frac{1}{C_{\mathrm{cov}}(n)} \frac{(x_i-1)^{n+8}}{(x_i+1)^{n+6}}.
	\end{align*}
	If \(1<x_i\leq3\), then \(x_i+1\leq4\), and hence
	\begin{align*}
		\frac{1}{C_{\mathrm{cov}}(n)}\frac{(x_i-1)^{n+8}}{(x_i+1)^{n+6}} \geq\frac{(x_i-1)^{n+8}}{2^{2n+12}C_{\mathrm{cov}}(n)}>M
	\end{align*}
	when \(x_i>1+(2^{2n+12}C_{\mathrm{cov}}(n)M)^{1/(n+8)}\).
	If \(x_i\geq3\), then \(x_i+1\leq2(x_i-1)\), so
	\begin{align*}
		\frac{1}{C_{\mathrm{cov}}(n)}\frac{(x_i-1)^{n+8}}{(x_i+1)^{n+6}} \geq\frac{(x_i-1)^2}{2^{n+6}C_{\mathrm{cov}}(n)}>M
	\end{align*}
	when \(x_i>1+(2^{n+6}C_{\mathrm{cov}}(n)M)^{1/2}\).
	Thus, with \(C=2^{2n+12}C_{\mathrm{cov}}(n)\), \(x_i>1+\max\{(CM)^{1/(n+8)},(CM)^{1/2}\}\) implies \eqref{0326-11}.

\end{proof}

\begin{proof}[Proof of Proposition \ref{prop-0421-02}]
	In this proof, we write \(x'=x_i\mathbf e_i'\).
		Recall that \(G(x',t_i^*(x'))=0\).

	Applying Lemma \ref{lem-0325-01} as in the proof of Proposition \ref{prop-0324-02} gives \(\partial_tG(x',t_i^*(x'))<0\) for every \(x_i>1\).

	We next claim that \(\partial_{x_i}G(x',t_i^*(x'))>0\) when \(x_i\) satisfies \eqref{0511-01}.
	We write
	\begin{align*}
		\partial_{x_i}G(x',t_i^*) & =\int_0^{t_i^*}\phi'(s)\partial_{x_i}\ln K(x',t_i^*-s) K(x',t_i^*-s)\,ds.
	\end{align*}
	Then Lemma \ref{lem-0326-01}, applied with \(M=M_2/M_0\), shows that \(\partial_{x_i}\ln K(x',t)\) is strictly increasing in \(t\) when \eqref{0511-01} holds.
	For \(0<s<t_1\), we have \(\phi'(s)>0\) and \(\partial_{x_i}\ln K(x',t_i^*-s)>\partial_{x_i}\ln K(x',t_i^*-t_1)\).
	For \(t_1<s<t_i^*\), we have \(\phi'(s)<0\) and \(\partial_{x_i}\ln K(x',t_i^*-s)<\partial_{x_i}\ln K(x',t_i^*-t_1)\).
	Hence
	\begin{align*}
		\partial_{x_i}G(x',t_i^*) & >\partial_{x_i}\ln K(x',t_i^*-t_1) \int_0^{t_i^*}\phi'(s)K(x',t_i^*-s)\,ds =\partial_{x_i}\ln K(x',t_i^*-t_1)G(x',t_i^*)=0.
	\end{align*}
	Lemma~\ref{lem-0424-20}, with \(f=\phi'\) and \(\sigma=1/2<\beta\), now gives joint \(C^1\) regularity of \(G(x_i\mathbf e_i',t)\).

	Since \(\partial_tG(x',t_i^*(x'))<0\), the implicit function theorem and the uniqueness in Proposition \ref{prop-0324-02} show that \(t_i^*(x_i\mathbf e_i')\) is \(C^1\).
	Hence
\begin{align*}
\frac{d}{dx_i}t_i^*(x_i\mathbf e_i')=-\frac{\partial_{x_i}G(x',t_i^*(x'))}{\partial_tG(x',t_i^*(x'))}>0,
\end{align*}
which implies that \(t_i^*(x_i\mathbf e_i')\) is strictly increasing under \eqref{0511-01}.

	For \(t_i^*(-r\mathbf{e}_i')\), let \(R_i y'=y'-2y_i\mathbf e_i'\) and set \(\psi_i^-(y')=\psi(R_i y')\).
	This function satisfies {\rm\textbf{(S)}} and has the same \(M_0\) and \(M_2\) as \(\psi\).
	Changing variables in \eqref{0320-27} gives \(G_\psi(-r\mathbf e_i',t)=-G_{\psi_i^-}(r\mathbf e_i',t)\).
	Applying the result for \(t_i^*(r\mathbf{e}_i')\) to \(\psi_i^-\) implies that \(t_i^*(-r\mathbf e_i')\) is strictly increasing in \(r\) under \eqref{0511-01}.
\end{proof}

\begin{remark}\label{rmk-0808-01}
Under {\rm\textbf{(T1)}}--{\rm\textbf{(T4)}}, the proof of Proposition \ref{prop-0324-02} gives \(\sgn(x_i)\partial_tG(x',t_i^*(x'))<0\) for every \(x'\in\Sigma\) with \(|x_i|>1\).
	By \eqref{0630-07}, \(\partial_t\partial_{x_n}w_i(x',0,t_i^*(x'))\) has sign \(-\sgn(x_i)\).
	Thus \(\sigma=1\) in Proposition \ref{prop-0512-02}{\rm(i)}.
\end{remark}

\subsubsection{Asymptotic behavior of the separation point}
We determine the asymptotic behavior of \(t_i^*(x_i\mathbf e_i')\) as \(|x_i|\to\infty\).

Recall \(T_k(t)\) from \eqref{0323-02}.
	For \(k\geq1\), the integrand \((t-s)^{\frac{2k-1}{2}}\phi'(s)\) vanishes at \(s=t\), so differentiation under the integral gives the relation \(T_k'(t) = \frac{2k-1}{2}T_{k-1}(t)\).
The next lemma gives the expansion used below.
\begin{lemma}\label{lem-0731-10}
	Let \(n\geq2\) and let \(\psi\) satisfy {\rm\textbf{(S)}}.
	Then, for \(0<r<1\) and \(0<\ep<\frac12\),
	\begin{align}\label{0731-10}
		 & \int_\Sigma\frac{1-\ep y_i}{|\mathbf e_i'-\ep y'|^n}\psi(y') Q\left(\frac{|\mathbf e_i'-\ep y'|^2}{4\ep^2r}\right)\,dy' =M_0+(n-1)m_i\ep\nonumber\\
		 & \quad+\left(\frac{n^2}{2}\int_\Sigma y_i^2\psi(y')\,dy'-\frac n2M_2-nM_0r\right)\ep^2\nonumber \\
		 & \quad+\left[\frac{n(n+1)}6\int_\Sigma \left((n+2)y_i^3-3y_i|y'|^2\right)\psi(y')\,dy' -n(n+1)m_ir\right]\ep^3+R(r,\ep),
	\end{align}
	where \(|R(r,\ep)|\leq C(n)M_0\ep^4\).
\end{lemma}
\begin{proof}
	Since \(|y'|\leq1\) and \(0<\ep<\frac12\), the triangle inequality gives \(\frac12<|\mathbf e_i'-\ep y'|<\frac32\).
	Thus \(\frac{|\mathbf e_i'-\ep y'|^2}{4\ep^2r}>\frac{1}{16\ep^2r}\), and \eqref{0322-12} of Lemma \ref{lem-0424-01} gives
	\begin{align}\label{0731-11}
		Q\left(\frac{|\mathbf e_i'-\ep y'|^2}{4\ep^2r}\right) =1-\frac{n\ep^2r}{|\mathbf e_i'-\ep y'|^2}+O(\ep^4),
	\end{align}
	uniformly for \(0<r<1\), \(0<\ep<\frac12\).
	A Taylor expansion in \(\ep\) gives
	\begin{align*}
		\frac{1-\ep y_i}{|\mathbf e_i'-\ep y'|^n} &=1+(n-1)y_i\ep+\frac{n}{2}\left(n y_i^2-|y'|^2\right)\ep^2 +\frac{n(n+1)}{6}\bigl((n+2)y_i^3-3y_i|y'|^2\bigr)\ep^3+O(\ep^4), \\
		\frac1{|\mathbf e_i'-\ep y'|^2} &=1+2y_i\ep+O(\ep^2),
	\end{align*}
	uniformly for \(|y'|\leq1\).
	Multiplying and collecting powers of \(\ep\) and then integrating against \(\psi\) gives \eqref{0731-10}.
\end{proof}

The next proposition gives the asymptotic expansion of \(t_i^*(x_i\mathbf e_i')\) as \(|x_i|\to\infty\).
This proves the asymptotic statement in Theorem \ref{thm-0324-01}{\rm(iv)}.
\begin{proposition}\label{prop-0731-11}
	Let \(1\leq i\leq n-1\), and let \(w\) solve the Stokes system \eqref{0308-01} with boundary data \eqref{0324-13}, where \(j=n\).
	Suppose that \(\psi\) satisfies {\rm\textbf{(S)}}, and that \(\phi\) satisfies {\rm\textbf{(T1)}}--{\rm\textbf{(T4)}}.
	Then
	\begin{align}\label{0731-12}
		t_i^*(x_i\mathbf e_i')=t_\infty+\frac{nT_1(t_\infty)}{T_0'(t_\infty)}x_i^{-2}+O(|x_i|^{-1-2\beta}) \qquad (|x_i|\rightarrow\infty).
	\end{align}
	If, in addition, \(\phi\) is \(C^{2,\alpha_1}\) near \(t_\infty\) for some \(\alpha_1\in(1/2,1)\), then
	\begin{align}\label{0731-13}
		t_i^*(x_i\mathbf e_i')=t_\infty+\frac{nT_1(t_\infty)}{T_0'(t_\infty)}x_i^{-2} +\frac{2nT_1(t_\infty)}{T_0'(t_\infty)}\frac{m_i}{M_0}x_i^{-3}+O(x_i^{-4}),\qquad (|x_i|\rightarrow \infty).
	\end{align}
\end{proposition}

\begin{proof}
It suffices to prove \eqref{0731-12} and \eqref{0731-13} for \(x_i>0\).
Indeed, the reflection \(R_i y'=y'-2y_i\mathbf e_i'\) and the function \(\psi_i^-(y')=\psi(R_i y')\) from the proof of Proposition \ref{prop-0421-02} give \(G_\psi(-r\mathbf e_i',t)=-G_{\psi_i^-}(r\mathbf e_i',t)\).
The function \(\psi_i^-\) has the same \(M_0\) and \(M_2\), while its \(i\)-th first moment is \(-m_i\).
Applying these expansions for \(x_i>0\) to \(\psi_i^-\) therefore yields them for \(x_i<0\), including the sign of the \(x_i^{-3}\) term.

Write \(\tau^*(x_i):=t_i^*(x_i\mathbf e_i')\), which is well-defined for large \(x_i\) by Propositions \ref{prop-0504-03} and \ref{prop-0324-02}, and satisfies \(t_1<\tau^*(x_i)<t_\infty\) by Lemma \ref{lem-0828-01}.
Applying \eqref{0731-10} with \(\ep=x_i^{-1}\) in \eqref{0320-27}, and using \eqref{0511-02}, gives
	\begin{align}\label{0731-14}
		x_i^{n-1}G(x_i\mathbf e_i',t) & =\left(M_0+\frac{(n-1)m_i}{x_i}+O(x_i^{-2})\right)T_0(t)-\left(\frac{nM_0}{x_i^2}+\frac{n(n+1)m_i}{x_i^3}\right)T_1(t) +O(x_i^{-4}),
	\end{align}
	uniformly for \(t\) in compact subsets of \((0,1)\).
	Since \(G(x_i\mathbf e_i',\tau^*(x_i))=0\) and \(M_0>0\), \eqref{0731-14} gives
	\begin{align}\label{0731-16}
		T_0(\tau^*(x_i)) =\left(\frac{n}{x_i^2}+\frac{2nm_i}{M_0x_i^3}+O(x_i^{-4})\right) T_1(\tau^*(x_i))+O(x_i^{-4})=O(x_i^{-2}).
	\end{align}
	As \(\tau^*(x_i)\in[t_1,t_\infty]\) and \(t_\infty\) is the only zero of \(T_0\) there (see Remark \ref{rmk-0730-03}), every limit point of \(\tau^*(x_i)\) is \(t_\infty\).
	Hence \(\tau^*(x_i)\to t_\infty\).
	Remark~\ref{rmk-0730-03} gives \(T_0'(t_\infty)<0\) under {\rm\textbf{(T1)}}--{\rm\textbf{(T4)}}.
	By \eqref{0731-16} and the nonzero derivative at \(t_\infty\), \(h:=\tau^*(x_i)-t_\infty=O(x_i^{-2})\).
	Lemma~\ref{lem-0424-20} with \(K(u)=u^{-1/2}\) and \(f=\phi'\) gives \(T_0\in C^{1,\beta-1/2}_{\mathrm{loc}}(0,1)\). Hence
	\begin{align*}
		T_0(t_\infty+h)=T_0'(t_\infty)h+O(|h|^{\beta+1/2}).
	\end{align*}
	In \eqref{0731-16}, \(T_1(t_\infty+h)=T_1(t_\infty)+O(|h|)\) because \(T_1'=T_0/2\).
	The two errors are therefore \(O(x_i^{-1-2\beta})\) from the Taylor expansion of \(T_0\) and \(O(x_i^{-3})\) from \eqref{0731-14}.
	Since \(\beta<1\), this proves \eqref{0731-12}.

	Assume now that \(\phi\) is \(C^{2,\alpha_1}\) on a neighborhood of \(t_\infty\) in \((0,1)\), with \(\alpha_1\in(1/2,1)\).
	Choose \(\eta\in C_c^\infty(0,1)\) supported in that neighborhood and equal to one near \(t_\infty\), and split \(T_0\) into the convolutions of \(u^{-1/2}\) with \(\eta\phi'\) and with \((1-\eta)\phi'\).
	For \(t\) near \(t_\infty\), the second integrand vanishes for \(s\) near \(t\), so the second convolution is smooth in \(t\).
	Since \((\eta\phi')(0)=0\), the derivative of the first convolution is the convolution of \(u^{-1/2}\) with \((\eta\phi')'\in C^{\alpha_1}([0,1])\), and Lemma~\ref{lem-0424-20} gives \(T_0\in C^2\) near \(t_\infty\).
	Moreover, \(T_1'=\frac12T_0\), so \(T_1'(t_\infty)=0\) and \(T_1\in C^2\).
	Consequently,
	\begin{align*}
		T_0(\tau^*(x_i))=T_0'(t_\infty)(\tau^*(x_i)-t_\infty)+O(x_i^{-4}), \qquad T_1(\tau^*(x_i))=T_1(t_\infty)+O(x_i^{-4}).
	\end{align*}
	Substituting these expansions into \eqref{0731-16} gives
	\begin{align*}
		T_0'(t_\infty)(\tau^*(x_i)-t_\infty) =\frac{nT_1(t_\infty)}{x_i^2} +\frac{2nm_iT_1(t_\infty)}{M_0x_i^3}+O(x_i^{-4}),
	\end{align*}
	which is \eqref{0731-13}.
\end{proof}

\subsection{Normal velocity component}\label{260056}
\subsubsection{Formula for the second normal derivative of \(w_n\)}
The divergence-free condition and the normal boundary data give \(\partial_{x_n}w_n(x',0,t)=0\), so we consider \(\partial_{x_n}^2w_n(x',0,t)\).
Define
\begin{align}\label{0419-01}
	k_i(x,t) := \int_0^t\int_{\Sigma}\partial_{x_n}\Gamma(x-z',t-s)g_i(z',s)\,dz'ds, \quad 1\leq i\leq n.
\end{align}
\begin{proposition}\label{prop-0504-01}
	Let \(w\) solve the Stokes system \eqref{0308-01} with boundary data \eqref{0324-13}, where \(j=n\).
	Suppose that \(\psi\) satisfies {\rm\textbf{(S)}}, and that \(\phi\) satisfies {\rm\textbf{(T1)}}.
	Then, for every \(x'\in \Sigma\) and \(t>1\),
	\begin{align}\label{0428-01}
		\partial_{x_n}^2w_n(x',0,t) = \frac{\Gamma\left(\frac{n}{2}\right)}{(4\pi)^{\frac{n}{2}}\Gamma\left(\frac{n-1}{2}\right)}\int_{\Sigma}\psi(y')\int_0^1 \frac{\phi(s)}{(t-s)^{\frac{n+3}{2}}}{}_1F_1\left(\frac{n}{2}, \frac{n-1}{2}; - \frac{|x'-y'|^2}{4(t-s)}\right)\,dsdy'.
	\end{align}
	For every \(x'\notin\supp\psi\) and \(0<t<1\),
	\begin{align}\label{0428-02}
		\begin{split}
			\partial_{x_n}^2w_n(x',0,t) & = - \frac{\Gamma\left(\frac{n}{2}\right)}{2^{n-1}\pi^{\frac{n}{2}}\Gamma\left(\frac{n-1}{2}\right)}\int_{\Sigma}\psi(y')\int_0^t\frac{1}{(t-s)^{\frac{n+1}{2}}}\left( \phi'(s) {}_1F_1\left(\frac{n}{2}, \frac{n-1}{2}; - \frac{|x'-y'|^2}{4(t-s)}\right)\right. \\
			& \quad\left. + \frac{n}{2}\frac{\phi(s)}{t-s}{}_1F_1\left(\frac{n+2}{2}, \frac{n-1}{2}; -\frac{|x'-y'|^2}{4(t-s)}\right)\right)\,dsdy'.
		\end{split}
	\end{align}
\end{proposition}
\begin{proof}
	We first show that
	\begin{align}\label{0420-04}
		\partial_{x_n}^2 w_n(x',0,t)= - 4 \sum_{k=1}^{n-1}\mathrm{p.v.} \int_{\Sigma}\partial_{x_k}N(x'-y',0)\partial_{y_k}\partial_{x_n}k_n(y',0,t)\,dy',
	\end{align}
	where \(k_n\) is given in \eqref{0419-01}.
	Indeed, \eqref{0330-01}, \eqref{0308-02}, and \eqref{0419-01} give
	\begin{align*}
		w_n(x,t) & = -2k_n(x,t) +4\int_0^t\int_{\Sigma}L_{nn}(x-y',t-s)g_n(y',s)\,dy'ds +2\int_{\Sigma}\partial_{x_n}N(x-y')g_n(y',t)\,dy'.
	\end{align*}
	Here \eqref{0419-02} gives
	\begin{align*}
		\int_0^t\int_{\Sigma} & L_{nn}(x-y',t-s)g_n(y',s)\,dy'ds\\
		& = \partial_{x_n}\int_{\Sigma\times(0,x_n)}\int_0^t\int_{\Sigma}\partial_{z_n}\Gamma(z,t-s)\partial_{x_n}N(x-y'-z)g_n(y',s)\,dy'dsdz\\
		& = \partial_{x_n}\int_{\Sigma\times(0,x_n)}\partial_{x_n}N(y)k_n(x-y,t)\,dy.
	\end{align*}
	Thus
	\begin{align*}
		\partial_{x_n}\int_0^t\int_{\Sigma} & L_{nn}(x-y',t-s)g_n(y',s)\,dy'ds\\
		& = -\frac12\partial_{x_n}\int_{\Sigma}\partial_{x_n}N(y',x_n)g_n(x'-y',t)\,dy' + \int_{\Sigma}\partial_{x_n}N(y',x_n)\partial_{x_n}k_n(x'-y',0,t)\,dy'\\
		& \quad + \int_{\Sigma\times(0,x_n)}\partial_{y_n}N(y)\partial_{x_n}^2k_n(x-y,t)\,dy,
	\end{align*}
where we used \eqref{260269} and the boundary limit \(\lim_{x_n\rightarrow0^+}k_n(x,t)=- \frac12 g_n(x',t)\).
	Thus
	\begin{align}\label{0420-02}
		\begin{split}
			\partial_{x_n}w_n(x,t)
			& = \left( -2\partial_{x_n}k_n(x,t) +2\partial_{x_n}k_n(x',0,t) +4\int_{\Sigma\times(0,x_n)}\partial_{y_n}N(y)\partial_{x_n}^2 k_n(x-y,t)\,dy\right)\\
			& \quad + \left( 4\int_{\Sigma}\partial_{x_n}N(x-y')\partial_{x_n}k_n(y',0,t)\,dy' - 2\partial_{x_n}k_n(x',0,t)\right)\\
			& =: I_1(x,t)+ I_2(x,t).
		\end{split}
	\end{align}

	For \(I_1\), we write
	\begin{align*}
		I_1(x,t)
		& = 2\int_0^{x_n}\left( 2\int_{\Sigma}\partial_{y_n}N(y)\partial_{x_n}^2k_n(x-y,t)\,dy'-\partial_{x_n}^2k_n(x',x_n-y_n,t)\right)\,dy_n.
	\end{align*}
	For fixed \(t>0\), choose \(T>t\), and let \(\chi\in C^\infty([0,T])\) vanish near \(0\) and equal \(1\) near \(t\).
	Decompose \(k_n=k_n^{(1)}+k_n^{(2)}\), where
	\begin{align*}
		k_n^{(1)}(x,s) & :=\int_0^s\int_\Sigma \partial_{x_n}\Gamma(x-z',s-\tau)\psi(z')\chi(\tau)\phi(\tau)\,dz'd\tau, \\
		k_n^{(2)}(x,s) & :=\int_0^s\int_\Sigma \partial_{x_n}\Gamma(x-z',s-\tau)\psi(z')(1-\chi(\tau))\phi(\tau)\,dz'd\tau.
	\end{align*}
	Let \(\alpha_0\in(0,\min\{\alpha,2\beta\})\).
			Since \(\chi\phi\) vanishes near \(0\), Lemma \ref{lem-260232} of Appendix \ref{260231} applies to \(k_n^{(1)}\).
	There exists \(d>0\) such that \(1-\chi(\tau)=0\) for \(t-d<\tau<t\).
	Thus \(t-\tau\geq d\) in the formula for \(k_n^{(2)}(\cdot,t)\), and the bounds of the heat kernel give
	\begin{align}\label{0923-01}
		 & \norm{\partial_{x_n}^2k_n(\cdot,t)}_{C_{x'}^{\alpha_0}L_{x_n}^\infty(\R_+^n)} +\norm{\nabla'\partial_{x_n}k_n(\cdot,t)}_{C_{x'}^{\alpha_0}L_{x_n}^\infty(\R_+^n)} +\norm{\partial_{x_n}k_n(\cdot,t)}_{L^\infty(\R_+^n)}\leq C_t\norm{\psi}_{C^{2,\alpha}(\Sigma)} \norm{\phi}_{C^{1,\beta}([0,T])}.
	\end{align}

	Then the estimate
	\begin{align*}
		\left| 2\int_{\Sigma}\partial_{y_n}N(y)f(x-y)\,dy' - f(x',x_n-y_n)\right| \leq C\norm{f}_{C_{x'}^{\alpha_0} L_{x_n}^\infty(\R_+^n)}y_n^{\alpha_0},
	\end{align*}
	applied with \(f=\partial_{x_n}^2k_n(\cdot,t)\) gives \(|I_1(x,t)| \leq C_t\norm{\psi}_{C^{2,\alpha}(\Sigma)} \norm{\phi}_{C^{1,\beta}([0,T])}x_n^{\alpha_0+1}\).
	Hence
	\begin{align}\label{0420-03}
		\lim_{x_n\rightarrow 0^+}\frac{I_1(x,t)}{x_n}=0.
	\end{align}

	For \(I_2\), harmonicity of \(N\) and integration by parts in \(y'\) give, for \(x_n>0\),
	\begin{align*}
		\partial_{x_n}I_2 (x,t) 
		 & = -4\sum_{k=1}^{n-1}\int_{\Sigma}\partial_{x_k}N(x-y')\partial_{y_k}\partial_{x_n}k_n(y',0,t)\,dy'.
	\end{align*}
	Using \(|\partial_{x_k}N(y',x_n)|\leq C|y'|^{1-n}\), the H\"{o}lder estimate \eqref{0923-01} and \eqref{0713-29} give
	\begin{align}\label{0420-01}
		\lim_{x_n\rightarrow 0^+} \partial_{x_n}I_2(x,t)= -4\sum_{k=1}^{n-1}\mathrm{p.v.}\int_{\Sigma}\partial_{x_k}N(x'-y',0)\partial_{y_k}\partial_{x_n}k_n(y',0,t)\,dy'.
	\end{align}

The limit \eqref{260267} and \eqref{0713-29} also give \(I_2(x,t)\to0\).
	Hence the fundamental theorem of calculus and \eqref{0420-01} imply
	\begin{align*}
		\lim_{x_n\to0^+}\frac{I_2(x,t)}{x_n} 
		=-4\sum_{k=1}^{n-1}\mathrm{p.v.}\int_{\Sigma}\partial_{x_k}N(x'-y',0) \partial_{y_k}\partial_{x_n}k_n(y',0,t)\,dy'.
	\end{align*}
	It follows from \eqref{0420-02}, \eqref{0420-03}, and the limit of \(I_2\) that \(\partial_{x_n}w_n(x',0,t)=0\).
	Consequently,
	\begin{align}\label{261001}
		\partial_{x_n}^2w_n(x',0,t) =\lim_{x_n\to0^+}\frac{I_1(x,t)+I_2(x,t)}{x_n}.
	\end{align}
	Combining \eqref{0420-03} and \eqref{261001} gives \eqref{0420-04}.

	We now derive \eqref{0428-01} and \eqref{0428-02}.
	To compute \(\partial_{x_n}k_n(x',0,t)\), integration by parts and \(g_n(x',0)=0\) give
	\begin{align*}
		\partial_{x_n}k_n(x,t) & =\int_0^t\int_{\Sigma}\partial_{x_n}^2\Gamma(x-y',t-s)g_n(y',s)\,dy'ds\\
		& = \int_0^t\int_{\Sigma}\Gamma(x-y',t-s)(\partial_s-\Delta_{y'})g_n(y',s)\,dy'ds.
	\end{align*}
	The identity \eqref{0713-28} gives
	\begin{align*}
		\partial_{x_n}k_n(x',0,t) = \int_0^t \int_{\Sigma}\Gamma(x'-y',0,t-s)(\partial_s-\Delta_{y'})g_n(y',s)\,dy'ds,
	\end{align*}
	and therefore
	\begin{align*}
		\mathcal{F}_{x'}(\partial_{x_n}k_n)(\xi',0,t) 
		& = \frac{1}{\sqrt{4\pi}}\left( \int_0^t \frac{e^{-|\xi'|^2(t-s)}}{\sqrt{t-s}}\phi'(s)\,ds + \int_0^t |\xi'|^2\frac{e^{-|\xi'|^2(t-s)}}{\sqrt{t-s}}\phi(s)\,ds\right) \mathcal{F}_{x'}\psi(\xi').
	\end{align*}
	Taking the Fourier transform in \(x'\) in \eqref{0420-04} gives
	\begin{align}\label{0429-06}
		\mathcal{F}_{x'}(\partial_{x_n}^2w_n)(\xi',0,t) = -\frac{1}{\sqrt{\pi}}\left(|\xi'|\mathcal{F}_{x'}\psi(\xi')\int_0^t \frac{\phi'(s)}{\sqrt{t-s}}e^{-|\xi'|^2(t-s)}\,ds +|\xi'|^3 \mathcal{F}_{x'}\psi(\xi') \int_0^t \frac{\phi(s)}{\sqrt{t-s}} e^{-|\xi'|^2(t-s)}\,ds \right).
	\end{align}

	If \(t>1\), integration by parts gives
	\begin{align*}
		\int_0^1 \phi'(s) \frac{e^{-|\xi'|^2(t-s)}}{\sqrt{t-s}}\,ds = - \int_0^1 \frac{\phi(s)}{\sqrt{t-s}}|\xi'|^2 e^{-|\xi'|^2(t-s)}\,ds - \frac12 \int_0^1 \frac{\phi(s)}{(t-s)^{\frac32}}e^{-|\xi'|^2(t-s)}\,ds,
	\end{align*}
	and thus
	\begin{align}\label{0429-05}
		\mathcal{F}_{x'}(\partial_{x_n}^2w_n)(\xi',0,t) = \frac{|\xi'|}{\sqrt{4\pi}}\mathcal{F}_{x'}\psi(\xi')\int_0^1 \frac{\phi(s)}{(t-s)^{\frac32}}e^{-|\xi'|^2(t-s)}\,ds.
	\end{align}

By \eqref{0401-07} and Lemma \ref{lem-0331-01},
	\begin{align*}
		\mathcal{F}_{\xi'}^{-1}\left(|\xi'|^ke^{-|\xi'|^2(t-s)}\right)= \frac{\Gamma\left( \frac{n-1+k}{2}\right)}{(4\pi)^{\frac{n-1}{2}}\Gamma\left(\frac{n-1}{2}\right)} \frac{1}{(t-s)^{\frac{n-1+k}{2}}}{}_1F_1\left(\frac{n-1+k}{2}, \frac{n-1}{2}; - \frac{|x'|^2}{4(t-s)}\right),\quad k \geq0.
	\end{align*}
	For \(t>1\), \(t-s\geq t-1>0\) on the interval of integration in \eqref{0429-05}, so Fourier inversion under the integral gives \eqref{0428-01}.

	By \eqref{0514-02}, the inverse Fourier transforms of \(|\xi'|^ke^{-|\xi'|^2(t-s)}\), \(k=1,3\), in \eqref{0429-06} are uniformly bounded there for \(0<t-s\leq t\) and \(y'\in\supp\psi\), which proves \eqref{0428-02}.
\end{proof}

To begin our analysis of separation points of \(w_n\), we introduce \(G_1(x',t)\) for \(t>1\) and \(G_2(x',t)\) for \(x'\notin\supp\psi\) and \(0<t<1\) as follows.
\begin{align}
	 & G_1(x',t):= \int_{\Sigma}\psi(y')\int_0^1 \frac{\phi(s)}{(t-s)^{\frac{n+3}{2}}}{}_1F_1\left(\frac{n}{2}, \frac{n-1}{2};-\frac{|x'-y'|^2}{4(t-s)}\right)\,dsdy',\label{0504-10}\\
	 & G_2(x',t) := \int_\Sigma \psi(y')\int_0^t \frac{1}{(t-s)^{\frac{n+1}{2}}}\left[\phi'(s) {}_1F_1\left(\frac{n}{2}, \frac{n-1}{2}; -\frac{|x'-y'|^2}{4(t-s)}\right)\right.\nonumber \\
	 & \quad \left. + \frac{n\phi(s)}{2(t-s)}{}_1F_1\left(\frac{n+2}{2}, \frac{n-1}{2};- \frac{|x'-y'|^2}{4(t-s)}\right)\right]\,dsdy'.\label{0508-01}
\end{align}
Then \eqref{0428-01} and \eqref{0428-02} give
\begin{align}
	\partial_{x_n}^2 w_n(x',0,t) & = \frac{\Gamma\left(\frac{n}{2}\right)}{(4\pi)^{\frac{n}{2}}\Gamma\left(\frac{n-1}{2}\right)}G_1(x',t) \,\text{for } t>1,\label{0517-10}\\
	\partial_{x_n}^2w_n(x',0,t)& = -\frac{\Gamma\left(\frac{n}{2}\right)}{2^{n-1}\pi^{\frac{n}{2}}\Gamma\left(\frac{n-1}{2}\right)}G_2(x',t) \,\text{ for }x'\notin\supp\psi,\ 0<t<1.\label{0517-11}
\end{align}

In particular, for \(G_2(x',t)\) the following representation formula will be useful.
\begin{lemma}\label{lem-0517-01}
	Let \(w\) solve the Stokes system \eqref{0308-01} with boundary data \eqref{0324-13}, where \(j=n\).
	Suppose that \(\psi\) satisfies {\rm\textbf{(S)}}, and that \(\phi\) satisfies {\rm\textbf{(T1)}}.
	Then, for \(|x'|>1\) and \(0<t<1\), \(G_2(x',t)\) in \eqref{0508-01} satisfies
	\begin{align}\label{0628-02}
		G_2(x',t) = \int_0^t \frac{\phi'(s)}{\sqrt{t-s}}W_2(x',t-s)\,ds,
	\end{align}
	where
	\begin{align}\label{0512-04}
		W_2(x',t) & := \frac{1}{(n-1)t^{\frac{n}{2}}}\int_{\Sigma}\psi(y')\mathcal{K}\left( \frac{|x'-y'|^2}{4t}\right)\,dy',
	\end{align}
	with \(\mathcal{K}(z)\) given in \eqref{0827-12}.
\end{lemma}
\begin{proof}
	For \(x'\neq0\), \(t>0\), and positive integers \(m\) and \(k\) such that \(m-k\notin\{0,2,4,\ldots\}\), set
	\begin{align*}
		W_{m,k}(x',t) := \left(\frac{|x'|^2}{t}\right)^{\frac{m}{2}}{}_1F_1\left(\frac{m}{2}, \frac{k}{2}; -\frac{|x'|^2}{4t}\right).
	\end{align*}
Then \eqref{0508-01} gives
	\begin{align*}
		G_2(x',t) & = \int_0^t\frac{\phi'(s)}{\sqrt{t-s}}\int_\Sigma \frac{\psi(y')}{|x'-y'|^n}W_{n,n-1}(x'-y',t-s)\,dy'ds\\
		& \quad + \frac{n}{2}\int_0^t \frac{\phi(s)}{\sqrt{t-s}}\int_{\Sigma}\frac{\psi(y')}{|x'-y'|^{n+2}}W_{n+2, n-1}(x'-y',t-s)\,dy'ds.
	\end{align*}
	Set
	\begin{align*}
		F(x',t):= -\int_0^t s^{-\frac12}W_{n+2, n-1}(x',s)\,ds.
	\end{align*}
	For \(x'\neq0\), the limit \(\lim_{t\rightarrow 0^+}W_{m,k}(x',t) = \frac{2^m\Gamma\left(\frac{k}{2}\right)}{\Gamma\left(\frac{k-m}{2}\right)}\) holds by \eqref{0327-10}.
	Hence \(F(x',t)\to0\) as \(t\to0^+\), and integration by parts using \(\phi(0)=0\) gives
	\begin{align}\label{0512-09}
		G_2(x',t) = \int_0^t \phi'(s)\int_{\Sigma}\frac{\psi(y')}{|x'-y'|^n}\left( \frac{1}{\sqrt{t-s}}W_{n, n-1}(x'-y',t-s) - \frac{n}{2}\frac{1}{|x'-y'|^2}F(x'-y',t-s)\right)\,dy'ds.
	\end{align}

	It remains to find \(F(x',t)\).
	We first note that the relations \eqref{0327-07} and \eqref{0506-04} respectively give
	\begin{align}
		\partial_t\left(\sqrt{t} W_{m,k}(x',t)\right) = -\frac{m-1}{2\sqrt{t}}W_{m,k}(x',t) + \frac{m}{4k\sqrt{t}}W_{m+2,k+2}(x',t),\label{0512-05} \\
		W_{m,k}(x',t) = W_{m,k+2}(x',t) - \frac{m}{2k(k+2)}W_{m+2, k+4}(x',t).\label{0512-06}
	\end{align}
	Combining \eqref{0512-06} at \((m,k)=(n+2,n-1)\) and \((n+2,n+1)\) with \eqref{0512-05} at \((m,k)=(n+2,n+1)\) and \((n+2,n+3)\) gives
	\begin{align*}
		\frac{1}{\sqrt{t}}W_{n+2, n-1}(x',t) = -\frac{2}{n-1}\partial_t\left(\sqrt{t}W_{n+2, n+1}(x',t)\right) + \frac{4}{n^2-1}\partial_t\left(\sqrt{t}W_{n+2, n+3}(x',t)\right).
	\end{align*}
	Integrating in \(t\) gives
	\begin{align*}
		F(x',t) = \frac{2\sqrt{t}}{n-1}W_{n+2, n+1}(x',t) - \frac{4\sqrt{t}}{n^2-1}W_{n+2, n+3}(x',t).
	\end{align*}
	Then \eqref{0512-09} implies
	\begin{align}
		G_2(x',t) & = \int_0^t \frac{\phi'(s)}{\sqrt{t-s}}\int_{\Sigma}\frac{\psi(y')}{|x'-y'|^n}\left[W_{n,n-1}(x'-y',t-s)\right.\nonumber\\
		& \quad \left. + \frac{t-s}{|x'-y'|^2}\left( \frac{2n}{n^2-1}W_{n+2,n+3}(x'-y',t-s) - \frac{n}{n-1}W_{n+2, n+1}(x'-y',t-s)\right)\right]\,dy'ds\nonumber\\
		& = \frac{1}{n-1} \int_0^t \frac{\phi'(s)}{\sqrt{t-s}}\nonumber\\
		& \quad \times\int_{\Sigma}\frac{\psi(y')}{(t-s)^{\frac{n}{2}}}\left[{}_1F_1\left(\frac{n}{2}, \frac{n+1}{2}; -\frac{|x'-y'|^2}{4(t-s)}\right) -\frac{2}{n+1}{}_1F_1\left(\frac{n}{2}, \frac{n+3}{2} ; -\frac{|x'-y'|^2}{4(t-s)}\right)\right]\,dy'ds,\label{0514-01}
	\end{align}
		where we have used the contiguous relation \eqref{0327-02} twice in the last equality.
\end{proof}

\subsubsection{Nonexistence of separation points}
The next proposition proves the nonexistence statement in Theorem \ref{thm-0507-01}{\rm(i)}. 
\begin{proposition}\label{prop-0730-20}
	Let \(w\) solve the Stokes system \eqref{0308-01} with boundary data \eqref{0324-13}, where \(j=n\).
	Suppose that \(\psi\) satisfies {\rm\textbf{(S)}}, and that \(\phi\) satisfies {\rm\textbf{(T1)}}.
	Fix \(0<t<1\) such that \(\phi\) is nondecreasing on \([0,t]\) and \(\phi(t)>0\).
	Then \(\partial_{x_n}^2w_n(x',0,t)>0\) whenever \(|x'|\geq1+2\sqrt{z_3t}\).
\end{proposition}
\begin{proof}
		For \(y'\in\supp\psi\) and \(0<s<t\), \(\frac{|x'-y'|^2}{4(t-s)}\geq\frac{(|x'|-1)^2}{4(t-s)}\geq z_3\frac{t}{t-s}>z_3\).
	Lemma \ref{lem-0514-01} therefore gives \(\mathcal{K}(\frac{|x'-y'|^2}{4(t-s)})<0\), so \(W_2(x',t-s)<0\), where \(W_2\) is defined in \eqref{0512-04}.
	Moreover, \(\phi'(s)\geq0\) on \((0,t)\) and \(\int_0^t\phi'(s)\,ds=\phi(t)>0\).
	Then \eqref{0628-02} and \eqref{0517-11} yield \(G_2(x',t)<0\) and \(\partial_{x_n}^2w_n(x',0,t)>0\).
\end{proof}

\subsubsection{Existence of separation points}
The next proposition proves the existence statement in Theorem \ref{thm-0507-01}{\rm(i)}.
\begin{proposition}\label{prop-0504-02}
	Let \(w\) solve the Stokes system \eqref{0308-01} with boundary data \eqref{0324-13}, where \(j=n\).
	Suppose that \(\psi\) satisfies {\rm\textbf{(S)}}, and that \(\phi\) satisfies {\rm\textbf{(T1)}} and \(\phi\not\equiv0\).
	Then, for \(|x'|>1+2\sqrt{z_1}\), \(\partial_{x_n}^2w_n(x',0,\cdot)\) has a sign-changing zero \(\widetilde{t_n^*}(x')\) with
	\begin{align*}
	    \frac{(|x'|-1)^2}{4z_1}<\widetilde{t_n^*}(x')<1+\frac{(|x'|+1)^2}{4z_1}.
	\end{align*}
	The point \((x',\widetilde{t_n^*}(x'))\) is a separation point of \(w_n\).
	If \(\phi\) also satisfies {\rm\textbf{(T2)}}, then there is a separation point \((x',t_n^*(x'))\) of \(w_n\) with \(0<t_n^*(x')<1\) for each \(|x'|>1+2\sqrt{z_1}\).
	For all sufficiently large \(|x'|\), depending only on \(n\), \(\widetilde{t_n^*}(x')\) is the unique zero of \(\partial_{x_n}^2w_n(x',0,\cdot)\) in \((1,\infty)\), and \(\partial_t\partial_{x_n}^2w_n(x',0,\widetilde{t_n^*}(x'))>0\).
\end{proposition}
\begin{proof}
		We first prove the existence of a sign-changing zero \(\widetilde{t_n^*}(x')>1\) of \(\partial_{x_n}^2w_n(x',0,t)\).
	Fix \(x'\in \Sigma\) such that \(|x'|>1+2\sqrt{z_1}\), and set \(t_+=1+(|x'|+1)^2/(4z_1)\).
	For \(0<s<1\), \(\frac{|x'-y'|^2}{4(t_+-s)}<\frac{(|x'|+1)^2}{4(t_+-1)}=z_1\).
	It follows that \(G_1(x',t_+)>0\).

	On the other hand, if \(1<t\leq \frac{(|x'|-1)^2}{4z_1}\), then, for \(0<s<1\), \( \frac{|x'-y'|^2}{4(t-s)}>\frac{(|x'|-1)^2}{4t}\geq z_1. \)
Hence \({}_1F_1\left(\frac{n}{2}, \frac{n-1}{2}; -\frac{|x'-y'|^2}{4(t-s)}\right)<0\) and \(G_1(x',t)<0\).

	By \eqref{0504-10} and Lemma \ref{lem-260044}{\rm(v)}, applied with \(T=1\), \(g=\psi\), \(q=\phi\), \(a=\frac n2\), \(b=\frac{n-1}{2}\), and \(c=\frac{n+3}{2}\), \(G_1(x',t)\) is real-analytic in \(t>1\).
	Together with Lemma \ref{lem-260044}{\rm(iii)}, this shows that \(G_1(x',t)\) has only finitely many zeros between \((|x'|-1)^2/(4z_1)\) and \(t_+\).
			Since its signs at these times are opposite, at least one of these zeros, \(\widetilde{t_n^*}(x')\), changes sign, and so does \(\partial_{x_n}^2w_n(x',0,\cdot)\) by \eqref{0517-10}.
	Proposition \ref{prop-260046}{\rm(i)} therefore shows that \((x',\widetilde{t_n^*}(x'))\) is a separation point of \(w_n\).

		To prove simplicity and uniqueness for large \(r=|x'|\), write \(H(u)={}_1F_1(\frac n2,\frac{n-1}{2};-u)\) and \(b=\frac{n+3}{2}\).
	The preceding argument gives \(G_1(x',t)<0\) for \(1<t\leq (r-1)^2/(4z_1)\).
	For \(t\geq 1+(r+1)^2/(4z_1)\), every argument of \(H\) lies in \((0,z_1)\), so \(G_1(x',t)>0\).
	Hence every zero after time \(1\) lies strictly between these two endpoints.

	Throughout the closed interval between them, \(0<s<1\) and \(y'\in\supp\psi\) imply
	\begin{align*}
	    \frac{z_1(r-1)^2}{(r+1)^2+4z_1}
	    \leq u:=\frac{|x'-y'|^2}{4(t-s)}
	    \leq \frac{z_1(r+1)^2}{(r-1)^2-4z_1}.
	\end{align*}
	For sufficiently large \(r\), both denominators are positive and both bounds tend to \(z_1\) as \(r\to\infty\).
		By \eqref{0327-07}, the definition of \(z_2\), and Lemma \ref{lem-0730-50}, \(H'(z_1)=-\frac n{n-1}\,{}_1F_1(\frac{n+2}{2},\frac{n+1}{2};-z_1)<0\).
	Hence \(H(z_1)=0\) implies that \(-bH(u)-uH'(u)\) is bounded above and below by positive constants depending only on \(n\) for every \(u\) in this range once \(r\) is sufficiently large.

The time derivative of the kernel in \eqref{0504-10} equals \((t-s)^{-b-1}[-bH(u)-uH'(u)]\).
	On the same interval, \(t-s\) is comparable to \(r^2\) uniformly in \(s\).
	Hence \eqref{0517-10} gives
	\begin{align}\label{late-normal-slope}
	    c_n M_0\|\phi\|_{L^1(0,1)}r^{-n-5}
	    \leq\partial_t\partial_{x_n}^2w_n(x',0,t)
	    \leq C_n M_0\|\phi\|_{L^1(0,1)}r^{-n-5}
	\end{align}
	for \((r-1)^2/(4z_1)\leq t\leq 1+(r+1)^2/(4z_1)\), with \(0<c_n\leq C_n\) depending only on \(n\).
	Strict increase on the entire interval containing the zeros proves uniqueness and simplicity.

	Assume now {\rm\textbf{(T2)}} to find \(t_n^*(x')\in(0,1)\).
	By \eqref{0514-01} and Lemma \ref{lem-260044}{\rm(iv)}, applied with \(T=1\), \(g=\psi\), \(q=\phi'\), \(m=0\), \(a=\frac n2\), \(c=\frac{n+1}{2}\), and \(b=\frac{n+1}{2},\frac{n+3}{2}\), \(G_2\) is jointly continuous for \(|x'|>1\), \(0<t\leq1\).
		For each fixed \(\omega'\in\mathbb S^{n-2}\) and \(t\in(0,1]\), \(G_2(r\omega',t)\) is real-analytic in \(r>1\).

	Consider the time \(t_2=\frac12\min\{t_1,z_1/z_3\}\).
	For \(|x'|>1+2\sqrt{z_1}\), \(y'\in\supp\psi\), and \(0<s<t\leq t_2\), \(\frac{|x'-y'|^2}{4(t-s)} >\frac{z_1}{t_2}\geq2z_3\).
	Lemma \ref{lem-0514-01}{\rm(i)} and \eqref{0512-04} give \(W_2(x',t-s)<0\).
	Since \(\phi'(s)>0\), \eqref{0628-02} gives \(G_2(x',t)<0\) for \(0<t\leq t_2\).

	We now prove \(G_2(x',1)>0\).
	Since \(\phi(0)=\phi(1)=0\), integration by parts gives
	\begin{align*}
		G_2(x',1) & = \int_{\Sigma}\psi(y')\int_0^1 \phi(s) \left[ \frac{n}{2(1-s)^{\frac{n+3}{2}}}{}_1F_1\left(\frac{n+2}{2}, \frac{n-1}{2};- \frac{|x'-y'|^2}{4(1-s)}\right) \right.\\
		& \quad \left. -\partial_s\left(\frac{1}{(1-s)^{\frac{n+1}{2}}} {}_1F_1\left(\frac{n}{2}, \frac{n-1}{2}; -\frac{|x'-y'|^2}{4(1-s)}\right)\right)\right]\,dsdy'\\
		& \quad +\int_{\Sigma}\psi(y')\lim_{s\rightarrow 1^{-}} \frac{\phi(s)}{(1-s)^{\frac{n+1}{2}}}{}_1F_1\left(\frac{n}{2}, \frac{n-1}{2};- \frac{|x'-y'|^2}{4(1-s)}\right)\,dy'\\
		& = \int_{\Sigma}\psi(y')\int_0^1 \frac{\phi(s)}{(1-s)^{\frac{n+3}{2}}} \left[\frac{n}{2}{}_1F_1\left(\frac{n+2}{2},\frac{n-1}{2};z\right)\right.\\
		& \quad\left.-\frac{n+1}{2}{}_1F_1\left(\frac{n}{2},\frac{n-1}{2};z\right) -z \frac{d}{dz}{}_1F_1\left(\frac{n}{2},\frac{n-1}{2};z\right)\right]\,dsdy'\\
		& =- \frac12 \int_\Sigma \psi(y')\int_0^1 \frac{\phi(s)}{(1-s)^{\frac{n+3}{2}}}{}_1F_1\left(\frac{n}{2}, \frac{n-1}{2}; -\frac{|x'-y'|^2}{4(1-s)}\right)\,dsdy'.
	\end{align*}
	Here \(z=-|x'-y'|^2/(4(1-s))\), and we used the contiguous relation \eqref{0327-01}.
	The boundary term containing the limit as \(s\to1^-\) vanishes since
	\begin{align*}
		\left|\frac{\phi(s)}{(1-s)^{\frac{n+1}{2}}}{}_1F_1\left(\frac{n}{2}, \frac{n-1}{2}; - \frac{|x'-y'|^2}{4(1-s)}\right)\right| & \leq \frac{C\norm{\phi'}_{L^\infty(0,1)}}{(1-s)^{\frac{n-1}{2}}}\left(\frac{1-s}{|x'-y'|^2}\right)^{\frac{n}{2}} \\
		& \leq \frac{C\norm{\phi'}_{L^\infty(0,1)}(1-s)^{\frac12}}{|x'-y'|^n}.
	\end{align*}

	If \(|x'|>1+2\sqrt{z_1}\), then \(\frac{|x'-y'|^2}{4(1-s)}\geq\frac{(|x'|-1)^2}{4}>z_1\), so \({}_1F_1\left(\frac{n}{2}, \frac{n-1}{2};-\frac{|x'-y'|^2}{4(1-s)}\right)<0\).
	Hence \(G_2(x',1)>0\).

	For each \(0<t\leq1\), \(G_2(r\omega',t)\) is not identically zero in \(r>1\).
	Indeed, substituting \eqref{0504-06} into \eqref{0514-01} gives, as \(r\to\infty\),
	\begin{align}\label{1006-07}
		G_2(r\omega',t) & =-\frac{2^{n-1}\Gamma\left(\frac{n-1}{2}\right)}{\sqrt\pi} \left[T_0(t)\int_\Sigma\frac{\psi(y')}{|r\omega'-y'|^n}\,dy' -3nT_1(t)\int_\Sigma\frac{\psi(y')}{|r\omega'-y'|^{n+2}}\,dy' +O(r^{-n-4})\right].
	\end{align}
	The absolute value of the last term is bounded by \(Cr^{-n-4}\int_0^t|\phi'(t-s)|s^{\frac32}\,ds\).
	Since \(T_1(t)=\frac12\int_0^t\phi(s)(t-s)^{-\frac12}\,ds>0\), \eqref{1006-07} gives
	\begin{align*}
		\lim_{r\to\infty}r^nG_2(r\omega',t) & =-\frac{2^{n-1}\Gamma\left(\frac{n-1}{2}\right)}{\sqrt\pi} M_0T_0(t)\neq0 & & \text{if }T_0(t)\neq0, \\
		\lim_{r\to\infty}r^{n+2}G_2(r\omega',t) & =\frac{3n2^{n-1}\Gamma\left(\frac{n-1}{2}\right)}{\sqrt\pi} M_0T_1(t)>0& & \text{if }T_0(t)=0.
	\end{align*}
	Thus Lemma \ref{lem-260044}{\rm(iii)} excludes an interval of zeros in \(r>1\) at any fixed \(0<t\leq1\).

	Take \(\omega'=x'/|x'|\).
	The connectedness argument in the proof of Proposition \ref{prop-0504-03}, applied to \(G_2(r\omega',t)\) on \([t_2,1]\), gives \(t_n^*(x')\in(t_2,1)\) such that \(G_2\) takes both signs in every neighborhood of \((|x'|,t_n^*(x'))\).
	By continuity, this gives \(G_2(x',t_n^*(x'))=0\).

	By \eqref{0517-11}, \(\partial_{x_n}^2w_n\) then takes both signs in every neighborhood of \((x',t_n^*(x'))\).
	Proposition \ref{prop-260046}{\rm(i)} shows that \((x',t_n^*(x'))\) is a separation point of \(w_n\).
\end{proof}

\subsubsection{Temporal vanishing rates of \(\partial_{x_n}^2w_n\) near a zero}
The next proposition estimates the temporal vanishing rates of \(\partial_{x_n}^2w_n\) near the two zeros found in Proposition \ref{prop-0504-02}.
For \(0<\beta\leq1/2\), it also constructs temporal boundary data attaining the stated rates for either \(\partial_{x_n}w_i(x',0,\cdot)\), \(i<n\), or \(\partial_{x_n}^2w_n(x',0,\cdot)\).
This proves Theorem \ref{thm-0507-01}{\rm(ii)}.
\begin{proposition}\label{prop-0512-01}
	Let \(w\) solve the Stokes system \eqref{0308-01} with boundary data \eqref{0324-13}, where \(j=n\).
	Suppose that \(\psi\) satisfies {\rm\textbf{(S)}}, and that \(\phi\) satisfies {\rm\textbf{(T1)}} and \(\phi\not\equiv0\).
	Fix \(|x'|>1+2\sqrt{z_1}\), and let \(\widetilde{t_n^*}(x')\in\bigl((|x'|-1)^2/(4z_1),\,1+(|x'|+1)^2/(4z_1)\bigr)\) be a sign-changing zero of \(\partial_{x_n}^2w_n(x',0,\cdot)\).
	Then the following hold.
	\begin{enumerate}
		\item\label{item-0518-01} For \(t\) sufficiently close to \(\widetilde{t_n^*}(x')\), with \(\kappa=1\) if \(\partial_t\partial_{x_n}^2 w_n(x',0,\widetilde{t_n^*}(x'))\neq0\) and \(\kappa=3\) otherwise,
		\begin{align}\label{0919-02}
			|\partial_{x_n}^2w_n(x',0,t)|\leq C(n)M_0\norm{\phi}_{L^1(0,1)}
			\begin{cases}
				\displaystyle\frac{|t-\widetilde{t_n^*}(x')|^\kappa}{\bigl((|x'|-1)^2-4z_1\bigr)^{\frac{n+3+2\kappa}{2}} }, & 1+2\sqrt{z_1}<|x'|<\max\{2,1+2\sqrt{2z_1}\}, \\
				\displaystyle\frac{|t-\widetilde{t_n^*}(x')|^\kappa}{|x'|^{n+3+2\kappa}}, & |x'|\geq\max\{2,1+2\sqrt{2z_1}\}.
			\end{cases}
		\end{align}

		For sufficiently large \(|x'|\), \(\kappa=1\) and
		\begin{align}\label{0919-02a}
		    |\partial_{x_n}^2w_n(x',0,t)|\asymp_n M_0\|\phi\|_{L^1(0,1)}|x'|^{-n-5}|t-\widetilde{t_n^*}(x')|.
		\end{align}

		\item\label{item-0518-02} Suppose that \(\phi\) satisfies {\rm\textbf{(T2)}}, let \((x',t_n^*(x'))\) be a separation point of \(w_n\) with \(t_n^*(x')\in(0,1)\), and let \(\beta\in(0,1)\) be the H\"older exponent in {\rm\textbf{(T1)}}.
		If \(t\in(0,1)\) and \(h=|t-t_n^*(x')|\) satisfies
		\begin{align*}
			0<h<\min\left\{t_n^*(x'),1-t_n^*(x')\right\},
		\end{align*}
		then
		\begin{align}\label{0919-03}
			|\partial_{x_n}^2w_n(x',0,t)|\leq C\frac{|x'|^2}{(|x'|-1)^{n+2}}
			\begin{cases}
				h^{\beta+\frac12}, & 0<\beta<\frac12, \\
				h|\ln h|,& \beta=\frac12,\\
				h^\sigma,& \frac12<\beta<1,
			\end{cases}
		\end{align}
		where, for \(\beta>1/2\), \(\sigma=1\) if \(\partial_t\partial_{x_n}^2w_n(x',0,t_n^*(x'))\neq0\) and \(\sigma=\beta+\frac12\) otherwise, and \(C>0\) depends only on \(n\), \(\beta\), \(M_0\), \(t_1\), and \(\norm{\phi'}_{C^\beta}\).

		\item\label{item-0801-20} Under the hypotheses of {\rm(ii)}, suppose that \(0<\beta<\frac12\) and that, for some \(c_0>0\),
		\begin{align*}
			\phi'(s)-\phi'(t_n^*(x'))=-c_0(s-t_n^*(x'))|s-t_n^*(x')|^{\beta-1} +o(|s-t_n^*(x')|^\beta) \quad\text{as }s\to t_n^*(x').
		\end{align*}
		Then
		\begin{align}\label{0731-20}
			\lim_{\ep\rightarrow0^+}\ep^{-\beta-\frac12}\,\partial_{x_n}^2w_n(x',0,t_n^*(x')+\ep) =-\frac{\Gamma\left(\frac n2\right)}{\pi^{\frac{n+1}{2}}}\,c_0 \frac{\Gamma(\beta+1)\Gamma(\frac12)}{\Gamma(\beta+\frac32)} \bigl(1+\sec(\pi\beta)\bigr) \int_{\Sigma}\frac{\psi(y')}{|x'-y'|^{n}}\,dy',
		\end{align}
		and
		\begin{align}\label{0731-20a}
			\lim_{\ep\rightarrow0^-}|\ep|^{-\beta-\frac12}\,\partial_{x_n}^2w_n(x',0,t_n^*(x')+\ep) =\frac{\Gamma\left(\frac n2\right)}{\pi^{\frac{n+1}{2}}}\,c_0 \frac{\Gamma(\beta+1)\Gamma(\frac12)}{\Gamma(\beta+\frac32)} \tan(\pi\beta) \int_{\Sigma}\frac{\psi(y')}{|x'-y'|^{n}}\,dy'.
		\end{align}
		Thus the exponent \(\beta+\frac12\) in {\rm(ii)} is sharp, and \(\partial_{x_n}^2w_n(x',0,t)\) is not differentiable at \(t_n^*(x')\).
	\end{enumerate}
For every \(0<\beta\leq1/2\) and each \(\psi\) satisfying {\rm\textbf{(S)}}, there exists \(\phi\) satisfying {\rm\textbf{(T1)}} and {\rm\textbf{(T2)}} for which the rate in {\rm(ii)} is attained at a separation time in \((t_1,1)\) for some sufficiently large \(|x'|\).
The same holds for \(\partial_{x_n}w_i(x',0,\cdot)\), \(1\leq i\leq n-1\), and the rate in \eqref{0919-01} at any fixed \(x'\) with \(|x_i|>1\), with a possibly different \(\phi\).
\end{proposition}
\begin{proof}
	\underline{Proof of (i).}\quad
	Fix \(x'\) as in the statement and write \(\tau=\widetilde{t_n^*}(x')\).
	By Proposition \ref{prop-0504-02}, \(G_1(x',\cdot)\) is real-analytic near \(\tau\) and changes sign there.
	If its first derivative vanishes, its second derivative must also vanish, since a nonzero quadratic term would give the same sign on both sides.
	Taylor's theorem therefore gives, for some \(\xi\) between \(t\) and \(\tau\),
	\begin{align}\label{0517-12}
	    G_1(x',t)=\frac{(t-\tau)^\kappa}{\kappa!}\partial_t^\kappa G_1(x',\xi),\qquad \kappa\in\{1,3\}.
	\end{align}

		For \(|t-\tau|\leq(\tau-1)/2\) and \(0<s<1\), Proposition \ref{prop-0504-02} gives \(t-s\geq\frac{\tau-1}{2}>\frac{(|x'|-1)^2-4z_1}{8z_1}\).
	For \(|x'|\geq\max\{2,1+2\sqrt{2z_1}\}\), the last bound is at least \(|x'|^2/(64z_1)\).

	With \(H\) and \(b\) as in the preceding proof, the \(k\)th time derivative of \((t-s)^{-b}H(u)\) is \((t-s)^{-b-k}\) times a linear combination of \(u^jH^{(j)}(u)\), \(0\leq j\leq k\).
	Equations \eqref{0327-07} and \eqref{0514-02} bound these factors uniformly for \(u>0\) and \(k\leq3\).
	Using the lower bounds on \(t-s\) in \eqref{0504-10}, and then \eqref{0517-12} and \eqref{0517-10}, gives both cases of \eqref{0919-02}.
	For sufficiently large \(|x'|\), integrating \eqref{late-normal-slope} between \(t\) and \(\tau\) gives \eqref{0919-02a}.

	\underline{Proof of (ii).}\quad
We observe that \(W_2(x',t)\) in \eqref{0512-04} satisfies
	\begin{align}\label{0624-10}
		\norm{\partial_t^k W_2(x',\cdot)}_{L_t^\infty(0,1)}\leq \frac{CM_0}{(|x'|-1)^{n+2k}},\quad 0\leq k\leq 3.
	\end{align}
		Indeed, \eqref{0327-01} and \eqref{0514-02} give \eqref{0624-10}, and together with \eqref{0504-06}, \(W_2\in C^2([0,1])\).
	The estimates for \(t>t_i^*\) in the proof of Proposition~\ref{prop-0512-02}{\rm(i)} apply to \(G_2,W_2\) by \eqref{0624-10}, with \(|x'|^2/(|x'|-1)^{n+1}\) replaced by \(|x'|^2/(|x'|-1)^{n+2}\).
	For \(t<t_n^*(x')\) when \(\partial_t\partial_{x_n}^2w_n(x',0,t_n^*(x'))=0\), the same dyadic difference-quotient estimate gives the power \(h^{\beta+1/2}\) in \eqref{0919-03}.
	Moreover, \(t_n^*(x')>t_2\) with \(t_2\) from the proof of Proposition \ref{prop-0504-02}, so \(t_1\) may be replaced by \(t_2\) in these estimates.

	\underline{Proof of (iii).}\quad
	Write \(t_n^*=t_n^*(x')\).
	The representation \eqref{0628-02} has the same form as \eqref{0629-01}, and \eqref{0624-10} gives the bounds used in the proofs of \eqref{0730-71} and \eqref{0730-71a}.
	Applying those calculations with \(W\) replaced by \(W_2\) gives
	\begin{align*}
		\lim_{\ep\to0^+}\ep^{-\beta-\frac12}G_2(x',t_n^*+\ep) & =-c_0W_2(x',0)\frac{\Gamma(\beta+1)\Gamma(\frac12)} {\Gamma(\beta+\frac32)}\bigl(1+\sec(\pi\beta)\bigr),\\
		\lim_{\ep\to0^-}|\ep|^{-\beta-\frac12}G_2(x',t_n^*+\ep) & =c_0W_2(x',0)\frac{\Gamma(\beta+1)\Gamma(\frac12)} {\Gamma(\beta+\frac32)}\tan(\pi\beta).
	\end{align*}
	It remains to evaluate \(W_2(x',0)\).
		By \eqref{0504-06} and \eqref{0514-02}, applying the dominated convergence theorem to the \(y'\)-integral in \eqref{0512-04} gives
	\begin{align*}
		W_2(x',0)=-\frac{2^n\Gamma(\frac{n+1}{2})}{(n-1)\Gamma(\frac12)} \int_\Sigma\frac{\psi(y')}{|x'-y'|^n}\,dy'<0.
	\end{align*}
	Substituting this into \eqref{0517-11} gives \eqref{0731-20} and \eqref{0731-20a}.
	If the same local expansion holds with \(\beta=1/2\), Remark~\ref{rmk-critical-abel} applies to \eqref{0628-02} with \(W_2\) in place of \(W\).
	Since \(W_2(x',0)<0\), the logarithmic upper bound in \eqref{0919-03} is attained on both sides of the earlier zero.

	To construct examples, fix \(0<\beta\leq1/2\) and take \(\phi_0(s)=s^2(1-s)^2\) on \([0,1]\), extended by zero for \(s>1\).
	Then \(\phi_0\) satisfies {\rm\textbf{(T1)}}--{\rm\textbf{(T4)}} with \(t_1=1/2\).
	Choose either \(\partial_{x_n}w_i\), \(i<n\), or \(\partial_{x_n}^2w_n\). For temporal boundary data \(f\), let \(V[f](u)\) denote the chosen derivative at \((x',0,u)\).

	For \(\partial_{x_n}w_i\), Theorem~\ref{thm-0324-01}{\rm(iii)} gives a sign-changing zero of \(V[\phi_0]\) in \((1/2,1)\) whenever \(|x_i|>1\).
	For \(\partial_{x_n}^2w_n\), Theorem~\ref{thm-0507-01}{\rm(iii)}--{\rm(iv)} gives the same conclusion for sufficiently large \(|x'|\), since the earlier zero converges to \(t_\infty\in(1/2,1)\).
	In either case choose \(1/2<a<b<1\) such that \(V[\phi_0](a)V[\phi_0](b)<0\).

	Choose \(\chi\in C_c^\infty((1/2,1))\) equal to one near \([a,b]\).
	For \(q\in[a,b]\) and \(\varepsilon>0\), set
	\begin{align*}
		\phi_q(s)=\phi_0(s)-\frac{\varepsilon}{1+\beta}\chi(s)|s-q|^{1+\beta}.
	\end{align*}
	The \(C^1\) norm of the perturbation is at most \(C\varepsilon\), uniformly in \(q\), and \(\max_{\supp\chi}\phi_0'<0\).
	Thus, for sufficiently small \(\varepsilon\), all \(\phi_q\) satisfy {\rm\textbf{(T1)}} and {\rm\textbf{(T2)}} with \(t_1=1/2\).

	The map \(q\mapsto\phi_q\) is continuous in \(C^1\).
	By \eqref{0629-01}, \eqref{0630-07}, \eqref{0517-11}, and \eqref{0628-02}, in either case \(V[f]\) is a constant multiple of a convolution of \(f'\) with a kernel bounded by \(C_{x'}s^{-1/2}\) on \(0<s<b\).
		Consequently \(V[f]\) is continuous on \([a,b]\) and
\begin{align*}
\sup_{a\leq u\leq b}|V[f](u)-V[g](u)|\leq 2C_{x'}\sqrt b\,\|f'-g'\|_{L^\infty(0,b)}.
\end{align*}
	Hence \(q\mapsto V[\phi_q](q)\) is continuous and, for small \(\varepsilon\), its values at \(a\) and \(b\) retain their opposite signs.
	Choose \(q_*\in(a,b)\) at which it vanishes.

	Since \(\chi=1\) near \(q_*\) and \(\phi_0\) is smooth there,
	\begin{align*}
		\phi_{q_*}'(s)-\phi_{q_*}'(q_*)
		=-\varepsilon\,\sgn(s-q_*)|s-q_*|^\beta+O(|s-q_*|).
	\end{align*}
	The one-sided limits above and in Proposition~\ref{prop-0512-02} and Remark~\ref{rmk-critical-abel} use only this expansion and \(V[\phi_{q_*}](q_*)=0\).
	They give nonzero coefficients of opposite signs on the two sides, so Proposition~\ref{prop-260046}{\rm(i)} makes \((x',q_*)\) a separation point.
\end{proof}

\subsubsection{Uniqueness of the zeros}
The next proposition shows that both zeros are unique for sufficiently large \(|x'|\).
With Proposition \ref{prop-0504-02}, it proves Theorem \ref{thm-0507-01}{\rm(iii)}.
\begin{proposition}\label{prop-0505-01}
	Let \(w\) solve the Stokes system \eqref{0308-01} with boundary data \eqref{0324-13}, where \(j=n\).
	Suppose that \(\psi\) satisfies {\rm\textbf{(S)}}, and that \(\phi\) satisfies {\rm\textbf{(T1)}} and {\rm\textbf{(T2)}} and is log-concave on \((0,1)\).
	Then, for each \(|x'|>R_+\), \(\widetilde{t_n^*}(x')\) is the unique zero of \(\partial_{x_n}^2w_n(x',0,\cdot)\) in \((1,\infty)\), and, for each \(|x'|>R_-\), \(t_n^*(x')\) is the unique zero in \((0,1)\), where
	\begin{align}\label{0730-40}
		R_+ & :=\frac{\sqrt{z_2}+\sqrt{z_1}}{\sqrt{z_2}-\sqrt{z_1}},\qquad
		R_-:=1+2\sqrt{\max\{z_5,z_6\}}.
	\end{align}
	If \(|x'|>\max\{R_+,R_-\}\), then
	\begin{align}\label{0730-41}
		\sgn\bigl(\partial_{x_n}^2w_n(x',0,t)\bigr) =\begin{cases} 1,& 0<t<t_n^*(x'),\\
			-1, & t_n^*(x')<t<\widetilde{t_n^*}(x'), \\
			1,& t>\widetilde{t_n^*}(x').
		\end{cases}
	\end{align}
	If \(\phi\) also satisfies {\rm\textbf{(T3)}}, then \(\partial_t\partial_{x_n}^2w_n(x',0,t_n^*(x'))<0\) for \(|x'|>R_-\).
\end{proposition}
\begin{proof}
	By \eqref{0701-05}, \(R_+=1+2\sqrt{z_1}/(\sqrt{z_2}-\sqrt{z_1})>1+2\sqrt{z_1}\), and \eqref{0512-03} with \(a=\frac12\) and \(b=\frac{n-1}{2}\) gives \(z_1<z_5\), so \(R_->1+2\sqrt{z_1}\).
	Hence Proposition~\ref{prop-0504-02} applies whenever \(|x'|>R_+\) or \(|x'|>R_-\).

	We first prove the uniqueness of \(\widetilde{t_n^*}(x')\).
	For \(t>1\), recall \(G_1(x',t)\) from \eqref{0504-10}.
	Write \(G_1(x',t) = \int_0^1 \phi(s)K_0(x',t-s)\,ds\), where
	\begin{align*}
		K_0(x',t) := \frac{1}{t^{\frac{n+3}{2}}}\int_\Sigma \psi(y') {}_1F_1\left(\frac{n}{2}, \frac{n-1}{2}; -\frac{|x'-y'|^2}{4t}\right)\,dy'.
	\end{align*}
	We show that \(K_0(x',t)\) has exactly one zero in \(t\) on \((0,\infty)\).
	Lemma~\ref{lem-0514-01}{\rm(ii)} and the definition of \(z_1\) give \(K_0(x',t)<0\) for small \(t\) and \(K_0(x',t)>0\) for large \(t\).
			If \(K_0(x',t_*)=0\), then \(\frac{(|x'|-1)^2}{4t_*}<z_1\).
	Otherwise, \(\frac{|x'-y'|^2}{4t_*}\geq z_1\) for \(|y'|\leq1\).
	The function \({}_1F_1\left(\frac{n}{2}, \frac{n-1}{2}; -\frac{|x'-y'|^2}{4t_*}\right)\) is then nonpositive for \(|y'|\leq1\) and negative on a set of positive measure where \(\psi>0\).
	Hence \(K_0(x',t_*)<0\), a contradiction.

	Using \eqref{0327-07} and \(K_0(x',t_*)=0\) gives
	\begin{align*}
		\partial_t K_0(x',t_*)=\frac{n}{4(n-1)t_*^{\frac{n+7}{2}}}\int_{\Sigma}\psi(y')|x'-y'|^2{}_1F_1\left( \frac{n+2}{2}, \frac{n+1}{2}; -\frac{|x'-y'|^2}{4t_*}\right)\,dy'.
	\end{align*}
	Set \(r=|x'|\).
			Since \(r=|x'|>R_+\), \(\frac{(r+1)^2}{4t_*}= \frac{(r+1)^2}{(r-1)^2}\frac{(r-1)^2}{4t_*}<z_2\).
	Thus \({}_1F_1\left( \frac{n+2}{2}, \frac{n+1}{2}; -\frac{|x'-y'|^2}{4t_*}\right)>0\) for all \(|y'|\leq 1\), so \(\partial_t K_0(x', t_*)>0\).
		Hence \(K_0(x',\cdot)\) satisfies the hypothesis of Lemma \ref{lem-0504-01}.
		Since \(\phi(0)=\phi(1)=0\) and \(\phi\) is positive and log-concave on \((0,1)\), Lemma~\ref{lem-0504-01} implies that \(G_1(x',\cdot)\) has at most one zero in \((1,\infty)\).
	Together with the existence result in Proposition \ref{prop-0504-02}, this proves uniqueness of \(\widetilde{t_n^*}(x')\).

	We next prove the uniqueness of \(t_n^*(x')\).
	For \(0<t<1\), set
	\begin{align}
		A_1(t) & := - \frac{1}{t^{\frac{n+1}{2}}}\int_{\Sigma}\psi(y'){}_1F_1\left(\frac{n}{2}, \frac{n-1}{2}; -\frac{|x'-y'|^2}{4t}\right)\,dy',\label{0701-01} \\
		A_2(t) & := \frac{n}{2t^{\frac{n+3}{2}}}\int_\Sigma \psi(y'){}_1F_1\left(\frac{n+2}{2}, \frac{n-1}{2}; -\frac{|x'-y'|^2}{4t}\right)\,dy'.\label{0701-02}
	\end{align}
	Then \eqref{0508-01} gives
	\begin{align*}
		G_2(x',t)=\int_0^t (A_2(t-s)\phi(s)- A_1(t-s)\phi'(s))\,ds.
	\end{align*}

For \(z>z_5\), \({}_1F_1(\frac{n+2}{2},\frac{n-1}{2};-z)>0\).
Since \(|x'|>R_-\geq1+2\sqrt{z_5}\), \eqref{0701-02} gives \(A_2(t)>0\) for \(0<t<1\).
	To prove \(A_1'(t)<0\), set \(F(z)=z^{\frac{n+1}{2}}{}_1F_1(\frac n2,\frac{n-1}{2};-z)\).
	Lemma \ref{lem-260032} gives \(F'(z)<0\) for \(z>z_6\).
	Since
	\begin{align}\label{0506-06}
		A_1'(t) = \frac{2^{n-1}}{t^2}\int_{\Sigma}\frac{\psi(y')}{|x'-y'|^{n-1}}F'\left( \frac{|x'-y'|^2}{4t}\right)\,dy',
	\end{align}
	and \(\frac{|x'-y'|^2}{4t}>z_6\) for \(|x'|>R_-\geq1+2\sqrt{z_6}\) and \(0<t<1\), \eqref{0506-06} gives \(A_1'(t)<0\).

	Finally, \eqref{0514-02} with \((a,b)=(\frac n2,\frac{n-1}{2})\) gives \( |{}_1F_1(\frac n2,\frac{n-1}{2};-z)|\leq C(1+z)^{-\frac n2}\).
		Thus \(|A_1(t)|\leq \frac{CM_0}{t^{1/2}(|x'|-1)^n}\).
	Lemma~\ref{lem-260032} and \eqref{0506-06} give \(|A_1'(t)|\leq \frac{CM_0}{t^{3/2}(|x'|-1)^{n}}\).
	The bound \eqref{0514-02}, applied with \((a,b)=(\frac{n+2}{2},\frac{n-1}{2})\), also gives \(|A_2(t)|\leq\frac{CM_0}{t^{1/2}(|x'|-1)^{n+2}}\), so \(A_2\in L^1(0,1)\).
	Thus \(A_1\) and \(A_2\) satisfy Lemma \ref{lem-0506-01} with \(\sigma=\frac12\).
Then \eqref{0428-02} gives
	\begin{align*}
		\partial_{x_n}^2w_n(x',0,t) = \frac{\Gamma\left(\frac{n}{2}\right)}{2^{n-1}\pi^{\frac{n}{2}}\Gamma\left(\frac{n-1}{2}\right)}\int_0^t\left( \phi'(s) A_1(t-s) - \phi(s) A_2(t-s)\right)\,ds.
	\end{align*}
	Lemma \ref{lem-0506-01} and Proposition \ref{prop-0504-02} prove that \(t_n^*(x')\) is unique.

	At the earlier zero, Lemma~\ref{lem-0506-01} gives a change from positive to negative because \(\phi\) and the constant in front of the integral are positive. It also gives a negative time derivative there if {\rm\textbf{(T3)}} holds.
	At the later zero, Lemma~\ref{lem-0504-01} and \eqref{0517-10} give a positive time derivative.
	Each zero thus yields a separation point of \(w_n\) by Proposition~\ref{prop-260046}{\rm(i)}.
	Formula \eqref{0730-41} follows from uniqueness, continuity, \eqref{0517-10}--\eqref{0517-11}, and the signs established in the proof of Proposition \ref{prop-0504-02}.

\end{proof}

\subsubsection{Monotonicity of the separation points}\label{subsubsect-0826-01}
The next proposition shows that both \(t_n^*(x')\) and \(\widetilde{t_n^*}(x')\) are strictly increasing along every fixed ray when \(|x'|\) is sufficiently large.
This proves the monotonicity of \(t_n^*(x')\) in Theorem~\ref{thm-0507-01}{\rm(iv)}.
\begin{proposition}\label{prop-0514-01}
	Let \(w\) solve the Stokes system \eqref{0308-01} with boundary data \eqref{0324-13}, where \(j=n\).
	Suppose that \(\psi\) satisfies {\rm\textbf{(S)}}, and that \(\phi\) satisfies {\rm\textbf{(T1)}}--{\rm\textbf{(T4)}}.
	Let \((x', \widetilde{t_n^*}(x'))\) and \((x', {t_n^*}(x'))\) be the separation points of \(w_n\) with \(t_n^*(x')<1<\widetilde{t_n^*}(x')\).
Then there exist \(C_1=C_1(n)>0\) and \(C_2=C_2(n)>0\) such that, for
	\begin{align}\label{0626-10}
		\widetilde R_+ := \frac{z_1+z_2+2\sqrt{z_1z_2(1+z_2-z_1)}}{z_2-z_1}, \qquad \widetilde R_{-}:= \max\left\{1+2\sqrt{Z_0},R_-,1+C_1,1+C_2\right\},
	\end{align}
 with \(R_-\) given in \eqref{0730-40}, \(r\mapsto\widetilde{t_n^*}(r\omega')\) is strictly increasing for \(r>\widetilde R_+\), and \(r\mapsto t_n^*(r\omega')\) is strictly increasing for \(r>\widetilde R_-\), for every fixed \(\omega'\in\mathbb S^{n-2}\).
	If \(\psi\) is radial, then \(t_n^*(x')\) and \(\widetilde{t_n^*}(x')\) depend only on \(|x'|\).
\end{proposition}

To prove Proposition \ref{prop-0514-01}, we need the following monotonicity lemma.
\begin{lemma}\label{lem-0515-02}
	Let \(M>0\), let \(\psi\) satisfy {\rm\textbf{(S)}} with \(\frac{M_2}{M_0}\leq M\), and define
	\begin{align}\label{0514-04}
		K_2(r,t) = \frac{1}{(n-1)t^{\frac{n+1}{2}}} \int_\Sigma \psi(y')\mathcal{K}\left( \frac{|r\mathbf e_1'-y'|^2}{4t}\right)\,dy',
	\end{align}
	where \(\mathcal{K}(z)\) is given in \eqref{0827-12}.
		Then there exist \(C_1,C_2>0\), depending only on \(n\), such that, if \(0<t<1\) and \( r>1+\max\left\{2\sqrt{Z_0},C_1M^{\frac{1}{n+7}},C_2M^{\frac12}\right\}\), then
	\begin{align}\label{0515-11}
		\partial_{t}\partial_{r}\ln(-K_2(r,t))>0,
	\end{align}
	where \(Z_0>0\) is given in Lemma \ref{lem-0520-01}.
\end{lemma}
\begin{proof}
	For \(0<t<1\) and \(|y'|\leq1\), \(u:=\frac{|r\mathbf e_1'-y'|^2}{4t}\geq \frac{(r-1)^2}{4}>Z_0>z_3\), so \(\mathcal{K}(u)<0\) and \(-K_2(r,t)=\int_{\Sigma}h(r,y',t)\,dy'\) with \(h(r,y',t):=-\frac{1}{(n-1)t^{\frac{n+1}{2}}}\mathcal{K}(u)\psi(y')\).
	Let \(f(r,y',t):=\partial_{r}\ln(-t^{-\frac{n+1}{2}}\mathcal{K}(u))\) and \(g(r,y',t):=\partial_{t}\ln(-t^{-\frac{n+1}{2}}\mathcal{K}(u))\).
	The computation in Lemma \ref{lem-0326-01} gives
	\begin{align}\label{0515-011}
		\partial_t\partial_{r}\ln(-K_2(r,t)) & = -\frac{\int_{\Sigma}\partial_t f(r,y',t)h(r,y',t)\,dy'}{K_2(r,t)}\nonumber\\
		 & \quad + \frac{\iint_{\Sigma\times\Sigma}(f(y')-f(z'))(g(y')-g(z'))h(y')h(z')\,dy'dz'}{2K_2(r,t)^2},
	\end{align}
	where \(r\) and \(t\) are suppressed in the double integral.
	Since \(\partial_r u=\frac{r-y_1}{2t}\) and \(\partial_tu=-\frac ut\),
	\begin{align}\label{0515-012}
		f(r,y',t) = \frac{r-y_1}{2t} \mathcal{L}'(u), \qquad \partial_t f(r,y',t) = \frac{r-y_1}{2t^2}\mathcal{M}(u), \qquad
		g(r,y',t) = -\frac{n+1}{2t} - \frac{u}{t}\mathcal{L}'(u),
	\end{align}
	where \(\mathcal{L}\) and \(\mathcal{M}\) are given in Lemma \ref{lem-0520-01}.
	Since \(\mathcal{M}\) is positive and strictly decreasing and \((r-1)^2\leq4tu\leq(r+1)^2\),
	\begin{align}\label{0515-013}
		\int_{\Sigma}\partial_t f(r,y',t)h(r,y',t)\,dy' & \geq \frac{r-1}{2t^2}\mathcal{M}\left( \frac{(r+1)^2}{4t}\right) \int_{\Sigma}h(r,y',t)\,dy'.
	\end{align}
	With \(M_f(r,t) := \sup_{|y'|\leq 1}|\nabla' f(r,y',t)|\) and \(M_g(r,t) := \sup_{|y'|\leq1}|\nabla' g(r,y',t)|\), the mean value theorem, \(\mathcal{L}'<0\), and \eqref{0520-30} give
	\begin{align}\label{0515-014}
		\frac{\left|\iint_{\Sigma\times\Sigma}(f(y')-f(z'))(g(y')-g(z'))h(y')h(z')\,dy'dz'\right|}{K_2(r,t)^2}
		& \leq 2M_f M_g\frac{\int_{\Sigma}|y'|^2h(r,y',t)\,dy'}{-K_2(r,t)}\nonumber\\ & \leq 2M_f M_g\frac{M_2(r+1)^n}{M_0(r-1)^n}.
	\end{align}
	By \eqref{0515-012} and \eqref{0520-40}, \(M_f(r,t)\leq2C_{\mathcal{L}}/(r-1)^2\), and, since \(\mathcal{M}\) is decreasing, \eqref{0520-31} gives
	\begin{align*}
		M_g(r,t) \leq \frac{r+1}{2t^2}\mathcal{M}\left( \frac{(r-1)^2}{4t}\right) \leq\frac{(r+1)^5}{t^2(r-1)^4} \mathcal{M}\left( \frac{(r+1)^2}{4t}\right).
	\end{align*}
	Hence \eqref{0515-011}, \eqref{0515-013}, and \eqref{0515-014} give
	\begin{align*}
		\partial_t\partial_r \ln (-K_2(r,t)) & \geq \frac{1}{2t^2}\mathcal{M}\left(\frac{(r+1)^2}{4t}\right)\left(r-1- 4C_{\mathcal{L}}\frac{(r+1)^{n+5}}{(r-1)^{n+6}} \frac{M_2}{M_0}\right),
	\end{align*}
	which is positive if
	\begin{align}\label{0625-03}
		\frac{M_2}{M_0} < \frac{1}{4C_{\mathcal{L}}} \frac{(r-1)^{n+7}}{(r+1)^{n+5}}.
	\end{align}
	If \(1<r\leq2\), then \(r+1\leq3\), and the right-hand side of \eqref{0625-03} exceeds \(M\) for \(r> 1+ C_1M^{\frac{1}{n+7}}\) with \(C_1=(4\cdot3^{n+5}C_{\mathcal{L}})^{\frac{1}{n+7}}\).
	If \(r>2\), then \(r+1\leq 3(r-1)\), and it exceeds \(M\) for \(r> 1+ C_2M^{\frac12}\) with \(C_2=(4\cdot3^{n+5}C_{\mathcal{L}})^{\frac12}\).
\end{proof}

\begin{proof}[Proof of Proposition \ref{prop-0514-01}]
	Fix \(\omega'\in\mathbb S^{n-2}\) and an orthogonal map \(O:\Sigma\to\Sigma\) with \(O\mathbf e_1'=\omega'\).
	Then \(\psi_O(y')=\psi(Oy')\) satisfies {\rm\textbf{(S)}} with the same \(M_0\) and \(M_2\), \(\psi_O=\psi\) if \(\psi\) is radial, and \eqref{0504-10} and \eqref{0508-01} give \(G_{k,\psi}(r\omega',t)=G_{k,\psi_O}(r\mathbf e_1',t)\) for \(k=1,2\).
	Hence it suffices to consider the positive \(x_1\)-axis; we suppress \(\psi_O\) and write \(G_k(r,t)=G_k(r\mathbf e_1',t)\), \(t_n^*(r)\), and \(\widetilde{t_n^*}(r)\).
	By \eqref{0730-40} and \eqref{0626-10}, \(\widetilde R_+>R_+\) and \(\widetilde R_-\geq R_-\), so Proposition \ref{prop-0505-01} gives the uniqueness of both zeros below.
	By the implicit function theorem, it suffices to prove that \(\partial_tG_1>0\) and \(\partial_rG_1<0\) at \((r,\widetilde{t_n^*}(r))\) for \(r>\widetilde R_+\), and that \(\partial_tG_2>0\) and \(\partial_rG_2<0\) at \((r,t_n^*(r))\) for \(r>\widetilde R_-\).

	\underline{1. Monotonicity of \(\widetilde{t_n^*}(r)\).}\quad
	Let \(r>\widetilde R_+\).
	By \eqref{0504-10},
	\begin{align}\label{0513-03}
		G_1(r,t) = \int_0^1 \phi(s) K_1(r,t-s)\,ds,\qquad K_1(r,t) := \frac{1}{t^{\frac{n+3}{2}}}\int_{\Sigma}\psi(y'){}_1F_1\left(\frac{n}{2}, \frac{n-1}{2}; -\frac{|r\mathbf e_1'-y'|^2}{4t}\right)\,dy',
	\end{align}
	and, since \(\phi(0)=\phi(1)=0\), \(\partial_t G_1(r,t) = \int_0^1 \phi'(s)K_1(r,t-s)\,ds\).
	Since \(G_1(r,\widetilde{t_n^*}(r))=0\) and \(\phi>0\) on \((0,1)\), \(K_1(r,\cdot)\), which is \(K_0(r\mathbf e_1',\cdot)\) in the proof of Proposition \ref{prop-0505-01}, has a zero \(t_*\), and, since \(r>R_+\), that proof shows that \(t_*\) is its only zero, \(\partial_tK_1(r,t_*)>0\), and \((r-1)^2/(4t_*)<z_1\).
	Hence \(K_1(r, \widetilde{t_n^*}(r)-s)\) is positive for \(0<s<s^*\) and negative for \(s^*<s<1\), where \(s^*=\widetilde{t_n^*}(r)-t_*\in(0,1)\).
	Since \(\phi'/\phi\) is strictly decreasing by {\rm\textbf{(T4)}},
	\begin{align*}
		\partial_t G_1(r, \widetilde{t_n^*}(r)) = \int_0^1\phi(s) \frac{\phi'(s)}{\phi(s)}K_1(r, \widetilde{t_n^*}(r)-s)\,ds > \frac{\phi'(s^*)}{\phi(s^*)}G_1(r, \widetilde{t_n^*}(r)) =0.
	\end{align*}
	Moreover, \(\widetilde{t_n^*}(r)>t_*>(r-1)^2/(4z_1)\), and \(r>\widetilde R_+\) gives
	\begin{align}\label{0513-09}
		\frac{|r\mathbf e_1'-y'|^2}{4(\widetilde{t_n^*}(r)-s)}\leq\frac{(r+1)^2}{4(\widetilde{t_n^*}(r)-1)} < \frac{z_1(r+1)^2}{(r-1)^2-4z_1}<z_2,\qquad |y'|\leq1,\quad 0<s<1.
	\end{align}
	Hence, by \eqref{0327-07},
	\begin{align*}
		\partial_r K_1(r,\widetilde{t_n^*}(r)-s) =- \frac{n}{2(n-1)(\widetilde{t_n^*}(r)-s)^{\frac{n+5}{2}}}\int_{\Sigma}\psi(y')(r-y_1) {}_1F_1\left(\frac{n+2}{2}, \frac{n+1}{2}; -\frac{|r\mathbf e_1'-y'|^2}{4(\widetilde{t_n^*}(r)-s)}\right)\,dy'<0,
	\end{align*}
	and \(\partial_r G_1(r, \widetilde{t_n^*}(r)) = \int_0^1 \phi(s) \partial_r K_1(r, \widetilde{t_n^*}(r)-s)\,ds<0\).

	\underline{2. Monotonicity of \({t_n^*}(r)\).}\quad
	Let \(r>\widetilde R_-\).
	By \eqref{0514-01},
	\begin{align}\label{0515-05}
		G_2(r,t) = \int_0^t\phi'(s) K_2(r,t-s)\,ds,
	\end{align}
	where \(K_2\) is given in \eqref{0514-04}, and \(K_2(r,t)<0\) for \(0<t<1\) by the proof of Lemma \ref{lem-0515-02}.
	We first show that
	\begin{align}\label{0515-03}
		t_1<t_n^*(r)<t_\infty.
	\end{align}
	Since \(\phi'>0\) on \((0,t_1)\), \eqref{0515-05} gives \(G_2(r,t)<0\) for \(0<t\leq t_1\), so \(t_n^*(r)>t_1\).
	By Proposition \ref{prop-0505-01}, it remains to show \(G_2(r,t_\infty)>0\).
	With \(W_2(r,t)=\sqrt{t}K_2(r,t)\),
	\begin{align}\label{0515-02}
		\partial_t W_2(r,t) = - \frac{t^{-\frac{n}{2}-1}}{n-1} \int_\Sigma \psi(y') \left[ \frac{n}{2} \mathcal{K}\left(\frac{|r\mathbf e_1'-y'|^2}{4t}\right) + \frac{|r\mathbf e_1'-y'|^2}{4t}\mathcal{K}'\left(\frac{|r\mathbf e_1'-y'|^2}{4t}\right) \right]\,dy',
	\end{align}
	and, since \(r>1+2\sqrt{Z_0}>1+2\sqrt{z_4}\), Lemma \ref{lem-0514-01} makes the bracket negative, so \(W_2(r,\cdot)\) is strictly increasing.
	Splitting at \(t_1\) and applying \eqref{0324-17} give
	\begin{align*}
		G_2(r,t_\infty) = \int_0^{t_\infty}\frac{\phi'(s)}{\sqrt{t_\infty-s}} W_2(r,t_\infty-s)\,ds >W_2(r,t_\infty-t_1)\int_0^{t_\infty} \frac{\phi'(s)}{\sqrt{t_\infty-s}}\,ds=0.
	\end{align*}

	By Lemma \ref{lem-0618-01}, \(|\partial_t K_2(r,t)|\leq CM_0(r-1)^{-n}t^{-3/2}\) for \(0<t<1\), so \(K_2\) satisfies \eqref{0829-01} with \(\sigma=\frac12\).
	Since {\rm\textbf{(T1)}} and {\rm\textbf{(T3)}} give \(\phi'\in C^\beta([0,1])\) with \(\beta\in(\frac12,1)\), Lemma \ref{lem-0424-20} applies to \eqref{0515-05} with \(f=\phi'\) and gives the joint continuity of \(\partial_tG_2\) and \(\partial_rG_2\).
	Since \(r>R_-\), applying Lemma \ref{lem-0506-01} as in the proof of Proposition \ref{prop-0505-01} gives \(\partial_t\partial_{x_n}^2w_n(r\mathbf e_1',0,t_n^*(r))<0\), and \eqref{0517-11} gives \(\partial_tG_2(r,t_n^*(r))>0\).

	Let \(k(s):=\partial_r\ln(-K_2(r,t_n^*(r)-s))\).
	Since \(\supp\psi\subset\overline{B_1'}\) gives \(M_2/M_0\leq1\), Lemma \ref{lem-0515-02} applies with \(M=1\), and \(k\) is strictly decreasing.
	By \eqref{0515-03}, \(\phi'(s)K_2(r,t_n^*(r)-s)\) is negative on \((0,t_1)\) and positive on \((t_1,t_n^*(r))\).
	Hence \eqref{0515-05} gives
	\begin{align*}
		\partial_r G_2(r,t_n^*(r)) =\int_0^{t_n^*(r)} \phi'(s)K_2(r,t_n^*(r)-s)k(s)\,ds < k(t_1)\int_0^{t_n^*(r)}\phi'(s)K_2(r,t_n^*(r)-s)\,ds=0.
	\end{align*}
\end{proof}

\subsubsection{Asymptotic behavior of the separation points}
The next proposition gives the large-\(|x'|\) expansions of \(t_n^*(x')\) and \(\widetilde{t_n^*}(x')\), with \(t_n^*(x')\) approaching \(t_\infty\) and \(\widetilde{t_n^*}(x')\) growing quadratically.
This proves the asymptotics in Theorem \ref{thm-0507-01}{\rm(iv)} and the monotonicity of \(\widetilde{t_n^*}(x')\) there.
\begin{proposition}\label{prop-0421-4}
	Let \(w\) solve the Stokes system \eqref{0308-01} with boundary data \eqref{0324-13}, where \(j=n\).
	Suppose that \(\psi\) satisfies {\rm\textbf{(S)}}, and that \(\phi\) satisfies {\rm\textbf{(T1)}} and \(\phi\not\equiv0\).
	For sufficiently large \(|x'|\), let \(\widetilde{t_n^*}(x')\) be the unique zero of \(\partial_{x_n}^2w_n(x',0,\cdot)\) in \((1,\infty)\). Then the following hold.
	\begin{enumerate}
			\item There exists \(N>1\) such that, for \(|x'|>N\),
		\begin{align}\label{0502-01}
			\begin{split}
				\widetilde{t_n^*}(x') & = \frac{1}{4z_1}|x'|^2 - \frac{\mathbf m\cdot x'}{2z_1M_0} + \frac{\norm{s\phi}_{L^1(0,1)}}{\norm{\phi}_{L^1(0,1)}} + \frac{M_2}{4z_1M_0}\\
				& \quad -\frac{2z_1+n-1}{4z_1M_0}\int_{\Sigma}\left| \frac{x'}{|x'|}\cdot \left( y' - \frac{\mathbf m}{M_0}\right)\right|^2 \psi(y')\,dy' + O(|x'|^{-1}).
			\end{split}
		\end{align}
		Moreover, uniformly for \(\omega'\in\mathbb S^{n-2}\) as \(r\to\infty\),
		\begin{align*}
			\frac{d}{dr}\widetilde{t_n^*}(r\omega')=\frac{r}{2z_1}-\frac{\omega'\cdot\mathbf m}{2z_1M_0}+O(r^{-2}).
		\end{align*}
		\item\label{item-0731-90} Suppose that \(\phi\) satisfies {\rm\textbf{(T2)}}--{\rm\textbf{(T4)}}, and let \(t_n^*(x')\) be the unique zero of \(\partial_{x_n}^2w_n(x',0,\cdot)\) in \((0,1)\).
			Then there exists \(N>1\) such that, for \(|x'|>N\),
		\begin{align}\label{0504-05}
			t_n^*(x') = t_\infty+\frac{3nT_1(t_\infty)}{T_0'(t_\infty)}|x'|^{-2} +O(|x'|^{-1-2\beta}),
		\end{align}
		where \(T_0\) and \(T_1\) are defined in \eqref{0323-02}.
		If, in addition, \(\phi\) is \(C^{2,\alpha_1}\) near \(t_\infty\) for some \(\alpha_1\in(1/2,1)\), then
		\begin{align}\label{260057}
			t_n^*(x') & =t_\infty+\frac{3nT_1(t_\infty)}{T_0'(t_\infty)}|x'|^{-2} +\frac{6nT_1(t_\infty)}{T_0'(t_\infty)} \frac{x'\cdot\mathbf m}{M_0|x'|^4}+O(|x'|^{-4}).
		\end{align}
	\end{enumerate}
\end{proposition}

\begin{proof}
	
	\underline{Proof of (i).}\quad
	We first consider \(\widetilde{t_n^*}(x')\).
	Let \(\ep=|x'|^{-1}\), \(\omega'=x'/|x'|\), and \(\tau=\ep^2t\).
	Then \eqref{0504-10} gives
	\begin{align}\label{0503-02}
		G_1(x',t) = \ep^{n+3}\mathcal{G}_1(\tau, \ep, \omega'),
	\end{align}
	where \(\omega'\in \mathbb S^{n-2}\) and
	\begin{align}\label{0503-01}
		\mathcal{G}_1(\tau,\ep,\omega') = \int_{\Sigma}\psi(y')\int_0^1 \frac{\phi(s)}{(\tau-\ep^2 s)^{\frac{n+3}{2}}}e^{-\frac{1-2\ep \omega'\cdot y' + \ep^2 |y'|^2}{4(\tau-\ep^2 s)}}G\left( \frac{1-2\ep \omega'\cdot y' + \ep^2 |y'|^2}{4(\tau-\ep^2 s)}\right)\,dsdy',
	\end{align}
	where \(G(z) := {}_1F_1\left(-\frac12, \frac{n-1}{2};z\right)\).
	Setting \(\ep=0\) gives
	\begin{align}\label{0620-01}
		\mathcal{G}_1(\tau, 0,\omega' ) = \frac{e^{-\frac{1}{4\tau}}}{\tau^{\frac{n+3}{2}}}{}_1F_1\left(-\frac12, \frac{n-1}{2}; \frac{1}{4\tau}\right) M_0\norm{\phi}_{L^1(0,1)}.
	\end{align}
Set \(\tau_0=\frac{1}{4z_1}\); the definition of \(z_1\) gives \(\mathcal{G}_1(\tau_0,0,\omega')=0\).

	Replacing \(1-2\ep\omega'\cdot y'+\ep^2|y'|^2\) in \eqref{0503-01} by \(|v'-\ep y'|^2\), with \(v'\in\Sigma\) near \(\mathbb S^{n-2}\), extends \(\mathcal{G}_1\) smoothly in \((\tau,\ep,v')\) near \(\{\tau_0\}\times\{0\}\times\mathbb S^{n-2}\).

	Then \eqref{0620-01} and \eqref{0327-07} give
	\begin{align*}
		\partial_\tau \mathcal{G}_1\left(\tau_0,0,\omega'\right) 
		 & = \frac{e^{-z_1}}{4(n-1)\tau_0^{\frac{n+7}{2}}}{}_1F_1\left(\frac12, \frac{n+1}{2}; z_1\right)M_0\norm{\phi}_{L^1(0,1)} >0,
	\end{align*}
	where the last inequality follows from \eqref{0324-14}.
	Its positivity, continuity, and the compactness of \(\mathbb S^{n-2}\) give \(\delta>0\) and \(\ep_0>0\) such that \(\partial_\tau\mathcal G_1(\tau,\ep,\omega')>0\) whenever \(|\tau-\tau_0|<\delta\), \(|\ep|<\ep_0\), and \(\omega'\in\mathbb S^{n-2}\).
	The implicit function theorem and the compactness of \(\mathbb S^{n-2}\) give a smooth function \(\widetilde{\tau_*}:(-\ep_0,\ep_0)\times\mathbb S^{n-2}\to\R\) with \(\widetilde{\tau_*}(0,\omega')=\tau_0\) and \(\mathcal{G}_1(\widetilde{\tau_*}(\ep,\omega'),\ep,\omega')=0\).

A Taylor expansion in \(\ep\), uniform for \(\omega'\in\mathbb S^{n-2}\), gives \(\widetilde{\tau_*}(\ep,\omega') =\tau_0+\mu_1(\omega')\ep+\mu_0(\omega')\ep^2+O(\ep^3)\).
	For \(\ep>0\) sufficiently small, \(\ep^{-2}\widetilde{\tau_*}(\ep,\omega')>1\).
	Equation \eqref{0503-02} and the uniqueness for \(t>1\) in Proposition~\ref{prop-0504-02} identify it with \(\widetilde{t_n^*}(x')\).
	Consequently, there exists \(R_0>0\) such that for \(|x'|\geq R_0\),
	\begin{align}\label{0503-08}
		\widetilde{t_n^*}(x')=\frac{1}{4z_1}|x'|^2 + \mu_1(\omega')|x'| + \mu_0(\omega')+O(|x'|^{-1}).
	\end{align}

By \eqref{0503-08},
	\begin{align}\label{0503-03}
		\begin{split}
		\frac{|x'-y'|^2}{4(\widetilde{t_n^*}(x')-s)} & = z_1- z_1\left( \frac{2x'\cdot y'}{|x'|} + 4z_1\mu_1\right) \frac{1}{|x'|}\\
			 & \quad + z_1\left(16z_1^2\mu_1^2- 4z_1(\mu_0-s)+ \frac{8x'\cdot y'}{|x'|}z_1\mu_1 + |y'|^2\right) \frac{1}{|x'|^2} + O(|x'|^{-3}).
		\end{split}
	\end{align}
		We also have
	\begin{align*}
		\frac{1}{(\widetilde{t_n^*}(x')-s)^{\frac{n+3}{2}}} & = \left(\frac{4z_1}{|x'|^2}\right)^{\frac{n+3}{2}}\left(1-\frac{2(n+3)z_1\mu_1}{|x'|} + O(|x'|^{-2})\right).
	\end{align*}
	Finally, expanding in a Taylor series and using \eqref{0503-03} gives
	\begin{align*}
		e^{-z}G(z) 
		  = e^{-z_1}\left[G'\left(z_1\right) \left(z-z_1\right) + \frac12 G''\left(z_1\right) \left(z-z_1\right)^2 - G'\left(z_1\right)\left(z-z_1\right)^2+ O\left(\left(z-z_1\right)^3\right)\right].
	\end{align*}
Substituting these expansions into \(G_1(x',\widetilde{t_n^*}(x'))=0\) and dividing by \(e^{-z_1}(4z_1/|x'|^2)^{(n+3)/2}\) gives
	\begin{align}\label{0503-04}
		\begin{split}
			0 & =-z_1G'\left(z_1\right) \int_{\Sigma}\psi(y')\int_0^1 \phi(s)\left(\frac{2x'\cdot y'}{|x'|} + 4z_1\mu_1\right)\,dsdy' \frac{1}{|x'|}\\
			& \quad + z_1G'\left(z_1\right) \int_{\Sigma}\psi(y')\int_0^1 \phi(s) \left(\left(4z_1\mu_1\right)^2- 4z_1(\mu_0-s) + \frac{2x'\cdot y'}{|x'|}4z_1\mu_1 + |y'|^2\right)\,dsdy' \frac{1}{|x'|^2} \\
			& \quad - z_1^2G'\left(z_1\right) \int_\Sigma\psi(y')\int_0^1 \phi(s)\left(\frac{2x'\cdot y'}{|x'|}+4z_1\mu_1\right)^2\,dsdy'\frac{1}{|x'|^2}\\
			& \quad + \frac{z_1^2}{2} G''\left(z_1\right) \int_\Sigma\psi(y')\int_0^1 \phi(s)\left(\frac{2x'\cdot y'}{|x'|}+4z_1\mu_1\right)^2\,dsdy'\frac{1}{|x'|^2}\\
			& \quad + {2(n+3)z_1^2\mu_1}G'\left(z_1\right) \int_{\Sigma}\psi(y')\int_0^1 \phi(s)\left(\frac{2x'\cdot y'}{|x'|} + 4z_1\mu_1\right)\,dsdy'\frac{1}{|x'|^2} + O(|x'|^{-3}).
		\end{split}
	\end{align}

Comparing the coefficients of \(|x'|^{-1}\) in \eqref{0503-04} gives
	\begin{align}\label{0503-05}
		\mu_1 (\omega')= - \frac{\omega'\cdot\mathbf m}{2z_1M_0}.
	\end{align}
	With \eqref{0503-05}, the first and fifth terms on the right-hand side of \eqref{0503-04} vanish.
Then, comparing the coefficients of \(|x'|^{-2}\) in \eqref{0503-04} gives
	\begin{align}\label{0503-06}
		\begin{split}
			& 4z_1\left(1-\frac{G''\left(z_1\right)}{2G'\left(z_1\right)}\right) \int_{\Sigma}\psi(y')\int_0^1 \phi(s)\left( \frac{x'}{|x'|}\cdot \left(y'- \frac{\mathbf m}{M_0}\right)\right)^2\,dsdy'\\
			& \quad = M_2\norm{\phi}_{L^1(0,1)} - 4z_1M_0\left(\mu_0\norm{\phi}_{L^1(0,1)}-\norm{s\phi}_{L^1(0,1)}\right),
		\end{split}
	\end{align}
	where \eqref{0503-05} gives \(16z_1^2\mu_1^2 M_0+8z_1\mu_1\int_\Sigma (\omega'\cdot y')\psi(y')\,dy'=0\).

	Kummer's equation \eqref{0327-03} and \(G(z_1)=0\) give \(1-\frac{G''(z_1)}{2G'(z_1)}=\frac12+\frac{n-1}{4z_1}\).
	Dividing \eqref{0503-06} by \(M_0\norm{\phi}_{L^1(0,1)}\) and isolating \(\mu_0\) give
	\begin{align}\label{0503-07}
		\mu_0= \frac{\norm{s\phi}_{L^1(0,1)}}{\norm{\phi}_{L^1(0,1)}} + \frac{M_2}{4z_1M_0}-\frac{2z_1+n-1}{4z_1M_0}\int_{\Sigma}\left| \frac{x'}{|x'|}\cdot \left( y' - \frac{\mathbf m}{M_0}\right)\right|^2 \psi(y')\,dy'.
	\end{align}
	Combining \eqref{0503-05}, \eqref{0503-07}, and \eqref{0503-08} yields the desired result \eqref{0502-01}.

	Differentiating \(\widetilde{t_n^*}(r\omega')=r^2\widetilde{\tau_*}(1/r,\omega')\) gives, since the terms \(2\mu_0(\omega')/r\) cancel,
\begin{align*}
\frac{d}{dr}\widetilde{t_n^*}(r\omega')
=2r\widetilde{\tau_*}(1/r,\omega')-\partial_\ep\widetilde{\tau_*}(1/r,\omega')
=2\tau_0r+\mu_1(\omega')+O(r^{-2}).
\end{align*}
By \eqref{0503-05}, this is the stated derivative, which is positive for large \(r\) in every direction.

	\underline{Proof of (ii).}\quad
	For all sufficiently large \(|x'|\), the computation in \eqref{0515-02}, with \(r\mathbf e_1'\) replaced by \(x'\), shows that \(W_2(x',t)\) is strictly increasing in \(t\in(0,1)\).
		Since \(W_2(x',t)<0\) there, \eqref{0628-02} gives \(G_2(x',t)<0\) for \(0<t\leq t_1\).
	Splitting the integral in \eqref{0628-02} at \(t_1\), using the fact that \(W_2\) is strictly increasing, and applying \eqref{0324-17} give \(G_2(x',t_\infty)>W_2(x',t_\infty-t_1)T_0(t_\infty)=0\).
	Proposition \ref{prop-0505-01} therefore gives \(t_1<t_n^*(x')<t_\infty\) for all sufficiently large \(|x'|\).

	Applying \eqref{0504-06} in \eqref{0508-01} and integrating by parts gives, uniformly for \(t\in[t_1,t_\infty]\) and \(\omega'\in\mathbb S^{n-2}\),
	\begin{align}\label{0518-20}
		\begin{split}
			\ep^{-n}G_2(x',t) & =\frac{2^n\Gamma\left(\frac{n-1}{2}\right)} {\Gamma\left(-\frac12\right)} \left[ T_0(t)\int_\Sigma \frac{\psi(y')}{|\omega'-\ep y'|^n}\,dy' -3n\ep^2T_1(t)\int_\Sigma \frac{\psi(y')}{|\omega'-\ep y'|^{n+2}}\,dy' +O(\ep^4) \right],
		\end{split}
	\end{align}
where \(T_1(t)=\frac12\int_0^t\phi(s)(t-s)^{-1/2}\,ds\) by integration by parts.

	Set \(t^*=t_n^*(x')\).
	Since both spatial integrals in \eqref{0518-20} equal \(M_0+O(\ep)\), the identity \(G_2(x',t^*)=0\) gives \(T_0(t^*)=O(\ep^2)\).
	The uniqueness of the zero of \(T_0\) in \([t_1,t_\infty]\) then implies \(t^*\to t_\infty\).
	Remark~\ref{rmk-0730-03} gives \(T_0'(t_\infty)<0\) under {\rm\textbf{(T1)}}--{\rm\textbf{(T4)}}.
	Consequently, \(t^*-t_\infty=O(\ep^2)\).
	Lemma~\ref{lem-0424-20} yields \(T_0\in C^{1,\beta-1/2}_{\mathrm{loc}}(0,1)\).
	Using this estimate, \(T_1'=T_0/2\), and the uniform spatial expansion \eqref{0518-20}, we obtain
	\begin{align*}
		T_0'(t_\infty)(t^*-t_\infty)
		=3nT_1(t_\infty)\ep^2+O(\ep^3+|t^*-t_\infty|^{\beta+1/2}).
	\end{align*}
	The remainder is \(O(\ep^{1+2\beta})\), uniformly in \(\omega'\), which proves \eqref{0504-05}.

	Assume finally that \(\phi\) is \(C^{2,\alpha_1}\) near \(t_\infty\) for some \(\alpha_1\in(1/2,1)\).
	For \(p\in\{n,n+2\}\),
	\begin{align*}
		\int_\Sigma\frac{\psi(y')}{|\omega'-\ep y'|^p}\,dy' =M_0+p\ep\,\omega'\cdot\mathbf m+O(\ep^2).
	\end{align*}
	The cutoff argument in the proof of Proposition~\ref{prop-0731-11} gives \(T_0\in C^2\) near \(t_\infty\).
	Thus \(T_0(t^*)=T_0'(t_\infty)(t^*-t_\infty)+O(\ep^4)\), while \(T_1'=T_0/2\) gives \(T_1(t^*)=T_1(t_\infty)+O(\ep^4)\).
	Therefore \eqref{0518-20} gives
	\begin{align*}
		T_0'(t_\infty)(t^*-t_\infty) =3nT_1(t_\infty)\ep^2 +\frac{6nT_1(t_\infty)}{M_0} (\omega'\cdot\mathbf m)\ep^3+O(\ep^4).
	\end{align*}
	Substituting \(\ep=|x'|^{-1}\) gives \eqref{260057}.
\end{proof}

	\section{Separation points for tangential boundary data}\label{secttan}
In this section, we prove Theorems \ref{thm-260010}, \ref{thm-260014}, and \ref{thm-260017}.
Throughout this section, \(1\leq j\leq n-1\) in \eqref{0324-13}, so that \(g(x',t)=\psi(x')\phi(t)\mathbf e_j\), and we consider the following cases.
\begin{enumerate}
	\item \(i=j\) (Subsection \ref{260058}).
	\item \(i\neq j\), \(i\neq n\) (Subsection \ref{260092}).
	\item \(i=n\) (Subsection \ref{260114}).
\end{enumerate}

	\subsection{Tangential velocity component in the direction of the boundary data}\label{260058}
\subsubsection{Formula for the normal derivative of \(w_i\)}\label{260059}

\begin{proposition}\label{prop-0703-01}
Let \(1\leq i\leq n-1\), and let \(w\) solve the Stokes system \eqref{0308-01} with boundary data \eqref{0324-13}, where \(j=i\), \(\psi\) satisfies {\rm\textbf{(S)}}, and \(\phi\) satisfies {\rm\textbf{(T1)}}.
	Then, for every \(x'\in \Sigma\) and \(t>0\),
	\begin{align}\label{0331-11}
		\partial_{x_n}w_i(x',0,t) & = \frac{1}{\sqrt{4\pi}}\int_{-\infty}^t \frac{1}{(t-s)^{\frac32}}\int_{\Sigma}\Gamma'(x'-y',t-s)\bigl(\psi(y')\phi(s)-\psi(x')\phi(t)\bigr)\,dy'ds\nonumber \\
		 & \quad- \phi(t)\partial_{x_i}R_i'\psi(x').
	\end{align}
\end{proposition}

\begin{proof}
	Recall from \eqref{0330-01} and \eqref{0308-02} that
	\begin{align*}
			w_i(x,t) & = -2\int_0^t\int_\Sigma \partial_{x_n}\Gamma(x-y',t-s) \psi(y')\phi(s)\,dy'ds + 4\int_0^t\int_\Sigma L_{ii}(x-y',t-s) \psi(y')\phi(s)\,dy'ds\\
			 & =:w_i^G(x,t)+w_{ii}(x,t).
	\end{align*}

	We first consider \(w_i^G\).
	By \eqref{260269}, \( w_i^G(x',0,t)=\psi(x')\phi(t)\), and \( \int_{-\infty}^t\int_\Sigma \partial_{x_n}\Gamma(x-y',t-s)\,dy'ds =-\frac12. \)
	Consequently,
	\begin{align}\label{0711-02}
		w_i^G(x,t)-\psi(x')\phi(t) = -2\int_{-\infty}^t\int_\Sigma \partial_{x_n}\Gamma(x-y',t-s) \bigl(\psi(y')\phi(s)-\psi(x')\phi(t)\bigr) \,dy'ds.
	\end{align}
	For \(\tau>0\) we have the exact identity
	\begin{align}\label{0711-01}
		\frac{\partial_{x_n}\Gamma(x-y',\tau)}{x_n} = -\frac{1}{4\sqrt{\pi}\,\tau^{3/2}} e^{-\frac{x_n^2}{4\tau}} \Gamma'(x'-y',\tau).
	\end{align}

	Since \(\psi\in C_c^{2,\alpha}(\Sigma)\) and \(\phi\in C^1([0,\infty))\), \(\int_\Sigma \Gamma'(x'-y',t-s) \bigl(\psi(y')\phi(s)-\psi(x')\phi(t)\bigr)\,dy'=O(t-s)\) as \(s\to t\).
	Using \eqref{0711-01} in \eqref{0711-02} and the dominated convergence theorem, we conclude that
	\begin{align}\label{0330-03}
		\partial_{x_n}w_i^G(x',0,t) = \frac1{\sqrt{4\pi}} \int_{-\infty}^t \frac1{(t-s)^{3/2}} \int_\Sigma \Gamma'(x'-y',t-s) \bigl(\psi(y')\phi(s)-\psi(x')\phi(t)\bigr) \,dy'ds.
	\end{align}

	We now consider \(w_{ii}\).
	Using \eqref{0419-02}, we have
	\begin{align*}
		w_{ii}(x,t) = 4\int_0^t\int_\Sigma \partial_{x_i} \int_{\Sigma\times(0,x_n)} \partial_{z_n}\Gamma(z-y',t-s)\, \partial_{x_i}N(x-z)\,dz\, \psi(y')\phi(s)\,dy'ds.
	\end{align*}
	Recall the potential \(k=k_\psi\) defined in \eqref{0703-05}.

	Fix \(t>0\), and set \(P_r=2\partial_{x_n}N(\cdot,r)\) for \(r>0\).
	Formula \eqref{260257}, with \(f=\psi\) and \(j=i\), gives
	\begin{align}\label{0708-13}
	    w_{ii}(\cdot,x_n,t)=2\int_0^{x_n} P_{x_n-q}*_\Sigma R_i'\partial_{x_i}k(\cdot,q,t)\,dq .
	\end{align}
	Here the term \(-2k\) in \eqref{260257} is precisely \(w_i^G\).

	For \(r,q>0\), the multipliers of \(P_r\) and \(R_i'\) are \(e^{-r|\xi'|}\) and \(-i\xi_i/|\xi'|\), respectively.
	Hence
	\begin{align}\label{0712-05}
	    (\partial_r P_r)*_\Sigma R_i'\partial_{x_i}k(\cdot,q,t)=P_r*_\Sigma\partial_{x_i}^2k(\cdot,q,t).
	\end{align}
	Moving the two tangential derivatives onto \(\psi\) in \eqref{0703-05} and using \(\int_0^\infty|\partial_q\Gamma_1(q,s)|\,ds=1/2\) gives
	\begin{align}\label{0712-03}
	    \sup_{q>0}\|\partial_{x_i}^2 k(\cdot,q,t)\|_\infty
	    \leq \tfrac12\|\phi\|_{L^\infty(0,1)}\|\partial_{x_i}^2\psi\|_\infty .
	\end{align}

	Since \(\|P_r\|_1=1\), \eqref{0712-03} bounds the derivative in \eqref{0712-05} uniformly down to \(r=0\).
	Differentiating \eqref{0708-13} by dominated convergence yields
	\begin{align}\label{0708-14}
	    \partial_{x_n} w_{ii}(\cdot,x_n,t)
	    =2R_i'\partial_{x_i}k(\cdot,x_n,t)+2\int_0^{x_n} P_{x_n-q}*_\Sigma\partial_{x_i}^2k(\cdot,q,t)\,dq .
	\end{align}

	Substitute \eqref{0703-05} into the last integral, with its tangential derivatives on \(\psi\).
	For each \(s>0\), \eqref{0708-23} identifies the resulting spatial convolution with \(L_{ni}(\cdot,x_n,s)*_\Sigma\partial_{x_i}\psi\) after moving one tangential derivative onto the heat kernel.
	Keeping the factor \(P_{x_n-q}=2\partial_{x_n}N(\cdot,x_n-q)\) gives
	\begin{align}\label{0708-38}
	    \partial_{x_n}w_{ii}(x,t)=2\partial_{x_i}R_i' k(x,t)
	    +4\int_0^t\phi(t-s)\bigl(L_{ni}(\cdot,x_n,s)*_\Sigma\partial_{x_i}\psi\bigr)(x')\,ds .
	\end{align}

	We now let \(x_n\to0^+\) in \eqref{0708-38}. We first claim that
	\begin{align}\label{0708-15}
		\lim_{x_n\to0^+} \int_0^t\int_\Sigma L_{ni}(x-z',t-s) \partial_{z_i}\psi(z')\phi(s)\,dz'ds =0.
	\end{align}
	Substituting \(s\mapsto t-s\) in the time integral and using the cancellation identity \eqref{0708-22} in Lemma~\ref{lem-260040}, we obtain
	\begin{align*}
		\int_0^t\int_\Sigma L_{ni}(x-z',t-s) \partial_{z_i}\psi(z')\phi(s)\,dz'ds = \int_0^t\phi(t-s) \int_\Sigma L_{ni}(x-z',s) \bigl( \partial_{z_i}\psi(z')-\partial_{x_i}\psi(x') \bigr)\,dz'ds.
	\end{align*}
	Since \(\psi\in C_c^{2,\alpha}(\Sigma)\), \( |\partial_{z_i}\psi(z')-\partial_{x_i}\psi(x')| \leq C|z'-x'|\), which, together with Lemma \ref{lem-260040}\eqref{260042} and \eqref{260043} and the dominated convergence theorem, gives \eqref{0708-15}.

	We next claim that
	\begin{align}\label{0708-37}
		\lim_{x_n\to0^+} \partial_{x_i}R_i'k(x,t) = -\frac12\phi(t)\partial_{x_i}R_i'\psi(x').
	\end{align}
	Writing \(\tau=t-s\) in \eqref{0703-05}, we obtain
	\begin{align*}
		k(x,t) = \int_0^t \partial_{x_n}\Gamma_1(x_n,\tau) \bigl(\Gamma'(\cdot,\tau)*_\Sigma\psi\bigr)(x') \phi(t-\tau)\,d\tau.
	\end{align*}
	Hence
	\begin{align}\label{0708-39}
		\partial_{x_i}R_i'k(x,t) = \int_0^t \partial_{x_n}\Gamma_1(x_n,\tau) \bigl( \Gamma'(\cdot,\tau)*_\Sigma \partial_{x_i}R_i'\psi \bigr)(x') \phi(t-\tau)\,d\tau.
	\end{align}
	By Lemma \ref{lem-260265}, \(h:=\partial_{x_i}R_i'\psi=R_i'(\partial_{x_i}\psi)\) belongs to \(C_0(\Sigma)\).
			Consequently, \(q_{x',t}(\tau) := \bigl(\Gamma'(\cdot,\tau)*_\Sigma h\bigr)(x') \phi(t-\tau)\) is bounded on \((0,t)\) and extends continuously to \(\tau=0\), with \( q_{x',t}(0)=h(x')\phi(t)\).
	Since the kernels \(-2\partial_{x_n}\Gamma_1(x_n,\tau)\) form an approximate identity at \(\tau=0\),
	\begin{align}\label{0712-08}
		\lim_{x_n\to0^+} \int_0^t \left(-\partial_{x_n}\Gamma_1(x_n,\tau)\right) q_{x',t}(\tau)\,d\tau = \frac12q_{x',t}(0),
	\end{align}
	and \eqref{0708-39} together with \eqref{0712-08} gives \eqref{0708-37}.

	Finally, taking \(x_n\to0^+\) in \eqref{0708-38} and using \eqref{0708-15} and \eqref{0708-37}, we obtain
	\begin{align}\label{0330-02}
		\partial_{x_n}w_{ii}(x',0,t) = -\phi(t)\partial_{x_i}R_i'\psi(x').
	\end{align}
	Combining \eqref{0330-03} and \eqref{0330-02} gives \eqref{0331-11}.
\end{proof}

From now on, we take \(i=j=1\).
The general case follows by relabeling the tangential coordinates.
We continue to assume that \(\psi\) satisfies {\rm\textbf{(S)}} and that \(\phi\) satisfies {\rm\textbf{(T1)}}, in which case \eqref{0331-11} can be written as
\begin{align}\label{0331-10}
	\begin{split}
		\partial_{x_n}w_1(x',0,t) & = \frac{1}{\sqrt{4\pi}} \int_0^t \frac{\phi(s)}{(t-s)^{3/2}} \int_\Sigma \Gamma'(x'-y',t-s) \bigl(\psi(y')-\psi(x')\bigr)\,dy'ds \\
		 & \quad+ \frac{\psi(x')}{\sqrt{4\pi}} \int_{0}^t \frac{\phi(s)-\phi(t)}{(t-s)^{3/2}}\,ds -\frac{\psi(x')\phi(t)}{\sqrt{\pi t}} - \phi(t)\partial_{x_1}R_1'\psi(x').
	\end{split}
\end{align}

We restrict attention to the region \(|x'|>1\), where \(\psi(x')=0\) and \eqref{0331-10} takes the form
\begin{align}\label{0401-09}
	\partial_{x_n}w_1(x',0,t) = \frac{1}{\sqrt{4\pi}}\int_0^t \frac{\phi(s)}{(t-s)^{\frac32}} \int_{\Sigma}\Gamma'(x'-y',t-s)\psi(y')\,dy'ds - \phi(t)\mathcal S(x'),
\end{align}
where \(\mathcal S\) is given in \eqref{260009} with \(i=1\).

\subsubsection{Nonexistence of separation points}\label{260060}
The next proposition excludes separation points where \(\mathcal S\leq0\), after the boundary data vanish, and in the regions specified below for radial \(\psi\).
This proves Theorem \ref{thm-260010}{\rm(i)}.

\begin{proposition}\label{prop-260061}
	Let \(w\) solve the Stokes system \eqref{0308-01} with boundary data \eqref{0324-13}, where \(j=1\), \(\psi\) satisfies {\rm\textbf{(S)}}, and \(\phi\) satisfies {\rm\textbf{(T1)}} and {\rm\textbf{(T2)}}.
	Let \(\mathcal S\) be given in \eqref{260009} with \(i=1\).
	Then the following hold.
	\begin{enumerate}
		\item[\rm(i)]\label{item-0907-01} If \(|x'|>1\) and \(t\geq1\), then \(\partial_{x_n}w_1(x',0,t)>0\), and \((x',t)\) is not a separation point of \(w_1\).

		\item[\rm(ii)]\label{item-0907-02} If \(|x'|>1\) and \(\mathcal S(x')\leq0\), then, for every \(t>0\), \(\partial_{x_n}w_1(x',0,t)>0\) and \((x',t)\) is not a separation point of \(w_1\).

		\item[\rm(iii)]\label{item-0907-03} If \(n=2\), then, for every \(|x_1|>1\) and \(t>0\), \(\mathcal S(x_1)<0\), \(\partial_{x_n}w_1(x_1,0,t)>0\), and \((x_1,t)\) is not a separation point of \(w_1\).

		\item[\rm(iv)]\label{item-0907-04} Suppose that \(n\geq3\) and \(\psi\) is radial.
		Then no separation point of \(w_1\) lies on the \(x_1\)-axis with \(|x_1|>1\), and for every \(0<\varepsilon<\sqrt{n-1}\) there exists \(r_\varepsilon>1\) such that no separation point \((x',t)\) satisfies
		\begin{align*}
			|x'|\geq r_\varepsilon \qquad\text{and}\qquad |x''|\leq\bigl(\sqrt{n-1}-\varepsilon\bigr)\,|x_1|.
		\end{align*}

		\item[\rm(v)]\label{item-0907-05} Suppose that \(n\geq3\) and \(\psi\) is radial.
		Then there exists \(r_0=r_0(n,\psi)>0\), independent of \(t\), such that, if
\begin{align}\label{260062}
|x_1|>r_0 \quad\text{and}\quad |x''|^2-(n-1)x_1^2\leq\frac{M_2}{(n-1)M_0},
\end{align}
then \(\mathcal S(x')<0\), and \((x',t)\) is not a separation point of \(w_1\) for any \(t>0\).
	\end{enumerate}
\end{proposition}

\begin{proof}
	Write
	\begin{align}\label{260063}
		K_{\mathcal S}(x'):=\frac{|x'|^2-nx_1^2}{|x'|^{n+2}} =\frac{|x''|^2-(n-1)x_1^2}{|x'|^{n+2}},\qquad x'\neq0.
	\end{align}
	Then \eqref{260009} with \(i=1\) and \eqref{0419-021} give
	\begin{align}\label{260064}
		\mathcal S(x')=\frac{\Gamma\left(\frac n2\right)}{\pi^{\frac n2}} \int_\Sigma K_{\mathcal S}(x'-y')\psi(y')\,dy', \qquad x'\notin\supp\psi.
	\end{align}
	The first term on the right-hand side of \eqref{0401-09} is strictly positive for \(|x'|>1\) and \(t>0\).
		Also, {\rm\textbf{(T1)}} and {\rm\textbf{(T2)}} imply that \(\phi>0\) on \((0,1)\), \(\phi\geq0\) on \((0,\infty)\), and \(\phi=0\) on \([1,\infty)\).

	\underline{Proof of (i).}\quad
	Since \(\phi(t)=0\) for \(t\geq1\), by \eqref{0401-09}, \(\partial_{x_n}w_1(x',0,t)\) is equal to its strictly positive integral term.

	\underline{Proof of (ii).}\quad
	Since \(\phi(t)\geq0\) and \(\mathcal S(x')\leq0\), \eqref{0401-09} gives \(\partial_{x_n}w_1(x',0,t)>0\).

	\underline{Proof of (iii).}\quad
	When \(n=2\), \eqref{260063} and \eqref{260064} give, for \(|x_1|>1\), \(\mathcal S(x_1)=-\frac1\pi\int_{\R}\frac{\psi(y_1)}{(x_1-y_1)^2}\,dy_1<0\), since \(\psi\geq0\) is continuous and not identically zero.
	Then {\rm(ii)} gives the result.

	Assume now that \(\psi\) is radial.
	For \(|x'|>2\), we claim that \(\mathcal{S}\) satisfies the asymptotic expansions
	\begin{align}
		 & \left|\mathcal S(x')-\frac{\Gamma\left(\frac n2\right)M_0 \bigl(|x''|^2-(n-1)x_1^2\bigr)} {\pi^{\frac n2}|x'|^{n+2}}\right| \leq C_nM_2|x'|^{-n-2},\label{260065}\\
		 & \mathcal S(x')=\frac{\Gamma\left(\frac n2\right)}{\pi^{\frac n2}} \left[\frac{M_0\bigl(|x''|^2-(n-1)x_1^2\bigr)-\frac{M_2}{n-1}}{|x'|^{n+2}}+\frac{\frac{(n+2)M_2}{2(n-1)}\bigl(|x''|^2-(n-1)x_1^2\bigr)- \frac{3(n+2)M_4}{2(n^2-1)}}{|x'|^{n+4}}\right.\nonumber\\
&\left.\quad \quad \quad \quad + \frac{3(n+2)(n+4)M_4}{8(n^2-1)}\frac{|x''|^2-(n-1)x_1^2}{|x'|^{n+6}}\right] +O\bigl(|x'|^{-n-6}\bigr).\label{260066}
	\end{align}
	Indeed, the radiality of \(\psi\) gives
	\begin{align*}
		\int_\Sigma y_jy_k\psi(y')\,dy'=\frac{M_2}{n-1}\delta_{jk},\qquad \int_\Sigma y_jy_ky_ly_m\psi(y')\,dy' =\frac{M_4}{n^2-1} (\delta_{jk}\delta_{lm}+\delta_{jl}\delta_{km}+\delta_{jm}\delta_{kl}).
	\end{align*}
	Radiality also makes every odd moment vanish. Taylor's theorem, with terms through degree five, therefore gives
	\begin{align}\label{0916-01}
		\int_\Sigma K_{\mathcal S}(x'-y')\psi(y')\,dy' =M_0K_{\mathcal S}(x')+\frac{M_2}{2(n-1)}\Delta'K_{\mathcal S}(x')+ \frac{M_4}{8(n^2-1)}(\Delta')^2K_{\mathcal S}(x')+O(|x'|^{-n-6}).
	\end{align}
	Direct differentiation of \eqref{260063} gives
	\begin{align*}
		K_{\mathcal S}(x')=\frac{|x''|^2-(n-1)x_1^2}{|x'|^{n+2}},\quad \Delta'K_{\mathcal S}(x')=-\frac{2}{|x'|^{n+2}} +\frac{(n+2)\bigl(|x''|^2-(n-1)x_1^2\bigr)}{|x'|^{n+4}},
	\end{align*}
	and
	\begin{align*}
		(\Delta')^2K_{\mathcal S}(x')=-\frac{12(n+2)}{|x'|^{n+4}} +\frac{3(n+2)(n+4)\bigl(|x''|^2-(n-1)x_1^2\bigr)}{|x'|^{n+6}}.
	\end{align*}
	Substituting these identities into \eqref{0916-01} gives \eqref{260066}.
	Expanding to first order gives \eqref{260065}.
	Let \(R(x')\) be the remainder in \eqref{0916-01}. Taylor's theorem also applies to each first derivative of \(K_{\mathcal S}\).
	The derivatives of \(K_{\mathcal S}\) of order \(k\) are homogeneous of degree \(-n-k\), and \(|x'-v y'|\geq |x'|/2\) for \(|x'|>2\), \(|y'|\leq1\), and \(0\leq v\leq1\).
	Thus \(R\) and its gradient satisfy
	\begin{align}\label{0916-02}
		|R(x')|+|x'|\,|\nabla'R(x')|\leq C_n M_6 |x'|^{-n-6},\qquad |x'|>2.
	\end{align}

	\underline{Proof of (iv).}\quad
	By radiality, the Riesz potential \((-\Delta')^{-1/2}\psi(x')=u(r)\) depends only on \(r=|x'|\).
	For \(r>1\), its integral representation is
	\begin{align*}
		u(r)=c_n\int_\Sigma |x'-y'|^{-(n-2)}\psi(y')\,dy',\qquad
		c_n=\frac{\Gamma(n/2)}{(n-2)\pi^{n/2}}.
	\end{align*}
	At \(x'=r\mathbf e_1'\), the distance to \(\supp\psi\) is at least \(r-1>0\), so we may differentiate under the integral.
	The radial derivative and the tangential Laplacian are
	\begin{align}
		u'(r)&=-c_n(n-2)\int_\Sigma\frac{(r-y_1)\psi(y')}{|r\mathbf e_1'-y'|^n}\,dy'<0,\nonumber\\
		u''(r)+\frac{n-2}{r}u'(r)&=c_n(n-2)\int_\Sigma\frac{\psi(y')}{|r\mathbf e_1'-y'|^n}\,dy'>0.\label{riesz-laplacian}
	\end{align}
	Both signs are strict because \(r-y_1>0\) on \(\supp\psi\) and \(\psi\) is positive on a nonempty open set.
	Hence \(u''(r)>0\).
	Using \(R_1'=-\partial_{x_1}(-\Delta')^{-1/2}\) in \eqref{260009}, we obtain
	\begin{align*}
		\mathcal S(\pm r\mathbf e_1')=-u''(r)<0,\qquad r>1.
	\end{align*}
	Part {\rm(ii)} therefore excludes separation points on the \(x_1\)-axis for \(|x_1|>1\).

	If \(|x''|\leq(\sqrt{n-1}-\varepsilon)|x_1|\), then \(|x''|^2-(n-1)x_1^2\leq -(2\varepsilon\sqrt{n-1}-\varepsilon^2)x_1^2\).
	In this region \(x_1^2\geq c_{n,\varepsilon}|x'|^2\), so \eqref{260065} gives \(\mathcal S(x')<0\) once \(|x'|\geq r_\varepsilon\) for a sufficiently large \(r_\varepsilon\).
	Applying (ii) again gives (iv).

	\underline{Proof of (v).}\quad
	The set \(\{\psi>0\}\) is open and nonempty, and therefore \(M_2^2<M_0M_4\) by the strict Cauchy--Schwarz inequality.
	It follows that
	\begin{align*}
		\frac{(n+2)M_2^2}{2(n-1)^2M_0}-\frac{3(n+2)M_4}{2(n^2-1)} <\frac{(n+2)(2-n)M_4}{(n-1)^2(n+1)}<0.
	\end{align*}
	Choose \(\lambda>\frac{2\pi^{\frac n2}C_nM_2}{\Gamma(\frac n2)M_0} +\frac{M_2}{(n-1)M_0}+1\), where \(C_n\) is the constant in \eqref{260065}.
	If \(|x''|^2-(n-1)x_1^2\leq-\lambda\), then \eqref{260065} gives
	\begin{align*}
		\mathcal S(x') \leq \frac{\Gamma\left(\frac n2\right)M_0}{\pi^{\frac n2}|x'|^{n+2}} \left(|x''|^2-(n-1)x_1^2+\frac{\pi^{\frac n2}C_nM_2} {\Gamma\left(\frac n2\right)M_0}\right)<0.
	\end{align*}
	If \(-\lambda\leq|x''|^2-(n-1)x_1^2\leq\frac{M_2}{(n-1)M_0}\), then the first term in \eqref{260066} is nonpositive and the coefficient of \(|x'|^{-n-4}\) is at most its strictly negative value at \(|x''|^2-(n-1)x_1^2=\frac{M_2}{(n-1)M_0}\).
	Hence \(\mathcal S(x')<0\) provided \(|x_1|>r_0(n,\psi)\).
	At a separation point with \(0<t<1\), \eqref{0401-09} and the strict positivity of its integral term give \(\mathcal S(x')>0\), and therefore \(|x''|^2-(n-1)x_1^2>\frac{M_2}{(n-1)M_0}\).
\end{proof}

The next lemma gives the sign of \(\mathcal S\) on \(x_1=0\).
\begin{lemma}\label{lem-260067}
	Let \(n\geq3\), let \(\psi\) be radial and satisfy {\rm\textbf{(S)}}, and let \(\mathcal S\) be given in \eqref{260009} with \(i=1\).
	Then
	\begin{align}\label{260068}
		\mathcal S(r\mathbf e_2')>0\qquad\text{for every }r>1.
	\end{align}
\end{lemma}
\begin{proof}
	Let \(u\) be the radial Riesz potential in the proof of Proposition~\ref{prop-260061}{\rm(iv)}.
	Since \(\mathcal S=-\partial_{x_1}^2u\), the formula for the Hessian of a radial function gives, for \(r>1\) and \(\omega'\in\mathbb S^{n-2}\),
	\begin{align*}
		\mathcal S(r\omega')=-u''(r)(\omega'_1)^2-\frac{u'(r)}{r}\bigl(1-(\omega'_1)^2\bigr).
	\end{align*}
	On the hyperplane \(x_1=0\), this reduces to \(-u'(r)/r\).
	The same proposition gives \(u'(r)<0\), which proves \eqref{260068}.
\end{proof}

\subsubsection{Existence and uniqueness of separation points}\label{260071}
The next proposition gives separation points of \(w_1\) where \(\mathcal S(x')>0\) and, for radial \(\psi\), on \(x_1=0\) and at \((x_1,\rho_*(x_1,t)\omega'',t)\) for large \(|x_1|\).
This proves Theorem \ref{thm-260010}{\rm(ii)}.
\begin{proposition}\label{prop-260072}
	Let \(n\geq3\) and let \(w\) solve the Stokes system \eqref{0308-01} with boundary data \eqref{0324-13}, where \(j=1\), \(\psi\) satisfies {\rm\textbf{(S)}}, and \(\phi\) satisfies {\rm\textbf{(T1)}} and {\rm\textbf{(T2)}}.
	Let \(\mathcal S\) be given in \eqref{260009} with \(i=1\).
	Then the following hold.
	\begin{enumerate}
		\item[\rm(i)] If \(x'\in\Sigma\), \(|x'|>1\), and \(\mathcal S(x')>0\), then \(T_{\rm van}(x')\) is a nonempty compact subset of \((0,1)\), and there exists \(t_1^*(x')\in T_{\rm van}(x')\) such that \((x',t_1^*(x'))\) is a separation point of \(w_1\).
		Moreover,
		\begin{align*}
			\partial_{x_n}w_1(x',0,t)<0 \ \ \text{for }0<t<\min T_{\rm van}(x'), \qquad \partial_{x_n}w_1(x',0,t)>0 \ \ \text{for }t>\max T_{\rm van}(x').
		\end{align*}
		\item[\rm(ii)] If \(\psi\) is radial, then there exists \(R_0=R_0(n,\psi)>2\) such that,
for \(0<t_0<1\), every \(|x_1|>R_{t_0}\), and every \(0<t\leq t_0\), where
		\begin{align}\label{260073}
			R_{t_0}:=R_0+\frac2{\sqrt n} \sqrt{\ln\frac{\phi(t_1)}{\phi(\max\{t_0,t_1\})}},
		\end{align}
		 there is \(\rho_*(x_1,t)>0\), independent of \(\omega''\in\mathbb S^{n-3}\), such that
		\begin{align*}
			\partial_{x_n}w_1(x_1,\rho_*(x_1,t)\omega'',0,t)=0, \qquad \left|\rho_*(x_1,t)^2-(n-1)x_1^2-\frac{M_2}{(n-1)M_0}\right|\leq1.
		\end{align*}
		The point \((x_1,\rho_*(x_1,t)\omega'',t)\) is a separation point of \(w_1\).
		\item[\rm(iii)] If \(\psi\) is radial and \(x''\in(\mathbf e_1')^\perp\) satisfies \(|x''|>1\), then \(\mathcal S(0,x'')>0\) and \(T_{\rm van}(0,x'')\neq\emptyset\).
		There exists \(t_1^*(0,x'')\in T_{\rm van}(0,x'')\), depending only on \(|x''|\), such that \((0,x'',t_1^*(0,x''))\) is a separation point of \(w_1\).
		Every \(t_1^*(0,x'')\in T_{\rm van}(0,x'')\) satisfies
		\begin{align}\label{0708-11}
			t_1^*(0,x'')\to1^- \qquad\text{as }|x''|\to\infty.
		\end{align}
	\end{enumerate}
\end{proposition}

\begin{proof}
	\underline{Proof of (i).}\quad
	Fix \(x'\in\Sigma\) with \(|x'|>1\) and \(\mathcal S(x')>0\).
	For \(t\geq1\), \(\phi(t)=0\), so \eqref{0401-09} gives \(\partial_{x_n}w_1(x',0,t)>0\), and, by continuity of the integral in \eqref{0401-09}, also for \(t<1\) close to \(1\).
	For \(0<t\leq t_1\), \(\phi(s)\leq\phi(t)\) for \(0<s<t\), and, with \(d:=\dist(x',\supp\psi)>0\), the change of variables \(u=d^2/(4(t-s))\) gives
	\begin{align}\label{260074}
		0< & \frac{1}{\sqrt{4\pi}}\int_0^t\frac{\phi(s)}{(t-s)^{3/2}} \int_\Sigma\Gamma'(x'-y',t-s)\psi(y')\,dy'ds\leq C_nM_0\phi(t)d^{-n}e^{-\frac{d^2}{8t}} =o\bigl(\phi(t)\bigr)\qquad\text{as }t\to0^+.
	\end{align}
	Thus \eqref{0401-09} is negative for all sufficiently small \(t>0\).
	Hence \(T_{\rm van}(x')\) is a nonempty compact subset of \((0,1)\), and continuity gives the two sign inequalities in {\rm(i)}.

	Let \(I\) be the maximal open interval containing \(x_1\) on which \(r^2+|x''|^2>1\).
	For \(0<t<1\), \(\partial_{x_n}w_1(r,x'',0,t)\) is real-analytic in \(r\in I\): Lemma \ref{lem-260044}{\rm(ii)}, with \(T=t\), \(q(s)=\phi(t-s)\), \(b=\frac{n+2}{2}\), \(m=0\), and \(g=\psi\), applies to the first term in \eqref{0401-09}, and Lemma \ref{lem-260044}{\rm(i)}, with \(m=1\), \(a=\frac n2\), and \(g=\psi\), applies to \(R_1'\psi(\cdot,x'')\), whose derivative is \(\mathcal S(\cdot,x'')\).
	For \(r\in I\) with \(|r|>1+\frac{|x''|+1}{\sqrt{n-1}}\), \eqref{260063}--\eqref{260064} give
	\begin{align}\label{1007-02}
		\mathcal S(r,x'') & =\frac{\Gamma(\frac n2)}{\pi^{\frac n2}} \int_\Sigma \frac{|x''-y''|^2-(n-1)(r-y_1)^2}{|(r,x'')-y'|^{n+2}} \psi(y')\,dy'\nonumber\\
		 & \leq\frac{\Gamma(\frac n2)}{\pi^{\frac n2}} \bigl[(|x''|+1)^2-(n-1)(|r|-1)^2\bigr] \int_\Sigma\frac{\psi(y')}{|(r,x'')-y'|^{n+2}}\,dy'<0,
	\end{align}
	so \(\partial_{x_n}w_1(r,x'',0,t)>0\) there, and Lemma \ref{lem-260044}{\rm(iii)} excludes an interval of zeros in \(I\).

	Choose \(0<t_-<t_+<1\) such that \(\partial_{x_n}w_1(x',0,t_-)<0<\partial_{x_n}w_1(x',0,t_+)\).
	The joint continuity of \(\partial_{x_n}w_1\) in \eqref{260261} with \(f=\psi\) and the connectedness argument in the proof of Proposition \ref{prop-0504-03}, applied on \([t_-,t_+]\), give \(t_1^*(x')\in(t_-,t_+)\) such that \(\partial_{x_n}w_1(r,x'',0,s)\) takes both signs in every neighborhood of \((x_1,t_1^*(x'))\) in \(I\times(0,1)\).
	Then \(\partial_{x_n}w_1(x',0,t_1^*(x'))=0\), and Proposition \ref{prop-260046}{\rm(i)} shows that \((x',t_1^*(x'))\) is a separation point.

	\underline{Proof of (ii).}\quad
	By {\rm\textbf{(T2)}},
	\begin{align}\label{260075}
		\frac{\phi(s)}{\phi(t)}\leq \frac{\phi(t_1)}{\phi(\max\{t_0,t_1\})} \qquad(0<s<t\leq t_0<1),
	\end{align}
	since \(\phi(s)\leq\phi(t)\) when \(t\leq t_1\), and \(\phi(s)\leq\phi(t_1)\) and \(\phi(t)\geq\phi(t_0)\) when \(t_1\leq t\leq t_0\).
	For each fixed \(\lambda>0\), \eqref{260066} gives, uniformly for \(\bigl||x''|^2-(n-1)x_1^2\bigr|\leq\lambda\) and sufficiently large \(|x_1|\),
	\begin{align}\label{260076}
		\mathcal S(x')=\frac{\Gamma\left(\frac n2\right) \left[M_0\bigl(|x''|^2-(n-1)x_1^2\bigr)-\frac{M_2}{n-1}\right]} {\pi^{\frac n2}|x'|^{n+2}} +O_\lambda\bigl(|x'|^{-n-4}\bigr).
	\end{align}

	Write \(r=|x'|>2\).
	Since \(|x'-y'|\geq r-1>1\) on \(\supp\psi\) and \(\int_0^1 s^{-(n+2)/2}e^{-\frac{(r-1)^2}{4}(1/s-1)}\,ds\leq C_n\), the integral in \eqref{0401-09} gives, uniformly for \(0<t<1\),
	\begin{align}\label{260077}
		0<\partial_{x_n}w_1(x',0,t)+\phi(t)\mathcal S(x')\leq Ce^{-c|x'|^2},
	\end{align}
	and, by \eqref{260075},
	\begin{align}\label{260078}
		0<\frac{\partial_{x_n}w_1(x',0,t)}{\phi(t)}+\mathcal S(x') \leq C_nM_0\frac{\phi(t_1)}{\phi(\max\{t_0,t_1\})} e^{-\frac{(r-1)^2}{4}} \qquad(0<t\leq t_0).
	\end{align}
	By \eqref{260065}, there is \(c=c(n,\psi)\geq1\) such that \(|x''|^2-(n-1)x_1^2\geq-c\), and hence \(r\geq\sqrt n\,|x_1|-\sqrt c\), whenever \(\mathcal S(x')>0\) or \(\left||x''|^2-(n-1)x_1^2-\frac{M_2}{(n-1)M_0}\right|\leq1\).
	Choose \(R_0\) so that \(\sqrt nR_0-\sqrt c-1\geq4\).
	For \(|x_1|>R_{t_0}\), \eqref{260073} then gives
	\begin{align}\label{260079}
		\ln\frac{\phi(t_1)}{\phi(\max\{t_0,t_1\})} -\frac{(r-1)^2}{4}\leq-2r+6,
	\end{align}
	and hence, in these two regions, uniformly for \(0<t\leq t_0\),
	\begin{align}\label{260080}
		0<\frac{\partial_{x_n}w_1(x',0,t)}{\phi(t)}+\mathcal S(x') \leq C_nM_0e^{-r}.
	\end{align}
	If \(|x''|^2-(n-1)x_1^2=\frac{M_2}{(n-1)M_0}-1\), then \eqref{260076} gives \(\mathcal S(x')<0\), so \(\partial_{x_n}w_1(x',0,t)>0\).
	If \(|x''|^2-(n-1)x_1^2=\frac{M_2}{(n-1)M_0}+1\), then \eqref{260076} and \eqref{260080} give \(\partial_{x_n}w_1(x',0,t)<0\) after increasing \(R_0\).

	To obtain a sign-changing zero, rotate the \(x''\)-coordinates so that \(x''=|x''|\mathbf e_2'\), and keep \(x_1\) fixed.
	The first term on the right-hand side of \eqref{0401-09} is real-analytic in \(x_2\) by Lemma \ref{lem-260044}{\rm(ii)} with \(m=0\).
	By \eqref{260063}--\eqref{260064}, \(\mathcal S\) involves the kernels \(|x'-y'|^{-n}\) and \((x_1-y_1)^2|x'-y'|^{-n-2}\), which are real-analytic in \(x_2\) by Lemma \ref{lem-260044}{\rm(i)} with \(m=0\), exponents \(n/2\) and \((n+2)/2\), and weights \(\psi(y')\) and \((x_1-y_1)^2\psi(y')\), since \(|x_1|>1\).
	The difference in \eqref{0401-09} is not identically zero, so Lemma \ref{lem-260044}{\rm(iii)} gives finitely many zeros between the two radii above.
	As the signs at these radii differ, \(\partial_{x_n}w_1\) takes both signs in every neighborhood of one of these zeros, \(\rho_*(x_1,t)\), and Proposition \ref{prop-260046}{\rm(i)} shows that \((x_1,\rho_*(x_1,t)\omega'',t)\) is a separation point.
	Since \(\psi\) is radial, \(\rho_*(x_1,t)\) is independent of \(\omega''\).

	\underline{Proof of (iii).}\quad
	Take \(x'=(0,x'')\).
	By \eqref{260068}, \(\mathcal S(x')>0\) for \(|x''|>1\), so {\rm(i)} gives a separation point \((x',t_1^*(0,x''))\).
	Since \(\psi\) is radial, \eqref{0401-09} is invariant under orthogonal transformations of the last \(n-2\) coordinates, so \(T_{\rm van}(x')\) depends only on \(|x''|\), and Proposition \ref{prop-260046}{\rm(i)} allows \(t_1^*(0,x'')\) to depend only on \(|x''|\).
	Since \(|y'|\leq1\) on \(\supp\psi\), \eqref{260063}--\eqref{260064} give, as \(|x''|\to\infty\),
	\begin{align}\label{0706-07}
		\mathcal S(x')=\frac{\Gamma\left(\frac n2\right)M_0}{\pi^{\frac n2}} |x''|^{-n}\bigl(1+O(|x''|^{-1})\bigr), \qquad c_-|x''|^{-n}\leq\mathcal S(x')\leq c_+|x''|^{-n},
	\end{align}
	where \(c_-,c_+>0\) depend only on \(n\) and \(\psi\).
	At every \(t_1^*(0,x'')\in T_{\rm van}(x')\), \eqref{0401-09} gives
	\begin{align}\label{0706-06}
		 & \frac{1}{\sqrt{4\pi}}\int_0^{t_1^*(0,x'')} \frac{\phi(s)}{(t_1^*(0,x'')-s)^{3/2}} \int_\Sigma\Gamma'(x'-y',t_1^*(0,x'')-s)\psi(y')\,dy'ds =\phi(t_1^*(0,x''))\mathcal S(x'),
	\end{align}
	so \eqref{260077} and \eqref{0706-07} give \(\phi(t_1^*(0,x''))\leq Ce^{-c|x''|^2}\), and \(t_1^*(0,x'')\) can accumulate only at \(0\) or \(1\).
	For \(0<t\leq t_1\), \eqref{0401-09}, \eqref{0706-07}, and \eqref{260074} give
	\begin{align*}
		\partial_{x_n}w_1(x',0,t) \leq C\phi(t)|x''|^{-n}e^{-\frac{(|x''|-1)^2}{8t}} -c_-\phi(t)|x''|^{-n}<0
	\end{align*}
	for all sufficiently large \(|x''|\), uniformly for \(0<t\leq t_1\), which excludes accumulation at \(0\) and proves \eqref{0708-11}.
\end{proof}

The next proposition gives uniqueness of \(t_1^*(x')\) for \(x'\) such that \(\mathcal S>0\), and uniqueness of \(\rho_*(x_1,t)\) for large \(|x_1|\) and fixed \(t\in(0,1)\).
This proves Theorem \ref{thm-260010}{\rm(iii)}.
\begin{proposition}\label{prop-260081}
	Let \(n\geq3\) and let \(w\) solve the Stokes system \eqref{0308-01} with boundary data \eqref{0324-13}, where \(j=1\), \(\psi\) satisfies {\rm\textbf{(S)}}, and \(\phi\) satisfies {\rm\textbf{(T1)}} and {\rm\textbf{(T2)}}.
	Let \(\mathcal S\) be given in \eqref{260009} with \(i=1\).
	Then the following hold.
	\begin{enumerate}
		\item[\rm(i)] Suppose that \(x'\in\Sigma\), \(|x'|>1\), \(\mathcal S(x')>0\), and \(\phi\) is log-concave on \((0,1)\).
		Then \(T_{\rm van}(x')=\{t_1^*(x')\}\) and
		\begin{align*}
			\partial_{x_n}w_1(x',0,t)
			\begin{cases}
				<0, & 0<t<t_1^*(x'), \\[2pt]
				=0, & t=t_1^*(x'),\\[2pt]
				>0, & t>t_1^*(x').
			\end{cases}
		\end{align*}
		The point \((x',t_1^*(x'))\) is a separation point of \(w_1\).
		\item[\rm(ii)] Suppose that \(\psi\) is radial.
	Then, after increasing \(R_0\) in \eqref{260073} if necessary, for \(0<t_0<1\), every \(|x_1|>R_{t_0}\), and every \(0<t\leq t_0\), there is a unique \(\rho_*(x_1,t)>0\) such that \(\partial_{x_n}w_1(x_1,\rho_*(x_1,t)\omega'',0,t)=0\), and
		\begin{align*}
			\sgn\bigl(\partial_{x_n}w_1(x',0,t)\bigr) =\sgn\bigl(\rho_*(x_1,t)-|x''|\bigr).
		\end{align*}
		The point \((x_1,\rho_*(x_1,t)\omega'',t)\) is a separation point of \(w_1\).
		\item[\rm(iii)] Suppose that \(\psi\) is radial and \(\phi\) is log-concave on \((0,1)\).
		Then, for \(x''\in(\mathbf e_1')^\perp\) with \(|x''|>1\), \(T_{\rm van}(0,x'')\) is a singleton \(\{t_1^*(0,x'')\}\), and \((0,x'',t_1^*(0,x''))\) is a separation point of \(w_1\).
	\end{enumerate}
\end{proposition}

\begin{proof}
	\underline{Proof of (i).}\quad
Let \(d=\operatorname{dist}(x',\supp\psi)>0\).
After substituting \(s\mapsto t-s\) in \eqref{0401-09}, we see that the spatial convolution multiplied by \(s^{-3/2}\) is bounded by \(C_nM_0 s^{-(n+2)/2}e^{-d^2/(4s)}\), which is integrable at zero.
Since \(\phi(0)=0\), differentiation gives
	\begin{align*}
		 & \partial_t\left[\frac{1}{\sqrt{4\pi}\,\phi(t)} \int_0^t\frac{\phi(s)}{(t-s)^{3/2}} \int_\Sigma\Gamma'(x'-y',t-s)\psi(y')\,dy'ds\right]\nonumber \\
		 & \quad=\frac{1}{\sqrt{4\pi}\,\phi(t)}\int_0^t\phi(t-s) \left[\frac{\phi'(t-s)}{\phi(t-s)}-\frac{\phi'(t)}{\phi(t)}\right] \frac{1}{s^{3/2}}\int_\Sigma\Gamma'(x'-y',s)\psi(y')\,dy'ds>0.
	\end{align*}
	Since \(\phi\) is log-concave, the difference of logarithmic derivatives is nonnegative.
	For each \(t\), it is positive when \(t-s>0\) is sufficiently small.
The convolution \(\Gamma'(\cdot,s)*_\Sigma\psi\) is positive for every \(s>0\) by {\rm\textbf{(S)}}, so the integral is strictly positive.

	Therefore the expression in square brackets on the left is strictly increasing, and \eqref{0401-09} has exactly one zero in \((0,1)\) by Proposition \ref{prop-260072}{\rm(i)}, with the stated signs.
	By Proposition \ref{prop-260046}{\rm(ii)}, this zero is a separation point of \(w_1\).

	\underline{Proof of (ii).}\quad
	Since \(K_{\mathcal S}\) is homogeneous of degree \(-n\), \(|\nabla'^kK_{\mathcal S}(x'-v y')|\leq C_{n,k}|x'|^{-n-k}\) for \(0\leq k\leq7\), \(|y'|\leq1\), \(0\leq v\leq1\), and sufficiently large \(|x_1|\).
	Hence a Taylor expansion applied to \(\nabla'K_{\mathcal S}\) allows differentiation of the \(O(|x'|^{-n-6})\) term in \eqref{260066}.
	The identities
	\begin{align*}
		\frac{x''\cdot\nabla_{x'}}{2|x''|^2} \bigl(|x''|^2-(n-1)x_1^2\bigr)=1, \qquad \frac{x''\cdot\nabla_{x'}}{2|x''|^2}|x'|^{-m} =-\frac m2|x'|^{-m-2},
	\end{align*}
	give
	\begin{align}\label{260082}
		\frac{x''\cdot\nabla_{x'}\mathcal S(x')}{2|x''|^2} =\frac{\Gamma\left(\frac n2\right)M_0}{\pi^{\frac n2}} \left[\frac1{|x'|^{n+2}}- \frac{(n+2)\left(|x''|^2-(n-1)x_1^2-\frac{2M_2}{(n-1)M_0}\right)} {2|x'|^{n+4}}\right] +O_\lambda\bigl(|x'|^{-n-6}\bigr)
	\end{align}
	uniformly when \(\bigl||x''|^2-(n-1)x_1^2\bigr|\leq\lambda\).
	Directly from the integral in \eqref{0401-09}, uniformly for \(|x'|>2\) and \(0<t<1\),
	\begin{align}\label{260083}
		\left|\nabla_{x'}\bigl(\partial_{x_n}w_1(x',0,t)+\phi(t)\mathcal S(x')\bigr)\right| \leq Ce^{-c|x'|^2},
	\end{align}
	where we used \(|\nabla\Gamma'(z',\tau)|\leq C\tau^{-1/2}\Gamma'(z',2\tau)\).
Direct differentiation of the heat kernel and \eqref{260075} give, for \(r=|x'|>2\),
	\begin{align*}
		\left|\nabla_{x'}\left( \frac{\partial_{x_n}w_1(x',0,t)}{\phi(t)}+\mathcal S(x')\right)\right| \leq C_nM_0(1+r)\frac{\phi(t_1)}{\phi(\max\{t_0,t_1\})} e^{-\frac{(r-1)^2}{4}}.
	\end{align*}
	Indeed, the additional factor in the integral is at most \((r+1)/(2s)\), and \(\int_0^1s^{-(n+4)/2}e^{-\frac{(r-1)^2}{4}(1/s-1)}\,ds\leq C_n\).
	Thus \eqref{260079} gives
	\begin{align}\label{260084}
		\left|\nabla_{x'}\left( \frac{\partial_{x_n}w_1(x',0,t)}{\phi(t)}+\mathcal S(x')\right)\right| \leq C_nM_0e^{-|x'|}
	\end{align}
	when \(|x_1|>R_{t_0}\), \(0<t\leq t_0\), and either \(\mathcal S(x')>0\) or \(\left||x''|^2-(n-1)x_1^2-\frac{M_2}{(n-1)M_0}\right|\leq1\).

	In the latter region, \eqref{260082} and \eqref{260084} give, after increasing \(R_0\),
	\begin{align}\label{transversality}
		\frac{x''}{|x''|}\cdot\nabla_{x'} \frac{\partial_{x_n}w_1(x',0,t)}{\phi(t)} \leq-\frac{c|x''|}{|x'|^{n+2}}<0,
	\end{align}
	where \(c>0\) depends only on \(n\) and \(\psi\).
	Indeed,
	\begin{align*}
		\frac{x''\cdot\nabla_{x'}}{2|x''|^2} \frac{\partial_{x_n}w_1(x',0,t)}{\phi(t)} \leq Ce^{-|x'|} -\frac{\Gamma(\frac n2)M_0}{2\pi^{\frac n2}}|x'|^{-n-2} \leq-\frac{\Gamma(\frac n2)M_0}{4\pi^{\frac n2}}|x'|^{-n-2},
	\end{align*}
	which gives \eqref{transversality} after multiplying by \(2|x''|\).

	The zero obtained in Proposition \ref{prop-260072}{\rm(ii)} lies in the region where \eqref{transversality} holds and is therefore unique there.

	At any zero, \eqref{0401-09} and \eqref{260080} give \(0<\mathcal S(x')\leq Ce^{-|x'|}\).
	Combining this bound with \eqref{260065} first shows that \(\bigl||x''|^2-(n-1)x_1^2\bigr|\) is bounded independently of \(x_1\) and \(t\).
	Applying \eqref{260076} then gives \(\left||x''|^2-(n-1)x_1^2-\frac{M_2}{(n-1)M_0}\right|\leq C|x'|^{-2}<1\) after increasing \(R_0\).

	Hence every zero lies in the region where \eqref{transversality} gives uniqueness, so the zero already obtained is the only one.
	The strict sign change at this zero and Proposition \ref{prop-260046}{\rm(i)} prove that \((x_1, \rho_*(x_1,t)\omega'',t)\) is a separation point.
	Continuity and uniqueness of the zero give \(\partial_{x_n}w_1(x',0,t)>0\) for \(0<|x''|<\rho_*(x_1,t)\) and \(\partial_{x_n}w_1(x',0,t)<0\) for \(|x''|>\rho_*(x_1,t)\).

	\underline{Proof of (iii).}\quad
	Since \(\mathcal S(0,x'')>0\) for \(|x''|>1\) by Proposition \ref{prop-260072}{\rm(iii)}, the conclusion follows from (i).
\end{proof}

\subsubsection{Asymptotics and monotonicity of the separation point}\label{260085}
The next proposition gives the asymptotic expansion of \(\rho_*(x_1,t)\) as \(|x_1|\to\infty\) and the rate at which \(t_1^*(0,x'')\) approaches \(1\) as \(|x''|\to\infty\).
This proves the asymptotics in Theorem \ref{thm-260010}{\rm(iv)}.

\begin{proposition}\label{prop-260086}
	Let \(n\geq3\) and let \(w\) solve the Stokes system \eqref{0308-01} with boundary data \eqref{0324-13}, where \(j=1\), \(\psi\) is radial and satisfies {\rm\textbf{(S)}}, and \(\phi\) satisfies {\rm\textbf{(T1)}} and {\rm\textbf{(T2)}}.
	Then the following hold.
	\begin{enumerate}
		\item[\rm(i)] After increasing \(R_0\) in \eqref{260073} if necessary, the unique positive zero \(\rho_*(x_1,t)\) of \(|x''|\mapsto\partial_{x_n}w_1(x_1,|x''|\omega'',0,t)\) satisfies, uniformly for \(|x_1|>R_{t_0}\) and \(0<t\leq t_0<1\),
		\begin{align}\label{260087}
			\rho_*(x_1,t)=\sqrt{n-1}\,|x_1| +\frac{M_2}{2(n-1)^{3/2}M_0|x_1|}+O\bigl(|x_1|^{-3}\bigr).
		\end{align}
		\item[\rm(ii)] Suppose that \(\phi\in C^N([0,1])\) for some \(N\geq1\) and that some derivative of \(\phi\) of order at most \(N\) is nonzero at each of \(0\) and \(1\).
		For \(|x''|\) sufficiently large, let \(t_1^*(0,x'')\in T_{\rm van}(0,x'')\), and let \(k_0,k_1\leq N\) be the smallest positive integers such that \(\phi^{(k_0)}(0)\neq0\) and \(\phi^{(k_1)}(1)\neq0\).
		Then
		\begin{align}\label{0706-12}
			\ln\bigl(1-t_1^*(0,x'')\bigr) =-\frac{|x''|^2}{4k_1}+O(|x''|) \qquad\text{as }|x''|\to\infty.
		\end{align}
	\end{enumerate}
\end{proposition}

\begin{proof}
	\underline{Proof of (i).}\quad
		Fix \(\omega''\in\mathbb S^{n-3}\) and set \(x'=(x_1,\rho_*(x_1,t)\omega'')\).
	At the zero of \(\partial_{x_n}w_1(x',0,t)\), \eqref{0401-09} and \eqref{260080} give \(0<\mathcal S(x')\leq Ce^{-|x'|}\).
	The estimate \eqref{260065} then shows that \(|\rho_*^2-(n-1)x_1^2|\) is bounded.
	Proposition \ref{prop-260061}{\rm(v)} gives \(\rho_*^2-(n-1)x_1^2>\frac{M_2}{(n-1)M_0}\), and \eqref{260076} gives \(0<\rho_*^2-(n-1)x_1^2-\frac{M_2}{(n-1)M_0}\leq C|x'|^{-2}\).
	Since \(|x'|^2=nx_1^2+\rho_*^2-(n-1)x_1^2=nx_1^2+O(1)\),
	\begin{align*}
		\rho_*=\sqrt{(n-1)x_1^2+\frac{M_2}{(n-1)M_0}+O(x_1^{-2})} =\sqrt{n-1}\,|x_1| \left(1+\frac{M_2}{(n-1)^2M_0x_1^2}+O(x_1^{-4})\right)^{1/2}.
	\end{align*}
	This gives \eqref{260087}.

	\underline{Proof of (ii).}\quad
The limit \eqref{0708-11} gives \(t_1^*(0,x'')\to1^-\).
	Since \(\phi(1)=0\) and \(k_1\) is the order of the first nonvanishing derivative of \(\phi\) at \(1\), Taylor's theorem gives \(\phi(t)=\frac{|\phi^{(k_1)}(1)|}{k_1!}(1-t)^{k_1}(1+o(1))\) as \(t\to1^-\).

	We estimate the integral in \eqref{0706-06}.
	For the upper bound, choose \(|x''|\) sufficiently large so that for \(y'\in\supp\psi\), \(|(0,x'')-y'|^2\geq(|x''|-1)^2\geq2(n+2)\).
	Thus \(\tau^{-\frac{n+2}{2}}e^{-\frac{|(0,x'')-y'|^2}{4\tau}}\) is increasing in \(\tau \in (0,1]\).
	Since \(0<t_1^*(0,x'')-s\leq1\),
	\begin{align*}
		(t_1^*(0,x'')-s)^{-\frac{n+2}{2}} e^{-\frac{|(0,x'')-y'|^2}{4(t_1^*(0,x'')-s)}} \leq e^{-\frac{|(0,x'')-y'|^2}{4}}\leq Ce^{-\frac{|x''|^2}{4}+\frac{|x''|}{2}}.
	\end{align*}
	Substituting this bound into the left-hand side of \eqref{0706-06}, we obtain
	\begin{align}\label{0706-06a}
		\frac{1}{\sqrt{4\pi}}\int_0^{t_1^*(0,x'')} \frac{\phi(s)}{(t_1^*(0,x'')-s)^{3/2}} \int_\Sigma\Gamma'((0,x'')-y',t_1^*(0,x'')-s)\psi(y')\,dy'ds \leq CM_0\norm{\phi}_{L^1(0,1)}e^{-\frac{|x''|^2}{4}+\frac{|x''|}{2}}.
	\end{align}

		For the lower bound, the estimate \(\phi(t_1^*(0,x''))\leq Ce^{-c|x''|^2}\) proved earlier in Proposition \ref{prop-260072} and a Taylor expansion give \(1-t_1^*(0,x'')\leq Ce^{-c|x''|^2}\).
	Hence, for \(0<s<|x''|^{-2}\) and large \(|x''|\),
	\begin{align*}
		t_1^*(0,x'')-s\geq1-|x''|^{-2}-Ce^{-c|x''|^2} \geq1-\frac32|x''|^{-2}.
	\end{align*}
	Since \(|(0,x'')-y'|^2\leq(|x''|+1)^2\leq|x''|^2+3|x''|\) on \(\supp\psi\), we obtain \(e^{-\frac{|(0,x'')-y'|^2}{4(t_1^*(0,x'')-s)}} \geq e^{-\frac{|x''|^2}{4}-|x''|}\).
	Since \((4\pi(t_1^*(0,x'')-s))^{-\frac{n-1}{2}}\geq c>0\),
	\begin{align}\label{0706-06b}
		\int_\Sigma\Gamma'((0,x'')-y',t_1^*(0,x'')-s)\psi(y')\,dy' \geq cM_0e^{-\frac{|x''|^2}{4}-|x''|}, \qquad 0<s<|x''|^{-2}.
	\end{align}
	On the other hand, since \(k_0\) is the order of the first nonvanishing derivative of \(\phi\) at \(0\) and \(\phi\geq0\), we have \(\phi(s)\geq cs^{k_0}\) for \(0<s<|x''|^{-2}\) when \(|x''|\) is sufficiently large.
	Since \(t_1^*(0,x'')\to1^-\), we also have \(|x''|^{-2}<t_1^*(0,x'')\) for sufficiently large \(|x''|\).
	Using \eqref{0706-06b} on \((0,|x''|^{-2})\), we obtain
	\begin{align}\label{0706-06c}
		 & \frac{1}{\sqrt{4\pi}}\int_0^{t_1^*(0,x'')} \frac{\phi(s)}{(t_1^*(0,x'')-s)^{3/2}} \int_\Sigma\Gamma'((0,x'')-y',t_1^*(0,x'')-s)\psi(y')\,dy'ds\nonumber \\
		 & \qquad\geq cM_0|x''|^{-2k_0-2}e^{-\frac{|x''|^2}{4}-|x''|} \geq cM_0e^{-\frac{|x''|^2}{4}-2|x''|}.
	\end{align}
	The bounds \eqref{0706-06a} and \eqref{0706-06c}, together with \eqref{0706-07} and \eqref{0706-06}, give \(\ln\phi(t_1^*(0,x''))=-\frac{|x''|^2}{4}+O(|x''|)\).
	A Taylor expansion of \(\phi(t)\) at \(1\) now yields \(\ln(1-t_1^*(0,x''))=-\frac{|x''|^2}{4k_1}+O(|x''|)\), which is \eqref{0706-12}.
\end{proof}

The next proposition shows that \(\rho_*(x_1,t)\) is strictly increasing in \(|x_1|>R_{t_0}\) and gives the growth rate of \(R_{t_0}\) as \(t_0\to1^-\).
This proves the monotonicity statement in Theorem \ref{thm-260010}{\rm(iv)}.

\begin{proposition}\label{prop-260088}
	Let \(n\geq3\) and let \(w\) solve the Stokes system \eqref{0308-01} with boundary data \eqref{0324-13}, where \(j=1\), \(\psi\) is radial and satisfies {\rm\textbf{(S)}}, and \(\phi\) satisfies {\rm\textbf{(T1)}} and {\rm\textbf{(T2)}}.
	After increasing \(R_0\) in \eqref{260073} if necessary, for every \(0<t_0<1\), the unique positive zero \(\rho_*(x_1,t)\) of \(|x''|\mapsto\partial_{x_n}w_1(x_1,|x''|\omega'',0,t)\) satisfies, for \(x_1>R_{t_0}\) and \(0<t\leq t_0\),
	\begin{align}\label{0706-33}
		\partial_{x_1}\rho_*(x_1,t) =\sqrt{n-1}-\frac{M_2}{2(n-1)^{3/2}M_0x_1^{2}}+O\bigl(x_1^{-4}\bigr) \ \geq\ \frac{\sqrt{n-1}}2>0.
	\end{align}
	Moreover, \(\rho_*(x_1,t)\) is strictly increasing in \(|x_1|\) on \(\{|x_1|>R_{t_0}\}\), uniformly for \(0<t\leq t_0\), and \(R_{t_0}\) satisfies
	\begin{align}\label{260089}
		\lim_{t_0\to1^-} \frac{R_{t_0}}{\sqrt{\ln\bigl(\phi(t_1)/\phi(t_0)\bigr)}} =\frac2{\sqrt n}.
	\end{align}
\end{proposition}

\begin{proof}
	Fix \(0<t\leq t_0<1\), \(x_1>R_{t_0}\), and \(\omega''\in\mathbb S^{n-3}\).
	Evaluate at \(x'=(x_1,\rho_*(x_1,t)\omega'')\).
	By definition,
	\begin{align*}
		\partial_{x_n}w_1(x_1,\rho_*(x_1,t)\omega'',0,t)=0.
	\end{align*}
	The spatial derivatives of \eqref{0401-09} are continuous near \(x'\), and \eqref{transversality} shows that \(\omega''\cdot\nabla_{x'}\partial_{x_n}w_1(x',0,t)<0\).
	The implicit function theorem therefore gives \(\partial_{x_1}\rho_*(x_1,t)\), and differentiating this identity in \(x_1\) yields
	\begin{align*}
		0=\partial_{x_1}\partial_{x_n}w_1(x',0,t) +\partial_{x_1}\rho_*(x_1,t)\, \omega''\cdot\nabla_{x'}\partial_{x_n}w_1(x',0,t).
	\end{align*}
Dividing the last identity by \(\phi(t)\) and rearranging gives
	\begin{align}\label{260090}
		 & \left(\partial_{x_1}\rho_*(x_1,t) -\frac{(n-1)x_1}{\rho_*(x_1,t)}\right) \omega''\cdot\nabla_{x'} \frac{\partial_{x_n}w_1(x',0,t)}{\phi(t)}=-\left(\partial_{x_1} +\frac{(n-1)x_1}{\rho_*(x_1,t)} \omega''\cdot\nabla_{x'}\right) \frac{\partial_{x_n}w_1(x',0,t)}{\phi(t)}.
	\end{align}

	We estimate the right-hand side.
	The expansion \eqref{260087} gives \(\rho_*(x_1,t)^2-(n-1)x_1^2-\frac{M_2}{(n-1)M_0}=O(x_1^{-2})\) and \(|x'|\asymp x_1\), uniformly for \(0<t\leq t_0\).
	By \eqref{0916-02}, the error in each first partial derivative of \eqref{260066} is bounded by \(C_n M_6|x'|^{-n-7}\).
	Since
	\begin{align*}
		 & \left(\partial_{x_1}+\frac{(n-1)x_1}{|x''|^2} x''\cdot\nabla_{x'}\right) \bigl(|x''|^2-(n-1)x_1^2\bigr)=0,\\
		 & \left(\partial_{x_1}+\frac{(n-1)x_1}{|x''|^2} x''\cdot\nabla_{x'}\right)|x'|^{-m} =-mnx_1|x'|^{-m-2},
	\end{align*}
	differentiating \eqref{260066} gives, uniformly when \(\bigl||x''|^2-(n-1)x_1^2\bigr|\leq\lambda\) for any fixed \(\lambda>0\),
	\begin{align}\label{260091}
		\left(\partial_{x_1}+\frac{(n-1)x_1}{|x''|^2} x''\cdot\nabla_{x'}\right)\mathcal S(x') =-\frac{n(n+2)\Gamma\left(\frac n2\right)x_1 \left[M_0\bigl(|x''|^2-(n-1)x_1^2\bigr)-\frac{M_2}{n-1}\right]} {\pi^{\frac n2}|x'|^{n+4}} +O_\lambda\bigl(|x'|^{-n-5}\bigr).
	\end{align}
	At the zero, the leading term in \eqref{260091} is \(O(x_1^{-n-5})\).
	The coefficient \((n-1)x_1/\rho_*(x_1,t)\) is bounded, so \eqref{260084} and \eqref{0401-09} give
	\begin{align*}
		\left|\left(\partial_{x_1} +\frac{(n-1)x_1}{\rho_*(x_1,t)}\omega''\cdot\nabla_{x'}\right) \frac{\partial_{x_n}w_1(x',0,t)}{\phi(t)}\right| \leq Cx_1^{-n-5}.
	\end{align*}
	On the other hand, \eqref{transversality} gives
	\begin{align*}
		\left|\omega''\cdot\nabla_{x'} \frac{\partial_{x_n}w_1(x',0,t)}{\phi(t)}\right| \geq c\rho_*(x_1,t)|x'|^{-n-2}\geq cx_1^{-n-1}.
	\end{align*}
	By \eqref{260090}, \(\partial_{x_1}\rho_*(x_1,t) =(n-1)x_1/\rho_*(x_1,t)+O(x_1^{-4})\).
	Moreover,
	\begin{align*}
		\frac{(n-1)x_1}{\rho_*(x_1,t)} =\sqrt{n-1}\left(1+\frac{M_2}{(n-1)^2M_0x_1^2}+O(x_1^{-4})\right)^{-\frac12} =\sqrt{n-1}-\frac{M_2}{2(n-1)^{3/2}M_0x_1^2}+O(x_1^{-4}),
	\end{align*}
	which gives \eqref{0706-33}.
	Symmetry under \(x_1\mapsto-x_1\) and uniqueness imply that \(\rho_*(\cdot,t)\) is even.
	Increasing \(R_0\) gives the inequality in \eqref{0706-33} and the strict monotonicity of \(\rho_*(\cdot,t)\).

	Finally, \eqref{260089} follows immediately from \eqref{260073} and \(\phi(t_0)\to0\).
\end{proof}

The next proposition describes how the positivity regions of \(\partial_{x_n}w_1(\cdot,0,t)\) behave as \(t\to0^+\) and \(t\to1^-\), in the sense of pointwise convergence of their indicator functions.
This proves Theorem \ref{thm-260010}{\rm(v)}.

\begin{proposition}\label{prop-0703-08}
	Let \(w\) solve the Stokes system \eqref{0308-01} with boundary data \eqref{0324-13}, where \(j=1\), \(\psi\) satisfies {\rm\textbf{(S)}}, and \(\phi\) satisfies {\rm\textbf{(T1)}} and {\rm\textbf{(T2)}}.
	Then, for every \(x'\in\Sigma\) with \(|x'|>1\), the indicator functions of
	\begin{align*}
		\Omega_t := \left\{ x'\in\Sigma:\ |x'|>1,\ \partial_{x_n}w_1(x',0,t)>0 \right\}, \qquad \Omega_0 := \left\{ x'\in\Sigma:\ |x'|>1,\ \mathcal S(x')\leq0 \right\},
	\end{align*}
	where \(\mathcal S\) is given in \eqref{260009} with \(i=1\), satisfy
	\begin{align*}
		\lim_{t\to0^+}\mathbf 1_{\Omega_t}(x') = \mathbf 1_{\Omega_0}(x'), \qquad \lim_{t\to1^-}\mathbf 1_{\Omega_t}(x') =1.
	\end{align*}
\end{proposition}
\begin{proof}
	Fix \(x'\in\Sigma\) with \(|x'|>1\).
	We first consider \(t\to0^+\).
	By \eqref{260074}, \eqref{0401-09} gives \(\partial_{x_n}w_1(x',0,t)/\phi(t)\to-\mathcal S(x')\) as \(t\to0^+\).
	If \(\mathcal S(x')<0\), then \(\partial_{x_n}w_1(x',0,t)>0\) for all sufficiently small \(t>0\).
	If \(\mathcal S(x')>0\), then \(\partial_{x_n}w_1(x',0,t)<0\) for all sufficiently small \(t>0\).
	If \(\mathcal S(x')=0\), then \(\partial_{x_n}w_1(x',0,t)>0\) by the strict positivity of the integral in \eqref{0401-09} for all \(t>0\).
	Hence \( \mathbf 1_{\Omega_t}(x')\to\mathbf 1_{\Omega_0}(x') \) as \(t\to0^+\).

	As \(t\to1^-\), \(\phi(t)\to0\), while dominated convergence shows that the \(\tau\)-integral in \eqref{0401-09} tends to the corresponding strictly positive integral with \(t=1\).
	Hence \(\partial_{x_n}w_1(x',0,t)>0\) for all \(t<1\) sufficiently close to \(1\), and therefore \( \mathbf 1_{\Omega_t}(x')\to1 \) as \(t\to1^-\).
\end{proof}

	\subsection{Tangential velocity component in a direction different from that of the boundary data}\label{260092}
We first prove a formula for \(i\neq j\), including \(i=n\), and then assume that \(i\leq n-1\) for the rest of the subsection.

\subsubsection{Formula for the normal derivative of \(w_i\)}\label{260093}

\begin{proposition}\label{prop-0703-02}
Let \(1\leq i\leq n\) and \(1\leq j\leq n-1\) with \(i\neq j\), and let \(w\) solve the Stokes system \eqref{0308-01} with boundary data \eqref{0324-13}, where \(\psi\) satisfies {\rm\textbf{(S)}} and \(\phi\) satisfies {\rm\textbf{(T1)}}.
	Then, for every \(x'\in \Sigma\) and \(t>0\),
	\begin{align}\label{0421-01}
		\partial_{x_n}w_i(x',0,t) =
		\begin{cases}
			-\phi(t)R_i'(\partial_{x_j}\psi)(x'), & 1\leq i\leq n-1,\quad i\neq j, \\[2mm]
			-\phi(t)\partial_{x_j}\psi(x'), & i=n.
		\end{cases}
	\end{align}
\end{proposition}
\begin{proof}
	From \eqref{0330-01} and \eqref{0308-02}, we have
	\begin{align*}
		w_i(x,t) & = 4\int_0^t \int_{\Sigma}g_j(z',s) \partial_{x_j}\left( \int_{\Sigma\times (0, x_n)}\partial_{y_n}\Gamma(y,t-s)\partial_{x_i}N(x-z'-y)\,dy\right)\,dz'ds \\
		 & = 4\partial_{x_j}\int_{\Sigma\times(0,x_n)}\partial_{x_i}N(x-y)\int_0^t\int_{\Sigma}g_j(z',s)\partial_{y_n}\Gamma(y-z',t-s)\,dz'dsdy\\
		 & = 4\partial_{x_j}\int_{\Sigma\times(0,x_n)}\partial_{x_i}N(x-y)k(y,t)\,dy = 4\int_{\Sigma\times(0,x_n)}\partial_{y_i}N(y)\,\partial_{x_j}k(x-y,t)\,dy,
	\end{align*}
	where \(k\) is given in \eqref{0703-05}.
The equalities follow by exchanging the order of integration, translating the integration variable, and moving the tangential derivative onto \(k\).
	Applying \eqref{260269} to \(\partial_{x_j}(\psi\phi)\) gives
	\begin{align}\label{0414-01-t}
		\lim_{x_n\to 0^+}\partial_{x_j}k(x,t) = -\frac12 \partial_{x_j}\psi(x')\phi(t).
	\end{align}
	Differentiating the last representation of \(w_i\) with respect to \(x_n\) and using this boundary value gives the first equality below.
	Integrating the volume term by parts in \(y_j\) gives the second equality.
	\begin{align}\label{0421-04}
		\partial_{x_n}w_i(x,t) & = -2\phi(t) \int_{\Sigma}\partial_{i}N(y',x_n)\, \partial_{x_j}\psi(x'-y')\,dy' + 4\int_{\Sigma\times(0,x_n)}\partial_{y_i}N(y)\,\partial_{x_n}\partial_{x_j}k(x-y,t)\,dy \nonumber \\
		 & = -2\phi(t) \int_{\Sigma}\partial_{i}N(y',x_n)\,\partial_{x_j}\psi(x'-y')\,dy' + 4\int_{\Sigma\times(0,x_n)}\partial_{y_j}\partial_{y_i}N(y)\,\partial_{x_n}k(x-y,t)\,dy,
	\end{align}
	where \(\partial_i=\partial_{y_i}\) if \(1\leq i\leq n-1\) and \(\partial_i=\partial_{y_n}\) if \(i=n\).

	For \(\varepsilon>0\) and \(x_n>0\), define the region
	\begin{align}\label{0709-02}
		D_{\varepsilon,x_n} := \{y\in\R^n:0<y_n<x_n,\ |y|>\varepsilon\}.
	\end{align}
	The volume integral in the second line of \eqref{0421-04}, containing \(\partial_{y_j}\partial_{y_i}N\), is understood as the limit of the corresponding integral over \(D_{\varepsilon,x_n}\) as \(\varepsilon\to0^+\).

		Since we are interested only in the limit \(x_n\to0^+\), we may assume throughout the rest of the proof that \(0<x_n<1\).
	Since \(\phi\) and \(\phi'\) vanish on \([1,\infty)\), \eqref{0713-31} gives
	\begin{align}\label{0713-21}
	\|\partial_{x_n}k(\cdot,t)\|_{L^\infty(\R^n_+)} \leq \frac{1}{\sqrt{\pi}} \left( \|\psi\|_{L^\infty(\Sigma)} \|\phi'\|_{L^\infty(0,1)} + \|\Delta'\psi\|_{L^\infty(\Sigma)} \|\phi\|_{L^\infty(0,1)} \right),
	\end{align}
	and Lemma \ref{lem-260234}, applied with \(T=\max\{1,t\}\) and \(\phi\equiv0\) on \([1,\infty)\), gives
	\begin{align}\label{0713-10}
		\|\partial_{x_n}k(\cdot,t)\|_ {C^\alpha(\overline{\R^n_+})} \leq C_t \|\psi\|_{C^{2,\alpha}(\Sigma)} \|\phi\|_{C^1([0,1])}.
	\end{align}
		Since \(i\neq j\), oddness in \(y_j\), these bounds, and \( |\partial_{y_i}N(y)|\leq C|y|^{1-n} \) show that no boundary term appears in the integration by parts in \eqref{0421-04}.

Denote the volume integral in \eqref{0421-04} by \(I(x,t)\).
	Since \(i\neq j\) and \(j\leq n-1\), the kernel \(\partial_{y_j}\partial_{y_i}N(y)\) is odd in the tangential variable \(y_j\), while the set \(D_{\varepsilon,x_n}\) is symmetric with respect to \(y_j\).
	Thus \(\int_{D_{\varepsilon,x_n}} \partial_{y_j}\partial_{y_i}N(y)\,dy= 0 \) for each \(\varepsilon>0\) and \(x_n>0\), and therefore
	\begin{align}\label{0713-20}
		I(x,t) = \lim_{\varepsilon\to0^+} \int_{D_{\varepsilon,x_n}} \partial_{y_j}\partial_{y_i}N(y) \bigl(\partial_{x_n}k(x-y,t)-\partial_{x_n}k(x,t)\bigr)\,dy.
	\end{align}
	We split the integral in \eqref{0713-20} into the regions \( D_{\varepsilon,x_n}\cap B_1\) and \(D_{\varepsilon,x_n}\setminus B_1\).
	Since \( |\partial_{y_j}\partial_{y_i}N(y)| \leq C|y|^{-n} \) for \(y \neq 0\), the H\"{o}lder continuity of \(\partial_{x_n}k\) gives
	\begin{align}\label{0713-18}
		\left| \int_{D_{\varepsilon,x_n}\cap B_1} \partial_{y_j}\partial_{y_i}N(y) \bigl(\partial_{x_n}k(x-y,t)-\partial_{x_n}k(x,t)\bigr)\,dy \right| 
		 & \leq C[\partial_{x_n}k(\cdot,t)]_{C^\alpha}\,x_n^\alpha,
	\end{align}
	while, on the complementary region, the boundedness of \(\partial_{x_n}k\) gives
	\begin{align}\label{0713-19}
		\left| \int_{D_{\varepsilon,x_n}\setminus B_1} \partial_{y_j}\partial_{y_i}N(y) \bigl(\partial_{x_n}k(x-y,t)-\partial_{x_n}k(x,t)\bigr)\,dy \right| 
		 & \leq C\|\partial_{x_n}k(\cdot,t)\|_{L^\infty}\,x_n.
	\end{align}
	Since \(0<\alpha<1\) and \(0<x_n<1\), we have \(x_n\leq x_n^\alpha\).
	Thus, using \eqref{0713-18}, \eqref{0713-19}, \eqref{0713-10}, and \eqref{0713-21}, we obtain
	\begin{align*}
		|I(x,t)| \leq C_t\left( \|\psi\|_{C^{2,\alpha}(\Sigma)} \|\phi\|_{C^1([0,1])} \right) x_n^\alpha,
	\end{align*}
	and hence \eqref{0713-20} gives
	\begin{align}\label{0708-41}
		\lim_{x_n\to0^+} \mathrm{p.v.} \int_{\Sigma\times(0,x_n)} \partial_{y_j}\partial_{y_i}N(y)\, \partial_{x_n}k(x-y,t)\,dy =0.
	\end{align}
	Therefore \eqref{0421-04} and \eqref{0708-41} give
	\begin{align}\label{0709-01}
		\lim_{x_n\to0^+}\partial_{x_n}w_i(x,t) = -2 \phi(t)\lim_{x_n\to0^+} \int_{\Sigma} \partial_{i}N(y',x_n)\, \partial_{x_j}\psi(x'-y')\,dy'.
	\end{align}
	Using \eqref{260266} and \eqref{260267} in \eqref{0709-01}, we conclude \eqref{0421-01}.
\end{proof}

For the rest of this subsection we assume \(i\neq j\) and \(i\neq n\).
	Since \(i\) and \(j\) are then two distinct tangential indices, we necessarily have \(n\geq3\).
Relabeling the tangential coordinates by a permutation of \(\{1,\dots,n-1\}\) leaves the system \eqref{0308-01}, the representation \eqref{0330-01} and the hypothesis {\rm\textbf{(S)}} invariant, and carries the pair \((j,i)\) to any other pair of distinct tangential indices.
	We may therefore assume throughout that \(j=1\) and \(i=2\).

We use the notation in \eqref{260013}.
By Proposition \ref{prop-0703-02}, for \(x'\notin\supp\psi\) we have \( \partial_{x_n}w_2(x',0,t)=-\phi(t)R_2'(\partial_{x_1}\psi)(x'). \)
Using \( R_2'=-\partial_{x_2}(-\Delta')^{-1/2} \) and the kernel of \((-\Delta')^{-1/2}\), this becomes
\begin{align}\label{260094}
	\partial_{x_n}w_2(x',0,t) = \frac{n\Gamma\left(\frac n2\right)}{\pi^{\frac n2}}\,\phi(t) \int_{\Sigma} \frac{(x_1-y_1)(x_2-y_2)}{|x'-y'|^{n+2}} \psi(y')\,dy', \qquad x'\notin\supp\psi.
\end{align}
Since \(\phi(t)>0\) for \(0<t<1\), \eqref{260094} gives, for \(|x'|>1\) and \(0<t<1\),
\begin{align}\label{0713-22}
	\partial_{x_n}w_2(x',0,t)=0 \quad\Longleftrightarrow\quad R_2'(\partial_{x_1}\psi)(x')=0.
\end{align}
By \eqref{260094}, the zero set of \(\partial_{x_n}w_2(\cdot,0,t)\) in \(\{|x'|>1\}\) is independent of \(t\in(0,1)\).
We therefore omit the time argument when denoting these zeros.

Set \(\psi_\sigma:=\psi\circ\sigma\), and write \(w[\eta]\) for the solution with \(\psi\) replaced by \(\eta\).
Changing variables in \eqref{260094} gives
\begin{align}\label{260095}
	\partial_{x_n}w_2[\psi_\sigma](x',0,t) =\partial_{x_n}w_2[\psi](\sigma x',0,t), \qquad x'\notin\supp\psi_\sigma.
\end{align}

The next lemma identifies the first two normal derivatives of \(w_2\) at the boundary.
\begin{lemma}\label{lem-1007-07}
	Let \(n\geq3\), and let \(w\) solve the Stokes system \eqref{0308-01} with boundary data \eqref{0324-13}, where \(j=1\), \(\psi\) satisfies {\rm\textbf{(S)}}, and \(\phi\) satisfies {\rm\textbf{(T1)}} and {\rm\textbf{(T2)}}.
	Then the following hold.
	\begin{enumerate}
		\item \(\partial_{x_n}w_2\) and \(\partial_{x_n}^2w_2\) extend continuously to \(\overline{\R_+^n}\times(0,\infty)\).
		\item For every \(x'\in\Sigma\),
		\begin{align}\label{1007-08}
		    \partial_{x_n}^2w_2(x',0,t)&=\partial_{x_2}p(x',0,t) &&(t>0),\nonumber\\
		    \partial_{x_n}w_2(x',0,t)&=0 &&(t\geq1).
		\end{align}
	\end{enumerate}
\end{lemma}
\begin{proof}
\underline{Proof of (i).}\quad
Fix a compact interval \(J\Subset(0,\infty)\) and \(0<\gamma<\alpha\), and extend \(\phi\) by zero to negative times.
Lemma~\ref{lem-260265} with exponent \(\gamma\), applied to \(\partial_{x_1}\psi\) and its first derivatives, gives \(f:=R_2'(\partial_{x_1}\psi)\in C^{1,\gamma}(\Sigma)\cap C_0(\Sigma)\).
Since the heat kernel has zero first moments, subtracting the linear Taylor polynomial of \(f\) gives \(\|\Gamma'(\cdot,u)*_\Sigma f-f\|_{L^\infty(\Sigma)}\leq C[\nabla'f]_{C^\gamma(\Sigma)}u^{(1+\gamma)/2}\) for \(0<u\leq1\).
The boundary limits of heat potentials \cite[Chapter~IV, Sections~1--2]{LSU} give, in \(L^\infty(\Sigma)\) and uniformly for \(t\in J\),
\begin{align}
	\partial_{x_n}R_2'\partial_{x_1}k(\cdot,x_n,t)
	\longrightarrow-\frac1{4\sqrt\pi}\int_0^\infty
	\frac{\phi(t-u)\Gamma'(\cdot,u)*_\Sigma f-\phi(t)f}{u^{3/2}}\,du
	\qquad\text{as }x_n\to0^+,\label{1007-09}
\end{align}
with a limit continuous in \(t\) with values in \(C_0(\Sigma)\), and, for \(0<x_n\leq1\),
\begin{align}
	\sup_{t\in J}\left\|\partial_{x_n}^2R_2'\partial_{x_1}k(\cdot,x_n,t)\right\|_{L^\infty(\Sigma)}
	&\leq Cx_n\int_0^\infty u^{-5/2}\left(1+\frac{x_n^2}{u}\right)e^{-x_n^2/(4u)}\min\{u^{(1+\gamma)/2},1\}\,du\leq Cx_n^{\gamma-1}.\label{1007-10}
\end{align}

Applying \eqref{260258} to \(\psi\) with \(i=2\) and \(j=1\), we obtain
\begin{align*}
	\partial_{x_n}w_2(\cdot,x_n,t)={}&-2\phi(t)\partial_{x_n}N(\cdot,x_n)*_\Sigma f+4\int_0^{x_n}\partial_{x_n}N(\cdot,r)*_\Sigma\partial_{x_n}R_2'\partial_{x_1}k(\cdot,x_n-r,t)\,dr.
\end{align*}
The integral is \(O(x_n)\) in \(L^\infty(\Sigma)\) by \eqref{1007-09}, and \(2\partial_{x_n}N(\cdot,x_n)\) is an approximate identity, so \(\partial_{x_n}w_2(\cdot,x_n,t)\to-\phi(t)f\) in \(L^\infty(\Sigma)\), uniformly for \(t\in J\).
The Fourier multipliers of the Poisson kernel and the Riesz transform give
\begin{align*}
	\partial_{x_n}\bigl[2\partial_{x_n}N(\cdot,x_n)*_\Sigma f\bigr]=2\partial_{x_n}N(\cdot,x_n)*_\Sigma\partial_{x_1}\partial_{x_2}\psi.
\end{align*}
Differentiating the last formula in \(x_n\) and using \eqref{1007-09} and \eqref{1007-10} therefore give, in \(L^\infty(\Sigma)\) and uniformly for \(t\in J\),
\begin{align}\label{1007-12}
	\partial_{x_n}^2w_2(\cdot,x_n,t) \longrightarrow{}&-\phi(t)\partial_{x_1}\partial_{x_2}\psi-\frac1{2\sqrt\pi}\int_0^\infty
	\frac{\phi(t-u)\Gamma'(\cdot,u)*_\Sigma f-\phi(t)f}{u^{3/2}}\,du.
\end{align}
Both limits are continuous in \((x',t)\), and \(J\) is arbitrary, which proves {\rm(i)}.

\underline{Proof of (ii).}\quad
Since \(\frac1{2\sqrt\pi}\int_0^\infty(1-e^{-u|\xi'|^2})u^{-3/2}\,du=|\xi'|\) and \(f=R_2'(\partial_{x_1}\psi)=R_1'(\partial_{x_2}\psi)\),
\begin{align*}
    \frac1{2\sqrt\pi}\int_0^\infty
    \frac{\Gamma'(\cdot,u)*_\Sigma f-f}{u^{3/2}}\,du
    =\partial_{x_1}\partial_{x_2}\psi .
\end{align*}
Splitting the numerator of the integrand in \eqref{1007-12} as \([\phi(t-u)-\phi(t)]\,\Gamma'(\cdot,u)*_\Sigma f+\phi(t)\,[\Gamma'(\cdot,u)*_\Sigma f-f]\) turns the right-hand side of \eqref{1007-12} into the right-hand side of \eqref{tangential-pressure-boundary-value} with \(i=2\) and \(j=1\), so \(\partial_{x_n}^2w_2(x',0,t)=\partial_{x_2}p(x',0,t)\) for \(t>0\).
For \(t\geq1\), \(\phi(t)=0\) and \(\phi(t-u)=0\) for \(u\geq t\), so
\begin{align}
	\partial_{x_n}^2w_2(x',0,t)
	=-\frac1{2\sqrt\pi}\int_0^t\frac{\phi(t-u)}{u^{3/2}}
	\bigl(\Gamma'(\cdot,u)*_\Sigma f\bigr)(x')\,du
	=\partial_{x_2}p(x',0,t),\label{1007-11}
\end{align}
where, at \(t=1\), the integral converges absolutely since \(\phi(1)=\phi'(1)=0\) and {\rm\textbf{(T1)}} give \(|\phi(1-u)|\leq Cu^{1+\beta}\).
Since the boundary value of \(\partial_{x_n}w_2\) in {\rm(i)} is \(-\phi(t)f=0\) for \(t\geq1\), this proves \eqref{1007-08}.
\end{proof}

\subsubsection{Nonexistence of separation points}\label{260096}
For \(w_2\), let \(\ell_2(t)\) denote the smallest positive integer \(m\) for which \(\partial_{x_n}^m w_2(\cdot,0,t)\) is not identically zero on \(\Sigma\setminus\supp\psi\).
The next proposition identifies this order on both sides of \(t=1\) and excludes separation points of \(w_2\) on \(\{|x_1|\geq1,\ |x_2|\geq1\}\).
This proves Theorem \ref{thm-260014}{\rm(i)}.

\begin{proposition}\label{prop-260097}
	Let \(n\geq3\) and let \(w\) solve the Stokes system \eqref{0308-01} with boundary data \eqref{0324-13}, where \(j=1\), \(\psi\) satisfies {\rm\textbf{(S)}}, and \(\phi\) satisfies {\rm\textbf{(T1)}} and {\rm\textbf{(T2)}}.
	Then \(\ell_2(t)=1\) for \(0<t<1\) and \(\ell_2(t)=2\) for \(t\geq1\).
	For every \(t>0\) and every \(x'\in\Sigma\) with \(|x_1|\geq1\) and \(|x_2|\geq1\),
	\begin{align}\label{0714-020}
		\sgn\bigl(\partial_{x_n}^{\ell_2(t)}w_2(x',0,t)\bigr)=\sgn(x_1x_2).
	\end{align}
	If \(|x_1|\geq1\), \(|x_2|\geq1\), and \(t>0\), then \((x',t)\) is not a separation point of \(w_2\).
\end{proposition}
\begin{proof}
	If \(|x_k|\geq1\) for \(k=1,2\), then \(x_k(x_k-y_k)\geq0\) on \(\supp\psi\).
	Both inequalities are strict on the nonempty open set \(\{\psi>0\}\subset B_1'\).
	Since \(\phi(t)>0\) for \(0<t<1\), formula \eqref{260094} gives \eqref{0714-020} and \(\ell_2(t)=1\) for \(0<t<1\).
	We next consider \(t\geq1\).
	Since \(\phi(t)=0\), Lemma \ref{lem-1007-07}{\rm(ii)} and \eqref{0822-13} give
	\begin{align*}
		\partial_{x_n}^2w_2(x',0,t)
		={}&\frac{\Gamma(\frac{n+2}{2})}{2^{n+2}\pi^{\frac n2}\Gamma(\frac{n+3}{2})}
		\int_\Sigma(x_1-y_1)(x_2-y_2)\psi(y')\\
		&\quad\times\int_0^t\frac{\phi(t-s)}{s^{\frac{n+5}{2}}}
		{}_1F_1\left(\frac{n+2}{2},\frac{n+3}{2};-\frac{|x'-y'|^2}{4s}\right)\,dsdy'.
	\end{align*}
	The inner integral is strictly positive.
	Since \(x_1x_2(x_1-y_1)(x_2-y_2)\geq0\) on \(\supp\psi\), with strict inequality where \(\psi>0\), the outer integral has the sign of \(x_1x_2\).
	Thus
	\begin{align*}
		\sgn\bigl(\partial_{x_n}^2w_2(x',0,t)\bigr)=\sgn(x_1x_2).
	\end{align*}
	Consequently, \eqref{0714-020} holds for \(t\geq1\), and \(\partial_{x_n}w_2(\cdot,0,t)\equiv0\) while \(\partial_{x_n}^2w_2(\cdot,0,t)\not\equiv0\), so \(\ell_2(t)=2\).
	If \(0<t<1\), the necessary condition \eqref{1021-01} excludes separation points in the stated region.
	If \(t>1\), Proposition \ref{prop-260046}{\rm(iii)}, applied with \(\ell=2\) on an interval compactly contained in \((1,\infty)\), gives the same conclusion.

	It remains to consider \(t=1\).
	At each such \(x'\), \(-R_2'(\partial_{x_1}\psi)(x')\) is nonzero and has the sign of \(x_1x_2\) by \eqref{260094}.
	Continuity of \(R_2'(\partial_{x_1}\psi)\) gives a spatial neighborhood of \(x'\) on which it keeps that sign, including when \(|x_1|=1\) or \(|x_2|=1\).
	Multiplying by \(\phi(s)\), we see that the first normal derivative has that sign throughout that neighborhood when \(0<s<1\) and vanishes when \(s\geq1\).

	By \eqref{0714-020} at \(t=1\) and Lemma \ref{lem-1007-07}{\rm(i)}, \(\partial_{x_n}^2w_2\) has the same strict sign for sufficiently small \(x_n\geq0\) and \(s\) sufficiently close to \(1\).
	Taylor's theorem gives
	\begin{align}\label{1021-03}
		w_2(x,s)=x_n\partial_{x_n}w_2(x',0,s) +\int_0^{x_n}(x_n-r)\partial_{x_n}^2w_2(x',r,s)\,dr
	\end{align}
	Therefore \(w_2\) does not change sign in a sufficiently small neighborhood in \(\R_+^n\times(0,\infty)\).
	Hence, by Definition \ref{260000}, \((x',1)\) is not a separation point either.
\end{proof}

\subsubsection{Separation points for radial \(\psi\)}\label{260098}
For radial \(\psi\), the next proposition characterizes separation points by \(x_1x_2=0\) on \(|x'|>1\) for \(0<t<1\) and on \(|x'|\geq\sqrt2\) for \(t\geq1\).
This proves the statement for radial \(\psi\) in Theorem \ref{thm-260014}{\rm(ii)}.

\begin{proposition}\label{prop-0708-10}
	Let \(n\geq3\) and let \(w\) solve the Stokes system \eqref{0308-01} with boundary data \eqref{0324-13}, where \(j=1\), \(\psi\) satisfies {\rm\textbf{(S)}} and is radial, and \(\phi\) satisfies {\rm\textbf{(T1)}} and {\rm\textbf{(T2)}}.
	Then \(\ell_2(t)=1\) for \(0<t<1\) and \(\ell_2(t)=2\) for \(t\geq1\), and
	\begin{align}\label{0713-23}
		\sgn\bigl(\partial_{x_n}^{\ell_2(t)}w_2(x',0,t)\bigr)=\sgn(x_1x_2),
	\end{align}
	whenever \(0<t<1\) and \(|x'|>1\), or \(t\geq1\) and \(|x'|\geq\sqrt2\).
	If either \(0<t<1\) and \(|x'|>1\), or \(t\geq1\) and \(|x'|\geq\sqrt2\), then \((x',t)\) is a separation point of \(w_2\) if and only if \(x_1x_2=0\).
	Moreover, if \(|x'|>1\) and \(x_1x_2=0\), then \((x',1)\) is a separation point of \(w_2\).
\end{proposition}

\begin{proof}
	We first consider \(0<t<1\).
	The identity \(R_2'=-\partial_{x_2}(-\Delta')^{-1/2}\) gives
	\begin{align}\label{0713-24}
	    R'_2(\partial_{x_1}\psi)=-\partial_{x_2}\partial_{x_1}(-\Delta')^{-1/2}\psi .
	\end{align}
	Since \(\psi\) is radial, its Riesz potential is the radial function \(u(r)=((-\Delta')^{-1/2}\psi)(x')\), where \(r=|x'|\). The chain rule gives
	\begin{align*}
	    \partial_{x_2}\partial_{x_1}\bigl((-\Delta')^{-1/2}\psi\bigr)(x')
	    =x_1x_2\left(\frac{u''(r)}{r^2}-\frac{u'(r)}{r^3}\right).
	\end{align*}
	Hence \eqref{0713-24} reads
	\begin{align}\label{260099}
		R_2'(\partial_{x_1}\psi)(x') = -x_1x_2A(|x'|), \qquad A(r):=\frac{u''(r)}{r^2}-\frac{u'(r)}{r^3}.
	\end{align}

	We claim that \(A(r)>0\) for every \(r>1\).
	Indeed, for \(r>1\),
	\begin{align}\label{0708-050}
		u(r) = \frac{\Gamma\left(\frac n2\right)}{(n-2)\pi^{\frac n2}} \int_{\Sigma} \frac{\psi(y')}{|x'-y'|^{n-2}}\,dy', \qquad |x'|=r.
	\end{align}
Writing \(\psi(y')=\Psi(|y'|)\), we obtain
	\begin{align}\label{0708-05}
		u(r) = \frac{\Gamma\left(\frac n2\right)}{(n-2)\pi^{\frac n2}} \int_0^1 \Psi(s)s^{n-2} \int_{\mathbb S^{n-2}} (r^2+s^2-2rs\omega_1)^{-\frac{n-2}{2}} \,d\sigma(\omega)ds.
	\end{align}
	For \(r>1\) and \(0\leq s\leq1\) we have \(s/r<1\), so Lemma \ref{lem-0422-03} gives
	\begin{align}\label{0708-16}
		\int_{\mathbb S^{n-2}} (r^2+s^2-2rs\omega_1)^{-\frac{n-2}{2}}\,d\sigma(\omega) = |\mathbb S^{n-2}|\,r^{-(n-2)} \sum_{k=0}^{\infty} a_k\left(\frac{s}{r}\right)^{2k}, \qquad a_k := \frac{ \left(\frac{n-2}{2}\right)_k\left(\frac12\right)_k }{ \left(\frac{n-1}{2}\right)_k\,k! }>0.
	\end{align}
	For each \(\varepsilon>0\), the series in \eqref{0708-16} converges uniformly for \(s\in[0,1]\) and \(r\geq1+\varepsilon\), so integrating \eqref{0708-05} term by term gives
	\begin{align}\label{0708-01}
		u(r) = \widetilde C_n \sum_{k=0}^{\infty} a_k\mu_{2k}\,r^{-(n-2+2k)}, \qquad \mu_{2k}:=\int_0^1\Psi(s)s^{n-2+2k}\,ds, \qquad \widetilde C_n:= \frac{\Gamma\left(\frac n2\right)|\mathbb S^{n-2}|} {(n-2)\pi^{\frac n2}}>0.
	\end{align}

	Since \(\psi\geq0\) and \(\psi\not\equiv0\), we have \(\mu_{2k}>0\) for every \(k\geq0\).
	Moreover, from \eqref{0708-16},
	\begin{align*}
		\frac{a_{k+1}}{a_k} = \frac{\left(k+\frac{n-2}{2}\right)\left(k+\frac12\right)} {\left(k+\frac{n-1}{2}\right)(k+1)}<1,
	\end{align*}
	so \(0<a_k\leq a_0=1\), and \(0\leq\mu_{2k}\leq\mu_0\) since \(\Psi\geq0\) and \(0\leq s\leq1\).
	Consequently, for every \(\varepsilon>0\), the series \eqref{0708-01} and the two series obtained from it by differentiating once and twice are dominated on \(r\geq1+\varepsilon\) by a constant multiple of \( \sum_{k\geq0}(1+k)^2(1+\varepsilon)^{-2k}<\infty\), so \eqref{0708-01} may be differentiated twice term by term on \((1,\infty)\).
	This gives
	\begin{align*}
		u''(r)-\frac{u'(r)}{r} = \widetilde C_n \sum_{k=0}^{\infty} a_k\mu_{2k}(n-2+2k)(n+2k)\,r^{-(n+2k)},
	\end{align*}
	whose terms are all nonnegative and whose \(k=0\) term is strictly positive since \(\mu_0>0\).
	Therefore \( A(r)=r^{-2}\left(u''(r)-r^{-1}u'(r)\right)>0 \) for every \(r>1\), as claimed.

	Combining Proposition \ref{prop-0703-02} with \eqref{260099}, we obtain, for \(|x'|>1\) and \(0<t<1\),
	\begin{align*}
		\partial_{x_n}w_2(x',0,t) = -\phi(t)R_2'(\partial_{x_1}\psi)(x') = \phi(t)\,x_1x_2\,A(|x'|).
	\end{align*}
	Since \(\phi(t)>0\) and \(A(|x'|)>0\), the sign of \(\partial_{x_n}w_2(x',0,t)\) is that of \(x_1x_2\), which is \eqref{0713-23}.

	Thus \(\ell_2(t)=1\) for \(0<t<1\).
	We next consider \(t\geq1\).
	Since \(\phi(t)=0\), Lemma \ref{lem-1007-07}{\rm(ii)} and \eqref{0825-13} give, with \(r=|x'|\),
	\begin{align*}
		\partial_{x_n}^2w_2(x',0,t)
		={}&\frac{\Gamma(\frac{n+2}{2})}{2^{n+2}\pi^{\frac n2}\Gamma(\frac{n+3}{2})}\frac{x_1x_2}{r^2}
		\int_\Sigma\bigl((r-y_1)^2-y_2^2\bigr)\psi(y')\\
		&\quad\times\int_0^t\frac{\phi(t-s)}{s^{\frac{n+5}{2}}}
		{}_1F_1\left(\frac{n+2}{2},\frac{n+3}{2};-\frac{|r\mathbf e_1'-y'|^2}{4s}\right)\,dsdy'.
	\end{align*}
	The inner integral is strictly positive.
	Since \(\psi\) is continuous and vanishes outside \(\overline{B_1'}\), every point where \(\psi(y')>0\) satisfies \(|y'|<1\).
	Thus for \(r\geq\sqrt2\) and every such \(y'\),
	\begin{align}\label{radial-geometric-factor}
	    (r-y_1)^2-y_2^2
	    =2\left(y_1-\frac r2\right)^2+\frac{r^2}{2}-(y_1^2+y_2^2)>0.
	\end{align}
	Therefore
	\begin{align*}
		\sgn\bigl(\partial_{x_n}^2w_2(x',0,t)\bigr)=\sgn(x_1x_2),\qquad |x'|\geq\sqrt2.
	\end{align*}
	It follows that \eqref{0713-23} holds for \(t\geq1\) and that \(\ell_2(t)=2\).
	The joint continuity up to the boundary required in Proposition \ref{prop-260046}{\rm(iii)} follows from Lemma \ref{lem-1007-07}{\rm(i)}.

	For \(0<t<1\), apply that proposition with \(f(x')=x_1x_2\) and \(\ell=1\).
	For \(t>1\) and \(|x'|\geq\sqrt2\) with \(x_1x_2\ne0\), continuity gives a neighborhood on which the second normal derivative is nonzero and has the sign of \(x_1x_2\), and the first normal derivative vanishes throughout a time interval by \eqref{1007-08}.
	Equation \eqref{1021-03} then shows that \(w_2\) does not change sign near \((x',t)\).

	If \(x_1x_2=0\), every neighborhood contains points \(z'\) with \(|z'|>\sqrt2\) and \(z_1z_2>0\), as well as points with \(|z'|>\sqrt2\) and \(z_1z_2<0\).
	At each such point, \eqref{1021-03} shows that \(w_2\) has the sign of \(z_1z_2\) for sufficiently small \(x_n>0\).
	By continuity in the interior, every cylinder in Definition \ref{260000} contains open sets on which \(w_2\) is positive and open sets on which it is negative, including when \(|x'|=\sqrt2\).

	At \(t=1\), let \(|x'|>1\) and \(x_1x_2=0\).
	By \eqref{0713-23} for \(s<1\), \(\partial_{x_n}w_2\) takes both signs in every neighborhood of \((x',1)\).
	Hence \((x',1)\) is a separation point by Proposition \ref{prop-260046}{\rm(i)}.

	If \(|x'|\geq\sqrt2\) and \(x_1x_2\neq0\), then \eqref{0713-23} and \eqref{1007-08} show that, in a sufficiently small neighborhood of \((x',1)\), \(\partial_{x_n}w_2(\cdot,0,s)\) is either zero or has the sign of \(x_1x_2\).
	By \eqref{0713-23} at \(t=1\) and Lemma \ref{lem-1007-07}{\rm(i)}, \(\partial_{x_n}^2w_2\) has the sign of \(x_1x_2\) for sufficiently small \(x_n\geq0\) and \(s\) sufficiently close to \(1\).
	
	Equation \eqref{1021-03} then shows that \(w_2\) does not change sign in the half-space near \((x',1)\), so \((x',1)\) is not a separation point.
\end{proof}

\subsubsection{Existence of the separation points}\label{260100}
For general \(\psi\), the next proposition gives separation points whose spatial locations are independent of \(t\in(0,1]\) on every line with \(|x_2|>1\) and, after exchanging the first two coordinates, on every line with \(|x_1|>1\).
This proves the existence statement in Theorem \ref{thm-260014}{\rm(ii)}.

\begin{proposition}\label{prop-260101}
	Let \(n\geq3\) and let \(w\) solve the Stokes system \eqref{0308-01} with boundary data \eqref{0324-13}, where \(j=1\), \(\psi\) satisfies {\rm\textbf{(S)}}, and \(\phi\) satisfies {\rm\textbf{(T1)}} and {\rm\textbf{(T2)}}.
	Then the following hold.
	\begin{enumerate}
		\item[\rm(i)] For every \(x''\in(\mathbf e_1')^\perp\) with \(|x_2|>1\), there are \(x_1^*(x'')\in(-1,1)\) and \(\delta>0\), independent of \(t\), such that \(\partial_{x_n}w_2(x_1^*(x'')\mathbf e_1'+x'',0,t)=0\) for every \(0<t<1\).
		If \(x_2>1\), then, for every \(0<t<1\),
		\begin{align}\label{0919-10}
			\partial_{x_n}w_2(x',0,t)
			\begin{cases}
				<0, & x_1^*(x'')-\delta<x_1<x_1^*(x''), \\
				>0, & x_1^*(x'')<x_1<x_1^*(x'')+\delta.
			\end{cases}
		\end{align}
		If \(x_2<-1\), both inequalities are reversed.
		For every \(0<t\leq1\), \((x_1^*(x'')\mathbf e_1'+x'',t)\) is a separation point of \(w_2\).

		\item[\rm(ii)] For every \(x^\flat\in(\mathbf e_2')^\perp\) with \(|x_1|>1\), there are \(x_2^*(x^\flat)\in(-1,1)\) and \(\delta>0\), independent of \(t\), such that \(\partial_{x_n}w_2(x',0,t)=0\) at \(x_2=x_2^*(x^\flat)\) for every \(0<t<1\).
		If \(x_1>1\), then, for every \(0<t<1\),
		\begin{align}\label{0919-11}
			\partial_{x_n}w_2(x',0,t)
			\begin{cases}
				<0, & x_2^*(x^\flat)-\delta<x_2<x_2^*(x^\flat), \\
				>0, & x_2^*(x^\flat)<x_2<x_2^*(x^\flat)+\delta.
			\end{cases}
		\end{align}
		If \(x_1<-1\), both inequalities are reversed.
		For every \(0<t\leq1\), \((x_2^*(x^\flat)\mathbf e_2'+x^\flat,t)\) is a separation point of \(w_2\).

	\end{enumerate}
\end{proposition}

\begin{proof}
	\underline{Proof of (i).}\quad
	For every \(x_1\in\R\), \(|x'|\geq|x_2|>1\).
	Hence the line \(\{(x_1,x''):x_1\in\R\}\) does not meet \(\overline{B_1'}\supset\supp\psi\), and \eqref{260094} is valid along it.
	Assume first that \(x_2>1\).
	Proposition \ref{prop-260097} gives \(\partial_{x_n}w_2(-1,x'',0,t)<0\) and \(\partial_{x_n}w_2(1,x'',0,t)>0\), so the intermediate value theorem gives a zero in \((-1,1)\).
	Since \( \inf_{y'\in\supp\psi}|x''-y''|\geq|x''|-1>0\), Lemma \ref{lem-260044}{\rm(i)}, applied with \(I=\R\), \(m=1\), \(g(y')=(x_2-y_2)\psi(y')\), and \(a=(n+2)/2\), shows that
	\begin{align*}
		\frac{1}{\phi(t)}\partial_{x_n}w_2(x',0,t)
	\end{align*}
	is real-analytic in \(x_1\in\R\).
	Its values at \(-1\) and \(1\) have opposite signs, so by Lemma \ref{lem-260044}{\rm(iii)} it has finitely many zeros in \([-1,1]\).
	Therefore at least one zero has negative values immediately to its left and positive values immediately to its right.
	Denote this zero by \(x_1^*(x'')\).
	Since it is isolated, \eqref{0919-10} holds for some \(\delta>0\).

	If \(x_2<-1\), then \(x_2-y_2<0\) on \(\supp\psi\).
		Repeating the preceding construction with all signs reversed gives the required zero.
	For each \(0<t<1\), \eqref{0919-10} shows that \(\partial_{x_n}w_2\) takes both signs arbitrarily close to \(x_1^*(x'')\mathbf e_1'+x''\).
	Proposition \ref{prop-260046}{\rm(i)} therefore shows that \((x_1^*(x'')\mathbf e_1'+x'',t)\) is a separation point.
	At \(t=1\), taking \(s<1\) in \eqref{0919-10} shows that the first normal derivative still takes both signs in every neighborhood of \((x_1^*(x'')\mathbf e_1'+x'',1)\).
	The same proposition then shows that \((x_1^*(x'')\mathbf e_1'+x'',1)\) is a separation point.

	\underline{Proof of (ii).}\quad
	Applying (i) to \(\psi_\sigma=\psi\circ\sigma\), which satisfies {\rm\textbf{(S)}}, and using \eqref{260095} gives the zero and \eqref{0919-11} in (ii).
	Applying Proposition \ref{prop-260046}{\rm(i)} with the signs in \eqref{0919-11} proves separation for \(0<t\leq1\).

\end{proof}

\subsubsection{Uniqueness of the separation points for large \(|x''|\) and \(|x^\flat|\)}\label{260102}
The next proposition shows that the two zeros are unique when \(|x''|\) or \(|x^\flat|\) is larger than \(1+2\sqrt{n+1}\).
This proves Theorem \ref{thm-260014}{\rm(iii)}.

\begin{proposition}\label{prop-260103}
		Let \(n\geq3\) and let \(w\) solve the Stokes system \eqref{0308-01} with boundary data \eqref{0324-13}, where \(j=1\), \(\psi\) satisfies {\rm\textbf{(S)}}, and \(\phi\) satisfies {\rm\textbf{(T1)}} and {\rm\textbf{(T2)}}.
	Then, with \(R:=1+2\sqrt{n+1}\), the following hold for every \(0<t<1\).
	\begin{enumerate}
		\item[\rm(i)] For every \(x''\in(\mathbf e_1')^\perp\) with \(|x_2|>1\) and \(|x''|>R\), the point \(x_1^*(x'')\in(-1,1)\) in \eqref{0919-10} is the unique zero of \(\partial_{x_n}w_2(x',0,t)\) in \(x_1\in\R\), \((x_1^*(x'')\mathbf e_1'+x'',t)\) is the unique separation point of \(w_2\) on the line \(\{x_1\mathbf e_1'+x'':x_1\in\R\}\), and
		\begin{align}\label{0928-003}
			\sgn \left( \partial_{x_n}w_2(x',0,t) \right) = \sgn(x_2)\,\sgn\bigl(x_1-x_1^*(x'')\bigr).
		\end{align}

		\item[\rm(ii)] For every \(x^\flat\in(\mathbf e_2')^\perp\) with \(|x_1|>1\) and \(|x^\flat|>R\), the point \(x_2^*(x^\flat)\in(-1,1)\) in \eqref{0919-11} is the unique zero of \(\partial_{x_n}w_2(x',0,t)\) in \(x_2\in\R\), \((x_2^*(x^\flat)\mathbf e_2'+x^\flat,t)\) is the unique separation point of \(w_2\) on the line \(\{x_2\mathbf e_2'+x^\flat:x_2\in\R\}\), and
		\begin{align*}
		\sgn \left(\partial_{x_n}w_2(x',0,t)\right)=\sgn(x_1)\,\sgn\bigl(x_2-x_2^*(x^\flat)\bigr).
		\end{align*}

	\end{enumerate}
\end{proposition}
\begin{proof}
	\underline{Proof of (i).}\quad
	Differentiating \eqref{260094} in \(x_1\) and using \( |x'-y'|^2=(x_1-y_1)^2+|x''-y''|^2 \) gives
	\begin{align}\label{260104}
		\partial_{x_1}\partial_{x_n}w_2(x',0,t) = \frac{n\Gamma\left(\frac n2\right)}{\pi^{\frac n2}}\,\phi(t) \int_{\Sigma} (x_2-y_2)\, \frac{|x''-y''|^2-(n+1)(x_1-y_1)^2} {|x'-y'|^{n+4}} \,\psi(y')\,dy'.
	\end{align}
	Let \(|x_1|\leq1\) and \(y'\in\supp\psi\).
	Then \(|x_1-y_1|\leq2\), while \(|x''-y''|\geq|x''|-1>2\sqrt{n+1}\).
		
	The integrand in \eqref{260104} therefore has the sign of \(x_2-y_2\), which is \(\sgn(x_2)\) throughout \(\supp\psi\) since \(|x_2|>1\).
	Consequently, \(x_1\mapsto\partial_{x_n}w_2(x',0,t)\) is strictly increasing on \([-1,1]\) when \(x_2>1\) and strictly decreasing there when \(x_2<-1\), so it has at most one zero in this interval.
	
	Proposition \ref{prop-260097} excludes zeros with \(|x_1|\geq1\), while Proposition \ref{prop-260101} supplies a zero in \((-1,1)\).
	Hence \(x_1^*(x'')\) is the unique zero on \(\R\), and, by Propositions \ref{prop-260046}{\rm(i)} and \ref{prop-260101}, \((x_1^*(x'')\mathbf e_1'+x'',t)\) is the unique separation point of \(w_2\) on this line.
	The monotonicity on \([-1,1]\) and \eqref{0714-020} for \(|x_1|\geq1\) give \eqref{0928-003}.

	\underline{Proof of (ii).}\quad
	Applying {\rm(i)} to \(\psi_\sigma=\psi\circ\sigma\), which satisfies {\rm\textbf{(S)}}, and using \eqref{260095} gives {\rm(ii)}.
\end{proof}

\subsubsection{Asymptotics and monotonicity of the separation points}\label{260105}
The next lemma compares two integral quotients.
\begin{lemma}\label{lem-260106}
	Let \(E\subset \Sigma\) be a bounded measurable set and let \(g\in L^1(E)\) be a nonnegative function with \(\int_E g(y')\,dy' >0\).
	Suppose that \(f\colon E\to \R\) is measurable and takes values in an interval of length \(L\), and that \(\varpi\colon E\to \R_+\) is measurable and satisfies
	\begin{align*}
		0<\inf_E \varpi\leq \sup_E\varpi<\infty, \qquad \frac{\sup_{E}\varpi}{\inf_E \varpi}\leq e^{\Lambda}
	\end{align*}
	for some \(\Lambda\geq0\).
	Then
	\begin{align*}
		\left| \frac{\int_E f\varpi g\,dy'}{\int_E \varpi g\,dy'} - \frac{\int_E fg\,dy'}{\int_E g\,dy'}\right| \leq L \tanh \left( \frac{\Lambda}{4}\right).
	\end{align*}
\end{lemma}
\begin{proof}
	Denote
	\begin{align*}
		a:= \inf_E \varpi,\quad b:= \sup_E \varpi, \quad \overline f := \frac{\int_E fg\,dy'}{\int_E g\,dy'}, \quad \overline \varpi := \frac{\int_E \varpi g\,dy'}{\int_E g\,dy'}\in[a,b].
	\end{align*}
	We may assume \(a<b\), the case \(a=b\) being trivial.
	Since \(\int_E \varpi g\,dy' = \overline\varpi\int_E g\,dy'\), for every constant \(c\),
	\begin{align}\label{0815-08}
		\frac{\int_E f\varpi g\,dy'}{\int_E \varpi g\,dy'} - \overline f & = \frac{\int_E(f-c)(\varpi-\overline \varpi)g\,dy'}{\overline \varpi\int_E g\,dy'}.
	\end{align}
	Suppose that \(f\) takes values in the interval \([f_-, f_+]\) with \(f_+-f_-= L\).
	Choose \(c:= \frac{f_- + f_+}{2}\) so that \(|f-c|\leq \frac{L}{2}\).
	Then \eqref{0815-08} gives
	\begin{align*}
		\left| \frac{\int_E f\varpi g\,dy'}{\int_E \varpi g\,dy'} - \overline f\right| \leq \frac{L}{2}\,\frac{\int_E|\varpi-\overline \varpi|g\,dy'}{\overline \varpi\int_E g\,dy'}.
	\end{align*}
	Now, by the convexity of \(|u-\overline \varpi|\) in \(u\) and the fact that \(\varpi\) takes values in \([a,b]\), we have
	\begin{align*}
		|\varpi-\overline \varpi|\leq \frac{b-\varpi}{b-a}(\overline \varpi -a) + \frac{\varpi-a}{b-a}(b-\overline \varpi) \qquad\text{on }E.
	\end{align*}
	Integrating this inequality against \(g\) over \(E\) and using \(\int_E (\varpi-d)g\,dy' = (\overline \varpi-d)\int_Eg\,dy'\) for every constant \(d\) give
	\begin{align*}
		\int_E|\varpi-\overline \varpi|g\,dy' \leq \frac{2(b-\overline \varpi)(\overline \varpi-a)}{b-a}\int_E g\,dy'.
	\end{align*}
	Thus we have
	\begin{align*}
		\frac{\int_E |\varpi-\overline \varpi|g\,dy'}{\int_E g\,dy'} \leq \frac{2(b-\overline \varpi)(\overline \varpi-a)}{b-a}, \qquad\text{hence}\qquad \left| \frac{\int_E f\varpi g\,dy'}{\int_E \varpi g\,dy'} - \overline f\right| \leq \frac{L(b-\overline \varpi)(\overline \varpi-a)}{(b-a)\overline \varpi}.
	\end{align*}
	Finally, \(\max_{v\in[a,b]}\frac{(b-v)(v-a)}{v} = (\sqrt b -\sqrt a)^2\), with the maximum being attained at \(v=\sqrt{ab}\), so that
	\begin{align*}
		\left| \frac{\int_E f\varpi g\,dy'}{\int_E \varpi g\,dy'} - \overline f\right| \leq L\,\frac{\sqrt{b}-\sqrt{a}}{\sqrt{b}+\sqrt{a}} = L\,\frac{\sqrt{b/a}-1}{\sqrt{b/a}+1} \leq L \tanh \left(\frac{\Lambda}{4}\right),
	\end{align*}
	where in the last inequality we used that \(\frac{r-1}{r+1}\) is increasing and that \(\sqrt{b/a}\leq e^{\Lambda/2}\).
\end{proof}

The next proposition gives the asymptotic expansions of \(x_1^*\) and \(x_2^*\) and their monotonicity along the coordinate directions.
This proves Theorem \ref{thm-260014}{\rm(iv)}.
\begin{proposition}\label{prop-260107}
		Let \(n\geq3\) and let \(w\) solve the Stokes system \eqref{0308-01} with boundary data \eqref{0324-13}, where \(j=1\), \(\psi\) satisfies {\rm\textbf{(S)}}, and \(\phi\) satisfies {\rm\textbf{(T1)}} and {\rm\textbf{(T2)}}.
Then, with \(R:=1+2\sqrt{n+1}\), the following hold.
	\begin{enumerate}
		\item\label{item-0714-31} For every \(x''\in(\mathbf e_1')^\perp\) with \(|x_2|>1\) and \(|x''|>R\), and for every \(x^\flat\in(\mathbf e_2')^\perp\) with \(|x_1|>1\) and \(|x^\flat|>R\),
		\begin{align}\label{260108}
			\left|x_1^*(x'')-\frac{m_1}{M_0}\right| \leq \frac{n+2}{|x''|-1} + \frac{1}{|x_2|-1}, \qquad \left|x_2^*(x^\flat)-\frac{m_2}{M_0}\right| \leq \frac{n+2}{|x^\flat|-1} + \frac{1}{|x_1|-1}.
		\end{align}
		\item\label{item-0714-32} Fix \(\omega''\in\mathbb S^{n-3}\) with \(\omega_2''\neq0\).
		Then, as \(\rho\to\infty\),
		\begin{align}\label{260109}
			x_1^*(\rho\,\omega'') & = \frac{m_1}{M_0} +\frac1{\rho M_0}\int_\Sigma \left(y_1-\frac{m_1}{M_0}\right) \left((n+2)\omega''\cdot y''-\frac{y_2}{\omega_2''}\right) \psi(y')\,dy' + O\!\left(\frac1{\rho^{2}}\right), \\
			x_2^*(\rho\,\sigma\omega'') & = \frac{m_2}{M_0} +\frac1{\rho M_0}\int_\Sigma \left(y_2-\frac{m_2}{M_0}\right) \left((n+2)\sigma\omega''\cdot y^\flat-\frac{y_1}{\omega_2''}\right) \psi(y')\,dy' + O\!\left(\frac1{\rho^{2}}\right), \label{260110}
		\end{align}
		and
		\begin{align}\label{260111}
			\frac{d}{d\rho}x_1^*(\rho\omega'') & =-\frac1{\rho^2M_0}\int_\Sigma \left(y_1-\frac{m_1}{M_0}\right) \left((n+2)\omega''\cdot y''-\frac{y_2}{\omega_2''}\right) \psi(y')\,dy'+O\!\left(\frac1{\rho^3}\right), \\
			\frac{d}{d\rho}x_2^*(\rho\sigma\omega'') & =-\frac1{\rho^2M_0}\int_\Sigma \left(y_2-\frac{m_2}{M_0}\right) \left((n+2)\sigma\omega''\cdot y^\flat-\frac{y_1}{\omega_2''}\right) \psi(y')\,dy'+O\!\left(\frac1{\rho^3}\right).
		\end{align}
		If \(\int_\Sigma\left(y_1-\frac{m_1}{M_0}\right)y_2\psi(y')\,dy'\neq0\), then there exists \(R_1=R_1(n,\psi)\geq R\) such that \( x_1^*(x_2\mathbf e_2')\) and \( x_2^*(x_1\mathbf e_1')\) are either both strictly increasing or both strictly decreasing on each of \((R_1,\infty)\) and \((-\infty,-R_1)\).
	\end{enumerate}
\end{proposition}

\begin{proof}
	By Proposition \ref{prop-260103}, \(\partial_{x_n}w_2(x_1^*(x'')\mathbf e_1'+x'',0,t)=0\) may be rewritten as
	\begin{align}\label{260112}
		x_1^*(x'') = \frac{\displaystyle\int_{\Sigma}y_1\varpi(y')\psi(y')\,dy'} {\displaystyle\int_{\Sigma}\varpi(y')\psi(y')\,dy'}, \qquad \varpi(y') := \frac{\sgn(x_2)\,(x_2-y_2)}{|x'-y'|^{n+2}}.
	\end{align}
	The factor \(\sgn(x_2)\) makes \(\varpi>0\) on \(\supp\psi\).

	\underline{Proof of (i).}\quad
	We first prove \eqref{260108} for \(x_1^*\).
	Lemma \ref{lem-260106} applies with \(E:=\supp\psi\), \(f(y')=y_1\), \(g=\psi\), and \(\varpi\) as in \eqref{260112}.
	Since \(f\) takes its values in \([-1,1]\), we may take \(L=2\).
	On the other hand,
	\begin{align*}
		\frac{\sup_{\supp\psi}\varpi}{\inf_{\supp\psi}\varpi} \leq \frac{|x_2|+1}{|x_2|-1} \left(\frac{|x'|+1}{|x'|-1}\right)^{n+2} \leq \frac{|x_2|+1}{|x_2|-1} \left(\frac{|x''|+1}{|x''|-1}\right)^{n+2}.
	\end{align*}
	We may therefore take
	\begin{align*}
		\Lambda = \ln\frac{|x_2|+1}{|x_2|-1} + (n+2)\ln\frac{|x''|+1}{|x''|-1}
	\end{align*}
	in Lemma \ref{lem-260106}, which gives, using \(\tanh r\leq r\),
	\begin{align*}
		\left|x_1^*(x'')-\frac{m_1}{M_0}\right| \leq 2\tanh\left(\frac{\Lambda}{4}\right) \leq \frac{\Lambda}{2}.
	\end{align*}
	Since \(\ln\frac{r+1}{r-1}=\ln\left(1+\frac{2}{r-1}\right)\leq\frac{2}{r-1}\) for \(r>1\), this is exactly \eqref{260108}.

	\underline{Proof of (ii).}\quad
	Let \(x''=\rho\,\omega''\), so that \(x_2=\rho\,\omega_2''\), and write \(r:=|x'|=\sqrt{x_1^2+\rho^2}\).
	For \(y'\in\supp\psi\) we have \( |x'-y'|^2=r^2(1+\tau) \) with \( \tau:=\frac{|y'|^2-2x'\cdot y'}{r^2}\), and \(|\tau|\leq\frac{1+2r}{r^2}\).
	Hence, for \(\rho\) large,
	\begin{align*}
		\frac{1}{|x'-y'|^{n+2}} = \frac{1}{r^{n+2}} \left( 1+\frac{n+2}{2}\cdot\frac{2x'\cdot y'-|y'|^2}{r^2}+O\!\left(\frac1{\rho^2}\right) \right) = \frac{1}{r^{n+2}} \left( 1+\frac{(n+2)\,\omega''\cdot y''}{\rho}+O\!\left(\frac1{\rho^2}\right) \right),
	\end{align*}
	uniformly for \(y'\in\supp\psi\) and \(|x_1|\leq1\).
	Here we used \(r^{-2}\rho=\rho^{-1}+O(\rho^{-3})\) and \(\frac{x_1y_1}{r^2}=O(\rho^{-2})\).
	Since also \( x_2-y_2=\rho\,\omega_2''\left(1-\frac{y_2}{\rho\,\omega_2''}\right)\), the weight of \eqref{260112} factorizes as
	\begin{align*}
		\varpi(y') = \frac{|\omega_2''|\,\rho}{r^{n+2}} \left( 1+\frac{a(y')}{\rho}+O\!\left(\frac1{\rho^2}\right) \right), \qquad a(y'):=(n+2)\,\omega''\cdot y''-\frac{y_2}{\omega_2''},
	\end{align*}
	again uniformly for \(y'\in\supp\psi\) and \(|x_1|\leq1\).
	Since \(a\) is bounded on \(\supp\psi\) by \((n+2)+|\omega_2''|^{-1}\), we obtain
	\begin{align*}
		x_1^*(\rho\,\omega'') =\frac{m_1 +\rho^{-1}\int_\Sigma y_1a\psi\,dy'+O(\rho^{-2})} {\displaystyle M_0+\rho^{-1}\int_\Sigma a\psi\,dy'+O(\rho^{-2})}=\frac{m_1}{M_0} +\frac1{\rho M_0}\int_\Sigma \left(y_1-\frac{m_1}{M_0}\right) a(y')\psi(y')\,dy' +O\!\left(\frac1{\rho^{2}}\right).
	\end{align*}
	Substituting the definition of \(a\) gives \eqref{260109}.
	Applying the same expansion to \(\psi_\sigma=\psi\circ\sigma\) and using \eqref{260095} gives \eqref{260110}.

	We next prove \eqref{260111}.
	Fix \(\omega''\) as in the statement and write \(\varepsilon=1/\rho\).
	At \(x''=\rho\omega''\), the equation obtained by setting the integral in \eqref{260094} equal to zero becomes, after multiplying by the positive factor \(\rho^{n+1}\),
	\begin{align}\label{260113}
		\int_\Sigma (x_1-y_1)(\omega_2''-\varepsilon y_2) \bigl(\varepsilon^2(x_1-y_1)^2+|\omega''-\varepsilon y''|^2\bigr)^{-\frac{n+2}{2}} \psi(y')\,dy'=0.
	\end{align}
	For \(|x_1|\leq1\), \(y'\in\supp\psi\), and sufficiently small \(|\varepsilon|\), the integrand and all its \((x_1,\varepsilon)\)-derivatives are bounded by a constant multiple of \(\psi\).
	Hence the left-hand side is smooth near \((m_1/M_0,0)\).
	At \(\varepsilon=0\), it equals \(\omega_2''(M_0x_1-m_1)\), and its \(x_1\)-derivative is \(\omega_2''M_0\neq0\).
	The implicit function theorem and \eqref{260108} show that the solution \(x_1^*(\omega''/\varepsilon)\), initially defined for \(\varepsilon>0\), extends smoothly to \(\varepsilon=0\) with value \(m_1/M_0\).
	Differentiating \eqref{260113} at \(\varepsilon=0\) gives
	\begin{align*}
		\left.\frac{d}{d\varepsilon}x_1^*(\omega''/\varepsilon)\right|_{\varepsilon=0} =\frac1{M_0}\int_\Sigma \left(y_1-\frac{m_1}{M_0}\right) \left((n+2)\omega''\cdot y''-\frac{y_2}{\omega_2''}\right) \psi(y')\,dy'.
	\end{align*}
	The derivative of this smooth extension is Lipschitz near \(0\).
	Since \(d\varepsilon/d\rho=-\rho^{-2}\), the first formula in \eqref{260111} follows.
	Applying this calculation to \(\psi_\sigma=\psi\circ\sigma\), which satisfies {\rm\textbf{(S)}}, and using \eqref{260095} gives the second formula.
If the leading integral is nonzero, the corresponding derivative has a fixed sign for all sufficiently large \(\rho\), proving the stated monotonicity.
	Taking \(\omega''=\mathbf e_2'\) and using
	\begin{align*}
		\int_\Sigma\left(y_1-\frac{m_1}{M_0}\right)y_2\psi(y')\,dy' =\int_\Sigma\left(y_2-\frac{m_2}{M_0}\right)y_1\psi(y')\,dy'
	\end{align*}
	gives the monotonicity statement in Theorem \ref{thm-260014}{\rm(iv)} on the positive axes.
	On the negative axes, take \(\omega''=-\mathbf e_2'\) and write \(s=-\rho\).
	The chain rule gives the same sign for the derivatives with respect to \(s\), opposite to the sign of the integral in Theorem \ref{thm-260014}{\rm(iv)}.
	Increasing \(R_1\) if necessary gives the final claim.
\end{proof}

\subsection{Normal velocity component}\label{260114}
By \eqref{0421-01}, \(\partial_{x_n}w_n(x',0,t)=0\) for \(|x'|>1\), so we consider \(\partial_{x_n}^2w_n(x',0,t)\) instead.

\subsubsection{Formula for the second normal derivative of \(w_n\)}\label{260115}
We first derive the formula for \(\partial_{x_n}^2w_n(x',0,t)\).

\begin{proposition}\label{prop-0703-03}
Let \(1\leq j\leq n-1\), and let \(w\) solve the Stokes system \eqref{0308-01} with boundary data \eqref{0324-13}, where \(\psi\) satisfies {\rm\textbf{(S)}} and \(\phi\) satisfies {\rm\textbf{(T1)}}.
	Then, for every \(x'\in\Sigma\) and \(t>0\),
	\begin{align}\label{0703-11}
		\begin{split}
			\partial_{x_n}^2w_n(x',0,t) & = \phi(t)R_j'(\Delta'\psi)(x') + \frac1{\sqrt{\pi}} \int_0^t\frac{1}{\sqrt{t-s}} \int_{\Sigma} \partial_{x_j}\Gamma'(x'-y',t-s) \left[ \psi(y')\phi'(s)-\phi(s)\Delta'\psi(y') \right]\,dy'ds,
		\end{split}
	\end{align}
where the spatial integral in the second term is interpreted using the cancellation of \(\partial_{x_j}\Gamma'\).
	If \(|x'|>1\), then
	\begin{align}\label{0703-12}
		\begin{split}
			\partial_{x_n}^2w_n(x',0,t) & = \frac{n\Gamma\left(\frac n2\right)}{\pi^{\frac n2}} \phi(t) \int_{\Sigma} \frac{x_j-y_j}{|x'-y'|^{n+2}}\psi(y')\,dy' \\
			 & \quad+ \frac12 \int_0^t\frac{1}{(t-s)^2} \int_{\Sigma} (x_j-y_j)\Gamma(x'-y',0,t-s) \psi(y')\phi(s)\,dy'ds.
		\end{split}
	\end{align}
	For \(t>1\) and every \(x'\in\Sigma\), \(\partial_{x_n}^2w_n(x',0,t)\) is the second line of \eqref{0703-12}, with the time integral over \(0<s<1\).
\end{proposition}
\begin{proof}
	We first prove \eqref{0703-11}.
	Fix \(t>0\).
	By {\rm\textbf{(T1)}} and the fact that \(\phi=0\) on \([1,\infty)\) with \(\phi'(1)=0\), we have \( \phi\in C^{1,\beta}([0,\max\{1,t\}]). \)
	Lemma \ref{lem-260235} therefore applies with \( \alpha_0:=\min\{\alpha,2\beta\}\in(0,1)\), and yields
	\begin{align}\label{260116}
		\left[ \partial_{x_n}^2k(\cdot,t) \right]_{C^{\alpha_0}(\overline{\R^n_+})} \leq C_t \|\psi\|_{C^{2,\alpha}(\Sigma)} \|\phi\|_{C^{1,\beta}([0,\max\{1,t\}])}.
	\end{align}

	The cancellation argument used in the proof of Proposition \ref{prop-0703-02} applies with \( \partial_{x_n}^2k \) in place of \( \partial_{x_n}k. \)
		Using \eqref{260116}, the bound \( |\partial_{y_n}\partial_{y_j}N(y)| \leq Cy_n|y|^{-n-1} \), and \( \int_{\Sigma\times(0,x_n)} \frac{y_n}{|y|^{n+1-\alpha_0}}\,dy \leq Cx_n^{\alpha_0}\), the principal value over \(D_{\varepsilon,x_n}\) as in \eqref{0709-02} exists and satisfies
	\begin{align}\label{260117}
		\left| \mathrm{p.v.} \int_{\Sigma\times(0,x_n)} \partial_{y_n}\partial_{y_j}N(y)\, \partial_{x_n}^2k(x-y,t)\,dy \right|\leq C \left[ \partial_{x_n}^2k(\cdot,t) \right]_{C^{\alpha_0}(\overline{\R^n_+})} x_n^{\alpha_0}.
	\end{align}

		For \(0<\varepsilon<x_n\), the upper face \(\{y_n=x_n\}\) does not meet the removed ball \(\{|y|\leq\varepsilon\}\).
	Therefore differentiating the principal-value volume term in \eqref{0421-04} with respect to \(x_n\) gives the boundary integral \(\int_{\Sigma} \partial_{y_n}\partial_{y_j}N(y',x_n) \partial_{x_n}k(x'-y',0,t)\,dy'\) plus the principal-value integral involving \(\partial_{x_n}^2k\), where the cancellation of \(\partial_{y_n}\partial_{y_j}N\) and \eqref{260116} make the convergence locally uniform in \(x_n>0\).

	It follows that
	\begin{align}\label{0703-13}
		\begin{split}
			\partial_{x_n}^2w_n(x,t) & = -2\phi(t) \int_{\Sigma} \partial_{x_n}^2N(y',x_n)\, \partial_{x_j}\psi(x'-y')\,dy' +4\int_{\Sigma} \partial_{y_j}\partial_{y_n}N(y',x_n) \partial_{x_n}k(x'-y',0,t)\,dy' \\
			 & \quad +4\mathrm{p.v.} \int_{\Sigma\times(0,x_n)} \partial_{y_n}\partial_{y_j}N(y)\, \partial_{x_n}^2k(x-y,t)\,dy.
		\end{split}
	\end{align}
	By \eqref{0713-29}, \( \partial_{x_n}k(\cdot,0,t)\in L^\infty(\Sigma) \) and \( \partial_{x_j}\partial_{x_n}k(\cdot,0,t)\in C_0(\Sigma). \)
	Moreover, for each fixed \(x_n>0\), the kernels \( \partial_{y_n}N(\cdot,x_n) \) and \( \partial_{y_j}\partial_{y_n}N(\cdot,x_n) \) belong to \(L^1(\Sigma)\), so integration by parts in \(y'\) gives
	\begin{align}\label{0713-32}
		\int_{\Sigma} \partial_{y_j}\partial_{y_n}N(y',x_n) \partial_{x_n}k(x'-y',0,t)\,dy' = \int_{\Sigma} \partial_{y_n}N(y',x_n)\, \partial_{x_j}\partial_{x_n}k(x'-y',0,t)\,dy'.
	\end{align}
	Using \( \partial_{x_n}^2N(y',x_n) = -\Delta'N(y',x_n) \) away from the origin and integrating by parts in \(y'\), we obtain
	\begin{align*}
		-2\phi(t) \int_{\Sigma} \partial_{x_n}^2N(y',x_n)\, \partial_{x_j}\psi(x'-y')\,dy' = 2\phi(t) \int_{\Sigma} \partial_{y_j}N(y',x_n)\, \Delta'\psi(x'-y')\,dy'.
	\end{align*}
	Combining this identity with \eqref{0713-32}, \eqref{0703-13} becomes
	\begin{align}\label{0703-14}
		\begin{split}
			\partial_{x_n}^2w_n(x,t) & = 2\phi(t) \int_{\Sigma} \partial_{y_j}N(y',x_n)\, \Delta'\psi(x'-y')\,dy' +4\int_{\Sigma} \partial_{y_n}N(y',x_n)\, \partial_{x_j}\partial_{x_n}k(x'-y',0,t)\,dy' \\
			 & \quad +4\mathrm{p.v.} \int_{\Sigma\times(0,x_n)} \partial_{y_n}\partial_{y_j}N(y)\, \partial_{x_n}^2k(x-y,t)\,dy.
		\end{split}
	\end{align}

	By \eqref{260117}, the last term in \eqref{0703-14} converges to zero as \(x_n\to0^+\).
	For the first term, since \(\Delta'\psi\in C_c^\alpha(\Sigma)\) by {{\rm{\textbf{(S)}}}}, the limit \eqref{260266} gives
	\begin{align}\label{0709-04}
		2\phi(t) \lim_{x_n\to0^+} \int_{\Sigma} \partial_{y_j}N(y',x_n)\, \Delta'\psi(x'-y')\,dy' = \phi(t)R_j'(\Delta'\psi)(x'),
	\end{align}
	and, for the second term, since \(\partial_{x_j}\partial_{x_n}k(\cdot,0,t)\in C_0(\Sigma)\), the limit \eqref{260267} yields
	\begin{align}\label{0709-05}
		4\lim_{x_n\to0^+} \int_{\Sigma} \partial_{y_n}N(y',x_n)\, \partial_{x_j}\partial_{x_n}k(x'-y',0,t)\,dy' = 2\partial_{x_j}\partial_{x_n}k(x',0,t).
	\end{align}
	Taking \(x_n\to0^+\) in \eqref{0703-14}, we therefore obtain
	\begin{align}\label{0703-15}
		\partial_{x_n}^2w_n(x',0,t) = \phi(t)R_j'(\Delta'\psi)(x') + 2\partial_{x_j}\partial_{x_n}k(x',0,t).
	\end{align}
	
	By \eqref{0713-28} and \eqref{0709-07},
	\begin{align*}
		2\partial_{x_j}\partial_{x_n}k(x',0,t) = \frac1{\sqrt{\pi}} \int_0^t\frac1{\sqrt{t-s}} \int_{\Sigma} \partial_{x_j}\Gamma'(x'-y',t-s) \left[ \psi(y')\phi'(s) - \phi(s)\Delta'\psi(y') \right]\,dy'ds,
	\end{align*}
	where the convergence of the integral follows from \eqref{0713-25}.
	Substituting this identity into \eqref{0703-15} gives \eqref{0703-11}.

	We now establish \eqref{0703-12} under {\rm\textbf{(T1)}} alone.
	Fix \(x_0'\in\Sigma\) with \(|x_0'|>1\), set \( d:=\dist(x_0',\supp\psi)>0\), and choose a neighborhood \(U\) of \(x_0'\) such that \( \dist(U,\supp\psi) \geq\frac d2. \)
	For every spatial multi-index \(\gamma\) with \(|\gamma|\leq3\), the definition \eqref{0703-05} gives
	\begin{align}\label{260118}
		\sup_{\substack{x'\in U\\
		0\leq x_n\leq1}} |\partial_x^\gamma k(x,t)| & \leq C_{\gamma,d}\|\psi\|_{L^1}\|\phi\|_{L^\infty(0,t)} \int_0^t u^{-\frac{n+1+|\gamma|}{2}} e^{-\frac{d^2}{Cu}}\,du<\infty\qquad (|\gamma|\leq3).
	\end{align}
			In particular, \(\partial_{x_n}^2k(\cdot,t)\) is Lipschitz near \((x_0',0)\).
	Differentiating \eqref{0713-31} in \(x_n>0\) and integrating \(|\partial_{x_n}\Gamma|\) over \(\Sigma\times(0,t)\) give
	\begin{align}\label{260119}
		\sup_{x\in\R_+^n}|\partial_{x_n}^2k(x,t)| \leq\frac12\sup_{0\leq s\leq t} \|\psi\phi'(s)-\phi(s)\Delta'\psi\|_{L^\infty(\Sigma)}<\infty.
	\end{align}
	Choose \(\eta>0\) so that \(B_{2\eta}'(x_0')\subset U\), and let \(x=(x_0',x_n)\) with \(0<x_n<\eta\).
	In the part \(|y|<\eta\) of the principal-value integral in \eqref{0703-14}, symmetry in \(y_j\) allows us to subtract \(\partial_{x_n}^2k(x,t)\).
	By \eqref{260118},
	\begin{align*}
		\left|\partial_{y_n}\partial_{y_j}N(y) \bigl[\partial_{x_n}^2k(x-y,t)-\partial_{x_n}^2k(x,t)\bigr]\right| \leq C y_n|y|^{-n}.
	\end{align*}
	Its integral over \(0<y_n<x_n\) is at most \(Cx_n\).
	On \(|y|\geq\eta\), \eqref{260119} and \(|\partial_{y_n}\partial_{y_j}N(y)|\leq Cy_n|y|^{-n-1}\) bound the integral by \(C_\eta x_n^2\).
	Thus the principal-value volume integral in \eqref{0703-14} tends to zero as \(x_n\to0^+\).
	By \eqref{260118} and \eqref{260119}, the calculation leading to \eqref{0703-15} is valid at \(x_0'\), and differentiation under the integral gives
	\begin{align*}
		2\partial_{x_j}\partial_{x_n}k(x_0',0,t) 
		 & = \frac12 \int_0^t\frac1{(t-s)^2} \int_\Sigma (x_{0,j}-y_j) \Gamma(x_0'-y',0,t-s) \psi(y')\phi(s)\,dy'ds.
	\end{align*}
	Since \(|x_0'-y'|\geq d>0\) for every \(y'\in\supp\psi\), the integration by parts in the term containing the Riesz transform is also classical and gives
	\begin{align*}
		R_j'(\Delta'\psi)(x_0') = \frac{n\Gamma\left(\frac n2\right)}{\pi^{n/2}} \int_\Sigma \frac{x_{0,j}-y_j}{|x_0'-y'|^{n+2}} \psi(y')\,dy'.
	\end{align*}
	Combining the preceding identities gives \eqref{0703-12}.

	Finally, let \(t>1\) and \(x'\in\Sigma\). Since \(t-s\geq t-1>0\) on \(0\leq s\leq1\), the heat equation gives
	\begin{align*}
	    \frac{d}{ds}\bigl(\partial_{x_j}\Gamma'(\cdot,t-s)*_\Sigma\psi\bigr)(x')
	    =-\bigl(\partial_{x_j}\Gamma'(\cdot,t-s)*_\Sigma\Delta'\psi\bigr)(x').
	\end{align*}
	Integrate the term containing \(\phi'\) in \eqref{0703-11} by parts in \(s\).
	The endpoint terms vanish because \(\phi(0)=\phi(1)=0\), and the two terms containing \(\Delta'\psi\) cancel. Thus
	\begin{align*}
	    \partial_{x_n}^2w_n(x',0,t)
	    =-\frac1{2\sqrt\pi}\int_0^1
	    \frac{\phi(s)}{(t-s)^{3/2}}
	    \bigl(\partial_{x_j}\Gamma'(\cdot,t-s)*_\Sigma\psi\bigr)(x')\,ds .
	\end{align*}
	Substituting \(\partial_{x_j}\Gamma'(z',u)=-z_j\Gamma'(z',u)/(2u)\) gives the stated integral on the whole boundary.
\end{proof}

\subsubsection{Nonexistence of separation points}\label{260120}
The next proposition fixes the sign of \(\partial_{x_n}^2w_n\) on \(\{|x'|>1,\ |x_j|\geq1\}\), so every separation point in \(\{|x'|>1\}\) must lie in \(\{-1<x_j<1\}\).
This proves Theorem \ref{thm-260017}{\rm(i)}.

\begin{proposition}\label{prop-260121}
	Let \(1\leq j\leq n-1\) and let \(w\) solve the Stokes system \eqref{0308-01} with boundary data \eqref{0324-13}, where \(\psi\) satisfies {\rm\textbf{(S)}} and \(\phi\) satisfies {\rm\textbf{(T1)}} and {\rm\textbf{(T2)}}.
	Then, for every \(t>0\) and every \(x'\in\Sigma\) with \(|x'|>1\) and \(|x_j|\geq1\),
	\begin{align*}
		\sgn\bigl(\partial_{x_n}^2w_n(x',0,t)\bigr)=\sgn(x_j).
	\end{align*}
	Every zero \(x'\) of \(\partial_{x_n}^2w_n(\cdot,0,t)\) with \(|x'|>1\) satisfies \(-1<x_j<1\).
	If \(|x'|>1\) and \(|x_j|\geq1\), then \((x',t)\) is not a separation point of \(w_n\).
	If \(n=2\), then
	\begin{align*}
		\sgn\bigl(\partial_{x_2}^2w_2(x_1,0,t)\bigr)=\sgn(x_1), \qquad |x_1|>1,\quad t>0.
	\end{align*}
\end{proposition}
\begin{proof}
	If \(x_j\geq1\), then \(x_j-y_j\geq0\) on \(\supp\psi\), with strict inequality on a subset of positive measure where \(\psi>0\).
	The first term in \eqref{0703-12} is nonnegative, while the second is strictly positive since \(\phi(s)>0\) for \(0<s<\min\{t,1\}\).
	Thus \(\partial_{x_n}^2w_n(x',0,t)>0\) for every \(t>0\).
	When \(x_j\leq-1\), both signs reverse.
	Proposition \ref{prop-260046}{\rm(i)} then excludes separation points with \(|x_j|\geq1\).
\end{proof}

In the remaining subsections, we assume \(j=1\).

\subsubsection{Separation points for radial \(\psi\)}\label{260122}
For radial \(\psi\), the next proposition identifies the set of zeros with \(\{x_1=0,\ |x'|>1\}\) and shows that every point of it is a separation point.
This proves the statement for radial \(\psi\) in Theorem \ref{thm-260017}{\rm(ii)}.

\begin{proposition}\label{prop-260123}
	Let \(n\geq3\) and let \(w\) solve the Stokes system \eqref{0308-01} with boundary data \eqref{0324-13}, where \(j=1\), \(\psi\) satisfies {\rm\textbf{(S)}} and is radial, and \(\phi\) satisfies {\rm\textbf{(T1)}} and {\rm\textbf{(T2)}}.
	Then, for every \(t>0\) and \(x'\in\Sigma\) with \(|x'|>1\),
	\begin{align*}
		\sgn\bigl(\partial_{x_n}^2w_n(x',0,t)\bigr)=\sgn(x_1).
	\end{align*}
	Every point \((x',t)\) with \(|x'|>1\), \(x_1=0\), and \(t>0\) is a separation point of \(w_n\).
\end{proposition}
\begin{proof}
	For each spatial integral in \eqref{0703-12}, form the corresponding vector integral by replacing \(x_1-y_1\) with \(x'-y'\).
	Radial symmetry makes each vector parallel to \(x'\).
	Its scalar product with \(x'\) is positive since \(|x'|^2-x'\cdot y'\geq|x'|(|x'|-|y'|)>0\) on \(\supp\psi\) and \(\psi\not\equiv0\).
			Taking the first components of these vectors in \eqref{0703-12}, with \(j=1\), shows that both spatial integrals in \eqref{0703-12} have the sign of \(x_1\).
		Since \(\phi(t)\geq0\) and \(\phi(s)>0\) for \(0<s<\min\{t,1\}\), \eqref{0703-12} gives the stated sign.

	By \eqref{0421-01}, \(\partial_{x_n}w_n(\cdot,0,t)=0\) on \(|x'|>1\).
	Proposition \ref{prop-260046}{\rm(iii)}, with \(\ell=2\), \(f(x')=x_1\), and the boundary continuity in \eqref{260262}, proves that the points with \(x_1=0\) are separation points.
\end{proof}
\subsubsection{Existence of the separation point}\label{260124}
For general \(\psi\), the next proposition gives a separation point on every line \(\{(x_1,x''):x_1\in\R\}\) with \(|x''|>1\).
This proves the existence statement in Theorem \ref{thm-260017}{\rm(ii)}.

\begin{proposition}\label{prop-260125}
	Let \(n\geq3\) and let \(w\) solve the Stokes system \eqref{0308-01} with boundary data \eqref{0324-13}, where \(j=1\), \(\psi\) satisfies {\rm\textbf{(S)}}, and \(\phi\) satisfies {\rm\textbf{(T1)}} and {\rm\textbf{(T2)}}.
	For every \(t>0\) and every \(x''\in(\mathbf e_1')^\perp\) with \(|x''|>1\), there are \(x_1^*(x'',t)\in(-1,1)\) and \(\delta>0\) such that \(\partial_{x_n}^2w_n(x_1^*(x'',t),x'',0,t)=0\) and
	\begin{align*}
		\partial_{x_n}^2w_n(x',0,t)
		\begin{cases}
			<0, & x_1^*(x'',t)-\delta<x_1<x_1^*(x'',t), \\
			>0, & x_1^*(x'',t)<x_1<x_1^*(x'',t)+\delta.
		\end{cases}
	\end{align*}
The point \((x_1^*(x'',t)\mathbf e_1'+x'',t)\) is a separation point of \(w_n\).
\end{proposition}
\begin{proof}
	Fix \(t>0\) and \(x''\in(\mathbf e_1')^\perp\) with \(|x''|>1\).
	Proposition \ref{prop-260121} gives \(\partial_{x_n}^2w_n(-1,x'',0,t)<0< \partial_{x_n}^2w_n(1,x'',0,t)\).
	The intermediate value theorem therefore gives at least one zero in \((-1,1)\).

	Since \( \inf_{y'\in\supp\psi}|x''-y''|\geq|x''|-1>0\), formula \eqref{0703-12} and Lemma \ref{lem-260044}{\rm(i)}--{\rm(ii)}, applied with \(I=\R\), \(E=\supp\psi\), \(g=\psi\), \(a=\frac{n+2}{2}\), \(T=t\), \(b=2+\frac n2\), \(m=1\), and \(q(s)=\phi(t-s)\), show that \(\partial_{x_n}^2w_n(x',0,t)\) is real-analytic in \(x_1\).
	Its values at \(-1\) and \(1\) have opposite signs, so it is not identically zero.
	Lemma \ref{lem-260044}{\rm(iii)} shows that it has finitely many zeros in \([-1,1]\).
	
	Thus at least one zero has negative values on its left and positive values on its right.
	Denote it by \(x_1^*(x'',t)\).
	This zero is isolated, so there exists \(\delta>0\) with the stated signs.
Proposition \ref{prop-260046}{\rm(i)} therefore shows that \((x_1^*(x'',t)\mathbf e_1'+x'',t)\) is a separation point of \(w_n\).
\end{proof}

\subsubsection{Uniqueness of the separation point for large \(|x''|\)}\label{260126}
For each fixed \(t>0\), the next proposition proves uniqueness in \(x_1\) when \(|x''|\) is large and determines the sign of \(\partial_{x_n}^2w_n(x',0,t)\) on either side of \(x_1^*(x'',t)\).
This proves Theorem \ref{thm-260017}{\rm(iii)}.

\begin{proposition}\label{prop-260127}
	Let \(n\geq3\) and let \(w\) solve the Stokes system \eqref{0308-01} with boundary data \eqref{0324-13}, where \(j=1\), \(\psi\) satisfies {\rm\textbf{(S)}}, and \(\phi\) satisfies {\rm\textbf{(T1)}} and {\rm\textbf{(T2)}}.
	Then there exists \(R>1\), depending only on \(n\) and \(\phi\), such that, for every \(t>0\) and every \(x''\in(\mathbf e_1')^\perp\) with \(|x''|>R\), \(\partial_{x_n}^2w_n(\cdot,x'',0,t)\) has a unique zero \(x_1^*(x'',t)\in(-1,1)\), and
	\begin{align}\label{0928-001}
		\sgn \left( \partial_{x_n}^2w_n(x',0,t) \right) = \sgn\bigl(x_1-x_1^*(x'',t)\bigr).
	\end{align}
	Moreover, \((x_1^*(x'',t)\mathbf e_1'+x'',t)\) is the unique separation point of \(w_n\) on the line \(\{x_1\mathbf e_1'+x'':x_1\in\R\}\), and \(x_1^*(x'',t)\) is continuously differentiable in \((x'',t)\) for \(|x''|>R\) and \(t>0\).
	If \(t\geq3/2\), \(\phi\) satisfies only {\rm\textbf{(T1)}}, and \(\phi\not\equiv0\), then \eqref{0928-001} holds for every \(x''\).
	The unique zero of \(\partial_{x_n}^2w_n(x',0,t)\) in \(x_1\) is simple and gives a separation point whenever it lies outside \(\supp\psi\).
\end{proposition}
\begin{proof}
	We show that \(\partial_{x_1}\partial_{x_n}^2w_n(x',0,t)>0\) at every zero with \(|x''|>R\), where \(R\) is independent of \(t\); the intermediate value theorem then gives uniqueness.

	\underline{The case \(t\geq1\).}\quad
	For \(t>1\), and for \(t=1\) and \(|x''|>1\), Proposition~\ref{prop-0703-03} represents \(\partial_{x_n}^2w_n(x',0,t)\) by the second line of \eqref{0703-12}.
	At a zero, \(x_1\) is therefore the average of \(y_1\) with respect to the positive weight \((t-s)^{-2}\Gamma(x'-y',0,t-s)\psi(y')\phi(s)\), so \(|x_1|\leq1\), and, since \(y_1^2\leq1\),
	\begin{align}\label{0814-02}
		\int_0^1 \int_{\Sigma}\frac{(x_1-y_1)^2}{(t-s)^2}\Gamma(x'-y',0,t-s)\psi(y')\phi(s)\,dy'ds
		\leq \int_0^1\int_{\Sigma}\frac{1}{(t-s)^2}\Gamma(x'-y',0,t-s)\psi(y')\phi(s)\,dy'ds.
	\end{align}
	By \eqref{0703-12},
	\begin{align}\label{260128}
		2\partial_{x_1}\partial_{x_n}^2w_n(x',0,t) & = \int_0^1\int_{\Sigma}\frac{1}{(t-s)^2}\Gamma(x'-y',0,t-s)\psi(y')\phi(s)\,dy'ds\nonumber\\
		 & - \frac12 \int_0^1 \int_{\Sigma}\frac{(x_1-y_1)^2}{(t-s)^3}\Gamma(x'-y',0,t-s)\psi(y')\phi(s)\,dy'ds =: \mathrm{I} - \frac12 \mathrm{II}.
	\end{align}

	Let \(1\leq t<\frac32\), and split \(\mathrm{II}=\mathrm{II}_1+\mathrm{II}_2\) at \(s=\frac14\).
	Since \(t-s\geq\frac34\) for \(s<\frac14\), \eqref{0814-02} gives \(\mathrm{II}_1\leq\frac43\mathrm{I}\).
	For \(s\geq\frac14\), \((x_1-y_1)^2\leq4\), and, once \(R\) is large, \((t-s)^{-\frac n2-3}e^{-\frac{|x'-y'|^2}{4(t-s)}}\) is decreasing in \(s\).
	Comparing with the contribution of \(s\in(\frac1{16},\frac18)\) to \(\mathrm I\) gives
	\begin{align}\label{0814-06}
		\int_{\frac14}^1 \frac{1}{(t-s)^3}\Gamma(x'-y',0,t-s)\phi(s)\,ds
		\leq 2^{\frac{n}{2}+3}\frac{\norm{\phi}_{L^1(0,1)}}{\norm{\phi}_{L^1\left(\frac{1}{16},\frac18\right)}} e^{-\frac{|x'-y'|^2}{55}}\int_0^1\frac{1}{(t-s)^2}\Gamma(x'-y',0,t-s)\phi(s)\,ds.
	\end{align}
	Choosing \(R\) so that the factor in front of the last integral is at most \(\frac1{12}\) for \(|x''|>R\) and \(y'\in\supp\psi\), we obtain \(\mathrm{II}_2\leq\frac13\mathrm{I}\).
	Hence \(\mathrm{II}\leq\frac53\mathrm{I}\), and \eqref{260128} gives \(2\partial_{x_1}\partial_{x_n}^2w_n(x',0,t)\geq\frac16\mathrm{I}>0\).

	Let \(t\geq\frac32\), fix any \(x''\), and assume only {\rm\textbf{(T1)}} and \(\phi\not\equiv0\).
	Since \(1/(t-s)<2\) for \(0<s<1\) and the weight in \eqref{0814-02} puts positive mass on \(\{y_1\ne x_1\}\), it follows from \eqref{0814-02} that \(\mathrm{II}<2\mathrm{I}\) at each zero, and \eqref{260128} shows that the \(x_1\)-derivative is positive there.
	By \eqref{0703-12}, \(\partial_{x_n}^2w_n(x',0,t)\) is negative for \(x_1\leq-1\) and positive for \(x_1\geq1\).
	Hence, for \(t\geq1\), with \(|x''|>R\) if \(t<\frac32\), the zero is unique and lies in \((-1,1)\), and \eqref{0928-001} holds.
	If this zero lies outside \(\supp\psi\), Proposition~\ref{prop-260046}{\rm(i)} shows that it gives a separation point.

	\underline{The case \(0<t<1\).}\quad
	At a zero, \eqref{0703-12} shows that \(x_1\) is the average of \(y_1\) with respect to the sum of the two positive weights in \eqref{0703-12}, so \(|x_1|\leq1\), and, as for \eqref{0814-02},
	\begin{align}\label{0814-09}
		 & \frac{n\Gamma\left(\frac{n}{2}\right)}{\pi^{\frac{n}{2}}}\phi(t)\int_{\Sigma}\frac{(x_1-y_1)^2}{|x'-y'|^{n+2}}\psi(y')\,dy' + \frac12 \int_0^t \int_{\Sigma}\frac{(x_1-y_1)^2}{(t-s)^2}\Gamma(x'-y',0,t-s)\psi(y')\phi(s)\,dy'ds\nonumber \\
		 & \leq \frac{n\Gamma\left(\frac{n}{2}\right)}{\pi^{\frac{n}{2}}}\phi(t)\int_{\Sigma}\frac{1}{|x'-y'|^{n+2}}\psi(y')\,dy' + \frac12 \int_0^t \int_{\Sigma}\frac{1}{(t-s)^2}\Gamma(x'-y',0,t-s)\psi(y')\phi(s)\,dy'ds.
	\end{align}
	By \eqref{0703-12},
	\begin{align}\label{260129}
		\partial_{x_1}\partial_{x_n}^2w_n(x',0,t) & = \frac{n\Gamma\left(\frac{n}{2}\right)}{\pi^{\frac{n}{2}}}\phi(t)\int_{\Sigma}\frac{1}{|x'-y'|^{n+2}}\psi(y')\,dy' + \frac12 \int_0^t \int_{\Sigma}\frac{1}{(t-s)^2}\Gamma(x'-y',0,t-s)\psi(y')\phi(s)\,dy'ds\nonumber\\
		 & - \frac{n\Gamma\left(\frac{n}{2}\right)}{\pi^{\frac{n}{2}}} (n+2)\phi(t)\int_{\Sigma}\frac{(x_1-y_1)^2}{|x'-y'|^{n+4}}\psi(y')\,dy' \nonumber \\
		 & - \frac14 \int_0^t \int_{\Sigma}\frac{(x_1-y_1)^2}{(t-s)^3}\Gamma(x'-y',0,t-s)\psi(y')\phi(s)\,dy'ds =: A+\frac12 B - C - \frac14 D.
	\end{align}

	Let \(0<t<\frac{1+t_1}{2}\).
	By {\rm\textbf{(T2)}}, \(\phi(s)\leq K_\phi\phi(t)\) for \(0<s<t\), with \(K_\phi:=\phi(t_1)/\phi(\frac{1+t_1}{2})\geq1\).
	For \(|x''|\geq2\), with \(a=(|x''|-1)^2\) and \(u=t-s\),
	\begin{align}\label{0814-10}
		\left|\frac12 B- \frac14 D\right| & \leq \frac{K_\phi M_0\phi(t)}{(4\pi)^{\frac n2}} \int_0^t\left(\frac1{2u^{\frac n2+2}}+\frac1{u^{\frac n2+3}}\right)e^{-\frac{a}{4u}}\,du \leq C_nK_\phi M_0\phi(t)e^{-\frac a8}.
	\end{align}
	For \(|x''|>1+2\sqrt{n+2}\), since \(|x_1|\leq1\) and \(|x'-y'|\leq|x''|+2\) on \(\supp\psi\),
	\begin{align}\label{0814-11}
		A-C \geq\frac{n\Gamma\left(\frac{n}{2}\right)}{\pi^{\frac{n}{2}}}\phi(t) \left(1-\frac{4(n+2)}{(|x''|-1)^{2}}\right)\frac{M_0}{(|x''|+2)^{n+2}}.
	\end{align}
	Hence \(|B/2-D/4|\leq(A-C)/2\) for \(|x''|>R\) with \(R\) large, and \eqref{260129} gives \(\partial_{x_1}\partial_{x_n}^2w_n(x',0,t)\geq(A-C)/2>0\).

	Let \(\frac{1+t_1}{2}\leq t<1\).
	For sufficiently large \(R\), \(C\leq\frac{t_1}{2(2+t_1)}A\).
	Split \(D=D_1+D_2\) at \(s=\frac{t_1}4\).
	Since \(t-s\geq\frac{2+t_1}{4}\) for \(s<\frac{t_1}4\), \eqref{0814-09} gives \(D_1\leq\frac{4}{2+t_1}(2A+B)\).
	For \(s\geq\frac{t_1}{4}\), the monotonicity and comparison used for \eqref{0814-06}, now with \(s\in(\frac{t_1}{16},\frac{t_1}8)\), give
	\begin{align*}
		D_2 \leq 2^{\frac{n}{2}+5}\frac{\norm{\phi}_{L^1(0,1)}}{\norm{\phi}_{L^1\left(\frac{t_1}{16},\frac{t_1}{8}\right)}} e^{-\frac{t_1(|x''|-1)^2}{32}}B\leq\frac{t_1}{2+t_1}B
	\end{align*}
	for \(|x''|>R\) with \(R\) large.
	Hence \(C+\frac14D\leq\frac{t_1+4}{2(2+t_1)}A+\frac{4+t_1}{4(2+t_1)}B<A+\frac12B\), and \eqref{260129} gives \(\partial_{x_1}\partial_{x_n}^2w_n(x',0,t)>0\).
	In both cases, \eqref{0928-001} now follows from Proposition \ref{prop-260125}.

	\underline{Differentiability of \(x_1^*\).}\quad
	Let \(U\Subset\{|x'|>1\}\) and \(0<T<\infty\).
	For \(x'\in U\), \(y'\in\supp\psi\), \(0<u<T\), and \(k<n\),
	\begin{align}\label{1007-17}
		 & \left|\frac{x_1-y_1}{u^2}\Gamma(x'-y',0,u)\right| +\left|\partial_{x_k}\left(\frac{x_1-y_1}{u^2} \Gamma(x'-y',0,u)\right)\right|
		 \leq \frac{C_{U,T}}{u^{\frac n2+3}} \exp\left(-\frac{\dist(U,\supp\psi)^2}{8u}\right) \in L^1((0,T),du).
	\end{align}
	Writing the time integral in \eqref{0703-12} with \(u=t-s\), \eqref{1007-17} and \(\phi(0)=0\) give
	\begin{align}\label{1007-18}
		\partial_t\partial_{x_n}^2w_n(x',0,t) & =\frac{n\Gamma\left(\frac n2\right)}{\pi^{\frac n2}}\phi'(t) \int_\Sigma\frac{x_1-y_1}{|x'-y'|^{n+2}}\psi(y')\,dy'
		+\frac12\int_0^t\int_\Sigma \frac{x_1-y_1}{u^2}\Gamma(x'-y',0,u)\psi(y')\phi'(t-u)\,dy'du,
	\end{align}
	and dominated convergence shows that \(\partial_{x_n}^2w_n(x',0,t)\) is \(C^1\) on \(\{|x'|>1\}\times(0,\infty)\).
	Since \(\partial_{x_1}\partial_{x_n}^2w_n(x_1^*(x'',t),x'',0,t)>0\) for \(|x''|>R\) and \(t>0\), the implicit function theorem and the uniqueness of the zero show that \(x_1^*\) is \(C^1\).
\end{proof}
\subsubsection{Asymptotics of the separation point}\label{260130}

The next proposition estimates the unique separation point for large \(|x''|\) for \(0<t<1\) and \(t\geq1\).
The three statements prove Theorem \ref{thm-260017}{\rm(iv)(a)}, {\rm(iv)(b)}, and {\rm(iv)(c)}, respectively.
\begin{proposition}
	\label{prop-260131}
	Let \(n\geq3\) and let \(w\) solve the Stokes system \eqref{0308-01} with boundary data \eqref{0324-13}, where \(j=1\), \(\psi\) satisfies {\rm\textbf{(S)}}, and \(\phi\) satisfies {\rm\textbf{(T1)}} and {\rm\textbf{(T2)}}.
	Then there exists \(R\geq2\), depending only on \(n\), \(\psi\), and \(\phi\), such that, for \(|x''|>R\), the following hold.
	\begin{enumerate}
				\item\label{item-0817-12} For every \(0<t<1\), with an absolute constant \(c>0\),
		\begin{align}\label{260132}
			\left|x_1^*(x'',t)-\frac{m_1}{M_0}\right| \leq\frac{n+2}{|x''|} +C_n\frac{\phi(t_1)}{\phi(\max\{t,t_1\})}e^{-c|x''|^2},
		\end{align}
		and
		\begin{align}\label{260133}
			\left|x_1^*(x'',t)-\frac{m_1}{M_0} -\frac{n+2}{M_0|x''|^2}x''\cdot\int_\Sigma \left(y_1-\frac{m_1}{M_0}\right)y''\psi(y')\,dy'\right| \leq\frac{C_{n,\psi}}{|x''|^2} +C_n\frac{\phi(t_1)}{\phi(\max\{t,t_1\})}e^{-c|x''|^2}.
		\end{align}
		\item\label{item-0817-11} Let \(t\geq1\) and \(\omega''\in\mathbb S^{n-3}\), and set \(F_{\omega''}:=\left\{y'\in\supp\psi:\omega''\cdot y''=\max_{z'\in\supp\psi}\omega''\cdot z''\right\}\).
		If \(\{y_1:y'\in F_{\omega''}\}=\{p_1\}\) for some \(p_1\), then
		\begin{align}\label{0817-11}
			\lim_{r\to\infty}x_1^*(r\omega'',t)=p_1.
		\end{align}
		\item\label{item-0817-13} For every \(t\geq1\),
		\begin{align}\label{260134}
			\left|x_1^*(x'',t)-\frac{m_1}{M_0}\right| \leq2\min\left\{1,\frac{1+|x''|}{t}\right\}.
		\end{align}
	\end{enumerate}
\end{proposition}

The next lemma is used to prove \eqref{0817-11}.
\begin{lemma}\label{lem-0817-01}
	Let \(n\geq3\), \(t\geq1\), let \(\phi\colon[0,\infty)\to\R\) satisfy {\rm\textbf{(T1)}} and {\rm\textbf{(T2)}}, and define, for \(\lambda>0\),
	\begin{align*}
		 & \mathcal T(\lambda) :=\int_0^1 \frac{\phi(s)}{(t-s)^{\frac{n}{2}+2}}e^{-\frac{\lambda}{4(t-s)}}\,ds,\\
		 & X(\lambda):=e^{\frac{\lambda}{4t}}\mathcal T(\lambda)=\int_0^1\frac{\phi(s)}{(t-s)^{\frac n2+2}}e^{-\frac{s\lambda}{4t(t-s)}}\,ds, \qquad D(\lambda):=-\frac{d}{d\lambda}\ln X(\lambda).
	\end{align*}
	Then the following hold.
	\begin{enumerate}
		\item\label{item-0817-02} If, in addition, \(\psi\) satisfies {\rm\textbf{(S)}} and \(w\) solves the Stokes system \eqref{0308-01} with boundary data \eqref{0324-13}, where \(j=1\), then, for \(|x'|>1\),
		\begin{align*}
			\partial_{x_n}^2w_n(x',0,t)=\frac{1}{2(4\pi)^{\frac n2}}\int_{\Sigma}(x_1-y_1)\mathcal T\bigl(|x'-y'|^2\bigr)\psi(y')\,dy'.
		\end{align*}
		\item\label{item-0817-03} \(X\) is positive, finite, smooth, and strictly decreasing on \((0,\infty)\).
		\item\label{item-0817-04} \(X\) is log-convex on \((0,\infty)\).
		\item\label{item-0817-06} \(D(\lambda)\to0\) as \(\lambda\to\infty\).
		\item\label{item-0817-07} For every \(\lambda>0\) and \(\delta\geq0\), \(e^{-D(\lambda)\delta}X(\lambda)\leq X(\lambda+\delta)\leq X(\lambda)\).
	\end{enumerate}
\end{lemma}
\begin{proof}
	\underline{Proof of (i).}\quad
	Substituting the heat kernel into \eqref{0703-12} and using \(\phi(t)=0\) for \(t\geq1\) gives the stated formula.

	\underline{Proof of (ii).}\quad
	Differentiation under the integral gives, for \(k\geq0\),
	\begin{align*}
		|X^{(k)}(\lambda)|\leq \int_0^1\frac{\phi(s)}{(t-s)^{\frac n2+2}} \left(\frac{s}{4t(t-s)}\right)^k e^{-\frac{\lambda s}{4t(t-s)}}\,ds<\infty.
	\end{align*}
	Then {\rm(ii)} is immediate.

	\underline{Proof of (iii).}\quad
	H\"older's inequality gives, for \(\lambda_1,\lambda_2>0\) and \(\theta\in(0,1)\),
	\begin{align*}
		X(\theta \lambda_1 + (1-\theta)\lambda_2) & = \int_0^1 \frac{\phi(s)}{(t-s)^{\frac{n}{2}+2}}e^{-\frac{s\theta \lambda_1}{4t(t-s)}}e^{-\frac{s(1-\theta)\lambda_2}{4t(t-s)}}\,ds\\
		 & \leq \left( \int_0^1 \frac{\phi(s)}{(t-s)^{\frac{n}{2}+2}}e^{-\frac{s\lambda_1}{4t(t-s)}}\,ds\right)^\theta \left( \int_0^1 \frac{\phi(s)}{(t-s)^{\frac{n}{2}+2}}e^{-\frac{s\lambda_2}{4t(t-s)}}\,ds\right)^{1-\theta} = X(\lambda_1)^\theta X(\lambda_2)^{1-\theta},
	\end{align*}
	which gives {\rm(iii)}.

	\underline{Proof of (iv).}\quad
	By {\rm(ii)}, \(D\geq0\). Fix \(\varepsilon\in(0,1)\).
	Splitting the integral for \(-X'(\lambda)\) at \(\varepsilon\), and bounding \(X(\lambda)\) below by its integral over \((0,\varepsilon/2)\), gives, for \(\lambda\geq1\),
	\begin{align*}
		0\leq D(\lambda) & \leq\frac{\varepsilon}{4t(t-\varepsilon)} +\exp\left(-\frac{(\lambda-1)\varepsilon}{4t(t-\varepsilon)} +\frac{\lambda\varepsilon}{8t(t-\varepsilon/2)}\right) \frac{\displaystyle\int_\varepsilon^1 \frac{s\phi(s)}{4t(t-s)^{\frac n2+3}} e^{-\frac{s}{4t(t-s)}}\,ds} {\displaystyle\int_0^{\varepsilon/2} \frac{\phi(s)}{(t-s)^{\frac n2+2}}\,ds}.
	\end{align*}
	Since \(\frac{\varepsilon}{4t(t-\varepsilon)}> \frac{\varepsilon}{8t(t-\varepsilon/2)}\),
	\begin{align*}
		0\leq\limsup_{\lambda\to\infty}D(\lambda) \leq\frac{\varepsilon}{4t(t-\varepsilon)}\longrightarrow0 \qquad\text{as }\varepsilon\to0^+.
	\end{align*}

	\underline{Proof of (v).}\quad
	Since \(X\) is decreasing and log-convex by {\rm(ii)} and {\rm(iii)}, 
	\begin{align*}
		-\delta D(\lambda)\leq\ln X(\lambda+\delta)-\ln X(\lambda)\leq0.
	\end{align*}
	Exponentiating gives the two bounds.
\end{proof}

\begin{proof}[Proof of Proposition \ref{prop-260131}]
	We take \(E:=\supp\psi\) in Lemma \ref{lem-260106}.
	For \(t>0\), setting the right-hand side of \eqref{0703-12} equal to zero at \(x_1=x_1^*(x'',t)\) expresses \(x_1^*(x'',t)\) as a weighted average of \(y_1\) with a positive denominator, so
	\begin{align}\label{0817-19}
		|x_1^*(x'',t)|\leq1,
	\end{align}
	and, for \(0<s<t<1\), {\rm\textbf{(T2)}} gives
	\begin{align}\label{260016}
		\frac{\phi(s)}{\phi(t)}\leq\frac{\phi(t_1)}{\phi(\max\{t,t_1\})}.
	\end{align}

	\underline{Proof of (i).}\quad
	For \(k=0,1\), let
	\begin{align*}
		 & S_k(x',t) :=\frac{n\Gamma\left(\frac{n}{2}\right)}{\pi^{\frac{n}{2}}}\phi(t) \int_{\Sigma}\frac{y_1^k}{|x'-y'|^{n+2}}\psi(y')\,dy', \\
		 & E_k(x',t) :=\frac12 \int_0^t \int_{\Sigma}\frac{y_1^k}{(t-s)^2}\Gamma(x'-y',0,t-s)\psi(y')\phi(s)\,dy'ds.
	\end{align*}
	By \eqref{0703-12}, \(x_1^*(x'',t)=(S_1+E_1)/(S_0+E_0)\), with \(S_k\) and \(E_k\) evaluated at \(x_1=x_1^*(x'',t)\), so
	\begin{align}\label{0817-08}
		x_1^*(x'',t)-\frac{m_1}{M_0} & =\frac{S_0}{S_0+E_0}\left(\frac{S_1}{S_0}-\frac{m_1}{M_0}\right) +\frac{E_0}{S_0+E_0}\left(\frac{E_1}{E_0}-\frac{m_1}{M_0}\right).
	\end{align}
	Since \(|x'|\geq|x''|\geq2\) and \(|x'|-1\leq |x'-y'|\leq |x'|+1\) on \(\supp\psi\), Lemma \ref{lem-260106} applies with \(f(y')=y_1\), \(g=\psi\), \(\varpi(y') = |x'-y'|^{-n-2}\), \(L=2\), and \(\Lambda = (n+2)\ln\frac{|x''|+1}{|x''|-1}\), and, since \(\frac12\ln\frac{|x''|+1}{|x''|-1}=\tanh^{-1}\frac{1}{|x''|}\) and \(\tanh(\frac{(n+2)u}{2})\leq\frac{n+2}{2}\tanh u\) for \(u>0\), gives
	\begin{align}\label{0817-03}
		\left|\frac{S_1}{S_0}-\frac{m_1}{M_0}\right| \leq2\tanh\left(\frac{n+2}{4}\ln\left(\frac{|x''|+1}{|x''|-1}\right)\right)\leq \frac{n+2}{|x''|}.
	\end{align}
	Since \(|y_1|\leq1\) on \(\supp\psi\), \(|E_1/E_0-m_1/M_0|\leq2\).
	Since \(|x''|-1\leq|x'-y'|\leq |x''|+2\) on \(\supp\psi\), \eqref{260016} and \(e^{-A/(4s)}\leq e^{-A/8}e^{-A/(8s)}\) for \(A>0\) and \(0<s<1\) give
	\begin{align}\label{0817-05}
		\frac{E_0}{S_0} & \leq\frac{(|x''|+2)^{n+2}}{2^{n+1}n\Gamma\left(\frac n2\right)} \frac{\phi(t_1)}{\phi(\max\{t,t_1\})} e^{-\frac{(|x''|-1)^2}{8}} \int_0^\infty\frac{e^{-\frac{(|x''|-1)^2}{8s}}}{s^{\frac n2+2}}\,ds\nonumber\\
		& \leq2^{\frac{5n+10}{2}} \frac{\phi(t_1)}{\phi(\max\{t,t_1\})}e^{-\frac{(|x''|-1)^2}{8}}.
	\end{align}
	Hence \eqref{0817-08} and \((|x''|-1)^2\geq |x''|^2/4\) give \eqref{260132} with \(C_n=2^{\frac{5n+12}{2}}\) and \(c=\frac{1}{32}\).

	Since \(|x'-y'|^2=|x''|^2-2x''\cdot y''+O(1)\) uniformly for \(y'\in\supp\psi\) and \(|x_1|\leq1\), for \(k=0,1\),
	\begin{align*}
		S_k(x',t) = \frac{n\Gamma\left(\frac{n}{2}\right)}{\pi^{\frac{n}{2}}}\frac{\phi(t)}{|x''|^{n+2}} \left(\int_\Sigma y_1^k\psi(y')\,dy' + \frac{n+2}{|x''|^{2}}x''\cdot \int_{\Sigma}y_1^ky'' \psi(y')\,dy' + O\left(\frac{1}{|x''|^{2}}\right)\right),
	\end{align*}
	so
	\begin{align*}
		\frac{S_1(x',t)}{S_0(x',t)}-\frac{m_1}{M_0} & =\frac{n+2}{M_0|x''|^2}x''\cdot\int_\Sigma \left(y_1-\frac{m_1}{M_0}\right)y''\psi(y')\,dy'+O_{\psi,n}\left(\frac1{|x''|^2}\right).
	\end{align*}
	By \eqref{0817-08}, \(|x_1^*(x'',t)-S_1/S_0| =\frac{E_0}{S_0+E_0}|E_1/E_0-S_1/S_0| \leq 2E_0/S_0\), and \eqref{0817-05} proves \eqref{260133}.

	\underline{Proof of (ii).}\quad
	Let \(x''=r\omega''\) and \(h=\max_{z'\in\supp\psi}\omega''\cdot z''\).
	Writing \(y''=q\omega''+z\) with \(z\cdot\omega''=0\), \(u=u(y'):=h-q\), and \(\rho:=r-h\), we have \(|x'-y'|^2 = \rho^2+2\rho u + c\), where \(c=c(x_1, y'):= (x_1-y_1)^2 + u^2 + |z|^2\in[0,9]\) for \(|x_1|\leq1\), since \(0\leq u\leq2\) and \(|z|\leq1\).
	By Lemma \ref{lem-0817-01}{\rm(i)},
	\begin{align}\label{0817-12}
		\partial_{x_n}^2 w_n(x',0,t) = \frac{X(\rho^2)}{2(4\pi)^{\frac{n}{2}}}e^{-\frac{\rho^2}{4t}}\int_{\Sigma}e^{-\frac{2u\rho + c }{4t}}(x_1-y_1)\psi(y')R_\rho (x_1, y')\,dy',
	\end{align}
	where \(R_\rho(x_1, y') := X(\rho^2+2\rho u + c)/X(\rho^2)\), and Lemma \ref{lem-0817-01}{\rm(v)} gives \(e^{-D(\rho^2)(2\rho u + c)} \leq R_\rho(x_1, y') \leq 1\).

	Fix \(\delta\in(0,1)\).
	Since \(\{u<\delta/2\}\) is an open neighborhood of a point of \(F_{\omega''}\subset\supp\psi\), \(\int_{u<\delta/2}\psi(y')\,dy'>0\).
	By Lemma \ref{lem-0817-01}{\rm(iv)}, there is \(\rho_\delta>0\) such that \(D(\rho^2)\leq \frac{1}{8t}\) for \(\rho\geq\rho_\delta\), and the bounds on \(R_\rho\) give
	\begin{align}\label{0817-14}
		\frac{ \int_{u\geq\delta} e^{-\frac{2\rho u + c}{4t}}R_\rho(x_1, y') \psi(y')\,dy'}{\int_{u<\frac{\delta}{2}} e^{-\frac{2\rho u + c}{4t}}R_\rho(x_1, y') \psi(y')\,dy'}
		& \leq e^{\frac{9}{4t}+ 9 D(\rho^2) - \rho \delta \left(\frac{1}{4t}-D(\rho^2)\right) }\frac{M_0}{\int_{u<\frac{\delta}{2}}\psi(y')\,dy'}\nonumber\\
		& \leq C_\delta e^{-\frac{\rho\delta}{8t}},\qquad \rho\geq\rho_\delta,
	\end{align}
	where \(C_\delta\) depends on \(t,\delta,\psi,\omega''\), and \(\rho_\delta\) on \(\delta,n,t,\phi\).

	Fix \(\varepsilon \in (0,1)\).
	Since \(\supp\psi\) is compact and \(\{y_1:y'\in F_{\omega''}\}=\{p_1\}\), there is \(\delta\in (0,1)\) such that \(|y_1-p_1|<\frac{\varepsilon}{2}\) if \(y'\in \supp\psi\) and \(u(y')<\delta\).
	Let \(x_1\in[p_1+\varepsilon,1]\), and split the integral in \eqref{0817-12} into \(I_1\) over \(\{u<\delta\}\) and \(I_2\) over \(\{u\geq\delta\}\).
	Since \(x_1-y_1\geq\varepsilon/2\) on \(\supp\psi\cap\{u<\delta\}\) and \(|x_1-y_1|\leq2\), \eqref{0817-14} gives, for \(\rho\geq\rho_\varepsilon\geq\rho_\delta\) with \(2C_\delta e^{-\frac{\rho_\varepsilon\delta}{8t}}\leq \frac{\varepsilon}{4}\),
	\begin{align*}
		I_1 \geq \frac{\varepsilon}{2}\int_{u<\frac{\delta}{2}} e^{-\frac{2u\rho + c }{4t}}\psi(y')R_\rho (x_1, y')\,dy'\geq2|I_2|,
	\end{align*}
	and \(I_1>0\), so \(\partial_{x_n}^2 w_n(x',0,t)>0\).
	Similarly, \(\partial_{x_n}^2 w_n(x',0,t)<0\) for \(x_1\in[-1,p_1-\varepsilon]\) and \(\rho\geq\rho_\varepsilon\).
	By \eqref{0817-19}, \(|x_1^*(r\omega'',t)-p_1|<\varepsilon\) for \(r\geq h+\rho_\varepsilon\), which proves \eqref{0817-11}.

	\underline{Proof of (iii).}\quad
	At \(t=1\), \eqref{260134} follows from \eqref{0817-19} and \(|m_1|\leq M_0\).
	For \(t>1\), \eqref{0703-12} gives \(x_1^*(x'',t)=E_1(x',t)/E_0(x',t)\), and Lemma \ref{lem-260106} applies with \(f(y')=y_1\), \(g=\psi\), and
	\begin{align*}
		\varpi(y')=\frac1{(4\pi)^{\frac n2}}\int_0^1 \frac{\phi(s)}{(t-s)^{\frac n2+2}} e^{-\frac{|x'-y'|^2}{4(t-s)}}\,ds.
	\end{align*}
	Since \(\bigl||x'-y'|^2-|x'-z'|^2\bigr| \leq1+4\sqrt{1+|x''|^2}\) for \(y',z'\in\supp\psi\), Lemma \ref{lem-260106} and \(\tanh r\leq\min\{1,r\}\) give
	\begin{align*}
		\left|x_1^*(x'',t)-\frac{m_1}{M_0}\right| \leq2\min\left\{1, \frac{1+4\sqrt{1+|x''|^2}}{16(t-1)}\right\}.
	\end{align*}
	This gives \eqref{260134}, since its right-hand side is \(2\) if \(t\leq1+|x''|\), and, if \(t>1+|x''|\), then \((t-1)/t\geq2/3\) and the last fraction is at most \((1+|x''|)/t\).
\end{proof}

\section{Tangential pressure derivatives}\label{260135}
In this section, we prove Theorems \ref{thm-260019}, \ref{thm-260021}, and \ref{thm-260027}.

For the Stokes system \eqref{0308-01} with boundary data \(g=(g_1,\dots,g_n)\), we choose the pressure representative given by \cite{Sol231,KangLaiLaiTsai2023}
\begin{align}\label{0418-01}
	p(x,t) & = -4\sum_{j=1}^n \partial_{x_j}\int_{\Sigma}\partial_{x_n} N(x-y')g_j(y',t)\,dy' - 2\partial_t \int_{\Sigma}N(x-y')g_n(y',t)\,dy' \nonumber\\
	& \quad - 2 \sum_{j=1}^n \int_{-\infty}^t\int_{\Sigma}\frac{g_j(z',s)-g_j(z',t)}{t-s}\int_{\Sigma}\Gamma(x'-y'-z',0,t-s)\partial_{y_j}N(y',x_n)\,dy'dz'ds,
\end{align}
where the boundary data are extended by zero to negative times.

Throughout this section, we assume \(n\geq2\) and write \(p(x',0,t):=\lim_{x_n\to0^+}p(x,t)\) whenever the limit exists.
Replacing \(p\) by \(p+c(t)\) adds \(c(t)\) to the boundary pressure in the next proposition and leaves its tangential derivatives unchanged.
\begin{proposition}\label{prop-0825-01}
		Let \(n\geq2\), \(1\leq i\leq n-1\), and let \(p\) be the pressure of the Stokes system \eqref{0308-01} with boundary data \eqref{0324-13}.
	Suppose that \(\psi\) satisfies {\rm\textbf{(S)}}, and that \(\phi\) satisfies {\rm\textbf{(T1)}}.
	If \(j<n\), the limit \(\partial_{x_i}p(x',0,t):=\lim_{x_n\downarrow0}\partial_{x_i}p(x,t)\) exists for every \((x',t)\in\Sigma\times(0,\infty)\) and is continuous in \((x',t)\).
	For \(x'\notin\supp\psi\), the following hold.
	\begin{enumerate}
		\item Suppose that \(j=n\).
		If \(n\geq3\), then
		\begin{align}\label{0822-11}
			p(x',0,t) = \frac{1}{2\sqrt{\pi}}\int_{\Sigma}\psi(y')\int_0^\infty \frac{2s\phi'(t)-\phi(t)-\phi(t-s)}{s^{\frac32}} \Gamma'(x'-y',s)\,dsdy'.
		\end{align}
		If \(n=2\), then
		\begin{align}\label{260136}
			p(x_1,0,t) & =-\frac{2\phi(t)}{\pi}\int_{\R}\frac{\psi(y_1)}{(x_1-y_1)^2}\,dy_1 -\frac{\phi'(t)}{\pi}\int_{\R}\psi(y_1)\ln|x_1-y_1|\,dy_1\nonumber\\
			& \quad-\frac{1}{2\sqrt{\pi}}\int_{\R}\psi(y_1)\int_0^\infty \frac{\phi(t-s)-\phi(t)}{s^{\frac32}} \Gamma_1(x_1-y_1,s)\,dsdy_1.
		\end{align}
		For every \(n\geq2\),
		\begin{align}\label{0822-12}
			\partial_{x_i}p(x',0,t) =\frac{1}{4\sqrt{\pi}}\int_{\Sigma}(x_i-y_i)\psi(y') \int_0^\infty \frac{\phi(t)+\phi(t-s)-2s\phi'(t)}{s^{\frac52}} \Gamma'(x'-y',s)\,dsdy'.
		\end{align}
		\item If \(j<n\), then
		\begin{align}\label{0822-10}
			p(x',0,t)= \frac{\Gamma\left(\frac{n}{2}\right)}{2^{n+1}\pi^{\frac{n}{2}}\Gamma\left(\frac{n+1}{2}\right)}\int_{\Sigma}(x_j-y_j)\psi(y')\int_0^\infty \frac{\phi(t)-\phi(t-s)}{s^{\frac{n+3}{2}}}{}_1F_1\left(\frac{n}{2},\frac{n+1}{2};-\frac{|x'-y'|^2}{4s}\right)\,dsdy',
		\end{align}
		and
		\begin{align}\label{0822-13}
			\partial_{x_i}p(x',0,t) & =-\frac{\Gamma(\frac{n+2}{2})} {2^{n+2}\pi^{\frac n2}\Gamma(\frac{n+3}{2})} \int_\Sigma(x_i-y_i)(x_j-y_j)\psi(y')\nonumber \\
			& \quad\times\int_0^\infty \frac{\phi(t)-\phi(t-s)}{s^{\frac{n+5}{2}}} {}_1F_1\left(\frac{n+2}{2},\frac{n+3}{2}; -\frac{|x'-y'|^2}{4s}\right)\,dsdy'\nonumber\\
			& \quad+\delta_{ij}\frac{\Gamma(\frac n2)} {2^{n+1}\pi^{\frac n2}\Gamma(\frac{n+1}{2})} \int_\Sigma\psi(y')\nonumber\\
			& \qquad\times\int_0^\infty \frac{\phi(t)-\phi(t-s)}{s^{\frac{n+3}{2}}} {}_1F_1\left(\frac n2,\frac{n+1}{2}; -\frac{|x'-y'|^2}{4s}\right)\,dsdy'.
		\end{align}
	\end{enumerate}
\end{proposition}
\begin{proof}
	\underline{Proof of (i).}\quad
	First assume \(n\geq3\).
	If \(g(x',t)=\psi(x')\phi(t)\mathbf e_n\), then \eqref{0418-01} gives
	\begin{align}\label{0825-01}
		p(x,t) & = -4 \phi(t)\partial_{x_n}^2\int_{\Sigma}N(x-y')\psi(y')\,dy' - 2 \phi'(t) \int_{\Sigma}N(x-y')\psi(y')\,dy' \nonumber \\
		 & \quad-\frac{1}{\sqrt{\pi}}\int_0^\infty \frac{\phi(t-s)-\phi(t)}{s^{\frac32}}\int_{\Sigma}\psi(z') \int_{\Sigma}\partial_{x_n}N(y',x_n) \Gamma'(x'-y'-z',s)\,dy'dz'ds.
	\end{align}
	Since \(x'\notin \supp \psi\), the first two terms on the right-hand side of \eqref{0825-01} converge to
	\begin{align}\label{0825-02}
		-\frac{2\Gamma\left(\frac{n}{2}\right)}{\pi^{\frac{n}{2}}}\phi(t)\int_{\Sigma}\frac{\psi(y')}{|x'-y'|^n}\,dy', \quad \text{and} \quad \frac{\Gamma\left(\frac{n-2}{2}\right)}{2\pi^{\frac{n}{2}}}\phi'(t) \int_{\Sigma}\frac{\psi(y')}{|x'-y'|^{n-2}}\,dy',
	\end{align}
	respectively as \(x_n\to0^+\).
	For \(f\in L^\infty(\Sigma)\) and \(s>0\), the identity \(\|2\partial_{y_n}N(\cdot,x_n)\|_{L^1(\Sigma)}=1\) gives
	\begin{align}\label{261013}
		 & \left|\frac{\phi(t-s)-\phi(t)}{2\sqrt{\pi}s^{\frac32}} \bigl(2\partial_{y_n}N(\cdot,x_n)*_\Sigma \Gamma'(\cdot,s)*_\Sigma f\bigr)(x')\right|\nonumber\\
		 & \quad\leq\frac{|\phi(t-s)-\phi(t)|}{2\sqrt{\pi}s^{\frac32}} \norm{2\partial_{y_n}N(\cdot,x_n)}_{L^1(\Sigma)} \norm{\Gamma'(\cdot,s)*_\Sigma f}_{L^\infty(\Sigma)}\nonumber\\
		 & \quad\leq\frac{\norm{f}_{L^\infty(\Sigma)}}{2\sqrt{\pi}s^{\frac32}} \min\left\{s\norm{\phi'}_{L^\infty(0,1)}, 2\norm{\phi}_{L^\infty(0,1)}\right\}\nonumber\\
		 & \quad\leq \frac{\max\left\{\norm{\phi'}_{L^\infty(0,1)}, 2\norm{\phi}_{L^\infty(0,1)}\right\}}{2\sqrt{\pi}}\, \norm{f}_{L^\infty(\Sigma)}\min\{s,1\}s^{-\frac32}.
	\end{align}
	The right-hand side of \eqref{261013} is integrable in \(s\), while \eqref{260267} gives
	\begin{align*}
		\lim_{x_n\to0^+}\int_{\Sigma}\partial_{x_n}N(y',x_n) \Gamma'(x'-y'-z',s)\,dy' =\frac12\Gamma'(x'-z',s).
	\end{align*}
	Therefore \eqref{261013}, with \(f=\psi\), allows us to apply the dominated convergence theorem to the \(s\)-integral in the third term of \eqref{0825-01}.
	This term converges as \(x_n\to0^+\) to
	\begin{align}\label{0825-03}
		-\frac{1}{2\sqrt{\pi}}\int_0^\infty \frac{\phi(t-s)-\phi(t)}{s^{\frac32}} \int_{\Sigma}\psi(y')\Gamma'(x'-y',s)\,dy'ds.
	\end{align}
	Substituting the identities
	\begin{align}\label{260138}
		\int_0^\infty \frac{\Gamma'(x'-y',s)}{s^{k+\frac12}}\,ds =\frac{2\cdot4^{k-1}\Gamma\left(\frac n2+k-1\right)} {\pi^{\frac{n-1}{2}}|x'-y'|^{n+2k-2}}, \qquad k\geq0,\quad \frac n2+k-1>0,
	\end{align}
	into \eqref{0825-01} and combining the result with \eqref{0825-02} and \eqref{0825-03} gives \eqref{0822-11}.

	Differentiating \eqref{0822-11}, which is justified since \(x'\notin \supp \psi\), gives
	\begin{align*}
		\partial_{x_i}p(x',0,t) = \frac{1}{4\sqrt{\pi}}\int_\Sigma(x_i-y_i)\psi(y') \int_0^\infty \frac{\phi(t)+\phi(t-s)-2s\phi'(t)}{s^{\frac52}} \Gamma'(x'-y',s)\,dsdy',
	\end{align*}
	which is \eqref{0822-12}.

	For \(n=2\), \eqref{0825-01} remains valid.
	Taking \(x_2\to0^+\), using \eqref{260267}, and differentiating \(N(x)=(2\pi)^{-1}\ln|x|\) gives \eqref{260136}.

	For \(n=2\), differentiating \eqref{260136} gives
	\begin{align}\label{260139}
		\partial_{x_1}p(x_1,0,t) & =\frac{4\phi(t)}{\pi}\int_\R\frac{\psi(y_1)}{(x_1-y_1)^3}\,dy_1 -\frac{\phi'(t)}{\pi}\int_\R\frac{\psi(y_1)}{x_1-y_1}\,dy_1\nonumber \\
		 & \quad+\frac{1}{4\sqrt{\pi}}\int_\R(x_1-y_1)\psi(y_1) \int_0^\infty\frac{\phi(t-s)-\phi(t)}{s^{\frac52}} \Gamma_1(x_1-y_1,s)\,dsdy_1.
	\end{align}
	Here we used \(\partial_{x_1}\Gamma_1(x_1-y_1,s) =-\frac{x_1-y_1}{2s}\Gamma_1(x_1-y_1,s)\).
	By \eqref{260138} with \(k=1,2\), we have
	\begin{align*}
		\frac{4\phi(t)}{\pi(x_1-y_1)^3} =\frac{\phi(t)(x_1-y_1)}{2\sqrt{\pi}}\int_0^\infty \frac{\Gamma_1(x_1-y_1,s)}{s^{\frac52}}\,ds, \qquad -\frac{\phi'(t)}{\pi(x_1-y_1)} =-\frac{\phi'(t)(x_1-y_1)}{2\sqrt{\pi}}\int_0^\infty \frac{\Gamma_1(x_1-y_1,s)}{s^{\frac32}}\,ds.
	\end{align*}
		Substituting these identities into \eqref{260139} and combining the three resulting integrals gives \eqref{0822-12} for \(n=2\).

	\underline{Proof of (ii).}\quad
	Take \(j=1\) for simplicity.
	If \(g(x',t)=\psi(x')\phi(t)\mathbf e_1\), then \eqref{0418-01} gives
	\begin{align}\label{0825-04}
		p(x,t) & = -4 \phi(t)\partial_{x_1}\partial_{x_n} \int_{\Sigma}N(x-y')\psi(y')\,dy'\nonumber \\
		 & \quad-\frac{1}{\sqrt{\pi}}\int_0^\infty \frac{\phi(t-s)-\phi(t)}{s^{\frac32}}\int_{\Sigma}\psi(z') \int_{\Sigma}\Gamma'(x'-y'-z',s) \partial_{y_1}N(y',x_n)\,dy'dz'ds.
	\end{align}
	The limit \eqref{260267} sends the first term on the right-hand side of \eqref{0825-04} to \(-2\phi(t)\partial_{x_1}\psi(x')=0\) as \(x_n\to0^+\), since \(x'\notin\supp\psi\).
	The second term satisfies
	\begin{align*}
		 & -\frac{1}{\sqrt{\pi}}\int_0^\infty \frac{\phi(t-s)-\phi(t)}{s^{\frac32}}\int_{\Sigma}\psi(z') \int_{\Sigma}\Gamma'(x'-y'-z',s) \partial_{y_1}N(y',x_n)\,dy'dz'ds\\
		 & \quad=-\frac{1}{\sqrt{\pi}}\int_0^\infty \frac{\phi(t-s)-\phi(t)}{s^{\frac32}}\int_{\Sigma}\int_{\Sigma} \psi(z')\Gamma'(x'-y',s) \partial_{y_1}N(y'-z',x_n)\,dz'dy'ds\\
		 & \quad=-\frac{1}{\sqrt{\pi}}\int_0^\infty \frac{\phi(t-s)-\phi(t)}{s^{\frac32}}\int_{\Sigma} \Gamma'(x'-y',s)\int_{\Sigma}\psi(z') \partial_{y_1}N(y'-z',x_n)\,dz'dy'ds.
	\end{align*}
	Lemma \ref{lem-260265} gives \(R_1'\psi\in C_0(\Sigma)\), and \eqref{261031} gives \( \partial_{y_1}N(\cdot,x_n)*_\Sigma\psi =\partial_{y_n}N(\cdot,x_n)*_\Sigma R_1'\psi \).
	Since \(\|2\partial_{y_n}N(\cdot,x_n)\|_{L^1(\Sigma)}=1\), the bound \eqref{261013} and the limit \eqref{260267}, applied with \(f=R_1'\psi\in C_0(\Sigma)\), allow us to apply the dominated convergence theorem to the remaining integrals.
	The last term in \eqref{0825-04} therefore converges to
	\begin{align}\label{0825-11}
		-\frac{1}{2\sqrt{\pi}}\int_0^\infty \frac{\phi(t-s)-\phi(t)}{s^{\frac32}}\int_{\Sigma}\Gamma'(x'-y',s)R_1'\psi(y')\,dy'ds.
	\end{align}
	We now claim that
	\begin{align}\label{0825-10}
		\int_{\Sigma}\Gamma'(x'-y',t)R_1'\psi(y')\,dy' = -\int_t^\infty \frac{1}{\sqrt{\pi(s-t)}}\partial_{x_1}\int_{\Sigma}\Gamma'(x'-y',s)\psi(y')\,dy'ds.
	\end{align}
	Indeed, the Fourier transform in \(x'\) of the left-hand side of \eqref{0825-10} is given by \(\frac{-i\xi_1 e^{-t|\xi'|^2}}{|\xi'|}\mathcal{F}_{x'}\psi(\xi')\), while the Fourier transform of the right-hand side of \eqref{0825-10} is given by
	\begin{align*}
		-\int_t^\infty \frac{1}{\sqrt{\pi(s-t)}}i\xi_1 e^{-s|\xi'|^2}\mathcal{F}_{x'}\psi(\xi')\,ds & = - \frac{i\xi_1 e^{-t|\xi'|^2}}{\sqrt{\pi}}\mathcal{F}_{x'}\psi(\xi')\int_0^\infty \frac{e^{-u|\xi'|^2}}{\sqrt{u}}\,du= - \frac{i\xi_1}{|\xi'|}e^{-t|\xi'|^2}\mathcal{F}_{x'}\psi(\xi').
	\end{align*}
			Hence \eqref{0825-10} holds.

	Substituting \eqref{0825-10}, with \(t=s\), into \eqref{0825-11} and differentiating \(\int_\Sigma\Gamma'(x'-y',\tau)\psi(y')\,dy'\) in \(x_1\), we transform \eqref{0825-11} into
	\begin{align}\label{0825-12}
		 & \frac{1}{2\pi}\int_0^\infty \frac{\phi(t-s)-\phi(t)}{s^{\frac32}}\int_s^\infty \frac{1}{\sqrt{\tau-s}}\partial_{x_1}\int_{\Sigma}\Gamma'(x'-y',\tau)\psi(y')\,dy'd\tau ds\nonumber \\
		 & = -\frac{1}{4\pi}\int_\Sigma (x_1-y_1)\psi(y')\int_0^\infty \frac{\phi(t-s)-\phi(t)}{s^{\frac32}}\int_s^\infty \frac{\Gamma'(x'-y',\tau)}{(\tau-s)^{\frac12}\tau}\,d\tau dsdy'.
	\end{align}
By \eqref{0324-14}, the inner integral is
	\begin{align}\label{1012-01}
		\int_s^\infty \frac{\Gamma'(x'-y',\tau)}{(\tau-s)^{\frac12}\tau}\,d\tau & = \frac{1}{(4\pi)^{\frac{n-1}{2}}}\int_s^\infty \frac{e^{-\frac{|x'-y'|^2}{4\tau}}}{\tau^{\frac{n+1}{2}}(\tau-s)^{\frac12}}\,d\tau\nonumber\\
		 & = \frac{\Gamma\left(\frac{n}{2}\right)}{2^{n-1}\pi^{\frac{n-2}{2}}\Gamma\left(\frac{n +1}{2}\right)}s^{-\frac{n}{2}}{}_1F_1\left(\frac{n}{2}, \frac{n+1}{2}; - \frac{|x'-y'|^2}{4s}\right).
	\end{align}
		Substituting \eqref{1012-01} into \eqref{0825-12} gives \eqref{0822-10}.

	Differentiating \eqref{0822-10}, which is justified by \eqref{0514-02}, and using \eqref{0327-07} give \eqref{0822-13}.

	For tangential boundary data, we now take the limit of \(\partial_{x_i}p(x,t)\) as \(x_n\downarrow0\) at every \(x'\in\Sigma\).
	Lemma \ref{lem-260265} gives \(R_j'(\partial_{x_i}\psi)\in C_0(\Sigma)\), while \(2\partial_{x_n}N(\cdot,x_n)\) is a nonnegative approximate identity of mass one.
	Differentiating \eqref{0825-04}, with \(1\) replaced by \(j\), moving tangential derivatives onto \(\psi\), and using \eqref{261031}, we obtain
	\begin{align*}
		\partial_{x_i}p(x,t) & =-4\phi(t)\bigl(\partial_{x_n}N(\cdot,x_n)*_\Sigma\partial_{x_i}\partial_{x_j}\psi\bigr)(x')\\
		&\quad-\frac1{\sqrt\pi}\int_0^\infty\frac{\phi(t-s)-\phi(t)}{s^{3/2}}\bigl(\partial_{x_n}N(\cdot,x_n)*_\Sigma\Gamma'(\cdot,s)*_\Sigma R_j'(\partial_{x_i}\psi)\bigr)(x')\,ds.
	\end{align*}
	The norm of the integrand in the second term is at most \(C\|R_j'(\partial_{x_i}\psi)\|_\infty\min\{s,1\}s^{-3/2}\), uniformly in \(x_n\), \(x'\), and \(t>0\), and this bound is integrable in \(s\).
	Moreover, associativity and commutativity of convolution give
	\begin{align*}
		 & \sup_{s>0}\bigl\|2\partial_{x_n}N(\cdot,x_n)*_\Sigma\Gamma'(\cdot,s)*_\Sigma R_j'(\partial_{x_i}\psi)-\Gamma'(\cdot,s)*_\Sigma R_j'(\partial_{x_i}\psi)\bigr\|_\infty\\
		 & \qquad\leq\sup_{s>0}\|\Gamma'(\cdot,s)\|_{L^1(\Sigma)}
		 \bigl\|2\partial_{x_n}N(\cdot,x_n)*_\Sigma R_j'(\partial_{x_i}\psi)-R_j'(\partial_{x_i}\psi)\bigr\|_\infty\\
		 & \qquad=\bigl\|2\partial_{x_n}N(\cdot,x_n)*_\Sigma R_j'(\partial_{x_i}\psi)-R_j'(\partial_{x_i}\psi)\bigr\|_\infty\longrightarrow0.
	\end{align*}
	Consequently, the limit as \(x_n\to0\) exists uniformly in \(x'\) and locally uniformly in \(t>0\), and equals
	\begin{align}\label{tangential-pressure-boundary-value}
		\partial_{x_i}p(x',0,t) & =-2\phi(t)\partial_{x_i}\partial_{x_j}\psi(x')\nonumber\\
		&\quad-\frac1{2\sqrt\pi}\int_0^\infty\frac{\phi(t-s)-\phi(t)}{s^{3/2}}\bigl(\Gamma'(\cdot,s)*_\Sigma R_j'(\partial_{x_i}\psi)\bigr)(x')\,ds.
	\end{align}
	Applying the same argument to \eqref{0825-04}, with \(1\) replaced by \(j\), gives the boundary value \(p(x',0,t)\). Lemma \ref{lem-260265} then identifies the \(x_i\)-derivative of \(p(\cdot,0,t)\) with \eqref{tangential-pressure-boundary-value}.
	The integrable bound \(C\|R_j'(\partial_{x_i}\psi)\|_\infty\min\{s,1\}s^{-3/2}\) and the continuity of \((x',s)\mapsto\bigl(\Gamma'(\cdot,s)*_\Sigma R_j'(\partial_{x_i}\psi)\bigr)(x')\) and of \(\phi\) show that \eqref{tangential-pressure-boundary-value} is jointly continuous in \((x',t)\in\Sigma\times(0,\infty)\).

	Since \(\psi\in C_c^{2,\alpha}(\Sigma)\) vanishes outside the closed unit ball, \(\partial_{x_i}\partial_{x_j}\psi=0\) on \(|x'|=1\).
	As \(x'\) approaches the unit sphere from \(|x'|>1\), \eqref{0822-13} tends to the boundary value in \eqref{tangential-pressure-boundary-value}.
\end{proof}

\begin{remark}\label{260142}
	For fixed \(t>0\) and \(j<n\), Lemma \ref{lem-260044}{\rm(vi)}, applied to \eqref{0822-10} and \eqref{0822-13}, shows that \(p(z'+r\omega',0,t)\) and \(\partial_{x_i}p(z'+r\omega',0,t)\) are real-analytic in \(r\) on every open interval where \(z'+r\omega'\notin\supp\psi\).
\end{remark}

Let \(\phi\) satisfy {\rm\textbf{(T1)}} and extend it by zero to negative times.
For \(r,t>0\), define
\begin{align}
	H_0(r,t)& :=\int_0^\infty \frac{\phi(t)-\phi(t-s)}{s^{\frac{n+3}{2}}}{}_1F_1\left(\frac n2,\frac{n+1}{2};-\frac{r^2}{4s}\right)\,ds,\label{0831-05}\\
	H_2(r,t)& :=\int_0^\infty \frac{\phi(t)-\phi(t-s)}{s^{\frac{n+5}{2}}}{}_1F_1\left(\frac{n+2}{2},\frac{n+3}{2};-\frac{r^2}{4s}\right)\,ds,\label{0831-06} \\
	\mathcal W(r,t) & :=\frac{1}{(4\pi)^{\frac n2}}\int_0^\infty \frac{\phi(t)+\phi(t-s)-2s\phi'(t)}{s^{\frac{n+4}{2}}}e^{-\frac{r^2}{4s}}\,ds.\label{0904-04}
\end{align}
Since \(\phi\) is bounded, \(\phi(t)-\phi(t-s)=O(s)\) as \(s\to0^+\) by {\rm\textbf{(T1)}}, and the factors involving the confluent hypergeometric functions satisfy \eqref{0514-02}, the integrals converge absolutely for \(r,t>0\).

We now establish the formula for \(\partial_{x_i}p(x',0,t)\) when \(\psi\) is radial.

\begin{proposition}\label{prop-260145}
		Let \(n\geq2\), \(1\leq i\leq n-1\), and let \(p\) be the pressure of the Stokes system \eqref{0308-01} with boundary data \eqref{0324-13}.
		Suppose that \(\psi\) satisfies {\rm\textbf{(S)}} and is radial, and that \(\phi\) satisfies {\rm\textbf{(T1)}}.
	Then, for \(r:=|x'|>1\), the following hold.
	\begin{enumerate}
		\item If \(j=n\), then
		\begin{align}\label{0831-01}
			\partial_{x_i}p(x',0,t) = \frac{x_i}{2r}\int_{\Sigma}(r-y_1)\psi(y')\mathcal W(|r\mathbf e_1'-y'|,t)\,dy'.
		\end{align}
		\item If \(n\geq3\) and \(j<n\), then
		\begin{align}\label{0825-13}
			\partial_{x_i}p(x',0,t) & =-\delta_{ij}\frac{\Gamma(\frac{n+2}{2})} {2^{n+2}\pi^{\frac n2}\Gamma(\frac{n+3}{2})} \int_\Sigma y_2^2\psi(y')H_2(|r\mathbf e_1'-y'|,t)\,dy'\nonumber \\
			& \quad+\delta_{ij}\frac{\Gamma(\frac n2)} {2^{n+1}\pi^{\frac n2}\Gamma(\frac{n+1}{2})} \int_\Sigma\psi(y')H_0(|r\mathbf e_1'-y'|,t)\,dy'\nonumber\\
			& \quad-\frac{\Gamma(\frac{n+2}{2})} {2^{n+2}\pi^{\frac n2}\Gamma(\frac{n+3}{2})} \frac{x_ix_j}{r^2}\int_\Sigma((r-y_1)^2-y_2^2)\psi(y') H_2(|r\mathbf e_1'-y'|,t)\,dy'.
		\end{align}
		\item If \(n=2\) and \(j=1\), then
		\begin{align}\label{260146}
			\partial_{x_1}p(x_1,0,t) & =-\frac{1}{12\pi^{3/2}}\int_{\R}(x_1-y_1)^2\psi(y_1)H_2(|x_1-y_1|,t)\,dy_1\nonumber \\
			& \quad+\frac{1}{4\pi^{3/2}}\int_{\R}\psi(y_1)H_0(|x_1-y_1|,t)\,dy_1.
		\end{align}
		Equivalently, for \(|x_1|>1\),
		\begin{align}\label{260155}
			\partial_{x_1}p(x_1,0,t) & =\frac{1}{6\pi^{\frac{3}{2}}}\int_{\R}\psi'(y_1)(x_1-y_1) \int_0^t\frac{\phi'(t-s)}{s^{\frac{3}{2}}} {}_1F_1\left(1,\frac52;-\frac{(x_1-y_1)^2}{4s}\right)\,dsdy_1.
		\end{align}
	\end{enumerate}
\end{proposition}
\begin{proof}
	\underline{Proof of (i).}\quad
	By \eqref{0822-12} and the definition \eqref{0904-04} of \(\mathcal W\), we have
	\begin{align}\label{0909-02}
		\partial_{x_i}p(x',0,t) = \frac12 \int_{\Sigma} (x_i-y_i)\psi(y') \mathcal W(|x'-y'|,t)\,dy'=: \frac12 V_i(x',t).
	\end{align}
	We now consider the vector field \(V(x',t)=(V_i(x',t))_{1\leq i\leq n-1}\).
Since \(\psi\) is radial, \(V(Qx',t)=QV(x',t)\) for every \(Q\in O(n-1)\).
	Hence \(V_i(x',t)=\frac{x_i}{r}V_1(r\mathbf e_1',t)\), where \(r=|x'|\), and
	\begin{align}\label{0909-01}
		V_1(r\mathbf e_1',t) = \int_{\Sigma}(r-y_1)\psi(y')\mathcal W(|r\mathbf e_1'-y'|,t)\,dy'.
	\end{align}
	Then \eqref{0909-02} and \eqref{0909-01} give \eqref{0831-01}.

	\underline{Proof of (ii).}\quad
	Assume \(n\geq3\).
	By \eqref{0822-13},
	\begin{align}\label{0831-04}
		\partial_{x_i}p(x',0,t) & = -\frac{\Gamma\left(\frac{n+2}{2}\right)}{2^{n+2}\pi^{\frac{n}{2}}\Gamma\left(\frac{n+3}{2}\right)} \int_{\Sigma}(x_i-y_i)(x_j-y_j)\psi(y')H_2(|x'-y'|,t)\,dy'\nonumber\\
		& \quad + \delta_{ij} \frac{\Gamma\left(\frac{n}{2}\right)}{2^{n+1}\pi^{\frac{n}{2}}\Gamma\left(\frac{n+1}{2}\right)}\int_{\Sigma}\psi(y')H_0(|x'-y'|,t)\,dy'\nonumber\\
		& = -\frac{\Gamma\left(\frac{n+2}{2}\right)}{2^{n+2}\pi^{\frac{n}{2}}\Gamma\left(\frac{n+3}{2}\right)}M_{ij}(x',t) + \delta_{ij} \frac{\Gamma\left(\frac{n}{2}\right)}{2^{n+1}\pi^{\frac{n}{2}}\Gamma\left(\frac{n+1}{2}\right)} S(x',t).
	\end{align}
Since \(\psi\) is radial, \(S\) is radial in \(x'\), and \(M(x',t):=(M_{ij}(x',t))_{1\leq i,j\leq n-1}\) satisfies \(M(Qx',t)=QM(x',t)Q^T\) for every \(Q\in O(n-1)\).
	Therefore
	\begin{align*}
		M_{ij}(x',t)=\delta_{ij}M_{22}(r\mathbf e_1',t) +\frac{x_ix_j}{r^2}\left[M_{11}(r\mathbf e_1',t) -M_{22}(r\mathbf e_1',t)\right].
	\end{align*}
	Then the definition of \(M_{ij}\) gives the following formula for \(1\leq i,j\leq n-1\).
	\begin{align}\label{0831-03}
		M_{ij}(x',t) & = \delta_{ij}\int_{\Sigma}y_2^2\psi(y') H_2(|r\mathbf e_1'-y'|,t)\,dy'+ \frac{x_ix_j}{r^2}\int_{\Sigma}((r-y_1)^2-y_2^2)\psi(y')H_2(|r\mathbf e_1'-y'|,t)\,dy'.
	\end{align}
	Then \eqref{0831-04} and \eqref{0831-03} give \eqref{0825-13}.

	\underline{Proof of (iii).}\quad
	The identity \eqref{260146} is \eqref{0822-13} with \(n=2\) and \(i=j=1\).

	To prove \eqref{260155}, substituting \eqref{0831-05}--\eqref{0831-06} with \(n=2\) into \eqref{260146} gives
	\begin{align}\label{0924-04}
		\partial_{x_1}p(x_1,0,t) & =\frac{1}{12\pi^{\frac32}}\int_{\R}\psi(y_1) \int_0^\infty\bigl(\phi(t)-\phi(t-s)\bigr) \biggl[ \frac{3}{s^{\frac52}} {}_1F_1\left(1,\frac32;-\frac{(x_1-y_1)^2}{4s}\right)\nonumber \\
		 & \hspace{45mm} -\frac{(x_1-y_1)^2}{s^{\frac72}} {}_1F_1\left(2,\frac52;-\frac{(x_1-y_1)^2}{4s}\right) \biggr]\,dsdy_1.
	\end{align}
	By \eqref{0911-21}, the expression in brackets equals
	\begin{align}\label{0928-20}
		\frac25\partial_s\left\{ s^{-\frac52}\left[(x_1-y_1)^2 {}_1F_1\left(2,\frac72;-\frac{(x_1-y_1)^2}{4s}\right) -5s\,{}_1F_1\left(1,\frac52;-\frac{(x_1-y_1)^2}{4s}\right)\right]\right\}.
	\end{align}
	Integration by parts in \(s\) is valid since the expression in braces in \eqref{0928-20} is \(O(s^{-\frac12})\) as \(s\to0^+\) and \(O(s^{-\frac32})\) as \(s\to\infty\), by \eqref{0514-02}, while \(\phi(t)-\phi(t-s)=O(s)\) as \(s\to0^+\).
	Since \(\phi'(t-s)=0\) for \(s>t\), it gives
	\begin{align}\label{260158}
		\partial_{x_1}p(x_1,0,t) & =-\frac{1}{30\pi^{\frac32}}\int_\R\psi(y_1)\int_0^t \frac{\phi'(t-s)}{s^{\frac52}}\biggl[ (x_1-y_1)^2{}_1F_1\left(2,\frac72;-\frac{(x_1-y_1)^2}{4s}\right)\nonumber \\
		 & \hspace{45mm} -5s\,{}_1F_1\left(1,\frac52;-\frac{(x_1-y_1)^2}{4s}\right) \biggr]\,dsdy_1\nonumber\\
		 & =-\frac{1}{6\pi^{\frac32}}\int_\R\psi(y_1)\int_0^t \frac{\phi'(t-s)}{s^{\frac32}} \partial_{y_1}\left[(x_1-y_1){}_1F_1\left(1,\frac52; -\frac{(x_1-y_1)^2}{4s}\right)\right]\,dsdy_1,
	\end{align}
	where the second equality follows from \eqref{0327-07}.
	Integration by parts in \(y_1\) in the last line of \eqref{260158} gives \eqref{260155}; the boundary terms vanish since \(\psi\) is compactly supported.
\end{proof}

\subsection{Sign of \texorpdfstring{\(\partial_{x_i}p(x',0,t)\)}{pressure derivatives} for normal boundary data with radial \texorpdfstring{\(\psi\)}{psi}}
The following proposition proves Theorem \ref{thm-260019}.
\begin{proposition}\label{prop-260147}
	Let \(n\geq2\), \(1\leq i\leq n-1\), and let \(p\) be the pressure of the Stokes system \eqref{0308-01} with boundary data \eqref{0324-13}, where \(j=n\).
		Suppose that \(\psi\) satisfies {\rm\textbf{(S)}} and is radial, and that \(\phi\) satisfies {\rm\textbf{(T1)}} and {\rm\textbf{(T2)}}.
	Then the following hold.
	\begin{enumerate}
		\item For \(0<t<t_1\), there exists \(r_{\rm nor}(t)>0\) such that
		\begin{align*}
			\sgn(\partial_{x_i}p(x',0,t)) =
			\begin{cases}
				\sgn(x_i),& 1<|x'|\leq r_{\rm nor}(t)-1, \\
				-\sgn(x_i), & |x'|\geq r_{\rm nor}(t)+1,
			\end{cases}
		\end{align*}
		where the first region is empty if \(r_{\rm nor}(t)\leq2\).
		Moreover, \(r_{\rm nor}\) is continuous on \((0,t_1)\), and
		\begin{align}\label{261038}
			n<\frac{\phi'(t)r_{\rm nor}(t)^2}{\phi(t)}<2n \qquad(0<t<t_1),
		\end{align}
		and
		\begin{align}\label{261045}
			r_{\rm nor}(t) =\left(\frac{n\phi(t)}{\phi'(t)}\right)^{\frac12}(1+o(1)) \qquad\text{as }t\to t_1^-.
		\end{align}
		If \(\phi\) is log-concave on \((0,1)\), then \(r_{\rm nor}(t)\to0\) as \(t\to0^+\).
		\item For \(t<t_1\) sufficiently close to \(t_1\), there exists a unique \(\rho_{\rm nor}(t)>1\) such that
		\begin{align*}
			\sgn(\partial_{x_i}p(x',0,t))=-\sgn(x_i)\sgn(|x'|-\rho_{\rm nor}(t)),\qquad |x'|>1.
		\end{align*}
		Moreover, \(|\rho_{\rm nor}(t)-r_{\rm nor}(t)|<1\), and as \(t\to t_1^-\),
		\begin{align}\label{261039}
			\rho_{\rm nor}(t)^2-r_{\rm nor}(t)^2
			=\frac{n+3}{n-1}\frac{M_2}{M_0}+O\bigl(r_{\rm nor}(t)^{-2}\bigr).
		\end{align}
		\item If \(|x'|>1\) and \(t\geq t_1\), then \(\sgn(\partial_{x_i}p(x',0,t))=\sgn(x_i)\).
	\end{enumerate}
\end{proposition}

\begin{lemma}\label{lem-260033}
	Let \(\phi\) satisfy {\rm\textbf{(T1)}} and {\rm\textbf{(T2)}}, and let \(\mathcal W\) be given in \eqref{0904-04}.
Then the following hold.
	\begin{enumerate}
		\item\label{260034} If \(t\geq t_1\), then \(\mathcal W(r,t)>0\) for all \(r>0\).
		\item\label{260035} For each \(t\in(0,t_1)\), \(\mathcal W(\cdot,t)\) has a unique zero \(r_\#(t)>0\), and
		\begin{align*}
			\mathcal W(r,t)
			\begin{cases}
				>0, & 0<r<r_\#(t), \\
				<0, & r>r_\#(t).
			\end{cases}
		\end{align*}
		\item\label{260036} \(r_\#\) is continuous on \((0,t_1)\), and
		\begin{align}\label{260037}
			n<\frac{\phi'(t)r_\#(t)^2}{\phi(t)}<2n \qquad (0<t<t_1).
		\end{align}
		If \(\phi\) is log-concave on \((0,1)\), then \(r_\#(t)\to0\) as \(t\to0^+\).
		As \(t\to t_1^-\), \(r_\#(t)\to\infty\) and
		\begin{align}\label{261007}
			r_\#(t) =\left(\frac{n\phi(t)}{\phi'(t)}\right)^{\frac12}(1+o(1)) \qquad\text{as }t\to t_1^-.
		\end{align}
	\end{enumerate}
\end{lemma}
\begin{proof}
	\underline{Proof of (i).}\quad
	If \(t_1\leq t<1\), then \(\phi(t)>0\) and \(\phi'(t)\leq0\); if \(t\geq1\), then \(\phi(t)=\phi'(t)=0\) and \(\phi(t-s)>0\) for \(s\in(t-1,t)\).
In both cases \(\phi(t)+\phi(t-s)-2s\phi'(t)\geq0\) is not identically zero in \(s\), so \(\mathcal W(r,t)>0\) for all \(r>0\).

	\underline{Proof of (ii).}\quad
	Set
	\begin{align}\label{0908-03}
		B(s,t) := \phi(t)+\phi(t-s)-2s\phi'(t), \qquad s>0,
	\end{align}
	so that
	\begin{align}\label{0920-01}
		\mathcal{W}(r,t)= \frac{1}{(4\pi)^{\frac n2}} \int_0^\infty B(s,t) s^{-\frac{n+4}{2}} e^{-\frac{r^2}{4s}}\,ds.
	\end{align}
	By {\rm\textbf{(T2)}}, \(B'(s,t)=-\phi'(t-s)-2\phi'(t)<0\) for \(0<s<t\), and \(B'(s,t)=-2\phi'(t)<0\) for \(s>t\).
Since \(\phi(0)=0\), \(B(\cdot,t)\) is continuous, hence strictly decreasing on \((0,\infty)\).
	Since
	\begin{align*}
		\lim_{s\to0^+}B(s,t)=2\phi(t)>0, \qquad \text{and } \qquad \lim_{s\to\infty}B(s,t)=-\infty,
	\end{align*}
	there exists a unique \(\sigma_\#(t)>0\) such that
	\begin{align}\label{260223}
		B(s,t)
		\begin{cases}
			>0, & 0<s<\sigma_\#(t), \\
			<0, & s>\sigma_\#(t).
		\end{cases}
	\end{align}

	Fix \(0<r_1<r_2\).
	Then \(e^{-\frac{r_2^2-r_1^2}{4s}}\) is strictly increasing in \(s\) on \((0,\infty)\).
	Hence, using \eqref{0920-01} and \eqref{260223},
	\begin{align}
		 & (4\pi)^{\frac n2} \left( \mathcal W(r_2,t) - e^{-\frac{r_2^2-r_1^2}{4\sigma_\#(t)}}\mathcal W(r_1,t) \right)= \int_0^\infty B(s,t)s^{-\frac{n+4}{2}} e^{-\frac{r_1^2}{4s}} \left( e^{-\frac{r_2^2-r_1^2}{4s}}- e^{-\frac{r_2^2-r_1^2}{4\sigma_\#(t)}} \right)\,ds <0. \label{260224}
	\end{align}
	The last inequality holds because \(B(s,t)\) and the last bracket have opposite signs for \(s\neq\sigma_\#(t)\).
Thus \(\mathcal W(r_1,t)\leq0\) implies \(\mathcal W(r_2,t)<0\), so \(\mathcal W(\cdot,t)\) has at most one zero on \((0,\infty)\).
Substituting \(u=r^2/(4s)\) gives
	\begin{align*}
		r^{n+2}\mathcal W(r,t) = \frac{4}{\pi^{\frac{n}{2}}} \int_0^\infty \left[ \phi(t) + \phi\left(t-\frac{r^2}{4u}\right) - \frac{r^2}{2u}\phi'(t) \right] u^{\frac n2}e^{-u}\,du.
	\end{align*}
	By the dominated convergence theorem,
	\begin{align}
		\lim_{r\to0^+}r^{n+2}\mathcal W(r,t) & = \frac{8}{\pi^{\frac{n}{2}}} \phi(t) \int_0^\infty u^{\frac n2}e^{-u}\,du = \frac{8\Gamma\left(\frac{n+2}{2}\right)} {\pi^{\frac n2}}\phi(t)>0. \label{260225}
	\end{align}

	On the other hand, for large \(r\), we write
	\begin{align}\label{261011}
		(4\pi)^{\frac n2}\mathcal W(r,t) & = \phi(t)\int_0^\infty s^{-\frac{n+4}{2}}e^{-\frac{r^2}{4s}}\,ds - 2\phi'(t)\int_0^\infty s^{-\frac{n+2}{2}}e^{-\frac{r^2}{4s}}\,ds + \int_0^t \phi(t-s)s^{-\frac{n+4}{2}}e^{-\frac{r^2}{4s}}\,ds\\
		 & = 2^{n+2}\Gamma\left(\frac{n+2}{2}\right)r^{-n-2}\phi(t) - 2^{n+1}\Gamma\left(\frac{n}{2}\right)r^{-n}\phi'(t) + \int_0^t \phi(t-s)s^{-\frac{n+4}{2}}e^{-\frac{r^2}{4s}}\,ds.
	\end{align}
	Moreover,
	\begin{align*}
		0 & \leq r^n\int_0^t \phi(t-s)s^{-\frac{n+4}{2}}e^{-\frac{r^2}{4s}}\,ds \leq 2^{n+2}\norm{\phi}_{L^\infty(0,1)} r^{-2}\int_{\frac{r^2}{4t}}^\infty u^{\frac n2}e^{-u}\,du \longrightarrow0 \qquad\text{as }r\to\infty.
	\end{align*}
	Therefore
	\begin{align}
		\lim_{r\to\infty}r^n\mathcal W(r,t) = -\frac{2\Gamma\left(\frac n2\right)} {\pi^{\frac n2}}\phi'(t)<0, \qquad 0<t<t_1. \label{260226}
	\end{align}
	It follows from \eqref{260225}, \eqref{260226}, and the continuity of \(\mathcal W\) in \(r\) that \(\mathcal W(r,t)\) has a zero for every \(t\in(0,t_1)\).
	Its uniqueness follows from \eqref{260224}, which gives \eqref{260035}.

	\underline{Proof of (iii).}\quad
	By the dominated convergence theorem, \(\mathcal W\) is jointly continuous on \((0,\infty)\times(0,\infty)\).
Fix \(t_0\in(0,t_1)\) and \(0<\eta\leq r_\#(t_0)/2\).
By \eqref{260035} and the joint continuity, \(\mathcal W(r_\#(t_0)-\eta,t)>0>\mathcal W(r_\#(t_0)+\eta,t)\) for \(t\) near \(t_0\), so \(|r_\#(t)-r_\#(t_0)|<\eta\).
Hence \(r_\#\) is continuous on \((0,t_1)\).

	Log-concavity of \(\phi\) on \((0,1)\) implies
	\begin{align}\label{260227}
		\lim_{t\to0^+}\frac{\phi'(t)}{\phi(t)}=\infty.
	\end{align}

	To prove \eqref{260037} under {\rm\textbf{(T1)}} and {\rm\textbf{(T2)}} alone, evaluate \eqref{261011} at \(r=r_\#(t)\) and substitute \(v=r^2/(4s)\) in the last integral.
	This gives
	\begin{align}\label{1012-02}
		\frac{\phi'(t)r^2}{\phi(t)}-n
		=\frac{2}{\Gamma(n/2)}\int_0^\infty
		\frac{\phi(t-r^2/(4v))}{\phi(t)}\,v^{n/2}e^{-v}\,dv .
	\end{align}
	The ratio in the integrand is zero for \(v\leq r^2/(4t)\) and lies strictly between zero and one for \(v>r^2/(4t)\).
	Since \(2\Gamma(1+n/2)/\Gamma(n/2)=n\), the right-hand side lies strictly between zero and \(n\).
	This proves \eqref{260037}, and hence
	\begin{align}\label{260228}
		r_\#(t)\longrightarrow0\quad\Longleftrightarrow\quad
		\frac{\phi(t)}{\phi'(t)}\longrightarrow0 \qquad(t\to0^+).
	\end{align}
	If \(\phi\) is log-concave, then \(\phi(t)/\phi'(t)\to0\) as \(t\to0^+\) by \eqref{260227}.

	At the other endpoint, \(\phi(t)\to\phi(t_1)>0\) and \(\phi'(t)\to0\).
	The lower bound in \eqref{260037} thus gives \(r_\#(t)\to\infty\) as \(t\to t_1^-\).

Since the ratio in \eqref{1012-02} lies in \([0,1]\) and vanishes for \(v\leq r^2/(4t)\), the right-hand side of \eqref{1012-02} at \(r=r_\#(t)\) is \(O_n(e^{-r_\#(t)^2/8})\) as \(t\to t_1^-\).
This proves \eqref{261007}.
\end{proof}

The initial limit in Lemma \ref{lem-260033}{\rm(iii)} can fail under {\rm\textbf{(T1)}} and {\rm\textbf{(T2)}}.
For \(0<u\leq1/2\), take \(v(u)=(1-2u)u^2(1+\sin(1/u)+u^2)\), with \(v(0)=0\), and define the temporal boundary data by
\begin{align*}
	\phi(t)=\int_0^{\min\{t,1-t\}}v(u)\,du\quad(0\leq t\leq1),\qquad \phi(t)=0\quad(t>1).
\end{align*}
Here \(v>0\) on \((0,1/2)\), \(v(0)=v(1/2)=0\), and \(v'\) is bounded.
Thus \(\phi'\) is Lipschitz and {\rm\textbf{(T2)}} holds with \(t_1=1/2\).

As \(t\to0^+\), \(\phi(t)=t^3/3+O(t^4)\), since integration by parts using \(\sin(1/u)=u^2\frac{d}{du}\cos(1/u)\) gives \(\int_0^t u^2\sin(1/u)\,du=O(t^4)\).
Along sequences with \(\sin(1/t)=1\) and \(\sin(1/t)=-1\), respectively, \(\phi(t)/\phi'(t)\) is asymptotic to \(t/6\) and \(1/(3t)\).
By \eqref{260037}, \(r_\#(t)\) tends to zero along the first sequence and to infinity along the second.

\begin{proof}[Proof of Proposition \ref{prop-260147}]
	Let \(r=|x'|\).

	\underline{Proof of (i).}\quad
	For \(0<t<t_1\), set \(r_{\rm nor}(t):=r_\#(t)\).
If \(\psi(y')>0\), then \(|y'|<1\), and hence \(r-1<|r\mathbf e_1'-y'|<r+1\).
	Hence \(|r\mathbf e_1'-y'|<r_{\rm nor}(t)\) for \(1<r\leq r_{\rm nor}(t)-1\), and \(|r\mathbf e_1'-y'|>r_{\rm nor}(t)\) for \(r\geq r_{\rm nor}(t)+1\).
	Since \(r-y_1>0\), Lemma \ref{lem-260033}{\rm(ii)} and \eqref{0831-01} give the signs in {\rm(i)}.
	The continuity, \eqref{261038}, and the limit as \(t\to0^+\) follow from Lemma \ref{lem-260033}{\rm(iii)}, and \eqref{261007} gives \eqref{261045}.

	\underline{Proof of (ii).}\quad
	By rotational symmetry, it suffices to take \(i=1\) and \(x'=r\mathbf e_1'\) with \(r>1\); when \(n=2\), \(i=1\) is the only tangential index.
	The other components follow from \eqref{0831-01}.
	Part {\rm(i)} and Lemma \ref{lem-260033}{\rm(iii)} give
	\begin{align}\label{261090}
		r_{\rm nor}(t)\longrightarrow\infty,\qquad \frac{\phi'(t)r_{\rm nor}(t)^2}{\phi(t)}\longrightarrow n \qquad\text{as }t\to t_1^-.
	\end{align}
	For \(t<t_1\) sufficiently close to \(t_1\), we have \(r_{\rm nor}(t)>3\).
	Part {\rm(i)} then gives opposite signs at \(r=r_{\rm nor}(t)-1\) and \(r=r_{\rm nor}(t)+1\).
	Since \(r\mapsto\partial_{x_1}p(r\mathbf e_1',0,t)\) is continuous on \((1,\infty)\) by \eqref{0831-01} and the joint continuity of \(\mathcal W\) established in Lemma \ref{lem-260033}{\rm(iii)}, the intermediate value theorem gives a zero in \((1,\infty)\).
	Every zero in \((1,\infty)\) belongs to
	\begin{align}\label{261091}
		(r_{\rm nor}(t)-1,r_{\rm nor}(t)+1).
	\end{align}
	We prove the uniqueness by estimating the radial derivative uniformly throughout this interval.

	By \eqref{261011}, for every \(r>0\),
	\begin{align}\label{261092}
		\mathcal W(r,t) & =\frac{2\Gamma(\frac n2)}{\pi^{\frac n2}}\left(n\phi(t)r^{-n-2}-\phi'(t)r^{-n}\right) \nonumber\\
		& \quad+\frac1{(4\pi)^{\frac n2}}\int_0^t\phi(t-s)s^{-\frac{n+4}{2}}e^{-\frac{r^2}{4s}}\,ds.
	\end{align}

	For \(m\in\{n,n+2\}\), \(0\leq q\leq\frac12\), and \(|y'|\leq1\), set
	\begin{align*}
		F_m(q,y'):=\frac{1-qy_1}{|\mathbf e_1'-qy'|^m},
		\qquad J_m(q):=\int_\Sigma\psi(y')F_m(q,y')\,dy'.
	\end{align*}
Since \(|\mathbf e_1'-qy'|\geq\frac12\), direct differentiation gives
	\begin{align*}
		F_m(0,y')=1,\qquad \partial_qF_m(0,y')=(m-1)y_1,\qquad
		|\partial_q^2F_m(q,y')|\leq C(n)|y'|^2.
	\end{align*}
	Radiality gives
	\begin{align*}
		\int_\Sigma y_\ell\psi(y')\,dy'=0,\qquad 1\leq\ell\leq n-1.
	\end{align*}
	Taylor's theorem and the fundamental theorem of calculus therefore give
	\begin{align}\label{261092a}
		|J_m(q)-M_0|\leq C(n)M_2q^2,
		\qquad |J_m'(q)|\leq C(n)M_2q.
	\end{align}
	For \(r\geq2\), taking \(q=r^{-1}\) gives
	\begin{align*}
		\frac1r\int_\Sigma\frac{(r-y_1)\psi(y')}{|r\mathbf e_1'-y'|^m}\,dy'
		=r^{-m}J_m(r^{-1}).
	\end{align*}
	Since \(M_2\leq M_0\), differentiating in \(r\) and using \eqref{261092a} with \(m=n\) and \(m=n+2\), respectively, give, for \(k=0,1\),
	\begin{align}\label{261093}
		\left|\partial_r^k\left[\frac1r\int_\Sigma\frac{(r-y_1)\psi(y')}{|r\mathbf e_1'-y'|^n}\,dy'-M_0r^{-n}\right]\right| & \leq C(n)M_0r^{-n-2-k},\nonumber\\
		\left|\partial_r^k\left[\frac1r\int_\Sigma\frac{(r-y_1)\psi(y')}{|r\mathbf e_1'-y'|^{n+2}}\,dy'-M_0r^{-n-2}\right]\right| & \leq C(n)M_0r^{-n-4-k}.
	\end{align}

	To estimate the term containing the exponential, consider the spatial integral
	\begin{align*}
		\mathcal G_0(r,s)=\int_\Sigma(r-y_1)\psi(y')e^{-\frac{|r\mathbf e_1'-y'|^2}{4s}}\,dy'.
	\end{align*}
	Differentiation under the integral gives
	\begin{align*}
		\partial_r\mathcal G_0(r,s)=\int_\Sigma\psi(y')\left(1-\frac{(r-y_1)^2}{2s}\right)e^{-\frac{|r\mathbf e_1'-y'|^2}{4s}}\,dy'.
	\end{align*}
	If \(r\geq2\), \(y'\in\supp\psi\), and \(0<s\leq1\), then
	\begin{align*}
		|r\mathbf e_1'-y'|\geq\frac r2,\qquad |r-y_1|\leq\frac{3r}{2},\qquad e^{-\frac{|r\mathbf e_1'-y'|^2}{4s}}\leq e^{-\frac{r^2}{32}}e^{-\frac{r^2}{32s}}.
	\end{align*}
	Consequently,
	\begin{align}\label{261093b}
		|\mathcal G_0(r,s)|&\leq \frac32M_0r e^{-\frac{r^2}{32}}e^{-\frac{r^2}{32s}},\nonumber\\
		|\partial_r\mathcal G_0(r,s)|&\leq M_0\left(1+\frac{9r^2}{8s}\right)e^{-\frac{r^2}{32}}e^{-\frac{r^2}{32s}}.
	\end{align}
	For every \(a>1\), the change of variables \(u=r^2/(32s)\) gives
	\begin{align}\label{261093c}
		\int_0^1s^{-a}e^{-\frac{r^2}{32s}}\,ds
		=\left(\frac{r^2}{32}\right)^{1-a}\int_{\frac{r^2}{32}}^\infty u^{a-2}e^{-u}\,du
		\leq C(a)r^{2-2a}.
	\end{align}

	Using \eqref{261093c} with \(a=(n+4)/2\) and \(a=(n+6)/2\), and differentiating the factor \(r^{n+1}\) when \(k=1\), gives
	\begin{align}\label{261094}
		\left|\partial_r^k\left[\frac{r^{n+1}}{M_0}\int_0^t\phi(t-s)s^{-\frac{n+4}{2}}\mathcal G_0(r,s)\,ds\right]\right|\leq C(n)\|\phi\|_{L^\infty(0,1)}e^{-cr^2},\qquad k=0,1,
	\end{align}
	uniformly for \(r\geq2\) and \(0<t\leq1\), where \(c>0\) is absolute.

	Combining \eqref{0831-01}, \eqref{261092}, \eqref{261093}, and \eqref{261094}, we obtain, uniformly for \(r\geq2\) and \(0<t<t_1\),
	\begin{align}\label{261095}
		\frac{\pi^{\frac n2}r^{n+1}}{\Gamma(\frac n2)M_0}\partial_{x_1}p(r\mathbf e_1',0,t)=n\phi(t)-r^2\phi'(t)+O\!\left(\frac{\phi(t)}{r^2}+\phi'(t)\right)+O\!\left(\|\phi\|_{L^\infty(0,1)}e^{-cr^2}\right),
	\end{align}
	and, using the cases \(k=1\) of \eqref{261093} and \eqref{261094} separately,
	\begin{align}\label{261096}
		\partial_r\left[\frac{\pi^{\frac n2}r^{n+1}}{\Gamma(\frac n2)M_0}\partial_{x_1}p(r\mathbf e_1',0,t)\right]=-2r\phi'(t)+O\!\left(\frac{\phi(t)}{r^3}+\frac{\phi'(t)}r\right)+O\!\left(\|\phi\|_{L^\infty(0,1)}e^{-cr^2}\right).
	\end{align}
	All implicit constants in these estimates depend only on \(n\).

	By \eqref{261090}, uniformly for \(r\) in \eqref{261091},
	\begin{align*}
		r\longrightarrow\infty,\qquad \frac{r}{r_{\rm nor}(t)}\longrightarrow1,\qquad \frac{\phi'(t)r^2}{\phi(t)}\longrightarrow n.
	\end{align*}
	Moreover, \(\phi(t)\to\phi(t_1)>0\).
	Multiplying \eqref{261096} by \(r/\phi(t)\) therefore gives, uniformly throughout \eqref{261091},
	\begin{align}\label{261097}
		\frac r{\phi(t)}\partial_r\left[\frac{\pi^{\frac n2}r^{n+1}}{\Gamma(\frac n2)M_0}\partial_{x_1}p(r\mathbf e_1',0,t)\right]\longrightarrow-2n<0 \qquad\text{as }t\to t_1^-.
	\end{align}
	For the error terms, the same normalization gives \(O(r^{-2}+\phi'(t)/\phi(t))\) and \(O(r e^{-cr^2}/\phi(t))\), both tending to zero uniformly on this interval.

	Consequently, for \(t<t_1\) sufficiently close to \(t_1\), the left-hand side of \eqref{261095} is strictly decreasing on \eqref{261091}.
	It has at most one zero there.
	The existence and localization already proved show that it has exactly one zero in \((1,\infty)\), denoted by \(\rho_{\rm nor}(t)\).
	The signs outside the interval \eqref{261091} and the strict decrease inside it give
	\begin{align*}
		\sgn\bigl(\partial_{x_1}p(r\mathbf e_1',0,t)\bigr)=-\sgn(r-\rho_{\rm nor}(t)),\qquad r>1.
	\end{align*}
	By \eqref{0831-01}, this yields
	\begin{align*}
		\sgn\bigl(\partial_{x_i}p(x',0,t)\bigr)=-\sgn(x_i)\sgn(|x'|-\rho_{\rm nor}(t)),\qquad |x'|>1.
	\end{align*}
	Finally, \eqref{261091} implies \(|\rho_{\rm nor}(t)-r_{\rm nor}(t)|<1\).
	Together with \eqref{261090}, this gives \(\rho_{\rm nor}(t)\sim r_{\rm nor}(t)\).

	To compare their squares, write \(R^2=n\phi(t)/\phi'(t)\), \(r=\rho_{\rm nor}(t)\), and \(q=1/r\).
	Radiality makes \(J_m\) even in \(q\) and gives \(\int_\Sigma y_1^2\psi=M_2/(n-1)\).
	For \(|q|\leq1/2\), the denominator of \(F_m\) stays bounded away from zero, so its fourth \(q\)-derivative is bounded uniformly on \(|y'|\leq1\).
	The coefficient of \(q^2\) in \(F_m(q,y')\) is \((m^2y_1^2-m|y'|^2)/2\).
	Consequently, for \(m=n,n+2\),
	\begin{align*}
		J_m(q)=M_0+\frac{m(m-n+1)}{2(n-1)}M_2q^2+O_n(M_0q^4).
	\end{align*}

	Since \(\partial_{x_i}p\) vanishes at \(|x'|=r\), equations \eqref{0831-01}, \eqref{261092}, and \eqref{261094} give
	\begin{align*}
		r^2\phi'(t)=n\phi(t)\frac{J_{n+2}(1/r)}{J_n(1/r)}
		+O_n\bigl(\|\phi\|_{L^\infty(0,1)}e^{-cr^2}\bigr),
	\end{align*}
	since \(J_n(1/r)\to M_0>0\).
	As \(r\sim R\) and \(\phi(t)\to\phi(t_1)>0\), dividing by \(\phi'(t)\) and expanding the quotient gives
	\begin{align*}
		r^2-R^2=\frac{n+3}{n-1}\frac{M_2}{M_0}\frac{R^2}{r^2}+O(R^{-2}).
	\end{align*}
	In particular, \(r^2-R^2=O(1)\), so \(R^2/r^2=1+O(R^{-2})\).
	At \(r=r_{\rm nor}(t)\), \eqref{1012-02} and the proof of Lemma \ref{lem-260033} give \(r_{\rm nor}(t)^2-R^2=O(R^{-2})\).
	Subtracting proves \eqref{261039}, whose leading coefficient is positive by {\rm\textbf{(S)}}.

	\underline{Proof of (iii).}\quad
Using \(r-y_1>0\) for \(y'\in\supp\psi\), Lemma \ref{lem-260033}{\rm(i)} and \eqref{0831-01} give \(\sgn(\partial_{x_i}p(x',0,t))=\sgn(x_i)\).
\end{proof}

\subsection{Sign of \texorpdfstring{\(\partial_{x_i}p(x',0,t)\)}{pressure derivatives} for tangential boundary data with radial \texorpdfstring{\(\psi\)}{psi}}
Throughout this subsection, \(1\leq j\leq n-1\) in \eqref{0324-13}, and we first consider \(i=j\) and then \(i\neq j,n\).

\subsubsection{The \(i=j\) case: preliminaries for all \(n\geq2\)}\label{subsubsect-261080-tools}

By rotation, we take \(i=j=1\) in Subsections \ref{subsubsect-261080-tools}--\ref{subsubsect-261022-offaxis}.
We first define \(F_t\) and establish Lemmas \ref{lem-0911-02b} and \ref{lem-261080}, which apply for every \(n\geq2\).
The case \(n=2\) is treated in Subsection \ref{subsubsect-260154-n2}.
For \(n\geq3\), the signs and zero curves of \(\partial_{x_1}p\) are studied in Subsections \ref{subsubsect-261020-axis} and \ref{subsubsect-261022-offaxis}, respectively.

For \(\phi\) satisfying {\rm\textbf{(T1)}} and \(t,z>0\), set
\begin{align}\label{260038}
	F_t(z):=-\int_0^t\frac{\phi'(t-s)}{s^{\frac{n+1}{2}}} {}_1F_1\left(\frac n2,\frac{n+3}{2};-\frac{z}{4s}\right)\,ds.
\end{align}
The next lemma gives the properties of \(F_t\) used below.
\begin{lemma}\label{lem-0911-02b}
Let \(n\geq2\), and suppose that \(\phi\) satisfies {\rm\textbf{(T1)}}--{\rm\textbf{(T4)}}.
	Then the following hold.
	\begin{enumerate}
		\item \(F_t(z)>0\) for every \(t\geq t_\infty\) and \(z>0\).
		\item If \(\phi\) also satisfies {\rm\textbf{(T5)}}, then
		\begin{align}\label{0911-35}
			-zF_{t_\infty}'(z)\geq\frac{n-1}{2}F_{t_\infty}(z), \qquad z>0.
		\end{align}
	\end{enumerate}
\end{lemma}
\begin{proof}
	\underline{Proof of (i).}\quad
	For fixed \(z>0\), the kernel
	\begin{align*}
	    K(z,t)=t^{-(n+1)/2}{}_1F_1\left(\frac n2,\frac{n+3}{2};-\frac{z}{4t}\right)
	\end{align*}
	is positive for \(t>0\). Identity \eqref{0911-21} with \((a,b,c)=(n/2,(n+3)/2,z/4)\) gives \(\partial_t K(z,t)<0\).
	The expansion \eqref{0327-10} and the same derivative identity show that \(K(z,\cdot)\) satisfies \eqref{0829-01} with \(\sigma=1/2\).
		Thus, by Lemma \ref{lem-0325-01}, \(t\mapsto F_t(z)\) has at most one zero in \((0,1)\) and \(\partial_tF_t(z)>0\) at every zero.

	We have \(F_t(z)<0\) for \(0<t\leq t_1\) since \(\phi'(t-s)>0\) for all \(0<s<t\leq t_1\).
	We now determine the value of \(F_{t_\infty}(z)\).
	Let \(s_0=t_\infty-t_1\).
	Since \( s^{-\frac n2} {}_1F_1\left(\frac n2,\frac{n+3}{2};- \frac{z}{4s}\right) = \frac{2^{n+1}\Gamma\left(\frac{n+3}{2}\right)}{\sqrt{\pi}}z^{-\frac n2}Q(z/(4s))\), Lemma \ref{lem-0424-01}{\rm(i)} shows that it is strictly decreasing in \(s>0\).
	Since \(T_0(t_\infty)=0\), subtracting its value at \(s=s_0\) gives
	\begin{align*}
		F_{t_\infty}(z) & =-\int_0^{t_\infty}\frac{\phi'(t_\infty-s)}{\sqrt{s}} \left[s^{-\frac n2}{}_1F_1\left( \frac n2,\frac{n+3}{2};-\frac{z}{4s}\right) -s_0^{-\frac n2}{}_1F_1\left( \frac n2,\frac{n+3}{2};-\frac{z}{4s_0}\right)\right]\,ds>0.
	\end{align*}
	Indeed, the two factors in the integrand have opposite signs on both \((0,s_0)\) and \((s_0,t_\infty)\).
	Since \(t\mapsto F_t(z)\) has at most one zero in \((0,1)\) and \(\partial_tF_t(z)>0\) at every zero, the inequalities \(F_t(z)<0\) for \(0<t\leq t_1\) and \(F_{t_\infty}(z)>0\) imply that \(F_t(z)>0\) for \(t_\infty\leq t<1\).

	For \(t\geq1\), write \(F_t(z)\) as
	\begin{align*}
		F_t(z)=-\int_0^1\phi'(s)K(z,t-s)\,ds.
	\end{align*}
	Since \(\partial_t K<0\), the derivative of \(s\mapsto K(z,t-s)\) is \(\lvert(\partial_t K)(z,t-s)\rvert\) for \(s<t\).
	Integration by parts in \(s\), using \(\phi(0)=\phi(1)=0\), therefore gives
	\begin{align*}
	    F_t(z)=\int_0^1\phi(s)\lvert(\partial_t K)(z,t-s)\rvert\,ds>0.
	\end{align*}
	This proves {\rm(i)}.

	\underline{Proof of (ii).}\quad
	Consider \(D(s):=-\int_0^s\phi'(t_\infty-u)u^{-1/2}\,du\) for \(0\leq s\leq t_\infty\).
	Since \(T_0(t_\infty)=0\), \(D(0)=D(t_\infty)=0\), and {\rm\textbf{(T2)}} gives \(D(s)>0\) for \(0<s<t_\infty\).
	Moreover,
	\begin{align*}
		\frac{D(s)}{\sqrt s} =-\int_0^1\phi'(t_\infty-sv)v^{-\frac12}\,dv.
	\end{align*}
	If \(0<s\leq t_\infty-t_1\), then \(t_\infty-sv\in[t_1,t_\infty]\) for \(0\leq v\leq1\), so {\rm\textbf{(T5)}} implies that \(D(s)/\sqrt{s}\) is nonincreasing in \(s\).
	If \(t_\infty-t_1<s<t_\infty\), then \(D'(s)<0\).
	Hence
	\begin{align}\label{0911-37}
		\frac{D(s)}{\sqrt s} \quad\hbox{is nonincreasing for }0<s<t_\infty.
	\end{align}

	Applying \eqref{0324-14} in \eqref{260038}, using Fubini's theorem, setting \(\lambda=u/(4s)\), and integrating by parts in \(s\), we have
	\begin{align}\label{0911-38}
F_{t_\infty}(z)& = -\frac{\Gamma\left(\frac{n+3}{2}\right)}{\Gamma\left(\frac{n}{2}\right)\Gamma\left(\frac32\right)}\int_0^1\int_0^{t_\infty} \frac{\phi'(t_\infty-s)}{s^{\frac{n+1}{2}}}e^{-\frac{zu}{4s}}u^{\frac{n-2}{2}}(1-u)^{\frac12}\,ds du\nonumber\\
& = \frac{2^{n+1}\Gamma\left(\frac{n+3}{2}\right)}{\Gamma\left(\frac{n}{2}\right)\Gamma\left(\frac32\right)}\int_0^\infty e^{-z\lambda}\lambda^{\frac{n-3}{2}} q(\lambda)\,d\lambda,
	\end{align}
	where
	\begin{align}\label{0911-39}
		q(\lambda):=\lambda^{\frac32} \int_0^{\min\{t_\infty,(4\lambda)^{-1}\}} \frac{D(s)}{\sqrt{1-4\lambda s}}\,ds.
	\end{align}
	Extend \(D\) by zero for \(s\geq t_\infty\).
	Since \(D(t_\infty)=0\), \eqref{0911-37} shows that \(D(s)/\sqrt{s}\) is nonnegative and nonincreasing for all \(s>0\).
	Also, \(D(s)/\sqrt{s}\leq2\|\phi'\|_{L^\infty(0,1)}\) by the definition of \(D\).
	Setting \(v=4\lambda s\) in \eqref{0911-39} now gives a formula with fixed integration limits for every \(\lambda>0\),
	\begin{align*}
		q(\lambda)=\frac18\int_0^1\sqrt{\frac{v}{1-v}}\,
		\frac{D(v/(4\lambda))}{\sqrt{v/(4\lambda)}}\,dv.
	\end{align*}
	For each \(v\in(0,1)\), the integrand is nondecreasing in \(\lambda\).
	Hence \(q\) is nondecreasing and \(0\leq q(\lambda)\leq (\pi/8)\|\phi'\|_{L^\infty(0,1)}\).

	Set \(a=(n-1)/2\).
	Substituting \(v=z\lambda\) in the last line of \eqref{0911-38} gives
	\begin{align*}
		z^a F_{t_\infty}(z)=
		\frac{2^{n+1}\Gamma((n+3)/2)}{\Gamma(n/2)\Gamma(3/2)}
		\int_0^\infty e^{-v}v^{a-1}q(v/z)\,dv.
	\end{align*}
	Since \(q\) is nondecreasing, the right-hand side is nonincreasing in \(z\).
	Hence \(\frac{d}{dz}[z^a F_{t_\infty}(z)]\leq0\), which gives \eqref{0911-35}.
\end{proof}

The next lemma collects the elementary bounds and asymptotics used for large times.
\begin{lemma}\label{lem-261080}
	Let \(n\geq2\), and let \(p\) be the pressure of the Stokes system \eqref{0308-01} with boundary data \eqref{0324-13}, where \(i=j=1\).
	Suppose that \(\psi\) satisfies {\rm\textbf{(S)}} and is radial, and that \(\phi\) satisfies {\rm\textbf{(T1)}} and \(\phi\not\equiv0\).
	Then the following hold.
	\begin{enumerate}
		\item The function
		\begin{align*}
			z\longmapsto n\,{}_1F_1\left(\frac{n+2}{2},\frac{n+1}{2};-z\right)-(n-1)\,{}_1F_1\left(\frac n2,\frac{n+1}{2};-z\right),\qquad z\geq0,
		\end{align*}
		equals \(1\) at \(z=0\), is positive on \([0,z_7)\) and negative on \((z_7,\infty)\), and has a negative derivative at \(z_7\).
		As \(z\to\infty\), it equals \(-\frac{(n-1)\Gamma(\frac{n+1}{2})}{\sqrt\pi}z^{-\frac n2}\bigl(1+O(z^{-1})\bigr)\).
		Moreover, its value at \(z=\varrho^2/4\), \(\varrho\geq0\), equals \(\frac{d}{d\varrho}\bigl[\varrho\,{}_1F_1\bigl(\frac n2,\frac{n+1}{2};-\frac{\varrho^2}{4}\bigr)\bigr]\).
		\item There exists \(C(n)>0\) such that, for every multi-index \(\beta\) with \(|\beta|\leq3\) and every \(\xi'\in\Sigma\),
		\begin{align*}
			\left|\partial_{\xi'}^\beta\partial_{\xi_1}\left[\xi_1\,{}_1F_1\left(\frac n2,\frac{n+1}{2};-\frac{|\xi'|^2}{4}\right)\right]\right| \leq\frac{C(n)}{(1+|\xi'|)^{n+|\beta|}}.
		\end{align*}
		\item For \(t>1\) and \(r\geq1\),
		\begin{align}\label{261081}
			\partial_{x_1}p(r\mathbf e_1',0,t) & =-\frac{(n+1)\Gamma\left(\frac{n+2}{2}\right)} {n\,2^{n+1}\pi^{\frac n2}\Gamma\left(\frac{n+3}{2}\right)} \int_0^1\frac{\phi(u)}{(t-u)^{\frac{n+3}{2}}}\nonumber\\
			& \qquad\times\int_\Sigma\psi(y') \left.\partial_{x_1}\left[x_1\,{}_1F_1\left(\frac n2,\frac{n+1}{2};-\frac{|x'|^2}{4(t-u)}\right)\right]\right|_{x'=r\mathbf e_1'-y'}\,dy'du.
		\end{align}
		At \(r=1\), the left-hand side denotes \(\lim_{x_n\downarrow0}\partial_{x_1}p(\mathbf e_1',x_n,t)\). Moreover, \(\partial_{x_1}p(-r\mathbf e_1',0,t)=\partial_{x_1}p(r\mathbf e_1',0,t)\) for every \(r\geq1\).
	\item There exists \(\delta_0=\delta_0(n)\in(0,\sqrt{z_7})\) such that, for every \(\delta\in(0,\delta_0]\), there exists \(T_\delta>2\), depending only on \(n\) and \(\delta\), such that, for every \(t\geq T_\delta\),
		\begin{align}\label{261082}
			\begin{aligned}
				\partial_{x_1}p(r\mathbf e_1',0,t) & <0, && 1\leq r\leq(2\sqrt{z_7}-\delta)\sqrt{t-1},\\
				\partial_r\partial_{x_1}p(r\mathbf e_1',0,t) & >0, && (2\sqrt{z_7}-\delta)\sqrt{t-1}<r<(2\sqrt{z_7}+\delta)\sqrt t,\\
				\partial_{x_1}p(r\mathbf e_1',0,t) & >0, && r\geq(2\sqrt{z_7}+\delta)\sqrt t.
			\end{aligned}
		\end{align}
		For \(t\geq T_{\delta_0}\), there is exactly one zero \(\rho(t)\) on \([1,\infty)\).
		It is simple and satisfies
		\begin{align*}
			\bigl|\rho(t)-2\sqrt{z_7t}\bigr|\leq C(n)t^{-1/2}.
		\end{align*}
	\end{enumerate}
\end{lemma}
\begin{proof}
	Throughout the proof, \(b=(n+1)/2\).

	\underline{Proof of (i).}\quad
	By \eqref{0327-04} and \eqref{260030}, \(e^z\) times the function in {\rm(i)} equals
	\begin{align*}
		n\,{}_1F_1\left(-\frac12,b;z\right)-(n-1)\,{}_1F_1\left(\frac12,b;z\right) =\sum_{k=0}^\infty c_k\frac{z^k}{k!}, \qquad c_k:=\frac{n(-\frac12)_k-(n-1)(\frac12)_k}{(b)_k}.
	\end{align*}
	Here \(c_0=1\), and \(c_k<0\) for \(k\geq1\) since \((-\frac12)_k<0<(\frac12)_k\).
	Hence this series equals \(1\) at \(z=0\), has negative derivative on \([0,\infty)\), and is at most \(1+c_1z\), so its only positive zero is \(z_7\), positive before and negative after it, with negative derivative there.
Multiplying by \(e^{-z}\) preserves these signs and, since the series vanishes at \(z_7\), the sign of the derivative at \(z_7\).
	The asymptotics follow from the second line of \eqref{0327-10}, which gives \({}_1F_1(\frac n2,b;-z)=\frac{\Gamma(b)}{\sqrt\pi}z^{-\frac n2}(1+O(z^{-1}))\) and \({}_1F_1(\frac{n+2}{2},b;-z)=O(z^{-\frac{n+2}{2}})\).
	Finally, let \(z=\varrho^2/4\).
	Then \(\frac{d}{d\varrho}\bigl[\varrho\,{}_1F_1(\frac n2,b;-z)\bigr] ={}_1F_1(\frac n2,b;-z)+2z\frac{d}{dz}{}_1F_1(\frac n2,b;-z)\), and \eqref{0327-01} with \(a=\frac n2\) gives
	\begin{align*}
		z\frac{d}{dz}{}_1F_1\left(\frac n2,b;-z\right) =\frac n2\left[{}_1F_1\left(\frac{n+2}{2},b;-z\right)-{}_1F_1\left(\frac n2,b;-z\right)\right].
	\end{align*}
	Substituting this identity gives the last claim in {\rm(i)}.

	\underline{Proof of (ii).}\quad
	Take \(q:=|\xi'|^2/4\).
	By \eqref{0327-07} and \eqref{0514-02}, for \(j\geq0\),
	\begin{align*}
		\left|\frac{d^j}{dq^j}{}_1F_1\left(\frac n2,b;-q\right)\right| =\frac{(\frac n2)_j}{(b)_j}\left|{}_1F_1\left(\frac n2+j,b+j;-q\right)\right| \leq\frac{C_j}{(1+q)^{\frac n2+j}}.
	\end{align*}
	The function in {\rm(ii)} equals \({}_1F_1(\frac n2,b;-q)+\frac{\xi_1^2}{2}\frac{d}{dq}{}_1F_1(\frac n2,b;-q)\).
	Both terms have the form \(\xi'^\gamma\frac{d^j}{dq^j}{}_1F_1(\frac n2,b;-q)\) with \(2j-|\gamma|=0\).
	Each derivative \(\partial_{\xi_l}\) maps such a term to a combination of two terms of the same form, in which either \(|\gamma|\) decreases by one or both \(j\) and \(|\gamma|\) increase by one.

	Hence \(\partial_{\xi'}^\beta\) of the function in {\rm(ii)} is a finite combination, with coefficients depending only on \(n\) and \(\beta\), of terms of this form with \(2j-|\gamma|=|\beta|\).
	Each of them is bounded by \(C|\xi'|^{|\gamma|}(1+|\xi'|^2/4)^{-\frac n2-j}\leq C(1+|\xi'|)^{-n-|\beta|}\).

	\underline{Proof of (iii).}\quad
	Fix \(t>1\) and write \(c_n:=\Gamma(\frac{n+2}{2})/(n2^n\pi^{\frac n2}\Gamma(\frac{n+3}{2}))\).
	Integrate by parts in \eqref{260038} using \(\phi(0)=\phi(1)=0\), and apply \eqref{0911-21} with \((a,b,c)=(\frac n2,\frac{n+3}{2},\frac z4)\).
	For \(z>0\) this gives
	\begin{align*}
		F_t(z)=\frac{n+1}{2}\int_0^1\frac{\phi(u)}{(t-u)^{\frac{n+3}{2}}}
		{}_1F_1\left(\frac n2,\frac{n+1}{2};-\frac{z}{4(t-u)}\right)du .
	\end{align*}
	Since \(t-u\geq t-1>0\), this formula extends \(F_t\) smoothly to \(z=0\).

	In \eqref{0822-10}, use \(\phi(t)=0\) and change variables to \(u=t-s\).
	Substituting the preceding identity gives, for \(|x'|>1\) in every dimension \(n\geq2\),
	\begin{align*}
		p(x',0,t)=-c_n\int_\Sigma\psi(y')(x_1-y_1)F_t(|x'-y'|^2)\,dy'.
	\end{align*}
	Since \(t-u\geq t-1\), differentiation in \(x_1\) under the spatial integral, or equivalently integration by parts in \(y_1\), gives, on the \(x_1\)-axis,
	\begin{align}\label{261083}
		\partial_{x_1}p(r\mathbf e_1',0,t)=-c_n\int_\Sigma\partial_{y_1}\psi(y')(r-y_1)F_t(|r\mathbf e_1'-y'|^2)\,dy'.
	\end{align}
	Using the integral for \(F_t\) in \eqref{261083} gives \eqref{261081} for \(r>1\).
	On \(1\leq r\leq2\), the integrand is bounded by \(C(n,t)\phi(u)|\psi(y')|\).
	Dominated convergence as \(r\downarrow1\), together with continuity of \(x'\mapsto\partial_{x_1}p(x',0,t)\) from \eqref{tangential-pressure-boundary-value}, gives \eqref{261081} at \(r=1\).
	Finally, substituting \(y'\mapsto-y'\) in \eqref{0822-10} shows that \(p(\cdot,0,t)\) is odd under \(x'\mapsto-x'\) on \(\{|x'|>1\}\).
	Its \(x_1\)-derivative is therefore even there. Letting \(r\downarrow1\) gives equality at \(\mathbf e_1'\) and \(-\mathbf e_1'\).

	\underline{Proof of (iv).}\quad
	For \(\xi'\in\Sigma\), let
	\begin{align*}
		\Theta(\xi'):=\partial_{\xi_1}\left[\xi_1\,{}_1F_1\!\left(\frac n2,\frac{n+1}{2};-\frac{|\xi'|^2}{4}\right)\right]
	\end{align*}
	and restrict it to the first coordinate axis
	\begin{align*}
		\theta(\varrho):=\Theta(\varrho\mathbf e_1')= \frac{d}{d\varrho}\left[\varrho\,{}_1F_1\!\left(\frac n2,\frac{n+1}{2};-\frac{\varrho^2}{4}\right)\right].
	\end{align*}
Set \(\varrho_*:=2\sqrt{z_7}\).
By {\rm(i)}, \(\theta\) is positive on \([0,\varrho_*)\), negative on \((\varrho_*,\infty)\), satisfies \(\theta'(\varrho_*)<0\), and obeys \(-\theta(\varrho)\geq c(n)\varrho^{-n}\) for all sufficiently large \(\varrho\).
For \(\varrho\geq0\) and \(0<\varepsilon\leq1\), define
	\begin{align*}
		\Pi(\varrho,\varepsilon):=\int_\Sigma\psi(y')\Theta(\varrho\mathbf e_1'-\varepsilon y')\,dy'.
	\end{align*}
	Radiality gives \(\int_\Sigma y'\psi(y')\,dy'=0\).
	Taylor's theorem, {\rm(ii)}, and \(|y'|\leq1\) on \(\supp\psi\) therefore give, for \(j=0,1\),
	\begin{align}\label{261085}
		\left|\partial_\varrho^j\Pi(\varrho,\varepsilon)-M_0\theta^{(j)}(\varrho)\right|\leq C(n)M_2\varepsilon^2(1+\varrho)^{-n-2-j}.
	\end{align}
	Choose \(\delta_0\in(0,\varrho_*/4)\) so that \(\theta'(\varrho)<0\) whenever \(|\varrho-\varrho_*|\leq2\delta_0\); then \(\delta_0<\sqrt{z_7}\).
	Fix \(\delta\in(0,\delta_0]\).
	By these properties of \(\theta\) and continuity, there is \(c_\delta>0\), depending only on \(n\) and \(\delta\), such that
	\begin{align*}
		\theta(\varrho)&\geq c_\delta &&(0\leq\varrho\leq\varrho_*-\delta),\\
		\theta(\varrho)&\leq-c_\delta(1+\varrho)^{-n} &&(\varrho\geq\varrho_*+\delta),\\
		\theta'(\varrho)&\leq-c_\delta &&(|\varrho-\varrho_*|\leq2\delta_0).
	\end{align*}
	Choose \(T_\delta>2\), depending only on \(n\) and \(\delta\), so large that
	\begin{align*}
		(\varrho_*-\delta)\sqrt{T_\delta-1}&\geq1,\\
		\left[(\varrho_*-\delta)\sqrt{\frac{t-1}{t}},(\varrho_*+\delta)\sqrt{\frac{t}{t-1}}\right]&\subset[\varrho_*-2\delta_0,\varrho_*+2\delta_0] \qquad(t\geq T_\delta),\\
		\frac{C(n)}{T_\delta-1}&\leq\frac{c_\delta}{2},
	\end{align*}
where \(C(n)\) is chosen at least as large as the constant in \eqref{261085}.
	Fix \(t\geq T_\delta\), and, for \(0\leq u\leq1\), set \(s:=t-u\), \(\varepsilon:=s^{-\frac12}\), and \(\varrho_s:=r/\sqrt s\).
	Since \(\partial_{x_1}[x_1\,{}_1F_1(\frac n2,b;-\frac{|x'|^2}{4s})]=\Theta(x'/\sqrt s)\), \eqref{261081} becomes
	\begin{align*}
		\partial_{x_1}p(r\mathbf e_1',0,t)=-\frac{(n+1)c_n}{2}\int_0^1\frac{\phi(u)}{s^{\frac{n+3}{2}}}\Pi(\varrho_s,\varepsilon)\,du.
	\end{align*}
	Because \(M_2\leq M_0\), \eqref{261085} and the choice of \(T_\delta\) imply that \(\Pi(\varrho_s,\varepsilon)>0\) when \(1\leq r\leq(\varrho_*-\delta)\sqrt{t-1}\), and \(\Pi(\varrho_s,\varepsilon)<0\) when \(r\geq(\varrho_*+\delta)\sqrt t\).
	Since \(\phi\geq0\) and \(\phi\not\equiv0\), these two signs give the first and third lines of \eqref{261082}.
In the remaining range, the last two conditions defining \(T_\delta\), together with \eqref{261085}, give \(\partial_\varrho\Pi(\varrho_s,\varepsilon)<0\).
	Differentiating the preceding integral formula with respect to \(r\) therefore gives \(\partial_r\partial_{x_1}p(r\mathbf e_1',0,t)>0\), which is the second line of \eqref{261082}.

	Take \(\delta=\delta_0\).
	The opposite signs at the two endpoints in \eqref{261082}, strict increase between them, and the signs outside give exactly one zero \(r=\rho(t)>1\) for \(t\geq T_{\delta_0}\).
	Set \(v:=\rho(t)/\sqrt t\).
	The choice of \(T_{\delta_0}\) places \(v\) and every \(\varrho_s=\rho(t)/\sqrt{t-u}\) in \([\varrho_*-2\delta_0,\varrho_*+2\delta_0]\), and
	\begin{align*}
		|\varrho_s-v|=v\bigl|(1-u/t)^{-1/2}-1\bigr|\leq C(n)/t,\qquad 0\leq u\leq1.
	\end{align*}
	By \eqref{261085}, \(M_2\leq M_0\), and boundedness of \(\theta'\) on this interval,
	\begin{align*}
		\left|\frac{\Pi(\varrho_s,s^{-1/2})}{M_0}-\theta(v)\right|\leq C(n)/t .
	\end{align*}
	At the zero, the average of \(\Pi(\varrho_s,s^{-1/2})\) with weight \(\phi(u)(t-u)^{-(n+3)/2}\) vanishes by \eqref{261081}.
	This weight has positive integral, so \(|\theta(v)|\leq C(n)/t\).
	Since \(\theta(\varrho_*)=0\) and \(-\theta'\geq c_{\delta_0}>0\) throughout the interval, the mean value theorem gives \(|v-\varrho_*|\leq C(n)/t\).
	Multiplying by \(\sqrt t\) proves the claimed bound.
\end{proof}

\subsubsection{The \(i=j\) case when \(n=2\)}\label{subsubsect-260154-n2}

Proposition \ref{prop-260154} gives the signs and zeros of \(\partial_{x_1}p(x_1,0,t)\) in Theorem \ref{thm-260021} for \(n=2\).
\begin{proposition}\label{prop-260154}
		Let \(n=2\), and let \(p\) be the pressure of the Stokes system \eqref{0308-01} with boundary data \eqref{0324-13}, where \(i=j=1\).
	Suppose that \(\psi\) satisfies {\rm\textbf{(S)}} and is even, and that \(\phi\) satisfies {\rm\textbf{(T1)}} and {\rm\textbf{(T2)}}.
	Then the following hold.
	\begin{enumerate}
		\item Suppose that \(\phi\) satisfies {\rm\textbf{(T3)}} and {\rm\textbf{(T4)}}. Then, as \(|x_1|\to\infty\),
		\begin{align}\label{260156}
			\partial_{x_1}p(x_1,0,t) =-\frac{1}{\pi^{\frac{3}{2}}x_1^2}\left[M_0T_0(t) +3\bigl(M_2T_0(t)-2M_0T_1(t)\bigr)x_1^{-2} +O(|x_1|^{-4})\right],
		\end{align}
		uniformly for \(t\) in compact subsets of \((0,\infty)\).
		For every fixed \(t>0\) with \(t\ne t_\infty\),
		\begin{align}\label{1013-109}
			\partial_{x_1}p(x_1,0,t)
			\begin{cases}
				<0, & 0<t<t_\infty, \\
				>0, & t>t_\infty,
			\end{cases}
			\qquad |x_1|\text{ sufficiently large}.
		\end{align}
		At \(t=t_\infty\),
		\begin{align*}
			\partial_{x_1}p(x_1,0,t_\infty)
			=\frac{6M_0T_1(t_\infty)}{\pi^{3/2}x_1^4}+O(|x_1|^{-6}),
		\end{align*}
		and \(\partial_{x_1}p(x_1,0,t_\infty)>0\) for all sufficiently large \(|x_1|\).
		\item Suppose that \(\psi\) is nonincreasing on \((0,\infty)\) and that \(\phi\) satisfies {\rm\textbf{(T3)}} and {\rm\textbf{(T4)}}. Then \(t\mapsto\partial_{x_1}p(\pm1,0,t)\) is continuous on \((0,\infty)\).
		Suppose, in addition, that \(\phi\) satisfies {\rm\textbf{(T5)}}.
		Then, at \(t=t_\infty\),
		\begin{align}\label{260157}
			\partial_{x_1}p(x_1,0,t_\infty)>0, \qquad |x_1|\geq1.
		\end{align}
		\item There exists \(t_0=t_0(\psi,\phi)\in(0,t_1)\) such that \(\partial_{x_1}p(x_1,0,t)<0\) for all \(|x_1|\geq1\) and \(0<t<t_0\).
		Suppose, in addition, that \(\psi\) is nonincreasing on \((0,\infty)\) and that \(\phi\) satisfies {\rm\textbf{(T3)}}--{\rm\textbf{(T5)}}.
		Then \(\overline t\) in \eqref{0911-23} satisfies \(t_\infty<\overline t<\infty\), and
		\begin{align*}
		    \partial_{x_1}p(\pm1,0,\overline t)=0,\qquad
		    \partial_{x_1}p(\pm1,0,t)<0\quad(t>\overline t).
		\end{align*}
		\item Suppose that \(\psi\) is nonincreasing on \((0,\infty)\) and that \(\phi\) satisfies {\rm\textbf{(T3)}}--{\rm\textbf{(T5)}}.
		Then there exists \(\delta_\infty=\delta_\infty(\psi,\phi)\in(0,t_\infty-t_1)\) such that, for \(t_\infty-\delta_\infty<t<t_\infty\), there is exactly one zero \(\rho(t)>1\) on \((1,\infty)\) and \(-\rho(t)\) is the only zero on \((-\infty,-1)\).
		Moreover,
		\begin{align}\label{0927-005}
			\sgn(\partial_{x_1}p(x_1,0,t)) =\sgn(\rho(t)-|x_1|)
		\end{align}
		for \(|x_1|\geq1\) and
		\begin{align*}
			\lim_{t\to t_\infty^-}(t_\infty-t)\rho(t)^2 =-\frac{6T_1(t_\infty)}{T_0'(t_\infty)}.
		\end{align*}
		\item Suppose that \(\psi\) is nonincreasing on \((0,\infty)\) and that \(\phi\) satisfies {\rm\textbf{(T3)}}--{\rm\textbf{(T5)}}.
		Then, for every \(t>\overline t\), the zeros of \(\partial_{x_1}p(x_1,0,t)\) on \((1,\infty)\) have a minimum \(\rho_-(t)\) and a maximum \(\rho_+(t)\), which satisfy
		\begin{align}\label{261014}
			\partial_{x_1}p(x_1,0,t)<0 \quad & \text{for }1\leq|x_1|<\rho_-(t), & \partial_{x_1}p(x_1,0,t)>0 \quad & \text{for }|x_1|>\rho_+(t).
		\end{align}
		Moreover, there exists \(\widehat t>1\) such that \(\overline t<\widehat t\) and, for every \(t\geq\widehat t\), \(\rho(t):=\rho_-(t)=\rho_+(t)\) satisfies \(\partial_{x_1}^2p(\rho(t),0,t)>0\) and
		\begin{align}\label{261015}
			\sgn\bigl(\partial_{x_1}p(x_1,0,t)\bigr)=\sgn\bigl(|x_1|-\rho(t)\bigr)\quad(|x_1|\geq1), \qquad \lim_{t\to\infty}\frac{\rho(t)}{\sqrt t}=2\sqrt{z_7}.
		\end{align}
	\end{enumerate}
\end{proposition}

\begin{proof}[Proof of Proposition \ref{prop-260154}]
	Since \(\psi\) is even, substituting \(y\mapsto-y\) in \eqref{260146} gives \(\partial_{x_1}p(-x_1,0,t)=\partial_{x_1}p(x_1,0,t)\) for \(|x_1|>1\).
	By \eqref{tangential-pressure-boundary-value}, \(\partial_{x_1}p(x_1,0,t)\) is continuous at \(x_1=\pm1\), so the symmetry extends there, and it suffices to consider \(x_1\geq1\).

	\underline{Proof of (i).}\quad
	Subtracting the first two terms of the kernel in \eqref{260155} leaves
	\begin{align*}
		E(x_1-y,s)=\frac{x_1-y}{s^{3/2}}{}_1F_1\left(1,\frac52;-\frac{(x_1-y)^2}{4s}\right)-\frac{6}{\sqrt{s}(x_1-y)}+\frac{12\sqrt{s}}{(x_1-y)^3}.
	\end{align*}
	For \(P(z)={}_1F_1(1,5/2;-z)\) and \(z=(x_1-y)^2/(4s)\), the first two \(x_1\)-derivatives of \((x_1-y)P(z)\) are \(P+2zP'\) and \(2z(3P'+2zP'')/(x_1-y)\), so \eqref{0327-07}, \eqref{0504-06}, and \(|x_1-y|\geq|x_1|-1\) on \(\supp\psi\) give
	\begin{align*}
		|\partial_{x_1}^j E(x_1-y,s)|\leq C s^{3/2}|x_1-y|^{-5-j},\qquad j=0,1,2,
	\end{align*}
	uniformly for \(0<s\leq t\), \(t\) in a compact subset of \((0,\infty)\), and \(|x_1|\) sufficiently large.
	Since \(\int_\R\psi'(y)E(x_1-y,s)\,dy=\int_\R\psi(y)\partial_{x_1}E(x_1-y,s)\,dy\), the bound for \(j=1\) gives
	\begin{align}\label{1012-04}
		\partial_{x_1}p(x_1,0,t) & =\frac{1}{\pi^{\frac32}}\left[ T_0(t)\int_\R\frac{\psi'(y)}{x_1-y}\,dy -2T_1(t)\int_\R\frac{\psi'(y)}{(x_1-y)^3}\,dy \right]+O(|x_1|^{-6}).
	\end{align}
	Integrating by parts in \(y\) and using the vanishing of the odd moments of \(\psi\), we obtain
	\begin{align}\label{0927-0011}
		\begin{split}
			\int_\R\frac{\psi'(y)}{x_1-y}\,dy & =-\int_\R\frac{\psi(y)}{(x_1-y)^2}\,dy =-\frac{M_0}{x_1^2}-\frac{3M_2}{x_1^4}+O(x_1^{-6}), \\
			\int_\R\frac{\psi'(y)}{(x_1-y)^3}\,dy & =-3\int_\R\frac{\psi(y)}{(x_1-y)^4}\,dy =-\frac{3M_0}{x_1^4}+O(x_1^{-6}),
		\end{split}
	\end{align}
	and \eqref{1012-04} gives \eqref{260156}.
	By Remark \ref{rmk-0730-03} and {\rm\textbf{(T4)}}, \(T_0>0\) on \((0,t_\infty)\) and \(T_0<0\) on \((t_\infty,\infty)\), so \eqref{260156} gives \eqref{1013-109}.
	At \(t=t_\infty\), the coefficient of \(x_1^{-2}\) vanishes, and \(T_1(t_\infty)>0\) by the proof of Proposition \ref{prop-0504-03}.

	\underline{Proof of (ii).}\quad
	Since \(|\phi(t)-\phi(t-s)|\leq C\min\{s,1\}\), \eqref{0514-02} gives, for \(1\leq x_1\leq2\) and \(y<1\),
	\begin{align}\label{260295}
		 & \int_0^\infty|\phi(t)-\phi(t-s)| \left|\frac{3}{s^{\frac52}} {}_1F_1\left(1,\frac32;-\frac{(x_1-y)^2}{4s}\right) -\frac{(x_1-y)^2}{s^{\frac72}} {}_1F_1\left(2,\frac52;-\frac{(x_1-y)^2}{4s}\right)\right|\,ds\nonumber \\
		 & \quad\leq C\left[\int_0^1 \frac{s^{-\frac12}(s+2(x_1-y)^2)} {(s+(x_1-y)^2)^2}\,ds +\int_1^\infty s^{-\frac52}\,ds\right] \leq C\left(1+\frac1{1-y}\right).
	\end{align}
	Since \(\psi(y)=O((1-y)^{2+\alpha})\) as \(y\to1^-\), dominated convergence in \eqref{0924-04} extends \(\partial_{x_1}p(x_1,0,t)\) continuously to \(x_1=1\), and, by \eqref{tangential-pressure-boundary-value}, the extension is the boundary value of \(\partial_{x_1}p\).

	Assume {\rm\textbf{(T5)}}, and let \(F=F_{t_\infty}\) as in \eqref{260038}.
	By \eqref{260155},
	\begin{align}\label{1012-05}
		\partial_{x_1}p(x_1,0,t_\infty) =-\frac1{6\pi^{\frac32}}\int_\R \psi'(y)(x_1-y)F((x_1-y)^2)\,dy,\qquad |x_1|>1.
	\end{align}
	By Lemma \ref{lem-0911-02b}{\rm(i)} and \eqref{0514-02}, \(0<F(z)\leq C\norm{\phi'}_{L^\infty(0,1)}\int_0^{t_\infty}s^{-1/2}(s+z)^{-1}\,ds\leq C\norm{\phi'}_{L^\infty(0,1)}z^{-1/2}\) for \(z>0\), so, since \(\psi'\in L^1(\R)\), \eqref{1012-05} extends to \(x_1=1\).
	By \eqref{0911-35} with \(n=2\), \(\frac{d}{dv}(vF(v^2))=F(v^2)+2v^2F'(v^2)\leq0\) for \(v>0\).
	Moreover, \(vF(v^2)\) is positive, real-analytic by \eqref{0911-38}, and tends to \(0\) as \(v\to\infty\), since \eqref{0504-06}, \eqref{260038}, and \(T_0(t_\infty)=0\) give \(F(z)=12T_1(t_\infty)z^{-2}+O(z^{-3})\).
	Hence \(vF(v^2)\) is strictly decreasing, and pairing \(y\) and \(-y\) in \eqref{1012-05} gives, for \(x_1\geq1\),
	\begin{align*}
		\partial_{x_1}p(x_1,0,t_\infty) & =-\frac1{6\pi^{\frac32}}\int_0^1\psi'(y) \left[(x_1-y)F((x_1-y)^2) -(x_1+y)F((x_1+y)^2)\right]\,dy>0,
	\end{align*}
	since the bracket is positive, \(\psi'\leq0\) on \((0,1)\), and \(\psi'\not\equiv0\).

	\underline{Proof of (iii).}\quad
	Let
	\begin{align}\label{0928-40}
		k(s):=s^{-\frac{5}{2}}\int_{\R}\psi(y)\left[(1-y)^2 {}_1F_1\left(2,\frac72;-\frac{(1-y)^2}{4s}\right) -5s{}_1F_1\left(1,\frac52;-\frac{(1-y)^2}{4s}\right)\right]\,dy,\quad s>0.
	\end{align}
	For \(1\leq x_1\leq2\), the bound \eqref{0514-02} for small \(s\) and the bound \(Cs\) for the bracket in \eqref{260158} when \(s\geq1\) give the majorant \(C|\psi(y)|(1-y)^{-2}\min\{s^{-1/2},s^{-3/2}\}\), which is integrable since \(\psi(y)=O((1-y)^{2+\alpha})\) near \(y=1\).
	Hence dominated convergence in \eqref{260158} gives
	\begin{align}\label{0924-01}
		\partial_{x_1}p(1,0,t)=\lim_{x_1\to 1^+}\partial_{x_1}p(x_1,0,t) =-\frac1{30\pi^{\frac{3}{2}}}\int_0^t\phi'(t-s)k(s)\,ds.
	\end{align}
	For \(y<1\), \(s^{-2}\) times the bracket in \eqref{0928-40} tends to \(30(1-y)^{-2}\) as \(s\to0^+\) by \eqref{0327-10}, and is bounded by \(C(1-y)^{-2}\) by \eqref{0514-02}.
	Since \(\psi(y)(1-y)^{-2}\in L^1(\R)\) by {\rm\textbf{(S)}}, dominated convergence gives
	\begin{align}\label{0924-02}
		k(s)\sim 30s^{-\frac{1}{2}}\int_{\R}\frac{\psi(y)}{(1-y)^2}\,dy>0 \quad(s\to0^+).
	\end{align}
	Choose \(s_0\in(0,\min\{t_1,1\})\) so that \(k>0\) on \((0,s_0)\).
	Since \(\phi'(t-s)>0\) for \(0<s<t<s_0\), \eqref{0924-01} gives \(\partial_{x_1}p(1,0,t)<0\) for \(0<t<s_0\).

	For \(x_1>1\), integrating by parts in \(y\) in \eqref{260155} and applying \eqref{0327-07} give
	\begin{align}\label{261012}
		\partial_{x_1}p(x_1,0,t) =\frac{1}{6\pi^{\frac32}}\int_0^t\frac{\phi'(t-s)}{s^{\frac32}} \int_{\R}\psi(y)\biggl[ {}_1F_1\left(1,\frac52;-\frac{(x_1-y)^2}{4s}\right) -\frac{(x_1-y)^2}{5s} {}_1F_1\left(2,\frac72;-\frac{(x_1-y)^2}{4s}\right) \biggr]\,dyds.
	\end{align}
	By \eqref{0504-06}, \({}_1F_1(1,\frac52;-z) -\frac{4z}{5}{}_1F_1(2,\frac72;-z) =-\frac{3}{2z}+O(z^{-2})\) as \(z\to\infty\), so there are \(Z,c,C>0\) such that
	\begin{align}\label{2610131}
		{}_1F_1\left(1,\frac52;-z\right) -\frac{4z}{5}{}_1F_1\left(2,\frac72;-z\right) \leq-\frac{c}{z}\quad(z\geq Z),
	\end{align}
	and the left-hand side of \eqref{2610131} is at most \(C\) for \(0<z<Z\).
	By {\rm\textbf{(S)}}, there are \(a\in(0,1)\) with \(\int_{-a}^a\psi(y)\,dy>0\) and \(C_\psi>0\) with \(0\leq\psi(y)\leq C_\psi(1-y)^{2+\alpha}\) for \(0\leq y\leq1\).
	Let \(t<t_1\), \(t\leq(1-a)^2/(4Z)\), \(2\sqrt{Zt}\leq1\), \(1<x_1\leq1+2\sqrt{Zt}\), and \(0<s\leq t\), and split the spatial integral in \eqref{261012} over \([-a,a]\), \(\{y\notin[-a,a]:(x_1-y)^2\geq4Zs\}\), and \(\{y:(x_1-y)^2<4Zs\}\).
	On \([-a,a]\), \(1-a\leq x_1-y\leq2+a\) and \((x_1-y)^2/(4s)\geq Z\), so \eqref{2610131} bounds this part, divided by \(s^{3/2}\), by \(-c_0s^{-1/2}\) with \(c_0:=4c(2+a)^{-2}\int_{-a}^a\psi(y)\,dy\).
	On the second set, \eqref{2610131} makes the integrand nonpositive.
	On the third set, points of \(\supp\psi\) satisfy \(0\leq 1-y\leq x_1-y<2\sqrt{Zs}\), so this part, divided by \(s^{3/2}\), is at most \(CC_\psi s^{-3/2}\int_0^{2\sqrt{Zs}}q^{2+\alpha}\,dq \leq Cs^{\frac{\alpha}{2}}\).
	Since \(\phi'(t-s)>0\) for \(0<s<t<t_1\) and \(s^{\alpha/2}\leq t^{(1+\alpha)/2}s^{-1/2}\), \eqref{261012} gives
	\begin{align*}
		\partial_{x_1}p(x_1,0,t) \leq\frac{T_0(t)}{6\pi^{\frac32}} \left[-c_0+Ct^{\frac{1+\alpha}{2}}\right].
	\end{align*}
	Let \(t_0:=\frac12\min\{t_1,s_0,\frac{(1-a)^2}{4Z}, \frac1{4Z},(\frac{c_0}{C})^{\frac{2}{1+\alpha}}\}\).
	For \(0<t<t_0\), this gives \(\partial_{x_1}p(x_1,0,t)<0\) for \(1<x_1\leq1+2\sqrt{Zt}\), and, for \(x_1\geq1+2\sqrt{Zt}\), \((x_1-y)^2/(4s)\geq Z\) on \(\supp\psi\) for \(0<s\leq t\), so \eqref{261012} and \eqref{2610131} give the same sign.
	With \(x_1=1\) and evenness, this proves the first assertion.

	By \eqref{260030}, applied separately to the factors involving the confluent hypergeometric functions in \eqref{0928-40} and their \(s\)-derivatives,
	\begin{align}\label{0924-03}
		k(s)=-5M_0s^{-\frac{3}{2}}+O(s^{-\frac{5}{2}}),\qquad k'(s)=\frac{15M_0}{2}s^{-\frac{5}{2}}+O(s^{-\frac{7}{2}}).
	\end{align}
	Hence, for sufficiently large \(t>1\), \(k'(t-u)>0\) for \(0\leq u\leq1\), and, since \(\phi(0)=\phi(1)=0\), \(\int_0^1\phi'(u)k(t-u)\,du =\int_0^1\phi(u)k'(t-u)\,du>0\), so \eqref{0924-01} gives \(\partial_{x_1}p(1,0,t)<0\).
	Under {\rm\textbf{(T5)}}, by {\rm(ii)}, \(t\mapsto\partial_{x_1}p(1,0,t)\) is continuous and positive at \(t_\infty\), and it is negative for large \(t\).
	Hence its last zero \(\overline t\) satisfies \(t_\infty<\overline t<\infty\), and it is negative on \((\overline t,\infty)\); evenness gives the same at \(x_1=-1\).

	\underline{Proof of (iv).}\quad
	By {\rm(ii)} and \eqref{260156}, for \(t<t_\infty\) close to \(t_\infty\), \(\partial_{x_1}p(1,0,t)>0\) and \(\partial_{x_1}p(x_1,0,t)<0\) for large \(x_1\), so there is a zero on \((1,\infty)\).
	Fix \(L>2\) and a compact interval \(K\subset(0,\infty)\).
	As in {\rm(ii)}, \eqref{0514-02} and \(|\phi(t)-\phi(t-s)|\leq C\min\{s,1\}\) bound the integrand in \eqref{0924-04} by \(C_{K,L}|\psi(y)|\min\{s^{-1/2},s^{-5/2}\}\) for \(3/2\leq x_1\leq L\) and \(t\in K\), which, with {\rm(ii)}, gives joint continuity on \([1,L]\times K\).
	Hence \eqref{260157} and uniform continuity on \([1,L]\times K\) give
	\begin{align}\label{261052}
		0<\min_{1\leq x_1\leq L}\partial_{x_1}p(x_1,0,t_\infty), \qquad \sup_{1\leq x_1\leq L} \left|\partial_{x_1}p(x_1,0,t) -\partial_{x_1}p(x_1,0,t_\infty)\right|\longrightarrow0 \quad(t\to t_\infty^-),
	\end{align}
	so \(\partial_{x_1}p>0\) on \([1,L]\) for \(t<t_\infty\) close to \(t_\infty\), and all zeros in \((1,\infty)\) tend to infinity as \(t\to t_\infty^-\).

	Let \(\rho(t)\) be such a zero.
	By \eqref{260156},
	\(0=M_0T_0(t)\rho(t)^2 +3(M_2T_0(t)-2M_0T_1(t))+o(1)\), so \(T_0(t)\rho(t)^2\to6T_1(t_\infty)\), and, since \(T_0'(t_\infty)<0\) by Remark~\ref{rmk-0730-03},
	\begin{align*}
		(t_\infty-t)\rho(t)^2 =\frac{T_0(t)\rho(t)^2}{T_0(t)/(t_\infty-t)} \longrightarrow-\frac{6T_1(t_\infty)}{T_0'(t_\infty)}.
	\end{align*}
	The bounds on \(\partial_{x_1}E\) and \(\partial_{x_1}^2E\) from {\rm(i)} permit differentiation under the integral in \eqref{260155} at large \(x_1\), and integration by parts gives
	\begin{align}\label{260280}
		\left|\int_\R\psi'(y)\partial_{x_1} E(x_1-y,s)\,dy\right|
		=\left|\int_\R\psi(y)\partial_{x_1}^2 E(x_1-y,s)\,dy\right|
		\leq C M_0s^{3/2}(x_1-1)^{-7}.
	\end{align}
	Differentiating \eqref{0927-0011} gives
	\begin{align}\label{260281}
		\begin{split}
			\frac{d}{dx_1}\int_\R\frac{\psi'(y)}{x_1-y}\,dy & =2\int_\R\frac{\psi(y)}{(x_1-y)^3}\,dy =\frac{2M_0}{x_1^3}+\frac{12M_2}{x_1^5} +O(x_1^{-7}),\\
			\frac{d}{dx_1}\int_\R\frac{\psi'(y)}{(x_1-y)^3}\,dy & =12\int_\R\frac{\psi(y)}{(x_1-y)^5}\,dy =\frac{12M_0}{x_1^5}+O(x_1^{-7}),
		\end{split}
	\end{align}
	and \eqref{260280}--\eqref{260281} give, uniformly for \(t\) in compact subsets of \((0,\infty)\),
	\begin{align}\label{0928-41}
		\partial_{x_1}^2p(x_1,0,t) =\frac1{\pi^{\frac32}x_1^3} \left[2M_0T_0(t) +12\bigl(M_2T_0(t)-2M_0T_1(t)\bigr)x_1^{-2} +O(x_1^{-4})\right].
	\end{align}
	At \(x_1=\rho(t)\), \(\pi^{\frac32}\rho(t)^5 \partial_{x_1}^2p(\rho(t),0,t) =2M_0T_0(t)\rho(t)^2 +12(M_2T_0(t)-2M_0T_1(t))+o(1) \to-12M_0T_1(t_\infty)<0\), uniformly over the zeros in \((1,\infty)\), since they tend to infinity and the remainders in \eqref{260156} and \eqref{0928-41} are uniform near \(t_\infty\).
	Hence there are \(L_0>2\) and \(0<\delta_1<t_\infty-t_1\) such that \(\partial_{x_1}^2p(x_1,0,t)<0\) at every zero \(x_1>L_0\) for \(t_\infty-\delta_1<t<t_\infty\), and, by \eqref{261052} with \(L=L_0\), there is \(0<\delta_\infty<\delta_1\) such that \(\partial_{x_1}p>0\) on \([1,L_0]\) for \(t_\infty-\delta_\infty<t<t_\infty\).
	For such \(t\), every zero in \((1,\infty)\) lies in \((L_0,\infty)\) and has negative slope, so, by the intermediate value theorem, \(\rho(t)\) is unique, and evenness gives \eqref{0927-005}.

	\underline{Proof of (v).}\quad
	Fix \(t>\overline t\).
	By {\rm(iii)} and \eqref{260156}, \(\partial_{x_1}p(x_1,0,t)\) is negative at \(x_1=1\) and positive for large \(x_1\), so, by continuity,
	\(\varnothing\neq Z_t:=\{x_1>1:\partial_{x_1}p(x_1,0,t)=0\}\subset[1+\eta,R]\) for some \(\eta>0\) and \(R>1+\eta\).
	As a closed set, \(Z_t\) has a minimum \(\rho_-(t)\) and a maximum \(\rho_+(t)\), and evenness gives \eqref{261014}.
	For \(t\geq\widehat t:=T_{\delta_0}\), Lemma~\ref{lem-261080}{\rm(iii)}--{\rm(iv)} gives the unique simple zero \(\rho(t)\) in \((1,\infty)\), the zero \(-\rho(t)\), and \eqref{261015}, so \(\rho_-(t)=\rho_+(t)=\rho(t)\), and also \(\partial_{x_1}p(1,0,t)<0\).
	Since \(\partial_{x_1}p(1,0,\overline t)=0\) by {\rm(iii)}, \(\overline t<\widehat t\).
\end{proof}

\subsubsection{The \(i=j\) case when \(n\geq3\): signs of \(\partial_{x_1}p\)}\label{subsubsect-261020-axis}

In Subsections \ref{subsubsect-261020-axis} and \ref{subsubsect-261022-offaxis}, \(n\geq3\), \(i=j=1\), and \(\psi\) satisfies {\rm\textbf{(S)}} and is radial.
We specify radial monotonicity in the results that use it.

For \(r=|x'|\geq1\), \(s>0\), and \(m=1,2\), define
\begin{align}\label{0925-03}
	\mathcal G^{(m)}(x_1,|x''|,s) & =\int_\Sigma\psi(y')\Bigg[ \left(\frac{x_1^2}{r^2}(r-y_1)^2+ \frac{|y''|^2}{n-2}\left(1-\frac{x_1^2}{r^2}\right)\right) {}_1F_1\left(\frac{n+2}{2},\frac{n+2m+1}{2}; -\frac{|r\mathbf e_1'-y'|^2}{4s}\right)\nonumber \\
	 & \hspace{25mm}-\frac{2(n+2m-1)}n s\, {}_1F_1\left(\frac n2,\frac{n+2m-1}{2}; -\frac{|r\mathbf e_1'-y'|^2}{4s}\right)\Bigg]\,dy'.
\end{align}

The following identities relate \(\mathcal G^{(m)}\) on the two coordinate axes and connect the kernels for \(m=1\) and \(m=2\).
\begin{lemma}\label{lem-260297}
	Let \(n\geq3\), and suppose that \(\psi\) satisfies {\rm\textbf{(S)}} and is radial.
Then the following hold.
\begin{enumerate}
\item For \(r\geq1\), \(s>0\), and \(m=1,2\),
	\begin{align}
		\mathcal G^{(m)}(0,r,s) & =-\frac{2(n+2m-1)s}{nr}\int_\Sigma \psi(y')(r-y_1){}_1F_1\left(\frac n2,\frac{n+2m-1}{2}; -\frac{|r\mathbf e_1'-y'|^2}{4s}\right)\,dy',\label{0911-20c} \\
		\mathcal G^{(m)}(r,0,s) & =-\frac{2(n+2m-1)s}{n}\partial_r \int_\Sigma \psi(y')(r-y_1){}_1F_1\left(\frac n2,\frac{n+2m-1}{2}; -\frac{|r\mathbf e_1'-y'|^2}{4s}\right)\,dy'\nonumber\\
		 & =\partial_r\left(r\mathcal G^{(m)}(0,r,s)\right). \label{260161}
	\end{align}
\item For \(x'\in\Sigma\) with \(|x'|\geq1\) and \(s>0\),
	\begin{align}\label{0911-20}
		\partial_s\left(s^{-\frac{n+3}{2}} \mathcal G^{(2)}(x_1,|x''|,s)\right) =-\frac{n+3}{2}s^{-\frac{n+5}{2}} \mathcal G^{(1)}(x_1,|x''|,s).
	\end{align}
\end{enumerate}
\end{lemma}
\begin{proof}
	\underline{Proof of (i).}\quad
	Differentiate the integral on the left-hand side of \eqref{0911-20e} twice in the directions perpendicular to \(\mathbf e_1'\).
	By \eqref{0327-07},
	\begin{align}\label{0911-20e}
		 & \left.\frac1{n-2}\sum_{j=2}^{n-1}\partial_{x_j}^2 \int_\Sigma\psi(y'){}_1F_1\left(\frac{n-2}{2},\frac{n+2m-3}{2}; -\frac{|x'-y'|^2}{4s}\right)\,dy'\right|_{x'=r\mathbf e_1'} =\frac{n(n-2)\mathcal G^{(m)}(0,r,s)} {4s^2(n+2m-3)(n+2m-1)}.
	\end{align}
	Differentiating the same convolution once in \(r\) gives
	\begin{align}\label{0925-01}
		 & \partial_r\int_\Sigma\psi(y'){}_1F_1\left(\frac{n-2}{2},\frac{n+2m-3}{2}; -\frac{|r\mathbf e_1'-y'|^2}{4s}\right)\,dy'\nonumber\\
		 & \qquad=-\frac{n-2}{2s(n+2m-3)} \int_\Sigma(r-y_1)\psi(y'){}_1F_1\left(\frac n2,\frac{n+2m-1}{2}; -\frac{|r\mathbf e_1'-y'|^2}{4s}\right)\,dy'.
	\end{align}
	For a twice differentiable radial function \(f\) on \(\Sigma\), its second derivatives perpendicular to \(\mathbf e_1'\) satisfy
	\begin{align*}
	    \left.\sum_{j=2}^{n-1}\partial_{x_j}^2f\right|_{x'=r\mathbf e_1'}=(n-2)\frac{f'(r)}r.
	\end{align*}
	Applying this identity to the convolution in \eqref{0911-20e} and substituting \eqref{0925-01} proves \eqref{0911-20c}.

	For the value on the \(x_1\)-axis, \eqref{0327-07} gives
	\begin{align*}
		 & -\frac{2(n+2m-1)s}{n}\partial_r \left((r-y_1){}_1F_1\left(\frac n2,\frac{n+2m-1}{2}; -\frac{|r\mathbf e_1'-y'|^2}{4s}\right)\right)\\
		 & \quad=(r-y_1)^2{}_1F_1\left(\frac{n+2}{2},\frac{n+2m+1}{2}; -\frac{|r\mathbf e_1'-y'|^2}{4s}\right)-\frac{2(n+2m-1)s}{n}{}_1F_1\left(\frac n2,\frac{n+2m-1}{2}; -\frac{|r\mathbf e_1'-y'|^2}{4s}\right).
	\end{align*}
	Integrating against \(\psi\) gives the first equality in \eqref{260161} by \eqref{0925-03}.
	The second follows from \eqref{0911-20c}.

	\underline{Proof of (ii).}\quad 
 Apply \eqref{0911-21} to the two factors containing the confluent hypergeometric functions in \(s^{-(n+3)/2}\mathcal G^{(2)}\), taking \(c=|r\mathbf e_1'-y'|^2/4\) in both.
 The respective parameter pairs are
 \begin{align*}
     (a,b)=\left(\frac{n+2}{2},\frac{n+5}{2}\right),\qquad
     (a,b)=\left(\frac n2,\frac{n+3}{2}\right).
 \end{align*}
 Combining the two differentiated terms gives \eqref{0911-20}.
\end{proof}

Integration by parts in time expresses the tangential pressure derivative in terms of \(\phi'\) and \(\mathcal G^{(2)}\).
\begin{lemma}\label{lem-260159}
		Let \(n\geq 3\), and let \(p\) be the pressure of the Stokes system \eqref{0308-01} with boundary data \eqref{0324-13}.
	Suppose that \(\psi\) satisfies {\rm\textbf{(S)}} and is radial, and that \(\phi\) satisfies {\rm\textbf{(T1)}}.
	Then, for \(|x'| \geq 1\),
	\begin{align}\label{0911-02}
		\partial_{x_1}p(x',0,t) =-\frac{\Gamma\left(\frac{n+2}{2}\right)} {2^{n+2}\pi^{\frac{n}{2}}\Gamma\left(\frac{n+5}{2}\right)} \int_0^t \frac{\phi'(t-s)}{s^{\frac{n+3}{2}}} \mathcal{G}^{(2)}(x_1,|x''|,s)\,ds.
	\end{align}
\end{lemma}
\begin{proof}
	First fix \(t>0\) and \(r=|x'|>1\).
	Radial symmetry allows us to replace \(y_2^2\) in the integrals in \eqref{0825-13} by the average \(|y''|^2/(n-2)\).
	Thus \eqref{0925-03} gives
	\begin{align*}
		\partial_{x_1}p(x',0,t) =-\frac{\Gamma\left(\frac{n+2}{2}\right)} {2^{n+2}\pi^{\frac n2}\Gamma\left(\frac{n+3}{2}\right)} \int_0^\infty \frac{\phi(t)-\phi(t-s)}{s^{\frac{n+5}{2}}} \mathcal{G}^{(1)}(x_1,|x''|,s)\,ds.
	\end{align*}
	On \(\supp\psi\), we have \(|r\mathbf e_1'-y'|\geq r-1>0\).
	Hence \eqref{0514-02} for small \(s\) and continuity of the factors involving the confluent hypergeometric functions at zero for large \(s\) give
	\begin{align*}
		\mathcal{G}^{(2)}(x_1,|x''|,s) =O\bigl(s^{\frac{n+2}{2}}\bigr)\quad(s\to0^+), \qquad \mathcal{G}^{(2)}(x_1,|x''|,s)=O(s)\quad(s\to\infty).
	\end{align*}
	Extend \(\phi\) by zero to negative times.
	Since \(\phi(t)-\phi(t-s)=O(s)\) at \(s=0\) and is bounded at infinity, we obtain
	\begin{align*}
		\bigl(\phi(t)-\phi(t-s)\bigr)s^{-\frac{n+3}{2}} \mathcal{G}^{(2)}(x_1,|x''|,s) =
		\begin{cases}
			O(s^{\frac{1}{2}}),\quad s\to0^+, \\
			O(s^{-\frac{n+1}{2}}),\quad s\to\infty.
		\end{cases}
	\end{align*}
	The endpoint terms at \(s=0\) and \(s=\infty\) therefore vanish.
	Since \(\phi(0)=0\), the map \(s\mapsto\phi(t)-\phi(t-s)\) is continuous at \(s=t\), so the two terms there cancel.
	Using \eqref{0911-20} for the differentiated kernel gives
	\begin{align*}
		\int_0^\infty\frac{\phi(t)-\phi(t-s)}{s^{\frac{n+5}{2}}} \mathcal G^{(1)}(x_1,|x''|,s)\,ds =\frac{2}{n+3}\int_0^t \frac{\phi'(t-s)}{s^{\frac{n+3}{2}}} \mathcal G^{(2)}(x_1,|x''|,s)\,ds.
	\end{align*}
	This proves \eqref{0911-02} for \(r>1\).

	To pass to \(r=1\), write \(x'=r\omega'\) with \(\omega'\in\mathbb S^{n-2}\) fixed.
	The quadratic factor in \eqref{0925-03} is bounded by \(|r\mathbf e_1'-y'|^2\), and \(|r\mathbf e_1'-y'|\geq|\mathbf e_1'-y'|\) on \(\supp\psi\).
	Consequently, \eqref{0514-02} bounds the absolute value of the integrand in \eqref{0911-02}, before integrating in \(y'\), by
	\begin{align*}
	    C\|\phi'\|_{L^\infty(0,1)} s^{-1/2}|\psi(y')|\,|\mathbf e_1'-y'|^{-n},
	    \qquad 0<s<t,\quad 1<r\leq2.
	\end{align*}
	By {\rm\textbf{(S)}}, \(|\psi(y')|\leq C|\mathbf e_1'-y'|^{2+\alpha}\), so this majorant is integrable in \((y',s)\), and dominated convergence gives the right-hand side of \eqref{0911-02} at \(r=1\).
	By continuity of \(x'\mapsto\partial_{x_1}p(x',0,t)\) in \eqref{tangential-pressure-boundary-value}, the left-hand side tends to \(\partial_{x_1}p(\omega',0,t)\).
\end{proof}

At a fixed radius, the values of \(\partial_{x_1}p\) on the \(x_1\)- and \(x_2\)-axes determine its value in every direction.

\begin{lemma}\label{lem-0911-02}
	Suppose that \(\psi\) satisfies {\rm\textbf{(S)}} and is radial, and that \(\phi\) satisfies {\rm\textbf{(T1)}}.
	Then, for every \(r=|x'|\geq1\) and \(t>0\),
	\begin{align}\label{0911-18}
		\partial_{x_1}p(x',0,t) & =\frac{x_1^2}{r^2}\partial_{x_1}p(r\mathbf e_1',0,t) +\frac{|x''|^2}{r^2}\partial_{x_1}p(r\mathbf e_2',0,t).
	\end{align}
	In particular, \(\partial_{x_1}p(-x_1,x'',0,t) =\partial_{x_1}p(x',0,t)\).
	Moreover,
	\begin{align}\label{260160}
	    \partial_{x_1}p(\cdot,0,\cdot)\in C\bigl(\{x'\in\Sigma:|x'|\geq1\}\times(0,\infty)\bigr).
	\end{align}
\end{lemma}
\begin{proof}
	Equation \eqref{0925-03} with \(m=2\) and \(x_1^2+|x''|^2=r^2\) gives
	\begin{align*}
		\mathcal G^{(2)}(x_1,|x''|,s) =\frac{x_1^2}{r^2}\mathcal G^{(2)}(r,0,s) +\frac{|x''|^2}{r^2}\mathcal G^{(2)}(0,r,s).
	\end{align*}
	Substituting this identity into \eqref{0911-02} gives \eqref{0911-18} for \(r>1\).
	Formula \eqref{tangential-pressure-boundary-value} gives joint continuity of \((x',t)\mapsto\partial_{x_1}p(x',0,t)\), proving \eqref{260160}.
	Letting \(r\downarrow1\) along \(x'=r\omega'\) proves \eqref{0911-18} also at every \(\omega'\in\mathbb S^{n-2}\).

	We also record a bound used in the proof of the next lemma to pass to \(r=1\).
	Fix \(1\leq k\leq n-1\) and let \(x'=r\mathbf e_k'\).
	For \(1\leq r\leq2\) and \(|y'|\leq1\),
	\begin{align*}
	    |r\mathbf e_k'-y'|\geq|\mathbf e_k'-y'|.
	\end{align*}
	Fix \(y'\ne\mathbf e_k'\).
	In \eqref{0822-13}, apply \(|\phi(t)-\phi(t-s)|\leq C\min\{s,1\}\) and \eqref{0514-02}, then set \(u=|x'-y'|^2/s\).
	This gives
	\begin{align}\label{260284}
		 & \int_0^\infty|\phi(t)-\phi(t-s)| \biggl[ \frac{|x_1-y_1|^2}{s^{\frac{n+5}{2}}} \left|{}_1F_1\left(\frac{n+2}{2},\frac{n+3}{2}; -\frac{|x'-y'|^2}{4s}\right)\right| +\frac1{s^{\frac{n+3}{2}}} \left|{}_1F_1\left(\frac n2,\frac{n+1}{2}; -\frac{|x'-y'|^2}{4s}\right)\right|\biggr]\,ds\nonumber \\
		 & \quad\leq C\int_0^\infty \left[\frac{|x'-y'|^2}{s^{(n+3)/2}(1+|x'-y'|^2/s)^{(n+2)/2}}
		 +\frac{1}{s^{(n+1)/2}(1+|x'-y'|^2/s)^{n/2}}\right]ds\nonumber\\
		 & \quad=\frac{C}{|x'-y'|^{n-1}}\int_0^\infty
		 \left[\frac{u^{(n-1)/2}}{(1+u)^{(n+2)/2}}
		 +\frac{u^{(n-3)/2}}{(1+u)^{n/2}}\right]du
		 \leq\frac{C}{|\mathbf e_k'-y'|^{n-1}}.
	\end{align}
	The \(u\)-integral is finite for \(n\geq3\), and the final inequality uses \(|x'-y'|\geq|\mathbf e_k'-y'|\).

	Since \(\psi\in C_c^{2,\alpha}(\Sigma)\) and vanishes outside \(\overline{B_1'}\), Taylor's theorem gives
	\begin{align}\label{260285}
		|\psi(y')|\,|\mathbf e_k'-y'|^{1-n} & \leq C|\mathbf e_k'-y'|^{3+\alpha-n},\nonumber \\
		\int_{\Sigma\cap B_2'(\mathbf e_k')} |\mathbf e_k'-y'|^{3+\alpha-n}\,dy' & \leq C\int_0^2 u^{1+\alpha}\,du<\infty.
	\end{align}
	Thus the bound \eqref{260284}, multiplied by \(|\psi(y')|\), is integrable in \(y'\) uniformly in \(1\leq r\leq2\) and \(t>0\).
\end{proof}

We now express the pressure derivative on the \(x_1\)-axis through a radial derivative of its value on the \(x_2\)-axis.
\begin{lemma}\label{lem-0911-03}
	Suppose that \(\psi\) satisfies {\rm\textbf{(S)}} and is radial, and that \(\phi\) satisfies {\rm\textbf{(T1)}}.
	Then the following hold.
	\begin{enumerate}
		\item For \(r\geq1\) and \(t>0\),
		\begin{align}\label{0911-22f}
			\partial_{x_1}p(r\mathbf e_1',0,t) & =\frac{\Gamma\left(\frac{n+2}{2}\right)} {n\,2^n\pi^{\frac{n}{2}}\Gamma\left(\frac{n+3}{2}\right)} \int_\Sigma\partial_{y_1}\psi(y')(r-y_1)\int_0^t \frac{\phi'(t-s)}{s^{\frac{n+1}{2}}} {}_1F_1\left(\frac n2,\frac{n+3}{2};- \frac{|r\mathbf e_1'-y'|^2}{4s}\right)\,dsdy'.
		\end{align}
		For \(r>1\) and \(t>0\),
		\begin{align}\label{0911-20ab}
			\partial_{x_1}p(r\mathbf e_2',0,t) & =\frac{\Gamma\left(\frac n2\right)} {2^{n+1}\pi^{\frac n2}\Gamma\left(\frac{n+3}{2}\right)r} \int_\Sigma (r-y_1)\psi(y')\int_0^t \frac{\phi'(t-s)}{s^{\frac{n+1}{2}}} {}_1F_1\left(\frac n2,\frac{n+3}{2};- \frac{|r\mathbf e_1'-y'|^2}{4s}\right)\,dsdy'.
		\end{align}
		\item For \(r>1\) and \(t>0\),
		\begin{align}\label{0911-20f}
			\partial_{x_1}p(r\mathbf e_1',0,t) =\partial_r\left( r\partial_{x_1}p(r\mathbf e_2',0,t)\right).
		\end{align}
		Moreover, if \(r=|x'|>1\), then
		\begin{align}\label{0911-20g}
			\partial_{x_1}p(x',0,t) & =\partial_{x_1}p(r\mathbf e_2',0,t) +\frac{x_1^2}{r}\partial_r \partial_{x_1}p(r\mathbf e_2',0,t).
		\end{align}
	\end{enumerate}
\end{lemma}
\begin{proof}
	\underline{Proof of (i).}\quad
	The identity \eqref{0327-07} shows that the expression in brackets in \(\mathcal G^{(2)}(r,0,s)\) equals
	\begin{align*}
		\frac{2(n+3)s}{n}\partial_{y_1} \left((r-y_1){}_1F_1\left(\frac n2,\frac{n+3}{2};- \frac{|r\mathbf e_1'-y'|^2}{4s}\right)\right).
	\end{align*}
	Substituting this identity into \eqref{0911-02} and integrating by parts in \(y_1\) prove \eqref{0911-22f} for \(r>1\).
	For \(1\leq r\leq2\), \eqref{0514-02} and the change of variables \(s=|r\mathbf e_1'-y'|^2u\) give
	\begin{align*}
		|\partial_{y_1}\psi(y')|\,|r-y_1| \int_0^t\frac{|\phi'(t-s)|}{s^{\frac{n+1}{2}}} \left|{}_1F_1\left(\frac n2,\frac{n+3}{2};- \frac{|r\mathbf e_1'-y'|^2}{4s}\right)\right|\,ds \leq C|\mathbf e_1'-y'|^{3+\alpha-n}\mathbf1_{\overline{B_1'}}(y').
	\end{align*}
	By \eqref{260285}, \eqref{260160}, and dominated convergence, \eqref{0911-22f} extends to \(r=1\).

	For the value on the \(x_2\)-axis, substituting \eqref{0911-20c} with \(m=2\) into \eqref{0911-02} gives \eqref{0911-20ab}.

	\underline{Proof of (ii).}\quad
	Differentiating \eqref{0911-02} gives
	\begin{align}\label{260282}
		\partial_r\left(r\partial_{x_1}p(r\mathbf e_2',0,t)\right) & =-\frac{\Gamma(\frac{n+2}{2})} {2^{n+2}\pi^{\frac n2}\Gamma(\frac{n+5}{2})} \int_0^t\frac{\phi'(t-s)}{s^{\frac{n+3}{2}}} \partial_r\left(r\mathcal G^{(2)}(0,r,s)\right)\,ds.
	\end{align}
	The differentiation is justified since, by \eqref{0514-02} and \eqref{260161}, the integrand is bounded by \(Cs^{-1/2}\) near \(s=0\), locally uniformly in \(r>1\).
	Now \eqref{260161} identifies the differentiated kernel with \(\mathcal G^{(2)}(r,0,s)\), so \eqref{0911-02} gives \eqref{0911-20f}.

	Substituting \eqref{0911-20f} into \eqref{0911-18}, we obtain
	\begin{align*}
		\partial_{x_1}p(x',0,t) =\partial_{x_1}p(r\mathbf e_2',0,t) +\frac{x_1^2}{r^2} \left[\partial_{x_1}p(r\mathbf e_1',0,t) -\partial_{x_1}p(r\mathbf e_2',0,t)\right] =\partial_{x_1}p(r\mathbf e_2',0,t) +\frac{x_1^2}{r}\partial_r \partial_{x_1}p(r\mathbf e_2',0,t),
	\end{align*}
	which is \eqref{0911-20g}.

\end{proof}

For the large-\(r\) expansions on the two axes, define
\begin{equation}\label{260279}\begin{gathered}
		A(t)=M_0T_0(t),\qquad D(t)=\frac{nM_2}{2(n-1)}T_0(t)-nM_0T_1(t),\\
		E(t)=M_0\left(|T_0(t)|+|T_1(t)|+ \int_0^t\phi(t-s)s^{\frac{1}{2}}\,ds\right),\qquad R_0(t)=2+C_n\sqrt t.
	\end{gathered}\end{equation}
Here \(z_0(a,b)\) is the lower bound on the argument in \eqref{0504-06}, and
\begin{align*}
    C_n=\max\left\{4,\ 2\max_{0\leq j\leq3}\sqrt{z_0\left(\frac n2+j,\frac{n+3}{2}+j\right)}\right\}.
\end{align*}
Below, \(C(n)\) denotes a constant depending only on \(n\).

\begin{lemma}\label{lem-0911-04a}
	Suppose that \(\psi\) satisfies {\rm\textbf{(S)}} and is radial, and that \(\phi\) satisfies {\rm\textbf{(T1)}}.
	Then, as \(r\to\infty\), uniformly for \(t\) in compact subsets of \((0,\infty)\),
	\begin{align}
		\partial_{x_1}p(r\mathbf e_2',0,t) & =\frac{\Gamma(\frac{n+2}{2})}{n\pi^{\frac{n}{2}}\Gamma(\frac{3}{2})r^n} \bigl[A(t)+\frac{D(t)}{r^2}+O(r^{-4})\bigr],\label{0911-20w}\\
		\partial_{x_1}p(r\mathbf e_1',0,t) & =-\frac{\Gamma(\frac{n+2}{2})}{n\pi^{\frac{n}{2}}\Gamma(\frac{3}{2})r^n}\bigl[(n-1)A(t)+\frac{(n+1)D(t)}{r^2}+O(r^{-4})\bigr],\label{0911-56b} \\
		\partial_r\partial_{x_1}p(r\mathbf e_1',0,t) & =\frac{\Gamma(\frac{n+2}{2})} {n\pi^{\frac n2}\Gamma(\frac32)r^{n+1}}\left[n(n-1)A(t) +\frac{(n+1)(n+2)D(t)}{r^2}+O(r^{-4})\right].\label{1013-106}
	\end{align}
	In particular, for every \(t>0\),
	\begin{align}
		\lim_{r\to\infty}r^n\partial_{x_1}p(r\mathbf e_2',0,t) & =\frac{\Gamma(\frac{n+2}{2})A(t)}{n\pi^{\frac{n}{2}}\Gamma(\frac{3}{2})}.\label{0911-20m}
	\end{align}
	For \(r\geq R_0(t)\), the \(j\)-th derivative of the \(O(r^{-4})\) term inside the brackets in \eqref{0911-20w} is bounded by \(C(n)E(t)r^{-4-j}\) for \(j=0,1,2\).
	The same bound holds for the \(O(r^{-4})\) term in \eqref{0911-56b} for \(j=0,1\).
\end{lemma}
\begin{proof}
	For \(r\geq R_0(t)\), the argument of each kernel in \eqref{0911-20ab} satisfies
	\begin{align*}
		    \frac{|r\mathbf e_1'-y'|^2}{4s}\geq\frac{(r-1)^2}{4t}
		    \geq\frac{C_n^2}{4}\geq\max_{0\leq j\leq3}z_0\left(\frac n2+j,\frac{n+3}{2}+j\right),
		    \qquad 0<s<t,\quad y'\in\supp\psi.
	\end{align*}
	Thus \eqref{0504-06} applies also to the shifted parameters in the derivative formula \eqref{0327-07}.
	Consequently, for \(k=0,1,2\),
	\begin{align}\label{260166}
		 & \left|\frac{d^k}{dr^k}\left[ \frac{1}{s^{\frac{n+1}{2}}}{}_1F_1\left(\frac n2,\frac{n+3}{2};- \frac{|r\mathbf e_1'-y'|^2}{4s}\right)- \frac{2^n\Gamma\left(\frac{n+3}{2}\right)} {\Gamma\left(\frac32\right)} \left(\frac{1}{\sqrt{s}|r\mathbf e_1'-y'|^n} -\frac{n\sqrt{s}}{|r\mathbf e_1'-y'|^{n+2}}\right)\right]\right| \leq \frac{C(n)s^{\frac{3}{2}}}{|r\mathbf e_1'-y'|^{n+4+k}},\nonumber \\
		 & \left|\partial_s\frac{d^k}{dr^k}\left[ \frac{1}{s^{\frac{n+1}{2}}}{}_1F_1\left(\frac n2,\frac{n+3}{2};- \frac{|r\mathbf e_1'-y'|^2}{4s}\right)- \frac{2^n\Gamma\left(\frac{n+3}{2}\right)} {\Gamma\left(\frac32\right)} \left(\frac{1}{\sqrt{s}|r\mathbf e_1'-y'|^n} -\frac{n\sqrt{s}}{|r\mathbf e_1'-y'|^{n+2}}\right)\right]\right| \leq \frac{C(n)\sqrt{s}}{|r\mathbf e_1'-y'|^{n+4+k}}.
	\end{align}

	In the integration by parts against \(\phi'(t-s)\), the boundary term at \(s=0\) vanishes by the first bound in \eqref{260166}, and the one at \(s=t\) vanishes because \(\phi(0)=0\).
	The second bound in \eqref{260166} then gives
	\begin{align}\label{260167}
		 & \left|\partial_r^k\int_0^t\phi'(t-s) \left[\frac{1}{s^{\frac{n+1}{2}}} {}_1F_1\left(\frac n2,\frac{n+3}{2};-\frac{|r\mathbf e_1'-y'|^2}{4s}\right) -\frac{2^n\Gamma(\frac{n+3}{2})}{\Gamma(\frac32)} \left(\frac{1}{\sqrt{s}|r\mathbf e_1'-y'|^n} -\frac{n\sqrt{s}}{|r\mathbf e_1'-y'|^{n+2}}\right)\right]ds\right|\nonumber \\
		 & \qquad\leq \frac{C(n)}{|r\mathbf e_1'-y'|^{n+4+k}} \int_0^t\phi(t-s)s^{\frac12}\,ds.
	\end{align}

	We next expand the spatial integrals in \eqref{260168} for \(r\geq2\).
	Radial symmetry cancels the odd moments and gives \(\int_\Sigma y_1^2\psi(y')\,dy'=M_2/(n-1)\).
	A Taylor expansion in \(y'/r\), also after \(k=0,1,2\) radial derivatives, therefore gives
	\begin{align}\label{260168}
		 & \left|\frac{d^k}{dr^k}\left[ \frac1r\int_\Sigma\frac{(r-y_1)\psi(y')} {|r\mathbf e_1'-y'|^n}\,dy'-\frac{M_0}{r^n} -\frac{nM_2}{2(n-1)r^{n+2}}\right]\right| \leq \frac{C(n)M_0}{r^{n+4+k}},\nonumber \\
		 & \left|\frac{d^k}{dr^k}\left[ \frac1r\int_\Sigma\frac{(r-y_1)\psi(y')} {|r\mathbf e_1'-y'|^{n+2}}\,dy'-\frac{M_0}{r^{n+2}}\right]\right| \leq \frac{C(n)M_0}{r^{n+4+k}},
	\end{align}
	where the remainder uses \(|y'|\leq1\) on \(\supp\psi\).
	Multiplying the two estimates in \eqref{260168} by the respective moment coefficients from \eqref{0323-02} yields
	\begin{align*}
		 & \left|\partial_r^k\left[\frac1r\int_\Sigma(r-y_1)\psi(y') \left(\frac{T_0(t)}{|r\mathbf e_1'-y'|^n} -\frac{nT_1(t)}{|r\mathbf e_1'-y'|^{n+2}}\right)\,dy' -\frac{M_0T_0(t)}{r^n}-\frac{D(t)}{r^{n+2}}\right]\right| \leq\frac{C(n)M_0}{r^{n+4+k}} \bigl(|T_0(t)|+|T_1(t)|\bigr).
	\end{align*}
	We now substitute \eqref{260167} and \eqref{260168} into \eqref{0911-20ab}.
	For the pressure derivative on the \(x_2\)-axis and its first two radial derivatives, this gives
	\begin{align}
		\partial_{x_1}p(r\mathbf e_2',0,t)& =\frac{\kappa_n}{r^n}\left(A(t)+\frac{D(t)}{r^2}+e_0(r,t)\right),\label{0911-60}\\
		\partial_r\partial_{x_1}p(r\mathbf e_2',0,t)& =\frac{\kappa_n}{r^{n+1}}\left(-nA(t)-\frac{(n+2)D(t)}{r^2}+e_1(r,t)\right),\label{0911-601} \\
		\partial_r^2\partial_{x_1}p(r\mathbf e_2',0,t) & =\frac{\kappa_n}{r^{n+2}}\left(n(n+1)A(t) +\frac{(n+2)(n+3)D(t)}{r^2}+e_2(r,t)\right),\label{0911-602}
	\end{align}
	where \(\kappa_n=\frac{\Gamma(\frac{n+2}{2})}{n\pi^{\frac{n}{2}}\Gamma(\frac32)}\), and
	\begin{align}\label{0911-60c}
		|e_k(r,t)|\leq \frac{C(n)E(t)}{r^4},\qquad k=0,1,2.
	\end{align}
	Equations \eqref{0911-60} and \eqref{0911-60c} give \eqref{0911-20w} and its limit \eqref{0911-20m}.
	Comparing the derivatives of \eqref{0911-60} with \eqref{0911-601}--\eqref{0911-602} yields
	\begin{align*}
		r\partial_re_0=e_1+ne_0,\qquad r^2\partial_r^2e_0=e_2-n(n+1)e_0+2nr\partial_re_0,
	\end{align*}
	so \eqref{0911-60c} controls the first two radial derivatives of \(e_0\).

	By \eqref{0911-20f}, the remainder inside the brackets in \eqref{0911-56b} is \((n-1)e_0-r\partial_re_0\).
	Its radial derivative is therefore bounded in absolute value by \(C(n)E(t)r^{-5}\).
	Differentiating \eqref{0911-56b} with this bound gives
	\begin{align*}
		\partial_r\partial_{x_1}p(r\mathbf e_1',0,t) & =\frac{\kappa_n}{r^{n+1}}\left[n(n-1)A(t) +\frac{(n+1)(n+2)D(t)}{r^2}+O(r^{-4})\right],
	\end{align*}
	as asserted in \eqref{1013-106}.
\end{proof}
The next proposition gives the signs and zeros on the \(x_1\)-axis used in Theorem \ref{thm-260021} when \(n\geq3\).

\begin{proposition}\label{prop-261020}
	Let \(n\geq3\), and let \(p\) be the pressure of the Stokes system \eqref{0308-01} with boundary data \eqref{0324-13}, where \(i=j=1\).
	Suppose that \(\psi\) satisfies {\rm\textbf{(S)}} and is radial, and that \(\phi\) satisfies {\rm\textbf{(T1)}} and {\rm\textbf{(T2)}}.
	Then the following hold.
	\begin{enumerate}
		\item Suppose that \(\phi\) satisfies {\rm\textbf{(T3)}} and {\rm\textbf{(T4)}}. Then, for every fixed \(t>0\) with \(t\ne t_\infty\), there exists \(R(t)>1\) such that
		\begin{align}\label{261041}
			\partial_{x_1}p(r\mathbf e_1',0,t)
			\begin{cases}
				<0, & 0<t<t_\infty, \\
				>0, & t>t_\infty,
			\end{cases}
			\qquad |r|>R(t).
		\end{align}
		If \(t\in(0,\infty)\setminus\{t_\infty\}\) and \(\varepsilon\in(0,\sqrt{n-1})\), then there exists \(R_\varepsilon(t)>1\) such that
		\begin{equation}\label{261047}
			\sgn(T_0(t))\partial_{x_1}p(x',0,t)
			\begin{cases}
				>0, & |x''|\geq(\sqrt{n-1}+\varepsilon)|x_1|, \\
				<0, & |x''|\leq(\sqrt{n-1}-\varepsilon)|x_1|,
			\end{cases}
			\qquad |x'|>R_\varepsilon(t).
		\end{equation}
		\item There exists \(t_0=t_0(n,\psi,\phi)\in(0,t_1)\) such that \(\partial_{x_1}p(r\mathbf e_1',0,t)<0\) for every \(|r|\geq1\) and \(0<t<t_0\).
		\item Suppose that \(\psi\) is radially nonincreasing and that \(\phi\) satisfies {\rm\textbf{(T3)}} and {\rm\textbf{(T4)}}. If \(\varepsilon\in(0,\sqrt{n+1})\), then there exists \(R_\varepsilon(t_\infty)>1\) such that
		\begin{equation}\label{261049}
			\partial_{x_1}p(x',0,t_\infty)
			\begin{cases}
				<0, & |x''|\geq(\sqrt{n+1}+\varepsilon)|x_1|, \\
				>0, & |x''|\leq(\sqrt{n+1}-\varepsilon)|x_1|,
			\end{cases}
			\qquad |x'|>R_\varepsilon(t_\infty).
		\end{equation}
		Suppose, in addition, that \(\phi\) satisfies {\rm\textbf{(T5)}}.
		Then \(\overline t\) in \eqref{0911-23} is well-defined and satisfies
		\begin{align*}
			t_\infty<\overline t<\infty.
		\end{align*}
		Moreover,
		\begin{align}\label{260286}
			\partial_{x_1}p(\pm\mathbf e_1',0,t)<0,\qquad t>\overline t,
		\end{align}
		and \(\partial_{x_1}p(\pm\mathbf e_1',0,\overline t)=0\).
		At \(t=t_\infty\),
		\begin{align}\label{0911-41}
			\partial_{x_1}p(r\mathbf e_1',0,t_\infty)>0,\qquad |r|\geq1.
		\end{align}
		\item Suppose that \(\psi\) is radially nonincreasing and that \(\phi\) satisfies {\rm\textbf{(T3)}}--{\rm\textbf{(T5)}}.
		Then there exists \(\delta_\infty=\delta_\infty(n,\psi,\phi)\in(0,t_\infty-t_1)\) such that, for \(t_\infty-\delta_\infty<t<t_\infty\), \(r\mapsto\partial_{x_1}p(r\mathbf e_1',0,t)\) has exactly one zero \(\rho(t)\) on \((1,\infty)\), and \(-\rho(t)\) is the only zero on \((-\infty,-1)\).
		Moreover,
		\begin{align}\label{1013-107}
			\sgn\bigl(\partial_{x_1}p(r\mathbf e_1',0,t)\bigr) =\sgn(\rho(t)-|r|),\qquad |r|\geq1,
		\end{align}
		and
		\begin{align}\label{261042}
			\lim_{t\to t_\infty^-}(t_\infty-t)\rho(t)^2 =-\frac{n(n+1)T_1(t_\infty)} {(n-1)T_0'(t_\infty)}.
		\end{align}
		\item Suppose that \(\phi\) satisfies {\rm\textbf{(T5)}}.
		Then, for every \(t>\overline t\), the zeros of \(r\mapsto\partial_{x_1}p(r\mathbf e_1',0,t)\) on \((1,\infty)\) have a minimum \(\rho_-(t)\) and a maximum \(\rho_+(t)\), which satisfy
		\begin{equation}\label{1013-108}
			\begin{aligned}
				\partial_{x_1}p(r\mathbf e_1',0,t) & <0 && \text{for }1\leq|r|<\rho_-(t), \\
				\partial_{x_1}p(r\mathbf e_1',0,t) & >0 && \text{for }|r|>\rho_+(t).
			\end{aligned}
		\end{equation}
		Moreover, there exists \(\widehat t>1\), depending only on \(n\), such that \(\overline t<\widehat t\) and, for every \(t\geq\widehat t\), \(\rho(t):=\rho_-(t)=\rho_+(t)\) satisfies \(\partial_r\partial_{x_1}p(r\mathbf e_1',0,t)>0\) at \(r=\rho(t)\) and
		\begin{align}\label{261034}
			\sgn\bigl(\partial_{x_1}p(r\mathbf e_1',0,t)\bigr)=\sgn\bigl(|r|-\rho(t)\bigr)\quad(|r|\geq1), \qquad \lim_{t\to\infty}\frac{\rho(t)}{\sqrt t}=2\sqrt{z_7}.
		\end{align}
	\end{enumerate}
\end{proposition}
\begin{proof}
	For \(r=|x'|\), substituting \eqref{0911-60}--\eqref{0911-601} into \eqref{0911-20g} and bounding the terms containing \(e_0\) and \(e_1\) by \eqref{0911-60c} give, uniformly in \(x'/r\),
	\begin{align}\label{261048}
		\partial_{x_1}p(x',0,t) =\frac{\kappa_n}{r^n}\left[ A(t)\frac{|x''|^2-(n-1)x_1^2}{r^2} +\frac{D(t)}{r^2}\frac{|x''|^2-(n+1)x_1^2}{r^2} +O_t(r^{-4})\right].
	\end{align}

	\underline{Proof of (i).}\quad
	Fix \(t\in(0,\infty)\setminus\{t_\infty\}\), so that \(A(t)=M_0T_0(t)\neq0\).
	By \eqref{0911-56b} and \eqref{0911-60c}, there is \(R(t)\geq R_0(t)\) such that
	\begin{align}\label{261061}
		\frac{(n+1)|D(t)|}{r^2}+\frac{C(n)E(t)}{r^4} \leq\frac{n-1}{2}|A(t)|, \qquad r>R(t),
	\end{align}
	so the term \((n-1)A(t)\) determines the sign of the bracket in \eqref{0911-56b}, and \(\sgn(\partial_{x_1}p(r\mathbf e_1',0,t)) =-\sgn(T_0(t))\) for \(r>R(t)\).
	Evenness in \(r\) gives \eqref{261041}.
	On the two regions in \eqref{261047}, \((|x''|^2-(n-1)x_1^2)/r^2\) has a fixed sign and is bounded away from zero, so the \(A(t)\)-term in \eqref{261048} determines the sign for \(r>R_\varepsilon(t)\).

	\underline{Proof of (ii).}\quad
	Let \(0<t<t_1\).
	By {\rm\textbf{(S)}}, \(\partial_{x_1}^2\psi\) vanishes on \(\{|x'|\geq1\}\), so \eqref{tangential-pressure-boundary-value} gives
	\begin{align}\label{1013-101}
		\partial_{x_1}p(x',0,t) =\frac1{2\sqrt\pi}\int_0^\infty \frac{\phi(t)-\phi(t-s)}{s^{\frac32}} \bigl(\Gamma'(\cdot,s)*_\Sigma\mathcal S\bigr)(x')\,ds,
		\qquad |x'|\geq1,
	\end{align}
	where \(\mathcal S = R_1'\left(\partial_{x_1}\psi\right)\) and \(\phi(t-s)=0\) for \(s\geq t\).
	By \eqref{0323-02}, \(\int_0^\infty(\phi(t)-\phi(t-s))s^{-3/2}\,ds =2T_0(t)\).
	Lemma~\ref{lem-260265} gives \(\mathcal S\in C^{1,\alpha}(\Sigma)\cap C_0(\Sigma)\), so subtracting its linear Taylor polynomial in \(\Gamma'(\cdot,s)*_\Sigma\mathcal S\) gives \(\|\Gamma'(\cdot,s)*_\Sigma\mathcal S-\mathcal S\|_\infty \leq C_{n,\alpha,\psi}\min\{s^{(1+\alpha)/2},1\}\).
	Hence the integral below converges uniformly, and, since its Fourier multiplier is \(-|\xi'|\) by integration by parts in \(s\), and \(-|\xi'|\widehat{\mathcal S}=-\xi_1^2\widehat\psi\), Plancherel's theorem gives
	\begin{align*}
		\frac1{2\sqrt\pi}\int_0^\infty \frac{\Gamma'(\cdot,s)*_\Sigma\mathcal S-\mathcal S}{s^{3/2}}\,ds =\partial_{x_1}^2\psi,
	\end{align*}
	which vanishes on \(\{|x'|\geq1\}\).
	Thus, subtracting \(T_0(t)\mathcal S/\sqrt\pi\) in \eqref{1013-101},
	\begin{align*}
		\partial_{x_1}p(x',0,t)-\frac{T_0(t)}{\sqrt\pi}\mathcal S(x') =-\frac1{2\sqrt\pi}\int_0^t\frac{\phi(t-s)}{s^{3/2}} [\Gamma'(\cdot,s)*_\Sigma\mathcal S-\mathcal S](x')\,ds,
	\end{align*}
	and, since \(\phi(t-s)\leq\phi(t)\leq\sqrt t\,T_0(t)\) on \(0<s<t\),
	\begin{align}\label{1013-102a}
		\sup_{|x'|\geq1}\left| \frac{\partial_{x_1}p(x',0,t)}{T_0(t)}-\frac{\mathcal S(x')}{\sqrt\pi} \right|\leq C_{n,\alpha,\psi}\,t^{(1+\alpha)/2}.
	\end{align}

	For \(0<t<t_1\), positivity of \(\phi'\) on \((0,t)\) gives
	\begin{align}\label{1013-103}
		0<T_1(t)\leq tT_0(t),\qquad 0\leq\int_0^t\phi(t-s)s^{\frac12}\,ds =\frac23\int_0^t\phi'(t-s)s^{\frac32}\,ds \leq\frac23t^2T_0(t),
	\end{align}
	so, since \(M_2\leq M_0\), \eqref{260279} gives \(|D(t)|+E(t)\leq C(n)M_0T_0(t)\).
	Hence there is \(R\geq2+C_n\), depending only on \(n\), such that \eqref{0911-56b} gives
	\begin{align}\label{1013-105}
		\partial_{x_1}p(r\mathbf e_1',0,t) \leq-\frac{(n-1)\Gamma(\frac{n+2}{2})M_0T_0(t)} {2n\pi^{\frac n2}\Gamma(\frac32)r^n}<0, \qquad r\geq R,\quad 0<t<t_1.
	\end{align}
	For \(1<r\leq R\), let \(u\) be the radial Riesz potential from Proposition~\ref{prop-260061}{\rm(iv)}.
	Then \(u'(r)<0\), \(\mathcal S(r\mathbf e_1')=-u''(r)\), and \eqref{riesz-laplacian} gives
	\begin{align*}
		-\mathcal S(r\mathbf e_1')=u''(r) \geq c_n(n-2)\int_\Sigma\frac{\psi(y')}{|r\mathbf e_1'-y'|^n}\,dy' \geq \frac{c_n(n-2)M_0}{(R+1)^n}=:m_R>0,
	\end{align*}
	also at \(r=1\) by continuity.
	Choose \(t_0\in(0,t_1)\) so that \(C_{n,\alpha,\psi}t_0^{(1+\alpha)/2}\leq m_R/(2\sqrt\pi)\).
	Then \eqref{1013-102a} gives \(\partial_{x_1}p(r\mathbf e_1',0,t)\leq-\frac{m_R}{2\sqrt\pi}T_0(t)<0\) for \(1\leq r\leq R\) and \(0<t<t_0\), and \eqref{1013-105} and evenness in \(r\) prove {\rm(ii)}.

	\underline{Proof of (iii).}\quad
	At \(r=1\), \(z=-|\mathbf e_1'-y'|^2/(4s)\) satisfies \(|z|\leq s^{-1}\) on \(\supp\psi\), so, for both parameter pairs in \eqref{0925-03}, \eqref{260030} gives \({}_1F_1(a,b;z)=1+O_n(s^{-1})\), and \eqref{0327-07} and \(|\partial_s z|\leq s^{-2}\) bound its \(s\)-derivative by \(O_n(s^{-2})\), uniformly in \(y'\).
	Hence, as \(s\to\infty\),
	\begin{align*}
		\mathcal G^{(2)}(1,0,s)=-\frac{2(n+3)}nM_0s+O_n(M_0),\qquad \partial_s\mathcal G^{(2)}(1,0,s)=-\frac{2(n+3)}nM_0+O_n(M_0s^{-1}),
	\end{align*}
	and \(\partial_s[s^{-\frac{n+3}{2}}\mathcal G^{(2)}(1,0,s)] =\frac{(n+3)(n+1)}nM_0s^{-\frac{n+3}{2}} +O(s^{-\frac{n+5}{2}})\).
	For large \(t\), this is positive at \(s=t-u\) for \(0\leq u\leq1\), and, since \(\phi(0)=\phi(1)=0\),
	\begin{align*}
		\int_0^1\phi'(u)(t-u)^{-\frac{n+3}{2}} \mathcal G^{(2)}(1,0,t-u)\,du =\int_0^1\phi(u) \left.\partial_s\bigl[s^{-\frac{n+3}{2}} \mathcal G^{(2)}(1,0,s)\bigr]\right|_{s=t-u}\,du>0,
	\end{align*}
	so the negative prefactor in \eqref{0911-02} gives \(\partial_{x_1}p(\mathbf e_1',0,t)<0\).

	At \(t=t_\infty\), \(A(t_\infty)=0\) and \(D(t_\infty)=-nM_0T_1(t_\infty)<0\), since \(T_1(t_\infty)>0\) by the proof of Proposition \ref{prop-0504-03}.
	On the two sets of directions in \eqref{261049}, \((|x''|^2-(n+1)x_1^2)/r^2\) has opposite signs and is bounded away from zero, so the \(D(t_\infty)\)-term in \eqref{261048} gives \eqref{261049} for large \(r\).

	Assume {\rm\textbf{(T5)}} and let \(r\geq1\).
	By \eqref{0911-22f} and \eqref{260038},
	\begin{align}\label{1013-100}
		\partial_{x_1}p(r\mathbf e_1',0,t_\infty) =-\frac{\Gamma(\frac{n+2}{2})} {n2^n\pi^{\frac n2}\Gamma(\frac{n+3}{2})} \int_\Sigma\partial_{y_1}\psi(y')(r-y_1) F_{t_\infty}(|r\mathbf e_1'-y'|^2)\,dy'.
	\end{align}
	Pairing \((y_1,y'')\) with \((-y_1,y'')\) and writing the difference of the two values of \((r-y_1)F_{t_\infty}(|r\mathbf e_1'-y'|^2)\) as an integral in \(v\) gives
	\begin{align*}
		 & \int_\Sigma\partial_{y_1}\psi(y')(r-y_1) F_{t_\infty}(|r\mathbf e_1'-y'|^2)\,dy'\\
		 & \quad=\int_{\{y_1>0\}}\partial_{y_1}\psi(y') \int_{-y_1}^{y_1} \bigl[-F_{t_\infty}(r^2+|y'|^2-2rv) -2r(r-v)F_{t_\infty}'(r^2+|y'|^2-2rv)\bigr]\,dvdy'<0.
	\end{align*}
	Indeed, for \(|y'|<1\) and \(-y_1\leq v\leq y_1\), \(z=r^2+|y'|^2-2rv\) satisfies \(2r(r-v)-z=r^2-|y'|^2>0\), so Lemma \ref{lem-0911-02b}{\rm(i)} and \eqref{0911-35} with \(n\geq3\) give
	\begin{align*}
		-F_{t_\infty}(z)-2r(r-v)F_{t_\infty}'(z) \geq\left(\frac{(n-1)r(r-v)}z-1\right)F_{t_\infty}(z)>0,
	\end{align*}
	while \(\partial_{y_1}\psi\leq0\) on \(\{y_1>0\}\) and \(\partial_{y_1}\psi<0\) on a nonempty open set.
	Thus \eqref{1013-100} gives \eqref{0911-41} for \(r\geq1\), and evenness gives it for \(r\leq-1\).
	By \eqref{260160}, \(t\mapsto\partial_{x_1}p(\mathbf e_1',0,t)\) is continuous, positive at \(t_\infty\), and negative for large \(t\), so its last zero \(\overline t\) satisfies \(t_\infty<\overline t<\infty\), and it is negative on \((\overline t,\infty)\).
	Evenness gives \eqref{260286}.

	\underline{Proof of (iv).}\quad
	At a zero on the positive \(x_1\)-axis with \(r\) large, \eqref{0911-56b} gives \((n-1)A(t)+(n+1)D(t)r^{-2}+O(r^{-4})=0\).
	Since \(D(t)\to-nM_0T_1(t_\infty)<0\) as \(t\to t_\infty^-\), multiplying this equation by \(r^2/D(t)\) and using it to eliminate \(A(t)\) from \eqref{1013-106} give
	\begin{align}\label{261043}
		-\frac{A(t)r^2}{D(t)} =\frac{n+1}{n-1}+O(r^{-2}),\qquad \partial_r\partial_{x_1}p(r\mathbf e_1',0,t) =\kappa_nr^{-n-3}\bigl[2(n+1)D(t)+O(r^{-2})\bigr].
	\end{align}
	Choose \(\overline\delta\in(0,t_\infty-t_1)\) such that \(D(t)\leq-\frac n2M_0T_1(t_\infty)\) for \(t_\infty-\overline\delta<t<t_\infty\).
	Both \(O(r^{-2})\) terms in \eqref{261043} are uniform on this interval, and we choose \(R\geq R_0(t_\infty)\) such that, for \(r\geq R\), their absolute values are at most \(\frac{n+1}{2(n-1)}\) and \(\frac{n(n+1)}{2}M_0T_1(t_\infty)\), respectively.
	Then \(\partial_r\partial_{x_1}p(r\mathbf e_1',0,t)<0\) at every zero with \(r\geq R\) and \(t_\infty-\overline\delta<t<t_\infty\), so there is at most one zero on \([R,\infty)\).

	By \eqref{0911-41}, \(m_{\infty,R}:=\min_{1\leq r\leq R} \partial_{x_1}p(r\mathbf e_1',0,t_\infty)>0\), and the joint continuity in \eqref{260160} gives \(\delta_\infty\in(0,\overline\delta)\) such that
	\begin{align*}
		\sup_{1\leq r\leq R}\left| \partial_{x_1}p(r\mathbf e_1',0,t) -\partial_{x_1}p(r\mathbf e_1',0,t_\infty)\right|<\frac{m_{\infty,R}}{2},\qquad t_\infty-\delta_\infty<t<t_\infty.
	\end{align*}
	For such \(t\), \(\partial_{x_1}p(r\mathbf e_1',0,t)\) is positive on \([1,R]\) and, by \eqref{0911-56b}, negative for large \(r\), so it has exactly one zero \(\rho(t)\) on \((1,\infty)\), \(\rho(t)>R\), and evenness gives \eqref{1013-107}.
	At \(r=\rho(t)\), \eqref{261043} gives
	\begin{align}\label{261069}
		\frac{n+1}{2(n-1)} \leq-\frac{A(t)\rho(t)^2}{D(t)} \leq\frac{3(n+1)}{2(n-1)}.
	\end{align}
	Since \(A(t)\to0^+\) and \(D(t)\to-nM_0T_1(t_\infty)\), \eqref{261069} gives \(\rho(t)\to\infty\) as \(t\to t_\infty^-\).
	Hence, by \eqref{261043}, \(-A(t)\rho(t)^2/D(t)\to(n+1)/(n-1)\), and \(A(t)/(t_\infty-t)\to-M_0T_0'(t_\infty)\) proves \eqref{261042}.

	\underline{Proof of (v).}\quad
	Fix \(t>\overline t\).
	By \eqref{260286}, \eqref{261041}, and the continuity in \eqref{260160}, \(\partial_{x_1}p(r\mathbf e_1',0,t)\) is negative at \(r=1\) and positive for large \(r\), so \(\varnothing\neq Z_t:=\{r>1:\partial_{x_1}p(r\mathbf e_1',0,t)=0\}\subset[1+\eta,L]\) for some \(\eta>0\) and \(L>1+\eta\).
	As a closed set, \(Z_t\) has a minimum \(\rho_-(t)\) and a maximum \(\rho_+(t)\), and evenness gives \eqref{1013-108}.
	For \(t\geq\widehat t:=T_{\delta_0}\), Lemma~\ref{lem-261080}{\rm(iii)}--{\rm(iv)} gives the unique simple zero \(\rho(t)\) in \((1,\infty)\), the zero \(-\rho(t)\), and \eqref{261034}, so \(\rho_-(t)=\rho_+(t)=\rho(t)\), and also \(\partial_{x_1}p(\mathbf e_1',0,t)<0\).
	Since \(\partial_{x_1}p(\mathbf e_1',0,\overline t)=0\) by {\rm(iii)}, \(\overline t<\widehat t\).
\end{proof}
\subsubsection{The \(i=j\) case when \(n\geq3\): zero curves of \(\partial_{x_1}p\)}\label{subsubsect-261022-offaxis}

Throughout this subsection, write
\begin{align*}
    P_t(r):=\partial_{x_1}p(r\mathbf e_2',0,t),\qquad r>1,\quad t>0.
\end{align*}
Primes on \(P_t\) denote derivatives in \(r\).
By \eqref{0911-20g}, the equation for a zero at radius \(r=|x'|\) is \(P_t(r)+(x_1^2/r)P_t'(r)=0\).
The next lemma gives the radial monotonicity needed to solve this equation for \(x_1>0\).
\begin{lemma}\label{lem-0911-08}
	Let \(n\geq3\).
	Suppose that \(\psi\) satisfies {\rm\textbf{(S)}} and is radial, and that \(\phi\) satisfies {\rm\textbf{(T1)}}.
	Fix \(t>0\) with \(\norm{\phi}_{L^1(0,t)}>0\).
	Then there exists \(C(n)\geq1\) for which \(R_{\rm off}(t)\) may be chosen as follows.
	If \(A(t)D(t)<0\), let
	\begin{align}\label{0911-60a}
		R_{\rm off}(t)=\max\left\{R_0(t), \left(\frac{C(n)|D(t)|}{|A(t)|}\right)^{\frac{1}{2}}, \left(\frac{C(n)E(t)}{|A(t)|}\right)^{\frac{1}{4}}\right\}.
	\end{align}
	If \(A(t)D(t)\geq0\), let
	\begin{align}\label{0911-60b}
		R_{\rm off}(t)=\max\left\{R_0(t), \min\left\{ \left(\frac{C(n)E(t)}{|A(t)|}\right)^{\frac{1}{4}}, \left(\frac{C(n)E(t)}{|D(t)|}\right)^{\frac{1}{2}} \right\}\right\},
	\end{align}
where any quotient with a zero denominator is interpreted as infinity.
	For \(r>R_{\rm off}(t)\), \(P_t'(r)\neq0\), and
	\begin{equation}\label{260189}
		\begin{aligned}
			\partial_r\left(r^2+ \frac{rP_t(r)} {P_t'(r)}\right)>0,\qquad \partial_r\left( \frac{rP_t(r)} {P_t'(r)}\right)<0.
		\end{aligned}
	\end{equation}
\end{lemma}
\begin{proof}
	Fix \(t>0\), and let \(c_0(n)\) be the constant in the remainder bound \eqref{0911-60c}.
	Integration by parts in \eqref{0323-02} gives
	\begin{align*}
		T_1(t)=\frac12\int_0^t\frac{\phi(s)}{\sqrt{t-s}}\,ds>0,
	\end{align*}
	since \(\phi\geq0\) and its integral on \((0,t)\) is positive.
	Thus \eqref{260279} implies that \(A(t)\) and \(D(t)\) do not vanish simultaneously.
	If \(A(t)D(t)<0\), \eqref{0911-60a} gives
	\begin{align*}
		\frac{|D(t)|}{r^2}<\frac{|A(t)|}{C(n)},\qquad \frac{E(t)}{r^4}<\frac{|A(t)|}{C(n)}.
	\end{align*}
	If \(A(t)D(t)\geq0\), \eqref{0911-60b} gives
	\begin{align*}
		\frac{E(t)}{r^4}<\frac{|A(t)|}{C(n)} \quad\hbox{or}\quad \frac{E(t)}{r^4}<\frac{|D(t)|}{C(n)r^2}.
	\end{align*}
	Thus the choice of \(R_{\rm off}(t)\) gives the following bounds for every \(r>R_{\rm off}(t)\).
	\begin{equation}\label{260274}
		\begin{aligned}
			A(t)D(t)<0& \ \Longrightarrow\ \frac{|D(t)|}{r^2}+\frac{E(t)}{r^4} \leq\frac{2|A(t)|}{C(n)},\\
			A(t)D(t)\geq0 & \ \Longrightarrow\ \frac{E(t)}{r^4} \leq\frac{|A(t)|}{C(n)}+\frac{|D(t)|}{C(n)r^2}.
		\end{aligned}
	\end{equation}

	We first use these bounds to control the denominator \(P_t'(r)\).
	If \(A(t)D(t)<0\), then \eqref{0911-601}, \eqref{0911-60c}, and the first inequality in \eqref{260274} give
	\begin{align*}
		\frac{r^{n+1}}{\kappa_n} \left|P_t'(r)\right| & \geq n|A(t)|-\frac{(n+2)|D(t)|}{r^2}-|e_1(r,t)|\geq\left(n-\frac{2(n+2+c_0(n))}{C(n)}\right)|A(t)|>0
	\end{align*}
	after increasing \(C(n)\).
	If \(A(t)D(t)\geq0\), then the second inequality in \eqref{260274} gives
	\begin{align*}
		\frac{r^{n+1}}{\kappa_n} \left|P_t'(r)\right| & \geq n|A(t)|+\frac{(n+2)|D(t)|}{r^2}-|e_1(r,t)|\geq\left(n-\frac{c_0(n)}{C(n)}\right) \left(|A(t)|+\frac{|D(t)|}{r^2}\right)>0
	\end{align*}
	for \(C(n)\) sufficiently large, since \(A(t)\) and \(D(t)\) cannot both vanish.
	Thus \(P_t'(r)\ne0\) in both cases.

	We may now differentiate the quotients in \eqref{260189}.
	The quotient rule gives
	\begin{align}\label{0911-61}
		\partial_r\left(r^2+ \frac{rP_t(r)} {P_t'(r)}\right) & =3r+\frac{P_t(r)} {P_t'(r)} \left(1-\frac{rP_t''(r)} {P_t'(r)}\right).
	\end{align}
	Multiplying the derivative in \eqref{0911-61} by \(P_t'(r)^2\) gives \(3rP_t'(r)^2+P_t(r)P_t'(r)-rP_t(r)P_t''(r)\).
	Insert \eqref{0911-60}--\eqref{0911-602} into this expression and divide by \(\kappa_n^2r^{-2n-1}\).
	Completing the square in the terms independent of \(e_0,e_1,e_2\) gives, for \(n\geq3\),
	\begin{align}\label{260299}
		 & 2n(n-1)A(t)^2+ \frac{4(n-1)(n+2)A(t)D(t)}{r^2} +\frac{2(n+1)(n+2)D(t)^2}{r^4}\nonumber \\
		 & \qquad=\frac{2(n-1)}n \left(nA(t)+\frac{(n+2)D(t)}{r^2}\right)^2 +\frac{4(n+2)D(t)^2}{nr^4} \geq\frac12\left(A(t)^2+\frac{D(t)^2}{r^4}\right).
	\end{align}
	For the last inequality, write \(B=D(t)/r^2\) and \(X=nA(t)+(n+2)B\).
	Since \(A(t)=X/n-(1+2/n)B\), the inequality \((u-v)^2\leq2u^2+2v^2\) gives
	\begin{align*}
	    \tfrac12\bigl(A(t)^2+B^2\bigr)
	    \leq \frac{X^2}{n^2}+\bigl((1+2/n)^2+\tfrac12\bigr)B^2
	    \leq \frac{2(n-1)}nX^2+\frac{4(n+2)}nB^2,
	\end{align*}
	where \(n\geq3\) in the last step.

	The derivative of \(rP_t(r)/P_t'(r)\) differs from \eqref{0911-61} by \(-2r\).
	Hence its terms independent of \(e_0,e_1,e_2\), with the same normalization, are
	\begin{align*}
	    -2nA(t)^2-\frac{4(n+2)A(t)D(t)}{r^2}-\frac{2(n+2)D(t)^2}{r^4}.
	\end{align*}
	If \(A(t)D(t)<0\), increase \(C(n)\) so that \(C(n)\geq\max\{2,8(n+2)/n\}\).
	The first inequality in \eqref{260274} then gives
	\begin{align*}
		-2nA(t)^2-\frac{4(n+2)A(t)D(t)}{r^2} -\frac{2(n+2)D(t)^2}{r^4} \leq-nA(t)^2 \leq-\frac n2\left(A(t)^2+\frac{D(t)^2}{r^4}\right).
	\end{align*}
	If \(A(t)D(t)\geq0\), then \(-4(n+2)A(t)D(t)/r^2\leq0\) and
	\begin{align*}
		-2nA(t)^2-\frac{2(n+2)D(t)^2}{r^4} \leq-\frac n2\left(A(t)^2+\frac{D(t)^2}{r^4}\right).
	\end{align*}
	Hence
	\begin{align}\label{260300}
		-2nA(t)^2-\frac{4(n+2)A(t)D(t)}{r^2} -\frac{2(n+2)D(t)^2}{r^4} \leq-\frac n2\left(A(t)^2+\frac{D(t)^2}{r^4}\right).
	\end{align}

	Every other term in the two normalized derivatives contains at least one of \(e_0,e_1,e_2\).
	Bounding each such factor by \eqref{0911-60c} gives, for \(r\geq R_0(t)\) wherever \(P_t'(r)\ne0\),
	\begin{equation}\label{260275}
		\begin{aligned}
			 & \left| \frac{(P_t'(r))^2} {\kappa_n^2r^{-2n-1}} \partial_r\left(r^2+ \frac{rP_t(r)} {P_t'(r)}\right)\right. \\
			 & \qquad\left.- \left(2n(n-1)A(t)^2+ \frac{4(n-1)(n+2)A(t)D(t)}{r^2} +\frac{2(n+1)(n+2)D(t)^2}{r^4}\right)\right|\\
			 & \qquad\leq\frac{c_3(n)E(t)}{r^4} \left(|A(t)|+\frac{|D(t)|}{r^2}+\frac{E(t)}{r^4}\right),\\
			 & \left| \frac{(P_t'(r))^2} {\kappa_n^2r^{-2n-1}} \partial_r\left( \frac{rP_t(r)} {P_t'(r)}\right) + 2nA(t)^2+\frac{4(n+2)A(t)D(t)}{r^2} +\frac{2(n+2)D(t)^2}{r^4}\right|\\
			 & \qquad\leq\frac{c_3(n)E(t)}{r^4} \left(|A(t)|+\frac{|D(t)|}{r^2}+\frac{E(t)}{r^4}\right).
		\end{aligned}
	\end{equation}
	Here \(c_3(n)\) depends only on \(n\).
	For \(r>R_{\rm off}(t)\), increase \(C(n)\) in \eqref{260274} so that the common error bound satisfies
	\begin{align}\label{260301}
		\frac{c_3(n)E(t)}{r^4} \left(|A(t)|+\frac{|D(t)|}{r^2}+\frac{E(t)}{r^4}\right) \leq\frac14 \left(A(t)^2+\frac{D(t)^2}{r^4}\right).
	\end{align}
	Then \eqref{260299}--\eqref{260301} give
	\begin{align*}
		 & \frac{(P_t'(r))^2} {\kappa_n^2r^{-2n-1}} \partial_r\left(r^2+ \frac{rP_t(r)} {P_t'(r)}\right) \geq\frac14 \left(A(t)^2+\frac{D(t)^2}{r^4}\right)>0,\\
		 & \frac{(P_t'(r))^2} {\kappa_n^2r^{-2n-1}} \partial_r\left( \frac{rP_t(r)} {P_t'(r)}\right) \leq-\frac n4 \left(A(t)^2+\frac{D(t)^2}{r^4}\right)<0.
	\end{align*}
	Since the normalizing factor is positive, these are the two signs in \eqref{260189}.
\end{proof}

The next proposition gives the zero curve used in Theorem \ref{thm-260021}{\rm(vi)}.
\begin{proposition}\label{prop-261022}
	Let \(n\geq3\), and let \(p\) be the pressure of the Stokes system \eqref{0308-01} with boundary data \eqref{0324-13}, where \(i=j=1\).
	Suppose that \(\psi\) satisfies {\rm\textbf{(S)}} and is radial, and that \(\phi\) satisfies {\rm\textbf{(T1)}}.
	Fix \(t>0\) and suppose that \(\norm{\phi}_{L^1(0,t)}>0\).
	On each line parallel to the \(x_1\)-axis with \(|x''|\geq1\), \(x_1\mapsto\partial_{x_1}p(x_1,x'',0,t)\) changes sign at a positive zero.
	For \(|x''|>R_{\rm off}(t)\), with \(R_{\rm off}(t)\) given in \eqref{0911-60a}--\eqref{0911-60b}, this zero is unique and simple.
	Its coordinate \(\rho(|x''|,t)\) is a strictly increasing \(C^1\) function of \(|x''|\) for \(|x''|>R_{\rm off}(t)\).
	For fixed \(t\), as \(|x''|\to\infty\),
	\begin{align}\label{261036}
		\rho(|x''|,t)=
		\begin{cases}
			\dfrac{|x''|}{\sqrt{n-1}}+
			\dfrac{1}{\sqrt{n-1}\,|x''|}\!
			\left(\dfrac{T_1(t)}{T_0(t)}-\dfrac{M_2}{2(n-1)M_0}\right)
			+O_t(|x''|^{-3}),& T_0(t)\ne0,\\[6pt]
			\dfrac{|x''|}{\sqrt{n+1}}+O_t(|x''|^{-1}),&T_0(t)=0.
		\end{cases}
	\end{align}
	Suppose, in addition, that \(\psi\) is radially nonincreasing and that \(\phi\) satisfies {\rm\textbf{(T2)}}--{\rm\textbf{(T5)}}.
Let \(t<t_\infty\) be sufficiently close to \(t_\infty\), and let \(\rho(t)\) be the unique zero of \(r\mapsto\partial_{x_1}p(r\mathbf e_1',0,t)\) on \((1,\infty)\).
	For every \(x''\), there is exactly one positive zero \(\rho(|x''|,t)\) with \(|x'|\geq\rho(t)\).
	This is the largest positive zero for the given \(x''\), and, when \(|x''|> R_{\rm off}(t)\), it is the unique zero obtained above.
	These zeros are simple, and \(\rho(\cdot,t)\) has an even real-analytic extension from \([0,\infty)\) to \(\R\).
	At \(x''=0\) it has a nondegenerate local maximum, with \(\rho(0,t)=\rho(t)\). For all \(x''\),
	\begin{align}\label{261037}
		\rho(|x''|,t)^2+|x''|^2\geq\rho(t)^2.
	\end{align}
	Moreover, \(\rho(\cdot,t)\) is strictly decreasing on \([0,\eta(t)]\) for some \(\eta(t)>0\).
\end{proposition}
\begin{proof}
	Let \(|x''|\geq1\), and write \(r=(x_1^2+|x''|^2)^{1/2}\).
	For \(x_1>0\), \(r>1\), and \eqref{0911-20g} gives \(\partial_{x_1}p(x',0,t)=\frac{d}{dx_1}\bigl[x_1P_t(r)\bigr]\).
	As \(x_1\downarrow0\), \(x_1P_t(r)\to0\), since \(P_t\) is continuous at \(r=1\) by \eqref{260160}, and \eqref{0911-18} bounds the derivative there.
	As \(x_1\to\infty\), \(x_1P_t(r)\to0\) by \eqref{0911-60}--\eqref{0911-60c}.
	Hence
	\begin{align}\label{0911-58}
		\int_0^\infty\partial_{x_1}p(x',0,t)\,dx_1=0,\qquad |x''|\geq1.
	\end{align}
	By the proof of Lemma \ref{lem-0911-08}, \(A(t)\) and \(D(t)\) do not vanish together, so \eqref{0911-60} gives \(P_t(r)\ne0\) for large \(r\).
	Since the antiderivative \(x_1P_t(r)\) is then nonzero and tends to zero, \(x_1\mapsto\partial_{x_1}p(x',0,t)\) is not identically zero.
	By \eqref{0911-58}, it takes both signs on \((0,\infty)\).
	It is real-analytic there by \eqref{0822-13} and Lemma~\ref{lem-260044}{\rm(vi)}, so Lemma~\ref{lem-260044}{\rm(iii)} gives finitely many zeros between two points of opposite signs, and it changes sign at one of them.

	Now let \(|x''|>R_{\rm off}(t)\).
	At any positive zero, \(r>R_{\rm off}(t)\), so \(P_t'(r)\ne0\) by Lemma \ref{lem-0911-08}, and \eqref{0911-20g} with \(r^2=x_1^2+|x''|^2\) gives
	\begin{align}\label{0911-50}
		x_1^2 & =-\frac{rP_t(r)} {P_t'(r)},\qquad |x''|^2 =r^2+\frac{rP_t(r)} {P_t'(r)}.
	\end{align}
	By \eqref{260189}, the right-hand side of the second identity is strictly increasing in \(r\), so it determines \(r\), and \(x_1=(r^2-|x''|^2)^{1/2}\).
	Hence the positive zero is unique and depends on \(x''\) only through \(|x''|\).
	At this zero, differentiating \eqref{0911-20g} and using \eqref{0911-50} give
	\begin{align}\label{261062}
		\partial_{x_1}^2p(x',0,t) & =\frac{x_1}{r^2}P_t'(r) \partial_r\left(r^2+ \frac{rP_t(r)} {P_t'(r)}\right),
	\end{align}
	which is nonzero by \eqref{260189}.
	Thus the zero is simple, and the implicit function theorem gives \(\rho(\cdot,t)\in C^1((R_{\rm off}(t),\infty))\).
	Differentiating \eqref{0911-50} along the zero curve and eliminating the derivative of \(r\) give
	\begin{align*}
		\partial_{|x''|}\rho(|x''|,t)
		=-\frac{|x''|}{\rho(|x''|,t)}
		\frac{\displaystyle\partial_r\left(\frac{rP_t(r)}{P_t'(r)}\right)}
		{\displaystyle\partial_r\left(r^2+\frac{rP_t(r)}{P_t'(r)}\right)}>0
	\end{align*}
	by \eqref{260189}.

	If \(T_0(t)\ne0\), then \(A(t)\ne0\), and \eqref{0911-60}, \eqref{0911-601}, and \eqref{0911-60c} give
	\begin{align*}
		x_1^2=-\frac{rP_t(r)}{P_t'(r)}
		=\frac{r^2}{n}-\frac{2D(t)}{n^2A(t)}+O_t(r^{-2}).
	\end{align*}
	With \(r^2=x_1^2+|x''|^2\), this gives \((n-1)x_1^2=|x''|^2-\frac{2D(t)}{nA(t)}+O_t(|x''|^{-2})\), and \eqref{260279} gives the first line of \eqref{261036}.
	If \(T_0(t)=0\), then \(A(t)=0\), \(D(t)<0\), and the same quotient gives \(x_1^2=r^2/(n+2)+O_t(1)\), so \((n+1)x_1^2=|x''|^2+O_t(1)\), which is the second line.

	Assume now that \(\psi\) is radially nonincreasing and \(\phi\) satisfies {\rm\textbf{(T2)}}--{\rm\textbf{(T5)}}, and let \(\delta_\infty\) be as in Proposition \ref{prop-261020}{\rm(iv)}.
	Since \(D(t)\to-nM_0T_1(t_\infty)<0\), \(E\) is bounded on compact time intervals, and \(R_0(t)=2+C_n\sqrt t\), there are \(0<\delta<\delta_\infty\), \(d_0,d_1,e_0>0\), and \(R_*>1\) such that
	\begin{align}\label{261063}
		d_0\leq |D(t)|\leq d_1,\qquad E(t)\leq e_0,\qquad R_0(t)\leq R_*,\qquad t_\infty-\delta<t<t_\infty.
	\end{align}
	Since \(A(t)>0\) there and \(2A(t)|D(t)|r^{-2}\leq A(t)^2+D(t)^2r^{-4}\), \eqref{261063} bounds the error in \eqref{260275} by \(Cr^{-2}(A(t)^2+D(t)^2r^{-4})\) with \(C=C(n,\psi,\phi)\).
	After increasing \(R_*\), this is at most \(\tfrac14(A(t)^2+D(t)^2r^{-4})\) for \(r\geq R_*\), and, by \eqref{261042}, we may reduce \(\delta\) so that \(\rho(t)>R_*\).

	By \eqref{261043} and \(\rho(t)\to\infty\), \(\frac{\rho(t)^2}{|D(t)|}\bigl[-nA(t)+\frac{(n+2)|D(t)|}{\rho(t)^2}+\frac{C(n)E(t)}{\rho(t)^4}\bigr]\to-\frac{2}{n-1}\) as \(t\to t_\infty^-\), so, after reducing \(\delta\),
	\begin{align}\label{261071}
		-nA(t)+\frac{(n+2)|D(t)|}{\rho(t)^2} +\frac{C(n)E(t)}{\rho(t)^4}<0.
	\end{align}
	The left-hand side decreases when \(\rho(t)\) is replaced by \(r\geq\rho(t)\), so \eqref{0911-601} and \eqref{0911-60c} give \(P_t'(r)<0\) for \(r\geq\rho(t)\).
	As \(P_t(r)>0\) for large \(r\) by \eqref{0911-20m}, this gives
	\begin{align}\label{261064}
		P_t(r)>0,\qquad P_t'(r)<0,\qquad r\geq\rho(t).
	\end{align}
	By \eqref{260299}, \eqref{260275}, and the choice of \(R_*\),
	\begin{align*}
		\frac{(P_t'(r))^2} {\kappa_n^2r^{-2n-1}} \partial_r\left(r^2+ \frac{rP_t(r)} {P_t'(r)}\right) \geq\frac14\left(A(t)^2+\frac{D(t)^2}{r^4}\right)>0,\qquad r\geq\rho(t),
	\end{align*}
	so, with \eqref{261064},
	\begin{equation}\label{261067}
		\partial_r\left(r^2+ \frac{rP_t(r)} {P_t'(r)}\right)>0, \qquad r\geq\rho(t).
	\end{equation}
	At \(r=\rho(t)\), \eqref{0911-20f} gives \(P_t(r)+rP_t'(r)=0\), that is,
	\begin{align}\label{261068}
		-\frac{rP_t(r)} {P_t'(r)}=r^2,\qquad r^2+\frac{rP_t(r)} {P_t'(r)}=0,
	\end{align}
	and, since \(A(t)>0\), \eqref{0911-60}--\eqref{0911-601} give
	\begin{align}\label{261065}
		r^2+\frac{rP_t(r)} {P_t'(r)} =\frac{n-1}{n}r^2+O_t(1)\longrightarrow\infty.
	\end{align}
	Hence, for each \(|x''|\geq0\), the equation
	\begin{align}\label{261066}
		|x''|^2=r^2+ \frac{rP_t(r)} {P_t'(r)}
	\end{align}
	has a unique solution \(r\in[\rho(t),\infty)\).
	With \(x_1:=(-rP_t(r)/P_t'(r))^{1/2}\), well-defined by \eqref{261064}, \eqref{261066} gives \(x_1^2+|x''|^2=r^2\), and \eqref{0911-20g} gives \(\partial_{x_1}p(x',0,t)=0\).

	Let \(\widetilde x_1>0\) be another zero for the same \(x''\), and set \(\widetilde r=(\widetilde x_1^2+|x''|^2)^{1/2}\).
	If \(\widetilde r\geq\rho(t)\), then \(P_t'(\widetilde r)\ne0\) by \eqref{261064}, so \(\widetilde r\) solves \eqref{261066} and \(\widetilde x_1=x_1\).
	If \(\widetilde r<\rho(t)\leq r\), then \(\widetilde x_1<x_1\).
	Thus the constructed zero is the largest positive zero, and, for \(|x''|>R_{\rm off}(t)\), it is \(\rho(|x''|,t)\).

	Let \(G_t(r):=r^2+rP_t(r)/P_t'(r)\).
	By Remark~\ref{260142}, \(P_t\) is real-analytic for \(r>1\), and \(P_t'<0\) and \(G_t'>0\) on \([\rho(t),\infty)\) by \eqref{261064} and \eqref{261067}.
	Hence \eqref{261062} gives \(\partial_{x_1}^2p<0\) at these zeros, and the inverse function theorem gives a real-analytic inverse of \(G_t\) on \([0,\infty)\), with \(G_t^{-1}(0)=\rho(t)\).
	Therefore \(\rho(|x''|,t)=\bigl(G_t^{-1}(|x''|^2)^2-|x''|^2\bigr)^{1/2}\), with a positive radicand by \eqref{261064} and \eqref{261066}, which gives the even real-analytic extension, \(\rho(0,t)=\rho(t)\), and, since \(G_t^{-1}(|x''|^2)\geq\rho(t)\), \eqref{261037}.

	Finally, \(\theta(t):=|D(t)|/(A(t)\rho(t)^2)\to\frac{n-1}{n+1}\) by \eqref{261043}.
	Since \(A(t)D(t)<0\), at \(r=\rho(t)\) the terms \(-2nA(t)^2-4(n+2)A(t)D(t)r^{-2}-2(n+2)D(t)^2r^{-4}\) in the second estimate of \eqref{260275} equal \(A(t)^2\bigl(-2n+4(n+2)\theta(t)-2(n+2)\theta(t)^2\bigr)\), with coefficient tending to \(4(n^2-3)/(n+1)^2\geq3/2\).
	By the choice of \(R_*\), the error there is at most \(\tfrac14(1+\theta(t)^2)A(t)^2\leq\tfrac12A(t)^2\) when \(\theta(t)\leq1\).
	After reducing \(\delta\), we obtain \(\partial_r(rP_t(r)/P_t'(r))>0\) at \(r=r_0:=\rho(t)\).
	Differentiating \(\rho(|x''|,t)^2=G_t^{-1}(|x''|^2)^2-|x''|^2\) twice in \(|x''|\) at \(0\) gives
	\begin{align*}
	    \rho'(0,t)=0,\qquad
	    \rho''(0,t)
	    =\frac{2r_0-G_t'(r_0)}{r_0G_t'(r_0)}
	    =-\frac{(rP_t/P_t')'(r_0)}{r_0G_t'(r_0)}<0,
	\end{align*}
	where the primes on \(\rho\) denote derivatives in its first argument.
	Hence the maximum at \(x''=0\) is nondegenerate, and \(\rho(\cdot,t)\) is strictly decreasing on \([0,\eta(t)]\) for some \(\eta(t)>0\).
\end{proof}

	\subsubsection{Distinct tangential indices \texorpdfstring{\(i\neq j\)}{i and j}}
The next proposition determines the sign of \(\partial_{x_i}p(x',0,t)\) for \(i\neq j\), including two sign regions throughout \(t_1<t<t_\infty\), with the first region possibly empty.
Its four statements prove Theorem \ref{thm-260027}{\rm(i)}--{\rm(iv)}.

\begin{proposition}\label{prop-0831-01}
	Let \(n\geq3\), \(1\leq i\neq j\leq n-1\), and let \(p\) be the pressure of the Stokes system \eqref{0308-01} with boundary data \eqref{0324-13}.
	Suppose that \(\psi\) satisfies {\rm\textbf{(S)}} and is radial, and that \(\phi\) satisfies {\rm\textbf{(T1)}} and \(\phi\not\equiv0\).
	Then the following hold.
	\begin{enumerate}
		\item Suppose that \(\phi\) satisfies {\rm\textbf{(T2)}}. If \(|x'|\geq\sqrt2\) and \(0<t\leq t_1\), then \(\sgn(\partial_{x_i}p(x',0,t))=-\sgn(x_ix_j)\).
		\item Suppose that \(\phi\) satisfies {\rm\textbf{(T2)}}--{\rm\textbf{(T4)}}. If \(t_1<t<t_\infty\), then there exists \(r_{\rm tan}(t)>0\) such that
		\begin{align*}
			\sgn(\partial_{x_i}p(x',0,t)) =
			\begin{cases}
				\sgn(x_ix_j),& \sqrt2\leq |x'|\leq r_{\rm tan}(t)-1, \\
				-\sgn(x_ix_j), & |x'|\geq\max\{\sqrt2,r_{\rm tan}(t)+1\},
			\end{cases}
		\end{align*}
		where the first region is empty if \(r_{\rm tan}(t)<1+\sqrt2\).
		Moreover, \(r_{\rm tan}(t)\) is continuous and satisfies
		\begin{align}\label{261046}
			r_{\rm tan}(t) =\left(-\frac{(n+2)T_1(t_\infty)}{T_0'(t_\infty)}\right)^{\frac12} (t_\infty-t)^{-\frac12}(1+o(1)) \qquad\text{as }t\to t_\infty^-.
		\end{align}
		\item Suppose that \(\phi\) satisfies {\rm\textbf{(T2)}}--{\rm\textbf{(T4)}}. Then, for \(t<t_\infty\) sufficiently close to \(t_\infty\), there exists a unique \(\rho_{\rm tan}(t)>\sqrt2\) such that
		\begin{align*}
			\sgn(\partial_{x_i}p(x',0,t))=-\sgn(x_ix_j)\sgn(|x'|-\rho_{\rm tan}(t)),\qquad |x'|\geq\sqrt2.
		\end{align*}
		Along each ray with \(x_ix_j\ne0\), this zero is simple.
		As \(t\to t_\infty^-\), \(\rho_{\rm tan}(t)-r_{\rm tan}(t)=O(r_{\rm tan}(t)^{-1})\), and
		\begin{align*}
			\rho_{\rm tan}(t)=\left(-\frac{(n+2)T_1(t_\infty)}{T_0'(t_\infty)}\right)^{\frac12}(t_\infty-t)^{-\frac12}(1+o(1)) \qquad\text{as }t\to t_\infty^-.
		\end{align*}
		\item If \(|x'|\geq\sqrt2\) and \(t\geq1\), then \(\sgn(\partial_{x_i}p(x',0,t))=\sgn(x_ix_j)\).
		Moreover, if \(\phi\) satisfies {\rm\textbf{(T2)}}--{\rm\textbf{(T4)}}, then the same sign holds for \(t_\infty\leq t<1\).
	\end{enumerate}
\end{proposition}

The following lemma gives the properties of \(H_2\) used in the proof.
\begin{lemma}\label{lem-0903-01}
	Let \(n\geq2\), and suppose that \(\phi\) satisfies {\rm\textbf{(T1)}} and \(\phi\not\equiv0\). Then the following hold.
	\begin{enumerate}
		\item\label{item-0906-02} If \(t\geq1\), then \(H_2(r,t)<0\) for all \(r>0\).
		If \(\phi\) satisfies {\rm\textbf{(T2)}}, then \(H_2(r,t)>0\) for \(0<t\leq t_1\) and all \(r>0\).
		\item\label{item-0906-01} Suppose that \(\phi\) satisfies {\rm\textbf{(T2)}}, and let \(t_1<t<1\).
		If \(T_0(t)\leq0\), then \(H_2(r,t)<0\) for all \(r>0\).
		If \(T_0(t)>0\), then \(H_2(\cdot,t)\) has a unique simple zero \(r_*(t)>0\).
		It is negative for \(0<r<r_*(t)\) and positive for \(r>r_*(t)\).
		\item\label{item-0906-03} If \(\phi\) satisfies {\rm\textbf{(T2)}}--{\rm\textbf{(T4)}}, then \(r_*\) in {\rm(ii)} exists and is continuous on \((t_1,t_\infty)\), and
		\begin{align}\label{0904-03}
			\lim_{t\to t_1^+}r_*(t)=0,\qquad \lim_{t\to t_\infty^-}r_*(t)=\infty.
		\end{align}
		Moreover,
		\begin{align}\label{261008}
			r_*(t) =\left(-\frac{(n+2)T_1(t_\infty)}{T_0'(t_\infty)}\right)^{\frac12} (t_\infty-t)^{-\frac12}(1+o(1)) \qquad\text{as }t\to t_\infty^-.
		\end{align}
	\end{enumerate}
\end{lemma}
\begin{proof}[Proof of Lemma \ref{lem-0903-01}]
	We first establish three facts used in {\rm(i)}--{\rm(iii)}.
	Integration by parts in \eqref{0831-06}, using \eqref{0911-21} with \((a,b,c)=((n+2)/2,(n+5)/2,r^2/4)\), gives
	\begin{align}\label{0831-07}
		H_2(r,t)=\frac{2}{n+3}\int_0^t\frac{\phi'(t-s)}{s^{\frac{n+3}{2}}} {}_1F_1\left(\frac{n+2}{2},\frac{n+5}{2};-\frac{r^2}{4s}\right)\,ds.
	\end{align}
	The corresponding boundary terms vanish by {\rm\textbf{(T1)}} and \eqref{0514-02}.

	For \(r,s>0\), write
	\begin{align*}
	    K_r(s)=s^{-(n+3)/2}P(r^2/(4s)),\qquad
	    P(z)={}_1F_1\!\left(\tfrac{n+2}{2},\tfrac{n+5}{2};-z\right),
	    \qquad L=\ln P .
	\end{align*}
	Take \(0<r_1<r_2\), set \(z_j=r_j^2/(4s)\), and consider the ratio \(q(s)=K_{r_2}(s)/K_{r_1}(s)\).
	Lemma~\ref{lem-0324-01}{\rm(v)} in dimension \(n+2\) shows that \(zL'(z)\) decreases strictly. Hence
	\begin{align*}
	    (\ln q)'(s)=\frac{z_1L'(z_1)-z_2L'(z_2)}{s}>0.
	\end{align*}
	Under {\rm\textbf{(T2)}}, for \(t_1<t<1\), write \(s_0=t-t_1\).
	Then \(\phi'(t-s)\) and \(q(s)-q(s_0)\) have the same strict sign on each side of \(s_0\). Subtracting \(q(s_0)\) times \eqref{0831-07} with \(r=r_1\) from \eqref{0831-07} with \(r=r_2\) gives
	\begin{align}\label{0831-08}
	    H_2(r_2,t)-q(s_0)H_2(r_1,t)
	    =\frac{2}{n+3}\int_0^t\phi'(t-s)K_{r_1}(s)[q(s)-q(s_0)]\,ds>0 .
	\end{align}
	In particular, a nonnegative value at \(r_1\) implies strict positivity at every \(r_2>r_1\).

	For the limit at zero, set \(s=r^2/(4u)\) in \eqref{0831-07} and extend the resulting integrand by zero below \(r^2/(4t)\).
	By \eqref{0514-02}, this integrand multiplied by \(r^{n+1}\) is bounded by \(C\|\phi'\|_{L^\infty(0,1)}\min\{u^{(n-1)/2},u^{-3/2}\}\).
	For the limit at infinity, \(r^{n+2}K_r(s)\leq Cs^{-1/2}\) by \eqref{0514-02}.
	Thus dominated convergence and \eqref{0327-10} give
	\begin{align*}
		\lim_{r\to0^+}r^{n+1}H_2(r,t) & =\frac{2^{n+2}}{n+3}\phi'(t) \int_0^\infty u^{\frac{n-1}{2}}{}_1F_1\left(\frac{n+2}{2},\frac{n+5}{2};-u\right)\,du, \\
		\lim_{r\to\infty}r^{n+2}H_2(r,t) & =\frac{2^{n+3}\Gamma(\frac{n+5}{2})}{(n+3)\Gamma(\frac32)}T_0(t).
	\end{align*}
	By \eqref{0324-14} and Tonelli's theorem, the integral in the first line is a beta integral, and
	\begin{align*}
		\int_0^\infty u^{(n-1)/2}{}_1F_1\left(\tfrac{n+2}{2},\tfrac{n+5}{2};-u\right)\,du
		=\frac{\sqrt\pi\,\Gamma((n+1)/2)\Gamma((n+5)/2)}{\Gamma((n+2)/2)}>0.
	\end{align*}

	\underline{Proof of (i).}\quad
	If \(0<t\leq t_1\), then \eqref{0831-07} gives \(H_2(r,t)>0\) since \(\phi'(t-s)>0\) for \(0<s<t\).
	For \(t\geq1\), we use \eqref{0831-06} directly. In this case, \(\phi(t)-\phi(t-s)=-\phi(t-s)\leq0\), and \({}_1F_1(\frac{n+2}{2},\frac{n+3}{2};-\frac{r^2}{4s})>0\) by \eqref{0324-14}.
	Since \(\phi\) is nonzero and nonnegative, the integrand is strictly negative on a set of positive measure.

	\underline{Proof of (ii).}\quad
	Fix \(t_1<t<1\). For \(T_0(t)<0\), the limit as \(r\to\infty\) is strictly negative.
	By \eqref{0831-08}, a nonnegative value of \(H_2(r,t)\) would imply strict positivity at every larger radius. Thus \(H_2(\cdot,t)<0\).

	If \(T_0(t)=0\), \(s\mapsto s^{-(n+2)/2}P(r^2/(4s))\) decreases strictly, since it is a positive multiple of \(Q(r^2/(4s))\) from Lemma~\ref{lem-0424-01}{\rm(i)} in dimension \(n+2\).
	Since the weight \(s^{-1/2}\phi'(t-s)\) in \eqref{0831-07} has integral \(T_0(t)=0\), we may subtract the value of this function at \(s_0=t-t_1\) inside the integral.
	The product is strictly negative on each side of \(s_0\), so \(H_2(r,t)<0\).

	If \(T_0(t)>0\), the limits at zero and infinity have opposite signs since \(\phi'(t)<0\).
	By continuity of \(H_2(\cdot,t)\) and the intermediate value theorem, a zero exists, and \eqref{0831-08} gives uniqueness and the signs on either side.

	To verify simplicity at this zero, set \(z=r^2/(4s)\).
	Then \(\partial_r\ln K_r(s)=2zL'(z)/r\) increases strictly in \(s\), by the same Lemma~\ref{lem-0324-01}{\rm(v)}.
	At \(r=r_*(t)\), differentiating under the integral and subtracting the value at \(s_0\) give
	\begin{align*}
	    \partial_r H_2(r,t)=\frac{2}{n+3}\int_0^t\phi'(t-s)K_r(s)
	    [\partial_r\ln K_r(s)-\partial_r\ln K_r(s_0)]\,ds>0.
	\end{align*}

	\underline{Proof of (iii).}\quad
	By Remark~\ref{rmk-0730-03}, \(T_0>0\) on \((t_1,t_\infty)\), so part {\rm(ii)} gives a simple zero throughout this interval.
	Fix \(0<a<b\) and \(T>0\), and extend \(\phi'\) by zero to negative times.
For \(a\leq r\leq b\) and \(0<t\leq T\), \eqref{0514-02} bounds the integrand of \eqref{0831-07} by \(C(n,a)\|\phi'\|_{L^\infty(0,1)} s^{-1/2}\) on \(0<s<T\), so dominated convergence gives joint continuity of \(H_2\) on \((0,\infty)^2\).
	Fix \(t_0\in(t_1,t_\infty)\) and \(0<\varepsilon<r_*(t_0)\).
	By \eqref{item-0906-01},
	\begin{align*}
		H_2(r_*(t_0)-\varepsilon,t_0)<0< H_2(r_*(t_0)+\varepsilon,t_0).
	\end{align*}
By joint continuity, \(H_2(r_*(t_0)-\varepsilon,t)<0<H_2(r_*(t_0)+\varepsilon,t)\) for \(t\) sufficiently close to \(t_0\).
	Hence \eqref{item-0906-01} gives \(|r_*(t)-r_*(t_0)|<\varepsilon\), proving that \(r_*\) is continuous.
	For every fixed \(r>0\), joint continuity and \(H_2(r,t_1)>0>H_2(r,t_\infty)\) give \(r_*(t)<r\) for \(t>t_1\) close to \(t_1\), and \(r_*(t)>r\) for \(t<t_\infty\) close to \(t_\infty\).
	Since \(r>0\) is arbitrary, \eqref{0904-03} follows.

	We finally prove \eqref{261008}.
	By \eqref{0904-03}, \(r_*(t)\to\infty\) as \(t\to t_\infty^-\).
	Applying \eqref{0504-06} in \eqref{0831-07}, uniformly for \(0<s\leq t\) and \(t<t_\infty\) sufficiently close to \(t_\infty\), gives
	\begin{align}\label{H2-expansion}
		H_2(r,t) =\frac{2^{n+3}\Gamma(\frac{n+5}{2})} {(n+3)\Gamma(\frac32)}r^{-n-2} \left(T_0(t)-\frac{(n+2)T_1(t)}{r^2}+O(r^{-4})\right).
	\end{align}
	Setting \(r=r_*(t)\) yields
	\begin{align*}
		T_0(t)r_*(t)^2=(n+2)T_1(t)+O(r_*(t)^{-2}).
	\end{align*}
	Now \(T_1(t)\to T_1(t_\infty)>0\) and \(T_0(t)/(t_\infty-t)\to-T_0'(t_\infty)>0\), so \eqref{261008} follows.
\end{proof}

\begin{proof}[Proof of Proposition \ref{prop-0831-01}]
	Write \(r=|x'|\).
	For \(r\geq\sqrt2\), the factor \((r-y_1)^2-y_2^2\) is strictly positive wherever \(\psi(y')>0\) by \eqref{radial-geometric-factor}.
	These points all satisfy \(|y'|<1\), so
	\begin{align*}
	    r-1<|r\mathbf e_1'-y'|<r+1.
	\end{align*}

	\underline{Proof of (i).}\quad
	Lemma \ref{lem-0903-01}{\rm(i)} gives \(H_2(|r\mathbf e_1'-y'|,t)>0\).
	Hence \eqref{0825-13} gives \(\sgn(\partial_{x_i}p(x',0,t))=-\sgn(x_ix_j)\).

	\underline{Proof of (ii).}\quad
	For \(t_1<t<t_\infty\), set \(r_{\rm tan}(t):=r_*(t)\).
	If \(\sqrt2\leq r\leq r_{\rm tan}(t)-1\), then \(|r\mathbf e_1'-y'|<r_{\rm tan}(t)\) wherever \(\psi(y')>0\).
	If \(r\geq\max\{\sqrt2,r_{\rm tan}(t)+1\}\), then \(|r\mathbf e_1'-y'|>r_{\rm tan}(t)\) at all such points.
	Lemma \ref{lem-0903-01}{\rm(ii)} and \eqref{0825-13} give the signs in {\rm(ii)}.
	The continuity and \eqref{261046} follow from Lemma \ref{lem-0903-01}{\rm(iii)} and \eqref{261008}, respectively.

	\underline{Proof of (iii).}\quad
	By relabeling the tangential coordinates, it suffices to take \(j=1\) and \(i=2\).
	Write \(r=|x'|\), and use the quantities \(A(t)\), \(D(t)\), \(E(t)\), \(R_0(t)\), and \(\kappa_n\) from \eqref{260279} and Lemma \ref{lem-0911-04a}.
	Consider the vector whose \(k\)-th component is obtained from the right-hand side of \eqref{0822-10}, with \(j=1\), by replacing \(x_1-y_1\) with \(x_k-y_k\).
	Because \(\psi\) and the kernel multiplying \(x_1-y_1\) in \eqref{0822-10} are radial, this vector is equivariant under every orthogonal transformation of \(\Sigma\).
	It is therefore of the form \(x'g(r,t)\) for a scalar function \(g\), and its first component gives \(p(x',0,t)=x_1g(r,t)\) in the normalization of \eqref{0822-10}.
	For \(r>1\), the differentiations below are justified by Remark \ref{260142}.
	At \(x'=r\mathbf e_2'\),
	\begin{align*}
		\partial_{x_1}p(r\mathbf e_2',0,t)=g(r,t),
	\end{align*}
	whereas at an arbitrary \(x'\) with \(|x'|=r>1\),
	\begin{align}\label{261098}
		\partial_{x_2}p(x',0,t)=\frac{x_1x_2}{r}\partial_rg(r,t)
		=\frac{x_1x_2}{r}\partial_r\partial_{x_1}p(r\mathbf e_2',0,t).
	\end{align}

	Along the \(x_2\)-axis, the radial derivative of \(\partial_{x_1}p\) is
	\begin{align*}
		q_t(r):=\partial_r\partial_{x_1}p(r\mathbf e_2',0,t),\qquad r\geq\sqrt2.
	\end{align*}
	Evaluating \eqref{261098} at \(x'=r(\mathbf e_1'+\mathbf e_2')/\sqrt2\), where \(x_1x_2/r>0\), and applying part {\rm(ii)}, we obtain
	\begin{align}\label{261099}
		q_t(r)>0 & \quad\text{if }\sqrt2\leq r\leq r_{\rm tan}(t)-1,\nonumber\\
		q_t(r)<0 & \quad\text{if }r\geq\max\{\sqrt2,r_{\rm tan}(t)+1\}.
	\end{align}
	By Lemma \ref{lem-0903-01}{\rm(iii)}, \(r_{\rm tan}(t)\to\infty\) as \(t\to t_\infty^-\).
	Thus, for \(t<t_\infty\) sufficiently close to \(t_\infty\), both sign regions in \eqref{261099} are nonempty, and the continuity of \(q_t\) gives at least one zero in \((\sqrt2,\infty)\).
	Define
	\begin{align*}
		Z_t:=\{r\geq\sqrt2:q_t(r)=0\}.
	\end{align*}
	Equation \eqref{261099} shows that
	\begin{align}\label{261100}
		\varnothing\neq Z_t\subset(r_{\rm tan}(t)-1,r_{\rm tan}(t)+1).
	\end{align}
	Since \(q_t\) is continuous, \(Z_t\) is closed in \([\sqrt2,\infty)\), hence compact by \eqref{261100}, and
	\begin{align*}
		\inf_{r\in Z_t}r\geq r_{\rm tan}(t)-1\longrightarrow\infty
		\qquad\text{as }t\to t_\infty^-.
	\end{align*}

	Choose \(\delta_0\in(0,t_\infty-t_1)\).
	On the compact time interval \([t_\infty-\delta_0,t_\infty]\), \(E(t)\) is bounded and \(R_0(t)\leq2+C_n\sqrt{t_\infty}\).
	Hence, for \(t<t_\infty\) sufficiently close to \(t_\infty\), every \(r\in Z_t\) satisfies \(r\geq R_0(t)\), so \eqref{0911-601}--\eqref{0911-60c} apply at every such zero.

	For \(r\in Z_t\), \eqref{0911-601} gives
	\begin{align}\label{261101}
		nA(t)r^2=-(n+2)D(t)+r^2e_1(r,t).
	\end{align}
	Substituting \eqref{261101} into \eqref{0911-602} gives the exact identity
	\begin{align}\label{261102}
		\frac{r^{n+4}}{\kappa_n}q_t'(r)=2(n+2)D(t)+r^2\bigl[(n+1)e_1(r,t)+e_2(r,t)\bigr].
	\end{align}
	Integration by parts in the definition \eqref{0323-02} gives
	\begin{align*}
		T_1(t)=\frac12\int_0^t\frac{\phi(s)}{\sqrt{t-s}}\,ds>0,\qquad 0<t<1.
	\end{align*}
	Since \(T_0(t_\infty)=0\), it follows that
	\begin{align*}
		D(t)\longrightarrow-nM_0T_1(t_\infty)<0 \qquad\text{as }t\to t_\infty^-.
	\end{align*}
	Moreover, \eqref{0911-60c} and \eqref{261100} give
	\begin{align*}
		\sup_{r\in Z_t}\left|\frac{r^{n+4}}{\kappa_n}q_t'(r)-2(n+2)D(t)\right|
		\leq C(n)E(t)\left(\inf_{r\in Z_t}r\right)^{-2}\longrightarrow0.
	\end{align*}
	It follows that
	\begin{align*}
		\frac{r^{n+4}}{\kappa_n}q_t'(r)\longrightarrow-2n(n+2)M_0T_1(t_\infty)<0
	\end{align*}
	uniformly for \(r\in Z_t\) as \(t\to t_\infty^-\).
	Hence, for \(t<t_\infty\) sufficiently close to \(t_\infty\), \(q_t'(r)<0\) at every \(r\in Z_t\).

	We now prove that \(Z_t\) has only one element.
	Every zero is isolated because \(q_t'\) is negative there; since \(Z_t\) is compact, \(Z_t\) is therefore finite.
	If it had at least two elements, let \(a<b\) be consecutive zeros.
	The inequalities \(q_t'(a)<0\) and \(q_t'(b)<0\) imply that \(q_t<0\) immediately to the right of \(a\) and \(q_t>0\) immediately to the left of \(b\).
	The intermediate value theorem would then give another zero in \((a,b)\), contradicting the choice of \(a\) and \(b\).
	Thus \(Z_t=\{\rho_{\rm tan}(t)\}\).
	Together with \eqref{261099}, uniqueness gives
	\begin{align*}
		q_t(r)>0\quad\text{for }\sqrt2\leq r<\rho_{\rm tan}(t),
		\qquad q_t(r)<0\quad\text{for }r>\rho_{\rm tan}(t).
	\end{align*}
	Using \eqref{261098} and then relabeling the two tangential coordinates gives
	\begin{align*}
		\sgn\bigl(\partial_{x_i}p(x',0,t)\bigr)=-\sgn(x_ix_j)\sgn(|x'|-\rho_{\rm tan}(t)),\qquad |x'|\geq\sqrt2.
	\end{align*}
	On each ray with \(x_ix_j\ne0\), the factor \(x_ix_j/r\) is nonzero and \(q_t'(\rho_{\rm tan}(t))<0\), so the zero at \(|x'|=\rho_{\rm tan}(t)\) on this ray is simple.

	To estimate \(\rho_{\rm tan}(t)-r_{\rm tan}(t)\), abbreviate \(\rho=\rho_{\rm tan}(t)\) and \(r_*=r_{\rm tan}(t)\).
	The expansion \eqref{H2-expansion} and \eqref{261101}, with \eqref{260279} and \eqref{0911-60c}, give
	\begin{align*}
	    T_0(t)r_*^2&=(n+2)T_1(t)+O(r_*^{-2}),\\
	    T_0(t)\rho^2&=(n+2)T_1(t)
	      -\frac{(n+2)M_2}{2(n-1)M_0}T_0(t)+O(\rho^{-2}).
	\end{align*}
	As \(t\uparrow t_\infty\), \eqref{261100} gives \(|\rho-r_*|<1\) and both radii tend to infinity.
	Since \(T_1(t_\infty)>0\), the two identities give \(r_*^{-2}+\rho^{-2}=O(T_0(t))\).
	Subtracting and dividing by \(T_0(t)>0\) yields \(\rho^2-r_*^2=O(1)\).
	Dividing by \(\rho+r_*\) gives the claimed \(O(r_*^{-1})\) bound, and the asymptotic formula for \(\rho\) follows from \eqref{261046}.

	\underline{Proof of (iv).}\quad
	For \(t\geq1\), Lemma~\ref{lem-0903-01}{\rm(i)} gives \(H_2(\cdot,t)<0\) under {\rm\textbf{(T1)}}.
	Under {\rm\textbf{(T2)}}--{\rm\textbf{(T4)}}, Remark~\ref{rmk-0730-03} gives \(T_0(t)\leq0\) for \(t_\infty\leq t<1\), so Lemma~\ref{lem-0903-01}{\rm(ii)} gives strict negativity there.
	Then \eqref{0825-13} gives \(\sgn(\partial_{x_i}p(x',0,t))=\sgn(x_ix_j)\).
\end{proof}

\appendix

\section{Bessel functions and radial Fourier inversion}
We recall the formulas used below from \cite[\S\S10.2 and 10.6]{DLMF}.
\subsection{Bessel function of the first kind}
For \(\nu\in\mathbb C\), the Bessel function of the first kind is defined, after fixing a branch of \(z^\nu\), by
\begin{align}\label{0401-05}
	J_\nu(z) = \sum_{k=0}^{\infty}\frac{(-1)^k}{k!\Gamma(k+\nu+1)} \left(\frac{z}{2}\right)^{2k+\nu}.
\end{align}
At a nonpositive integer \(k+\nu+1\), the coefficient is zero by the reciprocal gamma convention \cite[\S5.2(i)]{DLMF}.
With the principal branch, \(J_\nu\) is analytic on \(\mathbb C\setminus(-\infty,0]\).
If \(\nu\in\mathbb Z\), it extends to an entire function.
Its value is unambiguous on \((0,\infty)\).
For \(r,\rho>0\),
\begin{align}\label{0401-06}
	\partial_r(r^{-\nu}J_{\nu}(r\rho))= -\rho r^{-\nu}J_{\nu+1}(r\rho).
\end{align}

\subsection{Radial Fourier inversion}
For \(f,h\in L^1(\R^d)\), we use the convention
\begin{align*}
	\widehat f(\xi)=\int_{\R^d}e^{-ix\cdot\xi}f(x)\,dx, \qquad \mathcal{F}^{-1}h(x)=\frac{1}{(2\pi)^d}\int_{\R^d}e^{ix\cdot\xi}h(\xi)\,d\xi.
\end{align*}
Polar coordinates give the following formula; see \cite[Sect.~B.5]{Grafakos} or \cite[Sect.~IV.1]{SteinWeiss}.
\begin{lemma}
	Let \(d\geq1\), and let \(h\in L^1(\R^d)\) be radial, with \(h(\xi)=h_0(|\xi|)\) almost everywhere.
	Then, for every \(x\in\R^d\setminus\{0\}\),
	\begin{align}\label{0401-07}
		\mathcal{F}^{-1}h(x)=\frac{1}{(2\pi)^{\frac d2}|x|^{\frac{d-2}{2}}} \int_0^\infty h_0(\rho)\rho^{\frac d2} J_{\frac{d-2}{2}}(|x|\rho)\,d\rho.
	\end{align}
\end{lemma}

\section{Kernels involving the confluent hypergeometric functions}\label{260211}

We prove the lemmas stated in Subsection \ref{260031}.

\subsection{Proof of Lemma \ref{lem-0424-01}}\label{260212}
Let \(Q\) be given in \eqref{0827-10}.

\underline{Proof of (i).}\quad
Since \({}_1F_1(a,b;0)=1\), definition \eqref{0827-10} gives \(Q(0)=0\).
	For \(u>0\), substitute \(s=uv\) in the integral representation \eqref{0324-14} to obtain
\begin{align}\label{0615-03}
	Q(u) = \frac{1}{\Gamma\left(\frac{n}{2}\right)} \int_0^u s^{\frac{n}{2}-1} \sqrt{1-\frac{s}{u}} e^{-s}\,ds.
\end{align}
Extending the integrand in \eqref{0615-03} by zero to \((0,\infty)\), dominated convergence gives \(Q(u)\to1\) as \(u\to\infty\).
The function \(Q\) is strictly increasing because both the positive integrand and the interval of integration in \eqref{0615-03} increase with \(u\).

\underline{Proof of (ii).}\quad
The integral representation \eqref{0615-03} gives
\begin{align*}
	Q(u) = \frac{u^{\frac{n}{2}}}{\Gamma\left(\frac{n}{2}\right)}\int_0^1 e^{-ut}t^{\frac{n-2}{2}}(1-t)^{\frac12}\,dt < \frac{u^{\frac{n}{2}}}{\Gamma\left(\frac{n}{2}\right)}\int_0^1 t^{\frac{n-2}{2}}(1-t)^{\frac12}\,dt= \frac{\sqrt{\pi}}{2\Gamma\left(\frac{n+3}{2}\right)}u^{\frac{n}{2}}.
\end{align*}

\underline{Proof of (iii).}\quad
Write
\begin{align*}
	1-\frac{n}{4u}-Q(u)=\frac{1}{\Gamma\left(\frac n2\right)} \int_0^\infty s^{\frac n2-1}h\left(\frac{s}{u}\right)e^{-s}\,ds,
\end{align*}
where \(h(x):=1-\frac x2-\sqrt{\max\{1-x,0\}}\).
A direct computation gives \(|h(x)|\leq\frac{x^2}{2}\) for every \(x\geq0\), and hence
\begin{align*}
	\left| 1 - \frac{n}{4u} - Q(u) \right| & \leq \frac{1}{\Gamma\left(\frac{n}{2}\right)} \int_0^\infty s^{\frac{n}{2}-1} \frac{s^2}{2u^2} e^{-s}\,ds =\frac{n(n+2)}{8u^2}.
\end{align*}

\subsection{Proof of Lemma \ref{lem-0324-01}}\label{260213}

Let \(P\) and \(L\) be given in \eqref{0827-11}.
\eqref{0327-03} and the identities \(P'=L'P\) and \(P''=(L''+(L')^2)P\) give
\begin{align}\label{0327-06}
	uL''(u)+u(L'(u))^2+\left(\frac{n+3}{2}+u\right)L'(u)+\frac n2=0.
\end{align}
\begin{lemma}\label{lem-0327-01}
	Let \(r\in(0,\infty)\), and let \(\{f_k\}_{k=0}^\infty\) and \(\{g_k\}_{k=0}^\infty\) be sequences of real numbers such that \(g_k>0\) for all \(k\geq0\) and \(\{f_k/g_k\}_{k=0}^\infty\) is strictly decreasing.
	Suppose that \(f(u)=\sum_{k=0}^\infty f_ku^k\) and \(g(u)=\sum_{k=0}^\infty g_ku^k\) have radii of convergence at least \(r\).
	Then \(f(u)/g(u)\) is strictly decreasing for \(0<u<r\).
\end{lemma}
\begin{proof}
	It suffices to show that \(f'(u)g(u)-f(u)g'(u)<0\) for \(0<u<r\).
	We have
	\begin{align*}
		f'(u)g(u) - f(u)g'(u) 
		 & = \sum_{m=1}^{\infty}\sum_{\substack{k+l=m\\
		 k>l\geq 0}}(k-l)g_kg_{l}\left( \frac{f_k}{g_k}- \frac{f_{l}}{g_{l}}\right)u^{m-1}.
	\end{align*}
	Every term on the right-hand side is negative because \(g_k>0\) and \(\frac{f_k}{g_k}-\frac{f_l}{g_l}<0\) for \(k>l\).
	Hence \((f/g)'<0\) on \((0,r)\).
\end{proof}

\underline{Proof of (i).}\quad
\(P(0)=1\) follows from \eqref{260030}.
On the other hand, \eqref{0327-07} gives
\begin{align}\label{0327-09}
	P'(u) = - \frac{n}{n+3}{}_1F_1\left( \frac{n+2}{2}, \frac{n+5}{2};-u\right)<0.
\end{align}

Moreover, \eqref{0327-01} and \eqref{0324-14} give
\begin{align}\label{0327-08}
	P'(u)= - \frac{n}{2u}P(u) + \frac{n}{2u}{}_1F_1\left( \frac{n+2}{2}, \frac{n+3}{2};-u\right) >- \frac{n}{2u}P(u).
\end{align}
Therefore \(\frac{n+1}{2}P(u)+uP'(u)>\frac12P(u)>0\), as required.

\underline{Proof of (ii).}\quad
Set \(f(u):=(1+2u)^{n/2}P(u)\).
Then \eqref{0327-09}, \eqref{0327-08}, and \eqref{0327-02} with \(a=\frac{n+2}{2}\), \(b=\frac{n+3}{2}\), and \(z=-u\) give
\begin{align*}
	\begin{split}
		f'(u) & = (1+2u)^{\frac{n-2}{2}}(nP(u) + (1+2u)P'(u))\\
		& = (1+2u)^{\frac{n-2}{2}}\frac{n(n+2)}{n+3}{}_1F_1\left( \frac{n+4}{2}, \frac{n+5}{2}; -u\right)>0.
	\end{split}
\end{align*}

By \eqref{0327-10}, \(\lim_{u\to\infty}f(u)= \frac{2^{(n+2)/2}\Gamma((n+3)/2)}{\sqrt\pi}\).
Together with \(f(0)=1\), this proves \eqref{0529-03}.

\underline{Proof of (iii).}\quad
Using \eqref{0327-09},
\begin{align*}
	L'(u) = \frac{P'(u)}{P(u)} = - \frac{n}{n+3}\frac{{}_1F_1\left(\frac{n+2}{2}, \frac{n+5}{2};-u\right)}{{}_1F_1\left(\frac{n}{2}, \frac{n+3}{2};-u\right)}<0.
\end{align*}
To prove the inequality \eqref{0529-02}, we further show that \((1+2u)L'(u)\) is strictly decreasing.
\eqref{0327-01}, \eqref{0327-07}, and \eqref{0327-02} give
\begin{align}\label{0615-12}
	(1+2u){}_1F_1\left(\frac{n+2}{2}, \frac{n+5}{2};-u\right) = -(n+2){}_1F_1\left(\frac{n+4}{2}, \frac{n+5}{2};-u\right)+(n+3){}_1F_1\left(\frac{n}{2}, \frac{n+3}{2};-u\right).
\end{align}

Then \eqref{0615-12} and \eqref{0327-04} give
\begin{align*}
	(1+2u)L'(u) & = n\left( \frac{n+2}{n+3}\frac{{}_1F_1\left(\frac{n+4}{2}, \frac{n+5}{2};-u\right)}{{}_1F_1\left( \frac{n}{2}, \frac{n+3}{2}; -u\right)}-1\right)\\
	 & = n\left( \frac{n+2}{n+3}\frac{{}_1F_1\left(\frac{1}{2}, \frac{n+5}{2};u\right)}{{}_1F_1\left( \frac{3}{2}, \frac{n+3}{2}; u\right)}-1\right)= n\left( \frac{n+2}{n+3}\frac{f(u)}{g(u)}-1\right),
\end{align*}
where \(f(u):= {}_1F_1(\frac12, \frac{n+5}{2};u)\) and \(g(u):= {}_1F_1(\frac32, \frac{n+3}{2};u)\).
Every coefficient of \(g\) is positive, and the ratio
\begin{align*}
	\frac{\left(\frac12\right)_k/\left(\frac{n+5}{2}\right)_k}{\left(\frac32\right)_k/\left(\frac{n+3}{2}\right)_k} =\frac{\frac12}{\frac12+k}\frac{\frac{n+3}{2}}{\frac{n+3}{2}+k}
\end{align*}
is strictly decreasing in \(k\).
Lemma \ref{lem-0327-01} therefore shows that \(f(u)/g(u)\) is strictly decreasing.
Applying \eqref{0327-10} separately to the numerator and denominator in the formula for \(L'\) gives
\begin{align*}
	\frac{{}_1F_1\left(\frac{n+2}{2},\frac{n+5}{2};-u\right)} {{}_1F_1\left(\frac n2,\frac{n+3}{2};-u\right)} =\frac{n+3}{2u}\left(1+O(u^{-1})\right), \qquad L'(u)=-\frac{n}{2u}\left(1+O(u^{-1})\right).
\end{align*}
Thus \(\lim_{u\to\infty}(1+2u)L'(u)=-n\).
This limit, \(L'(0)=-\frac{n}{n+3}\), and the strict monotonicity of \((1+2u)L'(u)\) give \eqref{0529-02}.

\underline{Proof of (iv).}\quad
We find that
\begin{align*}
	L''(u) & = \frac{P''(u)P(u)-(P'(u))^2}{P(u)^2}.
\end{align*}
Set \(w_0(v):=v^{(n-2)/2}(1-v)^{1/2}\).
The numerator on the right-hand side is positive since
\begin{align*}
	\int_0^1 v^2e^{-uv}w_0(v)\,dv\int_0^1 e^{-uv}w_0(v)\,dv - \left( \int_{0}^1 ve^{-uv}w_0(v)\,dv\right)^2>0,
\end{align*}
by the strict Cauchy--Schwarz inequality.

Since \((1+2u)L'(u)\) is strictly decreasing, \eqref{0529-02} gives
\begin{align*}
	L''(u)=\frac{((1+2u)L'(u))'-2L'(u)}{1+2u} \leq  -\frac{2L'(u)}{1+2u}<\frac{2n}{(1+2u)^2}.
\end{align*}
This gives \eqref{0602-01}.

\underline{Proof of (v).}\quad
Using \eqref{0327-07} and \eqref{0324-14}, we have
\begin{align*}
	L'(u) & = -\frac{n}{n+3}\frac{{}_1F_1\left(\frac{n+2}{2}, \frac{n+5}{2};-u\right)}{{}_1F_1\left(\frac{n}{2}, \frac{n+3}{2};-u\right)}= -\frac{\int_0^1 e^{-ut}t^{\frac{n}{2}}(1-t)^{\frac12}\,dt}{\int_0^1 e^{-ut}t^{\frac{n-2}{2}}(1-t)^{\frac12}\,dt}.
\end{align*}
By integration by parts,
\begin{align*}
	u\int_0^1 e^{-ut}t^{\frac{n}{2}}(1-t)^{\frac12}\,dt &= - \int_0^1 \partial_t(e^{-ut})t^{\frac{n}{2}}(1-t)^{\frac12}\,dt \\
	&= \frac{n}{2}\int_0^1 e^{-ut}t^{\frac{n-2}{2}}(1-t)^{\frac12}\,dt - \frac12 \int_0^1 e^{-ut}t^{\frac{n}{2}}(1-t)^{-\frac12}\,dt.
\end{align*}
Hence
\begin{align}\label{0622-04}
	uL'(u) = -\frac{n}{2} + \frac12 \frac{\int_0^1 e^{-ut}t^{\frac{n}{2}}(1-t)^{-\frac12}\,dt}{\int_0^1 e^{-ut}t^{\frac{n-2}{2}}(1-t)^{\frac12}\,dt},
\end{align}
and symmetrization gives
\begin{align}
	(uL'(u))' & = \frac12 \frac{d}{du}\left(\frac{\int_0^1 e^{-ut}t^{\frac{n}{2}}(1-t)^{-\frac12}\,dt}{\int_0^1 e^{-ut}t^{\frac{n-2}{2}}(1-t)^{\frac12}\,dt} \right) \nonumber\\
	& = \frac12 \frac{\int_0^1 -te^{-ut}t^{\frac{n}{2}}(1-t)^{-\frac12}\,dt \int_0^1 e^{-us}s^{\frac{n-2}{2}}(1-s)^{\frac12}\,ds - \int_0^1 e^{-ut}t^{\frac{n}{2}}(1-t)^{-\frac12}\,dt \int_0^1 -s e^{-us}s^{\frac{n-2}{2}}(1-s)^{\frac12}\,ds}{\left( \int_0^1 e^{-ut}t^{\frac{n-2}{2}}(1-t)^{\frac12}\,dt\right)^2}\nonumber \\
	& = -\frac14 \frac{\int_0^1\int_0^1 e^{-u(t+s)}(ts)^{\frac{n-2}{2}}(1-t)^{-\frac12}(1-s)^{-\frac12}(t-s)^2\,dtds}{\left(\int_0^1 e^{-ut}t^{\frac{n-2}{2}}(1-t)^{\frac12}\,dt\right)^2}=-\frac14 \frac{B(u)}{A(u)^2},\label{0619-09}
\end{align}
where
\begin{align}\label{0619-090}
	A(u):=\int_0^1 e^{-ut}t^{\frac{n-2}{2}}(1-t)^{\frac12}\,dt,\qquad B(u):=\int_0^1\!\!\int_0^1 e^{-u(t+s)}(ts)^{\frac{n-2}{2}}(1-t)^{-\frac12}(1-s)^{-\frac12}(t-s)^2\,dtds.
\end{align}
Since \(B(u)>0\), \eqref{0619-09} already gives
\begin{align}\label{0731-33}
	\theta'(u)<0\qquad(u>0),\qquad \theta(u):=uL'(u).
\end{align}

It remains to prove that \(h(u):=(u+a_n)^2\theta'(u)\) is strictly increasing on \((0,\infty)\).
The bounds \eqref{0529-04} then follow from the limits computed below.
Set \(a:=a_n=(n+3)(n+5)/6\).

Multiplying \eqref{0327-06} by \(u\) and using \(u^2L''=u\theta'-\theta\), \(u^2(L')^2=\theta^2\), and \(uL'=\theta\) gives
\begin{align}\label{0731-34}
	u\theta'=-\theta\left(\theta+\tfrac{n+1}{2}\right)-u\left(\theta+\tfrac n2\right).
\end{align}
Differentiating \eqref{0731-34} once and twice yields
\begin{align}
	u\theta'' & =-\left(2\theta+\tfrac{n+3}{2}+u\right)\theta'-\left(\theta+\tfrac n2\right),\qquad u\theta'''=-2\theta'(\theta'+1)-\left(2\theta+\tfrac{n+5}{2}+u\right)\theta''.\label{0731-36}
\end{align}

By \eqref{0529-02}, \(\theta=uL'<0\), while \eqref{0622-04} gives \(\theta+\frac n2>0\).
Moreover, \eqref{0602-01} shows that \(L'\) is increasing, so \(\theta(u)=uL'(u)>uL'(0)=-\frac{n}{n+3}u\).
Combining these inequalities gives
\begin{align}\label{0731-37}
	-\min\left\{\tfrac n2,\ \tfrac{n}{n+3}u\right\}<\theta(u)<0\qquad \textrm{for } u>0.
\end{align}

By \eqref{0324-14}, \( P(u)=\frac{A(u)}{A(0)} \) with \(A\) given in \eqref{0619-090}, so \(\frac{d^k}{du^k}L=\frac{d^k}{du^k}(\ln A)\) for every \(k\geq1\).
For \(k\geq0\), set
\begin{align*}
	\mu_k:=\frac{(-1)^kA^{(k)}(0)}{A(0)} =\frac{\displaystyle\int_0^1t^{\,k+\frac{n-2}{2}}(1-t)^{\frac12}\,dt} {\displaystyle\int_0^1t^{\frac{n-2}{2}}(1-t)^{\frac12}\,dt} =\frac{\left(\frac n2\right)_k}{\left(\frac{n+3}{2}\right)_k}.
\end{align*}
Direct computation gives \(L'(0)=-\kappa_1=-\mu_1<0\), \(L''(0)=\kappa_2=\mu_2-\mu_1^2>0\), and \(L'''(0)=-\kappa_3=-\mu_3+3\mu_1\mu_2-2\mu_1^3\).
Hence
\begin{align}\label{0731-38}
	\kappa_1=\frac{n}{n+3},\qquad \kappa_2=\frac{6n}{(n+3)^2(n+5)},\qquad \kappa_3=-\frac{24\,n\,(n-3)}{(n+3)^3(n+5)(n+7)}.
\end{align}
These values also give
\begin{align}\label{0731-39}
	\frac{\kappa_1}{\kappa_2}=\frac{(n+3)(n+5)}{6}=a, \qquad 2\kappa_2^2-\kappa_1\kappa_3=\frac{24\,n^2(n+2)}{(n+3)^3(n+5)^2(n+7)}>0.
\end{align}

We now define
\begin{align}\label{0731-40}
	\mathcal{A}(u):=2\theta'(u)+(u+a)\theta''(u),\qquad \mathcal{B}(u):=\frac{3u}{u+a}+\left(a-\tfrac{n+5}{2}\right)-2\theta(u)+(u+a)\theta'(u),
\end{align}
so that \(h'(u)=(u+a)\mathcal{A}(u)\).
By \eqref{0731-38}--\eqref{0731-39}, direct computation gives
\begin{align}\label{0731-41}
	\mathcal{A}(0)=-2\kappa_1+2a\kappa_2=0,\qquad \mathcal{A}'(0)=6\kappa_2-3a\kappa_3=\frac{3}{\kappa_2}\bigl(2\kappa_2^2-\kappa_1\kappa_3\bigr)>0,
\end{align}
and, using \(1-\kappa_1=\frac{3}{n+3}\) and \(a\cdot\frac{3}{n+3}=\frac{n+5}{2}\),
\begin{align}\label{0731-42}
	\mathcal{B}(0)=a(1-\kappa_1)-\tfrac{n+5}{2}=0,\qquad \mathcal{B}'(0)=\tfrac3a+\kappa_1+2a\kappa_2>0.
\end{align}

Since \(\mathcal{A}'(u)=3\theta''+(u+a)\theta'''\), multiplying by \(u\) and using \eqref{0731-36} gives
\begin{align}\label{0731-43}
	u\,\mathcal{A}'(u)=\theta''\left[\,3u-(u+a)\left(2\theta+\tfrac{n+5}{2}+u\right)\right]-2(u+a)\theta'(\theta'+1).
\end{align}
If \(\mathcal{A}(u_1)=0\) for some \(u_1>0\), then \(\theta''(u_1)=-\frac{2\theta'(u_1)}{u_1+a}\).
Substituting this identity into \eqref{0731-43} gives
\begin{align}\label{0731-44}
	u_1\,\mathcal{A}'(u_1)=-2\theta'(u_1)\,\mathcal{B}(u_1).
\end{align}

Multiplying \eqref{0731-40} by \(u(u+a)>0\) and eliminating \(u\theta'\) by \eqref{0731-34} gives \(u(u+a)\mathcal{B}(u)=\widetilde{\mathcal{B}}\bigl(u,\theta(u)\bigr)\), where
\begin{align*}
	\widetilde{\mathcal{B}}(u,y):=3u^2+\left(a-\tfrac{n+5}{2}\right)u(u+a)-2yu(u+a) -(u+a)^2\left[y\left(y+\tfrac{n+1}{2}\right)+u\left(y+\tfrac n2\right)\right],
\end{align*}
a quadratic polynomial in \(y\) with leading coefficient \(-(u+a)^2<0\).
Set
\begin{align*}
	G_0(u,y):=u\,\partial_u\widetilde{\mathcal{B}}(u,y) -\partial_y\widetilde{\mathcal{B}}(u,y)\left[y\left(y+\tfrac{n+1}{2}\right)+u\left(y+\tfrac n2\right)\right],
\end{align*}
so that, by the chain rule and \eqref{0731-34}, when \(y=\theta(u)\),
\begin{align}\label{0731-46}
	u\,\frac{d}{du}\Bigl[u(u+a)\mathcal{B}(u)\Bigr]=G_0\bigl(u,\theta(u)\bigr).
\end{align}

We now divide \(G_0\), a cubic in \(y\), by \(\widetilde{\mathcal{B}}\), a quadratic in \(y\).
Set
\begin{align}\label{0806-01}
	d_0(u) & :=n^2+8n+6u+15>0,\nonumber \\
	q(u,y) & :=(n+1)(n+3)(n+5)+2\left(n^2+11n-6\right)u+4(n+3)(n+5)\,y+12u(u+2y),\nonumber\\
	p_1(u) & :=36u^2+12\left(n^2+8n+33\right)u+(n+3)(n+5)\left(n^2+8n-21\right),\\
	p_0(u) & :=36\left(n^2+14n+18\right)u^2+12\left(n^4+22n^3+127n^2+264n+54\right)u\nonumber \\
	& \qquad+(n+3)(n+5)^2\left(n^3+17n^2+24n+54\right),\nonumber
\end{align}
and
\begin{align}\label{0806-02}
	r_1(u):=\frac{u\,p_1(u)}{36},\qquad r_0(u):=\frac{u^2p_0(u)}{36\,d_0(u)}.
\end{align}
Direct expansion gives
\begin{align*}
	G_0(u,y)=-\frac{q(u,y)}{2d_0(u)}\,\widetilde{\mathcal{B}}(u,y)+r_1(u)\,y+r_0(u).
\end{align*}

We claim that
\begin{align}\label{0731-48}
	r_0(u)>0\ \ \text{for }\,u>0,\qquad r_0(u)-\tfrac{n}{n+3}u\,r_1(u)>0\ \ \text{for }\,0<u\leq\tfrac{n+3}{2},\qquad r_0(u)-\tfrac n2\,r_1(u)>0\ \ \text{for }\,u\geq\tfrac{n+3}{2}.
\end{align}
The first is immediate from \eqref{0806-01}--\eqref{0806-02}, since every coefficient of \(p_0\) is a polynomial in \(n\) with positive coefficients and \(d_0>0\).
For the second, substituting \(u=\frac{n+3}{2}t\) with \(t\in(0,1]\) gives the identity
\begin{align*}
	r_0\!\left(\tfrac{n+3}{2}t\right)-\tfrac{n}{n+3}\cdot\tfrac{n+3}{2}t\; r_1\!\left(\tfrac{n+3}{2}t\right)=\frac{t^2(n+3)^2}{48\,(n+3t+5)}\,q_0(t),
\end{align*}
where
\begin{align}\label{0806-05}
	q_0(t):={} & 3(n+2)(n+3)(n+5)^2(1-t)^3 +2\left(4n^4+70n^3+409n^2+984n+729\right)t(1-t)^2\nonumber\\
	& +\left(7n^4+139n^3+886n^2+2208n+1728\right)t^2(1-t) +2\left(n^4+22n^3+151n^2+384n+360\right)t^3.
\end{align}

For the third, substituting \(u=\frac{n+3}{2}+w\) with \(w\geq0\) gives
\begin{align*}
	r_0\!\left(\tfrac{n+3}{2}+w\right)-\tfrac n2\,r_1\!\left(\tfrac{n+3}{2}+w\right) =\frac{\widetilde p_0(w)}{72\left(n^2+11n+24+6w\right)},
\end{align*}
where
\begin{align}\label{0806-07}
	\widetilde p_0(w):={} & 72(n+2)(n+9)\,w^4 +12\left(2n^4+47n^3+350n^2+897n+756\right)w^3\nonumber\\
	& +2(n+3)\left(n^5+36n^4+471n^3+2618n^2+5784n+5238\right)w^2\nonumber \\
	& +(n+3)^2\left(n^5+36n^4+471n^3+2600n^2+5856n+5616\right)w\\
	& +3(n+3)^3\left(n^4+22n^3+151n^2+384n+360\right).\nonumber
\end{align}
Every coefficient in \eqref{0806-05} and \eqref{0806-07} is positive for \(n\geq2\).
Hence \(q_0>0\) on \([0,1]\) and \(\widetilde p_0>0\) for \(w\geq0\), which gives the remaining inequalities in \eqref{0731-48}.

Suppose that \(u>0\), \(\widetilde{\mathcal{B}}(u,y)=0\), and \(y\in[-\min\{\frac n2,\frac{n}{n+3}u\},0]\).
At such a zero, \(G_0(u,y)=r_1(u)y+r_0(u)\).
The affine expression on the right-hand side is positive at both endpoints of the stated interval by \eqref{0731-48}, hence \(G_0(u,y)>0\).

We now prove \(\mathcal{B}>0\) on \((0,\infty)\).
Set \(\widehat{\mathcal{B}}(u):=u(u+a)\mathcal{B}(u)\), which is smooth on \([0,\infty)\).
By \eqref{0731-42}, \(\widehat{\mathcal{B}}>0\) on some interval \((0,\varepsilon)\).
If \(\widehat{\mathcal{B}}\) had a zero in \((0,\infty)\), let \(u_2\) be the smallest one.
Then \(\widehat{\mathcal{B}}(u_2)=0\) and \(\widehat{\mathcal{B}}>0\) on \((0,u_2)\), so \(\widehat{\mathcal{B}}{}'(u_2)\leq0\).
On the other hand \(\widehat{\mathcal{B}}(u_2)=0\) means \(\widetilde{\mathcal{B}}(u_2,\theta(u_2))=0\), and \eqref{0731-37} gives
\begin{align*}
	\theta(u_2)\in \left[-\min\left\{\frac n2,\frac{n}{n+3}u_2\right\},0\right].
\end{align*}
The preceding paragraph and \eqref{0731-46} therefore give \(\widehat{\mathcal{B}}{}'(u_2)= \frac{G_0(u_2,\theta(u_2))}{u_2}>0\), a contradiction.
Thus \(\mathcal{B}>0\) on \((0,\infty)\).

By \eqref{0731-41}, \(\mathcal{A}>0\) on some \((0,\varepsilon)\).
If \(\mathcal{A}\) had a zero in \((0,\infty)\), let \(u_1\) be the smallest one.
Then \(\mathcal{A}'(u_1)\leq0\), while \eqref{0731-44}, \eqref{0731-33}, and \(\mathcal{B}(u_1)>0\) give \(u_1\mathcal{A}'(u_1)=-2\theta'(u_1)\mathcal{B}(u_1)>0\), a contradiction.
Hence \(\mathcal{A}>0\) on \((0,\infty)\) and \(h'(u)=(u+a)\mathcal{A}(u)>0\).

Since \(\theta'(0)=L'(0)=-\frac{n}{n+3}\), we obtain \(h(0)=a^2\theta'(0)=-\frac{n(n+3)(n+5)^2}{36}\).
For \(u\to\infty\), \eqref{0327-07} and \eqref{0504-06} give
\begin{align*}
	P'(u)=-\frac{n}{n+3}{}_1F_1\left(\frac{n+2}{2},\frac{n+5}{2};-u\right), \qquad P''(u)=\frac{n(n+2)}{(n+3)(n+5)} {}_1F_1\left(\frac{n+4}{2},\frac{n+7}{2};-u\right).
\end{align*}
Applying \eqref{0504-06} separately to \(P\), \(P'\), and \(P''\) gives
\begin{align*}
	\frac{P'(u)}{P(u)} =-\frac{n}{2u}+\frac{n}{4u^2}+O(u^{-3}), \qquad \frac{P''(u)}{P(u)} =\frac{n(n+2)}{4u^2}-\frac{n(n+2)}{4u^3}+O(u^{-4}).
\end{align*}
Consequently,
\begin{align*}
	L'(u)+uL''(u) =\frac{P'(u)}{P(u)} +u\frac{P''(u)}{P(u)} -u\left(\frac{P'(u)}{P(u)}\right)^2 =-\frac{n}{4u^2}+O(u^{-3}).
\end{align*}
It follows that
\begin{align*}
	h(u)=\left(1+\frac au\right)^2u^2\bigl(L'(u)+uL''(u)\bigr) \longrightarrow-\frac n4.
\end{align*}
Since \(h\) is strictly increasing and continuous on \([0,\infty)\), we conclude \(-\frac{n(n+3)(n+5)^2}{36}<h(u)<-\frac n4\) for every \(u>0\), which is \eqref{0529-04}.
Applying \eqref{0529-04} to the positive quantities \(-(uL'(u))'\) and \(-(vL'(v))'\) gives \eqref{0731-32}.

\underline{Proof of (vi).}\quad
By \eqref{0731-33}, \(\theta'(u)<0\), and part {\rm(v)} gives \(\mathcal A(u)>0\).
Solving \eqref{0731-40} for \(\theta''\) therefore gives
\begin{align}\label{0327-12}
	\bigl(u^2L''(u)\bigr)'=u\theta''(u)
	=\frac{u}{u+a}\bigl(\mathcal A(u)-2\theta'(u)\bigr)>0,\qquad u>0.
\end{align}

\subsection{Proof of Lemma \ref{lem-0618-01}}\label{260214}
Let \(\mathcal{K}\) and \(\mathcal{H}\) be given in \eqref{0827-12} and \eqref{0827-13}.

Throughout the proof, we use the following inequality for \(m>0\) and \(k>-1\).
\begin{align}\label{261003}
	\int_0^{1/2}e^{-mu}u^k\,du\leq\Gamma(k+1)m^{-k-1}.
\end{align}
By \eqref{0324-14},
\begin{align*}
	\mathcal{K}(z)= \frac{\Gamma\left(\frac{n+1}{2}\right)} {\Gamma\left(\frac n2\right)\Gamma\left(\frac12\right)} \int_0^1e^{-zu}u^{\frac{n-2}{2}}(1-u)^{-\frac12}(2u-1)\,du.
\end{align*}
Split the integral at \(1/2\). On \((0,1/2)\), it is bounded for \(0<z\leq1\), and \eqref{261003} bounds it by \(Cz^{-n/2}\) for \(z>1\).
On \((1/2,1)\), it is bounded by \(Ce^{-z/2}\leq C(1+z)^{-n/2}\).
Hence \(|\mathcal{K}(z)|\leq C(1+z)^{-n/2}\).

Differentiating the preceding representation gives
\begin{align*}
	\mathcal{H}(z)= \frac{\Gamma\left(\frac{n+1}{2}\right)} {\Gamma\left(\frac n2\right)\Gamma\left(\frac12\right)} \int_0^1e^{-zu}u^{\frac{n-2}{2}}(1-u)^{-\frac12}(2u-1) \left(\frac n2-zu\right)\,du.
\end{align*}
On \((0,1/2)\), integration by parts yields
\begin{align*}
	-\int_0^{1/2}e^{-zu}u^{\frac n2} \frac{d}{du}\left((1-u)^{-\frac12}(2u-1)\right)\,du,
\end{align*}
whose absolute value is at most \(C(1+z)^{-(n+2)/2}\).
The integral over \((1/2,1)\) has absolute value at most \(C(n+z)e^{-z/2}\), and hence at most \(C(1+z)^{-(n+2)/2}\).
Thus \(|\mathcal{H}(z)|\leq C(1+z)^{-(n+2)/2}\).

\subsection{Zeros of \texorpdfstring{\({}_1F_1\)}{1F1}}\label{apdxzero}
We recall the following theorem on the number of positive zeros of \({}_1F_1(-a,b;z)\) from \cite[Section 6.1]{Slater} and \cite[Eq.~13.9.1]{DLMF}.
\begin{theorem}[\cite{Slater}]\label{thm-0430-01}
Let \(a,b>0\).
Then \({}_1F_1(-a,b;z)\) has exactly \(\lceil a\rceil\) positive zeros.
\end{theorem}

We also recall the following interlacing property of the zeros of \({}_1F_1\).
The proof, although written in terms of Laguerre functions, is found in \cite[Theorem 6]{Javier}.
\begin{theorem}\label{thm-0512-01}
	Let \(a,b>0\), and order the positive zeros of \({}_1F_1(-a,b;z)\) as \(0<\zeta_1(a,b)<\cdots<\zeta_{\lceil a\rceil}(a,b)\).
	Then
	\begin{align}\label{0512-03}
		\zeta_k(a+1,b)<\zeta_k(a,b)<\zeta_{k+1}(a+1,b),\qquad 1\leq k\leq\lceil a\rceil.
	\end{align}
\end{theorem}

We now discuss the relative position of \(z_1\) and \(z_2\), the positive zeros of \({}_1F_1\left(\frac{n}{2}, \frac{n-1}{2};-z\right)\) and of \({}_1F_1\left(\frac{n+2}{2}, \frac{n+1}{2};-z\right)\), respectively.
By Kummer's transformation \eqref{0327-04} these are the positive zeros of \({}_1F_1\left(-\frac12,\frac{n-1}{2};z\right)\) and \({}_1F_1\left(-\frac12,\frac{n+1}{2};z\right)\).
\begin{lemma}\label{lem-0730-49}
	Let \(b>0\).
	Then \({}_1F_1\left(-\frac12,b;z\right)\) is strictly decreasing in \(z\) on \([0,\infty)\) and has exactly one zero \(\zeta(b)\in(0,2b)\).
\end{lemma}
\begin{proof}
	By \eqref{260030}, \({}_1F_1(-\frac12,b;z)=\sum_{k\geq0}\frac{(-\frac12)_k}{(b)_k k!}z^k\), and \((-\frac12)_k=-\frac12(\frac12)_{k-1}<0\) for every \(k\geq1\).
	Hence all coefficients with \(k\geq1\) are strictly negative, so the function is strictly decreasing on \([0,\infty)\).
		For \(z>0\), the terms with \(k\geq2\) give \({}_1F_1(-\frac12,b;z)<1-z/(2b)\).
	The function equals one at zero and is strictly negative at \(z=2b\), so continuity and strict decrease give exactly one zero in \((0,2b)\).
\end{proof}
Thus \(z_1=\zeta\left(\frac{n-1}{2}\right)\) and \(z_2=\zeta\left(\frac{n+1}{2}\right)\).

\subsection{Proof of Lemma \ref{lem-0730-50}}\label{260215}

Write \(b:=\frac{n-1}{2}\).
For \(k\geq1\), define \(a_k:=-\frac{(-\frac12)_k}{(b)_k k!}>0\).
By \eqref{0327-04} and Lemma \ref{lem-0730-49}, \(z_1\) and \(z_2\) are the unique positive zeros of \({}_1F_1(-\frac12,b;z)\) and \({}_1F_1(-\frac12,b+1;z)\), respectively.
The series \eqref{260030} and \((b+1)_k=\frac{(b+k)(b)_k}{b}\) give
\begin{align}
	{}_1F_1\left(-\frac12,b;z\right)=1-\sum_{k=1}^{\infty}a_kz^k,
	\qquad {}_1F_1\left(-\frac12,b+1;z\right)=1-\sum_{k=1}^{\infty}\frac{b}{b+k}\,a_kz^k.\label{0801-11}
\end{align}
Since \({}_1F_1(-\frac12,b;z_1)=0\), we have \(\sum_{k\geq1}a_kz_1^k=1\), and \eqref{0801-11} gives \({}_1F_1(-\frac12,b+1;z_1)>0\).
Thus \(z_1<z_2\).
Set \(s:=\sqrt{z_1}\) and \(c:=(1+\frac1s)^2\).
The coefficients \(a_k\) for \(k\geq2\) in \eqref{0801-11} are strictly positive, so
\begin{align*}
    1=\sum_{k=1}^\infty a_kz_1^k>a_1z_1=\frac{z_1}{2b}.
\end{align*}
Thus \(s<\sqrt{2b}\leq2b\), since \(b\geq\tfrac12\).
It follows that \(c=1+2/s+1/s^2>1+1/b\).
Bernoulli's inequality therefore gives \(\frac{bc^k}{b+k}>1\) for every \(k\geq1\).
Hence
\begin{align*}
	{}_1F_1\left(-\frac12,b+1;\left(1+\sqrt{z_1}\right)^2\right) =1-\sum_{k=1}^\infty a_kz_1^k\,\frac{b}{b+k}c^{k} <1-\sum_{k=1}^\infty a_kz_1^k=0.
\end{align*}
Thus \(z_2<(1+\sqrt{z_1})^2\), which together with \(z_1<z_2\) gives \eqref{0701-05}.

\subsection{Proof of Lemma \ref{lem-0514-01}}\label{260216}

We need a consequence of a single sign change in the coefficients of a power series.
\begin{lemma}\label{lem-0619-01}
	Let \(f(z)=\sum_{k=0}^{\infty}c_kz^k\not\equiv 0\) have radius of convergence \(r>0\).
	Suppose that the nonzero coefficients of \(\{c_k\}\) change sign at most once.
	Then \(f\) has at most one zero in \((0,r)\).
\end{lemma}

\begin{proof}
	If all nonzero coefficients have the same sign, then the conclusion is immediate.
	Thus, after possibly replacing \(f\) by \(-f\), we may assume that there exists \(k_0\geq0\) such that \(c_k\geq0\) for \(k\leq k_0\) and \(c_k\leq0\) for \(k>k_0\), with a nonzero coefficient on each side.
	If \(d_k:=-c_k\geq 0\) for \(k>k_0\), then
	\begin{align*}
		z^{-k_0}f(z)=\sum_{k=0}^{k_0}c_kz^{k-k_0}- \sum_{k=k_0+1}^\infty d_kz^{k-k_0}
	\end{align*}
	is strictly decreasing on \((0,r)\) because the derivative of \(z^{-k_0}f(z)\) is then a series of nonpositive terms with at least one negative term.
	It therefore has at most one zero, and so does \(f\).
\end{proof}

\underline{Proof of (i).}\quad
First consider \(\mathcal{K}(z)\).
Kummer's transformation \eqref{0327-04} and the series representation \eqref{260030} give
\begin{align}\label{0623-10}
	\mathcal{K}(z) & = e^{-z}\left( {}_1F_1\left(\frac12, \frac{n+1}{2};z\right) -\frac{2}{n+1}{}_1F_1\left(\frac32,\frac{n+3}{2};z\right)\right)\nonumber\\
	& =e^{-z}\sum_{k=0}^\infty \left(\frac{\left(\frac12 \right)_k}{\left(\frac{n+1}{2}\right)_k} - \frac{2}{n+1}\frac{\left(\frac32\right)_k}{\left(\frac{n+3}{2}\right)_k}\right)\frac{z^k}{k!} = e^{-z}\sum_{k=0}^\infty \frac{\left(\frac12\right)_k}{\left(\frac{n+1}{2}\right)_k}\frac{n-2k-1}{n+2k+1}\frac{z^k}{k!}.
\end{align}
The coefficients in \eqref{0623-10} are positive for \(k<(n-1)/2\) and negative for \(k>(n-1)/2\).
Lemma \ref{lem-0619-01} therefore gives at most one positive zero.
The decreasing function \(z^{-\lfloor(n-1)/2\rfloor}e^z\mathcal K(z)\) tends to \(-\infty\) as \(z\to\infty\), while \(\mathcal K(0)=(n-1)/(n+1)>0\).
Hence \(\mathcal{K}\) has exactly one positive zero and has the signs in Lemma \ref{lem-0514-01}{\rm(i)}.

\underline{Proof of (ii).}\quad
Using \eqref{0423-10}, we have 
\begin{align}\label{0515-01}
	\begin{split}
		\frac{n+1}{2}\mathcal{K}(z)+z\mathcal{K}'(z) 
		 & = \frac{n-1}{2}{}_1F_1\left(\frac{n}{2}, \frac{n-1}{2};-z\right).
	\end{split}
\end{align}
Theorem \ref{thm-0430-01} gives exactly one positive zero.
Since the value at zero is \((n-1)/2>0\), the function is positive before its zero and negative after it.

\underline{Proof of (iii).}\quad
Using \eqref{0515-01} and \eqref{0327-02}, with \((a,b)=(\frac{n}{2},\frac{n-1}{2})\) and \((a,b)=(\frac{n}{2},\frac{n+1}{2})\), we have
\begin{align*}
	\frac{n}{2}\mathcal{K}(z) + z\mathcal{K}'(z)
	 & = \frac{n}{2}{}_1F_1\left(\frac{n+2}{2}, \frac{n+1}{2};-z\right) - \frac{n}{n+1}{}_1F_1\left(\frac{n+2}{2}, \frac{n+3}{2};-z\right).
\end{align*}
Kummer's transformation \eqref{0327-04} and the series representation \eqref{260030} show that the last expression equals
\begin{align*}
	 & e^{-z}\left(\frac{n}{2}{}_1F_1\left(-\frac12, \frac{n+1}{2};z\right) -\frac{n}{n+1}{}_1F_1\left(\frac12, \frac{n+3}{2};z\right)\right)\\
	 & \quad=e^{-z}\left[\frac{n(n-1)}{2(n+1)}
	 -\frac n2\sum_{k=1}^\infty \frac{(\frac12)_k}{(\frac{n+1}{2})_k}
	 \frac{n+6k-1}{(2k-1)(n+2k+1)}\frac{z^k}{k!}\right].
\end{align*}
Every coefficient in the sum starting at \(k=1\) is positive, and its first term grows linearly in \(z\).
Thus the expression in square brackets decreases strictly from \(n(n-1)/(2(n+1))\) to \(-\infty\).
Since \(e^{-z}>0\), \(\frac n2\mathcal{K}(z)+z\mathcal{K}'(z)\) has exactly one positive zero and has the signs in Lemma \ref{lem-0514-01}{\rm(iii)}.

We now prove \eqref{0623-02}.
At \(z=z_3\), both \(\frac{n+1}{2}\mathcal{K}(z)+z\mathcal{K}'(z)\) and \(\frac n2\mathcal{K}(z)+z\mathcal{K}'(z)\) equal \(z_3\mathcal{K}'(z_3)<0\).
Indeed, the proof of Lemma \ref{lem-0619-01}, applied to the series in \eqref{0623-10}, gives \(\frac{d}{dz}[z^{-k_0}e^z\mathcal K(z)]<0\) on \((0,\infty)\), where \(k_0=\lfloor(n-1)/2\rfloor\).
Evaluating this derivative at \(z=z_3\), where \(\mathcal K(z_3)=0\), gives \(\mathcal K'(z_3)<0\).

Since both functions are positive near \(0\), it follows that \(z_1<z_3\) and \(z_4<z_3\).
At \(z=z_4\), the first function equals \(\frac12\mathcal{K}(z_4)>0\).
Thus \(z_4<z_1<z_3\), which is \eqref{0623-02}.

\subsection{Proof of Lemma \ref{lem-0520-01}}\label{260217}

\underline{Proof of (i).}\quad
The derivative identity \eqref{0327-07} expresses \(\mathcal K^{(j)}\), \(j=0,1,2,3\), in terms of the confluent hypergeometric functions with parameters \((a,b)=(\frac n2+j,\frac{n+1}{2}+j)\) and \((a,b)=(\frac n2+j,\frac{n+3}{2}+j)\).
Applying \eqref{0504-06} to these functions gives \(Z_1>0\) such that the following expansion and its first three derivatives hold for \(z>Z_1\).
\begin{align*}
	\mathcal{K}(z) = - \frac{\Gamma(\frac{n+1}{2})}{\sqrt{\pi}}z^{-\frac{n}{2}} \left[ 1 - \frac{3n}{4z} - \frac{5n(n+2)}{32z^2}+ O(z^{-3})\right].
\end{align*}
Applying \(\ln(1-x) = -x- \frac{x^2}{2}- \frac{x^3}{3}-\cdots\) to this expansion gives, for some constant \(C\) and \(z>Z_1\),
\begin{align}\label{0520-33}
	\mathcal{L}(z) = C- \frac{n}{2} \ln z - \frac{3n}{4z} - \frac{14n^2 + 10n}{32z^2} +O(z^{-3}).
\end{align}
Differentiating \eqref{0520-33} gives
\begin{align}\label{0520-35}
	\mathcal{M}(z) = \frac{3n}{4z^2} + \frac{7n^2 + 5n}{4z^3} + O(z^{-4}).
\end{align}
Hence \(\mathcal{L}'(z)<0\), \(\mathcal{M}(z)>0\), and the bound \eqref{0520-31} holds for sufficiently large \(z>Z_1\).
Differentiating \eqref{0520-35} gives
\begin{align*}
	\mathcal{M}'(z) = -\frac{3n}{2z^3}-\frac{3(7n^2+5n)}{4z^4}+O(z^{-5}),
\end{align*}
so \(\mathcal{M}(z)\) is decreasing for sufficiently large \(z>Z_1\).

\underline{Proof of (ii).}\quad
Since \(\mathcal M(z)>0\) for all sufficiently large \(z\), \(z\mathcal L'(z)\) is eventually decreasing.
Together with \eqref{0520-33}, this yields \(Z_2>Z_1\) such that \(z\mathcal L'(z)>-n/2\) for \(z>Z_2\).
Integrating \(-\mathcal{L}'(z)<\frac{n}{2z}\) from \(\zeta_1\) to \(\zeta_2\) gives
\begin{align*}
\ln\frac{-\mathcal{K}(\zeta_1)}{-\mathcal{K}(\zeta_2)}=\mathcal{L}(\zeta_1)-\mathcal{L}(\zeta_2)=-\int_{\zeta_1}^{\zeta_2}\mathcal{L}'(z)\,dz<\frac n2\ln\frac{\zeta_2}{\zeta_1}.
\end{align*}
This gives the bound \eqref{0520-30}.

\underline{Proof of (iii).}\quad
Using \eqref{0520-33}, we have, for \(z>Z_2\),
\begin{align*}
	 & \mathcal{L}'(z)=-\frac{n}{2z} + \frac{3n}{4z^2} + O(z^{-3}),\quad \mathcal{L}''(z) = \frac{n}{2z^2} - \frac{3n}{2z^3} + O(z^{-4}).
\end{align*}
Hence the bound \eqref{0520-40} holds.
Choosing \(Z_0>\max\{z_3,Z_1,Z_2\}\) sufficiently large completes the proof of {\rm(iii)}.

\subsection{Proof of Lemma \ref{lem-260032}}\label{260218}

Write \(b=(n-1)/2\) and \(H(z)={}_1F_1(-\frac12,b;z)\).
By \eqref{0327-04}, \(F(z)=z^{\frac{n+1}{2}}e^{-z}H(z)\).
Moreover, \eqref{0327-07} and \eqref{0324-14} give
\begin{align*}
	H'(z)=-\frac{1}{n-1}{}_1F_1\left(\frac12,\frac{n+1}{2};z\right)<0,\qquad z>0.
\end{align*}
Lemma \ref{lem-0730-49} shows that \(z_1\) is the unique positive zero of \(H\).
Consequently, \(F'(z_1)=z_1^{\frac{n+1}{2}}e^{-z_1}H'(z_1)<0\).

On the other hand, the power series \eqref{260030} gives
\begin{align*}
	F'(z)=\frac{n+1}{2}z^{\frac{n-1}{2}}\left(1+O(z)\right)>0
\end{align*}
for all sufficiently small \(z>0\).
Hence \(F'\) has at least one positive zero.

Differentiating \(F\) and then using \eqref{0327-01} and \eqref{0327-04}, we obtain
\begin{align*}
	F'(z)=\frac12z^{\frac{n-1}{2}}e^{-z}
	\left[{}_1F_1\left(-\frac12,b;z\right)
	+n{}_1F_1\left(-\frac32,b;z\right)\right].
\end{align*}
Applying \eqref{0327-10} to the two confluent hypergeometric functions yields
\begin{align}\label{1026-01}
	F'(z)=\frac{\Gamma(\frac{n-1}{2})}{2\Gamma(-\frac12)}
	z^{-\frac12}\left(1+O(z^{-1})\right)
	\qquad\text{as }z\to\infty.
\end{align}
Since \(\Gamma(-\frac12)<0\), it follows that \(F'(z)<0\) for all sufficiently large \(z\).

Choose \(0<\delta<1<Z\) so that \(F'(z)>0\) for \(0<z\leq\delta\) and \(F'(z)<0\) for \(z\geq Z\).
Since \(F'\) is real-analytic and is not identically zero, it has only finitely many zeros in \([\delta,Z]\).
Thus \(F'\) has a largest positive zero \(z_6\), and \(F'<0\) on \((z_6,\infty)\).

Finally, the power series above gives \(|F'(z)|\leq C(n)z^{(n-1)/2}\leq C(n)z^{-1/2}\) for \(0<z\leq1\), while \eqref{1026-01} gives \(|F'(z)|\leq C(n)z^{-1/2}\) for large \(z\).
Since \(F'\) is continuous on \((0,\infty)\), the estimate holds for every \(z>0\).

\section{Integral identities}\label{260240}
\subsection{Proof of Proposition \ref{prop-0322-001}}\label{260241}
\begin{lemma}\label{lem-0331-01}
For \(\alpha>0\), \(r>0\), and \(\mu,\nu\in \R\) with \(\nu>-1\) and \(\mu+\nu>-1\), 
	\begin{align}\label{0401-01}
		\int_0^\infty\rho^{\mu}e^{-\alpha\rho^2}J_{\nu}(r \rho)\,d\rho= \frac{r^{\nu}\Gamma\left(\frac{\mu+\nu+1}{2}\right)}{2^{\nu+1}\alpha^{\frac{\mu+\nu+1}{2}}\Gamma(\nu+1)} {}_1F_1\left(\frac{\mu+\nu+1}{2}, \nu+1; - \frac{r^2}{4\alpha}\right).
	\end{align}
\end{lemma}
\begin{proof}

	Substituting \(u=\sqrt\alpha\,\rho\) reduces the proof to \(\alpha=1\).
	In the Bessel series \eqref{0401-05}, the integral of each term reduces to
	\(\int_0^\infty u^{\mu+\nu+2k}e^{-u^2}\,du=\frac12\Gamma(\frac{\mu+\nu+1}{2}+k)\), and integrating term by term gives \eqref{0401-01}.
\end{proof}
\begin{proof}[Proof of Proposition \ref{prop-0322-001}]
	By \eqref{0419-021}, the principal-value kernel \(2\partial_{x_i}N(\cdot,0)\) defines \(R_i'\), with Fourier multiplier \(-i\xi_i/|\xi'|\).
	Multiplying half this symbol by \(e^{-t|\xi'|^2}\), the Fourier transform of \(\Gamma'(\cdot,t)\), gives
	\begin{align}\label{0630-20}
		\mathcal{F}_{x'}\left( \mathrm{p.v.}\int_{\Sigma} \Gamma'(x'-z',t)\partial_{z_i}N(z',0)\,dz' \right) = -\frac{i\xi_i}{2|\xi'|}e^{-t|\xi'|^2}.
	\end{align}
	This multiplier is integrable on \(\Sigma\), so its inverse Fourier transform is continuous.
	The right-hand side of \eqref{0320-01} is also continuous at the origin.
	It therefore suffices to work at \(x'\ne0\) and restore the minus sign in \eqref{0630-20} after inverting
	\begin{align}\label{0630-10}
		\frac{i\xi_i}{2|\xi'|}e^{-t|\xi'|^2}.
	\end{align}

	Suppose first that \(n\geq3\).
	Then \(d=n-1\geq2\), so \((2|\xi'|)^{-1}e^{-t|\xi'|^2}\) is integrable at the origin as well as at infinity.
	Let \(f\) be its inverse Fourier transform.	Then \(\partial_{x_i}f\) is the inverse transform of \eqref{0630-10}.
	The radial inversion formula \eqref{0401-07} gives, for \(x'\ne0\),
	\begin{align}\label{0927-0010}
		f(x',t) = \frac{1}{2(2\pi)^{\frac{n-1}{2}}|x'|^{\frac{n-3}{2}}} \int_0^\infty \rho^{\frac{n-3}{2}}e^{-t\rho^2} J_{\frac{n-3}{2}}(|x'|\rho)\,d\rho.
	\end{align}
	Differentiating \eqref{0927-0010} with respect to \(x_i\), and using \eqref{0401-06} with \(\nu=(n-3)/2\), we obtain
	\begin{align}\label{0630-22}
		\partial_{x_i}f(x',t) = -\frac{x_i}{2(2\pi)^{\frac{n-1}{2}}|x'|^{\frac{n-1}{2}}} \int_0^\infty \rho^{\frac{n-1}{2}}e^{-t\rho^2} J_{\frac{n-1}{2}}(|x'|\rho)\,d\rho.
	\end{align}

	For \(n=2\), the scalar factor \(|\xi_1|^{-1}\) is singular at zero.
	We therefore invert the integrable multiplier \eqref{0630-10} directly, obtaining
\begin{align*}
-\frac1{2\pi}\int_0^\infty e^{-t\rho^2}\sin(x_1\rho)\,d\rho.
\end{align*}
	Since \(J_{1/2}(s)=\sqrt{2/(\pi s)}\sin s\) by \eqref{0401-05}, this sine integral equals the right-hand side of \eqref{0630-22}.
	Thus \eqref{0630-22} represents the inverse transform in every dimension \(n\geq2\).

	We now evaluate the integral in \eqref{0630-22}.
	Applying \eqref{0401-01} with \( \alpha=t, r=|x'|, \mu=\nu=\frac{n-1}{2}\), we obtain
	\begin{align*}
		\int_0^\infty \rho^{\frac{n-1}{2}}e^{-t\rho^2} J_{\frac{n-1}{2}}(|x'|\rho)\,d\rho = \frac{|x'|^{\frac{n-1}{2}}\Gamma\left(\frac n2\right)} {2^{\frac{n+1}{2}}t^{\frac n2}\Gamma\left(\frac{n+1}{2}\right)} {}_1F_1\left( \frac n2,\frac{n+1}{2}; -\frac{|x'|^2}{4t} \right).
	\end{align*}
	Substituting this into \eqref{0630-22} yields
	\begin{align}\label{0630-24}
		\mathcal{F}_{x'}^{-1} \left( \frac{i\xi_i}{2|\xi'|}e^{-t|\xi'|^2} \right)(x') & = -\frac{\Gamma\left(\frac n2\right)} {2^{n+1}\pi^{\frac{n-1}{2}}\Gamma\left(\frac{n+1}{2}\right)} \frac{x_i}{t^{\frac n2}} {}_1F_1\left( \frac n2,\frac{n+1}{2}; -\frac{|x'|^2}{4t} \right)\nonumber\\
		 & =\frac{\Gamma\left(\frac n2\right)} {2\Gamma\left(\frac{n+1}{2}\right)} \sqrt t\,\partial_{x_i}\Gamma'(x',t) {}_1F_1\left( \frac12,\frac{n+1}{2}; \frac{|x'|^2}{4t} \right),
	\end{align}
	where we used Kummer's transformation \eqref{0327-04} in the second equality.

	Finally, by \eqref{0630-20}, the Fourier transform of the left-hand side of \eqref{0320-01} is the negative of the multiplier in \eqref{0630-10}.
The identity \eqref{0630-24} therefore gives \eqref{0320-01}.
\end{proof}

\subsection{Proof of Lemma \ref{lem-0422-03}}\label{260242}
We reduce the spherical integral to one variable using \eqref{0907-01}.

If \(n=2\), then the integrand on the left-hand side of \eqref{0422-07} is one.
On the right-hand side, \((0)_k=0\) for every \(k\geq1\), so both sides equal \(|\mathbb S^0|\).

Assume \(n\geq3\) and write \(\rho=s/r\in[0,1)\).
Multiplying the left-hand side of \eqref{0422-07} by \(r^{n-2}\) and applying \eqref{0907-01} gives
\begin{align*}
    |\mathbb S^{n-3}|\int_{-1}^1(1-2\rho\tau+\rho^2)^{-(n-2)/2}(1-\tau^2)^{(n-4)/2}\,d\tau.
\end{align*}
Applying \eqref{0907-01} to the constant function instead gives the normalization
\begin{align*}
    \frac{|\mathbb S^{n-2}|}{|\mathbb S^{n-3}|}
    =\int_{-1}^1(1-\tau^2)^{(n-4)/2}\,d\tau
    =\frac{\sqrt\pi\,\Gamma((n-2)/2)}{\Gamma((n-1)/2)}.
\end{align*}
Thus \eqref{0422-07} is equivalent to
\begin{align}\label{0422-06}
	\int_{-1}^1(1-2\rho\tau+\rho^2)^{-\frac{n-2}{2}}(1-\tau^2)^{\frac{n-4}{2}}\,d\tau = \frac{\sqrt\pi\,\Gamma\left(\frac{n-2}{2}\right)} {\Gamma\left(\frac{n-1}{2}\right)} \sum_{k=0}^\infty \frac{\left(\frac{n-2}{2}\right)_k\left(\frac12\right)_k} {k!\left(\frac{n-1}{2}\right)_k}\,\rho^{2k}.
\end{align}

We next show the integral identity
\begin{align}\label{0422-08}
	\int_{-1}^1(1-2\rho\tau+\rho^2)^{-\frac{n-2}{2}}(1-\tau^2)^{\frac{n-4}{2}}\,d\tau = \int_{-1}^1(1-\rho^2\tau^2)^{-\frac{n-2}{2}}(1-\tau^2)^{\frac{n-4}{2}}\,d\tau.
\end{align}
To prove \eqref{0422-08}, set \(u=(\rho-\tau)/(1-\rho\tau)\).
Direct calculation gives
\begin{align*}
	1-\rho u=\frac{1-\rho^2}{1-\rho\tau}, \qquad 1+\rho u=\frac{1-2\rho\tau+\rho^2}{1-\rho\tau}, \qquad 1-u^2=\frac{(1-\rho^2)(1-\tau^2)}{(1-\rho\tau)^2},
\end{align*}
and \(du/d\tau=-(1-\rho^2)/(1-\rho\tau)^2<0\).
Multiplying the first two identities and combining with the formulas for \(1-u^2\) and \(du/d\tau\), we obtain
\begin{align*}
	 & (1-\rho^2u^2)^{-\frac{n-2}{2}} (1-u^2)^{\frac{n-4}{2}}\left|\frac{du}{d\tau}\right|=(1-2\rho\tau+\rho^2)^{-\frac{n-2}{2}} (1-\tau^2)^{\frac{n-4}{2}}.
\end{align*}
Changing variables in the integral on the right-hand side of \eqref{0422-08} gives its left-hand side.

It remains to expand the right-hand side of \eqref{0422-08}.
Since \(\rho^2\tau^2\leq\rho^2<1\) on \([-1,1]\), the binomial series
\begin{align*}
	(1-\rho^2\tau^2)^{-\frac{n-2}{2}} = \sum_{k=0}^\infty\frac{\left(\frac{n-2}{2}\right)_k}{k!}\rho^{2k}\tau^{2k}
\end{align*}
converges uniformly in \(\tau\in[-1,1]\).
The beta integral gives
\begin{align*}
	\int_{-1}^1\tau^{2k}(1-\tau^2)^{\frac{n-4}{2}}\,d\tau = \frac{\left(\frac12\right)_k}{\left(\frac{n-1}{2}\right)_k} \cdot \frac{\sqrt\pi\,\Gamma\left(\frac{n-2}{2}\right)} {\Gamma\left(\frac{n-1}{2}\right)}, \qquad k\geq 0.
\end{align*}
Substituting this identity into the series integrated term by term gives \eqref{0422-06}.

\section{Volterra-type operators}\label{app-volterra}
\subsection{Proof of Lemma \ref{lem-0424-20}}\label{app-volterra-differentiation}
\begin{proof}
	The first bound in \eqref{0829-01} gives \(K\in L^1(0,1)\).
	Writing \((K*f)(t)=\int_0^t K(s)f(t-s)\,ds\), dominated convergence and continuity of \(f\) give \(K*f\in C(0,1)\).
	For \(0<h<1-t\), subtract \(f(t)\) inside the difference quotient to obtain
	\begin{align*}
		\frac{(K*f)(t+h)-(K*f)(t)}{h} & = \frac{1}{h}\left( \int_0^t (K(t+h-s)-K(t-s))(f(s)-f(t))\,ds \right.\\
		 & \quad \left. + \int_0^h K(s) (f(t+h-s)-f(t))\,ds + f(t) \int_t^{t+h}K(s)\,ds\right)=:I_1+I_2+I_3.
	\end{align*}
	For \(I_1\), the derivative bound for \(K\) gives
	\begin{align*}
		\left|\frac{1}{h}(K(t+h-s)-K(t-s))(f(s)-f(t))\right|\leq C[f]_{C^\alpha}(t-s)^{\alpha-\sigma-1}.
	\end{align*}
	Since \(\alpha>\sigma\), dominated convergence in \(s\) gives
	\begin{align*}
		\lim_{h\to 0^+} I_1 & = \int_0^t (f(s)-f(t))K'(t-s)\,ds.
	\end{align*}
	For \(I_2\), the inequality \(|f(t+h-s)-f(t)|\leq [f]_{C^\alpha}(h-s)^\alpha\) for \(0<s<h\) gives
	\begin{align*}
		|I_2|\leq \frac{C}{h}[f]_{C^\alpha}\int_0^h \frac{(h-s)^\alpha}{s^{\sigma}}\,ds \leq C[f]_{C^\alpha}h^{\alpha-\sigma}.
	\end{align*}
	Thus \(I_2\to0\) as \(h\to0^+\).
	Finally, \(I_3\to K(t)f(t)\), since \(K\) is continuous at \(t>0\).
	Adding the three limits identifies the right derivative of \(K*f\) with \eqref{0325-02}.

	Write the integral term in \eqref{0325-02} as \(A(t)=\int_0^t K'(u)[f(t-u)-f(t)]\,du\).
	Fix \([a,b]\subset(0,1)\) and take \(a\leq t<t+h\leq b\) with \(0<h<a\).
	On the common interval \(0<u<t\), the difference of the two brackets is bounded by
	\begin{align*}
		|f(t+h-u)-f(t+h)-f(t-u)+f(t)|
		\leq 2[f]_{C^\alpha}\min\{u^\alpha,h^\alpha\}.
	\end{align*}
	On \(t<u<t+h\), the kernel derivative is bounded since \(u\geq a\).
	Hence splitting the common interval at \(u=h\) gives
	\begin{align*}
		|A(t+h)-A(t)|
		&\leq C\left(\int_0^h u^{\alpha-\sigma-1}\,du
		+h^\alpha\int_h^t u^{-1-\sigma}\,du+h\right)
		\leq C_{a,b}h^{\alpha-\sigma}.
	\end{align*}
Since \(Kf\) is \(\alpha\)-H\"older continuous on \([a,b]\), the right derivative in \eqref{0325-02} is locally \((\alpha-\sigma)\)-H\"older continuous; as \(K*f\) is continuous, it is the ordinary derivative.

	For the kernel depending on \(r\), substituting \(u=t-s\) in the integrals for the two partial derivatives gives the bounds
	\begin{align*}
		|\partial_rK(r,u)f(t-u)|&\leq C_I\|f\|_\infty u^{-\sigma},\\
		|\partial_uK(r,u)[f(t-u)-f(t)]|&\leq C_I[f]_{C^\alpha}u^{\alpha-\sigma-1}.
	\end{align*}
	Both powers are integrable at zero, which gives differentiation under the integral in \(r\) and joint continuity in \((r,t)\).
\end{proof}

\subsection{Proof of Lemma \ref{lem-0325-01}}\label{260236}

Extend \(\phi\) by zero on \((-\infty,0]\) and fix \(0<t_a<t_b<\tau<1\).
For fixed \(0<u<t_a\), the derivative of \(t\mapsto\ln\phi(t-u)-\ln\phi(t)\) on \((u,1)\) is
\((\ln\phi)'(t-u)-(\ln\phi)'(t)\geq0\) by log-concavity of \(\phi\).
Consequently,
\begin{align}\label{log-concave-comparison}
    \frac{\phi(t_a-u)}{\phi(t_a)}
    \leq\frac{\phi(t_b-u)}{\phi(t_b)},\qquad 0\leq u\leq\tau .
\end{align}
For \(t_a<u<t_b\), \(\phi(t_a-u)=0<\phi(t_b-u)\), and for \(u\geq t_b\), both vanish.

For \(0<t<\tau\), integrate by parts on \([\varepsilon,t]\) and \([t,\tau]\).
Since \(\phi(0)=\phi(t-\tau)=0\), the terms at \(u=t\) agree and the endpoint term at \(\tau\) gives \(K(\tau)\). Thus
\begin{align*}
    \frac{(K*\phi')(t)}{\phi(t)}
    =K(\tau)+\int_0^\tau[-K'(u)]
    \left(1-\frac{\phi(t-u)}{\phi(t)}\right)du .
\end{align*}
Indeed, \(1-\phi(t-u)/\phi(t)=O_t(u)\) as \(u\downarrow0\).
By \eqref{0829-01}, the boundary term at \(\varepsilon\) is \(O_t(\varepsilon^{1-\sigma})\) and the integrand is \(O_t(u^{-\sigma})\).
Since \(-K'\geq0\) and \(\int_{t_a}^{t_b}(-K')=K(t_a)-K(t_b)>0\), the strict comparison on \((t_a,t_b)\) gives strict decrease.
Positivity of \(\phi\) then gives the asserted signs at any zero.

For {\rm(ii)}, assume now that \(\phi'\in C^\alpha\) with \(\alpha>\sigma\). Lemma~\ref{lem-0424-20} gives \(K*\phi'\in C^1(0,1)\).
Fix a zero \(t_*\in(0,1)\), so that
\begin{align}\label{0325-01}
	\int_0^{t_*}K(t_*-s)\phi'(s)\,ds=0.
\end{align}
Writing \(K(t_*-s)=K(t_*)-\int_0^s K'(t_*-u)\,du\) and using \eqref{0325-01}, Fubini's theorem, and \(\phi(0)=0\), we obtain
\begin{align}\label{0325-03}
			0 = K(t_*)\phi(t_*)-\int_0^{t_*}K'(t_*-s)(\phi(t_*)-\phi(s))\,ds.
\end{align}
Eliminating \(K(t_*)\phi(t_*)\) from the derivative formula \eqref{0325-02} by \eqref{0325-03} gives
\begin{align}\label{0325-04}
    (K*\phi')'(t_*)
    &=\int_0^{t_*} K'(t_*-s)
       \left[\phi'(s)-\frac{\phi'(t_*)}{\phi(t_*)}\phi(s)\right]\,ds \nonumber\\
    &=\int_0^{t_*} K'(t_*-s)\phi(s)
       \bigl[(\ln\phi)'(s)-(\ln\phi)'(t_*)\bigr]\,ds.
\end{align}
Log-concavity of \(\phi\) makes the bracket in \eqref{0325-04} nonnegative.
It is positive on some \((0,s_0)\) with \(s_0<t_*\), since \(\phi(0)=0\).

Since \(K\) is strictly decreasing, \(\int_{t_*-s_0}^{t_*}(-K'(u))\,du>0\).
The other factors are positive for \(0<s<s_0\), so \eqref{0325-04} gives \((K*\phi')'(t_*)<0\).

\subsection{Proof of Lemma \ref{lem-0506-01}}\label{260237}

Fix \(0<t_a<t_b<\tau<1\) and extend \(\phi\) by zero to negative times.
The integration by parts in Appendix~\ref{260236} applies to \(K_1\) by \eqref{0829-01}. Subtracting the resulting expressions for \((T\phi)(t_a)/\phi(t_a)\) and \((T\phi)(t_b)/\phi(t_b)\) gives
\begin{align*}
    &\frac{(T\phi)(t_b)}{\phi(t_b)}-\frac{(T\phi)(t_a)}{\phi(t_a)}\\
    &\qquad=-\int_0^\tau[-K_1'(u)+K_2(u)]
    \left(\frac{\phi(t_b-u)}{\phi(t_b)}
         -\frac{\phi(t_a-u)}{\phi(t_a)}\right)du<0 .
\end{align*}
The weight \(-K_1'+K_2\) is positive. By \eqref{log-concave-comparison}, the difference in parentheses is nonnegative throughout and positive on \((t_a,t_b)\).
This gives strict decrease and the signs at a zero for \(\phi\in C^1([0,1])\).

Assume now that \(\phi'\in C^\alpha\) with \(\alpha>\sigma\).
By Lemma \ref{lem-0424-20} with \(K=K_1\) and \(f=\phi'\), the first integral defining \(T\phi\) belongs to \(C^1(0,1)\).
Since \(K_2\in L^1(0,1)\) and \(\phi(0)=0\), \((K_2*\phi)'=K_2*\phi'\).
Thus \(T\phi\in C^1(0,1)\).

Fix a zero \(t_*\).
Its defining equation is
\begin{align}\label{0506-03}
	\int_0^{t_*}K_1(t_*-s)\phi'(s)\,ds = \int_0^{t_*}K_2(t_*-s)\phi(s)\,ds.
\end{align}
Writing \(K_1(t_*-s)=K_1(t_*)-\int_0^s K_1'(t_*-u)\,du\) and interchanging the integrals as in \eqref{0325-03}, we obtain from \eqref{0506-03}
\begin{equation}\label{0506-01}
	K_1(t_*)\phi(t_*) = \int_0^{t_*} K_2(t_*-s)\phi(s)\,ds + \int_0^{t_*} K_1'(t_*-u)(\phi(t_*) - \phi(u))\,du.
\end{equation}
Using \eqref{0325-02}, we obtain
\begin{equation}\label{0506-02}
	(T\phi)'(t_*) = K_1(t_*)\phi'(t_*) + \int_0^{t_*} K_1'(t_*-s)(\phi'(s) - \phi'(t_*))\,ds - \int_0^{t_*} K_2(t_*-s)\phi'(s)\,ds.
\end{equation}
Multiply \eqref{0506-02} by \(\phi(t_*)\) and substitute \eqref{0506-01} for \(K_1(t_*)\phi(t_*)\).
Collecting the \(K_1'\) and \(K_2\) terms gives
\begin{align*}
	\phi(t_*)(T\phi)'(t_*)=
	\int_0^{t_*}\bigl(K_1'(t_*-s)-K_2(t_*-s)\bigr)
	\bigl(\phi(t_*)\phi'(s)-\phi'(t_*)\phi(s)\bigr)\,ds .
\end{align*}

As in the proof of Lemma~\ref{lem-0325-01}, the logarithmic derivative at \(s\) is at least its value at \(t_*\) for every \(0<s<t_*\), with strict inequality on an initial subinterval.
Thus the second factor in the integral is nonnegative and positive on that subinterval.
Since \(K_1'-K_2<0\), we obtain \((T\phi)'(t_*)<0\).

\subsection{Proof of Lemma \ref{lem-0504-01}}\label{260238}
It suffices to show that \((K*\phi)'(t_*)>0\) at every zero \(t_*>1\).
At such a zero, integration by parts gives
\begin{align}\label{260239}
	(K*\phi)'(t_*) = \int_0^1 K'(t_*-s)\phi(s)\,ds = \int_0^1 K(t_*-s)\phi'(s)\,ds,
	\end{align}
	since \(\phi(0)=\phi(1)=0\).

Since \((K*\phi)(t_*)=0\), \(\phi>0\) on \((0,1)\), and \(K\) vanishes only at \(\tau_0\), the factor \(K(t_*-s)\) takes both signs.
It follows that \(t_*-1<\tau_0<t_*\).
Thus \(s_*=t_*-\tau_0\in(0,1)\), and \(K(t_*-s)\) is positive for \(s<s_*\) and negative for \(s>s_*\).

Subtracting \(\frac{\phi'(s_*)(K*\phi)(t_*)}{\phi(s_*)}=0\) from \eqref{260239} gives
\begin{equation}\label{0504-09}
	(K*\phi)'(t_*) = \int_0^1 K(t_*-s)\left[ \frac{\phi'(s)}{\phi(s)} - \frac{\phi'(s_*)}{\phi(s_*)} \right] \phi(s)\,ds.
\end{equation}

Log-concavity of \(\phi\) makes the bracket in \eqref{0504-09} nonnegative for \(s<s_*\) and nonpositive for \(s>s_*\).
It is positive on an initial subinterval since \(\phi(0)=0\), as in the proof of Lemma~\ref{lem-0325-01}.
The integrand is therefore nonnegative everywhere and positive on that subinterval. Hence \((K*\phi)'(t_*)>0\).

\section{H\"{o}lder estimates}\label{260231}
We first specify the spatial and parabolic H\"older norms used for the heat potentials.

Let \(d\geq 1\), let \(\Omega\subset\R^d\) be open, let \(q\in\mathbb N_0\), and let \(0<\eta<1\).
We set
\begin{align*}
	[f]_{C^\eta(\Omega)}
	:=\sup_{\substack{x,y\in\Omega\\x\neq y}}
	\frac{|f(x)-f(y)|}{|x-y|^\eta},\quad
	\|f\|_{C^{q,\eta}(\Omega)}
	:=\sum_{|a|\leq q}\|\partial^a f\|_{L^\infty(\Omega)}
	+\sum_{|a|=q}[\partial^a f]_{C^\eta(\Omega)}.
\end{align*}
For noninteger \(p=q+\eta\), write \(C^p(\Omega)=C^{q,\eta}(\Omega)\).
We also denote by \(C_c^{q,\eta}(\Omega)\) the subspace of functions with compact support in \(\Omega\).

A function \(f\in C^{q,\eta}(\Omega)\) belongs to
\(C^{q,\eta}(\overline\Omega)\) if \(f\) and all its derivatives through
order \(q\) extend continuously to \(\overline\Omega\).

For the boundary data, take \(I=[0,T]\), \(T>0\), or \(I=[0,\infty)\), and
write \(Q_I'=\Sigma\times I\).
For bounded \(v\colon Q_I'\to\R\), the separate spatial and temporal
seminorms of order \(0<\eta<1\) are
\begin{align*}
	[v]_{x';\eta;Q_I'}
	:=\sup_{\substack{t\in I\\x',y'\in\Sigma,\ x'\neq y'}}
	\frac{|v(x',t)-v(y',t)|}{|x'-y'|^\eta},\quad
	[v]_{t;\eta;Q_I'}
	:=\sup_{\substack{x'\in\Sigma\\t,s\in I,\ t\neq s}}
	\frac{|v(x',t)-v(x',s)|}{|t-s|^\eta}.
\end{align*}
Given noninteger \(p>0\), let \(q=\lfloor p\rfloor\) and \(\eta=p-q\).
A tangential derivative has parabolic order one and a time derivative has
parabolic order two.
With this convention, the parabolic H\"older norm is
\begin{align*}
	\|v\|_{C_{x',t}^{p,\frac p2}(Q_I')}
	&:=\sum_{|a|+2r\leq q}
	\left\|\partial_{x'}^a\partial_t^r v\right\|_{L^\infty(Q_I')}+\sum_{|a|+2r=q}
	\left[\partial_{x'}^a\partial_t^r v\right]_{x';\eta;Q_I'}\\
	&\quad+\sum_{\substack{|a|+2r\leq q\\0<p-|a|-2r<2}}
	\left[\partial_{x'}^a\partial_t^r v
	\right]_{t;\frac{p-|a|-2r}{2};Q_I'}.
\end{align*}
The space \(C_{x',t}^{p,p/2}(Q_I')\) consists of the functions with the
indicated derivatives for which this norm is finite.

Finally, let \(F\colon\overline{\R_+^n}\times I\to\R\) be bounded, and let
\(\alpha,\beta,\gamma\in(0,1)\).
Define
\begin{align*}
	&[F]_{x';\alpha}
	:=\sup_{\substack{x',y'\in\Sigma,\ x'\neq y'\\x_n\geq0,\ t\in I}}
	\frac{|F(x,t)-F(y',x_n,t)|}{|x'-y'|^\alpha},\quad
	[F]_{x_n;\beta}
	:=\sup_{\substack{x'\in\Sigma,\ t\in I\\r,s\geq0,\ r\neq s}}
			\frac{|F(x',r,t)-F(x',s,t)|}{|r-s|^\beta},\\
	&[F]_{t;\gamma}
	:=\sup_{\substack{x'\in\Sigma,\ x_n\geq0\\t,s\in I,\ t\neq s}}
	\frac{|F(x,t)-F(x,s)|}{|t-s|^\gamma},
\end{align*}
and
\begin{align*}
	\|F\|_{C_{x'}^\alpha C_{x_n}^\beta C_t^\gamma
	(\overline{\R_+^n}\times I)}
	:=\|F\|_{L^\infty(\overline{\R_+^n}\times I)}
	+[F]_{x';\alpha}+[F]_{x_n;\beta}+[F]_{t;\gamma}.
\end{align*}

The following estimate for the double-layer heat potential uses the
normalization of the heat kernel in \eqref{0703-05}.

\begin{lemma}\label{lem-260232}
	Let \(\nu\in\mathbb N_0^{n-1}\), \(l,m\in\mathbb N_0\),
	\(\alpha,\beta\in(0,1)\), \(0<\gamma<1/2\), and \(T>0\), and set \(I=[0,T]\) and \(p=|\nu|+l+2m+\max\{\alpha,\beta,2\gamma\}\).
	Suppose that \(g\in C_{x',t}^{p,p/2}(Q_I')\) satisfies
	\(\partial_t^r g(x',0)=0\) for all \(x'\in\Sigma\) and
	\(0\leq r\leq\lfloor p/2\rfloor\), and let
	\begin{align*}
		k(x,t):=\int_0^t\int_\Sigma
		\partial_{x_n}\Gamma(x-y',t-s)g(y',s)\,dy'ds.
	\end{align*}
	Then \(\partial_{x'}^\nu\partial_{x_n}^l\partial_t^m k\) extends continuously to \(\overline{\R_+^n}\times I\), and
	\begin{align}\label{260233}
		\left\|
		\partial_{x'}^\nu\partial_{x_n}^l\partial_t^m k
		\right\|_{C_{x'}^\alpha C_{x_n}^\beta C_t^\gamma
		(\overline{\R_+^n}\times I)}
		\leq C
		\|g\|_{C_{x',t}^{p,\frac p2}(Q_I')},
	\end{align}
	where \(C\) depends only on
	\(n,\nu,l,m,\alpha,\beta,\gamma\), and \(T\).
\end{lemma}
\begin{proof}
	See
	\cite[Chapter~IV, Section~2, especially estimate~(2.4),
	p.~274]{LSU}.
\end{proof}

The next lemma estimates \(\partial_{x_n}k\) when \(g(x',s)=\psi(x')\phi(s)\).

\begin{lemma}\label{lem-260234}
	Let \(T>0\) and \(\alpha\in(0,1)\).
	Suppose that \(\psi\in C_c^{1,\alpha}(\Sigma)\),
	\(\phi\in C^1([0,T])\), and \(\phi(0)=0\).
	Then, for every \(t\in(0,T]\),
	\(\partial_{x_n}k(\cdot,t)\) extends continuously to
	\(\overline{\R_+^n}\), and
	\begin{align}\label{0928-0111}
		\left\|\partial_{x_n}k(\cdot,t)\right\|_
		{C^\alpha(\overline{\R_+^n})}
		\leq C_T\|\psi\|_{C^{1,\alpha}(\Sigma)}
		\|\phi\|_{C^1([0,T])},
	\end{align}
	where \(C_T\) depends only on \(n,\alpha\), and \(T\).
\end{lemma}
\begin{proof}
See \cite[Chapter~IV, Sections~1--2]{LSU}.
\end{proof}

For the second normal derivative, the behavior as \(t\downarrow0\) depends
on the initial value of \(\phi'\).

\begin{lemma}\label{lem-260235}
	Let \(T>0\), \(\alpha,\beta\in(0,1)\), and
	\(\alpha_0=\min\{\alpha,2\beta\}\).
	Suppose that \(\psi\in C_c^{2,\alpha}(\Sigma)\),
	\(\phi\in C^{1,\beta}([0,T])\), and \(\phi(0)=0\).
	Then, for \(0<t\leq T\), \(\partial_{x_n}^2k(\cdot,t)\) extends continuously to \(\overline{\R_+^n}\), with
	\begin{align*}
		\partial_{x_n}^2k(x',0,t)
		=-\tfrac12\bigl[
		\psi(x')\phi'(t)-\phi(t)\Delta'\psi(x')
		\bigr].
	\end{align*}
	Moreover,
	\begin{align*}
		\|\partial_{x_n}^2k(\cdot,t)\|_
		{C^{\alpha_0}(\overline{\R^n_+})}
		&\leq
		C_T\|\psi\|_{C^{2,\alpha}(\Sigma)}
		\|\phi\|_{C^{1,\beta}([0,T])}
		+C t^{-\alpha_0/2}|\phi'(0)|\,\|\psi\|_\infty,
	\end{align*}
	where \(C_T\) depends only on \(n,\alpha,\beta,T\) and \(C\) only on
	\(n,\alpha,\beta\).
\end{lemma}
\begin{proof}
	See \cite[Chapter~IV, Sections~1--2]{LSU}.
\end{proof}

\section{Boundary limits for Newtonian and heat potentials}
\label{260264}

\begin{lemma}\label{lem-260265}
	Let \(n\geq2\), \(0<\alpha<1\), and \(1\leq i,j\leq n-1\).
Then the following hold.
\begin{enumerate}
\item If \(f\in C_c^\alpha(\Sigma)\), then \(R_i'f\in C^\alpha(\Sigma)\cap C_0(\Sigma)\) and
\begin{align*}
\|R_i'f\|_{C^\alpha(\Sigma)}
\leq C_{n,\alpha}
\bigl([f]_{C^\alpha(\Sigma)}
+\|f\|_{L^1(\Sigma)}\bigr).
\end{align*}
\item If \(f\in C_c^{1,\alpha}(\Sigma)\), then \(\partial_{x_j}R_i'f=R_i'(\partial_{x_j}f)\in C^\alpha(\Sigma)\cap C_0(\Sigma)\).
\item For \(f\in L^2(\Sigma)\) and \(x_n>0\), pointwise on \(\Sigma\),
\begin{align}\label{261031}
\partial_{y_i}N(\cdot,x_n)*_\Sigma f
=\partial_{y_n}N(\cdot,x_n)*_\Sigma R_i'f.
\end{align}
\item For every \(f\in C_c^\alpha(\Sigma)\) and \(x'\in\Sigma\),
		\begin{align}\label{260266}
			\lim_{x_n\to0^+}
			\int_{\Sigma}
			\partial_{y_i}N(y',x_n)f(x'-y')\,dy'
			=\frac12R_i'f(x').
		\end{align}

		\item
		For every \(f\in C_0(\Sigma)\) and \(x'\in\Sigma\),
		\begin{align}\label{260267}
			\lim_{x_n\to0^+}
			\int_{\Sigma}
			\partial_{y_n}N(y',x_n)f(x'-y')\,dy'
			=\frac12f(x').
		\end{align}
	\end{enumerate}
\end{lemma}
\begin{proof}
See \cite[Chapter~II, Section~6.9 and Chapter~III,
	Sections~2.1--2.3]{Stein1970}.
\end{proof}

\begin{lemma}\label{lem-260268}
	Let \(g\in C_c(\Sigma\times\R)\) satisfy
	\(g(x',s)=0\) for \(s\leq0\).
	Then, locally uniformly in \((x',t)\in\Sigma\times(0,\infty)\),
	\begin{align}\label{260269}
		\lim_{a\to0^+}
		\int_{-\infty}^t\int_\Sigma
		\partial_a\Gamma(x'-y',a,t-s)g(y',s)\,dy'ds
		=-\frac12g(x',t).
	\end{align}
\end{lemma}
\begin{proof}
See \cite[Chapter~IV, Sections~1--2]{LSU}.
\end{proof}

\begin{lemma}\label{lem-260270}
	Let \(n\geq2\), \(0<\alpha<1\), and \(1\leq j\leq n-1\).
	Suppose that \(\psi\) satisfies {\rm\textbf{(S)}} with H\"{o}lder
	exponent \(\alpha\) and that \(\phi\) satisfies
	{\rm\textbf{(T1)}}, and let \(h(y',s):=\psi(y')\phi'(s)-\phi(s)\Delta'\psi(y')\). Then
	\begin{align}\label{0713-30}
		M_t
		:=\sup_{0\leq s\leq t}\|h(\cdot,s)\|_{C^\alpha(\Sigma)}\leq C\|\psi\|_{C^{2,\alpha}(\Sigma)}\|\phi\|_{C^1([0,1])}<\infty.
	\end{align}
	For \(x_n>0\),
	\begin{align}\label{0713-31}
		\partial_{x_n}k(x,t)
		=\int_0^t\int_\Sigma
		\Gamma(x-y',t-s)h(y',s)\,dy'ds.
	\end{align}
	For every \(x'\in\Sigma\) and \(t>0\), the limit \(\partial_{x_n}k(x',0,t)
		:=\lim_{x_n\to0^+}\partial_{x_n}k(x,t)\) 
	exists, and
	\begin{align}\label{0713-28}
		\partial_{x_n}k(x',0,t)
		=\frac{1}{2\sqrt{\pi}}
		\int_0^t\frac{1}{\sqrt{t-s}}
		\int_\Sigma
		\Gamma'(x'-y',t-s)h(y',s)\,dy'ds.
	\end{align}
	Moreover,
	\begin{align}\label{0709-07}
		\partial_{x_j}\partial_{x_n}k(x',0,t)
		=\frac{1}{2\sqrt{\pi}}
		\int_0^t\tau^{-\frac12}
		\int_\Sigma
		\partial_{x_j}\Gamma'(x'-y',\tau)
		h(y',t-\tau)\,dy'\,d\tau,
	\end{align}
	and, for \(0<\tau<t\),
	\begin{align}\label{0713-25}
		\left|
		\int_\Sigma
		\partial_{x_j}\Gamma'(x'-y',\tau)
		h(y',t-\tau)\,dy'
		\right|
		\leq C M_t\tau^{\frac{\alpha-1}{2}}.
	\end{align}
	In particular,
	\begin{align}\label{0713-29}
		\partial_{x_n}k(\cdot,0,t)\in C_0(\Sigma),
		\qquad
		\partial_{x_j}\partial_{x_n}k(\cdot,0,t)\in C_0(\Sigma).
	\end{align}
Finally, for \(a>0\),
	\begin{align*}
		\|\partial_{x_j}\partial_{x_n}k(\cdot,a,t)
		-\partial_{x_j}\partial_{x_n}k(\cdot,0,t)\|_\infty
		\leq C_{n,\alpha}M_ta^\alpha.
	\end{align*}
\end{lemma}
\begin{proof}
See \cite[Chapter~IV, Sections~1--2]{LSU}.
\end{proof}

\section*{Statements and Declarations}

\textbf{Funding.} The work of both authors was supported by the National Research Foundation of Korea (NRF), grant nos. RS-2024-00336346 and RS-2024-00406821.

\textbf{Conflict of interests.} The authors declare that they have no competing interests.

\textbf{Data availability.} Data sharing is not applicable to this article as no datasets were generated or analyzed during the current study.

\textbf{Use of Artificial Intelligence.} During the preparation of this manuscript, the authors used ChatGPT and Claude solely for language editing and polishing of the manuscript. The authors take full responsibility for the correctness and validity of all mathematical results presented in this paper.
\begingroup
\setlength{\itemsep}{0pt plus 1pt}
\setlength{\parskip}{0pt}

\endgroup

\end{document}